\documentclass[letterpaper,11pt,oneside,reqno]{amsbook}
\usepackage{bm}
\usepackage{mathrsfs}
\usepackage{wrapfig}
\usepackage{subcaption}
\usepackage{amsfonts,amsmath, amssymb,amsthm,amscd,stmaryrd,mathtools}
\usepackage[height=9.6in,width=5.95in]{geometry}
\usepackage[mpinclude=false,DIV=13]{typearea}
\usepackage{verbatim}
\usepackage{hyperref}
\newcommand{\hfff}[2]{\label{gl:#1}\hypertarget{gl:#1}{#2}}
\newcommand{\fff}[1]{&\hyperlink{gl:#1}{\pageref*{gl:#1}}\\}
\usepackage{graphicx}
\usepackage{float}

\graphicspath{{figures/}}
\usepackage{latexsym}
\usepackage{lscape}
\usepackage{epsfig}

 \usepackage[utf8]{inputenc}
\usepackage{enumitem}

\usepackage{tikz}
\usetikzlibrary{decorations.pathreplacing,arrows.meta}

\usepackage{setspace}

\usepackage[skip=0.3\baselineskip plus 1pt]{parskip}
\usepackage{needspace}
\makeatletter
\def\paragraph{\@startsection{paragraph}{4}%
  \normalparindent\z@{-\fontdimen2\font}%
  {\normalfont\itshape}}
\def\subsubsection{\@startsection{subsubsection}{3}%
  \normalparindent{.25\linespacing\@plus.25\linespacing}{-.5em}%
  {\normalfont\itshape}}
\makeatother
\newcommand{\e}[0]{\epsilon}
\newcommand{\Capt}{\operatorname{Capt}}
\newcommand{\cF}{\mathcal{F}}

\newcommand{\Accept}{\operatorname{Accept}}

\newcommand{\midpiece}{\mathrm{MidPiece}}
\newcommand{\midconnect}{\mathrm{MidConnect}}

\newcommand{\cS}{\mathcal{S}}

\newcommand{\cost}{\operatorname{cost}}

\newcommand{\Halo}{\operatorname{Halo}}

\newcommand{\Seal}{\operatorname{Seal}}

\newcommand{\avf}{\mathfrak a_{\rm vf}}
\newcommand{\af}{\mathfrak a_{\rm f}}
\newcommand{\Cfine}{\mathsf{C}_{\rm f}}
\newcommand{\Bfine}{B_{\rm f}}
\newcommand{\qbarrow}{\overline{q}_{\rm row}}

\newcommand{\prot}{\operatorname{ph}}

\newcommand{\Obox}[1]{\widehat B[#1]}
\newcommand{\OBoxLat}{\widehat{\mathscr B}_K}

\newcommand{\ivf}{\mathrm{IVF}}
\newcommand{\llf}{\mathrm{LLF}}

\newcommand{\Bbox}[1]{B[#1]}
\newcommand{\Vtx}{\mathrm V}
\newcommand{\Edg}{\mathrm E}
\newcommand{\BoxUnion}[1]{[#1]_K}

\newcommand{\one}{\mathbf{1}}

\newcommand{\cH}{\mathcal{H}}

\newcommand{\EE}{\ensuremath{\mathbb{E}}}
\newcommand{\PP}{\ensuremath{\mathbb{P}}}
\newcommand{\N}{\ensuremath{\mathbb{N}}}
\newcommand{\R}{\ensuremath{\mathbb{R}}}
\newcommand{\Rlaw}{\ensuremath{\mathsf{R}}}
\newcommand{\Z}{\ensuremath{\mathbb{Z}}}
\newcommand{\Q}{\ensuremath{\mathbb{Q}}}
\newcommand{\Qdyn}{\ensuremath{\mathbb{Q}^{\rightarrow}}}
\newcommand{\E}[0]{\mathbb{E}}

\newtheorem{theorem}{Theorem}[chapter]

\newtheorem{lemma}[theorem]{Lemma}
\newtheorem{proposition}[theorem]{Proposition}
\newtheorem{corollary}[theorem]{Corollary}

\newtheorem{remark}[theorem]{Remark}

\numberwithin{section}{chapter}

\theoremstyle{definition}

\theoremstyle{definition}

\theoremstyle{definition}
\newtheorem{definition}[theorem]{Definition}

\theoremstyle{definition}

\theoremstyle{definition}

\newcommand{\intint}[1]{\llbracket #1 \rrbracket}

\newcommand{\Boxsupp}{\operatorname{supp}_K}
\newcommand{\LayerBoxsupp}{\operatorname{supp}_{K,N}}

\newcommand{\mc}{\mathcal}

\newcommand{\downspace}{\mathbb{D}}
\newcommand{\upspace}{\mathbb{U}}

\newcommand{\dbt}{\mathrm{DBT}}

\newcommand{\base}{\mathrm{base}}
\newcommand{\head}{\mathrm{head}}
\newcommand{\drift}{\mathfrak d}
\newcommand{\upwell}{\mathrm{UpClust}}
\newcommand{\downwell}{\mathrm{DownClust}}

\newcommand{\downcluster}{\mathrm{DownClust}}
\newcommand{\upcluster}{\mathrm{UpClust}}

\newcommand{\hmax}{h_{\rm max}}

\newcommand{\hadj}{h_{\rm adj}}

\newcommand{\Bed}{\mathrm{Bed}}
\newcommand{\Top}{\mathrm{Top}}

\newcommand{\fine}{\mathsf{Fine}}
\newcommand{\success}{\mathsf{SimpleJoin}}

\newcommand{\veryfine}{\mathsf{VeryFine}}

\newcommand{\sting}{\mathsf{Str}}
\newcommand{\culdesac}{\mathsf{CulDeSac}}
\newcommand{\culdesacset}{\mathrm{CdS}}
\newcommand{\cds}{\mathsf{CulDeSac}}
\newcommand{\cdsset}{\mathrm{CdS}}
\newcommand{\icds}{\mathsf{InvCulDeSac}}
\newcommand{\icdsset}{\mathrm{ICdS}}
\newcommand{\stingset}{S}
\newcommand{\cluster}{\mathsf{Cluster}}

\newcommand{\clusteruthree}{\mathsf{Cluster}_{[u,u+3r]}}

\newcommand{\cset}{\mc{C}}
\newcommand{\csetuthree}{\mc{C}_{u,u+3r}}
\newcommand{\overlinecsetuthree}{{\overline{\cset}}_{u,u+3r}}
\newcommand{\ucs}{\mathrm{UCS}}

\newcommand{\boxbig}{{\rm BigBox}}

\newcommand{\const}{K}
\newcommand{\surfacek}{\mathrm{Surface}_{\const}}

\newcommand{\slide}{{\rm Slide}}
\newcommand{\iter}{{\rm IterateSlide}}

\newcommand{\slab}{{\rm Slab}}

\newcommand{\fconn}{\overset{\scriptscriptstyle\rm fin}{\longleftrightarrow}}
\newcommand{\fconnext}{\overset{\scriptscriptstyle\rm ext}{\longleftrightarrow}}
\newcommand{\sconn}{\overset{\scriptscriptstyle\rm str}{\longleftrightarrow}}

\newcommand{\sausageset}{\mathscr{S}}

\newcommand{\badcoll}{\mathsf{Bad}}

\newcommand{\kappasting}{\kappa_{\rm string}}
\newcommand{\kappacluster}{\kappa_{\rm clust}}

\newcommand{\fv}{J_v}

\newcommand{\rlevel}{{\rm RenLev}}

\newcommand{\vecell}{{\pmb{\ell}}}
\newcommand{\cvl}{\cdot \vecell}
\newcommand{\cvlzero}{\cdot\vecell_0}
\newcommand{\uclust}{[\mathfrak{u}]}
\newcommand{\dclust}{\mathfrak{d}}

\newcommand{\simvecell}{\stackrel{\vecell}{\sim}}
\newcommand{\approxvecell}{\stackrel{\vecell}{\approx}}
\newcommand{\closedint}[1]{[#1]}

\newcommand{\fubuthree}{\{z\in\Z^d:z\cvl\in[u,u+3r]\}}

\newcommand{\dproj}{d_{{\rm proj}}}

\newcommand{\vertbdry}{\mathrm{VertBdry}}
\newcommand{\edgebdry}{\mathrm{EdgeBdry}}

\newcommand{\rmax}{r_{\rm max}}
\newcommand{\rmin}{r_{\rm min}}
\newcommand{\hgt}{\operatorname{hgt}}
\newcommand{\cover}{{\rm Cover}}
\newcommand{\lat}{\mathrm{Lat}}

\newcommand{\nkm}{\overline{N}_{k,k+m}}

\newcommand{\Nzero}{\mathbb N_0}

\newcommand{\diam}{{\rm diam}}
 \newcommand{\Ext}{\operatorname{Ext}}

\makeatletter
\newcommand\xxleftrightarrow[2][]{%
  \ext@arrow 9999{\longleftrightarrowfill@}{#1}{#2}}
\newcommand\longleftrightarrowfill@{%
  \arrowfill@\leftarrow\relbar\rightarrow}
\makeatother

\newcommand{\cV}{\mathcal{V}}
\newcommand{\notconn}{\mathrel{\longleftrightarrow\llap{\hspace{2.4pt}/\hspace{11.4pt}}}}

\newcommand{\Dmesh}{2K}
\newcommand{\Mesh}{\mathscr M}
\newcommand{\MeshBox}[1]{\mathsf D_{#1}}
\newcommand{\MCap}{\operatorname{Capt}_{2K}}
\newcommand{\MExt}{\operatorname{Ext}_{2K}}
\newcommand{\MeshSupp}{\operatorname{supp}_{\Mesh}}
\newcommand{\MeshInt}{\operatorname{Int}_{\Mesh}}
\newcommand{\ECap}{\mathcal E_{\rm cap}}
\newcommand{\RetRec}{\mathfrak H}
\newcommand{\MidRec}{\mathfrak G}
\newcommand{\ProsSeal}{\mathscr P}
\newcommand{\ProsExc}{\mathsf E^{\rm pros}}
\newcommand{\Reach}{\mathsf R}
\newcommand{\InitExc}{\mathsf E^{\rm init}}
\newcommand{\BdCost}{\cost^{\rm bd}}
\newcommand{\GainCost}{\cost^{\rm gain}}
\newcommand{\PrepCost}{\cost^{\rm prep}}
\newcommand{\SealCost}{\cost^{\rm seal}}
\newcommand{\SlideCost}{\cost^{\rm slide}}
\newcommand{\PassPrep}{\mathsf P}
\newcommand{\PassExc}{\mathsf Q}

\begin{document}

\title{Mass gap for finite connections \\ in supercritical percolation}

\author{Alexander Fribergh and Alan Hammond}
\AtEndDocument{
  \bigskip
  \small
  \par
  \textsc{Alexander Fribergh} \\
  \textsc{D\'epartement de math\'ematiques et de statistique} \\
  \textsc{Universit\'e de Montr\'eal} \\
  \textit{Email:} \texttt{alexander.fribergh@umontreal.ca} \\
}

\AtEndDocument{
  \bigskip
  \small
  \par
  \textsc{Alan Hammond} \\
  \textsc{Departments of Mathematics and Statistics, U.C. Berkeley} \\
  \textsc{899 Evans Hall, Berkeley, CA, 94720-3840, U.S.A.} \\
  \textit{Email:} \texttt{alanmh@berkeley.edu} \\
}

\begin{abstract}

We consider finite open clusters based at the origin in Bernoulli bond
percolation on~$\Z^d$ for $d\geq3$ at a supercritical value $p>p_c(d)$, whose
maximal height in the direction of a unit vector~$\vecell$ lies in a fixed-width
window at distance~$u$.  Conditioned on this event, we prove that, except with
probability exponentially small in~$u$, such a cluster is populated at positive
density by renewal points.  Here, a renewal point~$w$ is the only vertex of the
cluster in the fixed-width band immediately below the height hyperplane
through~$w$.

The main technical estimate is that, under this conditioning, the absence of a renewal point is exponentially unlikely.  From it we derive an exponential mass gap for the irreducible pieces delimited by successive renewal points, and from the mass gap the stated renewal abundance. 
 The mass gap, together with the renewal decomposition and endpoint estimates established here, yields a sharp point-to-point Ornstein--Zernike asymptotic for finite connectivity, with its Gaussian-type power-law correction. The multiplicative prefactor in this asymptotic is positive and real-analytic, and the associated equidecay and Wulff boundaries are real-analytic and strictly convex.

The principal  innovation of the article is the use of a slide
resampling that preserves the conditioned law.  The slide acts by cutting the
long finite connection at a given $\vecell$-height and relocating the upper
piece by a random shift.  The mass gap is obtained by a delicate analysis of
an iterated version of the slide operation along a positive density of
attempted heights.
\end{abstract}

\maketitle

\vspace{-2mm}

\setcounter{tocdepth}{1}
\tableofcontents

\chapter{Introduction}
\label{s.introduction}

Ornstein and Zernike~\cite{OrnsteinZernike} presented a theory that offers a sharp asymptotic description of the density correlations in classical fluids by assuming the validity of a form of locality expressed via a renewal equation.
Finding rigorous manifestations of their theory amounts to giving probabilistic descriptions of objects such as long clusters in terms of essentially one-dimensional chains of irreducible components.
In subcritical models, the formation of lengthy connections is an exponential rarity, with conditioned connections formed by the concatenation of irreducible blocks
at successive bottlenecks or  {\em renewal} points.
 Self-avoiding walk with an exponential penalty for length is such a model. It was addressed for axial directions in~\cite{ChayesChayes}
and it offered the setting for an early non-perturbative analysis in~\cite{Ioffe}.
 Campanino and Ioffe obtained sharp asymptotic descriptions for subcritical
percolation on Euclidean lattices~\cite{CampaninoIoffe}, while Campanino,
Ioffe and Velenik developed corresponding theories for finite-range Ising
models above the critical temperature~\cite{CampaninoIoffeVelenik03} and the
subcritical random-cluster model~\cite{CampaninoIoffeVelenik}.

Taken together, these works provide an Ornstein--Zernike framework: irreducible decomposition reduces a long connection to an effectively one-dimensional chain, and quantitative control of that chain permits a multivariate local-limit analysis producing the factor \(|x|^{-(d-1)/2}\) alongside the exponential decay.
  The analysis also
requires control of the
vector increments and endpoint pieces, identification of the relevant
convergence domain, and sufficient uniformity in direction to deduce
regularity of the equidecay surface and its polar Wulff shape.

Chief among the framework's hypotheses is a \emph{mass gap} condition.
Irreducible blocks concatenate to form long conditioned connections in a given height coordinate. In a subcritical model, the law of the block rooted at the origin is a sub-probability measure. Tilting the law exponentially in the height with a  suitable rate  yields a probability measure~$\mu$, permitting the concatenation of independently sampled blocks into an infinite random string. The mass gap condition asserts that the height of the renewal block under the law~$\mu$ has an exponential tail.
In transfer-operator formulations, the analogous conclusion is the isolation
of the leading eigenvalue by a positive spectral gap.  In the present exact
renewal formulation, the operative statement is the exponential tail of one
irreducible increment.

The papers~\cite{CampaninoIoffe,CampaninoIoffeVelenik03,CampaninoIoffeVelenik} verify mass gap in the subcritical models they investigate  by employing a ball-and-point search algorithm.
  It marks points on a lengthy connection and
tracks where the connection leaves a ball of fixed radius, in a suitable norm,
centred at the marked point.  The resulting coarse-grained description is
used both in the mass-gap analysis and in the concatenation arguments that
lead from renewal structure to the geometric Ornstein--Zernike conclusions. In these analyses, the rare object, a long subcritical connection, is being detected and explored by a local procedure.

By contrast, for bond percolation on the nearest neighbour Euclidean lattice~$(\Z^d,\sim)$,
a pair of vertices is {\em finitely connected} if the two vertices lie in a common
finite open connected component.
The open connected component of a given vertex may be iteratively explored at
any parameter value, and opening further edges can only extend such an
exploration towards a distant point.  This favourable monotonicity of ordinary
connection events is absent when a finite connection is sought
supercritically.  In this phase, the search is conducted
under the continual exertion of pressure from a pair of opposing forces, each threatening to thwart the search for a finite connection.
If the live end of the exploration contains  merely a few points, it is liable to terminate imminently, doing so if the unrealised incident edges turn out to be closed;
 while if this live cardinality is high,  the next few search steps risk contact with the infinite open cluster that pervades space.

For finite connectivity in the highly supercritical phase, Campanino and
Gianfelice~\cite{CampaninoGianfelice} obtained Ornstein--Zernike conclusions
in every asymptotic direction, generalising earlier axial results of Braga,
Procacci and Sanchis~\cite{BragaProcacciSanchis}.
 In this regime, the volume of a long finite connection  very typically remains within a small multiple of its   ground-state value, so that
such a connection is populated to high density by renewal points. Similarly to
early rigorous progress that was made with perturbative arguments, such as~\cite{MinlosZhizhina} for the Ising model at very high temperature,  the needed mass gap condition follows relatively directly, though some care is needed in specifying the notion of renewal in general non-axial directions.

In this article, we establish the renewal mass gap for finite connections in nearest-neighbour bond percolation on~$(\Z^d,\sim)$ for $d \geq 3$, throughout the supercritical phase and in every direction $\vecell \in S^{d-1}$.  The exponential decay of the probability of a long finite connection was proved in~\cite{FriberghHammond}.  We prove that, when a percolation exploration in the supercritical phase reaches a distant height without contact with the infinite open cluster, the resulting finite cluster withstands the stated dual pressures---imminent perishing when small; escape along an infinite open path when large---by adopting the random geometry characteristic of the Ornstein--Zernike picture.  Except with exponentially small probability under the point-to-hyperplane conditioning, this connection contains a positive density of renewal points.

At these renewal points the long connection breaks into successive irreducible pieces, or \emph{sausages}.  Exponential tilting of the sausage occurrence weights by the connection decay rate normalises them to a probability law~$\mu$, under which independent sausages may be concatenated into a renewal chain.  The mass-gap theorem asserts that the height of a sausage under~$\mu$ has an exponential tail.  Thus a second exponential decay appears at the level of the irreducible building blocks: beyond the exponential rarity of the long finite connection itself, the scale of the pieces in its exact renewal decomposition is exponentially controlled.  Together with endpoint estimates and vector renewal analysis, this structure yields a sharp point-to-point Ornstein--Zernike asymptotic for finite connectivity in every asymptotic direction, as well as a positive real-analytic amplitude and strict convexity and real analyticity of the equidecay surface and its polar Wulff shape.

The absence of monotonicity makes the mass gap problem a substantial challenge.
Our main proof technique is an iterated slide resampling of the conditioned cluster.
Its application shows that the competing pressures of imminent perishing and escape to
infinity are typically resolved through a positive density of near-extinction bottlenecks,
or {\em renewal points}.  Obtaining mass gap in this way is the principal technical task
we face, and we will present renewal setup and the slide in overview in the two chapters following the introduction. But the absence of monotonicity will be felt even after mass gap is derived, because it prevents
a direct application of the subcritical Ornstein--Zernike theory.  The renewal asymptotic in a
fixed height direction does follow from the mass gap much as it does
subcritically.  For the geometric Ornstein--Zernike conclusions, we prove an exact vector
renewal representation, exponential moments for the pair of {\em culs-de-sac} found at the two ends of the long finite connection, and
a marked decomposition of the finite-connectivity susceptibility.  The resulting
independent-increment local-renewal step is the IID analogue of the
transfer-operator local-limit calculation of Campanino--Ioffe--Velenik; here it is
proved directly from the lattice local central limit theorem and a Cram\'er bound.

We now fix notation and state the main results.  We first introduce renewal points, sausages and the tilted sausage law, and state the mass-gap theorem and its renewal-abundance consequence.  We then give the fixed-direction point-to-plane asymptotic and the point-to-point Ornstein--Zernike conclusions.

\section{Setting, renewal structure and the mass gap}
\label{ss.renewal-intro}

Let $d \geq 3$ and $p \in (p_c(d),1)$.  Write \hfff{percolation-law}{$\PP$} for the law of
supercritical bond percolation on $(\Z^d,\sim)$ with parameter~$p$.
The \emph{open clusters} are the connected components of the open-edge
subgraph.  Under $\PP$ there is almost surely a unique infinite open cluster;
an open cluster with finite vertex set is called a \emph{finite cluster}.
For $x\in\Z^d$, write $\mc C(x)$ for the open cluster containing~$x$.
Write $x \fconn y$ for the event that $x$ and $y$ belong to the same finite
cluster.

Let $|\cdot|$ denote the Euclidean norm on $\R^d$ and write
$S^{d-1} = \{x \in \R^d : |x| = 1\}$ for the unit sphere.
Fix a unit vector \hfff{height-axis}{$\vecell \in S^{d-1}$} and call $x\cvl$
the \emph{height} of $x \in \Z^d$.  Put
$r=r(\vecell)=\max_{1\le i\le d}|e_i\cvl|$, the greatest height spanned by a
single lattice edge.

The direction $\vecell$ is \emph{lattice} if the $\Z$-linear span of
$\{e_i\cvl : 1\le i\le d\}$ is discrete in $\R$, hence of the form
$s\Z$ for some $s>0$, the \emph{span}; otherwise it is
\emph{non-lattice}.  The
\emph{cover} $\cover(\vecell)$ is $s\N$ in the lattice case and
$(0,\infty)$ in the non-lattice case.  Write $f \approxvecell g$ if
$|f/g - 1| \leq Ce^{-cu}$ as $u\to\infty$ along the cover, for some
$C, c > 0$, and $f \simvecell g$ if $f/g \to 1$ as $u\to\infty$ along
the cover.

The \hfff{base-head-cluster-event}{\emph{base}} of a finite cluster is its lexicographically minimal
vertex of minimum height, and its \emph{head} is its lexicographically
maximal vertex of maximum height, where the lexicographic order on each
height hyperplane is induced by a fixed orthonormal basis for
$\{\vecell\}^\perp$.  Write $\cluster_u = \cluster_u(0)$ for the event that
the origin is the base of a finite cluster whose maximum height equals $u$.
More generally, for an interval $I\subseteq[0,\infty)$, put
$\cluster_I=\bigcup_{u'\in I}\cluster_{u'}(0)$.  Thus $\cluster_{[u,\infty)}$
is the event that this maximum height is at least $u$.

The asymptotic results arise from a one-dimensional renewal decomposition of long finite connections.  A vertex $x$ of a finite connected
subgraph $U$ is a \emph{renewal point} if every $y \in V(U)$ satisfying
$y\cvl \in (x\cvl - r, x\cvl]$ equals $x$; that is, $x$ is the only vertex of
$U$ in the band of width $r$ topped by its own level.  Since no edge spans a
height exceeding $r$, this isolation separates the part of $U$ below the band
from the part above~$x$; when both parts are non-empty, $x$ is a cut vertex.
Such a vertex is a point of maximal narrowness along the connected cluster:
the lower and upper parts meet only there, providing the geometric separation
underlying the renewal decomposition.  The set of renewal points of $U$ is written $R(U)$.

A \emph{string} from $x$ to $y$ (with $x\cvl<y\cvl$) is a finite open connected component of the slab between heights $x\cvl$ and $y\cvl$ for which $x$ and $y$ are renewal points; it decomposes uniquely into \emph{sausages}, the substrings whose only renewal points are their two endpoints.  These sausages are the irreducible renewal blocks of the Ornstein--Zernike picture.  Write $\sting_u=\sting_u(0)$ for the event that a string of height $u$ is based at the origin, and use interval subscripts for unions over terminal heights.

The exponential decay rate
$\zeta = \zeta(p,d,\vecell) \in (0,\infty)$ is specified by the string limit
\[
  \lim_{u\to\infty}u^{-1}\log\PP\bigl(\sting_{[u,u+r)}\bigr)=-\zeta
\]
established in Proposition~\ref{p.zeta}, after~\cite{FriberghHammond}.  The
point-to-plane asymptotics below show that the same rate governs the corresponding cluster events.
The \emph{sausage measure}~$\mu$ is the probability measure on sausages based at the origin given by
\[
  \mu(S)\, = \, e^{\zeta\hgt(S)}\,\PP(\mathsf{O}(S))\, ,
\]
where $\hgt(S)$ is the height of~$S$.  Here
$\mathsf{O}(S)$ is the event that $S$ occurs---every edge of $S$ is
open, while neighbouring edges are closed.
(A boundary convention partitions edges, so that consecutive sausages prescribe the status of disjoint edge sets.)
Proposition~\ref{p.muprob} says precisely that tilting by the decay rate $\zeta$ normalises the total sausage weight to one.  Concatenating independent $\mu$-sausages yields the renewal process $\mu^{\otimes\N}$, whose renewal levels $\rlevel$ obey, as Lemma~\ref{l.essn} and the sausage decomposition give,
\[
  e^{\zeta u}\,\PP(\sting_u)\, = \, \mu^{\otimes\N}(u \in \rlevel)\, .
\]

Our central structural conclusion is that this exact renewal law has an exponential moment in the sausage height.

\begin{theorem}[Mass gap]
\label{t.mass-gap}
There exists $\theta > 0$ such that $\E_\mu[e^{\theta\hgt(S)}] < \infty$.
\end{theorem}

In transfer-operator terms, the mass gap plays the role of a positive spectral
gap; in the present renewal formulation it is the analyticity of the
moment-generating function of $\hgt(S)$ in a complex neighbourhood of the origin.
The underlying slide estimate, Theorem~\ref{t.intermediate-renewal}(2), is that a long conditioned cluster contains an interior renewal point except with exponentially small probability.  It yields the mass gap above.  The resulting string renewal asymptotic, together with a Cram\'er estimate and the endpoint comparison, then returns this exponential control to the conditioned finite cluster in the form of linear renewal abundance.

For this statement, recall that $\cluster_{[u,u+3r]}$ is the event that the origin is the base of a finite cluster with maximum height in $[u,u+3r]$.  The fixed window of width $3r$ accommodates the general-direction resampling and allows the principal renewal estimates to be stated under one conditioning in both the lattice and non-lattice cases.

\begin{theorem}[Renewal abundance theorem]
\label{t.renewal}
There exist constants $c_0, C, c > 0$ such that for all $u$ large,
\[
  \PP\bigl(|R(\mc C(0))|\, \geq \, c_0 u
     \;\bigm|\; \clusteruthree\bigr)
\, \geq \, 1 - Ce^{-cu}\, .
\]
\end{theorem}

The two theorems record the renewal structure at different levels: the mass gap controls the height of an irreducible increment under the exact IID law~$\mu$, while renewal abundance says that, in a long conditioned finite cluster, such renewal cuts occur at linear density except with exponentially small probability.

\section{Point-to-plane asymptotics}

The renewal description also determines the sharp probability that a finite cluster reaches a distant height in the fixed direction~$\vecell$.

\begin{theorem}[Sharp point-to-plane asymptotics: lattice and non-lattice]
\label{t.stingwelltrap.cluster}
Let $d \geq 3$, $p \in (p_c,1)$, and $\vecell \in S^{d-1}$.  There exists
$\kappacluster = \kappacluster(d,p,\vecell) \in (0,\infty)$ such that the
following hold, the two cases differing according to whether the height
$z\cvl$ takes values in a lattice as $z$ ranges over $\Z^d$.
\begin{enumerate}
\item \textup{(Lattice case.)}  If $\vecell$ is lattice, with span $s$,
then for $u\in s\N$
\[
  \PP\bigl(\cluster_u\bigr) \;\approxvecell\; \kappacluster e^{-\zeta u}\, .
\]
\item \textup{(Non-lattice case.)}  If $\vecell$ is non-lattice, then
\[
  \PP\bigl(\cluster_{[u,\infty)}\bigr) \;\simvecell\; \kappacluster e^{-\zeta u}\, .
\]
\end{enumerate}
\end{theorem}

The probability that the cumulative $\mu$-sausage heights hit $u$ is the renewal mass at $u$.  In the lattice case these heights lie in $s\Z$, and the lattice renewal theorem gives convergence of this mass to a positive constant.  The mass gap upgrades this to an exponential rate of convergence; after the endpoint contributions are included, this gives the strong asymptotic $\approxvecell$ and determines $\kappacluster$.  In the non-lattice case no individual height has comparable mass, so the non-lattice renewal theorem is applied to the cumulative event $\cluster_{[u,\infty)}$.

\section{Ornstein--Zernike consequences}
\label{ss.oz-intro}

Theorem~\ref{t.stingwelltrap.cluster} gives a point-to-hyperplane asymptotic
in a fixed direction.  The sharper question is the asymptotic of the
two-point finite connection probability itself.  Because finite connection
is not a monotone event, the usual positive-association argument does not
furnish the supermultiplicative estimate needed to establish a concatenation
limit.  We therefore define the inverse correlation length through the
convergence domain.  Write

$$
 \tau^{\rm f}(x)\, = \, \PP(0\fconn x),\qquad
 \mathcal D\, = \, \left\{t:\sum_{x\in\Z^d}e^{t\cdot x}\tau^{\rm f}(x)<\infty\right\}\, ,
 \qquad \mathcal W\, = \, \overline{\mathcal D}\, ,
$$
and define
\begin{equation}
\label{e.xi-def}
\xi(x) \, = \, \sup_{t\in\mathcal W}t\cdot x,\qquad
K_\xi \, = \, \big\{ x\in\R^d:\xi(x)  \le  1 \big\} \,  .
\end{equation}
Thus \hfff{xi-bodies}{$\xi$} is homogeneous and convex by construction,
$K_\xi$ is the equidecay body, and $\mathcal W$ is its polar Wulff body.
The renewal analysis later identifies $\xi$ with the directional exponential
rate and establishes the regularity asserted below.

\begin{theorem}[Ornstein--Zernike asymptotic]
\label{t.oz-intro}
Let \(d\ge3\) and \(p\in(p_c,1)\).  There is a positive real-analytic function
\(A_{\rm f}\) on \(S^{d-1}\) such that, with \(\hat x=x/|x|\),

$$
  \PP\bigl(0\fconn x\bigr)
\, = \,
  A_{\rm f}(\hat x)\,|x|^{-(d-1)/2}\,e^{-\xi(x)}\,\bigl(1+o(1)\bigr)
  \qquad (|x|\to\infty),
$$

where the error tends to zero uniformly in the direction~$\hat x$.
The boundaries of the equidecay body $K_\xi$ and its polar $\mathcal W$ are
real analytic and strictly convex.
\end{theorem}
Under point-to-hyperplane conditioning, a finite cluster reaches a large height in a prescribed direction~\(\vecell\).  Write \hfff{drift}{\(\drift(\vecell)\)} for the cluster's typical direction.  Corollary~\ref{c.conditioned-head-drift} makes this interpretation precise by giving exponential concentration of the cluster head direction.  We show that
$\vecell\longmapsto\drift(\vecell)$
is a real-analytic diffeomorphism of \(S^{d-1}\).  Thus every point-to-point direction is obtained from a unique height direction, and the local asymptotics furnished by the corresponding height-directed analyses patch together to yield Theorem~\ref{t.oz-intro}.

We also present, in Theorem~\ref{t.oz-rooted}, the accompanying sharp
asymptotic when the two prescribed vertices are further required to be the
extremal vertices of their common finite cluster in a natural height
coordinate.

\section{Related work and further directions}
\label{ss.related}

Every edge in the external edge boundary of a finite open cluster is closed,
and the corresponding dual cells form an enclosing surface.  In two
dimensions this surface is a dual circuit.  On choosing two marked dual
vertices, the circuit divides into two disjoint connections in subcritical
dual percolation.  This reduction brings the planar problem within the
Ornstein--Zernike theory of subcritical bond percolation developed
in~\cite{CampaninoIoffe}.  Campanino, Ioffe and
Louidor~\cite{CampaninoIoffeLouidor} used it to obtain sharp asymptotics for
finite point-to-point connections in supercritical bond percolation on $\Z^2$.
In the case $d\ge3$ that concerns us, the enclosing surface is a genuinely
$(d-1)$-dimensional collection of dual plaquettes rather than a
one-dimensional object joining two marked points, and, absent such a reduction, the opposing pressures of extinction and escape must be reconciled.
Fribergh and
Hammond~\cite{FriberghHammond} established the exponential decay
rate $\zeta$ for long finite connections in this setting.  That first-order
rate is an input to the renewal analysis and sharp asymptotics developed here.

For subcritical percolation the renewal structure is
close to that of an independent chain~\cite{CampaninoIoffe}, whereas for
finite-range Ising models above the critical temperature and for the
subcritical random-cluster model the analysis is organised by a Ruelle transfer
operator over a countable
alphabet~\cite{CampaninoIoffeVelenik03,CampaninoIoffeVelenik}.  The problem we
treat here is the supercritical counterpart of these subcritical theories---in the random-cluster correspondence, the low-temperature phase of the Ising
model---in which the rare object is not a long connection but a long
connection sealed off from the infinite cluster.

This enclosing-surface perspective also underlies the large-deviation theory
of finite clusters.  Pisztora~\cite{Pisztora96} developed a coarse-graining
theory for supercritical random-cluster models in $d\ge3$, and
Deuschel--Pisztora~\cite{DeuschelPisztora1996} used such methods at high density
to prove surface-order large deviations for the volume of the largest cluster
and of the largest region that it surrounds.  In three dimensions,
Cerf~\cite{CerfLD3D} proved the corresponding Wulff crystal theorem: conditioned
to be finite and to have volume at least $n^3$, a suitably coarse-grained and
rescaled cluster concentrates, up to translation, around a dilation of the
Wulff crystal, the fixed-volume minimiser of the anisotropic surface energy of
the enclosing dual-plaquette interface.

The finite-connection Wulff body $\mathcal W$, defined above as the polar of
$K_\xi$, comes instead from the inverse correlation length of a long thin
point-to-point connection and is encoded by the effectively one-dimensional
renewal chain $\mu^{\otimes\N}$.  In two dimensions planar duality makes the
relation to the surface-tension Wulff crystal exact up to scale:
Campanino--Ioffe--Louidor~\cite{CampaninoIoffeLouidor} identify the
finite-connection rate with twice the inverse correlation length of the dual
subcritical model, equivalently twice the corresponding surface tension,
reflecting the two disjoint dual connections that form the enclosing circuit.
Thus $\mathcal W$ is a dilation of the surface-tension Wulff crystal, and the
two coincide after normalisation to fixed volume.  No analogous identification
is presently available in $d\ge3$.

Renewal is also a basic theme in many polymer-type models, including
self-avoiding walk.  For example, the proof of sub-ballisticity for
self-avoiding walk in $d\ge2$~\cite{DuminilCopinHammond} proceeds by
contradiction from a ballisticity hypothesis and proves that walks with a
positive density of renewal points occur at subexponential cost.

The iterated slide resampling with which we will locate renewal points in the conditioned cluster is designed as a law-preserving dynamic that uses retained data to diagnose geometric obstructions to renewal, and fresh resampling of discarded data to alter the offending geometry.
  It would be interesting to develop the slide technique in settings such as self-avoiding walk.
  It is also natural to ask how the
method might operate as $p\searrow p_c$, and what it might say about finite
connections in this regime.
The near-critical picture is far better understood in
two dimensions.  Garban, Pete and Schramm~\cite{GPS18} constructed the scaling
limit of near-critical (and dynamical) percolation on the triangular lattice, and
Duminil-Copin~\cite{DCWulff13} used it to show that the normalised subcritical
Wulff crystal tends to a disk as $p\nearrow p_c$.  On the triangular lattice,
planar duality transfers this conclusion to the normalised supercritical
finite-cluster crystal as $p\searrow p_c$; the scalar factor in the
finite-connection rate does not affect the normalised shape.  Near and at
criticality in $d\ge3$ it is the behaviour of these
enclosing dual surfaces that is at stake; Duminil-Copin, Kozma and
Tassion~\cite{DCKT19} proved upper bounds on the correlation length near
criticality in all dimensions, and Cerf~\cite{CerfTwoArm15} obtained lower
bounds on the two-arm exponent at $p_c$.  How the finite-connection Wulff shape $\mathcal W$
behaves as $p\searrow p_c$ in $d\ge3$ is an open
question for which our techniques may offer a point of entry.

Our interest in finite connections arises in an ongoing study of random walk
on the infinite open cluster of supercritical percolation on $(\Z^d,\sim)$
that is biased towards a preferred direction.  This model undergoes
a phase transition, as the strength of the bias increases, from a regime
of linear long-run progress to a sublinear regime in which traps detain the
walk for lengthy periods.  Building
on~\cite{BGP,Sznitman}, \cite{FriberghHammond} proved the existence of
a critical bias for this transition.  In the sublinear regime, traps
resemble lengthy finite percolation connections opened at one end, into
which the walker may enter to become waylaid near the opposing end.
Arguments in~\cite{FriberghHammond} used the first-order exponential decay
$e^{-\zeta u(1+o(1))}$ established there.  The mass gap derived in the present paper speaks to a far more precise
description of trap geometry in a broader study of biased walk in disorder.

{\bf Organisation.} This paper is presented in four parts.
The action of a single slide is the principal subject of Part~\ref{part.one}, which ends with a roadmap indicating how the operation will be iterated to establish the existence of renewal.
Part~\ref{part.two} develops the tools needed to control the resulting slide dynamics.  In Part~\ref{part.three}, renewal production and survival are proved using these dynamics.  Part~\ref{part.four} converts the resulting renewal estimate into the mass gap and develops its renewal-abundance and Ornstein--Zernike consequences.
A glossary of commonly used terms appears in Appendix~\ref{s.notation-guide}.

{\bf The preprint~\cite{FriberghHammondPreprint}.}
The present work largely subsumes the arXiv preprint~\cite{FriberghHammondPreprint}.  Much of the technical analysis in that preprint concerns the action of a single slide, developed in a form suitable for its later use along an iterated schedule.  This material is incorporated into Chapter~\ref{s.single-slide-positive-density}, in places essentially unchanged and elsewhere revised or developed.  The iterated slide schedule is conceptually central to the preprint, although the analysis specific to the iteration is concentrated near its end, in Sections~4.3 and~4.4 of~\cite{FriberghHammondPreprint}.  By making of order \(u^{1/2}\) slide attempts, separated on the scale \(u^{1/2}\), the argument gives an upper bound of the form \(Ce^{-cu^{1/2}}\) on the conditional probability of a long connection without renewal, and hence a corresponding stretched-exponential tail for the tilted sausage height.  Allied with the single-slide analysis, these two sections thus provide what is needed for the principal renewal conclusion of that paper, recorded in particular in the sausage-height estimate of Lemma~2.18.

While most of the preprint is absorbed into the present work, the arguments of Sections~4.3 and~4.4 are not reproduced.  The present paper instead develops a substantially more involved slide dynamics on a positive-density schedule, controlling the interactions between successive slides in order to replace the stretched-exponential conclusion by the exponential mass gap of Theorem~\ref{t.mass-gap}.  The distinction between the two schedules, and the additional difficulties created by positive-density iteration, are discussed in Section~\ref{ss.iterated-slide-overview}.

{\bf A recent preprint by Aoun and Khettabi.}
Aoun and Khettabi~\cite{AounKhettabi} have obtained, by a different route, the point-to-point Ornstein--Zernike asymptotic appearing in Theorem~\ref{t.oz-intro}.  Their argument uses the tilted sausage law of~\cite{FriberghHammondPreprint}, and in particular the just-mentioned stretched-exponential height estimate~\cite[Lemma~2.18]{FriberghHammondPreprint}, to obtain an a priori volume bound for long finite clusters.   
They then develop an ingenious adaptation of the Campanino--Ioffe--Velenik coarse-graining scheme,
in which the radii of coarse-graining regions are selected so that their boundaries meet relatively few edges incident to the cluster, and the volume bound is invoked to control the aggregate cost of separating the regions. 
 This yields an exponential estimate~\cite[Theorem~$3.8$]{AounKhettabi} for cone-point scarcity among connections whose volume is bounded by a constant multiple of their length, and thus a form of mass gap with a volume restriction that holds, under finite-connection conditioning, with a stretched-exponential failure probability admitting an upper bound of exponent $\tfrac{d-1}{2d}$.
 When this restriction is removed, Aoun and Khettabi obtain sufficient control to derive a cone-point decomposition and an effective directed random walk with stretched-exponential increment tails. This leads to the sharp point-to-point asymptotic~\cite[Theorem $1.2$]{AounKhettabi} with a positive continuous amplitude.
This argument shows, in particular, that this point-to-point asymptotic may be obtained without establishing an exponential mass gap for the tilted sausage law itself: stretched-exponential control after a suitable coarse graining suffices.  The central new conclusion of the present work is precisely such an exponential mass gap for the sausage law~$\mu$. 
This exact exponential renewal structure permits a direct IID multivariate renewal analysis yielding the point-to-point asymptotic and, in addition, the positive real-analytic amplitude and the real-analytic, strictly convex equidecay and Wulff boundaries asserted in Theorem~\ref{t.oz-intro}.

{\bf Formation of the paper and A.I. usage.}
This paper was written in two main stages.  The first, ending in 2024, produced the preprint~\cite{FriberghHammondPreprint}; the second, in 2026, developed the mass-gap argument and the present manuscript.   Part~\ref{part.one} descends substantially  from the 2024 preprint, though it has been revised and reorganised.  During the later phase we made extensive use of large language models, principally GPT-5.6 Sol, Claude Fable 5.1 and Claude Opus 4.8.  Claude was often used to critique chapters or sections, while GPT was used extensively in revising exposition, checking notation, and working through finer proof details.  The models sometimes proposed several paragraphs of text at a time, but such material was closely reviewed and repeatedly rewritten, often sentence by sentence, for both mathematical content and style.  In a long technical project, this assistance was practically valuable: it helped to locate gaps, stabilise notation, test arguments, and keep many local pieces of the proof moving coherently.

The authors developed the principal random-geometric constructions and directed the overall mathematical programme and proof architecture, including much of the fine-grained organisation of arguments. Much of the model assistance took place within this framework, entailing not only exposition and checking but also the formulation, development and revision of subsidiary mathematical arguments.

Moreover, there is a part of the mass-gap argument in which GPT’s assistance was more substantially conceptual: the assembly of the local geometric and accrual estimates into the multi-layer analytic framework used to prove the IVF-to-LLF transition.
Because the division between the underlying random-geometric construction and this later analytic assembly is most naturally described using the concepts and terminology developed in the argument itself, a more precise account is given at the end of the opening discussion of Chapter~\ref{s.shadow-bounds-delivery-fine}.

Across a project involving many such exchanges, it is difficult to characterise succinctly what the finished work owes to model assistance.  Our necessarily rough assessment is that model assistance typically contributed to the detailed development, checking and exposition of the proofs, while, in the multi-layer accrual analysis just mentioned, GPT supplied a conceptually valuable part of the argument that secures a positive density of renewal-production opportunities in the mass-gap proof.
  The authors take full responsibility for all statements, proofs and errors in the paper.

\section{Acknowledgements}

The second author is supported by the National Science Foundation under
DMS grants~1855550 and~2153359 and by the Simons Foundation as a 2021
Simons Fellow.  He thanks Hugo Duminil-Copin and Yvan Velenik for
helpful discussions around 2012 about finite connections in dimension
at least three.

\part{The Slide Method}\label{part.one}

\chapter{The renewal abundance theorem}
\label{s.renewal-abundance}

In Part~\ref{part.one}, we develop the slide method that will be used to locate renewal points.  Chapter~\ref{s.heuristics} presents the slide in overview, Chapter~\ref{s.single-slide-positive-density} analyses the single-slide operation rigorously, and Chapter~\ref{s.iterated-slide} introduces its iteration and concludes with a roadmap for the central slide argument.

In the present chapter, we define the renewal objects precisely and state the core
renewal estimates.  Theorem~\ref{t.intermediate-renewal}(2), which asserts that
a conditioned finite cluster fails to contain a renewal point only with
probability exponentially small in its height, is the outcome of the slide argument and the paper's principal technical estimate.
After passage to the independent-sausage process $\mu^{\otimes\N}$ and
renewal-theoretic post-processing, it yields the stronger renewal abundance
Theorem~\ref{t.renewal}, in which linearly many renewal points are present with
overwhelming probability.

In the subcritical Ornstein--Zernike analyses
of~\cite{CampaninoIoffe,CampaninoIoffeVelenik}, the rare object whose geometry
must be controlled is a long open path.  In the supercritical phase such paths
are commonplace; what is rare is a long open path that remains sealed off from
the infinite cluster.  The object that here plays the role of the subcritical
open path is therefore the surface of dual plaquettes that wraps a finite
cluster and separates it from the infinite one.  The event of a long finite
connection is neither increasing nor decreasing: opening an edge may lengthen
a finite cluster or, by joining it to the infinite cluster, destroy its
finiteness.  Where the subcritical theory
draws its control from disjoint occurrence and decay of correlations---the BK
inequality bounding disjoint connection probabilities by a product of
marginals~\cite[(2.2)]{CampaninoIoffe}, and exponential decay of correlations
supplying approximate independence between separated regions---it is the
iterated slide resampling behind
Theorem~\ref{t.intermediate-renewal}(2) that, by demonstrating renewal,
supplies control of the rare dual surface tightly wrapping the open path
within.

The height geometry, including slabs and the formal definition of a renewal
point, is recorded in Section~\ref{s.basics}, and the base, the head and the
string in Section~\ref{s.maincharacters}.  The sausage decomposition is developed later in
Section~\ref{s.sausage-law-mass-gap} and the renewal measure in
Section~\ref{s.provingsting}.  The route from
the core renewal estimates to the theorems stated above, and with it the
organisation of the remainder of the paper, is then set out in
Section~\ref{ss.principal-estimate-road-ahead}.

\section{Height geometry, renewal, and slabs}\label{s.basics}

 We regard $\vecell$ as pointing vertically in the positive direction; so $x \in \R^d$
is upwards from $z \in \R^d$ if $x \cvl > z \cvl$.

We adopt the \emph{axes ordering convention}: without loss of generality,
by reflecting and permuting coordinate axes, assume
\begin{equation}\label{e.aoc}
  e_i\cvl\, \geq \, 0 \quad\text{for }1\leq i\leq d,
  \qquad
  e_d\cvl\, = \, \max\{e_i\cvl : 1\leq i\leq d\}\, .
\end{equation}
Thus $r=e_d\cvl\geq d^{-1/2}$.

For $h \in \R$, the \hfff{height-spaces}{{\em forward} and {\em backward} half-spaces} of $h$ are given by
$$
 F_h\, = \, \big\{ u \in \R^d: u \cdot \vecell\, > \, h \big\}
 \, \, \, \textrm{and} \, \, \, B_h\, = \, \big\{ u \in \R^d: u \cdot \vecell \leq h \big\}\, .
$$

For $e = (x_1,x_2) \in E(\Z^d)$, write
$\closedint{e} = \big\{ \lambda x_1 + (1-\lambda) x_2:
\lambda \in [0,1] \big\} \subset \R^d$
for the endpoint-interpolating interval.
For $h \in \R$, we define the \emph{upspace} $\upspace_h$ and \emph{downspace} $\downspace_h$ of~$h$ to be subgraphs of $(\Z^d,\sim)$ whose edge sets are given by
$$
 E(\upspace_h)\, = \, \big\{ e \in E(\Z^d): \closedint{e} \cap F_h \not= \emptyset \big\} \, \, \, \textrm{and} \, \, \,  E(\downspace_h)\, = \, \big\{ e \in E(\Z^d): \closedint{e} \subset B_h   \big\}\, ,
$$
and whose vertex sets $V(\upspace_h)$ and $V(\downspace_h)$
are the sets of endpoints of elements of the respective edge sets.
Equivalently, if \(h^+(e)\) is the greater endpoint height of~\(e\), then
\[
 e\in E(\upspace_h)\ \Longleftrightarrow\ h^+(e)>h,
 \qquad
 e\in E(\downspace_h)\ \Longleftrightarrow\ h^+(e)\leq h.
\]
Thus these edge sets partition \(E(\Z^d)\); in particular, an edge contained
in the hyperplane of height~\(h\) is assigned to \(\downspace_h\).

\begin{figure}[htbp]
\centering
\includegraphics[width=0.75\textwidth]{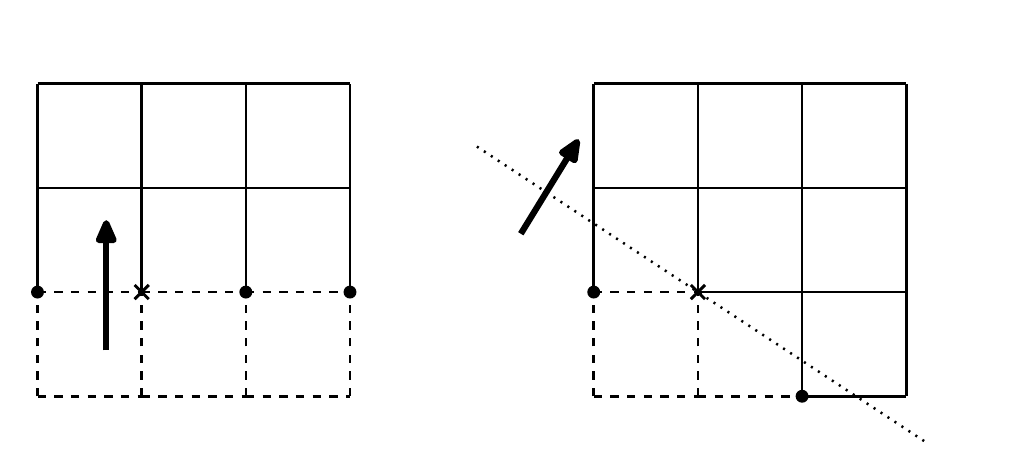}
\caption{The partition $E(\upspace_h) \cup E(\downspace_h)$ of edges in $\Z^2$ for an axial and a more general choice of $\vecell$. The arrow represents $\vecell$, while the cross marks the origin, which is an element of $\vertbdry_h$ for $h=0$. The discs mark the other elements of $\vertbdry_0$. Elements of $E(\upspace_0)$ are drawn with solid lines; those of $E(\downspace_0)$ are dashed. In accordance with the boundary convention, edges contained in the height-zero hyperplane are dashed.}\label{f.fblusketch}
\end{figure}

The {\em vertex boundary} $\vertbdry_h$ of $h$ equals the intersection of vertex sets $V(\upspace_h) \cap V(\downspace_h)$. The various characters appear in Figure~\ref{f.fblusketch}.

Next is the important definition of a renewal point.

\vspace*{3pt}

\begin{definition}\label{d.renewal}
Let $U$ be a finite connected subgraph of $(\Z^d,\sim)$
and let $x \in V(U)$.
The element $x$ is called a \hfff{renewal-point}{\em renewal point} of~$U$ if $V(U) \cap \vertbdry_{x \cvl}$
equals the singleton set $\{ x \}$.
We write $R(U)$ for the set of renewal points of~$U$ and
$\rlevel(U)=\{y\cvl:y\in R(U)\}$ for the corresponding set of renewal
levels.
\end{definition}

We record two elementary height facts.

\begin{lemma}\label{l.geometrybasics}
Let $h \in \R$.
\begin{enumerate}
\item $\vertbdry_h = \{v \in \Z^d : v\cvl \in (h-r,h]\}$.
\item If $S$ is a finite connected subgraph of $(\Z^d,\sim)$, the
heights of distinct elements of $R(S)$ differ by at least $r$.
\end{enumerate}
\end{lemma}
\noindent\textbf{Proof: (1).}  If $v \in \vertbdry_h$, then an edge incident to $v$ belongs to
$\downspace_h$, so $v\cvl \leq h$.  Another incident edge belongs to
$\upspace_h$; since one lattice step changes height by at most~$r$, this
gives $v\cvl>h-r$.
Conversely, if $v\cvl \in (h-r,h]$, then $(v,v+e_d)$ belongs to
$\upspace_h$ and $(v-e_d,v)$ belongs to $\downspace_h$. \\
{\bf (2).} Let distinct $x \in R(S)$ and  $y\in \Z^d$ satisfy $x\cvl\le y\cvl<x\cvl+r$.
Then (1) places $x$ in $\vertbdry_{y\cvl}$, so that~$y \not\in R(S)$. \qed

Thus Definition~\ref{d.renewal} is equivalent to the formulation in
Section~\ref{ss.renewal-intro}: $x$ is the only vertex of $U$ in the band of
width~$r$ topped by its own level.  Moreover, for every finite connected~$U$,
the map $x\mapsto x\cvl$ is a bijection from $R(U)$ to $\rlevel(U)$, and
distinct renewal levels are at least~$r$ apart.

For $m,n \in \Z$ with $m \leq n$, we denote by $\llbracket m,n \rrbracket$ the integer interval $\big\{ k \in \Z: m \leq k \leq n \big\}$, and write~$\intint{n}$ for $\llbracket 1,n\rrbracket$.

Let $G_1=\bigl(V(G_1),E(G_1)\bigr)$ and
$G_2=\bigl(V(G_2),E(G_2)\bigr)$ be two subgraphs of $(\Z^d,\sim)$. Their
intersection and union are the subgraphs
$$
G_1\cap G_2=\bigl(V(G_1)\cap V(G_2),\,E(G_1)\cap E(G_2)\bigr) \, \, \, \textrm{and} \, \, \, 
G_1\cup G_2=\bigl(V(G_1)\cup V(G_2),\,E(G_1)\cup E(G_2)\bigr) \, .
$$
Thus $G_1\cap G_2$ may contain isolated vertices, since this graph contains $v \in V(G_1) \cap V(G_2)$ even when~$v$ is incident to no element of $E(G_1) \cap E(G_2)$.
For $y,z\in\R$ with $y\le z$, we define the slab $\slab_{y,z}$ from $y$ to $z$ to be the subgraph
$\upspace_y\cap \downspace_z$.
The axes ordering convention and Lemma~\ref{l.geometrybasics}(1) give
\[
 V(\slab_{y,z})=\{v\in\Z^d:v\cvl\in(y-r,z]\}\, .
\]

We name the direction of increasing height \emph{north} and that of decreasing
height \emph{south}.  A slab is bounded only in the height direction, so its
boundary consists merely of its lower and upper sides; the slide arguments depend
on this simple feature, which we now record with the sides named.

\begin{lemma}[Entry to a slab is through its sides]
\label{l.slab-entry-faces}
For $y\le z$ the slab $\slab_{y,z}=\upspace_y\cap\downspace_z$ has
\emph{southside} $\vertbdry_y$ and \emph{northside} $\vertbdry_z$,
namely its lower and upper boundaries.  Any element of $V(\slab_{y,z})$
adjacent in $(\Z^d,\sim)$ to a vertex outside $V(\slab_{y,z})$ lies in
$\vertbdry_y\cup\vertbdry_z$.  Consequently every path with
one endpoint outside the slab and the other inside it meets $\vertbdry_y$ or
$\vertbdry_z$.
\end{lemma}

\noindent\textbf{Proof.} Let $v\in V(\slab_{y,z})$ be adjacent to $w\notin V(\slab_{y,z})$.
If $w\cvl\leq y-r$, then $v\cvl\leq w\cvl+r\leq y$, while
$v\cvl>y-r$, so Lemma~\ref{l.geometrybasics}(1) places
$v\in\vertbdry_y$.  If $w\cvl>z$, then
$v\cvl>w\cvl-r>z-r$, while $v\cvl\leq z$, so the same lemma places
$v\in\vertbdry_z$.  These alternatives exhaust the complement of
$V(\slab_{y,z})$.  On a path running from outside the slab to inside, the
first vertex in $V(\slab_{y,z})$ has its predecessor outside, proving the
last assertion.
\qed

\section{The base, the head, and strings}\label{s.maincharacters}

\subsection{The base and the head}\label{s.baseandhead}
For certain choices of $\vecell$, the hyperplane
$H_s=\{z\in\R^d:z\cvl=s\}$ may meet a finite set $A\subset\Z^d$ in more
than one point.  Since $H_0=\vecell^\perp$, we use the fixed orthonormal basis
of $\vecell^\perp$ chosen above to order $H_0$ lexicographically.  The map
$x\mapsto x-s\vecell$ identifies $H_s$ with $H_0$ and transports this order
to a total order on $H_s$.

Let $G$ be a finite subgraph of $(\Z^d,\sim)$. The base of $G$, denoted by $\base(G)$, is the lexicographically minimal element of $V(G)$ among those of minimum height.
For the head $\head(G)$, we replace `minimal' by `maximal' and `minimum' by `maximum'.
The usage of the term `height' is extended so that the height of $G$ is given by $\max \{ (v-u) \cdot \vecell : u,v \in V(G) \}$.  When the base of $G$ is the origin, this equals the maximum vertex height.

\subsection{The string}\label{s.tubestingcluster}
\begin{figure}[htbp]
\centering
\includegraphics[width=0.35\textwidth]{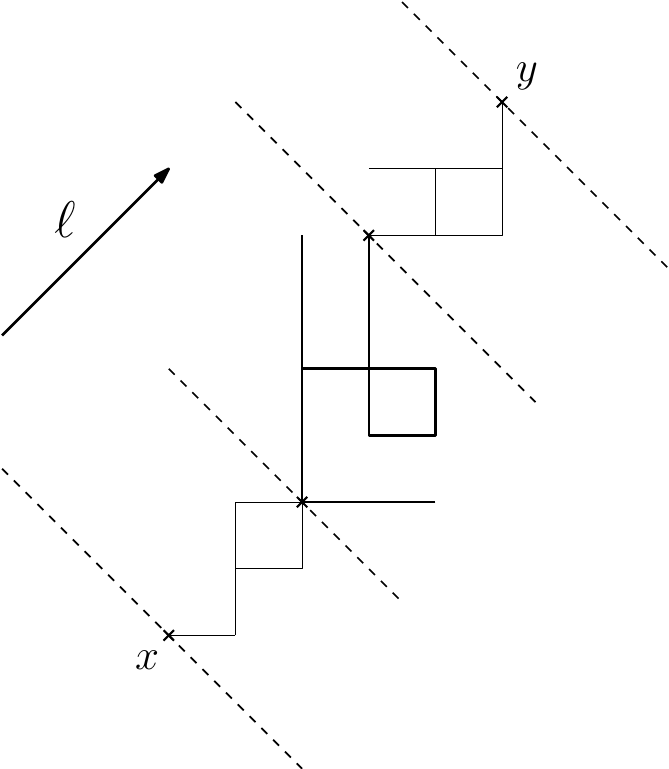}
\caption{Four crosses mark the renewal points in this string from $x$ to $y$. The middle of its three sausages is drawn in bold.}\label{f.sting}
\end{figure}

Let $x,y \in \Z^d$, $x \cvl < y \cvl$.
A \hfff{string}{string $S$ from $x$ to $y$}
is a finite open connected component in $\slab_{x \cvl,y \cvl}$
such that $x,y$ are renewal points of~$S$.
The term `string' is shorthand for `string of sausages': as anticipated above, a string breaks into sausages delimited by consecutive renewal points; see Figure~\ref{f.sting}.

\begin{lemma}[Cluster segment between renewal points]
\label{l.cluster-renewal-string}
Let $C$ be a finite open cluster and let $x,y\in R(C)$ satisfy
$x\cvl<y\cvl$.  The open component of $x$ in
$\slab_{x\cvl,y\cvl}$ is a string from $x$ to~$y$.
\end{lemma}
\noindent\textbf{Proof.} Choose a self-avoiding open path in $C$ from $x$ to~$y$.  If the path leaves
the slab through its southside after departing from~$x$, then
Lemma~\ref{l.slab-entry-faces} forces it to revisit $x$ when it re-enters,
because $x$ is the only cluster vertex in $\vertbdry_{x\cvl}$.  Nor may it
leave through the northside before reaching~$y$, because $y$ is the only
cluster vertex in $\vertbdry_{y\cvl}$.  Thus the path remains in the slab
until it reaches~$y$.

The open component~$S$ of $x$ in the slab therefore contains~$y$.  Since
$C$ is the full-lattice open cluster of~$x$, every vertex of~$S$ belongs
to~$C$.  Hence $x$ and $y$ remain renewal points of~$S$, and $S$ is a string
from $x$ to~$y$. \qed

For $u > 0$ and $x \in \Z^d$, let $\sting_u(x)$ denote the event
there exists a string from $x$ to $z$ for some $z \in \Z^d$ with $(z - x) \cvl = u$.
Calling a height \emph{attainable} if it is achieved by a lattice vertex,
note that $\sting_u(x)=\emptyset$ unless $u$ is attainable.
When $\sting_u(x)$ occurs, the slab is fixed by
$x$ and~$u$, so the string is the open component of~$x$ in this slab and
$z$ is its unique highest vertex.  Thus both $z$ and the string are unique.
We denote this string by $\stingset_u(x)$.
When $x = 0 \in \Z^d$, write $\sting_u = \sting_u(0)$ and $\stingset_u = \stingset_u(0)$.

The interval notation already introduced for cluster events extends to strings
as follows.

\vspace*{3pt}

\begin{definition}\label{d.enotation}
For an interval $I\subseteq(0,\infty)$, put
$\sting_I=\bigcup_{u\in I}\sting_u$.  We use this notation when $I$ is finite
or semi-infinite.
\end{definition}

\section{Seeking renewal by resampling: once and iteratively}
\label{ss.principal-estimate-road-ahead}

We state Theorem~\ref{t.intermediate-renewal} and outline the route by which it yields the results stated in the introduction.

Write \hfff{conditioned-law}{$\Q = \PP(\,\cdot\mid\clusteruthree)$} for the conditional law under which
the origin is the base of a finite cluster of height in $[u,u+3r]$.  The
dependence $\Q=\Q_u$ on~$u$ is left implicit.  The
principal aim is to establish Theorem~\ref{t.intermediate-renewal}(2) for this
conditioned law.  The bounded height offset $3r$ secures a sufficiently rich
set of admissible heights for the base of the translated upper piece in the
slide resampling described below when $\vecell$ is non-lattice.  We say that an event holds \hfff{overwhelming-probability}{\emph{with overwhelming probability}}
under a given probability law if its probability is at least $1-Ce^{-cu}$ for
some $C,c>0$ and all large $u$.

\begin{theorem}[Core renewal estimates]
\label{t.intermediate-renewal}
Let $d\geq3$ and $p\in(p_c,1)$.
\begin{enumerate}
\item \textup{(Positive-probability renewal density.)}
There exist $a,c>0$ and $u_0\in\N$ such that, for all $u\geq u_0$,
\[
   \Q\bigl(|R(\mc C(0))|\geq au\bigr)\, \geq \, c\, .
\]
\item \textup{(A renewal point is overwhelmingly likely.)}
There exist $C,c>0$ such that, for all $u\geq r$,
\[
   \Q\bigl(\rlevel(\mc C(0))\cap(0,u)\neq\emptyset\bigr)
\, \geq \, 1-Ce^{-cu}\, .
\]
\end{enumerate}
\end{theorem}

Part~(1) says that a positive density of renewal points occurs with
probability bounded away from zero.  Part~(2)---that a renewal-free conditioned
cluster is exponentially rare in its height---is the central estimate of the
paper.  It is the input from which the mass gap and, through renewal theory and
endpoint estimates, the renewal abundance Theorem~\ref{t.renewal} and the
asymptotic results stated in the introduction are derived.

Both estimates are obtained by slide resampling.  One cuts
the conditioned cluster across a chosen height into a lower piece and an upper
piece joined by a midsection; the midsection is discarded, the upper piece is
slid to a new position relative to the lower piece, and their junction is
reconstructed with fresh randomness.  The operation is built to preserve the
conditional law $\Q$.
The single-slide analysis shows that suitable heights offer a uniformly
positive chance to produce a renewal point; together with the abundance of
such heights, this yields Theorem~\ref{t.intermediate-renewal}(1).

An increasing schedule containing a positive density of heights provides
order $u$ opportunities for renewal production.  With this abundance comes
the danger that slides interact: a later slide may undo the renewal status of
a point created earlier.  To prove Theorem~\ref{t.intermediate-renewal}(2),
the iterated argument must show that, with probability $1-Ce^{-cu}$, at least
one renewal produced along the way survives to the terminal state.

The derivation of this result through a delicate analysis of the slide dynamics on a positive-density schedule will occupy us until the end of Part~\ref{part.three}.  Deferring a heuristic overview of this argument, we instead comment briefly on what follows once Theorem~\ref{t.intermediate-renewal}(2) has been obtained.

The analysis in  Part~\ref{part.four} passes through the independent-sausage renewal
process.  By Lemma~\ref{l.cluster-renewal-string}, a finite cluster, cut at
its lowest and highest renewal points, is a middle string between two endpoint
culs-de-sac; the string decomposes into irreducible sausages.  Theorem~\ref{t.intermediate-renewal}(2),
combined with the comparison between sausages and full clusters, supplies an
exponentially decaying tilted sausage-height tail.  Together with the string
decay rate, this normalises
$\mu(S)=e^{\zeta\hgt(S)}\PP(\mathsf O(S))$ and proves the mass gap
Theorem~\ref{t.mass-gap}.  Concatenation of independent $\mu$-sausages then
gives a one-dimensional renewal process.  The renewal theorem yields the sharp
string asymptotic, and the cul-de-sac--string decomposition and endpoint
estimates give the sharp point-to-plane cluster asymptotic.  A Cram\'er estimate
shows that strings with too few renewals have an additional exponential cost;
the endpoint estimates and the point-to-plane asymptotic transfer this bound to
conditioned clusters, yielding the renewal abundance Theorem~\ref{t.renewal}.
Finally, control of the full vector increment of a sausage and its behaviour
under small rotations of~$\vecell$ gives the Ornstein--Zernike conclusions.


\chapter{The slide resample in overview}\label{s.heuristics}

The core renewal estimate underlying Theorem~\ref{t.renewal} will be proved by
an iterative {\em slide resampling} operation.  In this chapter, we offer a rough conceptual guide to the motivation
for and the form of this operation.

The iterative procedure is composed of a sequence of slides along a schedule of heights.
The slide at a given height
 itself requires some care to set up, and it may be harnessed to derive
the positive-probability renewal density of
Theorem~\ref{t.intermediate-renewal}(1).
The present chapter gives a conceptual account of the single slide;
Chapter~\ref{s.single-slide-positive-density} implements it rigorously and proves that
estimate.

In the informal treatment of the present chapter, we take $\vecell =e_d$ to be vertical.
In this axial case, we may dispense with the nicety of the bounded offset~$3r$ in conditioned cluster height.
Thus we will consider the typical form of the cluster $\mc{C}$ under percolation $\PP$
given the event $\cluster_u$ that a finite cluster based at the origin runs to height equal to $u \in \N$. We will write $\mc{C}_u$ for $\mc{C}$ to emphasise this conditioning.

\section{Dynamic formation of the cluster}

Under the law $\PP \big( \cdot \big\vert \cluster_u \big)$,
the cluster $\mc{C}_u$ is finite and contains the origin, with the interval of heights occupied by its vertices being exactly the integer interval~$\llbracket 0, u \rrbracket$.

Take the unconditioned percolation law~$\PP$ and seek to form the cluster~$\mc{C}_u$ in stages indexed by a given increasing sequence $\big\{ h_i: i \in \intint{L} \big\}$  of integer heights lying in~$[0,u]$.
At the first stage, percolation is realised in the downspace~$\downspace_{h_1}$.
Write $O \subset \downspace_{h_1}$ for the open component containing the origin of this percolation, and denote by $\Top(O)$ the set of elements of~$V(O)$ of maximum height.
 If the event $\cluster_u$ is to occur,
 then $O$ must be finite, with the origin its base, and the common height of elements of $\Top(O)$ must equal~$h_1$, so that $O$
 reaches the upper bounding hyperplane at height~$h_1$. When $\cluster_u$ occurs, this component is a subgraph of $\mc{C}_u$
that we will call the {\em downcluster}~$\downcluster_{h_1}$ associated to height~$h_1$.
At the $i$\textsuperscript{th} stage, the presently realised percolation is extended to~$\downspace_{h_i}$. On the event~$\cluster_u$, this produces a nested downcluster sequence
$$
  \downcluster_{h_1} \subset \downcluster_{h_2}  \subset \cdots \subset \downcluster_{h_L}\, ,
$$
with percolation in $\downspace_{h_L}$ sampled by the end of the $L$\textsuperscript{th} stage. In a final stage, the configuration in the remainder of the space, namely the upspace~$\upspace_{h_L}$, is realised, and the occurrence of $\cluster_u$ is checked, by examining whether $\downcluster_{h_L}$ forms part of a finite whole-space cluster whose head has height~$u$.

In this way, we may view the cluster~$\mc{C}_u$ under $\PP \big( \cdot \big\vert \cluster_u \big)$ as being formed dynamically, in a nested increasing sequence of downclusters, and a completing step.
Consider a height schedule $\{h_i: i \in \intint{L}\}$ with bounded positive increments between consecutive terms and with $L$ of order~$u$.
Under the unconditioned law~$\PP$, the formation of~$\cluster_u$ may fail at any given stage. Take~$i \in \intint{L}$, and suppose that the nested sequence has formed until entry
$\downcluster_{h_{i-1}}$. Elements of $\Top(\downcluster_{h_{i-1}})$ lie in the southside border of the slab $\slab_{h_{i-1},h_i}$ on which percolation will be sampled at the~$i$\textsuperscript{th} stage. These elements may be viewed as seeds from which open structure may form.
If this stage is to be successful, with~$\downcluster_{h_i}$ being added to the forming nested sequence,
then open growth from the seeds must occur in a controlled way. One threat to such growth is represented by the possibility that every edge in $\upspace_{h_{i-1}}$ incident to $\Top(\downcluster_{h_{i-1}})$ happens to be closed. In this event, the whole-space cluster of the origin is simply~$\downcluster_{h_{i-1}}$, since all avenues for addition to this downcluster are closed off.
A second threat is posed by the possibility that open paths run from sites in $\Top(\downcluster_{h_{i-1}})$ across the slab $\slab_{h_{i-1},h_i}$. By the general thickening form of Grimmett and Marstrand~\cite[Theorem A]{GrimmettMarstrand}, sufficiently wide slabs of any given orientation are supercritical. Thus the constant increment $h_i-h_{i-1}$ may be chosen sufficiently large, as a function of $p>p_c$, that this slab contains an infinite open connected component. The forming open paths may connect to this cluster in an event that would render~$\downcluster_{h_i}$ infinite, incompatibly with the occurrence of~$\cluster_u$.

The dynamic formation is thus harried at each of linearly many stages by this pair of threats: if the present downcluster top is small, a few local closures may terminate exploration; and, whatever the value of this top cardinality, but all the more so if it is high, growth into the upcoming slab risks {\em open escape} and a failure of finite connection. The mantra
\[
   \hbox{{\bf t}hin: {\bf i}nstantly  {\bf p}erish;}
   \qquad
   \hbox{ {\bf t}hick:  {\bf o}penly  {\bf e}scape}\, .
\]
represents the twin threats, which are issued on an order of~$u$ occasions.
Under  the law $\PP \big( \cdot \big\vert \cluster_u \big)$,
 in which percolation is conditioned on an event whose probability has exponential decay in~$u$,
 the exploration survives all these threats; it may be viewed as advancing on \hfff{tiptoe}{TipToe},
harnessing a good fortune that is exponential in~$u$
in an ongoing geometric balance at the forming downcluster top  that continually manages both threats.


\section{Studying a probability measure by resampling}
\label{s.resampling-overview}

Dynamic cluster formation has offered a prospect on the TipToe dichotomy that guides our analysis, but it is not the method we will use to prove renewal abundance.

Our tool is instead resampling. In order to study a probability measure of interest called~$\Q$, a resampling operation $\omega \mapsto \omega'$ is specified that maps the state space of~$\Q$
to itself in such a way as to preserve this law. Think of representing the input configuration or sample point~$\omega$
as an ordered list of pieces of information that characterise it, and then drawing a line under a certain item in this list. Information above the line is retained; the remainder is lost. The remaining information is reconstructed according to the requisite conditional probability, given the form of the retained data, and the reconstructed data is recorded below the line, to produce an altered list, and an output~$\omega'$ whose law is~$\Q$.

Such a tool is used to study~$\Q$ via its output~$\omega'$. Perhaps we wish to show that a given event is typical for~$\Q$. The input~$\omega$ may lack the requisite property. If the retained data is well chosen, the information imparted by the input~$\omega$ may act in concert with sufficiently powerful independent randomness in the reconstruction step that a favourable output~$\omega'$ is typically produced for most choices of the input~$\omega$.

Of course care is needed if this method is to work. To give the pair of trivial extreme examples, the operation may act as the identity $\omega \mapsto \omega$, so that the output always equals the input;
or the output~$\omega'$ may be selected as an independent copy of the law~$\Q$, without reference to the form of the input~$\omega$. In the first case, the input is informative but the action is uselessly weak, offering no prospect for undoing features of~$\omega$ that prevent the occurrence of the putatively typical event. In the latter, the action is too broad or powerful, rendering the input totally uninformative.

\section{Undoing the capital \texorpdfstring{$N$}{N} by a single slide resample}\label{s.slideintorenewal}

The single slide resample illustrates this framework. The input~$\omega$ may lack a desired feature, but, in the presence of a geometric condition that is known to be typical,
the resampling action $\omega \mapsto \omega'$ may serve to show that the output satisfies this feature typically, or with positive probability.

Our discussion of dynamic cluster formation may suggest that, under a typical sample of $\PP$ given $\cluster_u$, the cluster $\mc C_u$ has a bounded cross-section at typical heights.  The {\em global surface bound} Proposition~\ref{p.nobigcluster} delivers a rigorous version of this typical boundedness through control of the exterior boundary of a coarse graining of the cluster, and this feature of the conditioned random geometry plays a role in demonstrating that the single slide resample output has a renewal point near the slide height with positive probability.  Proposition~\ref{p.nobigcluster} will be proved using the coarse-graining method of Antal and Pisztora~\cite{AntalPisztora}; we will develop this coarse-graining framework considerably in the slide analysis.

Here we briefly signpost some ideas in the derivation of the global surface bound. Dispense with the dynamic formation and view the cluster $\mc{C}_u$ in
its entirety, as a static object. It is a finite open subgraph of vertical
extension $u$ that is encased away from the infinite open cluster by
dual plaquettes through whose midpoints run closed edges in the cluster
boundary.  This dual plaquette
surface is in essence subcritical, as a coarse graining serves to show.  The endpoints of the closed edge through a
constituent plaquette lie in different open clusters and thus have infinite
open chemical distance.  The coarse graining uses a finite local witness
of this separation: failure of the two sides to connect inside a suitably
enlarged box.  Such failures are represented by bad boxes, large connected collections of which are exponentially rare.

We will develop these ideas in proving Proposition~\ref{p.nobigcluster} and in deriving from it bounded cluster intersection in typical height bands.  It is worth comparing this with the highly
supercritical phase, where Ornstein--Zernike asymptotics were obtained by
Campanino and Gianfelice~\cite{CampaninoGianfelice}.  There the global volume
very typically lies within a factor $1+\varepsilon$ of its minimum---in an axial
direction the cluster is essentially a straight line of length~$u$---so that
renewal points already sit at high density.  Our global surface bound is reached by the
same coarse-graining, and is in this sense a reduction to that phase; but the
reduction is only partial, with the volume contributed by the typical bounded
cross-sections amounting to a \emph{large constant} multiple of~$u$, rather
than the near-minimal $(1+\varepsilon)u$ geometry of the highly supercritical
phase.  The single-slide argument must therefore handle,
at a typical height, cross-sections of bounded cardinality that are nevertheless
spatially spread, a feature absent in the highly supercritical regime.  We now
illustrate this geometry.

\begin{figure}[htbp]
\centering
\includegraphics[width=0.75\textwidth]{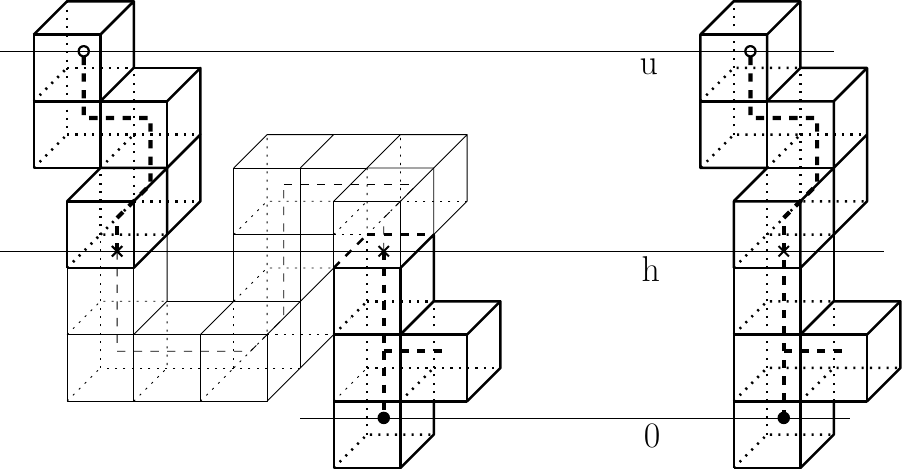}
\caption{A slide eliminates from a capital $N$ its midpiece.}\label{f.Nslide}
\end{figure}

The cluster $\mc{C}_u$ runs between heights zero and $u$. The single slide operation is indexed by a height $h \in [0,u]$
that we take to be typical, say one around level $u/2$. The resulting typical-band control indicates that the cross-section of $\mc{C}_u$ at height $h$ is likely to have a cardinality that is bounded above independently of $u$.
If this cross-section occurs at a unique vertex, that vertex would be a renewal point for the cluster; so in seeking to derive renewal abundance, we wish to argue that this circumstance is typical under the conditioned law. As it stands, though we know that the cross-section is typically bounded in cardinality, it may be high in diameter. The cluster could resemble a capital `$N$', as depicted in Figure~\ref{f.Nslide}(left).
Elements of $E(\mc{C}_u)$ are dashed, and have the rough shape of a path that runs upwards from the base, doubling back across level $h \approx u/2$, and then running upwards again to the head. Slicing through at height $h$ decomposes the cluster into three pieces: emboldened in the depiction are the downcluster and the upcluster,
the former already seen in the discussion of dynamic formation, these being the respective open connected components of the cluster's base or head within the half-spaces delimited by the slice;
 and a remaining part, the midpiece, which runs between the two vertices marked with crosses that lie one apiece in the downcluster and the upcluster.

The capital $N$ example is instructive as we consider a useful definition of the single slide operation. The midpiece is a costly item in the example, because it is sheathed by a subcritical surface  of dual plaquettes away from the two surfaces associated to the upcluster and the downcluster.
Not only expensive but also, it would seem, needless: the downcluster covers a height of $h$; the upcluster, one of $u-h$; so together, all the distance needed to realise the event $\cluster_u$ is covered.
If specifying a resample operation $\omega \mapsto \omega'$ is a matter of gathering the data of~$\omega$ into a list, and drawing a line to delimit retained from lost data, a natural first effort would be to record the
downcluster and upcluster and then draw the line, with a view to jettisoning the wasteful midpiece in~$N$ by resampling.
What exactly is the lost data in this proposal? To answer this, note that the proposed retained data is constituted of a certain collection of edges in the half-spaces $\downspace_h$ and $\upspace_h$, comprising the down- and upclusters, which are declared open; alongside edges neighbouring these clusters in the respective half-spaces, declared closed in order to determine that the two clusters are no larger. In this sense, the lost data is the status of every other edge, those whose status have not been thus declared.

However, in the terms of Section~\ref{s.resampling-overview}, this specification
is too rigid: it lies too close to the identity $\omega\mapsto\omega$ and too
far from independent resampling from $\Q$.  The upcluster base and downcluster
top are trapped in place at a high separation, making any connecting midpiece
unnecessarily costly. What is needed is a more nuanced list specification,
which will enable an argument gauging the rarity of the $N$ picture in the
tempting way, permitting us to shift the upcluster horizontally so that its
base coincides with the head of the downcluster; which is to say, so that the
two crosses coincide in Figure~\ref{f.Nslide}, whose right image shows the
shifted picture. There is no expensive surface enclosing a midpiece and so
the right picture seems of the two to be much more probable. In this resample,
the law $\Q$ of the input may naturally be taken to be
$\Q=\PP\bigl(\,\cdot\mid\cluster_u\bigr)$.
The retained data is the downcluster; and also the upcluster, but only up to a horizontal shift, so that this data includes only the equivalence class of the upcluster under identification by translations valued in~ $\Z^{d-1} \times \{ 0\}$.
The line is drawn here, with all further data lost.
The lost data may be specified explicitly by recording a vector $V \in \Z^{d-1}$ that specifies the value of the shift that determines the upcluster from its shift-equivalence class; when this value is known, the available data is the collection of declared open and closed edges in the two half spaces (as in the first discussed resample), which means that lost data is filled out by the status of the further, undeclared, edges throughout the space.

In the reconstruction step in this slide operation,  the lost data is randomly resampled according to the law of $\PP$ given $\cluster_u$ and also given the two pieces of retained data. In the depicted case, this reconstruction job involves two tetris-like pieces, the lower one held in place, and the upper one resting on the level-$h$ hyperplane, and capable of gliding over it, with any given spot being chosen with
what we may call the {\em connection chance}, namely with
 a weight proportional to the probability that, were the rest of space to be completed according to the ordinary percolation law, a suitable (and perhaps trivial) finite midpiece would form to connect the two pieces and reconstitute a cluster~$\mc C$ that realises $\cluster_u$.

When the upcluster is shifted so that the two crosses coincide, there is unit connection chance; whereas the connection chance for distant displacements may be expected to be of rapid decay.

This form of resample shows some promise, but there are objections. One  problem is that, if the slice height is such that the downcluster has more than one vertex at the height, no horizontal shift of the upcluster will lead to a renewal point at that height. To overcome this difficulty, we may simply cut the downcluster off from above at height $h$, while cutting the upcluster from below at height $h+2$; then, if we shift the upcluster so that its base hovers two units above the head of the downcluster, a resampling in which the two intervening vertical edges are declared open will, alongside a suitable set of other edges being declared closed, result in a resampled cluster~$\mc{C}_u$ that has a renewal point at height $h+1$; we say that a {\em simple join} occurs when the resampling acts in this way.

For exposition, we have so far considered the axial case $\vecell=e_d$.  For general $\vecell\in S^{d-1}$, a second difficulty arises.  When the cluster is sliced at a given height, the possible locations at which the slide may place the upcluster can be sparse, because the corresponding hyperplane may meet~$\Z^d$ only on a widely spaced mesh, or even at a single point.  The bounded offset~$3r$ in cluster height addresses this difficulty: it allows the slide to move the upcluster base through a bounded range of heights, thereby making a sufficiently rich set of shifts available.

\section{Cardinality upper bounds facilitate the slide}
These heuristics appear to offer in outline an argument that, when the input has a capital~$N$ form, the slide output will have a renewal point near the slide height with positive probability.
Some caution is needed, however. The downcluster and upcluster---the lower and upper halves of the left and right vertical struts in the figure~$N$---could in principle broaden out significantly at least at some heights either below or above the resample height~$h$. By so doing, they may form crevices that nearly enclose substantial regions in undeclared edges. The midpiece that forms in the reconstruction step could pass long distances through such nearly protected territory. Such a possibility presents a challenge to the notion represented by the TipToe mantra:
a midpiece could grow to a substantial size, so that it is unlikely to instantly perish;
by edging along the outside of voluminous downcluster or upcluster structure, it may be also shielded from open escape in its passage.

Such a scenario raises the prospect that for many potential shifts, the connection chance may not be small.
The favourable outcome of the slide operation is that the vector associated to the simple join event is selected, with the requisite finite energy decisions then being made.
This shift has a unit-order connection chance; but this is insufficient to ensure that the favourable outcome has positive probability. Indeed,
the selected shift  has probability given by the normalised connection chance;
 the connection-chance sum (or partition function)  appearing in the denominator must also be bounded above.

\begin{figure}[htbp]
\centering
\includegraphics[width=0.75\textwidth]{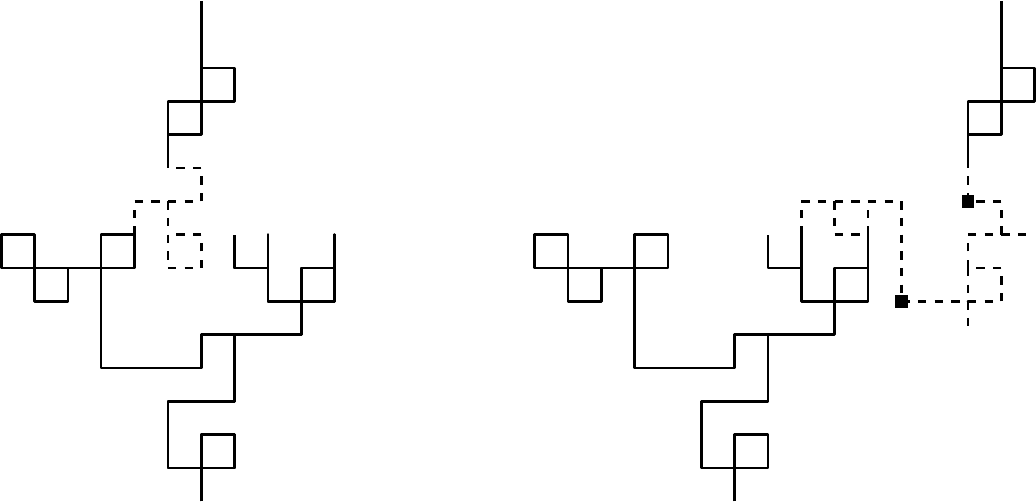}
\caption{The slide resample pushes the upcluster some way to the right in this depiction. The downcluster has some presence near the slide height and this permits the midpiece in the right sketch to travel for a while around the downcluster, but the slide is big enough to force the midpiece to make a journey in virgin territory between the two square vertices. The disjoint surface that encloses this part of the midpiece reflects the low connection chance (and the high connection cost!) of the proposed slide location.}\label{f.slideschematic}
\end{figure}

The remedy for this difficulty is illustrated by
Figure~\ref{f.slideschematic}, which shows a slide resample with a gap of two at the slide height. The illustration is schematic, depicting a two-dimensional example for simplicity.
Simple join would occur were the upcluster on the left resampled by a shift to the left by one unit, with two open vertical midpiece edges running down from the resulting upcluster base, and suitable neighbouring edges being closed.
The resampling on the right shows a dashed midpiece whose surface (which is not depicted)
shares several facets with that of the downcluster. That the downcluster has several routes up to around level~$h$ broadens the ways in which the midpiece and downcluster surfaces may coincide in part. Note however between the two square vertices runs a midpiece section whose surface is disjoint from those of the upcluster and the downcluster. In general, we may draw a path in the midpiece from the downcluster head to the upcluster base, and take the subpath between its last visit to a neighbour of the downcluster and its first visit to a neighbour of the upcluster.  This subpath runs through virgin territory, where a connection is costly: we will identify a path, called the \hfff{tether}{\emph{Tether}}, through a suitable collection of coarse-grained boxes.  The exponential cost of its length holds the shifted upcluster close to the downcluster beneath it.  For a given proposed upcluster shift, the connection chance is of order $e^{-c\kappa}$, where $\kappa$ is the minimum Tether length over the viable midpieces for that shift.

If the downcluster and upcluster cardinalities are bounded above suitably in a neighbourhood of the slice, then the number of shift locations for which the $\kappa$-value is bounded above is limited. The $\kappa$-value grows linearly in the distance away from such special sites. This means that a given shift that may lead to simple join and a renewal point (such as one unit to the left in Figure~\ref{f.slideschematic}(left))
has a connection chance that is bounded below, while the sum over all shifts of the connection chance is bounded above, with a bounded population of sites making a unit-order contribution and others making one that is exponentially small as a function of distance between downcluster and shifted upcluster.

We will infer bounds on cluster cardinality in typical height bands from Proposition~\ref{p.nobigcluster}, and argue consequently that downcluster and upcluster cardinalities at heights~$h'$ close to the slide height~$h$ typically grow at a suitably low rate as $\vert h' - h \vert$ rises.
The slide resample may then operate with an input that typically satisfies a geometric condition under which the slide produces renewal with positive conditional probability.
This typical condition is a suitably specified collection of upper bounds on the cardinality of the intersection of the downcluster and the upcluster with bands at variable displacement from the resample height. Such `fine' conditions serve to secure the connection chance sum or partition function upper bound and to deliver the positive renewal density  Theorem~\ref{t.intermediate-renewal}(1), the principal direct consequence of single slide analysis.

In the next chapter, we present single-slide results and tools.
When we later iterate the slide to prove Theorem~\ref{t.intermediate-renewal}(2), we apply order~$u$ resamples on an increasing positive-density schedule.  Each slide offers a positive probability of producing renewal, but the slides interact:
when an upcluster shifts in a slide, one midpiece exits and another joins the cluster, and these changes may degrade the fine geometry needed for  production later along the schedule, or destroy a renewal point already created.  Controlling these effects is the substance of the iterated-slide analysis, and ultimately leads to the main renewal result, Theorem~\ref{t.renewal}.

\chapter{The single slide operation and positive renewal density}
\label{s.single-slide-positive-density}

This chapter defines the single slide operation rigorously and proves Theorem~\ref{t.intermediate-renewal}(1), while also recording resampling tools needed when the operation is iterated.  We begin by proving the global surface bound, Proposition~\ref{p.nobigcluster}, which supplies the typical cardinality control used by the slide analysis.

\section{Coarse-grained percolation}\label{s.coarse}

This section introduces the required coarse graining and proves the global surface bound, Proposition~\ref{p.nobigcluster}, on the size of the coarse
boundary of a lengthy finite cluster.

The argument is local.  A box adjacent to the coarse exterior boundary must
exhibit one of two failures of the usual supercritical connectivity; such a box
will be called \emph{bad}.  The Antal--Pisztora estimate makes the bad-box field
subcritical when the box scale is large.  A large coarse boundary therefore
forces a large subcritical box cluster, giving the exponential estimate in
Proposition~\ref{p.nobigcluster}.

Two distinct vertices $v,w \in \Z^d$ are {\em $*$-neighbours} if  $w_i - v_i \in \{-1,0,1\}$ for $i \in \intint{d}$. A $*$-path is a $\Z^d$-valued path whose pairs of consecutive elements are $*$-neighbours.
The $*$-distance of $v,w \in \Z^d$ equals $n - 1$, where $n$ is the number of elements in the shortest $*$-path that connects $v$ and $w$.

Let $K \in \N$ be a given {\em even} integer.
For $v\in\Z^d$, the \hfff{k-box}{$K$-box $B(v)$} with central vertex $2Kv$ is the
subgraph of $(\Z^d,\sim)$ induced by the vertex set
$2Kv+\llbracket-K,K\rrbracket^d$.
For $A$ a subgraph of $(\Z^d,\sim)$, $\overline{A}$ will denote the set of $K$-boxes that contain an edge in $A$.

The family $\{B(v):v\in\Z^d\}$ is indexed by $\Z^d$; two $K$-boxes are
$*$-neighbours when their indices are $*$-neighbours.

An unqualified path in a lattice makes nearest-neighbour steps.  Likewise,
a path in a box lattice uses nearest-neighbour adjacency of its index lattice;
a $*$-path instead uses $*$-adjacency.
\begin{definition}\label{d.exterior}
Let \(\mathscr C\) be a collection of \(K\)-boxes.  Its coarse exterior
boundary \(\partial_{\rm ext}\mathscr C\) is the set of \(K\)-boxes \(B\)
such that \(B\) is a \(*\)-neighbour of an element of \(\mathscr C\) and
there exists an infinite path in \(K\)-boxes, starting at \(B\), that is
disjoint from \(\mathscr C\).  For a subgraph \(G\) of \((\Z^d,\sim)\), this
definition gives the notation \(\partial_{\rm ext}\overline G\).
\end{definition}

\begin{lemma}[Exterior boundary connectivity]
\label{l.extbdy-halfspace}
If \(\mathscr C\) is a finite \(*\)-connected collection of \(K\)-boxes,
then its coarse exterior boundary is nearest-neighbour connected.
\end{lemma}

\noindent\textbf{Proof.} Through the box index, this is the classical nearest-neighbour connectivity
of the exterior boundary of a finite \(*\)-connected subset of \(\Z^d\), by
Kesten~\cite[Lemma~2.23]{Kesten1986}. \qed

Here is the promised result.

\begin{proposition}\label{p.nobigcluster}
For $d \geq 2$, there exist positive constants $C$, $c$ and $D_0$ such that, when $u\geq K$ and $D \geq D_0$,
 $$
 \PP \Big(  \big\vert \partial_{\rm ext} \overline{\cset} \big\vert  \geq Du   \, \Big\vert \, \cluster_{[u,u+3r]} \Big)\, \leq \,
 C \exp \big\{ - c D u \big\}\, .
$$
The constants may be chosen uniformly in~$\vecell\in S^{d-1}$.
\end{proposition}

\subsection{Box percolation}
As tools we now present will serve to show, the above proposition will be derived by arguing that the $K$-box elements of  $\partial_{\rm ext} \overline{\cset}$ are unusual in a way that can be locally diagnosed.

The {\em boxlet} ${\rm boxlet}(v)$ and the {\em big box} $\boxbig(v)$
with central vertex $2Kv$ are the induced subgraphs with respective vertex sets
$2Kv+\llbracket-K/2,K/2\rrbracket^d$ and
$2Kv+\llbracket-5K,5K\rrbracket^d$.

For $G$ a subgraph of $(\Z^d,\sim)$,
we write $\partial G = \big\{ v \in V(G): \exists \, w \in \Z^d \setminus V(G) \, , \, v \sim w \big\}$.

Fix $v\in\Z^d$.  For the percolation configuration
sampled under $\PP$, the $K$-box $B(v)$ is called {\em bad} if at least one of
the following two random events occurs:
\begin{enumerate}
\item there exists no open path from ${\rm boxlet}(v)$ to $\partial B(v)$---we say that `${\rm boxlet}(v)$ is isolated'; or
\item there exists an open path $P$ from  ${\rm boxlet}(v)$ to $\partial B(v)$ (so that $P$ lies in $B(v)$ and thus also in $\boxbig(v)$), and an open path $Q$, also contained in $\boxbig(v)$, of Euclidean diameter at least~$K/2$ that starts in a $K$-box that $*$-neighbours $B(v)$, such that $P$ and $Q$
lie in different open connected components of the $\PP$-percolation restricted
to~$\boxbig(v)$.
\end{enumerate}
 \begin{lemma}\label{l.bdbad}
 Let \(\cset\) be the finite open cluster of the origin and suppose that its
 Euclidean diameter is at least~\(K\).  Let $B = B(v)$ be a $K$-box, with
 $v\in\Z^d$, and let $B'$ be a $K$-box $*$-neighbour of $B$. If
 $B\notin\overline{\cset}$ and $B'\in\overline{\cset}$, then $B$ is bad.
 \end{lemma}
 \noindent{\bf Proof.} We will show that if the first condition that would qualify $B$ as bad is not met, then the second one is.  Suppose then that ${\rm boxlet}(v)$ is not isolated, and let $P$ denote an open path from this set to $\partial B$. Let $x$ be a vertex in both $\cset$ and $B'$.  Since the Euclidean diameter of $\cset$ is at least~$K$,
 $\max \big\{ d(x,y): y \in V(\cset) \big\} \geq K/2$ (where $d$ is Euclidean distance). Follow an open path in $\cset$ from $x$ to such a vertex and stop when it first exits the Euclidean ball of radius $K/2$ about~$x$. The resulting path~$Q$ has diameter at least $K/2$, and every one of its vertices lies within distance $K/2+1$ of~$x$.  Thus
 $Q$ lies in $\boxbig(v)$, because $x$ lies in this big box at distance at least $2K$ from its boundary.  No vertex of $P$ lies in $\cset$: otherwise the open path $P$ would belong to the open cluster $\cset$, so that one of its edges, which lies in $B$, would place $B$ in $\overline{\cset}$. Hence there is no open path from $P$ to $Q$ in $\boxbig(v)$.
 We have shown what we sought:  if ${\rm boxlet}(v)$ is not isolated, then $B$ is bad via the second condition specifying this notion.  This completes the proof of  Lemma~\ref{l.bdbad}. \qed

In the next result, we again rely on the indexing of $K$-boxes by $\Z^d$.
\begin{lemma}\label{l.percodominate}
Define a random process on $\Z^d$ by declaring that $v \in \Z^d$ is bad if the $K$-box $B(v)$ is bad. For any $\e >0$, there exists $K \in \N$ such that the process of bad indices---a random subset of~$\Z^d$---is stochastically dominated by site percolation of parameter~$\e$.
\end{lemma}
\noindent{\bf Proof.} For a given $v\in\Z^d$, the big box $\boxbig(v)$ can be contained in the set $B'_{\bf i}(N)$ specified in \cite[Equation~$(2.7)$]{AntalPisztora}, where the parameter $N$ in this reference is set to be minimal to permit this containment.
It is readily verified that, if $B(v)$ is bad, then the complement of the event $R_{\bf i}^{(N)}$ in~\cite[Equation~$(2.9)$]{AntalPisztora} occurs. It is~\cite[Equation~(2.24)]{AntalPisztora} which establishes the rarity of this complementary event. Indeed, this bound shows that any given index in $\Z^d$ is bad with probability at most $\e$, where $\e > 0$ tends to zero in the limit of high $K$. The process of bad indices is $5$-dependent in the sense of~\cite{LSS}. Thus Lemma~\ref{l.percodominate} follows from \cite[Theorem~$0.0$(i)]{LSS}. \qed

To compare the bad-box field under $\PP$ with its law after conditioning on
edges of a subgraph $A$, we simply exclude the two-box neighbourhood of $A$.  The following
notation records such neighbourhoods.

\begin{definition}\label{d.plustwominustwo}
Let $A$ denote a subgraph of $(\Z^d,\sim)$.
For $s\in\N$, we denote by $\overline{A}[s]$ the set of $K$-boxes whose
$*$-distance from $\overline{A}$ is at most~$s$.  We write
$\overline{A}[-s]$ for the set of $K$-boxes in $\overline{A}$ whose
$*$-distance from the complement of $\overline{A}$ is at least $s+1$.
Thus any $*$-path of $K$-boxes from $\overline{A}[-s]$ to the complement of
$\overline{A}$ contains at least $s$ boxes in
$\overline{A}\setminus\overline{A}[-s]$.
\end{definition}
\begin{lemma}\label{l.condbadperc}
Let $A$ be a subgraph of $(\Z^d,\sim)$. Let $\mathsf{E}$ denote a positive-probability percolation event that is measurable with respect to the edges in $A$.
Let $\badcoll$ denote the random set of bad $K$-boxes that do not lie in~$\overline{A}[2]$. Then the distribution of $\badcoll$ is equal under the laws $\PP$ and $\PP ( \cdot \vert \mathsf{E} )$.
\end{lemma}
\noindent{\bf Proof.}
For each \(K\)-box \(B(v)\notin\overline A[2]\), whether \(B(v)\) is bad is determined by the percolation configuration in \(\boxbig(v)\).  Every edge of \(\boxbig(v)\) is contained in a \(K\)-box at \( * \)-distance at most two from \(B(v)\).  Hence \(\boxbig(v)\) contains no edge of \(A\), since otherwise some \(K\)-box in \(\overline A\) would lie at \( * \)-distance at most two from \(B(v)\), contrary to \(B(v)\notin\overline A[2]\).  Thus \(\badcoll\) is determined by edge states outside \(A\), which under \(\PP\) are independent of the edge states in~\(A\).  Conditioning on the event \(\mathsf E\) therefore leaves the law of \(\badcoll\) unchanged. \qed

\subsection{Bounding the cardinality of the cluster boundary}

Here we prove Proposition~\ref{p.nobigcluster}.

\begin{lemma}[Exterior boundary of a large finite cluster]
\label{l.large-cluster-boundary}
Suppose that $\cset$, the open cluster of the origin, satisfies $\infty > \diam(\cset)\ge K$.  Put $\mathcal B=\partial_{\rm ext}\overline{\cset}$
and $n=|\mathcal B|$.  Then $\mathcal B$ is connected,
every box in $\mathcal B$ is bad, and, for some $C=C(d,K)<\infty$,
$\mathcal B$ contains a box whose index lies within distance $C n^{1/(d-1)}$ of the origin.
\end{lemma}
\noindent{\bf Proof.} The box family $\overline{\cset}$ is finite and
$*$-connected.  Lemma~\ref{l.extbdy-halfspace} therefore gives the
connectivity of~$\mathcal B$, while Lemma~\ref{l.bdbad} shows that every
box in $\mathcal B$ is bad.  Since $0\in V(\cset)$, if the boundary
were at index distance greater than~$R$ from the origin, it would separate
the origin from infinity outside the index ball of radius~$R$ and so
contain at least $cR^{d-1}$ boxes.  Thus $n\geq cR^{d-1}$, and the
asserted distance bound follows. \qed

\begin{lemma}\label{l.bigwell}
For $d\geq2$, there exist positive constants $C$ and~$c$ such that, for
$m\geq1$,
\[
 \PP\Big(\cset\text{ is finite},\ \diam(\cset)\geq K,\
      |\partial_{\rm ext}\overline{\cset}|\geq m\Big)
 \leq Ce^{-cm}\, .
\]
In particular, for $u\geq K$ and $D\geq1$,
\[
 \PP\Big(\cluster_{[u,u+3r]},\
      |\partial_{\rm ext}\overline{\cset}|\geq Du\Big)
 \leq Ce^{-cDu}\, .
\]
\end{lemma}
\noindent{\bf Proof.} Suppose that the cluster~$\cset$ is
finite, has Euclidean diameter at least~$K$, and has an exterior boundary
of cardinality~$n$.  Lemma~\ref{l.large-cluster-boundary} supplies a
connected set of $n$ bad boxes meeting the index ball
of radius~$C n^{1/(d-1)}$ about the origin.

Choose $K$ so that the parameter~$\e$ of the dominating site
percolation in Lemma~\ref{l.percodominate} is subcritical.  The event
under consideration therefore entails that Bernoulli site percolation of parameter~$\e$
contains a connected bad set of at least $n$ sites meeting the index ball
of radius~$C n^{1/(d-1)}$.  By the exponential tail for the size of a subcritical
site-percolation cluster~\cite{AizenmanBarsky}, this has probability at most
$C n^{d/(d-1)} e^{-cn}$.  Summing over $n\geq m$ and reducing~$c$ gives the first
assertion.  On $\cluster_{[u,u+3r]}$ with $u\geq K$, the cluster is finite
and has Euclidean diameter at least~$u\geq K$; taking $m=Du$ gives the
second. \qed

\noindent{\bf Proof of Proposition~\ref{p.nobigcluster}.} That there exists $c > 0$ such that $\PP(\cluster_{[u,u+3r]}) > c^u$
is demonstrated by considering an explicit cluster, such as a concatenation of $\lceil u \cdot (e_d \cvl)^{-1} \rceil$ $e_d$-directed edges rising from~$0$.  Proposition~\ref{p.nobigcluster} then follows from Lemma~\ref{l.bigwell}
by making a suitably high choice for the lower bound $D_0$ on $D$.  The estimate in Lemma~\ref{l.bigwell} is direction-free, while the explicit-cluster lower bound is uniform in~$\vecell$ because $e_d \cvl \geq d^{-1/2}$ by the axes ordering convention~(\ref{e.aoc}).  Thus $C,c$ and~$D_0$ may all be chosen uniformly in~$\vecell\in S^{d-1}$. \qed

\subsection{Latitude and isoperimetry for $K$ boxes}

It is useful to know that when an object such as the downcluster is large in a given height interval, the enclosing surface is also large there.
Here we state and prove Proposition~\ref{p.dn.gen} to this effect, via Lemma~\ref{l.lift}.

\begin{definition}\label{d.boxlatitude}
Recall that a $K$-box takes the form $B = x + \llbracket -K,K \rrbracket^d$ with $x \in 2K \cdot \Z^d$. The latitude of $B$ equals $\bigl \lfloor \frac{x \cvl}{2K e_d \cvl} \bigr  \rfloor$.
\end{definition}
 A $K$-box analogue of height, latitude measures box location relative to the hyperplane $z\cvl=0$ in integer units.

 For $k \in \Z$ and $m \in \N$, let $\overline{N}_{k,k+m}$ denote the set of $K$-boxes whose latitude  belongs to $\llbracket k,k+m \rrbracket$.

\begin{proposition}\label{p.dn.gen}
Let $d\geq3$. For some $k \in \Z$, suppose that $\overline{A}$ and $\overline{B}$
are finite subsets of $\nkm$
such that $\overline{B}$ separates $\overline{A}$
from infinity in $\nkm$. Then
$$
 \vert \overline{A} \vert\, \leq \, (m+1) \vert \overline{B} \vert^{\tfrac{d-1}{d-2}}\, .
$$
\end{proposition}

Writing $\mc{B}$ for the set of $K$-boxes, we introduce the projection $P:\mc{B} \to \Z^{d-1}$
that sends $2Ku + \llbracket -K,K \rrbracket^d$ to $(u_1,\cdots,u_{d-1})$: namely, $P$ sends any $K$-box to the projection on the first $d-1$ coordinates of its natural index in $\Z^d$.
\begin{lemma}\label{l.lift}
Suppose that $m \geq 1$. Let $\overline{C} \subset \overline{N}_{k,k+m}$. Suppose that $x \in \Z^{d-1}$ is not separated from infinity by $P(\overline{C})$ in~$\Z^{d-1}$.
Then nor is any $B \in \overline{N}_{k,k+m}$ with $P(B) = x$
separated from infinity by~$\overline{C}$ in $\nkm$.
\end{lemma}
\noindent{\bf Proof.} Let $y$ and $z$ be nearest neighbours in $\Z^{d-1}$, with $z = y + e_i$ for $i \in \intint{d-1}$.
Let $\mc{B}[y]$ denote the set of $K$-boxes in $\nkm$ that project to $y$ under $P$.
Then $\vert \mc{B}[y] \vert$ equals $m +1$; label the lowest and highest elements of $\mc{B}[y]$ in the form
 $B_y$ and $B'_y = B_y + 2K m e_d$.
Similarly for~$z$. We {\em claim} that some element of $\mc{B}[z]$ is obtained from some element of $\mc{B}[y]$ by translation by the vector $2K e_i$.
Indeed, if this were not the case, then, since $m \geq 1$, either $B_z$ would be displaced from $B'_y + 2K e_i$ by a positive multiple of $2K e_d$; or $B'_z$ would be displaced from $B_y + 2K e_i$
by a negative multiple of $2K e_d$. In the former case, $B'_y + 2K e_i$ would have latitude less than $k$; in the latter, $B_y + 2K e_i$ would have latitude more than $k+m$. However, $0 \leq e_i \cvl \leq e_d \cvl$
by the axes ordering convention; and thus the difference in latitude between boxes $B$ and $\hat{B}$ satisfying $\hat{B} = B + 2K e_i$ is seen to be at most one in absolute value.
Since $B'_y$ has latitude $k+m$ and $B_y$ has latitude $k$, we have found a contradiction and proved the claim.

Whenever elements of $\nkm$ respectively project to nearest neighbours
$y,z\in\Z^{d-1}$, the claim permits us to connect them by a path of
$K$-boxes in $\nkm$ whose elements project to one or other of $y$ and~$z$.
Thus paths in $\Z^{d-1}$ may be lifted to paths in $\nkm$.

Indeed, since $x\in\Z^{d-1}$ is not separated from infinity by
$P(\overline C)$, there exists a path $Q$ from $x$ to infinity in
$\Z^{d-1}$ that avoids $P(\overline C)$.  Choosing an element of $\nkm$
above each successive vertex of~$Q$ and interpolating consecutive choices
by the preceding local construction yields a path to infinity in $\nkm$
that avoids $\overline C$.  This proves the lemma. \qed

\noindent{\bf Proof of Proposition~\ref{p.dn.gen}.}
The first assertion of \cite[Lemma~$2.3$]{GLA}
is a simple consequence of the Loomis-Whitney inequality~\cite{LoomisWhitney}:
if $R \subseteq \Z^{d-1}$
denotes the collection of sites that a finite set $W \subseteq \Z^{d-1}$
separates from
infinity, then $\vert R \vert \leq \vert W \vert^{(d-1)/(d-2)}$.

Since $\overline{A}$ is separated from infinity in $\nkm$ by $\overline{B}$,
 Lemma~\ref{l.lift}
 implies that $P(\overline A)$ is separated from infinity in $\Z^{d-1}$ by
 $P(\overline B)$.  Taking $R = P( \overline{A})$ in the recalled bound, we  find that
 $\big\vert P ( \overline{A} ) \big\vert \leq \big\vert P ( \overline{B} ) \big\vert^{\tfrac{d-1}{d-2}}$.
 Using $\vert\overline{A} \vert \leq (m +1)    \big\vert P ( \overline{A} ) \big\vert$
 and $\vert P (\overline{B}) \vert \leq \vert \overline{B} \vert$, we find that
 $\vert \overline{A} \vert \, \leq \, (m +1) \big\vert \overline{B}  \big\vert^{\tfrac{d-1}{d-2}}$
 and thus obtain Proposition~\ref{p.dn.gen}. \qed

\subsection{Height intervals and box latitudes}

The next result will permit us to interpolate between box latitudes and heights in the underlying lattice.

\begin{lemma}\label{l.latint}
For $q \in \R$, let $\Gamma_q$ denote the set of latitudes of $K$-boxes that contain a vertex with $v \cvl \in [q,q+1]$. Then
\[
 \Gamma_q\subseteq\llbracket a_q,b_q\rrbracket\, ,
 \qquad
 a_q\, = \, \left\lfloor\frac{q-Kd^{1/2}}{2Kr}\right\rfloor\, ,
 \quad
 b_q\, = \, \left\lfloor\frac{q+Kd^{1/2}+1}{2Kr}\right\rfloor\, .
\]
In particular, $|\Gamma_q|$ is bounded uniformly in~$q$ and~$\vecell$, and there is
$L_*=L_*(K,d)<\infty$, independent of~$\vecell$, such that $|q-q'|>L_*$ implies that
$\Gamma_q$ and $\Gamma_{q'}$ are at distance at least two.
\end{lemma}
\noindent{\bf Proof.} Suppose that the $K$-box $B = 2Ku + \llbracket -K,K \rrbracket^d$ contains $v$. Then $\vert v \cvl - 2K u \cvl \vert \leq K d^{1/2}$, so that $v \cvl \in [q,q+1]$ implies $q - K d^{1/2} \leq 2K u \cvl \leq q + K d^{1/2} +1$.
Taking floors after division by $2Kr$ gives the asserted containment.  By the
axes ordering convention, $d^{-1/2}\leq r\leq1$.  The length of the containing
interval is therefore bounded in terms of $d$ and~$K$, uniformly in both~$q$
and~$\vecell$.  Moreover, changing $q$ by $\Delta$ shifts the unfloored
endpoints by $\Delta/(2Kr)$, whose absolute value is at least
$|\Delta|/(2K)$.  Thus a choice of $L_*$ depending only on $d$ and~$K$
separates the two containing intervals, and hence $\Gamma_q$ and
$\Gamma_{q'}$, by at least two whenever $|q-q'|>L_*$. \qed

\subsection{Height bands and external boundary}

We note a consequence of Proposition~\ref{p.dn.gen}
that offers control on the number of height bands having large intersection with $\csetuthree$
in terms of the cardinality of the bad box boundary~$\partial_{\rm ext}\overlinecsetuthree$.

\begin{lemma}[Height-band boundary estimate]\label{l.heightbandboundary}
Let \(\alpha=(d-2)/(d-1)\). Whenever \(\clusteruthree\) occurs, the
following holds for every $s\in\R$, $L\geq1$ and $M>0$. Set
\(J_h=[h+s,h+s+L]\). If
\(F_L\subseteq\N\cap[0,u]\) is the set of heights \(h\) for which
$|\{x\in V(\csetuthree):x\cdot\vecell\in J_h\}|>M$, then
 $|F_L| \leq C_L M^{-\alpha}
   \big\vert \partial_{\rm ext}\overlinecsetuthree \big\vert$,
where \(C_L<\infty\) depends only on \(L,d,K\), uniformly in~\(\vecell\).
\end{lemma}

\noindent{\bf Proof.} We first prove the assertion for unit intervals, with $L=1$. Let
$L_*=L_*(K,d)$ be supplied by Lemma~\ref{l.latint}. Choose
\(G\subseteq F_L\) so that the distance between any pair of elements of
\(G\) exceeds \(L_*\), with \(|G|\geq C^{-1}|F_L|\).

For \(h\in G\), the vertices of \(\csetuthree\) whose heights lie in
\(J_h\) are contained in \(K\)-boxes whose latitudes belong to
\(\Gamma_{h+s}\). By Lemma~\ref{l.latint}, this set has cardinality
at most a constant depending only on \(d,K\). Since each
\(K\)-box contains at most \((2K+1)^d\) vertices, there exists
\(q_h\in\Gamma_{h+s}\) such that at least \(cM\) distinct elements of
\(\overlinecsetuthree\) have latitude \(q_h\), where \(c>0\) depends
only on \(d,K\).

The choice of \(L_*\) ensures that the slabs
\(\overline N_{q_h,q_h+1}\), \(h\in G\), are pairwise disjoint. For
each \(h\in G\), the set
\(\overlinecsetuthree\cap\overline N_{q_h,q_h+1}\) is separated from
infinity in \(\overline N_{q_h,q_h+1}\).  More precisely, the separating set
may be taken to be
\(\partial_{\rm ext}\overlinecsetuthree\cap
\overline N_{q_h,q_h+1}\).  Proposition~\ref{p.dn.gen} therefore implies that
\[
   \big\vert
      \partial_{\rm ext}\overlinecsetuthree
      \cap \overline N_{q_h,q_h+1}
   \big\vert
\, \geq \, c M^\alpha\, .
\]
Summing over \(h\in G\), and using the disjointness of the selected
latitude slabs, gives
\[
   \big\vert \partial_{\rm ext}\overlinecsetuthree \big\vert
\, \geq \, c M^\alpha |G|
\, \geq \, c M^\alpha |F_L|\, ,
\]
after adjusting \(c>0\). This proves the first assertion.

Now consider $L \geq 1$ general. Cover \(J_h\) by at most
\(\lceil L\rceil+1\) unit intervals. If \(J_h\) contains more than \(M\)
vertices of \(\csetuthree\), then one of these unit intervals contains
more than \(M/(\lceil L\rceil+1)\) such vertices. Applying the
unit-interval assertion to each of the finitely many unit offsets, and
summing the resulting bounds, proves the stated estimate for \(F_L\).
\qed

\section{The single slide in a general direction}\label{s.onetime}

We now implement the slide resample of Section~\ref{s.slideintorenewal} for a general direction~$\vecell$.
This single slide analysis will yield Theorem~\ref{t.intermediate-renewal}(1), or linearly many renewal points with uniformly positive probability.

  In the axial case, we could speak of proposed slides as horizontal, and indexed clearly by a copy of~$\Z^{d-1}$.  The vector~$\vecell$ may however be non-lattice, and we have to accept that there is no natural substitute for the lattice~$\Z^{d-1}$ as an index of possible shifts.  It is inevitable that the height of the cluster head will change, albeit perhaps only by a bounded quantity, under the actions of reasonable shifts.  This consideration finds expression in the conditioning~$\clusteruthree$, where we condition on a cluster of height in the short interval~$[u,u+3r]$, rather than of the specific value~$u$; the latter choice is a luxury that is unavailable in the non-lattice case.

Let \(u,h \in (0,\infty)\) satisfy \(u>h\).  We introduce the two
half-spaces that will determine the slide decomposition at height \(h\).
The downspace is simply the lower half-space \(\downspace_h\).  The upspace is
specified only on the event \(\clusteruthree\), and is defined using an
adjusted height.

Recall that
\(r=\max\{|e_i\cdot\vecell|:i\in\intint d\}=e_d\cdot\vecell\), where the
second equality is the axes ordering convention~\eqref{e.aoc}.  When
\(\clusteruthree\) occurs, there is a finite open connected component
\(\mc C\), with base \(0\), whose height lies in \([u,u+3r]\); for
emphasis, we denote this cluster by \(\csetuthree\).  We write
$\hmax=\max\{z\cdot\vecell:z\in V(\csetuthree)\}$ for its maximum height.  Thus \(\head(\csetuthree)\) is the
lexicographically maximal vertex in \(\csetuthree\) of height \(\hmax\).
The upspace is the upper half-space \(\upspace_{\hadj}\), where the
\hfff{slide-decomposition}{adjusted height} is
\[
  \hadj\, = \, h+(\hmax-u)\, .
\]
Since \(\hmax\in[u,u+3r]\), we have \(\hadj\in[h,h+3r]\).

Throughout the remainder of Chapter~\ref{s.single-slide-positive-density}, the
single slide height is fixed at~\(h\), and we suppress it from the notation for the downcluster and
upcluster.  From Chapter~\ref{s.iterated-slide} onwards, the height is again
displayed because it varies along the slide schedule.  On \(\clusteruthree\), the
downcluster, denoted \(\downcluster\), is the connected component of the
subgraph induced by \(\csetuthree\) in the downspace \(\downspace_h\) that
contains the origin.  Thus \(v\in\Z^d\) lies in \(V(\downcluster)\) if
there exists an open path in \(\downspace_h\) from \(0\) to \(v\).  The edge set
of \(\downcluster\) is the union of the edges belonging to all such open
paths.  The top of the downcluster, \(\Top(\downwell)\), is the set of
vertices in \(\downwell\) that lie in \(\vertbdry_h\).

The upcluster, denoted \(\upcluster\), is the open connected component
containing \(\head(\csetuthree)\) of the subgraph induced by
\(\csetuthree\) in the upspace \(\upspace_{\hadj}\).  Thus \(v\in V(\upcluster)\)
if and only if an open path in \(\upspace_{\hadj}\) runs from
\(\head(\csetuthree)\) to \(v\).  The edge set of \(\upcluster\) is the
union of the edges in such paths.  Note that
\(\head(\upwell)=\head(\csetuthree)\).  The \emph{bed} of the upcluster,
\(\Bed(\upwell)\), is the set of vertices in \(\upwell\) that lie in
\(\vertbdry_{\hadj}\).

The definition of the downcluster is largely consistent with the heuristics of Section~\ref{s.slideintorenewal}, but the upcluster is specified a little differently: the hyperplane that delimits the upspace, and hence the upcluster, from below is set at a height that is $u-h$ units below that of the upcluster head.  This definition is made to face the challenge of setting up the slide resample when $\vecell$ may be non-lattice and some counterpart must be found for the space~$\Z^{d-1}$ of admissible shifts of the upcluster head.  By specifying the upcluster relative to its head, we obtain a definition that is invariant under a suitable class of shifts.

\subsection{Geometry of the slide decomposition}

The construction of strings in a slab repeatedly calls for a
path that rises from the base to just below a prescribed height without
meeting either bounding hyperplane prematurely.  The following lemma supplies
such a path; it is a one-dimensional ordering statement about the heights of
the unit steps.  Its numerical threshold is three times the maximal height
change~$r$ of one lattice step.

\begin{lemma}\label{l.heightordering}
Let $x \in \Z^d$ with $x\cvl=h$, and suppose $h \geq 3r$.  There is a
path from $0$ to $x-e_d$ all of whose vertices, other than
the initial vertex $0$, have height in $(0,h-r]$.
\end{lemma}
\noindent\textbf{Proof.} Write $x-e_d=\sum_{j=1}^{d}m_je_j$ with $m_j\in\Z$.  A path
from $0$ to $x-e_d$ is obtained by listing, in any order, $|m_j|$ steps equal
to $\operatorname{sgn}(m_j)e_j$ for each $j$; every such ordering reaches
$x-e_d$, so it suffices to order the steps so that each partial sum has height
in $(0,h-r]$.  The axes-ordering convention~\eqref{e.aoc} gives
$0\le e_j\cvl\le r$ for every~$j$.  Thus each signed step changes
the height by $\operatorname{sgn}(m_j)\,e_j\cvl\in[-r,r]$, and the heights of
all the steps sum to $(x-e_d)\cvl=h-r$.  Set $B=h-r$, so that $B\geq 2r$.

Order the steps greedily, maintaining the running height $H$, initially $0$.
If $H\leq r$, append a step of positive height: one remains, because the
heights of the unused steps sum to $B-H\geq B-r>0$ and so cannot all be
non-positive, and the new height lies in $(0,2r]\subseteq(0,B]$.  If $H>r$,
append a step of negative height when one remains, giving a new height that is
positive (as $H>r$) and at most $H\leq B$; if none remains, the unused steps
have non-negative height and sum to $B-H$, so appending them in any order keeps
the running height in $(H,B]\subseteq(0,B]$.  Steps of zero height may be
placed at any point after the first.  The first step has positive height, since
$B>0$, so every non-initial partial sum lies in $(0,B]=(0,h-r]$. \qed

\begin{lemma}\label{l.updowncluster}
Let \(H_d\) and \(H_u\) denote the heights of \(\downcluster\) and
\(\upcluster\).  Almost surely under \(\PP\) given \(\clusteruthree\),
\[
  h-r<H_d\, \le \, h
  \qquad\textrm{and}\qquad
  u-h\, \le \, H_u<u-h+r\, .
\]
Moreover, \(\head(\downcluster)\in\Top(\downcluster)\) and
\(\base(\upcluster)\in\Bed(\upcluster)\).
\end{lemma}

\noindent{\bf Proof.} Because \(0=\base(\csetuthree)\), every vertex of \(\csetuthree\) has
non-negative height.  Every vertex of \(\downspace_h\) has height at most~\(h\).
Moreover, a path in \(\csetuthree\) from \(0\) to its head must leave
\(\downspace_h\), and its last vertex in the downcluster lies in
\(\vertbdry_h\).  Lemma~\ref{l.geometrybasics}(1) therefore gives, in fact,
\(h-r<H_d\le h\), as asserted.  Since \(\head(\downcluster)\) has height
\(H_d\), the same lemma gives \(\head(\downcluster)\in\vertbdry_h\), and
hence \(\head(\downcluster)\in\Top(\downcluster)\).

For the upcluster, follow such a path in the opposite direction, from the head
towards~\(0\).  Its portion before it leaves \(\upspace_{\hadj}\) lies in
\(\upcluster\), and its last vertex there belongs to
\(\Bed(\upcluster)\); in particular, the bed is non-empty.  Put
\(b=\base(\upcluster)\).  Any vertex of \(\upspace_{\hadj}\) has height strictly
greater than \(\hadj-r\), because an incident edge whose interior meets
\(F_{\hadj}\) changes height by at most~\(r\).  Since the bed contains a vertex
of height at most~\(\hadj\), the minimum-height vertex \(b\) also has height at
most~\(\hadj\).  Lemma~\ref{l.geometrybasics}(1) thus shows that
$b\cvl\in(\hadj-r,\hadj]$ and $b\in\Bed(\upcluster)$.  The head of
\(\upcluster\) is the head of \(\csetuthree\), of height~\(\hmax\).
Consequently $H_u=\hmax-b\cvl$ and
$\hmax-\hadj\le H_u<\hmax-\hadj+r$.  Substituting \(\hadj=h+(\hmax-u)\) cancels the head-height offset and gives
\(u-h\le H_u<u-h+r\), as asserted.
\qed

We next decompose the part of the cluster between the downcluster and the
upcluster.  This middle graph will be called the midpiece.

\begin{definition}
Under \(\clusteruthree\), let \(\midpiece\) be
the graph whose edge set consists of the edges of \(\csetuthree\) that
belong to neither \(\downcluster\) nor \(\upcluster\).  Its vertex set is
the set of endpoints of these edges.
\end{definition}

The midpiece may have several connected components, and it may be empty.
If the retained pieces meet, their contact already joins the two retained
parts and no connecting midpiece edge is needed.

\begin{lemma}\label{l.midpiece}
\leavevmode
\begin{enumerate}
\item If \(V(\downwell)\cap V(\upwell)\ne\varnothing\), then
\[
 V(\downwell)\cap V(\upwell)
 \subseteq \Top(\downwell)\cap\Bed(\upwell).
\]
\item If \(V(\downwell)\cap V(\upwell)=\varnothing\), then \(\midpiece\)
contains at least one connected component whose vertex set intersects both
\(\Top(\downwell)\) and \(\Bed(\upwell)\).
\item The vertex set of any connected component of \(\midpiece\)
intersects the union of the vertex sets of \(\downwell\) and \(\upwell\).
Any vertex in this intersection lies in
\(\Top(\downwell)\cup\Bed(\upwell)\).
\end{enumerate}
\end{lemma}

\noindent{\bf Proof: (1).}
Let \(x\in V(\downwell)\cap V(\upwell)\).  Since \(x\) is a vertex of the
downspace, \(x\cvl\le h\).  Since it is a vertex of the upspace,
\(x\cvl>\hadj-r\).  As \(\hadj\ge h\), this gives
\(x\cvl\in(h-r,h]\), and hence \(x\in\vertbdry_h\) by
Lemma~\ref{l.geometrybasics}(1).  Thus \(x\in\Top(\downwell)\).  The same
bounds give \(x\cvl\in(\hadj-r,\hadj]\), so
\(x\in\vertbdry_{\hadj}\) and hence \(x\in\Bed(\upwell)\).

{\bf (2).}
Choose a path in \(\csetuthree\) from \(\base(\csetuthree)\) to
\(\head(\csetuthree)\).  It begins in \(\downwell\) and ends in
\(\upwell\).  Since the two retained vertex sets are disjoint, the subpath
from the final vertex in \(\downwell\) to the first vertex in \(\upwell\)
has at least one edge.  None of its edges belongs to either retained graph,
so all belong to \(\midpiece\).  Its endpoints lie respectively in
\(\Top(\downwell)\) and \(\Bed(\upwell)\), and the component containing the
subpath has the required property.

{\bf (3).}
Let \(H\) be a connected component of \(\midpiece\).  Since
\(\csetuthree\) is connected, a path in \(\csetuthree\) runs from any
given vertex of \(H\) to \(\base(\csetuthree)\).  When this path leaves
\(H\) for the final time, it does so at a vertex \(x\) which also belongs to
either \(\upwell\) or \(\downwell\).  Suppose first that \(x\in\upwell\).
Since \(x\in V(H)\), an edge of \(H\subseteq\midpiece\) is incident to~\(x\).
This open edge cannot belong to \(E(\upspace_{\hadj})\), since it would then
belong to the upcluster; hence it belongs to \(E(\downspace_{\hadj})\).
As \(x\in V(\upspace_{\hadj})\), we have
\(x\in\vertbdry_{\hadj}\), and thus \(x\in\Bed(\upwell)\).  The case
\(x\in\downwell\) is symmetric: an incident midpiece edge cannot belong to
\(E(\downspace_h)\), and the edge partition at height~\(h\) gives
\(x\in\vertbdry_h\), hence \(x\in\Top(\downwell)\).
\qed

We use one connecting component only when a middle connection is present.

\begin{definition}
If \(V(\downwell)\cap V(\upwell)\ne\varnothing\), put
\(\midconnect=\varnothing\).  Otherwise, by Lemma~\ref{l.midpiece}(2), let
\(V_*\) be the union of the vertex sets of the components of \(\midpiece\)
which intersect both \(\Top(\downwell)\) and \(\Bed(\upwell)\), and let
\(\midconnect\) be the component containing the lexicographically minimal
vertex of~\(V_*\).
\end{definition}

\subsection{Shift equivalence and the fine conditions}

We next prepare to shift the upcluster.  The resampling retains its shape but
not its placement, so we introduce its configuration space and then identify
configurations that differ by a translation.

\begin{definition}
The {\em upcluster configuration space} \(\ucs\) consists of all values that can be adopted by the subgraph \(\upwell\) in a percolation configuration that realises the event~\(\clusteruthree\).
\end{definition}

The height of the head of an element \(U\in\ucs\) lies in \([u,u+3r]\).  In a configuration with \(\upcluster\) taking the value~\(U\), the adjusted height
\(\hadj\) equals \(h+\bigl(\head(U)\cvl-u\bigr)\) and lies in \([h,h+3r]\).

\begin{definition}\label{d.equivalence}
Two elements \(U\) and \(U'\) of \(\ucs\) are {\em shift equivalent} if there exists \(v\in\Z^d\) such that \(U'\) is the graph obtained from \(U\) by shifting by \(v\), so that the vertex and edge sets of \(U'\) are obtained from the counterparts for \(U\) by shifting by this vector.

We will write \([\upwell]\) for the shift equivalence class of \(\upwell\). Thus \([\upwell]\) is a random variable that is defined whenever \(\clusteruthree\) occurs.
\end{definition}

In this way, a height band of width \(3r\) is made available to accommodate the shift vector under a slide, permitting a rich enough class of shifts.  The use of the adjusted height in the definition of the upspace is important here.  Equivalent elements of \(\ucs\) are equal when viewed in the adjusted coordinates associated with their heads.

Our direct aim in this chapter is
Theorem~\ref{t.intermediate-renewal}(1).  The discussion in
Section~\ref{s.slideintorenewal} identifies the local input needed at a
proposed slide height: the retained pieces must not be too populous in height
bands near the slide, so that the connection-chance sum remains controlled.
Given such input, the particular two-edge reconstruction called a simple join
in that section has uniformly positive conditional probability and produces a
renewal point.  We first formalise the two families of cardinality bounds that
enter this one-slide argument.  We then add some further regularity, not for
renewal production by a single slide, but so that the same notion of a fine
height can later be used along the iterated schedule.

\begin{definition}\label{d.fvf}
The subgraph of \((\Z^d,\sim)\) induced by the vertices whose heights lie in a
given unit interval will be called a \emph{unit height band}.

Let \(\af:\Nzero\to(0,\infty)\) be a {\em fine allowance} function.  For
\(k\in\Nzero\), let
\begin{itemize}
\item[{\bf A}${}^{-}_{\rm f}(k)$:]
the number of vertices in \(\downcluster\) with
height in \([h-k-1,h-k]\) is at most \(\af(k)\);

\item[{\bf A}${}^{+}_{\rm f}(k)$:]
the number of vertices in \(\upcluster\) with height
in \([\hadj+k-1,\hadj+k]\) is at most \(\af(k)\).
\end{itemize}
Thus \({\bf A}^{-}_{\rm f}\) is the southside profile of the retained
downcluster, measured downwards from the fixed slide height~\(h\), while
\({\bf A}^{+}_{\rm f}\) is the northside profile of the retained upcluster,
measured upwards from the adjusted height~\(\hadj\).  When the upcluster is
translated, this latter reference level moves with it.  These are precisely
the band-cardinality bounds used below to control the connection-chance sum
and obtain a positive conditional chance of the simple join.

For the later iteration, the fine event also records a thin profile some
distance north of the adjusted cut.  Beginning at height \(\hadj+b_0\), we
mark \(3i_{\rm br}\) successive, evenly spaced heights, and from each take the
short slab of width \(2K\) immediately to its north.  In south-to-north order,
these slabs are divided into three consecutive rows of \(i_{\rm br}\) slabs
each.  This arrangement is convenient in the iterated analysis, where the
slabs serve as test layers used as breakers against midpieces created by
other slides.

Fix an even integer \(i_{\rm br}\geq2\).  The breaker parameters
\(b_0,C_0,C_1\) are chosen so that
$g_{\rm br}:=C_0/(3i_{\rm br})\geq2K$.  For
\(\rho\in\{0,1,2\}\) and \(j\in\intint{i_{\rm br}}\), put
\begin{align}
\label{e.breaker-profile-band}
 a_{\rho,j}(h)
   &=\hadj+b_0+(\rho i_{\rm br}+j-1)g_{\rm br},
   &J_{\rho,j}(h)&=[a_{\rho,j}(h),a_{\rho,j}(h)+2K],\\
\label{e.breaker-enlarged-interval}
 J^\sharp_{\rho,j}(h)
   &=[a_{\rho,j}(h)-r,a_{\rho,j}(h)+2K].
\end{align}
Thus the \(a_{\rho,j}(h)\) are the marked heights, in south-to-north
order, and the indices with fixed~\(\rho\) form one of the three rows.  The
interval \(J_{\rho,j}(h)\) is the nominal height interval of the associated
short slab.  The dependence of \(a_{\rho,j}(h)\), \(J_{\rho,j}(h)\) and
\(J^\sharp_{\rho,j}(h)\) on the ambient configuration through \(\hadj\)
is suppressed from the notation.  Since membership of an edge in the slab is
determined by its upper endpoint height, the lower endpoint may lie as much as
one edge-height below~\(J_{\rho,j}(h)\).  We therefore also use the slightly
enlarged interval \(J^\sharp_{\rho,j}(h)\) when bounding the vertices of the
slab.  Since \(b_0>r\) in the eventual parameter choice, this enlarged
interval remains strictly north of the adjusted cut.

The condition \({\bf B}^{+}_{\rm f}\) requires that, for each
\(\rho\in\{0,1,2\}\), at least \(i_{\rm br}/2\) indices
\(j\in\intint{i_{\rm br}}\) satisfy
\[
 |V(\upcluster)\cap\{x:x\cvl\in J^\sharp_{\rho,j}(h)\}|\, \leq \, C_1\, .
\]
The event \hfff{fine-very-fine}{$\fine_h=\fine_h(\af,b_0,C_0,C_1)$} occurs when
$\clusteruthree$ does so alongside \({\bf A}^{-}_{\rm f}(k)\) and
\({\bf A}^{+}_{\rm f}(k)\) for every \(k\in\Nzero\), and
\({\bf B}^{+}_{\rm f}\).  The standing row length~\(i_{\rm br}\) is suppressed
from this notation.

For the iterated slide we also introduce a stronger condition on the input
configuration, before the schedule begins.  Let
\(\avf:\Nzero\to(0,\infty)\) be a {\em very fine allowance} function.  For
\(k\in\Nzero\), the conditions \({\bf A}^{-}_{\rm vf}(k)\) and
\({\bf A}^{+}_{\rm vf}(k)\) require respectively that
\begin{align*}
 \bigl|V(\csetuthree)\cap\{x:x\cvl\in[h-k-1,h-k]\}\bigr|
   &\, \leq \, \avf(k)\, ,\\
 \bigl|V(\csetuthree)\cap\{x:x\cvl\in[h+k,h+k+1]\}\bigr|
   &\, \leq \, \avf(k)\, .
\end{align*}
Thus \({\bf A}^{-}_{\rm vf}\) is the whole-cluster counterpart of
\({\bf A}^{-}_{\rm f}\). 
For the northside conditions there is also a change of coordinates:
\({\bf A}^{+}_{\rm vf}\) controls the whole input cluster in bands indexed
from~\(h\), whereas \({\bf A}^{+}_{\rm f}\) follows the retained upcluster
through the adjusted height~\(\hadj\).  Since
\(\hadj-h\in[0,3r]\), these two band systems differ by only a bounded shift;
this yields the northside part of the inclusion
\(\veryfine_h\subseteq\fine_h\) in Lemma~\ref{l.vf-implies-fine}.

The condition \({\bf B}^{+}_{\rm vf}\) uses the same enlarged intervals but
counts the whole input cluster: for each \(\rho\in\{0,1,2\}\), at least
\(i_{\rm br}/2\) indices~\(j\) satisfy
\[
 |V(\csetuthree)\cap\{x:x\cvl\in J^\sharp_{\rho,j}(h)\}|\leq C_1.
\]
The event $\veryfine_h=\veryfine_h(\avf,b_0,C_0,C_1)$ is the subset of
$\clusteruthree$ on which \({\bf A}^{-}_{\rm vf}(k)\) and
\({\bf A}^{+}_{\rm vf}(k)\) hold for every \(k\in\Nzero\), and
\({\bf B}^{+}_{\rm vf}\) holds.
\end{definition}

For Theorem~\ref{t.intermediate-renewal}(1), the very-fine event has a modest
role: in Proposition~\ref{p.veryfinelevels}, many very-fine heights will be
supplied, with Lemma~\ref{l.vf-implies-fine} turning these into heights at
which the one-slide renewal-production argument applies.  The stronger
whole-cluster bounds, the breaker profile, and the slack between the very-fine
and fine allowances are chosen in this form because they will also be needed
when slides are iterated: earlier slides may add material before the slide at a
later target height is attempted.

We now specify the allowances used throughout the paper.  Their form is fixed
here; the constants entering them will be chosen in stages as the argument
develops.
 Put
\begin{equation}\label{e.allowance-design}
   \avf(k)\, = \, A_0(1+k)^3\, ,
   \qquad
   \af(0)=A_1\, ,
   \qquad
   \af(k)\, = \, \Cfine A_0^2(1+k)^{\Bfine}\quad(k\geq1),
   \qquad \Bfine=3(d+2).
\end{equation}
Thus the fine allowance is specified separately at depth zero and at positive
depth.  The need for this distinction, and for the extra room at positive
depth, will become clear in the iterated-slide analysis.

\begin{lemma}[Very fine implies fine]
\label{l.vf-implies-fine}
There are fixed lower bounds on \(A_1/A_0\) and on \(\Cfine\), independent
of~\(\vecell\), such that, for the allowances in
\eqref{e.allowance-design}, \(\veryfine_h\subseteq\fine_h\).
\end{lemma}

\noindent{\bf Proof.} The southside assertion follows from
\(\downcluster\subseteq\csetuthree\) once \(A_1\geq A_0\) and
\(\Cfine A_0^2(1+k)^{\Bfine}\geq\avf(k)\) for \(k\geq1\).
Put \(\delta=\hmax-u\in[0,3r]\).  The interval
$[\hadj+k-1,\hadj+k]=[h+\delta+k-1,h+\delta+k]$ meets at most two members of the families
\([h-q-1,h-q]\) and \([h+q,h+q+1]\), \(q\in\Nzero\).  Their indices satisfy
\(q\leq k+3r+1\leq k+4\), using~\(r\leq1\).  A southside member can occur only for \(k=0\): for
\(k\geq1\) the left endpoint is at least \(h\).  Thus, for \(k=0\), the
very-fine band bounds give a fixed estimate \(CA_0\), with an absolute
constant~\(C\), which is absorbed by choosing \(A_1\geq CA_0\).  For
\(k\geq1\), the same adjacent-band comparison gives
\[
 \bigl|V(\upcluster)\cap
       \{x:x\cvl\in[\hadj+k-1,\hadj+k]\}\bigr|
\, \leq \, C\avf(k)\, ,
\]
and this is absorbed by \(\Cfine A_0^2(1+k)^{\Bfine}\) after a fixed lower choice of
\(\Cfine\).  Enlarge \(\Cfine\), if necessary, so that
\(A_1\leq\af(1)\); the resulting \(\af\) is non-decreasing.  The breaker
implication follows from \(\upcluster\subseteq\csetuthree\). \qed

\begin{remark}[Finite adjacent-band summation]
\label{r.new-adjacent-band-form}
Fix \(k\in\Nzero\).  Let \(J\subset\R\) be covered by at most \(m_0\) bands of the form
\([h-q-1,h-q]\) or \([h+q,h+q+1]\).  Suppose that every southside index
satisfies \(|q-k|\leq m_0\), while a northside index can occur only when
\(k,q\leq m_0\).  On \(\veryfine_h\),
\[
 |V(\csetuthree)\cap\{x:x\cvl\in J\}|
\, \leq \, C_{m_0}\avf(k)\, .
\]
Indeed, sum the applicable unit-band allowances and use
\(\avf(q)=A_0(1+q)^3\).  No widened-band allowance is imposed.
\end{remark}

The connection-chance analysis for a single slide uses only the two fine
allowance families ${\bf A}^{-}_{\rm f}$ and ${\bf A}^{+}_{\rm f}$.  For that
purpose the fine allowance need only have at most polynomial growth: for some
\(C_{\rm f},q<\infty\),
\begin{equation}\label{e.fine-polynomial-growth}
   \af(k)\, \leq \, C_{\rm f}(1+k)^q\, ,
   \qquad k\in\Nzero\, .
\end{equation}

\subsection{The slide operator and simple join}

Section~\ref{s.slideintorenewal} described the favourable reconstruction in
which a translate of the upcluster is joined to the downcluster by a two-edge
bridge and no further midpiece.  We now formulate this \emph{simple join}
exactly.  In the line--list picture of Chapter~\ref{s.heuristics}, the retained
data, above the line, is \(\bigl(\downwell_h,[\upwell_h]\bigr)\).  The placement
of the upcluster, represented by its shift vertex \(V\), and the midpiece lie
below the line and are resampled.

Take \(v=\head(\downcluster)+e_d\).

The event \hfff{simple-join}{\(\success\)} occurs when \(\clusteruthree\) does alongside the
conditions that
\begin{enumerate}
\item \(\midconnect\) has merely three vertices and two edges, the pair of
edges having the form \([v,v+e_d]\) and \([v-e_d,v]\); and
\item the graph \(\midpiece\setminus\midconnect\) is empty.
\end{enumerate}
On \(\success\), the vertex \(v=\head(\downcluster)+e_d\) will be called
the {\em midvertex}.

\begin{lemma}[Fine is retained-data measurable]
\label{l.fine-retained}
On \(\clusteruthree\), the event \(\fine_h\) is determined by the retained pair
\(\bigl([\upwell_h],\downwell_h\bigr)\).  Equivalently, two configurations in
\(\clusteruthree\) with the same retained pair either both belong to
\(\fine_h\) or both lie outside it.
\end{lemma}

\noindent{\bf Proof.} The southside conditions \({\bf A}^{-}_{\rm f}(k)\) depend only on the
retained downcluster.  For the northside conditions, suppose that a
representative \(\upwell\) is translated by a vector~\(t\).  The corresponding
maximum height and adjusted height increase by \(t\cdot\vecell\), so
\[
 (z+t)\cdot\vecell-(\hadj+t\cdot\vecell)
 =z\cdot\vecell-\hadj
 \qquad\text{for every }z\in V(\upwell).
\]
Thus every northside band count \({\bf A}^{+}_{\rm f}(k)\) is unchanged by
changing representative of \([\upwell]\).  The same identity translates each
interval \(J^\sharp_{\rho,j}(h)\) onto its counterpart for
the translated representative, without changing its upcluster occupation.
Hence \({\bf B}^{+}_{\rm f}\) is unchanged as well.  On \(\clusteruthree\),
these conditions together determine \(\fine_h\), proving the claim. \qed

We write \(\PP^{\fine_h}=\PP(\,\cdot\,\vert\fine_h)\), and use
$\PP_{[\upwell],\downwell}^{\fine_h}
=\PP^{\fine_h}(\,\cdot\,\mid[\upwell],\downwell)$ for a regular conditional law given \(([\upwell],\downwell)\).  The
configuration space and the space in which the retained pair takes values are
standard Borel spaces, so this law may be chosen as a measurable probability
kernel of the retained pair.

\begin{definition}[The single-slide operator]
\label{d.slide-operator}
For \(h\in(0,u)\), define the Markov kernel \hfff{single-slide}{\(\slide_h\)} by
\[
   \slide_h(\omega,A)\, = \,
   \begin{cases}
   \PP^{\fine_h}_{[\upwell_h(\omega)],\downwell_h(\omega)}(A)\, ,
      & \omega\in\fine_h\, ,\\
   \mathbf 1_A(\omega),&\omega\notin\fine_h\, ,
   \end{cases}
   \qquad A\in\mc B\, .
\]
A slide is called  \emph{active} or  \emph{dormant} accordingly.  Thus an active slide retains \(([\upwell_h],\downwell_h)\) and resamples
the configuration from its conditional law on \(\fine_h\), while a dormant
slide leaves the configuration unchanged.  Every active output therefore
remains in \(\fine_h\), and hence in \(\cluster_{[u,u+3r]}\).  The explicit reconstruction of this
conditional law is developed in Section~\ref{s.simplejoin}.
\end{definition}

The next result gives a lower bound on the probability that the slide
resample realises a simple join at height~\(h\).

\begin{proposition}\label{p.success}
Let the fine allowance \(\af\) satisfy~\eqref{e.fine-polynomial-growth}.
There exists \(c>0\), depending on \(\af\) but uniform in~\(\vecell\) and in
the breaker-profile constants \(b_0,C_0,C_1,i_{\rm br}\), such that, whenever
\(u\ge K\), \(h>2r\) and \(u-h>2r\),
\[
  \PP^{\fine_h}\big(\success\,\big\vert\,[\upwell],\downwell\big)
  \,\geq\, c .
\]
For the particular choice of \(\af\) in~\eqref{e.allowance-design}, write
\(c_{\rm SJ}\) for such a lower bound.  It may be chosen so that
\begin{equation}
\label{e.simplejoin-Cfine-lower}
   c_{\rm SJ}\,\geq\,c_*(2+\Cfine)^{-D_*}
\end{equation}
for fixed \(c_*>0\) and \(D_*<\infty\), uniformly in~\(\vecell\).
\end{proposition}

\begin{figure}[htbp]
\centering
\includegraphics[width=0.6\textwidth]{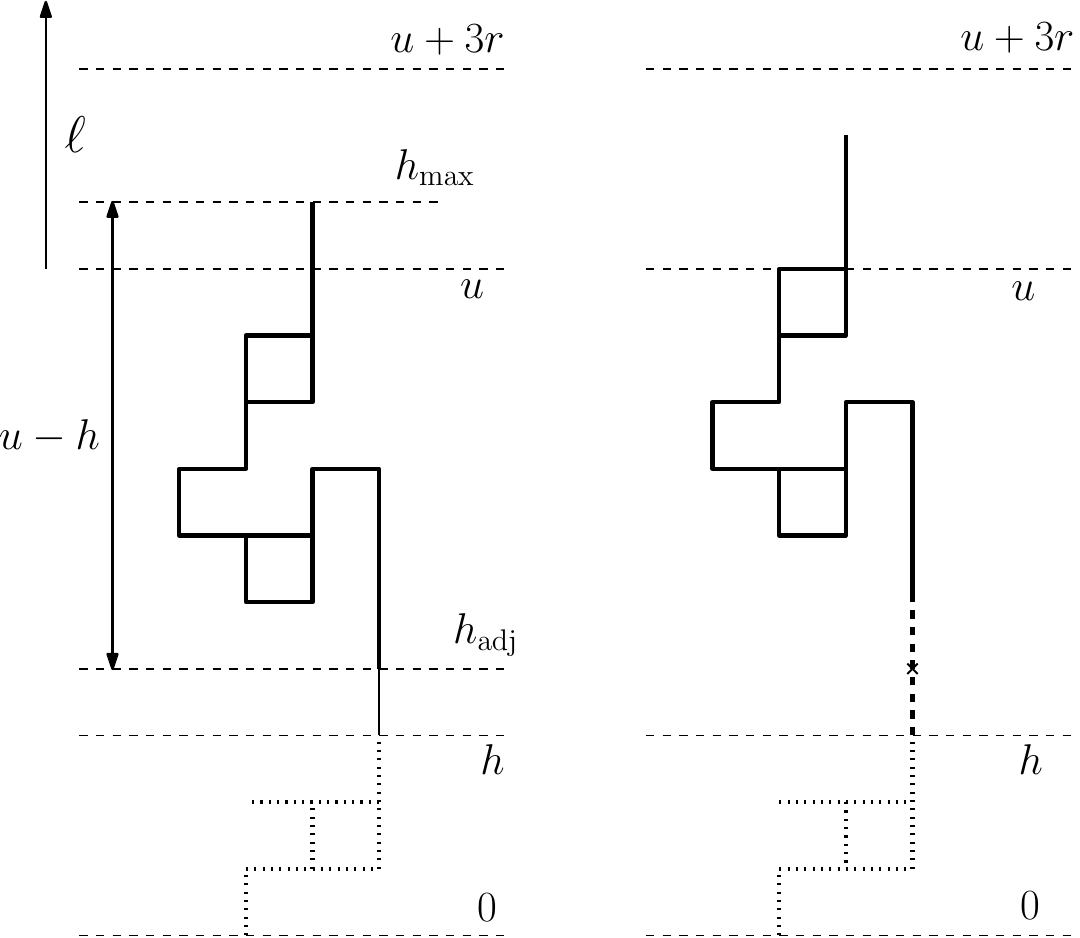}
\caption{The slide resample illustrated, with \(\success\) occurring as a
result. {\em Left}: \(\upwell\) is bold; its head lies
\(u-h=h_{\rm max}-h_{\rm adj}\) above the adjusted cut; \(\downwell\) is bold dotted.
{\em Right}: \(\upwell\) moves up by one unit, to accommodate the two bold
dashed edges that meet at the midvertex, marked with a cross.}
\label{f.resample}
\end{figure}

For \(u\ge K\), Proposition~\ref{p.success} gives an active slide a uniformly
positive conditional probability of \(\success\), while Lemma~\ref{l.success} identifies
the midvertex of every such outcome as a renewal vertex.  Figure~\ref{f.resample}
illustrates this reconstruction.

\begin{lemma}\label{l.success}
When \(\success\) occurs, the midvertex is a renewal vertex of the
cluster~\(\csetuthree\).
\end{lemma}

For \(t\in\R\), let the {\em edge boundary} \(\edgebdry_t\) be the set of
edges incident to the vertex boundary \(\vertbdry_t\) of
Section~\ref{s.basics}.

\begin{lemma}\label{l.twointernal}
Let \(v\in\Z^d\).  No edge in
\(E(\upspace_{(v+e_d)\cdot\vecell})\) is incident to
\(\vertbdry_{v\cdot\vecell}\); equivalently,
\[
   E(\upspace_{(v+e_d)\cdot\vecell})
      \cap\edgebdry_{v\cdot\vecell}\, = \, \varnothing\, .
\]
\end{lemma}

\noindent{\bf Proof.} Take \(e\in E(\upspace_{(v+e_d)\cdot\vecell})\), and denote its higher and lower endpoints by
\(w\) and \(x\).  Put \(r=e_d\cdot\vecell\).  The axes-ordering convention
gives \(w\cdot\vecell-x\cdot\vecell\leq r\), while the definition of the
upspace gives
\(w\cdot\vecell>(v+e_d)\cdot\vecell=v\cdot\vecell+r\).  Hence
\(x\cdot\vecell\geq w\cdot\vecell-r>v\cdot\vecell\).  By
Lemma~\ref{l.geometrybasics}(1), neither endpoint
lies in \(\vertbdry_{v\cdot\vecell}\). \qed

\noindent{\bf Proof of Lemma~\ref{l.success}.}
On \(\success\), the cluster \(\csetuthree\) is the concatenation of
\(\downcluster\), the two-edge middle bridge, and \(\upcluster\).  More
precisely, it fragments at \(v\) into two subgraphs whose edge sets are
the union of \([v,v+e_d]\) and \(E(\upcluster)\), and the union of
\([v-e_d,v]\) and \(E(\downcluster)\).  To prove the lemma, it is enough
to argue that \(e\in E(\upcluster)\cup E(\downcluster)\) implies
\(e\notin\edgebdry_{v\cdot\vecell}\).  For the upcluster, first note that \(v+e_d\in\Bed(\upcluster)\).
Indeed, \(\midconnect\) is the two-edge bridge with vertices
\(v-e_d,v,v+e_d\), and it meets \(\Bed(\upcluster)\).  If \(v-e_d\) were
in the bed, then the open edge \([v-e_d,v]\) would belong to
\(E(\upspace_{\hadj})\) and hence to \(E(\upcluster)\), contrary to its
membership of \(\midpiece\).  The same argument with \([v,v+e_d]\) rules
out \(v\in\Bed(\upcluster)\).  Thus \(v+e_d\in\Bed(\upcluster)\), so
\((v+e_d)\cdot\vecell\le\hadj\).  It follows that
\(E(\upcluster)\subset E(\upspace_{\hadj})\subset
E(\upspace_{(v+e_d)\cdot\vecell})\), and Lemma~\ref{l.twointernal}
applies.  For the downcluster we argue by vertex
height.  Every vertex of \(\downcluster\) has height at most
\((v-e_d)\cdot\vecell=v\cdot\vecell-r\), its head being the highest, and so
lies outside
\(\vertbdry_{v\cdot\vecell}=\{w:w\cdot\vecell\in(v\cdot\vecell-r,\,v\cdot\vecell]\}\);
hence no edge of \(\downcluster\) is incident to \(\vertbdry_{v\cdot\vecell}\).
This proves the lemma. \qed

\subsection{Positive renewal density}\label{ss.positive-renewal-density}

We now make the choices of constants in the standing allowances
\eqref{e.allowance-design} needed for the renewal-density result,
Theorem~\ref{t.intermediate-renewal}(1), while retaining the freedom required
by the later iterated analysis.  Proposition~\ref{p.success}, which gives an
active slide a uniformly positive chance to produce renewal by a simple join,
requires only that \(\af\) have at most polynomial growth.  We will next
establish, in Proposition~\ref{p.veryfinelevels}, an abundance estimate for
very-fine input heights; its proof requires
\(\sum_{k\geq0}\avf(k)^{-\alpha}<\infty\), where
\(\alpha=(d-2)/(d-1)\).  The later iterated argument also requires the fine
allowance to absorb \(\avf^{1/\alpha}\).

With the coarse scale~\(K\) fixed as in Section~\ref{s.coarse},
choose \(A_0\geq1\) sufficiently large for the very-fine allowance
part of the abundance argument below and its later sampling in
Section~\ref{ss.sampling-ivf-points}.  The depth-zero allowance \(A_1\)
is then fixed sufficiently large for Lemma~\ref{l.vf-implies-fine} and for
the deterministic host estimate in the depth-zero delivery argument of
Chapter~\ref{s.shadow-bounds-delivery-fine}.  The local breaker occupation
threshold \(C_1\) is fixed next, sufficiently large for the breaker-profile
part of the same argument and sampling step; both bounds are uniform in the
row length and its separation.  The remaining parameters---the positive-depth
coefficient \(\Cfine\), the breaker-row strength and the random-mesh gap---are
chosen simultaneously in Lemma~\ref{l.delivery-constant-compatibility}.
That lemma also enforces \(A_1\leq\af(1)\), so that \(\af\) is
non-decreasing.

The split at zero is deliberate.  The finite-energy cost of prescribing the
two-edge join depends on \(\af(0)=A_1\).  The probability of selecting this
favourable reconstruction is then normalised by the total connection chance
over possible shifts; its polynomial dependence on~\(\Cfine\) is recorded in
\eqref{e.simplejoin-Cfine-lower}.  The extra slack at positive depth serves a
different purpose in the iterated analysis, where earlier slides may add
material before a later target height is reached.  Meanwhile, since
\(3\alpha>1\) for every \(d\geq3\), the very-fine allowance is summable
to the required power.  The fixed exponent~\(\Bfine\) places the
positive-depth fine allowance well above~\(\avf^{1/\alpha}\), providing the
room needed to absorb this later accrual.

The following estimate supplies the promised abundance of input heights.

\begin{proposition}\label{p.veryfinelevels}
Let \(N\) be the number of \(h\in\N\cap[0,u]\) for which
\(\veryfine_h\) occurs.  With \(A_0\) and \(C_1\) chosen sufficiently large,
there exist \(C,c>0\), uniformly in the offset \(b_0\), separation \(C_0\),
and row length \(i_{\rm br}\), such that
\[
   \PP\Big( N < \tfrac12 u \,\Big\vert\, \clusteruthree \Big)
\, \leq \, C \exp\{-cu\}\, .
\]
The choices of \(A_0,C_1\) and the constants \(C,c\) may be made uniformly in~\(\vecell\). 
\end{proposition}
 The factor of \(u\) is taken to equal one-half since this suffices for our purpose, though it could be increased to any value strictly less than one.

\noindent{\bf Proof of Theorem~\ref{t.intermediate-renewal}(1).} We may assume
\(u\ge K\), since the theorem concerns all sufficiently large~\(u\).  For
$h \in [0,u]$, let $v = \head(\downcluster) + e_d$. Denote by
$\mathscr{R}_u(h) \subseteq \clusteruthree$ the event that $v$ is a renewal
point in $\csetuthree$.  We claim that there exists $c>0$ such that, when
$h>2r$ and $u-h>2r$,
\begin{eqnarray*}
& & \PP \big( \mathscr{R}_u(h) \big)\, \geq \,
\PP \big( \clusteruthree , \success_h \big)\, \geq \,
\PP \big( \fine_h , \success_h \big) \\
&\, = \, &
 \E \, \Big[ \, \PP_{[\upwell_h],\downwell_h}^{\fine_h} \big(  \success \big) \cdot {\bf 1}_{\fine_h} \, \Big]\, \geq \,
c \, \PP \big( \fine_h \big)\, ,
\end{eqnarray*}
using the conditional-law notation above.
Indeed, Lemma~\ref{l.success}, $\fine_h \subseteq \clusteruthree$ and Proposition~\ref{p.success} deliver the respective inequalities.

For each $h \in [0,u]$, the associated value of $\head(\downcluster)$ lies in $\vertbdry_h$. By Lemma~\ref{l.geometrybasics}(1), if two such values of $h$ differ by at least $r = e_d \cvl$, the associated vertices~$v$ cannot be equal.
Since $\vecell\in S^{d-1}$ we have $r = e_d\cvl \leq 1$, so consecutive integer heights, being one apart, differ by at least~$r$; thus for distinct $h\in\N\cap[0,u]$ the vertices~$v$ are distinct, and ascertaining the occurrence of a certain number of the events $\mathscr{R}_u(h)$ for $h \in \N \cap [0,u]$ guarantees at least that many elements of~$R(\csetuthree)$.

We may thus divide the bound derived above by $\PP(\clusteruthree)$ and sum
the result over $h \in \N \cap [0,u]$ for which $h>2r$ and $u-h>2r$, to find that
$\E\big[\vert R(\csetuthree)\vert \,\big\vert\, \clusteruthree\big]$ is at least
$c\,\E[N_{\rm f} \,\vert\, \clusteruthree]$, where $N_{\rm f}$ denotes the number of integer
$h \in \N \cap [0,u]$, with $h>2r$ and $u-h>2r$, for which $\fine_h$ occurs.
Here the stronger input condition is used only through
$\veryfine_h\subseteq\fine_h$; it is the condition needed later for the
iterated slide.  Accordingly, the count $N_{\rm f}$ is at least the very-fine
count $N$ of Proposition~\ref{p.veryfinelevels} less a bounded number of
endpoint heights.  That proposition gives
$\E[N\,\vert\,\clusteruthree]\geq c_0u$ for some $c_0>0$ and all large~$u$,
and hence $\E[N_{\rm f} \,\vert\, \clusteruthree] \geq c_1 u$ for some
$c_1 > 0$.
Whence, using the almost-sure bound
$\vert R(\csetuthree)\vert\leq 1+(u+3r)r^{-1}\leq5ur^{-1}$ for $u\geq r$,
\begin{eqnarray*}
 & & 5ur^{-1} \, \PP \big(  \vert R(\csetuthree) \vert\, \geq \, c_2 u \big\vert  \clusteruthree \big) \, + \, c_2 u \, \PP \big(   \vert R(\csetuthree) \vert < c_2 u \big\vert  \clusteruthree \big) \\
 &\, \geq \, & \E \big[   \vert R(\csetuthree) \vert  \big\vert \clusteruthree \big]\, \geq \, c c_1 u\, ,
\end{eqnarray*}
so that $\PP \big( \vert R(\csetuthree) \vert \geq c_2 u \big\vert \clusteruthree \big) \geq \tfrac r5(cc_1 - c_2)$ for any $c_2\in(0,cc_1)$; since the axes-ordering convention~\eqref{e.aoc} gives $r = e_d\cvl\geq d^{-1/2}$, this is at least $\tfrac15d^{-1/2}(cc_1-c_2)$.  Taking $c_2 = c c_1/2$ and relabelling the density and probability constants as $a$ and $c$, respectively, we obtain Theorem~\ref{t.intermediate-renewal}(1). \qed

We also record the form in which the preceding proposition and lemma will be
used later.

\begin{corollary}[Uniform renewal chance at a fine active slide]
\label{c.uniform-renewal-chance-active-slide}
There exists \(p_{\rm ren}>0\) such that, whenever \(u\ge K\), the slide at
height \(h\) is active, \(h>2r\), and \(u-h>2r\), the conditional probability
of \(\success\) is at least \(p_{\rm ren}\), uniformly over the value of
\((\downwell,[\upwell])\) and over~\(\vecell\in S^{d-1}\).  In particular, with conditional probability at
least \(p_{\rm ren}\), the resample creates a renewal vertex whose height lies
in \((h,h+r]\).
\end{corollary}

\noindent{\bf Proof.} By Proposition~\ref{p.success}, conditionally on \((\downwell,[\upwell])\)
and on \(\fine_h\), the event \(\success\) has probability at least a
positive constant.  By Lemma~\ref{l.success}, occurrence of \(\success\)
makes the midvertex a renewal vertex.  Write this vertex as~\(m\).  Since
\(m=\head(\downwell)+e_d\), Lemma~\ref{l.geometrybasics}(1) gives
\(h<m\cdot\vecell\leq h+r\).  Taking \(p_{\rm ren}\) to be the preceding
positive constant proves both assertions. \qed

\subsubsection{Abundance of very-fine heights}

The proof of Theorem~\ref{t.intermediate-renewal}(1) given above uses
Propositions~\ref{p.veryfinelevels} and~\ref{p.success}.  We prove the former
in the remainder of this section and the latter in
Section~\ref{s.simplejoin}.  For the allowance fixed above, the summability condition
$\sum_{k\geq0}\avf(k)^{-\alpha}<\infty$ holds.  The prefactor \(A_0\) and the local breaker threshold \(C_1\) will be
fixed so that the deterministic estimate below reduces failure of very-fineness
to the event controlled by the global surface bound,
Proposition~\ref{p.nobigcluster}.  The
number and separation of the breaker tests do not enter this estimate.

Proposition~\ref{p.veryfinelevels} is reduced to the following deterministic
estimate.

\begin{lemma}[Very-fine failure count]\label{l.nhigh}
There is a constant \(C<\infty\), depending only on \(d,K\), with
the following property.  Put
$Q=C(\sum_{k\geq0}\avf(k)^{-\alpha}+C_1^{-\alpha})$.  For every admissible breaker profile, uniformly in \(b_0,C_0\) and
\(i_{\rm br}\), and for every configuration that realises
\(\clusteruthree\), the number of \(h\in\N\cap[0,u]\) for which
\(\veryfine_h\) fails is at most
$Q|\partial_{\rm ext}\overlinecsetuthree|$.
\end{lemma}

Let \(D_0\) be the threshold in Proposition~\ref{p.nobigcluster}.  In making
the choices of \(A_0\) and \(C_1\) above, we require
\begin{equation}\label{e.vf-abundance-choice}
   QD_0<\frac1{8K}\, .
\end{equation}
This is compatible with the order above: first increase \(A_0\) until
$\sum_{k\geq0}\avf(k)^{-\alpha}$ is sufficiently small; after \(A_1\) has
been fixed, increase \(C_1\) until the second term in~\(Q\) is sufficiently
small.  The uniformity in Lemma~\ref{l.nhigh} is also uniform in~$\vecell$; it leaves the row length and its
separation to be chosen later.

\noindent{\bf Proof of Proposition~\ref{p.veryfinelevels}.}
On \(\clusteruthree\), the event \(N<u/2\) implies
$|\partial_{\rm ext}\overlinecsetuthree|>(2Q)^{-1}u$.  By
\eqref{e.vf-abundance-choice}, \((2Q)^{-1}>D_0\), so Proposition
\ref{p.nobigcluster}, applied with \(D=(2Q)^{-1}\), gives
\[
   \PP\bigl(N<u/2\mid\clusteruthree\bigr)
\, \leq \, C\exp\{-c(2Q)^{-1}u\}\, .
\]
Since \(Q\) is fixed, this is the asserted estimate. \qed

\noindent{\bf Proof of Lemma~\ref{l.nhigh}.}
We consider the three possible ways in which \(\veryfine_h\) may fail.

First consider the southside allowance condition \({\bf A}^{-}_{\rm vf}\).
For \(k\in\Nzero\), let \(F^-_k\) be the set of \(h\in\N\cap[0,u]\) for
which the number of vertices of \(\csetuthree\) whose height lies in
\([h-k-1,h-k]\) exceeds \(\avf(k)\). Applying
Lemma~\ref{l.heightbandboundary} with \(M=\avf(k)\) gives
\[
   |F^-_k|
\, \leq \, C\avf(k)^{-\alpha}
      \big\vert \partial_{\rm ext}\overlinecsetuthree \big\vert\, .
\]
The northside allowance condition \({\bf A}^{+}_{\rm vf}\) is
handled in the same way.  Let \(F^+_k\) be the set of
\(h\in\N\cap[0,u]\) for which the interval \([h+k,h+k+1]\) contains more
than \(\avf(k)\) vertices of \(\csetuthree\).  Lemma~\ref{l.heightbandboundary}
gives
\[
   |F^+_k|
\, \leq \, C\avf(k)^{-\alpha}
      \big\vert \partial_{\rm ext}\overlinecsetuthree \big\vert\, .
\]

It remains to consider the breaker profile.  For
\(\rho\in\{0,1,2\}\) and \(j\in\intint{i_{\rm br}}\), let
\(F^{\rm br}_{\rho,j}\) be the set of heights~\(h\) for which the interval
\(J^\sharp_{\rho,j}(h)\) contains more than~\(C_1\) vertices of the
whole cluster.  On \(\clusteruthree\), the adjusted offset \(\hmax-u\) is
independent of \(h\).  Thus, for fixed \((\rho,j)\), these intervals have fixed
width~\(2K+r\) and are translates of~\(h\) by one fixed offset.  Lemma
\ref{l.heightbandboundary}, with \(L=2K+r\) and \(M=C_1\), therefore gives
\[
 |F^{\rm br}_{\rho,j}|
 \leq C C_1^{-\alpha}
       \big\vert\partial_{\rm ext}\overlinecsetuthree\big\vert,
\]
where the constant is independent of \(\rho,j,i_{\rm br},b_0,C_0\) and~\(\vecell\): indeed $2K+r\leq2K+1$, and the constants in Lemma~\ref{l.heightbandboundary} are uniform in~$\vecell$.

Let \(F^{\rm br}\) be the set of heights at which
\({\bf B}^{+}_{\rm vf}\) fails.  At every \(h\in F^{\rm br}\), one of the
three rows has fewer than \(i_{\rm br}/2\) good indices and
hence at least \(i_{\rm br}/2\) bad ones.  Double-counting the bad incidences
and using the preceding bound gives
\[
 \frac{i_{\rm br}}2|F^{\rm br}|
 \leq \sum_{\rho=0}^2\sum_{j=1}^{i_{\rm br}}|F^{\rm br}_{\rho,j}|
 \leq 3i_{\rm br} C C_1^{-\alpha}
       \big\vert\partial_{\rm ext}\overlinecsetuthree\big\vert.
\]
Consequently
$|F^{\rm br}|\leq6C C_1^{-\alpha}
|\partial_{\rm ext}\overlinecsetuthree|$, with no dependence on the row length or on the number of test layers.

The set of \(h\in\N\cap[0,u]\) for which \(\veryfine_h\) fails is
contained in $F^{\rm br}\cup\bigcup_{k\geq0}(F^-_k\cup F^+_k)$.  Combining the three estimates proves the lemma. \qed

For later sampling from prescribed sets of candidate heights, we record the
same pointwise estimate after restriction to a deterministic subset.

\begin{corollary}\label{c.nhigh-subset}
For any deterministic \(J\subseteq\N\cap[0,u]\), on \(\clusteruthree\) the number
of \(h\in J\) for which \(\veryfine_h\) fails is at most
\(Q\,\big\vert\partial_{\rm ext}\overlinecsetuthree\big\vert\), the bound of
Lemma~\ref{l.nhigh}; equivalently, at least
\(|J|-Q\,\big\vert\partial_{\rm ext}\overlinecsetuthree\big\vert\) of the heights
in \(J\) are very fine.
\end{corollary}

\noindent{\bf Proof.} The set of \(h\in J\) at which \(\veryfine_h\) fails is contained in the set of
all \(h\in\N\cap[0,u]\) at which it fails, and Lemma~\ref{l.nhigh} bounds the
cardinality of the latter by \(Q\,\big\vert\partial_{\rm
ext}\overlinecsetuthree\big\vert\); the complementary count is the second
assertion. \qed

\section{Simple joins are often found at fine levels}\label{s.simplejoin}

We now derive Proposition~\ref{p.success}, the remaining input needed to prove Theorem~\ref{t.intermediate-renewal}(1).

Section~\ref{ss.explicit-slide-law} gives the conditional law of the slide.
Alongside it, we record product-field and stopping-domain facts used later in
the iterated-slide analysis.
The law contains a normalising sum over possible head locations~\(v\) satisfying
\(v\cdot\vecell\in[u,u+3r]\).  Each summand is the connection chance obtained after
fixing the downcluster and translating the upcluster so that its head is~\(v\).

The space of available head locations  is infinite and summability of the connection chances must be
proved.  As discussed in Chapter~\ref{s.heuristics}, a shift that places the
head far from the downcluster top near height~\(h\) forces a long Tether through
undiscovered material.  Its connection chance falls exponentially in the
Tether length.  Proposition~\ref{p.finesum} makes this comparison uniform and
shows that the normalising sum is finite.  The same estimates also show that
the resulting head-placement law carries little mass on placements of very
small individual probability.  Section~\ref{s.reduction} reduces
Proposition~\ref{p.success} to the partition-function bound, whose proof
occupies Section~\ref{s.finesum}.

\subsection{An explicit description of the conditional law in the slide resample}
\label{ss.explicit-slide-law}

Two definitions concerning the resample permit us to present
Lemma~\ref{l.probform}, which rigorously expresses the notion advanced in
Section~\ref{s.slideintorenewal} that any given shift value is selected
under resampling with probability proportional to its connection chance.

Recall the random shift-equivalence class \([\upwell]\) specified under
\(\PP\) given \(\clusteruthree\) by Definition~\ref{d.equivalence}.

\begin{definition}\label{d.displace}
Let \(v\in\Z^d\) satisfy \(v\cdot\vecell\in[u,u+3r]\).
On the event \(\clusteruthree\), we denote by \(\upwell^v\) the unique
translate of \(\upwell\) whose head equals~\(v\).  This translated graph
depends only on the shift-equivalence class \([\upwell]\).  We set
\[
   V=\head(\upwell)=\head(\csetuthree),
   \qquad V\cdot\vecell\in[u,u+3r]\, .
\]
Thus \(\upwell=\upwell^V\).
\end{definition}

\begin{definition}[Fixed-shift reconstruction]\label{d.prime}
Let \(\Omega\) be the configuration space of bond percolation on
\((\Z^d,\sim)\): a point \(\omega\in\Omega\) is a set of {\em open} edges, a
subset of \(\Edg(\Z^d)\), equivalently a point of \(\{0,1\}^{\Edg(\Z^d)}\)
through its indicator.  We take \(\mc B\) to be the product
\(\sigma\)-field---generated by the coordinate events \(\{e\text{ open}\}\),
\(e\in\Edg(\Z^d)\)---and \(\PP\) to be the product Bernoulli(\(p\)) measure,
under which distinct edges are open independently with probability~\(p\).  Bond
percolation on \((\Z^d,\sim)\) is treated as defined on the probability space
\((\Omega,\mc B,\PP)\).

We select a second copy \((\Omega',\mc B',\PP')\) of this probability
space \((\Omega,\mc B,\PP)\).  The product space
\((\Omega\times\Omega',\mc B\times\mc B',\PP\times\PP')\) describes two
independent bond percolations on \(\Z^d\) of parameter~\(p\).  The first
may be called the {\em original} percolation and the second, associated to
the prime variable, the {\em auxiliary} percolation.

Let \(\Omega''\) denote a further copy of the configuration space~\(\Omega\).  For each
\(v\in\Z^d\) with \(v\cdot\vecell\in[u,u+3r]\), we
specify a further configuration \(\omega''_v\in\Omega''\) as a function of
\((\omega,\omega')\in\Omega\times\Omega'\).

An edge is called \(\omega\)-open or \(\omega\)-closed if it is, or is not,
an element of \(\omega\); naturally the same usage may be made for
\(\omega'\) and \(\omega''_v\).  Set \(\omega''_v\) to be a definite
choice---the identically closed configuration, say---if
\(\omega\notin\clusteruthree\).  Suppose now that
\(\omega\in\clusteruthree\).  In what follows, every random variable, such
as \(\upwell^v\) or \(\downwell\), is understood to be specified by the
configuration~\(\omega\).

To specify the configuration \(\omega''_v\), we classify edges in
\(\Z^d\).  For the proposed placement \(v\), the cluster head has height
\(v\cdot\vecell\).  Its adjusted cut height is therefore
$\hadj(v):=h+(v\cdot\vecell-u)=v\cdot\vecell-(u-h)$.
Put
\[
 \begin{split}
 O_v&=E(\upwell^v)\cup E(\downwell),\\
 C_v&=\bigl(E(\downspace_h)\cap\partial_E(\downwell)\bigr)
       \cup
       \bigl(E(\upspace_{\hadj(v)})
                    \cap\partial_E(\upwell^v)\bigr) \, , 
 \end{split}
\]
where \(\partial_E(G)\) is the edge boundary of a subgraph \(G\), namely the set
of edges in \(E(\Z^d)\setminus E(G)\) incident to a vertex of \(G\).

Edges in \(O_v\) are called \emph{discovered open}, and those in \(C_v\)
\emph{discovered closed}; remaining edges, in
\(E(\Z^d)\setminus(O_v\cup C_v)\), are \emph{undiscovered}.

The sets \(O_v\) and \(C_v\) are disjoint.
 Indeed, \(v\cdot\vecell\geq u\) gives
\(\hadj(v)\geq h\).  The two downcluster terms lie in
\(E(\downspace_h)\), while the two upcluster terms lie in
\(E(\upspace_{\hadj(v)})\); their cross-intersections are therefore
excluded by the edge partition at height~\(\hadj(v)\), while the
same-cluster intersections are excluded by
\(E(A)\cap\partial_E(A)=\varnothing\).  Every discovered-open edge is
\(\omega''_v\)-open, every discovered-closed edge is \(\omega''_v\)-closed,
and each undiscovered edge is \(\omega''_v\)-open precisely when it is
\(\omega'\)-open.
\end{definition}

Let \(\uclust\) and \(\dclust\) denote possible values for
\([\upwell]\) and \(\downwell\) as \(\omega\) ranges over
\(\clusteruthree\).  For \(v\in\Z^d\) with
\(v\cdot\vecell\in[u,u+3r]\), write
\(\mathfrak u^v\) for the unique translate determined by the fixed class
\(\uclust\) whose head is~\(v\), in parallel with
Definition~\ref{d.displace}.  Set \(R^v_{\uclust,\dclust}(\omega')\) equal to the
configuration produced by Definition~\ref{d.prime} from retained data
\((\uclust,\dclust)\), proposed head~\(v\), and auxiliary configuration
\(\omega'\).  On these fixed values, the sets \(O_v\) and \(C_v\) of
Definition~\ref{d.prime} are determined by
\((\uclust,\dclust,v)\).  Their states are fixed, while every undiscovered
edge is assigned the state of the corresponding edge of \(\omega'\).

Set the \hfff{fixed-shift-law}{fixed-shift reconstruction law}
\[
   \Rlaw^v_{\uclust,\dclust}:\mc B\to[0,1]\, ,
   \qquad
   \Rlaw^v_{\uclust,\dclust}(A)
\, = \,
   \PP'\big(\omega'\in\Omega'\, :
      R^v_{\uclust,\dclust}(\omega')\in A\big)\, .
\]

For the slide height~\(h\) under consideration and \(v\in\Z^d\) with
\(v\cdot\vecell\in[u,u+3r]\), define
the \hfff{accept-event}{fixed-shift acceptance event \(\Accept_v\)} by
\begin{equation}
\label{e.fixed-shift-acceptance}
 \Accept_v
 =\fine_h\cap\clusteruthree\cap\{V=v\}.
\end{equation}
We view this event under \(\Rlaw^v_{\uclust,\dclust}\).  On
\(\Accept_v\), the reconstructed configuration has downcluster~\(\dclust\)
and upcluster~\(\mathfrak u^v\) at height~\(h\): the discovered closed
boundaries in Definition~\ref{d.prime} prevent either retained graph from
extending in its defining half-space, \(\clusteruthree\) places the retained
pieces in the finite origin cluster, and \(V=v\) identifies its upper piece.
We call \(\Rlaw^v_{\uclust,\dclust}(\Accept_v)\) the \emph{connection chance}
of the proposed head~\(v\).

Since \(\fine_h\subseteq\clusteruthree\), the middle factor in
\eqref{e.fixed-shift-acceptance} is redundant; it is retained to display the
finite-connection requirement entering the reconstruction.

In the next result, the left-hand side is an ordinary conditional probability
under \(\PP\).  The events \(B\) and \(H\) are subsets of the original
configuration space \(\Omega\); on the right-hand side,
\(\Rlaw^v_{\uclust,\dclust}\) evaluates these events at the reconstructed
configuration \(R^v_{\uclust,\dclust}(\omega')\).

\begin{lemma}\label{l.probform}
Let \(B,H\subseteq\Omega\) denote two \(\mc B\)-measurable events, with
\(\PP(B\vert\clusteruthree)>0\).  Let \((\uclust,\dclust)\) be a value
that \(([\upwell],\downwell)\) assumes with positive probability under
\(\PP\) given \(B\cap\clusteruthree\).  Then
\begin{eqnarray*}
& & \PP\Big(
       H \,\Big\vert\, [\upwell]\, = \, \uclust,\downwell\, = \, \dclust\, ,
       B,\clusteruthree
    \Big) \\
&\, = \, &
\frac{
\sum_{v\in\fubuthree}
\Rlaw^v_{\uclust,\dclust}
\big(H\cap B\cap\clusteruthree\cap\{V\, = \, v\}\big)
}{
\sum_{v\in\fubuthree}
\Rlaw^v_{\uclust,\dclust}
\big(B\cap\clusteruthree\cap\{V\, = \, v\}\big)
}\, ,
\end{eqnarray*}
where \(V\) is the random vector specified under
\(\PP(\cdot\vert\clusteruthree)\) in Definition~\ref{d.displace}.
\end{lemma}

\noindent{\bf Proof.} Let \(G\subseteq\Omega\) be an arbitrary \(\mc B\)-measurable event.  Fix
\(v\in\Z^d\) with \(v\cdot\vecell\in[u,u+3r]\), and consider the event
$\clusteruthree\cap G\cap\{\upwell=\mathfrak u^v\}
\cap\{\downwell=\dclust\}$.  Does this event occur? We may find the answer in two stages.  First, we
check a necessary condition, an event that we label~\(N\): namely, that
elements of \(E(\mathfrak u^v)\cup E(\dclust)\) are \(\omega\)-open and
that elements of
\[
   E(\downspace_h)\cap\partial_E(\dclust)
   \quad\textrm{and}\quad
   E\big(\upspace_{\hadj(v)}\big)\cap\partial_E(\mathfrak u^v)
\]
are \(\omega\)-closed.  If this test is passed, we move to the second
stage, and realise all further edges in~\(E(\Z^d)\).  Check that the event
\(\clusteruthree\cap G\) occurs.  If this happens, then \(\downwell\) must
equal~\(\dclust\), but \(\upwell\) is not yet assuredly \(\mathfrak u^v\).
We need to confirm the latter condition, and this is done by checking, at
the end of the second stage, that the random vector \(V\) equals \(v\).

Writing \(\varphi_v=\PP(N)>0\), note that \(\varphi_v\) and
\(\Rlaw^v_{\uclust,\dclust}\big(G\cap\clusteruthree\cap\{V=v\}\big)\) equal
the respective probabilities that the checks at the first and second
stages are passed.  The event~\(N\) prescribes the edges of~\(O_v\) open and
those of~\(C_v\) closed.  By the disjointness established in
Definition~\ref{d.prime} and the product Bernoulli(\(p\)) law,
\[
   \varphi_v=\PP(N)=p^{|O_v|}(1-p)^{|C_v|}.
\]
The pair \((\mathfrak u^v,\upspace_{\hadj(v)})\) is a translate, in adjusted
coordinates, of a fixed pair.  Thus the cardinalities of
\(E(\mathfrak u^v)\) and
\(E(\upspace_{\hadj(v)})\cap\partial_E(\mathfrak u^v)\) do not depend on~\(v\).
The corresponding downcluster sets are fixed.  Hence \(|O_v|\) and \(|C_v|\)
are independent of~\(v\), and \(\varphi_v\) may be written
\(\varphi=p^{|O_v|}(1-p)^{|C_v|}\).  What we have learnt is that, for every \(v\in\Z^d\) satisfying
\(v\cdot\vecell\in[u,u+3r]\),
\begin{equation}\label{e.resampleclaim}
\PP\Big(
   G,\upwell\, = \, \mathfrak u^v,\downwell\, = \, \dclust,\clusteruthree
\Big)
\, = \,
\varphi\cdot
\Rlaw^v_{\uclust,\dclust}
\big(G\cap\clusteruthree\cap\{V\, = \, v\}\big)\, .
\end{equation}
Now set \(G=B\) and sum~\eqref{e.resampleclaim} over all such~\(v\).
The resulting value is positive by our assumptions.
Then take \(G=H\cap B\) and sum similarly. Dividing the latter of the
resulting equalities by the former, we obtain Lemma~\ref{l.probform}. \qed

Taking \(B=\fine_h\) in Lemma~\ref{l.probform} gives, for every admissible
retained pair,
\begin{equation}
\label{e.fixed-shift-head-law}
 \PP^{\fine_h}\bigl(V=v\mid[\upwell]=\uclust,\downwell=\dclust\bigr)
 =
 \frac{\Rlaw^v_{\uclust,\dclust}(\Accept_v)}
 {\sum_{w\in\fubuthree}\Rlaw^w_{\uclust,\dclust}(\Accept_w)}.
\end{equation}
Thus the active slide chooses its head with probability proportional to the
fixed-shift connection chance.

\medskip
Identity~\eqref{e.resampleclaim} also yields the following single-slide invariance.

\begin{proposition}[Single-slide invariance]
\label{p.single-slide-invariance}
For every height $h\in(0,u)$ the single-slide operator $\slide_h$ of
Definition~\ref{d.slide-operator} preserves $\Q=\PP(\,\cdot\mid\clusteruthree)$:
if $\omega\sim\Q$, then $\slide_h\omega\sim\Q$.
\end{proposition}

\noindent{\bf Proof.} By Lemma~\ref{l.fine-retained}, on $\clusteruthree$ the event $\fine_h$ is
determined by the retained pair $(\downwell_h,[\upwell_h])$.  On a retained-pair
fibre contained in $\fine_h$, the output is obtained by drawing the shift $V$
according to the connection-chance law~\eqref{e.fixed-shift-head-law} and
reconstructing the midpiece.  Identity~\eqref{e.resampleclaim}, summed over all \(v\) satisfying
\(v\cdot\vecell\in[u,u+3r]\) with the retained data fixed, shows that this is the
conditional law of $\Q$ given the retained pair; conditioning additionally on
$\fine_h$ adds no information on such a fibre.  On every complementary
retained-pair fibre the slide is dormant and leaves the configuration
unchanged.  Conditioning on the retained pair therefore shows that
$\slide_h$ preserves $\Q$. \qed

\medskip

The preceding one-slide calculation uses a deterministic division between
retained and resampled edge states.  In the iterated-slide analysis, a
prospective midpiece will be explored in stages and paused at selected layers,
so the set of auxiliary edges inspected before a pause is random.  We
record two elementary product-measure facts that will preserve the product law
on the uninspected edges in these later explorations.

For an edge set \(A\subseteq\Edg(\Z^d)\) write
$\cF_A=\sigma(\{e\text{ open}\}:e\in A)\subseteq\mc B$ for the
\(\sigma\)-field generated by the states of the edges of \(A\); a random
variable is \(\cF_A\)-measurable exactly when it is a function of \(\omega\)
through \(A\) alone.  Under \(\PP\), the \(\sigma\)-fields associated with
pairwise disjoint edge sets are independent.  We write \(\sqcup\) for
disjoint union and \(\bigsqcup\) for an indexed disjoint union.  Thus a
deterministic partition
\(\Edg(\Z^d)=\bigsqcup_i A_i\) yields independent product
Bernoulli(\(p\)) restrictions \(\{\omega|_{A_i}\}_i\).

\medskip

Call a random edge set \(\Xi=\Xi(\omega)\subseteq\Edg(\Z^d)\) a \emph{stopping
domain} if it is almost surely finite and
\[
   \{\Xi\, = \, A\}\in\cF_A
   \qquad\text{for every finite edge set }A\, ;
\]
this is the usual situation when a percolation exploration inspects edges one
at a time, each choice of the next edge being determined by the states already
seen.  Thus whether the exploration has inspected exactly \(A\) is decided by
the states of the edges of \(A\) alone.  Its \emph{stopped \(\sigma\)-field} is
\[
   \cF_\Xi\, = \, \bigl\{G\in\mc B:\ G\cap\{\Xi\, = \, A\}\in\cF_A\ \text{for every }A\bigr\}\, .
\]

\begin{lemma}[Freshness of the unexplored complement]
\label{l.freshness}
Let \(\Xi\) be a stopping domain.  For every finite edge set \(A\), on the
event \(\{\Xi=A\}\) the conditional law given \(\cF_\Xi\) gives
\(\omega|_{A^c}\) its product Bernoulli(\(p\)) law.  More precisely, for every
bounded measurable function \(f_A\) of the edge states in~\(A^c\),
\[
 \mathbf 1_{\{\Xi=A\}}
 \E\bigl[f_A(\omega|_{A^c})\mid\cF_\Xi\bigr]
 =
 \mathbf 1_{\{\Xi=A\}}
 \E \bigl[f_A(\omega|_{A^c})\bigr]
 \qquad\PP\text{-almost surely.}
\]
\end{lemma}

\noindent{\bf Proof.} Fix a finite edge set \(A\) and an event \(G\in\cF_\Xi\).  By the definition
of the stopped \(\sigma\)-field,
\(G\cap\{\Xi=A\}\in\cF_A\).  Independence of the edge states under the product law gives
\[
 \E\bigl[\mathbf 1_G\mathbf 1_{\{\Xi=A\}}
          f_A(\omega|_{A^c})\bigr]
 =
 \PP\bigl(G\cap\{\Xi=A\}\bigr)
 \E_{\PP}\bigl[f_A(\omega|_{A^c})\bigr] .
\]
This is exactly the displayed conditional-expectation identity. \qed

Throughout, an edge is \emph{unread}, relative to the exploration history in force, if its state has not been revealed in that history.

The just-stated freshness conclusion applies on the finite part of a possibly infinite exploration stopping time~\(T\).  Let \((\cF_n)\) denote the exploration filtration, with \(\cF_n\) generated by the edge states revealed through the first \(n\) queries.  Apply the lemma to the finite explored domain at time \(T\wedge n\) and restrict to \(\{T=n\}\).  Taking the union over finite~\(n\) gives the product law for the unread edges given \(\cF_T\), on \(\{T<\infty\}\).  This argument also applies conditionally on fixed data when the remaining edge states initially have a product law; it requires no conditioning on eventual termination of the whole exploration.

\subsection{Reducing Proposition~\ref{p.success} to a partition function upper bound}\label{s.reduction}

Here \(\fine\) denotes \(\fine_h\).  We first record the elementary translation
underlying the construction.  The reduction then uses the SimpleJoin lower
bound of Lemma~\ref{l.qfine} and the partition-function upper bound of
Proposition~\ref{p.finesum}.

\begin{lemma}\label{l.qfine-translate}
For any \(a\in\Bed(\upwell)\) and \(b\in\Top(\downwell)\) there is a unique
translate \(\upwell^x\) of \(\upwell\) in which the translate of \(a\)
equals \(b+2e_d\); its head location \(x\) satisfies
\(x\cdot\vecell\in[u,u+3r]\).
\end{lemma}

\noindent{\bf Proof.} A translate of \(\upwell\) is determined by the image of any single vertex, so
requiring the translate of \(a\) to lie at \(b+2e_d\) fixes the translate, and
with it the head location \(x\), uniquely.  The height \(x\cdot\vecell\) equals
\(b\cdot\vecell+2r+\hmax-a\cdot\vecell\).  Since
\(b\in\vertbdry_h\) and \(a\in\vertbdry_{\hadj}\),
Lemma~\ref{l.geometrybasics}(1) gives
\[
 b\cdot\vecell\in(h-r,h]\, ,
 \qquad
 \hmax-a\cdot\vecell\in[u-h,u-h+r)\, .
\]
It follows that \(x\cdot\vecell\in(u+r,u+3r)\subset[u,u+3r]\). \qed

\begin{lemma}\label{l.qfine}
There exists \(c>0\), depending only on \(p,d\) and \(\af(0)\), such
that
\[
   \sum_{v\in\fubuthree}
      \Rlaw^v_{[\upwell],\downwell}(\Accept_v)\, \geq \, c
\]
whenever \(\omega\in\fine\). Indeed, for such \(\omega\) there exists
\(x\in\Z^d\) satisfying \(x\cdot\vecell\in[u,u+3r]\) for which
\[
   \Rlaw^x_{[\upwell],\downwell}
      \big(\success\cap\Accept_x\big)\, \geq \, c\, .
\]
\end{lemma}

We next state the partition-function upper bound. Its proof is the subject of
the next subsection.

\begin{proposition}\label{p.finesum}
Suppose that the fine allowance satisfies~\eqref{e.fine-polynomial-growth}.
For any such choice of \(\af\), and for any choice of the breaker-profile
constants \(b_0,C_0,C_1,i_{\rm br}\) in \({\bf B}^{+}_{\rm f}\), there exists
\(C<\infty\), independent of these breaker-profile constants and uniform in~\(\vecell\), such that\[
   \sum_{\substack{v\in\Z^d:\\ v\cdot\vecell\in[u,u+3r]}}
      \Rlaw^v_{[\upwell],\downwell}(\Accept_v)
   \,\leq\, C .
\]
whenever \(\omega\in\fine\), \(u\ge K\), \(h>2r\) and \(u-h>2r\).  For the allowance in
\eqref{e.allowance-design}, once \(A_0,A_1\) and the coarse geometry are
fixed, the upper-bound constant may be chosen so that
\begin{equation}
\label{e.finesum-Cfine-polynomial}
 C\,\leq\, C_*(2+\Cfine)^{D_*}
\end{equation}
for fixed \(C_*,D_*<\infty\), which may be chosen uniformly in~\(\vecell\).
\end{proposition}

\noindent{\bf Proof of Proposition~\ref{p.success}.}
Take \(H=\success\) and \(B=\fine_h\) in Lemma~\ref{l.probform}.
By \eqref{e.fixed-shift-acceptance}, its denominator is the
connection-chance sum
\(\sum_{v\in\fubuthree}\Rlaw^v_{[\upwell],\downwell}(\Accept_v)\), and hence
is bounded above by Proposition~\ref{p.finesum}.  Its numerator is bounded
below by the second assertion of Lemma~\ref{l.qfine}. We obtain
$\PP^{\fine_h}(\success\mid[\upwell],\downwell)\geq c$ whenever
\(u\ge K\), \(h>2r\) and \(u-h>2r\), by
Proposition~\ref{p.finesum}.  For the allowance in~\eqref{e.allowance-design}, the lower bound in
Lemma~\ref{l.qfine} depends only on the already fixed \(A_1=\af(0)\), while
\eqref{e.finesum-Cfine-polynomial} bounds the denominator by
\(C_*(2+\Cfine)^{D_*}\).  The numerator lower bound from Lemma~\ref{l.qfine} and the denominator
bound from Proposition~\ref{p.finesum} are both uniform in~\(\vecell\).
Decreasing \(c_*\) if necessary gives
\eqref{e.simplejoin-Cfine-lower} uniformly in direction.  This proves Proposition~\ref{p.success}.
\qed

To make the prescribed two-edge SimpleJoin with uniformly positive
probability, we need a uniform bound on the number of possible bridge
endpoints.

\begin{lemma}\label{l.finecard}
The event \(\fine\) entails that the vertex sets \(\Bed(\upwell)\) and
\(\Top(\downwell)\) have cardinality at most \(C_{\rm bd}\), where
\(C_{\rm bd}<\infty\) depends only on \(d\) and \(\af(0)\), uniformly in~\(\vecell\).
\end{lemma}

\noindent{\bf Proof.} By Lemma~\ref{l.geometrybasics}(1), every element of
\(\Top(\downwell)\subset\vertbdry_h\) is a vertex in \(\csetuthree\) whose
height lies in \([h-r,h]\), where \(r=e_d\cdot\vecell\). Since \(r\leq1\),
this interval is contained in \([h-1,h]\), and \({\bf A}^{-}_{\rm f}(0)\) gives \(|\Top(\downwell)|\leq\af(0)\).

Similarly, every element of \(\Bed(\upwell)\subset\vertbdry_{\hadj}\) has
height in \([\hadj-r,\hadj]\), and this interval is contained in
\([\hadj-1,\hadj]\). \({\bf A}^{+}_{\rm f}(0)\) gives
\(|\Bed(\upwell)|\leq\af(0)\). This proves the lemma. \qed

\noindent{\bf Proof of Lemma~\ref{l.qfine}.}
We prescribe a SimpleJoin.  Take \(b=\head(\downwell)\).  By
Lemma~\ref{l.updowncluster}, \(b\in\Top(\downwell)\).  Take
\(a=\base(\upwell)\), which belongs to \(\Bed(\upwell)\) by the same lemma.
By Lemma~\ref{l.qfine-translate} there is a unique translate \(\upwell^x\)
in which the translate of \(a\) equals \(b+2e_d\), with head location
\(x\) satisfying \(x\cdot\vecell\in[u,u+3r]\).  Its adjusted cut height is
$\hadj(x)=h+(x\cdot\vecell-u)=x\cdot\vecell-(u-h)$.  Put \(m=b+e_d\); this is the proposed SimpleJoin midvertex.
The two bridge edges are undiscovered in Definition~\ref{d.prime}.  Indeed,
\(b\cdot\vecell>h-r\), so their greater endpoint heights exceed~\(h\), while
\(b+2e_d\in\Bed(\upwell^x)\) has height at most~\(\hadj(x)\).  Thus both edges lie in
\(E(\upspace_h)\cap E(\downspace_{\hadj(x)})\), and hence in neither~\(O_x\) nor~\(C_x\).

Let \(\mathcal J\) be the event that the edges \([b,m]\) and
\([m,b+2e_d]\) are \(\omega'\)-open, while all other edges incident to
\(m\) are \(\omega'\)-closed, as are all other undiscovered edges that
are incident to \(\Bed(\upwell^x)\cup\Top(\downwell)\). For each
\(\omega\in\fine\), Lemma~\ref{l.finecard} gives
\[
   \R^x_{[\upwell],\downwell}(\mathcal J)
\, \geq \,
      p^2\,(1-p)^{2(2d-1)C_{\rm bd}+2d-2}\, .
\]
Indeed, the right-hand side opens the two \(e_d\)-oriented edges \([b,m]\) and
\([m,b+2e_d]\)---the factor \(p^2\)---closes the other edges incident to the
shared endpoint \(m\), and closes the other needed undiscovered edges incident to
\(\Bed(\upwell^x)\cup\Top(\downwell)\), of which there are at most
\(2(2d-1)C_{\rm bd}\), since
\(|\Bed(\upwell)|,|\Top(\downwell)|\leq C_{\rm bd}\) by
Lemma~\ref{l.finecard}.  The first group contains \(2d-2\) edges, giving
the exponent displayed above.  Absorbing this positive quantity
into \(c\), we obtain a lower bound depending only on \(p,d\) and
\(\af(0)\), uniformly in~\(\vecell\).

When the configuration \(R^x_{[\upwell],\downwell}(\omega')\) realises
\(\mathcal J\), the only open edges leaving the union of the retained
downcluster and upcluster are the two bridge edges.  Indeed, an edge leaving
this union from a vertex of \(\downwell\) is either discovered closed or
undiscovered.  In the latter case it cannot belong to \(E(\downspace_h)\),
so the edge partition places it in \(E(\upspace_h)\); its endpoint in
\(\downwell\) consequently lies in \(\vertbdry_h\), and hence in
\(\Top(\downwell)\), where \(\mathcal J\) closes it unless it is the bridge
edge \([b,m]\).  The corresponding statement for \(\upwell^x\) holds with
\(\Bed(\upwell^x)\), and all other edges incident to~\(m\) are closed by
\(\mathcal J\).  Hence the origin cluster is precisely the union of
\(\downwell\), \(\upwell^x\), and the two bridge edges.  It is finite and has
head~\(x\), with \(x\cdot\vecell\in[u,u+3r]\); its component \(\midconnect\) consists precisely of
the two bridge edges and \(\midpiece\setminus\midconnect\) is empty.  Thus
\(R^x_{[\upwell],\downwell}(\omega')\) realises the event \(\success\), with
midvertex \(m\).

We next confirm that \(V=x\).  It is enough to check that the upcluster under
\(R^x_{[\upwell],\downwell}(\omega')\) is equal to \(\upwell^x(\omega)\).
There is a unique open connected component of \(\csetuthree\) in
\(\upspace_{\hadj(x)}\) under this configuration. Indeed, this
component contains \(\upwell^x(\omega)\), and no open edge of the
resampled configuration joins \(\upwell^x(\omega)\) inside
\(\upspace_{\hadj(x)}\) to any further vertex of \(\csetuthree\).
For every edge in
$E(\upspace_{\hadj(x)})\cap\partial_E(\upwell^x)$, this follows from the specification of discovered closed edges in
Definition~\ref{d.prime}; and every undiscovered edge incident to
\(\Bed(\upwell^x)\) is closed under \(\mathcal J\), except for the
single edge \([m,b+2e_d]\), whose upper endpoint already lies in
\(\upwell^x(\omega)\). Thus, the unique open connected component of
\(\csetuthree\) in \(\upspace_{\hadj(x)}\) equals
\(\upwell^x(\omega)\), and hence \(V=x\).

The resampled configuration therefore has the same retained pair
\(\bigl([\upwell],\downwell\bigr)\) as \(\omega\): its downcluster is unchanged and
its upcluster is the translate \(\upwell^x\).  Since \(\omega\in\fine\),
Lemma~\ref{l.fine-retained} now shows that the resampled configuration also
lies in \(\fine\).  Thus it belongs to \(\fine\cap\success\cap\{V=x\}\).

We learn that, for \(\omega\in\fine\),
$\R^x_{[\upwell],\downwell}(\fine\cap\success\cap\{V=x\})\geq c$, where \(c>0\) depends only on \(p,d\) and \(\af(0)\), uniformly in~\(\vecell\). The latter
assertion of Lemma~\ref{l.qfine} has been derived; the former is a
trivial consequence. \qed

\subsection{The Tether and the partition function}\label{s.finesum}

Here we prove Proposition~\(\ref{p.finesum}\), the bound on the sum of connection chances that keeps the shift at the end of a short Tether.
 To argue that the connection chance
\(\Rlaw^v_{[\upwell],\downwell}(\Accept_v)\) is small when the shift \(v\) is
large, we work with the coarse-grained surface associated to the middle piece.
Only the two allowance families ${\bf A}^{-}_{\rm f}$ and~${\bf A}^{+}_{\rm f}$ are used below.  Dropping ${\bf B}^{+}_{\rm f}$ enlarges
the event in every summand, so the resulting constant is uniform in all
breaker-profile parameters.
A connection between the shifted upcluster and the downcluster must either
bring their coarse neighbourhoods together, or else trail a Tether---a path of
bad boxes through territory undiscovered under the resampling law---whose length
grows with the reach of the shift. The
fine allowance bounds the number of possible endpoint boxes at each
latitude, while a long Tether has exponentially small probability.

Set
\[
   A_{\rm f}(n)
\, = \, 1+\sum_{m=0}^{\lceil n\rceil}\af(m)\, ,
      \qquad n\, \geq \, 0\, .
\]
The polynomial-growth condition~\eqref{e.fine-polynomial-growth} gives
\[
   A_{\rm f}(n)
\, \leq \, 1+C_{\rm f}\sum_{m=0}^{\lceil n\rceil}(1+m)^q
\, \leq \, C_{\rm f}'(1+n)^{q+1}\, .
\]

Recall from the coarse-grained setup preceding
Definition~\ref{d.exterior} that, for a subgraph \(A\),
\(\overline A\) is the collection of \(K\)-boxes containing an edge of
\(A\).  If \(\midconnect=\varnothing\), put
\(\mc H=\surfacek=\varnothing\).  
 Otherwise, since $\midconnect$ is connected,
$\overline{\midconnect})$ is nearest-neighbour connected in the box lattice,
and hence $*$-connected.
  It therefore lies in a single \(*\)-connected
component of~\(\overline{\midpiece}\); call this component~\(\mc H\), and
set
\[
   \surfacek\, = \, \partial_{\rm ext}\mc H\, .
\]

\begin{lemma}\label{l.surfacek}
If \(\midconnect\ne\varnothing\), then
\begin{enumerate}
\item the collection \(\surfacek\) is a nearest-neighbour connected set of
\(K\)-boxes;
\item no element of \(\surfacek\) contains an edge of \(\midpiece\);
\item every infinite \(K\)-box path beginning at a box that intersects
\(\midconnect\) visits \(\surfacek\).
\end{enumerate}
\end{lemma}

\noindent{\bf Proof: (1).} The collection \(\surfacek\) is the box-lattice
exterior boundary of the finite \(*\)-connected box collection \(\mc H\).
Lemma~\ref{l.extbdy-halfspace} applies.

{\bf (2).}  Suppose that \(B\in\surfacek\) contains an edge of
\(\midpiece\).  Then \(B\in\overline{\midpiece}\).  Since \(B\) is a
\(*\)-neighbour of \(\mc H\), and \(\mc H\) is a \(*\)-connected component
of \(\overline{\midpiece}\), this would give \(B\in\mc H\), contrary to
\(B\in\partial_{\rm ext}\mc H\).

{\bf (3).}  A \(K\)-box that intersects \(\midconnect\) either belongs to
\(\mc H\) or is a \(*\)-neighbour of an element of \(\mc H\).  If the given
infinite path never visits \(\mc H\), its initial box therefore belongs to
\(\partial_{\rm ext}\mc H\).  Otherwise, since \(\mc H\) is finite and the
path runs to infinity, take the first box after its last visit to \(\mc H\).
It is a nearest neighbour of \(\mc H\), and the infinite tail beginning there
is disjoint from \(\mc H\); hence this box belongs to
\(\partial_{\rm ext}\mc H=\surfacek\). \qed

\begin{lemma}\label{l.surfacebad}
Suppose that \(\clusteruthree\) occurs and \(u\ge K\). Let \(B\) be a
\(K\)-box that lies in \(\surfacek\) and whose vertex set is disjoint from
the union of the vertex sets of \(\upwell\) and \(\downwell\). Then \(B\)
is a bad \(K\)-box.
\end{lemma}

\noindent{\bf Proof.}
If \(\midconnect=\varnothing\), then \(\surfacek=\varnothing\) and there is
nothing to prove.  Otherwise the hypothesis excludes edges of \(\downwell\)
and \(\upwell\) from~\(B\),
while Lemma~\ref{l.surfacek}(2) excludes edges of \(\midpiece\).  Thus
\(B\notin\overlinecsetuthree\).  Let \(B'\) denote a \(K\)-box
\(*\)-neighbour of \(B\) that lies in~\(\mc H\); it contains an edge of
\(\midpiece\subseteq\csetuthree\), and hence
\(B'\in\overlinecsetuthree\).  Since \(\clusteruthree\) occurs and
\(u\ge K\), the cluster \(\csetuthree\) has Euclidean diameter at least~\(K\).
Lemma~\ref{l.bdbad} now implies that \(B\) is a bad \(K\)-box. \qed

It is useful to record the way in which the fine allowance controls the
number of possible endpoint boxes at a given latitude. The next result and
its proof utilise the notions of \(K\)-box latitude from
Definition~\ref{d.boxlatitude}, and of the \(K\)-box fattening
\(\overline A[s]\) from Definition~\ref{d.plustwominustwo}.  Fix an integer-valued 
 \emph{Tether buffer}
\(b_{\rm T}\geq 2d+3\) large enough that,
whenever a box lies outside \(\overline A[b_{\rm T}]\), its associated big
box contains no edge of~\(A\) and no edge incident to~\(A\).  Such a fixed
choice is possible because the big-box radius and the box overlaps are fixed
multiples of~\(K\).

Throughout this lemma we reckon latitude relative to the height-\(h\) layer.
For a \(K\)-box \(B\), put
\[
   \lat_h(B)=\lat(B)-\Big\lfloor\frac{h}{2Kr}\Big\rfloor,
\]
where \(\lat(B)\) is the absolute latitude of
Definition~\ref{d.boxlatitude}.  If \(B=2Kw+\llbracket-K,K\rrbracket^d\),
then
\begin{equation}\label{e.relative-latitude-slack}
   \big|2Kw\cdot\vecell-h-2Kr\,\lat_h(B)\big|<2Kr.
\end{equation}

\begin{lemma}\label{l.latcard}
There exist \(C<\infty\) and \(C_{\rm lat}<\infty\), depending only on
\(d,K\), uniformly in~\(\vecell\), such that, whenever \(\fine\) occurs,
\begin{enumerate}
\item for every integer \(i\geq -C_{\rm lat}\),
\[
   \big\vert\{B\in\overline{\downwell}[b_{\rm T}]:\lat_h(B)=-i\}\big\vert
\, \leq \, C A_{\rm f}\big(C(1+|i|)\big)\, ;
\]
\item for every integer \(j\geq -C_{\rm lat}\),
\[
   \big\vert\{B\in\overline{\upwell}[b_{\rm T}]:\lat_h(B)=j\}\big\vert
\, \leq \, C A_{\rm f}\big(C(1+|j|)\big)\, .
\]
\end{enumerate}
\end{lemma}
\textbf{Proof: (1).}
Let $B = 2Kw + \llbracket -K,K\rrbracket^d$ be a $K$-box with
\(\lat_h(B)=-i\) lying in $\overline{\downwell}[b_{\rm T}]$.  By definition of this fixed
fattening, there exists a $K$-box $B_0 \in \overline{\downwell}$ at
$*$-distance at most $b_{\rm T}$ from $B$ in the $K$-box graph.  Since
$B_0 \in \overline{\downwell}$, it contains an edge of $\downwell$, hence a
vertex $x \in V(\downwell)$.

We bound $h - x\cdot\vecell$.  Write
\(B_0=2Kw_0+\llbracket-K,K\rrbracket^d\), whose centre is \(2Kw_0\).
Since $B_0$ is at $*$-distance at most $b_{\rm T}$ from
$B$, we have $\|w_0 - w\|_\infty \leq b_{\rm T}$, so
\[
   |2K(w_0 - w)\cdot\vecell|
\, \leq \, 2K\sum_{q=1}^d |w_0^q - w^q|\,|e_q\cdot\vecell|
\, \leq \, 2Kb_{\rm T}d\,r
\, \leq \, C_0 K\, ,
\]
where $C_0$ depends only on $d$ and the fixed value of
$b_{\rm T}$, since $r\leq1$.  Since $x \in B_0$, we
also have $|x\cdot\vecell - 2Kw_0\cdot\vecell| \leq K\|\vecell\|_1 \leq C_0 K$
after enlarging $C_0$.  Hence
\[
   |x\cdot\vecell - 2Kw\cdot\vecell|\, \leq \, C_1 K\, ,
\]
where $C_1$ depends only on $d$ and $b_{\rm T}$, uniformly in~$\vecell$.

By~\eqref{e.relative-latitude-slack}, $|2Kw\cdot\vecell-h+2Kri|<2Kr$.
Consequently
\[
   h - x\cdot\vecell
\, \leq \, |h-2Kw\cdot\vecell|+C_1K
\, \leq \, 2Kr(1+|i|)+C_1K
\, \leq \, C_2K(1+|i|)\, ,
\]
for a constant $C_2$ depending only on $d$ (and the fixed box geometry), uniformly in~$\vecell$.  Since
$x \in V(\downwell) \subset V(\downspace_h)$, we have $x\cdot\vecell \leq h$,
so $h - x\cdot\vecell \geq 0$.

We also verify the lower bound on $i$.  Any $B \in \overline{\downwell}[b_{\rm T}]$
has a neighbouring $B_0 \in \overline{\downwell}$ containing a vertex of
$\downwell$, which has height at most $h$.  The centre height of $B_0$
is then at most $h + C_0 K$, and so the centre height of $B$ is at most
$h + 2C_0 K$.  Equation~\eqref{e.relative-latitude-slack} and
$r\geq d^{-1/2}$ therefore give \(\lat_h(B)\leq C_{\rm lat}\) for a constant
$C_{\rm lat}$ depending only on $d,K$.  Hence the condition
$i \geq -C_{\rm lat}$ accounts for all $K$-boxes that can occur in
$\overline{\downwell}[b_{\rm T}]$.

Thus $x$ lies in the height slab $[h - C_2 K(1+|i|),\, h]$, which is
covered by the bands $[h - m - 1, h - m]$ for $0 \leq m \leq C_2 K(1+|i|)$.
The fine allowance condition ${\bf A}^{-}_{\rm f}(m)$ bounds the number of vertices
of $\downwell$ in the band $[h-m-1,h-m]$ by $\af(m)$.  Summing over
admissible $m$,
\[
   \bigl|\{x \in V(\downwell) : x\cdot\vecell \in [h - C_2K(1+|i|),\, h]\}\bigr|
\, \leq \, \sum_{m=0}^{\lceil C_2 K(1+|i|)\rceil} \af(m)
\, \leq \, A_{\rm f}\bigl(C_2 K(1+|i|)\bigr)\, .
\]
For a fixed vertex~$x$, only $C(d,b_{\rm T})$ such boxes~$B$ are compatible
with~$x$ in the construction above.  Indeed, at most \(2^d\) \(K\)-boxes contain
\(x\), and each has at most \((2b_{\rm T}+1)^d\) \(K\)-boxes at
\(*\)-distance at most $b_{\rm T}$.
Combining, the number of $K$-boxes in
$\overline{\downwell}[b_{\rm T}]$ with relative latitude~$-i$ is at most
$C(d,b_{\rm T})A_{\rm f}(C_2K(1+|i|))\leq C A_{\rm f}(C(1+|i|))$, since $K$ is fixed and $A_{\rm f}$ is non-decreasing.  This proves~(1).

\medskip
\noindent{\bf (2).}
The argument mirrors~(1) with $\downwell$ replaced by $\upwell$ and heights
measured in the adjusted coordinate $x\cdot\vecell - \hadj$.

Let $B = 2Kw + \llbracket -K,K\rrbracket^d$ be a $K$-box with
\(\lat_h(B)=j\) lying in $\overline{\upwell}[b_{\rm T}]$.  By the same fixed-fattening argument,
there exists $B_0 \in \overline{\upwell}$ at $*$-distance at most $b_{\rm T}$ from
$B$, containing a vertex $x \in V(\upwell)$, and the centre-height estimate
gives
\[
   |x\cdot\vecell - 2Kw\cdot\vecell|\, \leq \, C_1 K\, .
\]
By Lemma~\ref{l.geometrybasics}(1), every vertex $x \in V(\upwell)$
satisfies $x\cdot\vecell \geq \hadj - r$, since $\upwell$ is a subgraph of
$\upspace_{\hadj}$ and any vertex of $\upspace_{\hadj}$ has height at
least $\hadj - r$.  Hence $x\cdot\vecell - \hadj \geq -r \geq -1$, which
is bounded below uniformly in $u$ and $h$.

Equation~\eqref{e.relative-latitude-slack} gives
\[
   |2Kw\cdot\vecell-h-2Krj|<2Kr.
\]
Since \(\hadj-h=\hmax-u\in[0,3r]\), combining this with
\(\lvert x\cdot\vecell-2Kw\cdot\vecell\rvert\leq C_1K\) gives
\[
   x\cdot\vecell - \hadj
\, \leq \, C_2K(1+|j|)
\]
for a constant $C_2$ depending only on the fixed geometry, uniformly in~$\vecell$.
For $x\in\upspace_{\hadj}$ the adjusted height satisfies
$x\cdot\vecell-\hadj\geq -r$, so, by the same centre-height argument as in
part~(1), \(\lat_h(B)\) is bounded below uniformly.  Enlarging
\(C_{\rm lat}\) if necessary, the condition $j \geq -C_{\rm lat}$ therefore
accounts for all boxes that arise in $\overline{\upwell}[b_{\rm T}]$.

\begin{samepage}
Put \(M_j=\lceil C_2K(1+|j|)\rceil\).  Hence \(x\) lies in one of the
adjusted-coordinate bands
\([\hadj+m-1,\hadj+m]\), \(m\in\llbracket0,M_j\rrbracket\).  The condition
\({\bf A}^{+}_{\rm f}(m)\) bounds the number of upcluster vertices in this
band by \(\af(m)\).  Summing and applying the bounded-multiplicity argument
from part~(1) gives
\[
 \big|\{B\in\overline{\upwell}[b_{\rm T}]:\lat_h(B)=j\}\big|
 \leq C A_{\rm f}(C(1+|j|)).
\]
This proves~(2), and hence the lemma. \qed
\end{samepage}

We now realise the Tether of Chapter~\ref{s.heuristics} on coarse-grained
boxes.  When the retained graphs are disjoint, it is extracted from the
surface \(\surfacek\) as a chain of bad, undiscovered boxes lying between
them.  When the retained graphs meet, no such chain is needed and the Tether
will be empty.

Continue for the moment in the disjoint case.  Choose an open path \(Q\) in
\(\midconnect\) from \(w\in\Top(\downwell)\) to
\(w'\in\Bed(\upwell)\), in some definite manner. Let \(Q^-\) be the lexicographically minimal open path in
\(\downwell\) from \(\base(\csetuthree)=0\) to \(w\), and let \(Q^+\) be
the similarly minimal open path in \(\upwell\) from \(w'\) to
\(\head(\csetuthree)\). Denote by \(P_1\) the concatenation
\(Q^-\circ Q\circ Q^+\).

Under the axes ordering convention, \(e_d\cdot\vecell\) is the largest
positive coordinate increment.  Prefix to \(P_1\) an infinite path \(P_0\)
of consecutive positively oriented \(e_d\)-parallel edges that runs directly
into \(\base(\csetuthree)\), and postfix such a path~\(P_2\) starting at
\(\head(\csetuthree)\). The resulting path
\(P=P_0\circ P_1\circ P_2\) is a bi-infinite path in \(\Z^d\).

To remove the ambiguity caused by overlapping \(K\)-boxes, assign each
vertex \(z\in\Z^d\) to the box \(B(u)\), where \(u\in\Z^d\) is the unique
index satisfying $z-2Ku\in\llbracket-K,K-1\rrbracket^d$.  Adjacent vertices are assigned either to the same box or to boxes with
nearest-neighbour indices.  Replace the vertices of \(P\) by their assigned
boxes and delete consecutive repetitions; the resulting sequence is the
\(K\)-box path \(P[K]\).  Let \(B_1\) be the last element
of \(\surfacek\) at or before the first box of \(P[K]\) that intersects
\(\midconnect\), and let \(B_2\) be the first element of \(\surfacek\) at or after
the last box of \(P[K]\) that intersects \(\midconnect\).

\begin{lemma}\label{l.surfacepath}
Suppose that \(u\ge K\), \(h>2r\) and \(u-h>2r\), and
that \(\Accept_v\) occurs. There is a deterministic choice of
\(D[\mathcal{T}]\in\overline{\downwell}[b_{\rm T}]\), \(U[\mathcal{T}]\in\overline{\upwell}[b_{\rm T}]\),
and a \(K\)-box path \(\mathcal{T}\), the \emph{Tether}, possibly empty, such that
\begin{enumerate}
\item every element of \(\mathcal{T}\) lies in
\(\surfacek\setminus\big(\overline{\downwell}[b_{\rm T}]\cup
\overline{\upwell}[b_{\rm T}]\big)\);
\item the first element of \(\mathcal{T}\), if it exists, is a nearest neighbour
of \(D[\mathcal{T}]\), and the final element of \(\mathcal{T}\), if it exists, is a
nearest neighbour of \(U[\mathcal{T}]\);
\item if \(\mathcal{T}\) is empty, then \(D[\mathcal{T}]\) and \(U[\mathcal{T}]\) are equal or
nearest-neighbour \(K\)-boxes.
\end{enumerate}
When non-empty, \(\mathcal T\) is self-avoiding.  Every box of
\(\mathcal{T}\) is bad, and every edge in its associated big box is
undiscovered.
\end{lemma}

\noindent{\bf Proof.} Suppose first that \(V(\downwell)\cap V(\upwell)\ne\varnothing\), so that
\(\midconnect=\varnothing\).  Take the lexicographically minimal common
vertex~\(x\).  In the fixed box order choose
\(B^-\in\overline{\downwell}\) and \(B^+\in\overline{\upwell}\)
containing retained edges incident to~\(x\).  Both boxes contain~\(x\), so
their indices are equal or \(*\)-neighbours.  As \(b_{\rm T}\ge1\), the
box~\(B^-\) belongs to
\(\overline{\downwell}[b_{\rm T}]\cap
\overline{\upwell}[b_{\rm T}]\).  Take
\(D[\mathcal T]=U[\mathcal T]=B^-\) and \(\mathcal T=\varnothing\).
All assertions of the lemma follow, the last one vacuously.

Assume now that the retained graphs are disjoint.  The boxes \(B_1\) and
\(B_2\) lie in \(\surfacek\). The segment of \(P[K]\) preceding its first
box that intersects \(\midconnect\) is the box trace of \(P_0\circ Q^-\);
the segment following its last such box is the box trace of \(Q^+\circ P_2\).

If \(B_1\) occurs on the box trace of \(Q^-\), then it is at
\(*\)-distance at most one from \(\overline{\downwell}\).  Suppose instead
that it occurs on \(P_0\), and let \(z=-ne_d\), \(n\geq0\), be a vertex
of this ray assigned to \(B_1\).  Since \(B_1\in\surfacek\), it is a
\(*\)-neighbour of a box \(D\in\mc H\).  The box \(D\) contains an edge of
\(\midpiece\); choose an endpoint \(y\) of this edge.  The assigned-box
geometry gives \(\|y-z\|_\infty\leq4K\).  Since every cluster vertex has
non-negative height,
\[
   nr\, = \, -z\cdot\vecell
\, \leq \, |(y-z)\cdot\vecell|
\, \leq \, 4K\|\vecell\|_1
\, \leq \, 4K\sqrt d\, ,
\]
and \(r\geq d^{-1/2}\) gives \(n\leq4Kd\).  Thus \(B_1\) is at
\(*\)-distance at most \(2d+1\) from the box assigned to the origin.
Because \(h>2r\), the path \(Q^-\) contains an edge, so this assigned box
is at distance at most one from \(\overline{\downwell}\).  Hence
\(B_1\in\overline{\downwell}[b_{\rm T}]\).

The same argument at the upper endpoint uses that every cluster vertex has
height at most \(\hmax\): if \(z=\head(\csetuthree)+ne_d\) is assigned to a
surface box on \(P_2\), then \(nr\leq4K\sqrt d\).  Since \(u-h>2r\), the
path \(Q^+\) contains an edge, and consequently
\(B_2\in\overline{\upwell}[b_{\rm T}]\).

Since \(\surfacek\) is connected, there is a path in
\(\surfacek\) from \(B_1\) to \(B_2\). Choose one such path \(\Gamma\) of
minimal length, in some definite manner. If
\(\overline{\downwell}[b_{\rm T}]\) and \(\overline{\upwell}[b_{\rm T}]\) meet along
\(\Gamma\), choose the first box \(B\) of their intersection in the fixed box
order and set
\(D[\mathcal{T}]=U[\mathcal{T}]=B\) and \(\mathcal{T}=\emptyset\).

Otherwise, let \(D[\mathcal{T}]\) be the last element of
\(\Gamma\cap\overline{\downwell}[b_{\rm T}]\), and let \(U[\mathcal{T}]\) be the first
subsequent element of \(\Gamma\cap\overline{\upwell}[b_{\rm T}]\). We take
\(\mathcal{T}\) to be the portion of \(\Gamma\) strictly between these two boxes.
If this portion is empty, \(D[\mathcal{T}]\) and \(U[\mathcal{T}]\) are
consecutive boxes of \(\Gamma\).  The three displayed properties follow, and
minimality of \(\Gamma\) makes every non-empty \(\mathcal T\) self-avoiding.

It remains to verify the final assertion. If \(B\in\mathcal{T}\), then
\(B\in\surfacek\), and \(B\) is disjoint from both
\(\overline{\downwell}[b_{\rm T}]\) and \(\overline{\upwell}[b_{\rm T}]\). In particular,
the vertex set of \(B\) is disjoint from those of \(\downwell\) and
\(\upwell\), and
Lemma~\ref{l.surfacebad} implies that \(B\) is bad. Write \(B=B(w)\).
By the defining choice of $b_{\rm T}$, no edge of the associated big box
\(\boxbig(w)\) belongs to either retained graph or is incident to either one.
Definition~\ref{d.prime} therefore declares every edge in \(\boxbig(w)\) undiscovered
for the law \(\Rlaw^v_{[\upwell],\downwell}\). \qed

The projected distance \(\dproj\) between two \(K\)-boxes \(B_1\) and
\(B_2\) is given by translating the two boxes in the \(e_d\)-direction so
that they have latitude zero, and setting \(\dproj(B_1,B_2)\) to be the
\(\ell_1\)-distance between these translates.

Let \(i,j\in\Z\) be at least \(-C_{\rm lat}\), and let \(k\in\Nzero\). We
set \(\fv[i,j,k]\) to be the event that the deterministic construction in
Lemma~\ref{l.surfacepath} produces \(D[\mathcal{T}]\) with
\(\lat_h(D[\mathcal T])=-i\), \(U[\mathcal{T}]\) with
\(\lat_h(U[\mathcal T])=j\), and
\(\dproj(D[\mathcal{T}],U[\mathcal{T}])=k\).

\begin{lemma}\label{l.rfprob}
Let \(\epsilon\in(0,1/(2d-1))\) denote a parameter that may be chosen
arbitrarily small by making a sufficiently high choice of \(K\). For
\(\omega\in\fine\) with \(u\ge K\), the following hold.
\begin{enumerate}
\item For \(v\in\Z^d\) with \(v\cdot\vecell\in[u,u+3r]\),
\(i,j\geq-C_{\rm lat}\), and \(k\in\Nzero\),
\[
\Rlaw^v_{[\upwell],\downwell}\big(\fv[i,j,k]\big)
\, \leq \,
   C A_{\rm f}\big(C(1+|i|)\big)
   \big((2d-1)\epsilon\big)^{c(|i+j|+k)}\, .
\]
\item For such \((i,j,k)\), the set
\[
 \bigl\{v\in\Z^d : v\cdot\vecell\in[u,u+3r],
   \Rlaw^v_{[\upwell],\downwell}\big(\fv[i,j,k]\big)>0\bigr\}
\]
has cardinality at most
\[
   C K^d
   A_{\rm f}\big(C(1+|i|)\big)
   A_{\rm f}\big(C(1+|j|)\big)
   \big(1+|i+j|+k\big)^d\, .
\]
\end{enumerate}
\end{lemma}

\noindent{\bf Proof: (1).}
For given \(\omega\in\fine\) and \(v\in\Z^d\) satisfying
\(v\cdot\vecell\in[u,u+3r]\), the set of
nearest-neighbour edges in \(\Z^d\) is partitioned in
Definition~\ref{d.prime} into discovered open, discovered closed and
undiscovered edges. When
\(R^v_{[\upwell],\downwell}(\omega')\in\fv[i,j,k]\), the path \(\mathcal{T}\)
supplied by Lemma~\ref{l.surfacepath} is the Tether: a \(K\)-box path of
bad boxes whose associated big boxes are undiscovered.

The box \(D[\mathcal{T}]\) lies in \(\overline{\downwell}[b_{\rm T}]\) with
relative latitude \(-i\), and hence may be chosen in at most
\(C A_{\rm f}(C(1+|i|))\) ways by Lemma~\ref{l.latcard}. Once \(D[\mathcal{T}]\)
is fixed, the vertical and projected displacement between
\(D[\mathcal{T}]\) and \(U[\mathcal{T}]\) gives
\[
   |\mathcal{T}|\, \geq \, c\bigl(|i+j|+k\bigr)-C\, .
\]
Indeed, adjoin the endpoint boxes to the Tether.  Each nearest-neighbour step
changes latitude by at most one, while, after translating each box to latitude
zero, it changes the normalised box index by a bounded amount, uniformly
in~\(\vecell\).  The endpoint differences in these two quantities are
\(|i+j|\) and~\(k\), respectively, and the displayed bound follows after an
additive constant for the two endpoint steps.  If the Tether is empty,
Lemma~\ref{l.surfacepath}(3) instead says that its endpoint boxes are equal or
adjacent, so \(|i+j|+k\) is bounded and the same estimate holds.

Put $s=|i+j|+k$ and $n_0=\max\{0,\lfloor cs-C\rfloor\}$.  By
Lemma~\ref{l.surfacepath}, a non-empty Tether is a self-avoiding bad-box path,
starting next to $D[\mathcal T]$, with at least $n_0$ boxes, and every edge in
the associated big boxes is undiscovered.  Badness of a box is determined by
the edge states in its associated big box.  Under
\(\Rlaw^v_{[\upwell],\downwell}\), these undiscovered edge states have their
original Bernoulli product law.  Thus, for any fixed candidate path whose
associated big boxes are undiscovered, the joint badness field along the path
has the same law as under~\(\PP\).  Lemma~\ref{l.percodominate}
therefore gives stochastic domination by Bernoulli site percolation of
parameter~$\epsilon$.  There are at most
$2d(2d-1)^{n-1}$ self-avoiding box paths of length~$n$ from a specified
starting box.  A union bound therefore gives
\[
 \sum_{n\, \ge \, n_0}2d(2d-1)^{n-1}\epsilon^n
\, \le \, C\big((2d-1)\epsilon\big)^{n_0}
\, \le \, C\big((2d-1)\epsilon\big)^{cs}\, .
\]
Multiplying by the number of possible choices of \(D[\mathcal{T}]\) proves the
first assertion.

{\bf (2).}
If \(\Rlaw^v_{[\upwell],\downwell}\big(\fv[i,j,k]\big)\) is non-zero, then
there must exist a pair \((D[\mathcal{T}],U[\mathcal{T}])\), the former in
\(\overline{\downwell}[b_{\rm T}]\) with relative latitude \(-i\), and the
latter in \(\overline{\upwell}[b_{\rm T}]\) with relative latitude \(j\), whose
projected distance is \(k\).  Decompose the upcluster translation vector
modulo $2K\Z^d$, giving at most $(2K)^d$ residue classes.  For each class
admitting a proposed head in the height window, choose one such placement
as representative.  The other placements in that class translate its
boxes by vectors in $2K\Z^d$.  Since the head height changes by at most
$3r$, the output latitude~$j$ restricts the box in the representative to a
bounded number of latitudes near~$j$, by~\eqref{e.relative-latitude-slack}.
Lemma~\ref{l.latcard} thus bounds the number of endpoint pairs in the
downcluster and this representative by
$C A_{\rm f}(C(1+|i|))A_{\rm f}(C(1+|j|))$.  Once these boxes have been chosen,
the requirements that their output relative latitudes be $-i$ and $j$ and that their
projected distance be~$k$ confine the translation index to an $\ell_1$-shell
contained in a box of side $C(1+|i+j|+k)$.  It therefore has at most
$C(1+|i+j|+k)^d$ possible values. Multiplication by the two endpoint counts
and the number of residue classes proves the second assertion. \qed

\smallskip\noindent
Part~(1) gives the exponential Tether estimate.  A head placement represented
by \((i,j,k)\) forces a Tether of length at least \(c(|i+j|+k)-C\).
Its bad-box path is controlled by Bernoulli site percolation of parameter
\(\epsilon\).  Writing \(q=(2d-1)\epsilon<1\), the connection chance therefore
falls like \(q^{c(|i+j|+k)}\), the \(e^{-c\kappa}\) of
Chapter~\ref{s.heuristics} with \(\kappa\) the Tether length.  A short Tether
allows only bounded separation of the endpoint boxes; as their separation
increases, the bound decreases geometrically.

\noindent{\bf Proof of Proposition~\ref{p.finesum}.}
For \(v\in\Z^d\) satisfying \(v\cdot\vecell\in[u,u+3r]\), occurrence
of \(\Accept_v\) implies occurrence
of \(\fv[i,j,k]\) for some \(i,j\geq-C_{\rm lat}\) and \(k\in\Nzero\), by
Lemma~\ref{l.surfacepath}. Thus
\[
\begin{aligned}
\sum_{\substack{v\in\Z^d:\\ v\cdot\vecell\in[u,u+3r]}}
   \Rlaw^v_{[\upwell],\downwell}(\Accept_v)                    
&\quad\, \leq \, \,
   \sum_{i,j,k}
   \sum_{\substack{v\in\Z^d:\\ v\cdot\vecell\in[u,u+3r]}}
      \Rlaw^v_{[\upwell],\downwell}\big(\fv[i,j,k]\big)        \\
&\quad\, \leq \, \, 
   C K^d\sum_{i,j,k}
      A_{\rm f}\big(C(1+|i|)\big)^2
      A_{\rm f}\big(C(1+|j|)\big)                              \\
&\qquad\qquad{} \, \times
      \big(1+|i+j|+k\big)^d
      \big((2d-1)\epsilon\big)^{c(|i+j|+k)}\, .
\end{aligned}
\]
The sum is over \(i,j\geq-C_{\rm lat}\) and \(k\geq0\).  Put
\(s=|i+j|+k\).  The lower bounds on \(i\) and~\(j\), together with
\(i+j\leq s\), imply \(|i|,|j|\leq s+C_{\rm lat}\).  Thus all allowance and
counting factors are bounded by \(C(1+s)^D\) for a fixed~\(D\), while the
number of triples with a given value of~\(s\) is polynomial in~\(s\).
Consequently the sum is finite when \(\epsilon\) is chosen sufficiently
small.  For the allowance in~\eqref{e.allowance-design},
$A_{\rm f}(n)\leq C(2+\Cfine)(1+n)^{\Bfine+1}$, so every allowance-dependent factor in the last summation is bounded by a
fixed power of \(2+\Cfine\); its exponential Tether factor is unchanged.
All constants in this estimate are uniform in~\(\vecell\): the endpoint
counts are uniform by Lemma~\ref{l.latcard}, the Tether displacement bound is
uniform by the axes-ordering inequalities used in Lemma~\ref{l.rfprob}, and
the bad-box domination is independent of the height direction.  This proves
both the uniform bound and \eqref{e.finesum-Cfine-polynomial}. \qed

We record this further consequence for later use in the iterated-slide
analysis.  For a retained pair \((\uclust,\dclust)\) arising
under~\(\fine_h\), write
\[
 \nu(v)=
 \PP^{\fine_h}\bigl(V=v\mid[\upwell]=\uclust,\downwell=\dclust\bigr),
 \qquad v\in\Z^d,\quad v\cdot\vecell\in[u,u+3r].
\]
Thus~\(\nu\) is the head-placement law in~\eqref{e.fixed-shift-head-law}.

\begin{corollary}[Tether summability of the head law]
\label{c.tether-head-summability}
Under the hypotheses of Proposition~\ref{p.finesum}, there are
\(C_{\rm T}<\infty\) and \(\gamma_{\rm T}>0\) such that, uniformly in the
retained pair,
\[
   \sum_{v:\,\nu(v)<t}\nu(v)
\, \le \, C_{\rm T}t^{\gamma_{\rm T}}\, ,
   \qquad 0<t\, \le \, 1\, .
\]
For the allowance in~\eqref{e.allowance-design}, \(\gamma_{\rm T}\) may be
chosen independently of \(\Cfine\), and
\begin{equation}
\label{e.tether-Cfine-dependence}
   C_{\rm T}\, \leq \, C(2+\Cfine)^{D_{\rm T}}
\end{equation}
for constants \(C,D_{\rm T}<\infty\) independent of~\(\Cfine\).  These
constants and~\(\gamma_{\rm T}\) may also be chosen independently of the
breaker-profile parameters~\(b_0,C_0,C_1,i_{\rm br}\).
\end{corollary}

\noindent{\bf Proof.}
Put
\[
 Z=\sum_{\substack{w\in\Z^d:\\ w\cdot\vecell\in[u,u+3r]}}
       \Rlaw^w_{\uclust,\dclust}(\Accept_w).
\]
By Lemma~\ref{l.qfine}, \(Z\geq c>0\).  For a placement~\(v\) with
\(\Rlaw^v_{\uclust,\dclust}(\Accept_v)>0\), let
\[
 q(v)=\min\bigl\{|i+j|+k:
 \Rlaw^v_{\uclust,\dclust}(\fv[i,j,k])>0\bigr\}.
\]
Occurrence of~\(\Accept_v\) entails one of the events~\(\fv[i,j,k]\), so
the minimum exists.  The two assertions of Lemma~\ref{l.rfprob}, together
with the polynomial growth of the allowance, give constants \(C,c,D>0\)
such that, for \(L\geq0\),
\[
 \sum_{q(v)>L}\Rlaw^v_{\uclust,\dclust}(\Accept_v)
 \leq Ce^{-cL},
 \qquad
 \#\{v:q(v)\leq L\}\leq C(1+L)^D.
\]
Dividing by~\(Z\) and splitting the sum according as \(q(v)>L\) or
\(q(v)\leq L\), the preceding two bounds give
\[
 \sum_{\nu(v)<t}\nu(v)
 \leq Ce^{-cL}+Ct(1+L)^D.
\]
Taking \(L=B\log(1/t)\), with fixed sufficiently large~\(B\), proves the
first assertion.

For the allowance in~\eqref{e.allowance-design},
\(A_{\rm f}(n)\leq C(2+\Cfine)(1+n)^{\Bfine+1}\).  The prefactors in the
two bounds from Lemma~\ref{l.rfprob} are therefore bounded by a fixed power
of~\(2+\Cfine\), while the exponential rate is unchanged.  The lower bound
on~\(Z\) depends only on \(\af(0)=A_1\).  The two bounds from
Lemma~\ref{l.rfprob} and this lower bound are independent of the
breaker-profile parameters~\(b_0,C_0,C_1,i_{\rm br}\).  Thus
\(\gamma_{\rm T}\) may be fixed independently of~\(\Cfine\) and of these
parameters, while the remaining prefactor satisfies
\eqref{e.tether-Cfine-dependence} with the same uniformity. \qed

\chapter{The iterated slide: concepts and setup}
\label{s.iterated-slide}


The first section gives a general view of the iterated slide schedule and of
the initially-very-fine to limelight-fine (IVF-to-LLF) transition on which renewal production
depends.  The second sets out the iterated slide formally, with the operator,
the randomised schedule, and the IVF and LLF index sets defined precisely.
The third offers a roadmap for the control of dynamic cluster accrual
and for the route from this control to the proof of
Theorem~\ref{t.intermediate-renewal}(2) in
Chapter~\ref{s.renewal-survival}.

\section{Overview of the iterated slide}
\label{ss.iterated-slide-overview}

In seeking to prove that renewal absence along the entire height range~$[0,u]$ has probability that decays exponentially in~$u$
under $\Q = \PP \big(  \cdot \big\vert \clusteruthree \big)$, we will attempt slide after slide, along an increasing schedule of heights $\big\{ h_i : i \in\intint{L} \big\}$, with the total number $L$ of attempted slides being deterministic of order~$u$. The operation
$$
 \omega_0 \mapsto \omega_1 \mapsto \cdots \mapsto \omega_L
$$
preserves the law of~$\Q$ at each single-slide step, so that the composition $\omega_0 \mapsto \omega_L$ sending input to output also preserves~$\Q$.
When precisely specified, this map will be the iterated slide, our principal resampling tool.

When a single slide was performed at height~$h$, the event $\fine_h$ acted as a test of input regularity and served to ensure a uniformly positive conditional probability of producing a renewal point, via simple join. As we plan the slide operation in its iterated form, we are concerned not only with securing such a favourable outcome for individual slides, but also with containing the damage when simple join fails to occur. Indeed, any active slide will relocate the concerned upcluster by a shift, with a new midpiece forming in order to connect the upcluster bed with the downcluster top. If for example this midpiece travels far to the south of the present slide height, it may undo renewal secured earlier as the slide height rose up the system.

This need for damage control leads us to insist that, when the slide reaches a given height in the schedule, the proposed slide will occur only if, in its present configuration, the fine event is verified at that height. In the terminology of Definition~\ref{d.slide-operator}, the step is then active; otherwise it is dormant and leaves the configuration unchanged. Proposition~\ref{p.finesum} bounds the connection-chance sum of every active slide by a uniform constant. The Tether estimate underlying this bound sharply limits the effective range of possible head placements. Consequently, the later damage analysis may be carried out in fixed-shift experiments---where the remaining undiscovered edges form an independent percolation field---and then summed over the possible placements.

A schedule of order~$u$ heights is proposed in order to offer such an order of opportunities to produce renewal. If renewal is produced often enough, and is not often later undone, the probability that the output~$\omega_L$ has no renewal may be shown to be at most $C e^{-cu}$; since~$\omega_L$ has law~$\Q$, this gives the estimate required for Theorem~\ref{t.intermediate-renewal}(2). By performing slides only on auspicious occasions, we may however face a supply problem for renewal production opportunity. A given level in the proposed height schedule awaits the arrival of the rising slide height. That moment of arrival is the level's  {\em limelight} time, when the level may index an active slide. But this moment will be enjoyed only if the level is {\em limelight fine}, with the configuration at the moment of the proposed slide verifying the fine condition at that height. By Definition~\ref{d.fvf}, this means that every condition ${\bf A}^{-}_{\rm f}(k)$ and ${\bf A}^{+}_{\rm f}(k)$, $k\in\Nzero$, and the upcluster breaker condition ${\bf B}^{+}_{\rm f}$ hold at that moment.

Renewal production opportunities are in adequate supply only if they number~\(\Theta(u)\) with overwhelming probability. How to ensure this? Recall that a resampling \(\omega\mapsto\omega'\) seeks to retain geometric information from the input while using fresh randomness to produce favourable output geometry.  Here this balance must persist dynamically.  Although \(\omega_0\mapsto\omega_L\) is the full composite resampling, for the \(i\)\textsuperscript{th} scheduled height \(h_i\) the relevant evolution is \(\omega_0\mapsto\omega_{i-1}\): the useful regularity present at \(h_i\) in the input must survive, in weakened form, until its limelight time.

The event $\veryfine_h$ of Definition~\ref{d.fvf} is the strengthened input condition.  For the allowances fixed in Section~\ref{s.onetime},
\(
   \avf(k)=A_0(1+k)^3
\)
and the split fine allowance is given by~\eqref{e.allowance-design}.  Thus the
very-fine allowance is cubic, the depth-zero fine allowance is the separately
fixed constant~\(A_1\), and the positive-depth fine allowance grows as a
higher power of the depth, with the coefficient~\(\Cfine\) left to be fixed in
the iterated argument.

By Proposition~\ref{p.veryfinelevels}, the very-fine event occurs at a positive density of integer heights with overwhelming probability under~$\Q$.
We will apply this form of regularity to the input~$\omega_0$, identifying a
positive density of attempted slide heights that are {\em initially very
fine} (IVF).  The renewal production supply problem is then reduced to
arguing that, with overwhelming probability, a positive fraction of these
heights are {\em limelight fine} (LLF), verifying the fine condition at the
unique moment that matters, when the slide arrives.

Establishing this \hfff{fine-delivery}{IVF-to-LLF transition or \emph{fine delivery}} is delicate, because it entails tracking many cluster-bound conditions at each element in a positive-density set of scheduled heights: for a given such height, active slides earlier in the schedule may individually or collectively threaten these conditions by adding downcluster material during the run-up to limelight.  We control this accrual in two ways.  Wavebreakers attenuate the influence of slides far below the target height, while single-layer accrual bounds control the material contributed by nearby slides.  The separation between the very-fine and fine allowances is chosen to absorb the resulting accrual with high probability.

This suggests taking \(\Cfine\) large, thereby giving the positive-depth fine conditions more room.  The extra room has a cost when one takes account of the acceptance conditioning along the slide dynamics: the resulting estimates become less favourable as \(\Cfine\) grows, a dependence already visible in the SimpleJoin estimate~\eqref{e.simplejoin-Cfine-lower}.  We compensate by strengthening each breaker row, increasing the number \(i_{\rm br}\) of thin test layers it contains; the wavebreaker analysis will show that finite passage through such a row then becomes exponentially unlikely in \(i_{\rm br}\).  Proposition~\ref{p.veryfinelevels} is uniform in the row length and separation, so lengthening the rows does not reduce the proved density of very-fine integer heights.  It does require the attempted-slide heights to be spaced farther apart so that the resulting breaker windows remain disjoint.  Since these scales are fixed independently of \(u\), this leaves a smaller but still positive density of renewal-production opportunities.  The quantitative balance is carried out when fine delivery is proved in Chapter~\ref{s.shadow-bounds-delivery-fine}.

\medskip
\noindent{\em The height schedule at positive or lower density.}  Our schedule has positive density, whereas~\cite{FriberghHammondPreprint} considered one of order~$u^{1/2}$. 
In the approach of the arXiv preprint, midpieces are forced to intersect typical integer bands in bounded cross-sections; the requisite control is furnished at most heights by the typical-band estimate derived from the global surface bound, Proposition~\ref{p.nobigcluster}.
The passage of a midpiece between consecutive schedule heights
is then shown to have probability at most of order~$e^{-cu^{1/2}}$, in a manifestation
of the {\em thin, instantly perish} aspect of TipToe (in the language of the
present work).  In this sense, any interaction between slides is treated as a
global catastrophe contributing to the upper bound on the probability of a
renewal-free conditioned connection.  It seems plausible that, with some
refinement, such an argument could treat a schedule whose gaps are of
order~$\log u$, with a resulting upper bound of the form~$e^{-cu/\log u}$ on the
$\Q$-probability of renewal absence; and a further bootstrapping argument might
even replace the logarithm by a function of much slower growth.  Theorem~\ref{t.intermediate-renewal}(2), and the
resulting mass gap, would however require that this function be replaced by a
uniformly bounded one.  The detailed study that we undertake of the
interaction between slides on the positive-density schedule is needed in order to reach this conclusion.

\section{Specifying the procedure}

We now make precise the IVF-to-LLF transition described in
Section~\ref{ss.iterated-slide-overview}.  An IVF height will be a sampled
height at which the input state is very fine.  The same height is LLF when the
state immediately before its scheduled turn is fine; it then indexes an active
slide.  Thus fine delivery is a statement about the evolution between two
specified times.  We next define the iterated operator, the randomised schedule,
and these two conditions on its sampled heights.

\subsection{Candidate heights and the iterated operator}
\label{ss.candidate-heights}

Let \(\overline h=(h_i:i\in\intint{L})\) be an increasing schedule of attempted slide heights, with $0 \leq h_1 < h_L \leq u$.  At height
\(h\), write \(\slide_h\) for the single-slide Markov kernel.  If the
current configuration is not in \(\fine_h\), this kernel acts as the
identity.  Write \(\slide_h(\omega)\) for a random configuration whose law is
given by \(\slide_h(\omega,\cdot)\).  When this notation is iterated,
independent randomness is used in successive applications.

For such a sequence, set the \emph{iterated operator}
\[
   \iter_{\overline h}
\, = \,
   \slide_{h_L}\circ\cdots\circ\slide_{h_1}\, ,
\]
where composition is understood in this stochastic sense.  More explicitly,
if the input configuration is \(\omega_0\in\clusteruthree\), the procedure
constructs configurations \(\omega_1,\ldots,\omega_L\) recursively.  Given
\(\omega_{i-1}\), the \(i\)\textsuperscript{th} step applies the single-slide
operator of Definition~\ref{d.slide-operator} at height \(h_i\), so that
\(\omega_i=\slide_{h_i}(\omega_{i-1})\).  If
\(\omega_{i-1}\notin\fine_{h_i}\), the step is \emph{dormant} and
\(\omega_i=\omega_{i-1}\).  If \(\omega_{i-1}\in\fine_{h_i}\), the step is
\emph{active}: the downcluster \(\downwell_{h_i}(\omega_{i-1})\) and the
shift-equivalence class \([\upwell_{h_i}(\omega_{i-1})]\) of the upcluster are
held fixed, in the sense of Lemma~\ref{l.probform}, and the remaining
configuration is reconstructed by the connection-chance law there.  Because
\(\slide_{h_i}\) is a Markov kernel, given the schedule and conditionally on
\(\omega_{i-1}\), the law of \(\omega_i\) is
\(\slide_{h_i}(\omega_{i-1},\cdot)\), depending on the earlier configurations
only through \(\omega_{i-1}\).  Finally,
\(\iter_{\overline h}(\omega_0)=\omega_L\).

For the formal iteration, the cluster may act as the state variable.  Let
\(\mathfrak C_u\) be the space of finite clusters arising from configurations
in~\(\clusteruthree\), and write \(\pi(\omega)=\csetuthree(\omega)\).  The cluster
\(\pi(\omega)\) determines its open edges and incident closed boundary, and
hence the downcluster, the shift class of the upcluster, and the events
\(\fine_h\) and \(\veryfine_h\).  For \(C\in\mathfrak C_u\) and measurable \(A\subseteq\mathfrak C_u\), set
\begin{equation}
\label{e.cluster-slide-kernel}
 \overline\slide_h(C,A)
 =\slide_h\bigl(\omega,\pi^{-1}(A)\bigr),
 \qquad \pi(\omega)=C.
\end{equation}
The right-hand side is independent of the choice of \(\omega\) satisfying
\(\pi(\omega)=C\): Lemma~\ref{l.fine-retained} determines the active or dormant
status from the retained pair, and Definition~\ref{d.slide-operator} uses an
active input only through the retained downcluster and translated upcluster.
At an active slide the old undiscovered edge states are discarded, and the
reconstruction resamples them from a fresh product field.  Thus
the law of the next cluster depends on the present configuration only through
its present cluster, and \(\overline\slide_h\) is a Markov kernel on
\(\mathfrak C_u\).  Moreover, if
\(\overline\Q=\Q\circ\pi^{-1}\), Proposition~\ref{p.single-slide-invariance}
gives
\[
 \int\overline\Q(dC)\,\overline\slide_h(C,A)
 =\overline\Q(A).
\]

Put \(S_i=\pi(\omega_i)\) and
\[
 \mathcal K_i
 =\sigma\bigl(S_0,\overline h,S_1,\ldots,S_i\bigr).
\]
Conditionally on \(\mathcal K_{i-1}\),
\begin{equation}
\label{e.cluster-trajectory-markov}
 \operatorname{Law}(S_i\mid\mathcal K_{i-1})
 =\overline\slide_{h_i}(S_{i-1},\cdot).
\end{equation}
Thus the projected clusters form a Markov chain in their own right.

\subsection{The random schedule and dynamic law}
\label{ss.random-mesh}

The attempted height schedule used in the iterative procedure will be randomised independently of the input configuration and other randomness introduced in resampling during its steps.
This randomisation ensures that a macroscopic number of very fine input heights
enter the attempted slide schedule.

With the sufficiently large even coarse scale~\(K\) fixed in
Section~\ref{s.coarse}, set
\[
   \beta_n=K(2n+1/2),\qquad n\in\mathbb Z.
\]
Thus \(\beta_{n+1}-\beta_n=2K\). These heights delimit the layers in the
partition  introduced  in  Section~\ref{s.slab-setting} for use in monitoring cluster accrual. 
Since attempted slides occur
only at such heights, every physical slide height is a face of the partition and
no physical slide cut lies in the interior of a partition layer.  

Let \(m,g\in\N\) be schedule parameters, to be chosen independently of~\(u\),
with \(m/(m+g)>1/2\).  They are fixed in
Chapter~\ref{s.shadow-bounds-delivery-fine}.  Let \(n_0\) be the least integer
such that \(\beta_{n_0}>2r\).  For \(i\geq1\), set
\[
  \mathsf I_i
\, = \, \{\beta_{n_0+(i-1)(m+g)},\ldots,
             \beta_{n_0+(i-1)(m+g)+m-1}\}\, .
\]
Let \(L\) be the greatest positive integer for which the last member of
\(\mathsf I_L\) is less than \(u-2r\).  Thus each
\(\mathsf I_i\) contains exactly \(m\) consecutive candidate heights and exactly
\(g\) candidate heights lie strictly between two consecutive intervals.  From
each interval \(\mathsf I_i\), choose one height \(H_i\) uniformly and
independently.  The random vector \(\overline H=(H_i:i\in\intint{L})\) is
the sequence of attempted slide heights, and the slides are tested in the
order \(H_1,\ldots,H_L\).

The integer \(g\) records the separation of successive sampling blocks:
exactly \(g\) candidate heights are left unused between them.  We require
\begin{equation}
\label{e.random-mesh-gap}
   g\, \geq \, \big\lceil\max\bigl(3r + b_0 + C_0,\;4r\bigr)\big\rceil + 1\, .
\end{equation}
Here \(r=e_d\cdot\vecell\), while \(b_0\) and \(C_0\) are the offset and
window-length parameters in the northside very-fine breaker condition of
Definition~\ref{d.fvf}.  This requirement will serve to keep the breaker
windows associated with distinct attempted heights separate.  Since \(r\leq1\),
a single direction-independent choice of \(g\) is possible, for example by
imposing the stronger bound obtained by replacing \(r\) by~\(1\) in~\eqref{e.random-mesh-gap}.

In a dynamic version of the resampling viewpoint of
Section~\ref{s.resampling-overview}, we record the successive cluster states
along the iterated slide.

\begin{definition}[Dynamic law]
\label{d.qdyn}
Let \(\overline H=(H_1,\ldots,H_L)\) be the random schedule just specified,
independently of \(S_0\sim\overline\Q\).  Write
\hfff{dynamic-law}{\(\Qdyn\)} for the joint law of the random attempted-height
schedule \(H_1,\ldots,H_L\) and the successive cluster states
\(S_0,S_1,\ldots,S_L\), where the latter are formed by successive application
of the kernels \(\overline\slide_{H_i}\).
\end{definition}

Thus under~\(\Qdyn\), \(S_i\) is the cluster state after the
\(i\)\textsuperscript{th} step.  Expressions such as
\(\downcluster_h(S_i)\), \(\fine_h(S_i)\) and \(\csetuthree(S_i)\) are
interpreted in this state.

\subsection{Sampling IVF heights}
\label{ss.sampling-ivf-points}

Proposition~\ref{p.veryfinelevels} supplies a macroscopic abundance of very
fine heights in the input. The random schedule converts this abundance into
a macroscopic supply of attempted heights that are very fine at time zero.

Set $\ivf\, = \, \{i\in\intint{L}:S_0\in\veryfine_{H_i}\}$.

For \(i\in\ivf\), the attempted height \(H_i\) is called {\em initially very fine},
or {\em IVF}.

Let \(\gamma\in(0,1)\).  Given the sampled schedule, we say that the input is
\(\gamma\)-{\em satisfactory} if $|\ivf|\, \geq \, \gamma L$..

\begin{lemma}
\label{l.sampling-ivf-points}
There exist positive constants \(\gamma,c\) and \(C<\infty\), depending only
on the fixed schedule parameters and on the constants in
Lemma~\ref{l.nhigh} and Proposition~\ref{p.nobigcluster}, but not on the
realised schedule, such that
\[
   \Qdyn\big(\hbox{the input is not \(\gamma\)-satisfactory}\big)
\, \leq \, C\exp\{-cu\}\, .
\]
One may take \(\gamma=1/24\) for every fixed choice of schedule parameters
satisfying \(m/(m+g)>1/2\).
\end{lemma}

\noindent{\bf Proof.}  Let \(Q\) be the constant of Lemma~\ref{l.nhigh} and \(D_0\) the
threshold of Proposition~\ref{p.nobigcluster}.  By
\eqref{e.vf-abundance-choice}, \(QD_0<1/(8K)\), while the schedule parameters
were chosen above with \(m/(m+g)>1/2\).

Successive sampling blocks begin \(2K(m+g)\) height units apart, while the
endpoint exclusions have bounded size.  Hence
\[
   \lambda_u:=\frac{mL}{u}
   \, = \, \frac{m}{2K(m+g)}
        +O\!\left(\frac1u\right)
   \, > \, \frac1{4K}
\]
for all sufficiently large~\(u\).  Thus \(QD_0<\lambda_u/2\).

By Proposition~\ref{p.nobigcluster}, there is an event \(\Omega_u\) of input
clusters such that
\[
 \overline\Q(\Omega_u^c)\, \leq \, C\exp\{-cu\}
 \quad\hbox{and}\quad
 |\partial_{\rm ext}\overlinecsetuthree|\, \leq \, D_0 u
 \quad\hbox{on }\Omega_u\, .
\]
On \(\Omega_u\), Corollary~\ref{c.nhigh-subset}, applied to the deterministic
set \(J=\bigcup_i\mathsf I_i\), shows that the number of non-very-fine
candidate heights in the sampling union is less than
\(mL/2\).  Hence at least \(mL/2\) of the candidate heights in
\(\bigcup_i\mathsf I_i\) are very fine.

Call \(\mathsf I_i\) {\em grim} if fewer than \(m/4\) of its \(m\)
candidate heights are very fine.  If a proportion \(p\) of the \(L\)
intervals were grim, the total number of very fine heights in the union
would be at most \((1-3p/4)mL\).  Since this total is at least \(mL/2\),
we have \(p\leq2/3\).  Thus at least \(L/3\) intervals are not grim.

Given any input cluster state in~\(\Omega_u\), the choices from different
sampling blocks are independent, and in each non-grim block the selected
height is very fine with probability at least \(1/4\).  Therefore
\(|\ivf|\) stochastically dominates
\(\mathrm{Bin}(\lceil L/3\rceil,1/4)\).  A standard lower-tail estimate
consequently shows that, with \(\gamma=1/24\) and suitable \(c>0\), the
conditional probability of \(|\ivf|<\gamma L\) is at most
\(\exp\{-cL\}\).  Adding the contribution from \(\Omega_u^c\), and using
\(L=\Theta(u)\), completes the proof. \qed

\subsection{Limelight fine heights}
\label{ss.breaker-llf-points}

Set
\[
   \llf\, = \, \{i\in\intint{L}:S_{i-1}\in\fine_{H_i}\}\, .
\]
For \(i\in\llf\), the attempted height \(H_i\) is called {\em limelight fine}, or {\em LLF}.
Thus membership of \(\ivf\) records the very fine condition at time zero,
whereas membership of \(\llf\) records the fine condition at the index's
\hfff{limelight}{limelight time}, immediately before its slide is attempted.
\label{gl.ivf-to-llf-transition}

\subsection{The output law}
\label{ss.qdyn}

It remains to verify the cluster marginal of the dynamic law.

\begin{proposition}
\label{p.qdyn-output}
Under~$\Qdyn$ the marginal of \(S_i\) at every stage is
\(\overline\Q=\Q\circ\pi^{-1}\), the cluster law induced by~\(\Q\).
In particular \(S_L\) has this law.
\end{proposition}

\noindent{\bf Proof.} Fix a realisation \(\overline H=(h_1,\ldots,h_L)\).  By independence, the
conditional law of \(S_0\) is~\(\overline\Q\).  Suppose inductively that the
conditional law of \(S_{i-1}\) is~\(\overline\Q\).  For every bounded
measurable function \(f\) on \(\mathfrak C_u\),
\[
 \EE_{\Qdyn}\bigl[f(S_i)\mid\overline H\bigr]
 =\int\overline\Q(dC)
      \int\overline\slide_{h_i}(C,dC')f(C')
 =\int f\,d\overline\Q .
\]
The second equality follows from
\eqref{e.cluster-slide-kernel} and
Proposition~\ref{p.single-slide-invariance}; the kernel already includes
the active and dormant cases.  Induction and averaging over
\(\overline H\) prove the assertion.
\qed

\section{Roadmap for the iterated slide argument}
\label{ss.roadmap-iterated}

Part~\ref{part.two} develops tools for understanding  the slide dynamics and Part~\ref{part.three} applies them. Together, these parts form the lengthy central passage of the paper  by
which we will obtain the renewal existence estimate of
Theorem~\ref{t.intermediate-renewal}(2) in
Chapter~\ref{s.renewal-survival}. We present a guide now.

By Proposition~\ref{p.qdyn-output} the output cluster \(S_L\) has the
cluster law induced by~$\Q$.  Theorem~\ref{t.intermediate-renewal}(2) is a
statement about this law---that the cluster contains a renewal point with probability at
least $1-Ce^{-cu}$---so the question before us is how to find a renewal
in~\(S_L\) with overwhelming $\Qdyn$-probability.

A target height~$\ell$ is limelight fine when the fine condition holds there at
the moment the schedule reaches it.  The part of this condition which may be
lost during the run-up to limelight is the set of southside band-count
constraints in Definition~\ref{d.fvf}.  
For each \(k\geq0\), the number of vertices of the target-height downcluster \(\downcluster_\ell\) in the unit band a distance~\(k\) below~\(\ell\) must be at most the fine allowance \(\af(k)\) in \({\bf A}^{-}_{\rm f}(k)\). When~\(\ell\) is initially very fine, the corresponding stronger input condition \({\bf A}^{-}_{\rm vf}(k)\) imposes the bound \(\avf(k)\).
 The danger is that material added to
$\downcluster_\ell$ by active slides below~$\ell$ exhausts the 
slack  \(\af(k) - \avf(k)\) before the schedule reaches~$\ell$.  We must therefore control the
accrual slide by slide.

Fix \(k\geq0\), and place the unit band at depth~\(k\) below~\(\ell\) inside a parallel coarse layer.
The layer has bounded width in the $\vecell$-direction, is infinite in the
transverse directions, and is chosen wide enough that percolation restricted
to it is supercritical.  The band is where accrual is counted, while the containing layer
is the fixed environment in which percolation-based arguments control it.  We
monitor the restriction of the whole conditioned cluster to this layer as the
slide schedule rises up to and then beyond the layer, stopping when the target
height~$\ell$ is reached.  This monitored graph contains the restriction of the
target-height downcluster to the layer as a subgraph.  If the material
accumulated after the input remains within the separation
between the very-fine and fine allowances, the southside fine condition at
depth~$k$ holds at limelight.

We first need a precise account of how one slide changes the cluster seen at a
different height.  In the single-slide analysis of Section~\ref{s.onetime}, the
cluster was divided at the slide height into a downcluster, a midpiece and an
upcluster.  Chapter~\ref{s.dmu} formalises this DMU decomposition and compares
the decompositions at two heights.  This comparison identifies the material
that a lower slide may add to the target-height downcluster.
It also identifies a powerful \emph{cut-it-and-lose-it} mechanism.  If a slide
south of the layer deposits midpiece material there, the next active slide
still south of the layer cuts through that material and removes it from the
layer.

During the approach of the rising slide schedule to the target height~$\ell$,
there is a distinguished moment of {\em passage}, at which the slide height
crosses the given layer (the physical attempted slide heights lie on the
boundaries of the partition layers).  The evolution in the layer has two phases: an advance phase before
passage, and a post-passage phase in which the slide height runs north of the
layer up to~$\ell$.

In the advance phase, the part of the input cluster lying in the upcluster above
the rising slide height is iteratively transported by the successive slide
shifts.  As the height schedule rises towards the layer, the intersection of
this transported material with the layer shifts from one active slide to the
next; at some slides, part of it may be lost as material exits at the rising
cut.  We call this transported part of the live cluster the \emph{host}.  It
may be accompanied by a contribution from the currently sampled midpiece, the
\emph{guest}; each such guest is transient in the advance phase, since
cut-it-and-lose-it removes it at the next active slide.  At passage, the
cluster material present in the layer therefore 
consists of the surviving host
together with the guest, if any, created by the last active slide in the advance phase.

This cluster material evolves one slide at a time.  Consider one step.  The
open-and-closed information already determined in the fixed layer consists of
finitely many discovered open edges together with the closed boundary edges
exposed in determining that structure.  The slide acts at a given height,
above or below the layer.  Condition on the present downcluster and upcluster shift class, and
fix a proposed head location~$v$.  Under the fixed-shift law
$\Rlaw^v_{\uclust,\dclust}$ of Section~\ref{s.onetime}, the undiscovered edges are
sampled from an independent percolation field. The new midpiece may be
explored in stages delimited by the southside and northside of the layer:
outside, then through the layer, and outside again, until the proposed
connection has been decided.  Monitoring only the portions inside the layer
gives a sequence of percolation explorations whose starting points may depend
on the information revealed earlier in the same slide.

The slide is ultimately accepted only when this fixed-shift exploration
satisfies the event \(\Accept_v\).  It is natural,
however, first to ask how accrual may arise before this conditioning is
imposed.  Significant finite accrual can occur only
to the extent that the discovered open and closed record shields further
starting points from the infinite open cluster of the layer.  Some territory
may be enclosed completely; elsewhere the record may merely make connection to
the infinite component unlikely.  An adaptive exploration could in principle
visit such regions successively and accumulate substantial finite volume
without encountering open escape.  We therefore need an upper bound on the
shielding capacity of an arbitrary finite discovered record in the layer.

The means by which we prove this shielding estimate is the \emph{seal}, whose
properties are developed in Chapters~\ref{s.overlapping-quilts-seals}
and~\ref{s.seed-exploration-seal-extension}. 
In the first of these chapters, we work within the physical layer, the slab orthogonal to~\(\vecell\), and assign its edges to connected boxlets cut out by transverse axial columns.
 Each coarse box is the union of \(2^{d-1}\) neighbouring
boxlets, and neighbouring boxes overlap by one boxlet-width, or fifty per cent,
in each transverse coordinate.  With the box width fixed sufficiently large,
most boxes satisfy an Antal--Pisztora--type condition.  A connected collection
of such boxes forms a \emph{quilt}; the macroscopic open pieces in their shared
boxlets join to form a connected open fabric.  An outward search from the
finite discovered record locates a surrounding quilt collection without
examining the edges beyond it.  When this collection is selected as a stopping
domain and captures the record in its finite interior, it is a seal.

A fundamental property is that the uninspected part of the layer beyond a
seal retains its conditional product law.  When the fixed-shift exploration
enters the unread exterior, it pauses at a new seed vertex before querying
the first unread edge.  That edge determines a wholly untouched connected
boxlet whose remaining edges are unread.  From such a seed there is a uniformly positive
probability of connection to the infinite open component of the layer.  On
the complementary finite-cluster outcome, the volume newly explored in the
layer has a stretched-exponential tail.  A slide may visit the layer several
times and choose later seeds adaptively; the same alternative is applied at
each visit.

Chapter~\ref{s.seed-exploration-seal-extension} develops the sealing machinery
used as the slide schedule advances.  Before passage, the retained upcluster
record shifts as the cut rises, while the current midpiece contribution is cut
and resampled at each subsequent active slide during the advance phase.  At
the last such slide, the surviving retained record together with the current
midpiece record forms the finite open-and-closed passage record in the layer.
Two sealing operations enter here.  A one-off, comparatively expensive
\emph{enclosure} places an initial seal around the retained record.  \emph{Seal
enlargement}, within the same fixed-shift configuration, then deforms the seal
outward as needed to incorporate whatever part of the midpiece record lies
beyond its present capture.  The resulting passage seal traps the finite record
with product exterior. It forms the intiial condition for a process of post-passage seal growth.

Indeed, after passage the slide height rises towards the target~\(\ell\) and each active
slide finds a seal enclosing the cluster material already discovered in the layer.
This material is retained, but the fresh resampling may disturb the woven fabric on
which the seal depends. The enclosing task remains, but the old seal no longer
reliably discharges it. In its place, we build a fresh seal farther out, in the
virgin product environment supplied by the new fixed-shift field, surrounding the
previously captured region and taking over the old seal's role. Exploration in the
layer can then proceed, and any further finite cluster material discovered is
absorbed by enlarging the new seal outward. This pattern of replacement and
enlargement repeats as the schedule advances towards limelight.

The advance phase may contain order~\(u\) attempted heights, so an
unattenuated budget for every slide in this phase is unavailable.
Cut-it-and-lose-it reduces the actual passage record to the last active slide
of the advance phase, but the index of that slide is not a stopping time.
Chapter~\ref{s.seal-sequence-layer} therefore prepares a possible passage seal
at each active advance slide whose midpiece reaches the layer.  At passage,
the candidate from the last active advance slide is selected if one was
prepared; otherwise an initial seal is constructed around the surviving host.

The required attenuation is supplied by the wavebreakers developed at the end
of Part~\ref{part.two}, in Chapter~\ref{s.wavebreakers}.  The IVF heights
furnish a positive density of breaker windows, and intervening windows make it
geometrically unlikely that the fresh midpiece from a far-below candidate for
the last active slide reaches the monitored layer.  This decay makes it
possible to prepare and control candidate passage seals throughout the
advance phase: the farther below the layer the slide lies,
the less likely its midpiece is to reach the layer, and the resulting
probabilities are summable over the advance phase.

The elements treated so far concern one layer and one target to its north.
The southside conditions must instead be controlled over many layer--target
pairs at once, and the layers cannot be analysed separately: a single
proposed slide may visit several of them, and its effects in different layers
are dependent through the shared shift and midpiece.
Chapter~\ref{s.seal-sequence-layer} therefore gathers the contributions of one
proposed slide across the monitored layers---fresh finite volume, new seal
boundary created when finite exploration in a layer is enclosed, and any
prospective passage preparation at that slide---into a single non-negative
random variable, and proves a conditional exponential-moment bound that can be
iterated along the schedule.

Chapter~\ref{s.model-problems-accrual} then
follows the successive fresh seals through the post-passage slides.
Enlargement within one fixed-shift configuration may move part of the seal
boundary outward as new finite material is incorporated.  At the next active
slide the seal must be replaced in a fresh field around the inherited capture,
and the new seal may lie farther out again.  The succession of captures
resembles a Russian doll: each fresh seal encloses the region captured at the
preceding stage, while newly incorporated finite material may push its boundary
farther out.  Chapter~\ref{s.model-problems-accrual} shows that the resulting
outward front advances at linear speed through the replacements, up to an
exponentially tailed excess; the new boundary sites created by finite
enlargement have already been charged to the one-slide cost in
Chapter~\ref{s.seal-sequence-layer}.  
For the layer containing the depth-\(k\) band, only order~\(k\) such replacements intervene before the target reaches limelight, giving the required depth-dependent control of the captured region.
 A single proposed slide may
affect several monitored layers, with the shared shift and midpiece coupling
the resulting contributions.  The analysis therefore treats these
contributions jointly at each slide.  Attenuation from intervening breaker
windows, together with weights decreasing with target depth, makes the
resulting estimates summable over all monitored layers.

The remaining probabilistic organisation is carried out through the reference
dynamics of Chapter~\ref{s.model-problems-accrual}, which uses the same
proposal law as the accepted dynamics but runs each corresponding fixed-shift
experiment without conditioning on acceptance.  Two principles govern this
organisation.  First, the auxiliary sealing construction must be invisible to
the cluster dynamics.  Sealing may read edges that the cluster exploration
itself would not read.  For a fixed proposed head~\(v\),
Chapter~\ref{s.seal-sequence-layer} therefore proves that projecting away the
auxiliary sealing data recovers the ordinary fixed-shift cluster experiment.
Proposals in the reference dynamics are chosen from the marked cluster history
independently of the auxiliary sealing coordinate, so these extra revelations
do not affect later proposal choices.  Conditioning the fixed-shift cluster
experiment on acceptance then gives the accepted transition.

Second, the accumulated acceptance correction is applied globally only once.
The one-slide estimates used along the schedule are proved before conditioning on acceptance and are then iterated under the reference dynamics to obtain the joint multi-layer accrual bound of Corollary~\ref{c.reference-multilayer-accrual}.
Chapter~\ref{s.shadow-bounds-delivery-fine} then uses this bound to control
positive-depth southside failure through a tower-and-shadow scheme under the
reference dynamics.  At this stage, the trajectory-level change of measure,
controlled by the Tether, transfers the resulting aggregate failure estimate
to~\(\Qdyn\).

The same chapter completes fine delivery.  A positive-density set of
well-spaced targets is selected from the IVF heights, each with a linear supply
of earlier breaker windows behind it.  At depth zero, where no post-passage
phase exists, southside failure is treated separately by tracking directly the
reach of the guest supplied by the last active slide.  The constants of the
positive-depth and depth-zero arguments are chosen together, and at least half
of the selected targets satisfy every fine condition at limelight.  This is
the IVF-to-LLF transition, Theorem~\ref{t.vf-to-fine-delivery}.
Every limelight-fine height indexes an active slide.  The one-slide
simple-join estimate gives a uniformly positive conditional probability of
renewal production there.  Applying this conditional bound successively
along the attempted schedule gives a positive density of produced renewal
points with overwhelming probability.

It remains to show that at least one produced renewal survives the later
slides.  At each active slide, a thin three-row breaker window lies just to its
north.  The slide's own midpiece first encounters the lower row and is
likely to be arrested there, before reaching the middle and upper rows.
Those two rows are then available as breaker units standing between the
renewal currently being followed and the schedule as it advances to the
north.  The protection is asymmetric: later
attacks come from the north, and no subsequent cut is guaranteed to remove a
midpiece that has passed a breaker, so every tested unit is treated as
single-use.  As the schedule advances, southgoing midpieces remove units
from the northern end of the queue, while new active slides typically add a
fresh pair farther north.  Chapter~\ref{s.renewal-survival} makes this moving
queue precise and shows that it has positive drift.  If the queue is
exhausted, the tracked renewal is abandoned and the construction restarts at
a later renewal-producing slide.  With overwhelming probability, some
restarted queue remains non-empty to the terminal cluster~\(S_L\), leaving its
associated renewal intact.  Since \(S_L\) has the cluster law induced by~\(\Q\),
Proposition~\ref{p.qdyn-output} completes the proof of
Theorem~\ref{t.intermediate-renewal}(2).

This completes Part~\ref{part.three} and the central slide argument, delivering to Part~\ref{part.four} the renewal-existence estimate that leads to the mass gap and the Ornstein--Zernike conclusions
by means of renewal analysis.

\part{Tools for the Slide Dynamics}\label{part.two}

\chapter{The down-mid-up decomposition}
\label{s.dmu}

Quilts, seals and wavebreakers are tools to analyse the slide dynamics that are the respective objects of attention in the later chapters of Part~\ref{part.two}. Understanding the dynamics  requires
a basic study in the present chapter of how a pair of slides may interact.

The \hfff{dmu}{\emph{down--mid--up}, or DMU, decomposition} splits the static cluster at a
given height into a downcluster, a midpiece and an upcluster. In a first section,
we set up the decomposition and compare it across heights to obtain the
cut-it-and-lose-it mechanism. In a second, we study DMU along a slide
sequence, noting that higher upclusters are rigidly transported by lower slides,
and how successive downcluster increments assemble the output edge set of an
active increasing slide schedule.

\section{DMU geometry across heights}

\subsection{The decomposition at one height}
\label{ss.dmu-one-height}

On the event \(\clusteruthree\), the finite open connected component
\(\csetuthree\) has base \(0\) and height $\hmax := \max \big\{ v
\cvl: v \in V(\csetuthree) \big\}$ belonging to \([u,u+3r]\), with
\(\head(\csetuthree)\) the lexicographically maximal vertex among those of
height \(\hmax\).

Let $z \in [0,u]$ be a given height, and recall that
\(\hadj(z)=z+(\hmax-u)\).  The DMU decomposition at height~$z$ uses the
lower space \(\downspace_z\) and the upper space \(\upspace_{\hadj(z)}\).
Recall that the height-$z$ downcluster \(\downcluster_z\) is the
component of \(0\) in the subgraph of \(\csetuthree\) induced
by~\(\downspace_z\); the upcluster \(\upcluster_z\) is the component of
\(\head(\csetuthree)\) in the subgraph of \(\csetuthree\) induced
by~\(\upspace_{\hadj(z)}\); and the midpiece \(\midpiece_z\) is the
graph whose edge set consists of the edges of \(\csetuthree\) that lie
in neither \(\downcluster_z\) nor \(\upcluster_z\), and whose vertex set
is formed from the endpoints of these edges.

\subsection{Backtrack characterisations}
\label{ss.dmu-backtrack}

For \(v\in V(\csetuthree)\), define its
\hfff{backtrack-protection}{\emph{down-backtrack}} by
$$
\dbt(v)
\, = \,
\min \Bigl\{ \, \hmax(P): P \textrm{ is an open path in } \csetuthree
\textrm{ from } 0 \textrm{ to } v \Big\}\, ,
$$
where \(\hmax(P)=\max\{x\cvl:x\in V(P)\}\); thus \(\dbt(v)\) is the
lowest ceiling beneath which
\(v\) connects to the base.

For an edge \(e\in E(\csetuthree)\), write \(\hmax(e)\) for the greater of
its two endpoint heights.  Recalling Section~\ref{s.basics}, the edge set of the
upspace~\(\upspace_h\) consists of the edges \(e\) for which
\(\hmax(e)>h\).  The upcluster at parameter~\(t\) is the component containing
\(\head(\csetuthree)\) in
\(\csetuthree\cap\upspace_{\hadj(t)}\), with \(\hadj(t)=\hmax-(u-t)\) and
\(\hmax=\hmax(\csetuthree)\).  Thus membership in the upper space is decided
by the greater endpoint height of an edge.  In particular, an edge for which
\(\hmax(e)\in(\hadj(t)-r,\hadj(t)]\) does not belong to the upper space; if it
does not belong to the downcluster, it is consequently assigned to the
midpiece.

\begin{proposition}[DMU decomposition at two heights]
\label{p.dmutwoheights}
Let \(0<\ell<h<u\).  Then:
\begin{enumerate}
\item for every \(t\in(0,u)\), \(V(\downcluster_t)=\{v:\dbt(v)\le t\}\),
  and \(\downcluster_t\) is the subgraph of \(\csetuthree\) induced on
  this vertex set;
\item \(\downcluster_\ell\subseteq \downcluster_h\) and
  \(\upcluster_h\subseteq \upcluster_\ell\);
\item
  $\downcluster_h = \downcluster_\ell \cup \bigl((\midpiece_\ell\cup
  \upcluster_\ell)\cap\downcluster_h\bigr)$;
\item
  $E(\midpiece_h) =
  E(\csetuthree)\setminus\bigl(E(\downcluster_h)\cup
  E(\upcluster_h)\bigr)$.
\end{enumerate}
\end{proposition}

\noindent{\bf Proof: (1).}
A vertex \(v\) belongs to \(\downcluster_t\) precisely when there is an
open path in \(\csetuthree\cap \downspace_t\) from \(0\) to \(v\).
This is equivalent to the existence of an open path \(P\) in
\(\csetuthree\) from \(0\) to \(v\) with \(\hmax(P)\leq t\), and hence
to \(\dbt(v)\leq t\); so \(V(\downcluster_t)=\{v:\dbt(v)\le t\}\).  Since
\(v\cdot\vecell\leq\dbt(v)\) for every vertex \(v\), an edge of
\(\csetuthree\) with both endpoints of down-backtrack at most~\(t\)
lies in \(\downspace_t\) and joins two vertices of
\(\downcluster_t\), hence lies in \(\downcluster_t\); so
\(\downcluster_t\) is the subgraph induced on \(\{v:\dbt(v)\le t\}\), as
stated.

{\bf (2).}
For the downcluster, (1) identifies \(\downcluster_\ell\) and
\(\downcluster_h\) as the subgraphs of \(\csetuthree\) induced on the
nested vertex sets
\(\{v:\dbt(v)\leq\ell\}\subseteq\{v:\dbt(v)\leq h\}\).
Thus \(\downcluster_\ell\subseteq\downcluster_h\).  For the upcluster,
\(\hadj(\ell)<\hadj(h)\) gives \(\upspace_{\hadj(h)}\subseteq
\upspace_{\hadj(\ell)}\), so the head-component within the
smaller space is contained in that within the larger:
\(\upcluster_h\subseteq \upcluster_\ell\).

{\bf (3).}
By (2), \(\downcluster_\ell\subseteq\downcluster_h\), so the
right-hand side of the displayed identity is a subgraph of
\(\downcluster_h\).

For the reverse inclusion, first note that
\(0\in V(\downcluster_\ell)\).  Let \(e\in E(\downcluster_h)\).
If \(e\in E(\downcluster_\ell)\), then \(e\) belongs to the right-hand
side.  Otherwise, since \(e\in E(\csetuthree)\), the definition of
\(\midpiece_\ell\) gives $e\in E(\midpiece_\ell)\cup E(\upcluster_\ell)$.
Thus,
\[
 e\in
 E\bigl((\midpiece_\ell\cup\upcluster_\ell)\cap\downcluster_h\bigr).
\]
Hence every edge of \(\downcluster_h\) belongs to the right-hand side.
Since \(\downcluster_h\) is connected and contains~\(0\), every vertex
other than~\(0\) is incident to one of these edges; the vertex~\(0\)
already belongs to \(\downcluster_\ell\).  The reverse subgraph
inclusion follows.

{\bf (4).}
This is the definition of the midpiece at height \(h\). \qed

\subsection{Protection height}
\label{ss.protection-height}

For a finite vertex set \(A\subseteq V(\csetuthree)\), define its {\em
protection height} by
$$
\prot(A)\, = \, \max\{\dbt(v):v\in A\}\, .
$$
We use the convention \(\prot(\varnothing)=-\infty\).
For a finite subgraph \(G\subseteq\csetuthree\), set
\(\prot(G)=\prot(V(G))\); for an edge \(e=\{x,y\}\), this gives
\(\prot(e)=\max\{\dbt(x),\dbt(y)\}\).  Thus
\(A\subseteq V(\downcluster_t)\) if and
only if \(\prot(A)\leq t\).  Since \(\downcluster_t\) is induced in
\(\csetuthree\), the same criterion holds for a finite subgraph
\(G\subseteq\csetuthree\).

The protection height records the threshold at which material has entered the
downcluster.  Thus, when the cut rises from $\ell$ to $h$,
Proposition~\ref{p.dmutwoheights}(3) identifies the material newly added to
$\downcluster_h$ as the part of
$\midpiece_\ell\cup\upcluster_\ell$ whose protection height is at most
$h$.  More precisely, for edges,
\[
 E(\downcluster_h)\setminus E(\downcluster_\ell)
 =\{e\in E(\midpiece_\ell)\cup E(\upcluster_\ell):
       \ell<\prot(e)\leq h\}.
\]

\begin{corollary}[Protection height and the downcluster]
\label{c.protection-height-lower-record}
Let \(S\in\clusteruthree\), and define \(\dbt\) and
\(\downcluster_t\) in the state \(S\).  For an edge
\(e=\{x,y\}\in E(\csetuthree(S))\), set \(\prot_S(e)\) equal to the
maximum of \(\dbt_S(x)\) and \(\dbt_S(y)\).  If
\(F\subseteq E(\csetuthree(S))\), write \(F_{\leq t}\) for the set of edges
\(e\in F\) such that \(\prot_S(e)\leq t\).  Then, for every \(t\in(0,u)\),
\[
   F\cap E(\downcluster_t(S))\, = \, F_{\leq t}\, .
\]
\end{corollary}

\noindent{\bf Proof.} By Proposition~\ref{p.dmutwoheights}(1), \(\downcluster_t(S)\) is the
subgraph of \(\csetuthree(S)\) induced on the vertices
\(v\) such that \(\dbt_S(v)\leq t\).  Thus an edge
\(e=\{x,y\}\in E(\csetuthree(S))\) belongs to
\(E(\downcluster_t(S))\) if and only if both \(\dbt_S(x)\leq t\) and
\(\dbt_S(y)\leq t\).  This is equivalent to \(\prot_S(e)\leq t\).
Intersecting with \(F\) proves the claim. \qed

\begin{figure}[htbp]
\centering
\includegraphics[width=0.7\textwidth]{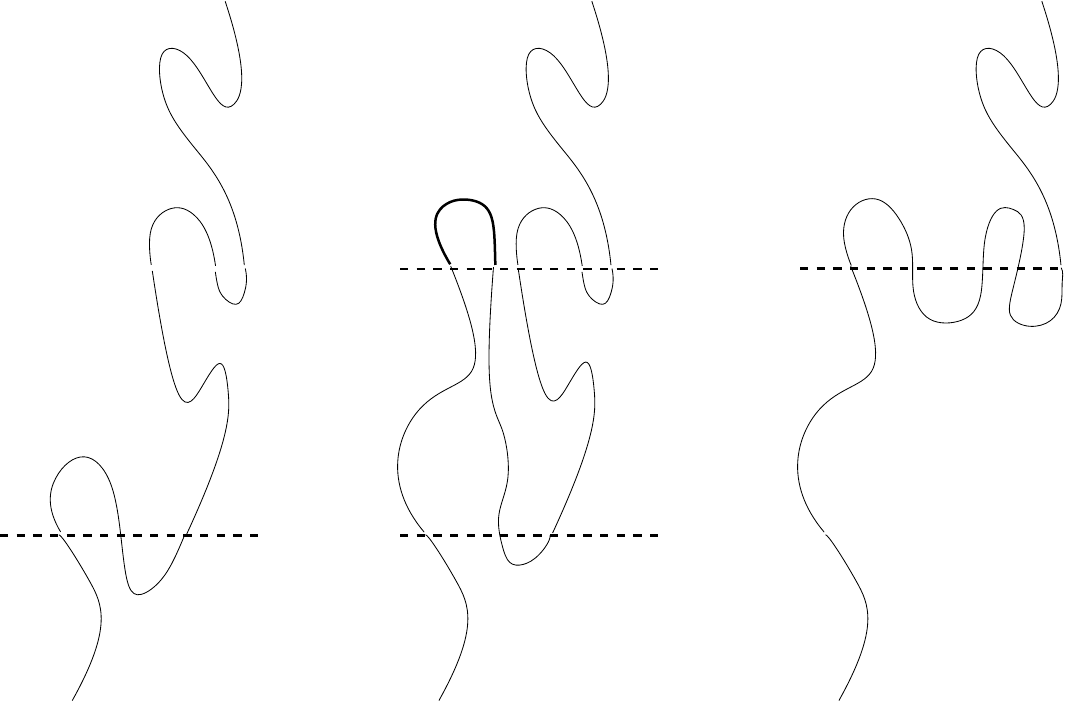}
\caption{A lower slide occurs before an upper. The midpiece introduced at the lower slide has a provisional status. It may survive the upper slide only by forming part of the downcluster at the higher slide height. Otherwise, it forms part of the midpiece for the latter slide, and exits when this slide takes place. In particular, any part of the lower slide midpiece that rises above the higher level will thus retire. The bold curve in the middle sketch illustrates this cut-it-and-lose-it mechanism.}\label{f.cutitandloseit}
\end{figure}

Figure~\ref{f.cutitandloseit} illustrates the basic
\hfff{cut-lose}{cut-it-and-lose-it mechanism}.  For
\(0<\ell<h<u\), the DMU decomposition at height~\(\ell\) informs our
understanding of the action of a slide at height~\(h\).  In the state
immediately before this slide,
\(\upcluster_h\subseteq\upcluster_\ell\).  Consider an edge of
\(\midpiece_\ell\cup\upcluster_\ell\) that does not belong to
\(\upcluster_h\).  If the edge  lies in
\(\downcluster_h\), it will automatically be retained at the slide.
To meet this condition, the edge must have an open route back to
the origin that remains in~\(\downspace_h\).
If instead every such route rises above
the height-$h$ cut, the edge belongs to \(\midpiece_h\)---since it lies in neither
 \(\downcluster_h\) nor  \(\upcluster_h\)---and
exits at the slide cut;  the edge rejoins the cluster on completion of the slide only if the shifted higher
upcluster or the new midpiece runs through it.

\section{DMU under the slide dynamics}

\subsection{Higher upclusters shift with lower slides}
\label{ss.higher-ucs-shift}

We next record the upper-side rigidity of a slide.  The first lemma is
a deterministic consequence of how the upcluster sits in the upspace.
Although it will be used to analyse slides, its proof depends only on nested
upspaces and translation; we state it separately so that it may be invoked
whenever an upper piece is translated as a whole.

\begin{lemma}[Translation of nested upclusters]
\label{l.sub-ucs-carried}
Let \(S\) and \(S'\) be two configurations in \(\clusteruthree\).
Suppose that, for some \(\ell\in(0,u)\), there is a translation \(\tau\)
such that
\[
   \upcluster_\ell(S')\, = \, \tau \cdot \upcluster_\ell(S)\, ,
\]
with \(\tau\) carrying \(\head(S)\) to \(\head(S')\).  Then, for every
\(h\in(\ell,u)\),
\[
   \upcluster_h(S')\, = \, \tau \cdot \upcluster_h(S)\, .
\]
\end{lemma}
\noindent{\bf Proof.} Write $\hadj(S,t)=t+\bigl(\hmax(S)-u\bigr)$,
and likewise for \(S'\).  For \(h>\ell\), let \(K_h(S)\) be the subgraph
of \(\upcluster_\ell(S)\) with edge set
\[
 E(K_h(S))
 =
 E(\upcluster_\ell(S))\cap E(\upspace_{\hadj(S,h)}),
\]
and with vertex set consisting of the endpoints of these edges together
with \(\head(S)\).  Since \(h>\ell\), the upper space at adjusted height
\(\hadj(S,h)\) is contained in that at adjusted height
\(\hadj(S,\ell)\).  Hence every path from the head in the former lies in
\(\upcluster_\ell(S)\), and therefore in \(K_h(S)\).  Conversely, every
path in \(K_h(S)\) lies in
\(\csetuthree(S)\cap\upspace_{\hadj(S,h)}\).  Thus
\(\upcluster_h(S)\) is precisely the component of \(\head(S)\) in
\(K_h(S)\).  The same statement holds for \(S'\).

Put \(a=\tau\cvl\).  Since
\(\tau\cdot\head(S)=\head(S')\),
$\hmax(S')=\hmax(S)+a$ and
$\hadj(S',h)=\hadj(S,h)+a$.  Moreover, for every edge \(e\),
$\hmax(\tau e)=\hmax(e)+a$.
Together with $\upcluster_\ell(S')=\tau\cdot\upcluster_\ell(S)$,
this gives $K_h(S')=\tau\cdot K_h(S)$.
Since \(\tau\) maps \(\head(S)\) to \(\head(S')\), it maps the component
of the former in \(K_h(S)\) onto the component of the latter in
\(K_h(S')\).  These components are \(\upcluster_h(S)\) and
\(\upcluster_h(S')\), respectively, which proves the claim. \qed
\begin{lemma}[Higher upclusters shift with lower slides]
\label{l.higher-ucs-shift}
Let \(S\in\clusteruthree\), and suppose that an active slide is
performed at height \(\ell\), producing the state \(S'\).  Let \(\tau\)
be the translation selected by the slide, so that
\(\upcluster_\ell(S') = \tau \cdot \upcluster_\ell(S)\).  Then, for
every \(h\in(\ell,u)\),
\[
   \upcluster_h(S')\, = \, \tau \cdot \upcluster_h(S)\, .
\]
\end{lemma}
\noindent{\bf Proof.} Translation preserves height order and the lexicographic order within each
height hyperplane, so \(\tau\) carries \(\head(S)\) to \(\head(S')\).
Lemma~\ref{l.sub-ucs-carried} therefore yields the claim. \qed

Thus a lower slide applies the same translation to the upcluster at the slide
height and to every higher upcluster contained in it.  This is particularly
simple in adjusted-height coordinates.  A translation of height component
\(a\) increases both \(x\cvl\) and \(\hmax(S)\) by \(a\), so
\[
 (\tau x)\cvl-\bigl(\hmax(S')-u\bigr)
 =x\cvl-\bigl(\hmax(S)-u\bigr).
\]
Hence a lower slide preserves the adjusted heights of the vertices of every
higher upcluster.  In particular, its vertex count in any interval specified
relative to the adjusted height is unchanged.  We record this for the
intervals \(J^\sharp_{\rho,j}\) appearing in the breaker condition.

\begin{corollary}[Upcluster transport of the breaker profile]
\label{c.breaker-points-persist}
Let \((H_i:i\in\intint L)\) and
\((S_i:i\in\llbracket0,L\rrbracket)\) be the attempted heights and successive
states of the iterated slide procedure of
Section~\ref{ss.candidate-heights}.  Fix \(q\in\intint L\) and put
\(h=H_q\).  For a state \(S\), let \(J^\sharp_{\rho,j}(S,h)\) denote the
interval in~\eqref{e.breaker-enlarged-interval}, with \(\hadj\) evaluated
in~\(S\).

For every \(i\in\llbracket0,q-1\rrbracket\), there is a translation
\(\tau_i\) such that
\[
 \upcluster_h(S_i)=\tau_i\cdot\upcluster_h(S_0),
\]
and, for every
\((\rho,j)\in\{0,1,2\}\times\intint{i_{\rm br}}\),
\[
 \bigl|V(\upcluster_h(S_i))
        \cap\{x:x\cvl\in J^\sharp_{\rho,j}(S_i,h)\}\bigr|
 =
 \bigl|V(\upcluster_h(S_0))
        \cap\{x:x\cvl\in J^\sharp_{\rho,j}(S_0,h)\}\bigr|.
\]
Consequently \({\bf B}^{+}_{\rm f}(h)\) holds in \(S_i\) if and only if it
holds in \(S_0\).
\end{corollary}

\noindent{\bf Proof.} Every active slide before time \(q\) is performed at a height below \(h\).
Iterating Lemma~\ref{l.higher-ucs-shift} over the active steps up to \(i\),
and using that dormant steps leave the state unchanged, gives a translation
\(\tau_i\) such that
\(\upcluster_h(S_i)=\tau_i\cdot\upcluster_h(S_0)\); for \(i=0\), take
\(\tau_0\) to be the identity.

Put \(a_i=\tau_i\cvl\).  Since \(\tau_i\) carries
\(\head(S_0)\) to \(\head(S_i)\),
\(\hmax(S_i)=\hmax(S_0)+a_i\), and hence
\(\hadj(S_i,h)=\hadj(S_0,h)+a_i\).
It follows directly from
\eqref{e.breaker-enlarged-interval} that
\[
 J^\sharp_{\rho,j}(S_i,h)
 =
 J^\sharp_{\rho,j}(S_0,h)+a_i .
\]
Translation by \(\tau_i\) therefore gives a bijection between the two
upcluster vertex sets counted in the displayed equality of the statement.
This proves the equality of their cardinalities and hence the final
assertion. \qed

\subsection{Downcluster increments along active slides}
\label{ss.production-rounds-protection}

Let \(N\geq1\), let \(h_1<\cdots<h_N\) be the active slide heights, and let the
\(i\)\textsuperscript{th} active slide send \(S_{i-1}\) to \(S_i\).  The
slide at \(h_i\) fixes its downcluster, so
\(\downcluster_{h_i}(S_i)=\downcluster_{h_i}(S_{i-1})\).  For \(i<N\),
define the \(i\)\textsuperscript{th} downcluster increment by
\begin{equation}
\label{e.downcluster-increment}
 \Lambda_i\, = \,
 E(\downcluster_{h_{i+1}}(S_i))
 \setminus E(\downcluster_{h_i}(S_i))\, ,
\end{equation}
and put \(\Lambda_N=E(\midpiece_{h_N}(S_N))\).  Thus
\(\Lambda_1,\ldots,\Lambda_{N-1}\) describe successive changes in
the downcluster, while \(\Lambda_N\) denotes the midpiece remaining after
the last slide.

The exits and entrances at successive slides assemble into an explicit record
of how the slide dynamics transform the cluster from input to output.
Although we make no later appeal to the resulting description, we end the
chapter with this application of the DMU decomposition.

\begin{proposition}[Downcluster increments and the output cluster]
\label{p.production-rounds-protection}
For every \(i<N\),
\begin{align}
 \Lambda_i
 &\, = \, \bigl\{e\in E(\csetuthree(S_i))\, :
       h_i<\prot_{S_i}(e)\, \leq \, h_{i+1}\bigr\} \label{e.increment-protection}\\
 &\, = \, E\Bigl(
   (\midpiece_{h_i}(S_i)\cup\upcluster_{h_i}(S_i))
   \cap\downcluster_{h_{i+1}}(S_i)
   \Bigr). \label{e.increment-dmu}
\end{align}
Moreover,
\begin{equation}
\label{e.final-downcluster-decomposition}
 E(\downcluster_{h_N}(S_N))
\, = \, E(\downcluster_{h_1}(S_0))
  \sqcup\Lambda_1\sqcup\cdots\sqcup\Lambda_{N-1}\, ,
\end{equation}
and the output edge set has the disjoint decomposition
\begin{equation}
\label{e.output-dmu-decomposition}
 E(\csetuthree(S_N))
\, = \, E(\downcluster_{h_1}(S_0))
  \sqcup\Lambda_1\sqcup\cdots\sqcup\Lambda_N
  \sqcup E(\upcluster_{h_N}(S_N))\, .
\end{equation}
\end{proposition}

\noindent\textbf{Proof.} Apply Corollary~\ref{c.protection-height-lower-record} at the two heights
\(h_i<h_{i+1}\) in state \(S_i\).  The difference in
\eqref{e.downcluster-increment} is exactly the protection-height set in
\eqref{e.increment-protection}.  Proposition~\ref{p.dmutwoheights}(3),
again in state \(S_i\), gives \eqref{e.increment-dmu}.

The slide at \(h_{i+1}\) fixes \(\downcluster_{h_{i+1}}\), and therefore
\[
 E(\downcluster_{h_{i+1}}(S_{i+1}))
\, = \, E(\downcluster_{h_{i+1}}(S_i))
\, = \, E(\downcluster_{h_i}(S_i))\sqcup\Lambda_i\, .
\]
The first active slide also fixes
\(\downcluster_{h_1}\), so
\(E(\downcluster_{h_1}(S_1))=E(\downcluster_{h_1}(S_0))\).
Starting from this equality and iterating the displayed identity proves
\eqref{e.final-downcluster-decomposition}; at each stage the union is disjoint
because the two downclusters in \eqref{e.downcluster-increment} are nested.
Finally, Proposition~\ref{p.dmutwoheights}(4) at height \(h_N\) in state
\(S_N\) gives
\[
 E(\csetuthree(S_N))
\, = \, E(\downcluster_{h_N}(S_N))
  \sqcup E(\midpiece_{h_N}(S_N))
  \sqcup E(\upcluster_{h_N}(S_N))\, .
\]
Substitution of \eqref{e.final-downcluster-decomposition} and the definition
of \(\Lambda_N\) proves \eqref{e.output-dmu-decomposition}. \qed

\chapter{Quilts and uniform open escape in a layer}
\label{s.overlapping-quilts-seals}

This chapter and the next furnish seals and provide the shielding estimates used
to prove the IVF-to-LLF transition.  To keep track of cluster accrual during
this transition, we work inside a \emph{layer}: a fixed-width slab graph
transverse to~\(\vecell\).  In the later application, such a layer contains the unit allowance band whose accrual is being monitored, in an alignment that allows the increasing slide schedule to pass over it without cutting through its interior.
The monitored layers form an edge partition of
\(E(\mathbb Z^d)\); consecutive layers share their common boundary vertices but no edges.

These shielding estimates play different roles for slides originating near
and far from a given layer.  A nearby slide may reach the layer, but the volume
that it adds there is controlled.  An exploration arising from a slide far
below must instead traverse a long line of wavebreakers.  At each wavebreaker
the shielding estimate gives the exploration a uniformly positive chance of
failing to cross intact; the exploration therefore typically never reaches
the layer to contribute to accrual.  In both cases the underlying question is
the same: within a given layer,
to what extent can a prescribed finite collection of discovered open and closed edges shield subsequent growth from open escape?

The answer is obtained by encircling the discovered structure with a fabric of open paths in the layer.  Such a fabric is built from overlapping coarse-grained boxes that satisfy an Antal--Pisztora-type {\em aptness} condition; a connected collection of apt boxes is called a \emph{quilt}.  A seal is a suitably selected surrounding quilt collection that cuts off a finite part of the layer containing the discovered structure while leaving the exterior uninspected.  Conditionally on the information thereby captured, the exterior edges retain their Bernoulli product law.
The principal conclusion of the present chapter,
 Proposition~\ref{p.uniform-open-escape-exterior-seed}, shows that
  an exploration in the layer that grows from a given seed location bordering this product exterior consequently runs a uniformly positive risk of open escape. In Chapter~\ref{s.seed-exploration-seal-extension}, we will find that,  on the complementary finite-cluster event, the finite open volume gained in the layer by the exploration has a stretched-exponential tail.

The Antal--Pisztora coarse-graining of Section~\ref{s.coarse} was used to control the enclosing surface of the single-slide cluster.  Here it must do more.  The sealing construction must capture every discovered open and closed edge without inspecting the exterior; and the resulting conditioning must not make the boxes immediately outside the boundary atypically hostile to open passage.  The seal will therefore be both a stopping domain and, in the terminology introduced below, \emph{welcoming}: after the captured information has been exposed, each nearby exterior box still has a high conditional probability of being apt.

In the first of four sections, we set up the local coarse geometry and the
apt-box fabric from which quilts are formed.  We then introduce the disjoint
mesh, whose edge partition makes capture and exterior exact edgewise notions,
together with the boundary reduction needed to remove buried pieces.  A
welcoming condition controls boxes that straddle the resulting boundary after
the captured information has been exposed.  Finally, boundary connectivity and
high-density exterior percolation produce an infinite apt route from the
boundary fabric, to which a bounded finite-energy opening attaches a prescribed
exterior seed, yielding the uniform open-escape estimate of
Proposition~\ref{p.uniform-open-escape-exterior-seed}.

\section{Local coarse geometry and apt boxes}
\label{s.overlapping-boxes-boxlets}

We partition space into layers bounded by hyperplanes orthogonal to~$\vecell$.
This is the natural geometry for the slide dynamics, whose steps offer
upcluster shifts that preserve the height coordinate~$x\cvl$ up to a bounded
offset. Within each layer we coarse-grain in the ordinary coordinate
directions.
 The resulting axial columns cut the tilted
slab into bounded connected boxlets, and neighbouring groups of boxlets form the
overlapping boxes from which quilts will be built.

\subsection{Boxlets and overlapping boxes}

Recall the axes-ordering convention~\eqref{e.aoc}, and write
\(\ell_i=e_i\cvl\), \(i\in\intint d\).  Thus
\(0\le\ell_i\le\ell_d\) and \(\ell_d\ge d^{-1/2}\).  In the box geometry below, we call
\(x_1,\ldots,x_{d-1}\) the \emph{transverse coordinates}.  Every lattice edge has
a unique representation \(e=\{x,x+e_i\}\).  We call
\(e^\uparrow=x+e_i\) its \emph{upper endpoint}; thus
\(e^\uparrow\cvl=h^+(e)\).  When \(\ell_i=0\), the two endpoints have the
same height, and this convention simply selects one of them.

For \(k\in\Z\), put
\[
   a_k=K\bigl(k-\tfrac32\bigr),\qquad
   b_k=K\bigl(k+\tfrac12\bigr),\qquad
   \mathsf A_k=\{x\in\Z^d:a_k<x\cvl\le b_k\},
\]
and let the \hfff{layer}{\emph{layer}} of index \(k\) be
\(N_k=\slab_{a_k,b_k}\).  Since membership of an edge in a slab is decided by its greater endpoint
height,
\begin{equation}
\label{e.layer-upper-endpoint}
       E(N_k)=\{e\in E(\Z^d):e^\uparrow\in\mathsf A_k\}.
\end{equation}
The two sides of \(N_k\) are its \emph{southside}
\(\partial_-N_k=\vertbdry_{a_k}\) and its \emph{northside}
\(\partial_+N_k=\vertbdry_{b_k}\).
The even-indexed layers $\big\{ N_{2n}: n \in \mathbb{Z} \big\}$ form the partition of edge space that we will use.

We retain the usual Antal--Pisztora terminology \emph{box}: the present boxes
are the layer-adapted counterparts of the overlapping axial boxes of the
standard coarse graining, and need not be axis-aligned boxes in physical
space.  We first define \emph{boxlets}, their smaller axial-column constituents.

\begin{figure}[htbp]
\centering
\includegraphics[width=0.8\textwidth]{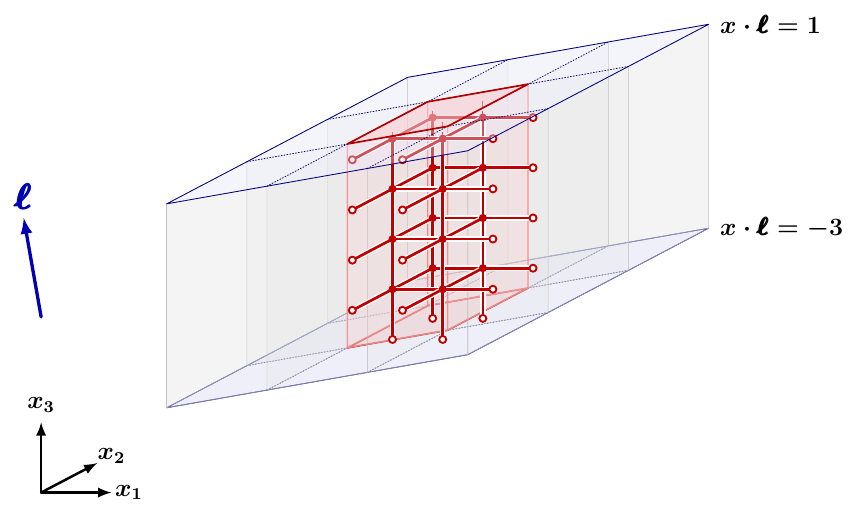}
\caption{A boxlet is drawn in red in $d=3$ with $K=2$, for a choice of $\vecell$ lying in the $(x_1,x_3)$-plane at a small angle from the positive vertical direction. The layer has width $2K=4$ in the $\vecell$-direction, and is bounded by parallel tilted planes on which lines mark the intersection with the axial vertical hyperplanes that specify the boxlet edge-partition.  (Figure prepared with Claude.)}\label{f.boxlet}
\end{figure}

Write a transverse index as
\(u=(u_1,\ldots,u_{d-1})\in\Z^{d-1}\), and set
\[
 \mathsf C_k(u)=
 \Bigl\{x\in\mathsf A_k:
       K(u_i-1)<x_i\le Ku_i,\quad i\in\intint{d-1}\Bigr\}.
\]
The piece of the layer cut out by this transverse \(K\)-column is
the \emph{boxlet} \(\Bbox{k,u}\), the subgraph with
\begin{equation}
\label{e.cell-edge-partition}
 E(\Bbox{k,u})
   =\{e\in E(N_k):e^\uparrow\in\mathsf C_k(u)\},
\end{equation}
and with vertex set given by the endpoints of these edges.  In other words, each layer edge is assigned to the unique transverse
\(K\)-column containing its upper endpoint.  The boxlet edge sets therefore
partition \(E(N_k)\).  A lower endpoint may lie outside \(\mathsf C_k(u)\), but every boxlet edge has its upper endpoint in \(\mathsf C_k(u)\).
See Figure~\ref{f.boxlet} for a depiction.

The elementary geometry behind this construction is as follows.  For
\(x'=(x_1,\ldots,x_{d-1})\in\Z^{d-1}\), let
\[
 J_k(x')
 =\left(
   \frac{a_k-\sum_{i<d}\ell_i x_i}{\ell_d},
   \frac{b_k-\sum_{i<d}\ell_i x_i}{\ell_d}
   \right].
\]
The vertices of \(\mathsf A_k\) above \(x'\) are exactly the points
\((x',m)\) with \(m\in J_k(x')\cap\Z\); we call this set the
\(e_d\)-\emph{fibre} above~\(x'\).  The interval \(J_k(x')\) has length
\(2K/\ell_d\ge2K\), while
\(J_k(x'+e_i)=J_k(x')-\ell_i/\ell_d\), with
\(0\le\ell_i/\ell_d\le1\).  Hence the \(e_d\)-fibres above adjacent
transverse points have at least \(2K-1\) consecutive integer levels in common.

\begin{lemma}[Connected boxlets]
\label{l.connected-layer-cells}
Let \(A\subseteq\Z^{d-1}\) be finite and nearest-neighbour connected.
The union \(\bigcup_{u\in A}\mathsf C_k(u)\) is nearest-neighbour connected
in \(\Z^d\), and the union of the boxlet graphs \(\Bbox{k,u}\), \(u\in A\),
is connected.  In particular every
boxlet is connected.  Moreover there is \(C_d<\infty\) such that
\[
       |E(\Bbox{k,u})|\le C_dK^d,
\]
and, from every vertex incident to an edge of \(\Bbox{k,u}\), that boxlet
contains a path of Euclidean diameter at least \(K-1\).
\end{lemma}

\noindent{\bf Proof.} For fixed \(x'\), the fibre \(J_k(x')\cap\Z\) is a non-empty interval of at
least \(2K\) consecutive integers.  If \(x'\) and \(x'+e_i\) are neighbours,
their fibres have a common integer level, so a horizontal edge joins them.
Inside one transverse \(K\)-column, and then across neighbouring columns,
these horizontal joins together with the vertical \(e_d\)-fibres prove the
first assertion.  Every nearest-neighbour edge between two vertices of \(\mathsf C_k(u)\) belongs to the boxlet edge set.  Every remaining boxlet edge has its upper endpoint in \(\mathsf C_k(u)\), so adjoining these edges preserves connectedness.  This proves connectedness of each boxlet graph and of the stated connected unions.

There are \(K^{d-1}\) possible transverse integer coordinates in one boxlet
and at most \(2K/\ell_d+1\le 2d^{1/2}K+1\) vertices of \(\mathsf C_k(u)\)
above each.  At most \(d\) boxlet edges have any one of these vertices as upper
endpoint, which gives the cardinality bound.  Finally, after traversing one
edge if necessary, start from a vertex of \(\mathsf C_k(u)\).  Its \(e_d\)-fibre
contains at least \(2K\) consecutive vertices, so one of the two vertical
directions supplies a path in the boxlet of at least \(K-1\) edges. \qed

For \(v\in\Z^{d-1}\), the \hfff{overlapping-layer-box}{\emph{overlapping box}
\(\Obox{k,v}\)} is the union of the \(2^{d-1}\) boxlets
\(\Bbox{k,v+\sigma}\), \(\sigma\in\{-1,0\}^{d-1}\).  Its upper endpoints are exactly the lattice points in the half-open
parallelotope
\begin{equation}
\label{e.layer-box-region}
 \Lambda_{k,v}
 =\Bigl\{x\in\R^d:
       a_k<x\cvl\le b_k,\quad
       K(v_i-2)<x_i\le Kv_i,\ i\in\intint{d-1}\Bigr\}.
\end{equation}
Thus the box has width \(2K\) in each transverse coordinate, while increasing
one coordinate of \(v\) by one shifts the box by \(K\) in that direction.
Neighbouring boxes therefore overlap by one boxlet-width in each transverse coordinate in which their indices differ.
  As a graph,
\(\Obox{k,v}\) may also contain the lower endpoints of boundary edges whose
upper endpoints lie in \(\Lambda_{k,v}\).

Let \(\OBoxLat=\{\Obox{k,v}:(k,v)\in\Z\times\Z^{d-1}\}\).
Two boxes whose full indices have \(\ell^\infty\)-distance greater than one
have disjoint edge sets.  In a fixed layer, two boxes whose transverse
indices have \(\ell^\infty\)-distance at most one share at least one whole
boxlet.

If \(\mathscr A\subseteq\OBoxLat\), set
\(\Edg(\mathscr A)=\bigcup_{B\in\mathscr A}\Edg(B)\) and
\(\Vtx(\mathscr A)=\bigcup_{B\in\mathscr A}\Vtx(B)\), and write
\(\BoxUnion{\mathscr A}=(\Vtx(\mathscr A),\Edg(\mathscr A))\).
For a vertex set \(A_V\), put
\[
 \Boxsupp(A_V)=\{B\in\OBoxLat:\Vtx(B)\cap A_V\ne\emptyset\}.
\]
For a subgraph \(A\), use the same notation with
\[
       \Boxsupp(A)=\Boxsupp(V(A)).
\]
Thus this coarse support is intentionally vertex based.  Exact edgewise
support will be recorded below using the disjoint mesh.  For a finite subgraph
\(A\), write \(|A|_*=|\Boxsupp(A)|\).  Since \(K\) is fixed in the
applications, microscopic cardinality is bounded by a fixed multiple of
coarse cardinality on unions of coarse boxes.

\subsection{The aptness condition}
\label{s.apt-boxes}

An apt box carries a macroscopic open piece through each of its boxlets and absorbs every other macroscopic open path that it contains.  The first property supplies macroscopic open pieces in shared boxlets; the second forces these pieces to join and excludes a competing macroscopic cluster.

\begin{definition}[Apt boxes]
\label{d.apt-box}
An overlapping box \(B\in\OBoxLat\) is
\hfff{apt-woven}{\emph{apt}}\footnote{The adjective is offered as an
Antal--Pisztora tribute.} if there is an open cluster
\(\mathcal W(B)\) in \(B\) such that
\begin{enumerate}
\item for every boxlet \(L\) of \(B\),
      \(\mathcal W(B)\cap L\) contains an open path of Euclidean diameter at
      least \(3K/4\);
\item every open path in \(B\) of Euclidean diameter at least \(K/2\) is
      contained in \(\mathcal W(B)\).
\end{enumerate}
The cluster \(\mathcal W(B)\) is called the \emph{woven cluster} of \(B\).
\end{definition}
The absorption property of aptness makes the woven cluster unique whenever it exists.

We next state the layer-box Antal--Pisztora estimate derived in
Appendix~\ref{s.nonaxial-ap-proof}.  In order to state it in the general form
that we will need, we first extend the layer-box notation to arbitrary height
offsets.

For \(k\in\Z\), write
\[
   \overline N_k=\{\Obox{k,v}:v\in\Z^{d-1}\}.
\]
By \eqref{e.cell-edge-partition} and the definition of overlapping boxes,
\(\BoxUnion{\overline N_k}=N_k\).

For wavebreakers, we will also have cause to consider layers at height
offsets which do not belong to the fixed even-indexed partition.  For
\(a\in\R\), put
\[
 \mathsf A(a)=\{x\in\Z^d:a<x\cvl\le a+2K\},
 \qquad N(a)=\slab_{a,a+2K}.
\]
Use the same transverse \(K\)-columns, replacing \(\mathsf A_k\) by
\(\mathsf A(a)\) in \eqref{e.cell-edge-partition}; define boxlets and
overlapping boxes by the same formulas, and write \(\overline N(a)\) for the
resulting overlapping-box lattice.  When a transverse mesh is used on
\(N(a)\), take the boxes with even transverse index~\(2z\), as for a
partition layer.  A lattice translation~\(t\) sends the physical slab
\(N(a)\) to \(N(a+t\cvl)\); this does not assert invariance of the chosen
transverse mesh.

A field \((X_v)_{v\in\Z^m}\) is \(r\)-dependent if the families
\((X_v)_{v\in I}\) and \((X_v)_{v\in J}\) are independent whenever
\(\min\{\|v-w\|_\infty:v\in I,w\in J\}>r\).

\begin{proposition}[Layer-box Antal--Pisztora input]
\label{p.nonaxial-ap}
There are constants \(C,c>0\) and \(K_0\in\N\) such that, for \(K\ge K_0\),
with \(\varepsilon_K=Ce^{-cK}\), every \(B\in\OBoxLat\) satisfies
\[
       \PP(B\text{ is not apt})\le\varepsilon_K.
\]
Aptness of \(B\) is measurable with respect to the edge configuration on
\(E(B)\), and the inapt-box field on the full index lattice
\(\Z\times\Z^{d-1}\) is \(1\)-dependent.  For every \(a\in\R\), the same
estimate and measurability hold for the boxes of \(\overline N(a)\), uniformly
in~\(a\), and the corresponding inapt-box field on \(\Z^{d-1}\) is
\(1\)-dependent.  The constants \(C,c,K_0\) may be chosen uniformly in the
height direction~\(\vecell\in S^{d-1}\).
\end{proposition}

\subsection{Quilts}

For boxes $B,B'$ in a common layer-box lattice, write $B\sim B'$ for
nearest-neighbour adjacency of their transverse indices and $B\sim_*B'$ when
the distinct transverse indices have $\ell^\infty$-distance one.

\begin{definition}[Quilts and quilt collections]
A \hfff{quilt}{quilt} is a finite nearest-neighbour connected collection of
apt boxes in one layer-box lattice.  A \emph{quilt collection} is a finite
collection of apt boxes; its nearest-neighbour connected components, each a
quilt, are called its \emph{quilt components}.
\end{definition}

\begin{lemma}[Neighbouring apt boxes and quilt clusters]
\label{l.woven-fabric-note}
\leavevmode
\begin{enumerate}
\item If $B,B'$ are distinct apt boxes in a common layer-box lattice with
      $B\sim_*B'$, then
      $\mathcal W(B)\cap\mathcal W(B')\ne\emptyset$.
\item If $\mathscr Q$ is a quilt, its \emph{woven fabric}
      $\mathcal W(\mathscr Q)=\bigcup_{B\in\mathscr Q}\mathcal W(B)$
      is a finite connected open subgraph of the layer.  If an open path has
      a subpath of diameter at least $K/2$ contained in a box of
      $\mathscr Q$, then this subpath belongs to $\mathcal W(\mathscr Q)$,
      and the full path belongs to the open cluster containing this fabric.
\end{enumerate}
\end{lemma}

\noindent{\bf Proof: (1).}
The boxes $B$ and $B'$ have a common boxlet $L$.  By aptness,
$\mathcal W(B)\cap L$ contains an open path of diameter at least $K/2$.
This path lies in $B'$, so the absorption property of aptness for $B'$ places it in
$\mathcal W(B')$.

{\bf (2).}
Chain part~(1) along a nearest-neighbour path in the quilt to obtain the
connectedness of the union.  The absorption property puts the indicated
subpath in the woven cluster of its box; the full open path is consequently
connected to the woven fabric. \qed

\begin{figure}[!htbp]
\centering
\includegraphics[width=0.75\textwidth]{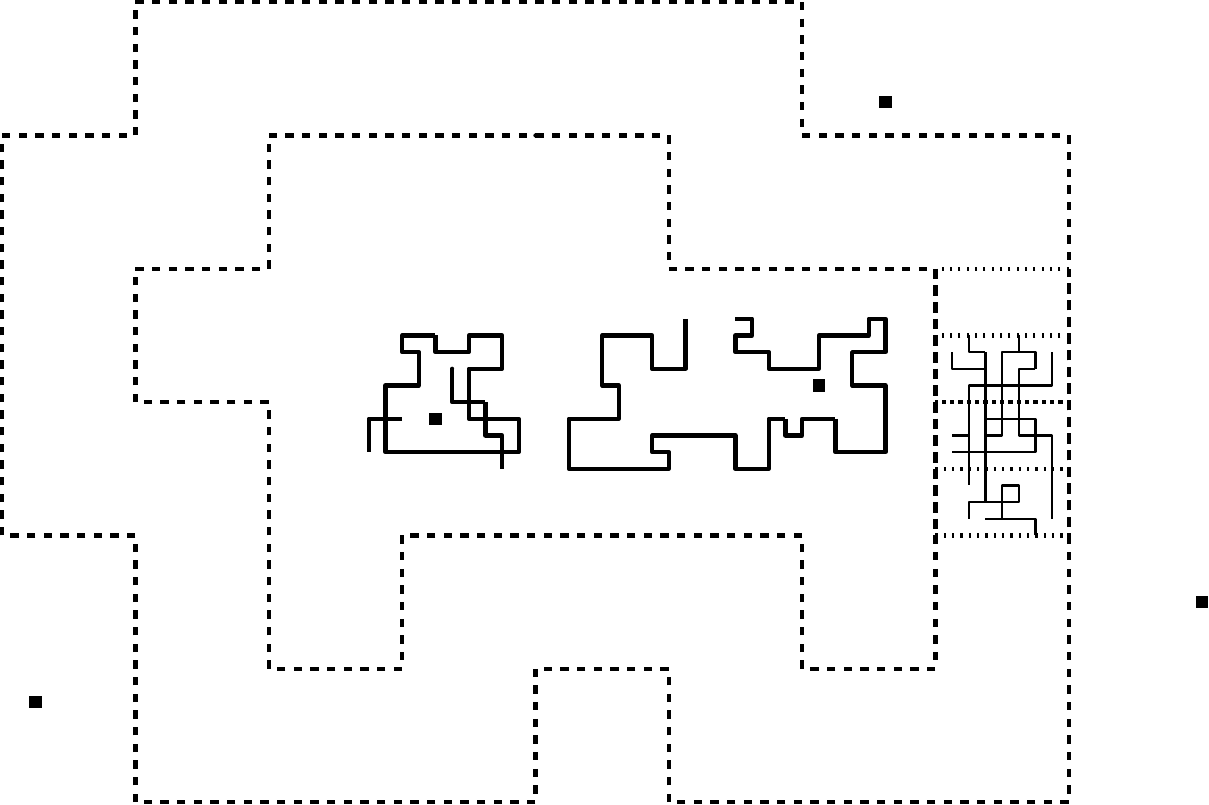}
\caption{A schematic of the geometry in a layer.  The dashed contours mark an
enclosing quilt around discovered material.  At right, the woven clusters in
two adjacent apt boxes join across their overlap.  The square markers indicate
seeds from which later open-cluster explorations in the layer may start.  The
discovered material shields growth from the two interior seeds, completely for
the left seed and partially for the right one.  Growth from any of the three
exterior seeds cannot be shielded in this way, since any open approach to the
discovered structure would first connect to the quilt that captures it.}
\label{f.quilt}
\end{figure}

Figure~\ref{f.quilt} gives a schematic indication of the role quilts will play as enclosing fabric in a layer.

In each application below, \(K\) is chosen once subject to the stated lower
bounds and is then unchanged throughout
Chapters~\ref{s.overlapping-quilts-seals} and
\ref{s.seed-exploration-seal-extension}.  The notation \(\varepsilon_K\)
refers to this fixed one-box error bound.  The single-layer constructions and
estimates below are stated for members of the fixed partition.  When one is
invoked on \(N(a)\), it is understood with the offset box structure just
defined; the same constants apply uniformly in~\(a\).

\section{Mesh boundaries and capture in a fixed layer}
\label{s.slab-setting}

We first record the partition property of the even-indexed layers and then fix
one such layer.  The boxes with even transverse indices form a disjoint mesh
whose edge sets partition the layer.  Finite mesh regions will specify the
edge information lying on either side of a boundary.  Along such a mesh
region we include every overlapping box that is contained in it, so that the
box geometry used for aptness remains available.

\subsection{The mesh and overlapping boxes}

\begin{lemma}[The layer partition]
\label{l.layer-partition}
The edge sets \(E(N_{2n})\), \(n\in\Z\), are pairwise disjoint and have union
\(E(\Z^d)\).  Consecutive partition layers have the common side
\[
       \partial_+N_{2n}
       =\partial_-N_{2n+2}
       =\vertbdry_{K(2n+1/2)}.
\]
\end{lemma}

\noindent{\bf Proof.} By~\eqref{e.layer-upper-endpoint}, \(E(N_{2n})\) consists of the edges whose
upper-endpoint heights lie in
\[
       \bigl(K(2n-\tfrac32),K(2n+\tfrac12)\bigr].
\]
These half-open intervals partition \(\R\).  The side identity follows from
\(b_{2n}=a_{2n+2}\). \qed

Fix henceforth a partition layer \(N=N_{2n}\), and write
\(\overline N=\overline N_{2n}\).  The transverse index identifies the box
graph of \(\overline N\) with \(\Z^{d-1}\); nearest-neighbour and
\( * \)-adjacency below refer to this lattice.  Define the \emph{mesh} of
\(N\) to be the disjoint subcollection
\[
   \Mesh=\Mesh_{2n}
   =\{\MeshBox z:=\Obox{2n,2z}:z\in\Z^{d-1}\}.
\]
In transverse coordinate \(i\), the mesh box \(\MeshBox z\) consists of the
two boxlet labels \(2z_i-1\) and \(2z_i\).  These pairs partition \(\Z\).
Thus a mesh box is one of the overlapping boxes, distinguished by its even
transverse index.

\begin{lemma}
\label{l.disjoint-twoK-mesh}
The edge sets \(E(\MeshBox z)\), \(z\in\Z^{d-1}\), are pairwise disjoint
and have union \(E(N)\).
\end{lemma}

\noindent{\bf Proof.}
The boxlet edge sets partition \(E(N)\), and the displayed pairs partition
their transverse indices. \qed

For an edge set \(H\subseteq E(N)\), define its \emph{mesh support} by
\[
  \MeshSupp(H)=\{M\in\Mesh:E(M)\cap H\ne\emptyset\}.
\]
For a subgraph or an overlapping box we use the same notation for the support
of its edge set; in particular,
\[
  \MeshSupp(B)=\{M\in\Mesh:E(M)\cap E(B)\ne\emptyset\},
  \qquad B=\Obox{2n,v}\in\overline N.
\]
Thus the mesh support records exactly which members of the disjoint mesh
partition contribute edges to the object in question.  Equivalently for an
overlapping box, these are the mesh boxes containing its boxlets.
In each transverse coordinate this support has one mesh index when the
corresponding coordinate of \(v\) is even and two adjacent indices when it
is odd.  Hence
\begin{equation}
\label{e.mesh-support-bound}
       |\MeshSupp(B)|\le2^{d-1}.
\end{equation}
For nearest-neighbour \(B,B'\in\overline N\), their mesh supports are nested:
in the changed coordinate a singleton is replaced by the adjacent pair, or
conversely, and every other coordinate is unchanged.

For \(A\subseteq\Mesh\), put
\[
       A^{\mathrm{ov}}
       =\{B\in\overline N:\MeshSupp(B)\subseteq A\}.
\]
Thus \(A^{\mathrm{ov}}\) is the collection of overlapping boxes contained,
edgewise, in the mesh region \(A\).  It contains every mesh box in \(A\),
and it contains the intervening overlapping box between any two neighbouring
mesh boxes of \(A\).  This enrichment does not change the underlying boxlet
region: \(\BoxUnion{A^{\mathrm{ov}}}=\BoxUnion A\).

\begin{lemma}
\label{l.overlap-mesh-region}
Let \(A\subseteq\Mesh\) be finite and nearest-neighbour connected.  Then
\(A^{\mathrm{ov}}\) is nearest-neighbour connected and
\(|A^{\mathrm{ov}}|\le 3^{d-1}|A|\).
\end{lemma}

\noindent{\bf Proof.} If
\(\MeshBox z\) and \(\MeshBox{z+e_i}\) are neighbours, the overlapping box
with transverse index \(2z+e_i\) has these two mesh boxes as its support and
bridges them.  These bridges connect the mesh boxes along a
path in \(A\).  Every remaining member of \(A^{\mathrm{ov}}\) is joined to a
supporting mesh box by successively resolving its odd coordinates.  A fixed
mesh box belongs to the support of at most \(3^{d-1}\) overlapping boxes,
which gives the cardinality estimate. \qed

\subsection{Capture and boundary reduction}

\hfff{exterior-capture}{For a finite $\mathscr A\subseteq\overline N$, let}
$\Ext_{\overline N}(\mathscr A)$ be the unbounded nearest-neighbour component
of $\overline N\setminus\mathscr A$ and put
\[
 \Capt_{\overline N}(\mathscr A)
   =\overline N\setminus\Ext_{\overline N}(\mathscr A),
 \qquad
 \partial_{\rm ext}^{\overline N}\mathscr A
   =\{B\in\Ext_{\overline N}(\mathscr A):B\sim_*A
        \text{ for some }A\in\mathscr A\}.
\]
We call $\Capt_{\overline N}(\mathscr A)$ the \emph{capture} of~$\mathscr A$
and $\Ext_{\overline N}(\mathscr A)$ its \emph{exterior}.  Thus the capture
consists of $\mathscr A$ together with the bounded components of its
complement, while the exterior is the component connected to infinity.  The
same notation and language are used on the mesh lattice.  A finite mesh set
$D$ is \emph{filled} if $\Capt_{\Mesh}(D)=D$, equivalently if
$\Mesh\setminus D$ has no bounded nearest-neighbour component; we call
$\Capt_{\Mesh}(D)$ the \emph{filling} of~$D$.

For a finite $S\subseteq\overline N$, define
\begin{equation}
\label{e.mesh-edge-partition}
 \MCap(S)=\Mesh\cap\Capt_{\overline N}(S),\qquad
 \MExt(S)=\Mesh\cap\Ext_{\overline N}(S),
\end{equation}
and put
\[
 \MeshInt(S)=\MCap(S)\setminus(S\cap\Mesh).
\]
Thus $\MeshInt(S)$, the \emph{mesh interior} of~$S$, consists of the captured mesh
boxes lying strictly behind~$S$.  Let
\[
 \ECap(S)=\bigcup_{M\in\MCap(S)}E(M),\qquad
 \mathcal U(S)=\bigcup_{M\in\MExt(S)}E(M).
\]
We call $\MCap(S)$ the captured mesh set and $\MExt(S)$ the exterior mesh
set, and $\ECap(S)$ and $\mathcal U(S)$ the corresponding captured and
exterior edge sets.  Since the mesh edge sets partition $E(N)$, these two edge
sets are disjoint and exhaust the layer.  Thus every layer edge has an
unambiguous captured or exterior status even though the boxes in $S$ may
belong to the full overlapping collection.

For a finite mesh set $A$, let $U^*(A)$ denote the unbounded
$*$-component of $\Mesh\setminus A$.  The next lemma relates the exterior of
$A^{\mathrm{ov}}$ in the overlapping-box lattice to this mesh exterior.

\begin{lemma}[Mesh transfer]
\label{l.mesh-transfer}
Let $A\subseteq\Mesh$ be finite and put $S=A^{\mathrm{ov}}$.  Then
\begin{equation}
\label{e.mesh-transfer-capture}
 \MExt(S)=U^*(A),\qquad
 \MeshInt(S)=(\Mesh\setminus A)\setminus U^*(A),\qquad
 \MCap(S)=A\,\sqcup\,\MeshInt(S).
\end{equation}
Moreover,
\begin{equation}
\label{e.mesh-transfer-full-exterior}
 \Ext_{\overline N}(S)
 =\{B\in\overline N:\MeshSupp(B)\cap U^*(A)\ne\varnothing\}.
\end{equation}
\end{lemma}

\noindent{\bf Proof.}
Let $M,M'\in\Mesh\setminus A$ be $*$-adjacent.  There is an overlapping
box whose mesh support contains both.  Resolving its odd transverse
coordinates first towards $M$ and then towards $M'$ gives a
nearest-neighbour box path from $M$ to $M'$ outside $S$: every box on the
path has one of $M,M'$ in its support and hence cannot belong to $A^{\mathrm{ov}}$.
Thus an infinite $*$-path in $\Mesh\setminus A$ lifts to a path to infinity
in $\overline N\setminus S$.

Conversely, let $B_0,B_1,\ldots$ be a nearest-neighbour path to infinity in
$\overline N\setminus S$.  For each $j$ choose
$M_j\in\MeshSupp(B_j)\setminus A$, which is possible because
$B_j\notin A^{\mathrm{ov}}$.  The supports of consecutive boxes are nested, so
$M_j$ and $M_{j+1}$ are equal or $*$-adjacent.  Since mesh supports have
uniformly bounded diameter, the resulting mesh path leaves every finite set.
If $B_0=M$ is a mesh box, this proves that $M$ is exterior for $S$ exactly
when $M\in U^*(A)$, and hence, by the definitions of $\MeshInt(S)$ and
$\MCap(S)$, gives~\eqref{e.mesh-transfer-capture}.

For a general $B\notin S$, a support member $M\in U^*(A)$ joins $B$ to
infinity outside $S$ by first resolving the odd coordinates towards $M$ and
then using the preceding construction.  Conversely, applying the preceding
selection to a path from $B$ to infinity produces a member of
$\MeshSupp(B)\cap U^*(A)$.  This proves
\eqref{e.mesh-transfer-full-exterior}. \qed

The remaining lattice-topological input is exterior-boundary connectivity,
already used in the single-slide analysis.

\begin{lemma}[Exterior boundary connectivity]
\label{l.extbdy-slab}
Let $d\ge3$, and let $\mathcal L$ denote either $\Mesh$ or $\overline N$,
identified through its transverse index with $\Z^{d-1}$.  If
$C\subseteq\mathcal L$ is finite and $*$-connected, then
$\partial_{\rm ext}^{\mathcal L}C$ is nearest-neighbour connected.
\end{lemma}

\noindent{\bf Proof.} Under the transverse identification, this is
Lemma~\ref{l.extbdy-halfspace}, with $\Z^d$ replaced by $\Z^{d-1}$. \qed

With graph-theoretic arguments, Tim\'ar~\cite{Timar2013} has obtained
boundary-connectivity results generalising those of Kesten and
Deuschel--Pisztora~\cite{DeuschelPisztora1996}.

For a finite mesh set $D$, write
\[
       A(D)=\partial_{\rm ext}^{\Mesh}D.
\]
Thus $A(D)$ is the exterior-facing mesh boundary of $D$: its members lie in
the unbounded nearest-neighbour component of $\Mesh\setminus D$ and are
$*$-adjacent to $D$.  The associated overlapping-box boundary is
$A(D)^{\mathrm{ov}}$.

Passing from $D$ to its exterior mesh boundary $A(D)$ can seal a narrow
entrance that was open in $\Mesh\setminus D$.  Consequently
$A(D)^{\mathrm{ov}}$ may capture mesh boxes outside $D$, and may even bury
components of $A(D)$ itself.  The lemma below first records the resulting
enclosure of $D$, and then introduces the reduced boundary $A^\circ$, which
discards these buried components without changing the capture or exterior.

\begin{lemma}[Boundary enclosure and reduction]
\label{l.mesh-boundary-capture}
Let $D\subseteq\Mesh$ be finite, and put $A=A(D)$ and $S=A^{\mathrm{ov}}$.
Then the following assertions hold.
\begin{enumerate}
\item
\label{l.mesh-boundary-capture.core}
Every mesh box of $D$ belongs to the mesh interior $\MeshInt(S)$.

\item
Define the reduction of $S$ by
\[
 A^\circ=\partial_{\rm ext}^{\Mesh}\MeshInt(S),
 \qquad
 S^\circ=(A^\circ)^{\mathrm{ov}}.
\]
Then $A^\circ\subseteq A$ and $S^\circ\subseteq S$.  Moreover,
\begin{equation}
\label{e.mesh-reduction-same-exterior}
 \Ext_{\overline N}(S^\circ)=\Ext_{\overline N}(S),
 \qquad
 \Capt_{\overline N}(S^\circ)=\Capt_{\overline N}(S).
\end{equation}
The mesh interior $\MeshInt(S^\circ)$ is obtained from $\MeshInt(S)$ by filling its bounded
nearest-neighbour complementary components, and
\[
 S^\circ=(S^\circ\cap\Mesh)^{\mathrm{ov}},
 \qquad
 S^\circ\cap\Mesh=\partial_{\rm ext}^{\Mesh}\MeshInt(S^\circ).
\]
\end{enumerate}
\end{lemma}

The last identities say, in particular, that the reduction clears out the
buried boundary in one pass: reducing $S^\circ$ again leaves $S^\circ$
unchanged.

\noindent{\bf Proof of Lemma~\ref{l.mesh-boundary-capture}: (1).}
Let $\widehat D$ be the filling of $D$.  The unbounded component of
$\Mesh\setminus D$ is also the unbounded component of
$\Mesh\setminus\widehat D$, and $D\subseteq\widehat D$, so
$A(D)\subseteq A(\widehat D)$.  For the reverse inclusion, let
$x\in A(\widehat D)$.  If $x$ is $*$-adjacent to a site of $D$, then
$x\in A(D)$.  Otherwise it is $*$-adjacent to a site $y$ in a bounded
component of $\Mesh\setminus D$.  Resolve the $*$-step from $y$ to $x$
one coordinate at a time.  The resulting nearest-neighbour path must meet
$D$, while every one of its vertices lies within $\ell^\infty$-distance one
of $x$.  Thus $x$ is $*$-adjacent to $D$, and again $x\in A(D)$.  Hence
$A(\widehat D)=A(D)$.  Replacing $D$ by $\widehat D$, we may and do assume
that $D$ is filled; this does not change $A$ or $S$.

By the definition of $\MeshInt(S)$ and Lemma~\ref{l.mesh-transfer},
\[
 \MeshInt(S)=\MCap(S)\setminus A
     =(\Mesh\setminus A)\setminus U^*(A).
\]
Suppose that $M\in D$ belonged to the unbounded $*$-component $U^*(A)$.
Follow a $*$-path in $\Mesh\setminus A$ from $M$ to infinity and take its
first site outside $D$.  Since $D$ is filled, this site lies in the unbounded
nearest-neighbour component of $\Mesh\setminus D$ and is $*$-adjacent to
$D$, so it belongs to $A$, a contradiction.  Hence $D\subseteq \MeshInt(S)$.

{\bf (2).}
Put $I=\MeshInt(S)$.  If $X\in A^\circ\setminus A$, then
$X\in\Mesh\setminus A$ is $*$-adjacent to some $Y\in I$.  The two sites
therefore belong to the same $*$-component of $\Mesh\setminus A$,
contradicting $X\notin I$.  Thus $A^\circ\subseteq A$, and hence
$S^\circ\subseteq S$.

Every member of $A$ is $*$-adjacent to the filled set $D\subseteq I$.
Consequently, on each nearest-neighbour component of $A$, either all of its
mesh boxes belong to the unbounded nearest-neighbour component of
$\Mesh\setminus I$, or none do.  Thus $A^\circ$ is a union of
nearest-neighbour components of $A$.

Let $F=\Capt_{\Mesh}(I)$, the filling of $I$.  Every site of $U^*(A)$ is
nearest-neighbour connected to infinity in
$A\cup U^*(A)=\Mesh\setminus I$: resolve an infinite $*$-path one coordinate
at a time.  Any inserted site outside $A$ is $*$-adjacent to a site of
$U^*(A)$ and therefore belongs to the same $*$-component.  Since every member
of $A$ is $*$-adjacent to $I$, it lies in the unbounded nearest-neighbour
complement of $I$ exactly when it belongs to $A^\circ$.  Hence
\begin{equation}
\label{e.mesh-reduction-decomposition}
 F=I\,\sqcup\,(A\setminus A^\circ),
 \qquad
 \Mesh\setminus F=A^\circ\,\sqcup\,U^*(A).
\end{equation}
Thus passing from $I$ to its filling $F$ absorbs exactly the boundary boxes
in $A\setminus A^\circ$; these are the parts of the original boundary that
have become buried.  The surviving set $A^\circ$ is the part of $A$ that
remains exposed to the unbounded complement.

No site of $U^*(A)$ is $*$-adjacent to $F$.  Adjacency to $I$ is impossible,
since it would join the unbounded $*$-component of $\Mesh\setminus A$ to one
of its bounded $*$-components.
If $u\in U^*(A)$ were instead $*$-adjacent to
$x\in A\setminus A^\circ$, resolve that step one coordinate at a time.  An
intermediate site in $I$ would be $*$-adjacent to $u$, again a contradiction;
if there is no such intermediate site, the resolved path lies in
$\Mesh\setminus I$ and joins $x$ to its unbounded component, which would put
$x$ in $A^\circ$.  It follows that
\[
 A^\circ=\partial_{\rm ext}^{\Mesh}F.
\]
Indeed, a site of $A^\circ$ has an infinite nearest-neighbour path in
$\Mesh\setminus I$.  This path cannot pass from $A^\circ$ to
$A\setminus A^\circ$ within $A$, since $A^\circ$ is a union of
nearest-neighbour components of $A$; nor can it enter
$A\setminus A^\circ$ from $U^*(A)$, since no site of $U^*(A)$ is
$*$-adjacent to $F$.  Thus the path lies in $\Mesh\setminus F$.
Conversely, a site of the exterior boundary of $F$ lies, by
\eqref{e.mesh-reduction-decomposition}, in $A^\circ\cup U^*(A)$, and the
non-adjacency just proved excludes the second alternative.

We next claim that $U^*(A^\circ)=U^*(A)$.  The inclusion
$U^*(A)\subseteq U^*(A^\circ)$ follows from $A^\circ\subseteq A$.
Conversely, $F$ is filled and has exterior mesh boundary $A^\circ$.  A
$*$-path beginning in $F$ and leaving $F$ must at its first exit visit
$A^\circ$, so no site of $F$ is $*$-connected to infinity in
$\Mesh\setminus A^\circ$.  Equation~\eqref{e.mesh-reduction-decomposition}
then gives the reverse inclusion.  Lemma~\ref{l.mesh-transfer}, applied to
$A$ and $A^\circ$, now yields~\eqref{e.mesh-reduction-same-exterior}.

The same lemma and~\eqref{e.mesh-reduction-decomposition} give
\[
 \MeshInt(S^\circ)=F,
 \qquad
 S^\circ\cap\Mesh=A^\circ=\partial_{\rm ext}^{\Mesh}F.
\]
Thus $\MeshInt(S^\circ)$ is the filling of $\MeshInt(S)$ and
$S^\circ\cap\Mesh=\partial_{\rm ext}^{\Mesh}\MeshInt(S^\circ)$.  Since a mesh box
belongs to $(A^\circ)^{\mathrm{ov}}$ exactly when it belongs to $A^\circ$, also
$S^\circ=(S^\circ\cap\Mesh)^{\mathrm{ov}}$. \qed

\begin{lemma}[Boundary contacts of a captured excursion]
\label{l.mesh-boundary-screening}
Let $D\subseteq\Mesh$ be finite, put $S=A(D)^{\mathrm{ov}}$, and let $S^\circ$ be its
reduction.  If a nearest-neighbour box path has a finite segment contained in
$\Capt_{\overline N}(S^\circ)$ whose first and last boxes belong to
$S^\circ$ and whose intermediate boxes belong to
$\Capt_{\overline N}(S^\circ)\setminus S^\circ$, then its two contacts with
$S^\circ$ belong to the same nearest-neighbour component of $S^\circ$.
\end{lemma}

\noindent{\bf Proof.}
The assertion is immediate if there is no intermediate box.  Otherwise put
$A^\circ=S^\circ\cap\Mesh$ and $F=\MeshInt(S^\circ)$.  By
Lemma~\ref{l.mesh-boundary-capture}(2), $F$ is filled and
$A^\circ=\partial_{\rm ext}^{\Mesh}F$.

From the mesh support of each intermediate box choose a member outside
$A^\circ$.  Such a member exists because the box is not in $S^\circ$, and it
must belong to $F$: a member in the mesh exterior would, after resolving odd
coordinates towards it, join the intermediate box to infinity without meeting
$S^\circ$.  For consecutive intermediate boxes the supports are nested, so
the chosen mesh boxes are equal or $*$-adjacent.  They therefore lie in one
$*$-component $C$ of $F$.

At either end, let $B\in S^\circ$ be the contact box and $B'$ the adjacent
intermediate box.  The support of $B'$ contains a member of $F$, whereas the
support of $B$ is contained in $A^\circ$.  Since the supports are nested,
and the reverse inclusion would put $\MeshSupp(B')$ inside $A^\circ$,
we have $\MeshSupp(B)\subseteq\MeshSupp(B')$.  Every member of $\MeshSupp(B)$ is
consequently $*$-adjacent to $C$.  Since $C$ is a $*$-component of $F$, its
exterior $*$-boundary is disjoint from $F$.  As $F$ is filled, every site of
this boundary lies in the unbounded complement of $F$; being $*$-adjacent to
$F$, it therefore belongs to $\partial_{\rm ext}^{\Mesh}F=A^\circ$.
Lemma~\ref{l.extbdy-slab} makes this boundary nearest-neighbour connected.
The mesh support of each contact box lies in this boundary, so
Lemma~\ref{l.overlap-mesh-region} joins the two contacts in one
nearest-neighbour component of $S^\circ$. \qed

\begin{lemma}[Exterior boxes meet the capture at the boundary]
\label{l.mesh-boundary-interface}
Let $D\subseteq\Mesh$ be finite, put $A=A(D)$ and $S=A^{\mathrm{ov}}$.  If
$B\in\Ext_{\overline N}(S)$ and
$M\in\MeshSupp(B)\cap\MCap(S)$, then $M\in A$.  The same assertion holds
for the reduced pair $(S^\circ,A^\circ)$.
\end{lemma}

\noindent{\bf Proof.}
If $M\notin A$, resolving the odd transverse coordinates of $B$ towards $M$
gives a path outside $S$ from $B$ to $M$.  Thus $M$ would belong to the
exterior component of $S$, contrary to $M\in\MCap(S)$.  The same argument
applies to $(S^\circ,A^\circ)$. \qed

\section{Welcoming boundaries and the halo}
\label{ss.visible-welcoming-exteriors}

Let $D$ be a finite set of mesh boxes.  Surrounding its
overlapping-box boundary $S=A(D)^{\mathrm{ov}}$ are exterior overlapping
boxes that nonetheless contain captured edges.  We call these the
\emph{halo boxes}, and put
\[
 \Halo(S)=\{B\in\Ext_{\overline N}(S):E(B)\cap\ECap(S)\ne\emptyset\}.
\]
For $B\in\Ext_{\overline N}(S)\setminus\Halo(S)$, no edge of $B$ is
captured, so conditioning on the captured edge states leaves the law in $B$
unchanged and hence does not affect its probability of being apt.  For a halo
box $B$, by contrast, part of $E(B)$ lies on the captured side, and
conditioning on these edges may make aptness substantially less likely.  To
locate useful quilt structure without querying the exterior, we will impose a
welcoming condition requiring that halo boxes are likely to be apt after this
captured-side conditioning.  For $B\in\Ext_{\overline N}(S)$ put
\[
 \mathcal I_B(S)=E(B)\cap\ECap(S),\qquad
 q_B(S)=\PP\bigl(B\text{ is not apt}\mid\cF_{\mathcal I_B(S)}\bigr).
\]
For a halo box $B$, we call $\mathcal I_B(S)$ its
\emph{captured-side conditioning set}.
Lemma~\ref{l.mesh-boundary-interface} gives
\begin{equation}
\label{e.halo-reads-seal}
   \mathcal I_B(S)\subseteq E(S),\qquad B\in\Halo(S).
\end{equation}
Thus $q_B(S)$ is measurable from captured-side boundary data, while the
exterior edge states retain their product law.

To formulate the welcoming condition uniformly, we need an allowance for
$q_B(S)$ that depends only on how many boxlets of $B$ lie on the captured
side, with enough reserve for further captured-side exposure.  The one-boxlet
case shows the mechanism.  Suppose first that the captured-side conditioning
set consists of one prescribed boxlet $L$ of an overlapping box $B$.  Write
\[
 q=\PP\bigl(B\text{ is not apt}\mid\cF_{E(L)}\bigr).
\]
Since $\EE q=\PP(B\text{ is not apt})\le\varepsilon_K$, Markov's
inequality gives
\[
 \PP\bigl(q>\sqrt{\varepsilon_K}\bigr)\le\sqrt{\varepsilon_K}.
\]
Thus, with probability at least $1-\sqrt{\varepsilon_K}$, the configuration
seen on this captured-side boxlet leaves $B$ likely to be apt: its conditional
probability of inaptitude is at most $\sqrt{\varepsilon_K}$.

In general the captured-side conditioning set of a halo box may contain
several boxlets, in a pattern determined by the interface geometry.  Let
$\mathscr L(B)$ be the set of boxlets of $B$, and write
\[
       m_d=|\mathscr L(B)|=2^{d-1}.
\]
For $\mathscr A\subseteq\mathscr L(B)$, put
\[
 I_B(\mathscr A)=\bigcup_{L\in\mathscr A}E(L),
 \qquad
 q_{B,\mathscr A}
   =\PP\bigl(B\text{ is not apt}\mid\cF_{I_B(\mathscr A)}\bigr).
\]
Because boxlet edge sets partition the layer, every captured-side conditioning
set $\mathcal I_B(S)$ is of this form.  Thus there are at most $2^{m_d}$
possible boxlet conditioning sets, independently of~$K$.

The one-boxlet estimate may now be iterated.  Take the boxlets in any
convenient order and successively adjoin their edge states to the conditioning
field.  At each enlargement, the tower property and Markov's inequality give
the next square-root allowance.  Since the tower property applies to every
nested pair of boxlet conditioning sets, the resulting allowances may be
graded simply by the number of boxlets present.  Put
\[
  \eta_{K,0}=\varepsilon_K,
  \qquad
  \eta_{K,j}=\varepsilon_K^{2^{-j}},\quad 1\le j\le m_d-1,
  \qquad
  \delta_K=\varepsilon_K^{2^{-(m_d-1)}}.
\]
For large $K$, these numbers lie in $(0,1)$ and
\[
 \eta_{K,0}<\eta_{K,1}<\cdots<\eta_{K,m_d-1}\le\delta_K,
 \qquad
 \frac{\eta_{K,i}}{\eta_{K,j}}\le\delta_K
       \quad(0\le i<j\le m_d-1).
\]
Thus the allowance at level $j$ depends only on how many boxlets enter the
conditioning set.  For any nested pair $\mathscr A\subsetneq\mathscr A'$, the
ratio bound controls passage from the smaller conditioning field to the larger
one.  Since only finitely many subsets are possible, simultaneous control costs
only a dimension-dependent factor, while $\delta_K$ still tends to zero
exponentially as $K\to\infty$.

For $B\in\Halo(S)$, let
\[
 \mathscr L_B(S)
   =\{L\in\mathscr L(B):E(L)\subseteq\mathcal I_B(S)\},
 \qquad r_B(S)=|\mathscr L_B(S)|.
\]
Then
\[
 \mathcal I_B(S)=I_B(\mathscr L_B(S)),
 \qquad 1\le r_B(S)\le m_d-1.
\]
The lower bound is the definition of the halo.  For the upper bound,
\eqref{e.mesh-transfer-full-exterior} gives an exterior mesh box in
$\MeshSupp(B)$, so at least one boxlet of $B$ lies outside
$\mathcal I_B(S)$.
We call the deterministic boundary $S$ \emph{welcoming} if every box of
$S$ is apt and
\begin{equation}
\label{e.welcoming-halo}
   q_B(S)\le\eta_{K,r_B(S)}
       \qquad\text{for every }B\in\Halo(S).
\end{equation}
Thus a halo box is required to remain apt with high conditional probability
given exactly its captured-side conditioning set.  The threshold depends only
on the number of boxlets in this set, not on their arrangement,
and retains reserve for further captured-side exposure while the box remains
exterior.
For deterministic $S=A(D)^{\mathrm{ov}}$, write $\mathrm{Wel}(S)$ for the event that
$S$ is welcoming.  Aptness of the boxes in $S$ is determined by their own
edges, while~\eqref{e.halo-reads-seal} places every halo conditioning set
within $E(S)$.  Thus whether $S$ is welcoming is determined by the edge states
in $E(S)$; equivalently, $\mathrm{Wel}(S)\in\cF_{E(S)}$.

\begin{lemma}[Finite-boxlet conditioning reserve]
\label{l.sqrt-welcoming-note}
Let $\mathscr A\subsetneq\mathscr A'\subsetneq\mathscr L(B)$, with
$i=|\mathscr A|$ and $j=|\mathscr A'|$.  Then
\begin{equation}
\label{e.welcome-martingale}
 \EE\bigl(q_{B,\mathscr A'}\mid\cF_{I_B(\mathscr A)}\bigr)
   =q_{B,\mathscr A}.
\end{equation}
Consequently, on
$\{q_{B,\mathscr A}\le\eta_{K,i}\}$,
\begin{equation}
\label{e.welcome-reserve-transition}
 \PP\bigl(q_{B,\mathscr A'}>\eta_{K,j}
      \mid\cF_{I_B(\mathscr A)}\bigr)
 \le\delta_K,
 \qquad
 \PP\bigl(B\text{ is not apt}\mid\cF_{I_B(\mathscr A)}\bigr)
 \le\delta_K.
\end{equation}
\end{lemma}

When $\mathscr A=\varnothing$, the hypothesis in the lemma is automatic:
$q_{B,\varnothing}=\PP(B\text{ is not apt})\le\varepsilon_K=\eta_{K,0}$.

\noindent{\bf Proof of Lemma~\ref{l.sqrt-welcoming-note}.}
Identity~\eqref{e.welcome-martingale} is the tower property.  Conditional
Markov inequality gives
\[
 \PP\bigl(q_{B,\mathscr A'}>\eta_{K,j}
      \mid\cF_{I_B(\mathscr A)}\bigr)
 \le \frac{q_{B,\mathscr A}}{\eta_{K,j}}
 \le \frac{\eta_{K,i}}{\eta_{K,j}}
 \le\delta_K.
\]
The second estimate in~\eqref{e.welcome-reserve-transition} is the definition
of $q_{B,\mathscr A}$, together with $\eta_{K,i}\le\delta_K$. \qed

The proof uses only the Bernoulli-$p$ product law on the as-yet
unconditioned edges of $B$.  Consequently, if additional information has been
conditioned on but the edges in $E(B)\setminus I_B(\mathscr A)$ retain this
product law, then \eqref{e.welcome-martingale} and
\eqref{e.welcome-reserve-transition} remain valid relative to the enlarged
conditioning field.  We use this observation later in the conditional
seal-completion argument.

For an overlapping box \(B\), define the local defect event
\begin{equation}
\label{e.absolute-seal-bad-box}
 \mathsf E_B
 =\{B\text{ is not apt}\}
   \cup
   \bigcup_{\substack{\varnothing\ne\mathscr A\subsetneq\mathscr L(B)}}
   \{q_{B,\mathscr A}>\eta_{K,|\mathscr A|}\}.
\end{equation}
By Proposition~\ref{p.nonaxial-ap}, Lemma~\ref{l.sqrt-welcoming-note} and a
union bound over at most \(2^{m_d}\) boxlet conditioning sets,
\begin{equation}
\label{e.failed-seal-test-density}
       \PP(\mathsf E_B)\le C_d\delta_K.
\end{equation}
The event \(\mathsf E_B\) is determined by \(E(B)\), so the field
\((\mathbf1_{\mathsf E_B})\) is one-dependent.

For \(S=A(D)^{\mathrm{ov}}\), write \(S^\circ\) for the reduction supplied by
Lemma~\ref{l.mesh-boundary-capture}(2).

\begin{lemma}[Reduction preserves welcome]
\label{l.reduction-preserves-welcome}
For deterministic \(S=A(D)^{\mathrm{ov}}\),
\[
 \ECap(S^\circ)=\ECap(S),\qquad
 \mathcal U(S^\circ)=\mathcal U(S),\qquad
 \Halo(S^\circ)=\Halo(S).
\]
If \(\mathrm{Wel}(S)\) holds, then so does \(\mathrm{Wel}(S^\circ)\).
\end{lemma}

\noindent{\bf Proof.}
Lemma~\ref{l.mesh-boundary-capture}(2) gives equality of the captured and
exterior box regions, hence the first two displayed identities.  The halo is
therefore unchanged; on this common halo, each set \(\mathcal I_B\) and its
conditioning rank \(r_B\) are unchanged as well, so the halo inequalities
are unchanged.  Finally \(S^\circ\subseteq S\), and therefore every box of
\(S^\circ\) is apt whenever every box of \(S\) is apt.
\qed

\section{Exterior percolation and open escape}
\label{s.exterior-percolation-open-escape}

We conclude the chapter with the uniform open-escape estimate.  Given
a welcoming boundary, any prescribed vertex incident to an unread
exterior edge has a uniformly positive conditional probability of connection
to infinity in the layer.  The first subsection constructs an open route to
infinity through the exterior and the boundary quilts; the second attaches the
prescribed vertex to this route.

\subsection{Exterior routes}
\label{s.viable-routes-slab}

For a welcoming boundary, conditioning on the captured edge states leaves a
highly supercritical field of apt boxes throughout the exterior,
including---by the welcoming condition---the surrounding halo.  By declaring
the captured boxes open, we extend this to an auxiliary field on the full box
lattice.  From a given exterior starting box, this field has a uniformly
positive chance to carry a path to infinity.  Such a path may make finite
excursions through captured boxes, which have been declared open only in the
auxiliary field.  Screening reroutes these excursions along the boundary
quilts, producing a path to infinity of apt boxes; their woven clusters then
yield an open route to infinity.

\subsubsection{Conditional percolation}

Fix deterministic $S=A(D)^{\mathrm{ov}}$ and work conditionally on the captured
edge states, on the event $\mathrm{Wel}(S)$.  For an exterior box $B$, its
captured edges are exactly those in $\mathcal I_B(S)$, while its remaining
edges belong to $\mathcal U(S)$.  Product independence therefore gives
\begin{equation}
\label{e.conditional-inapt-q}
 \PP(B\text{ is not apt}\mid\cF_{\ECap(S)})=q_B(S).
\end{equation}

\begin{lemma}[Conditional density of inapt boxes]
\label{l.low-density-unusable-exterior-boxes}
Let $\rho>0$ and take $K$ sufficiently large. Uniformly in $S$, on
$\mathrm{Wel}(S)$ and conditionally on $\cF_{\ECap(S)}$, the inapt-box
field on $\Ext_{\overline N}(S)$ is one-dependent, with every one-box
failure probability at most $\rho$.
\end{lemma}

\noindent{\bf Proof.} For a halo box, \eqref{e.conditional-inapt-q} and the halo condition in
\(\mathrm{Wel}(S)\) give failure probability at most \(\delta_K\).  Outside the
halo, no captured edge belongs to the box, so Proposition~\ref{p.nonaxial-ap}
gives the bound $\varepsilon_K$.  Boxes at transverse index distance at least
two use disjoint exterior edge sets and are therefore conditionally
independent.  Since $\delta_K\vee\varepsilon_K\to0$, the result follows. \qed

\begin{proposition}[LSS domination]
\label{p.lss-usable-boxes}
Let $q\in(0,1)$ and take $K$ sufficiently large.  Uniformly in $S$, on
$\mathrm{Wel}(S)$ and conditionally on $\cF_{\ECap(S)}$, the auxiliary field
on the full box lattice defined by
\[
 X_B=\begin{cases}
       1,&B\in\Capt_{\overline N}(S),\\
       \mathbf 1\{B\text{ is apt}\},&B\in\Ext_{\overline N}(S)
     \end{cases}
\]
stochastically dominates Bernoulli site percolation of parameter $q$ on
$\Z^{d-1}$.
\end{proposition}
\noindent{\bf Proof.} Captured boxes are deterministically open in the auxiliary field.  On the
exterior, Lemma~\ref{l.low-density-unusable-exterior-boxes} gives a
one-dependent field with one-box failure probability tending to zero with
$K$.  By Liggett--Schonmann--Stacey~\cite{LSS}, this gives the stated
stochastic domination. \qed

\subsubsection{Screening}

The auxiliary path may make finite excursions through boxes of the capture
which are open only in the auxiliary field.  Screening removes these auxiliary
portions at the box level: each excursion is rerouted along one quilt component
of the reduced boundary.  After loop erasure, this yields an infinite self-avoiding path
consisting entirely of apt boxes.  We then use their woven clusters to pass
from this coarse path to open passage in the lattice.

\begin{lemma}[Screening by the reduced boundary]
\label{l.screened-holes-site-path}
Let $S_0=A(D)^{\mathrm{ov}}$ and let $S=S_0^\circ$ be its reduction.  Suppose that
the boxes of $S$ are apt.  An infinite self-avoiding path in the box lattice
which starts in the common exterior of $S_0$ and $S$, and whose exterior boxes
are apt, may be modified, without changing its initial box, to an infinite self-avoiding
path every box of which is apt.  The modification replaces each
finite excursion into the common capture of $S_0$ and $S$ by a path in one quilt component
of $S$ and then erases loops.
\end{lemma}

\noindent{\bf Proof.}
By Lemma~\ref{l.mesh-boundary-capture}(2), $S_0$ and $S$ have the same
capture and exterior.  Consider a maximal finite segment of the path
in the common capture of $S_0$ and $S$.  Its first and last boxes belong to $S$: a box of
$\Capt_{\overline N}(S)\setminus S$ which were adjacent to the exterior
would itself be connected to infinity in $\overline N\setminus S$.  Break
the segment at all of its visits to $S$.  Each intervening portion is either
already a path in $S$ or has its two endpoint contacts in $S$ and all
intermediate boxes in $\Capt_{\overline N}(S)\setminus S$.  In the latter
case Lemma~\ref{l.mesh-boundary-screening} replaces it by a path in one
component of $S$.

Since the capture is finite and the original path is self-avoiding, only
finitely many such excursions occur.  After making all of the replacements,
the path agrees with the original from some exterior box onward.  The modified
initial portion is disjoint from this unchanged tail except at their common
endpoint: its new pieces lie in the capture, while the original path is
self-avoiding.  Loop-erasing this initial portion, starting from the original
initial box, therefore preserves both that box and the tail to infinity.
Every exterior box retained from the original path is apt by hypothesis, while
each replacement lies in a quilt component of $S$ and hence consists of apt
boxes.  The resulting infinite self-avoiding path has the asserted
properties. \qed

\begin{lemma}[Open passage along an apt box path]
\label{l.viable-route-open-passage}
Suppose a path to infinity in the box lattice begins at~$B_0$ and every box
on it is apt.  Then $\mathcal W(B_0)$ is connected to infinity by an open path
in~$N$.
\end{lemma}

\noindent{\bf Proof.} Consecutive boxes have intersecting woven clusters by
Lemma~\ref{l.woven-fabric-note}.  The union of these woven clusters is an
unbounded connected open subgraph of $N$ containing $\mathcal W(B_0)$, and
hence contains an open path to infinity from $\mathcal W(B_0)$. \qed

\subsection{Uniform open escape}
\label{s.open-escape-seed}

For a deterministic boundary $S$, condition on the captured edge states
$\cF_{\ECap(S)}$; the edges in $\mathcal U(S)$ are unread relative to this
conditioning.  The remaining step is local: an unread edge incident to the
prescribed vertex lies in a boxlet all of whose edges are unread.
Opening this boxlet at bounded finite-energy cost supplies the attachment to
the exterior route.

\begin{lemma}[The boxlet of an unread edge is fresh]
\label{l.unread-edge-launch-cell}
Let \(S=A(D)^{\mathrm{ov}}\) and let \(e\in\mathcal U(S)\).  Let \(L\) be the
unique boxlet to which \(e\) is assigned.  Then
\[
        E(L)\subseteq\mathcal U(S).
\]
Moreover every overlapping box in \(\overline N\) having \(L\) as a
boxlet belongs to \(\Ext_{\overline N}(S)\).
\end{lemma}

\noindent{\bf Proof.} Let \(M\in\Mesh\) be the unique mesh box containing the boxlet \(L\).  Since
\(e\in\mathcal U(S)\) and the mesh edge sets are disjoint, \(M\in\MExt(S)\).
The whole boxlet \(L\) belongs to this same mesh box, so
\(E(L)\subseteq E(M)\subseteq\mathcal U(S)\).

Now let \(B\) contain \(L\).  Then \(M\in\MeshSupp(B)\).
Certainly \(B\notin S=A(D)^{\mathrm{ov}}\), because every member of \(S\) has
all of its mesh support in \(A(D)\subseteq\MCap(S)\).  If \(B\) were in
\(\Capt_{\overline N}(S)\setminus S\), resolve its odd transverse
coordinates towards the mesh box \(M\).  Every box on the resulting
bounded path has \(M\) in its mesh support and hence does not belong to
\(S\).  The path joins \(B\) to the exterior mesh box \(M\) without meeting
\(S\), contradicting capture.  Thus \(B\) is exterior. \qed

The exterior route and the fresh boxlet supplied by the preceding lemma now
combine to give the promised uniform open-escape estimate.

\begin{proposition}[Uniform open escape from an exterior vertex]
\label{p.uniform-open-escape-exterior-seed}
Let \(d\ge3\).  For all sufficiently large fixed \(K\) there is
\(c_{\rm esc}=c_{\rm esc}(d,p,K)>0\) with the following property.
Let \(S=A(D)^{\mathrm{ov}}\) be deterministic, and let \(s\in V(N)\) be incident to at
least one edge of \(\mathcal U(S)\).  Let \(L\) be the boxlet to which the
first such edge is assigned in a fixed deterministic ordering.  Then, almost
surely on \(\mathrm{Wel}(S)\),
\[
 \PP\bigl(s\longleftrightarrow\infty\text{ in }N
      \mid\cF_{\ECap(S)}\bigr)\ge c_{\rm esc}.
\]
\end{proposition}

\noindent{\bf Proof.} By Lemma~\ref{l.unread-edge-launch-cell},
$E(L)\subseteq\mathcal U(S)$ and every overlapping box containing $L$ is
exterior.  Choose one such box $B_0$ by a fixed deterministic rule.

Choose
\(q\in(p_c^{\rm site}(\Z^{d-1}),1)\) and write
\(a_{\rm slab}=\PP_q(0\longleftrightarrow\infty)>0\).
Proposition~\ref{p.lss-usable-boxes} gives, with conditional probability at
least \(a_{\rm slab}\), a path to infinity from $B_0$ in the auxiliary open
field.  On this event choose such a path self-avoiding.  Let $S^\circ$ be the
reduction of $S$.  By Lemma~\ref{l.reduction-preserves-welcome}, it has the
same exterior and capture as $S$ and is welcoming.  Lemma~\ref{l.screened-holes-site-path}
modifies the auxiliary path to an infinite self-avoiding path of apt boxes,
still beginning at $B_0$.  Lemma~\ref{l.viable-route-open-passage} therefore
connects $\mathcal W(B_0)$ to infinity.  Since $B_0$ is apt and $L$ is one of
its boxlets, $\mathcal W(B_0)\cap L$ contains an open path of diameter at
least $3K/4$.  Hence a vertex of $L$ is connected to infinity.

Let
\[
 E_\infty^S(L)=\{\text{some vertex of }L\text{ is joined to infinity in }N
       \text{ using only edges of }E(S)\cup\mathcal U(S)\}.
\]
The route just constructed uses only these edges.  Indeed, its screened
portions use edges of $S^\circ\subseteq S$, while for every exterior box its
captured-side edges, when present, lie in
$\mathcal I_B(S)\subseteq E(S)$ by \eqref{e.halo-reads-seal}; all remaining
edges lie in $\mathcal U(S)$.  Thus
\(\PP(E_\infty^S(L)\mid\cF_{\ECap(S)})\ge a_{\rm slab}\).

Let \(F_L\) be the event that every edge of \(L\) is open.  By
Lemma~\ref{l.connected-layer-cells} and Lemma~\ref{l.unread-edge-launch-cell},
\[
       |E(L)|\le C_dK^d,
       \qquad
       \PP(F_L\mid\cF_{\ECap(S)})=p^{|E(L)|}\ge p^{C_dK^d}.
\]
Conditionally on $\cF_{\ECap(S)}$, the edge states on $\mathcal U(S)$ have
their Bernoulli-$p$ product law.  Both $F_L$ and $E_\infty^S(L)$ are increasing
events of this product field, so the Harris--FKG inequality gives
\[
 \PP(F_L\cap E_\infty^S(L)\mid\cF_{\ECap(S)})
 \ge p^{C_dK^d}a_{\rm slab}.
\]
On \(F_L\), all edges of the connected boxlet graph \(L\) are open.  Since
\(s\) is incident to the selected edge of \(L\), it is connected inside
\(L\) to every vertex of \(L\), in particular to a vertex that reaches
infinity on \(E_\infty^S(L)\).  Thus one may take
\(c_{\rm esc}=p^{C_dK^d}a_{\rm slab}\). \qed

\chapter{Seals: enclosure and enlargement}
\label{s.seed-exploration-seal-extension}

Fix one of the coarse layers used to monitor the target-height downcluster.
As the slide dynamics proceeds, the cluster in the layer changes slide by
slide.  We enclose the finite open-and-closed record left by the exploration in
a finite stopping-domain quilt collection without inspecting the exterior.
Conditionally on the revealed information, the exterior therefore retains its
Bernoulli product law.  Such an enclosure is a \emph{seal}.

Two forms of sealing are needed.  In \emph{enclosure sealing}, a seal is
placed so as to enclose a prescribed finite region in the layer.  In
\emph{enlargement sealing}, one starts from a seal already in force and builds
it outward to enclose additional finite material bordering its capture.

The slide dynamics gives both operations natural roles.  During the advance
phase the slide height lies south of the layer.  At passage, the cluster
material in the layer consists of the part of the transported upcluster
carried by the last active advance slide, the \emph{host}, together with the
\emph{guest}, if any, formed by its midpiece.  A one-off \emph{initial seal}
is placed around the host by enclosure sealing; any part of the guest outside
its capture is then absorbed by enlargement sealing.  As the slide heights
advance into the post-passage phase, a second reason for enclosure sealing
appears.  A seal constructed in one configuration need
not remain valid once the next active slide resampling takes effect.  The finite
region previously enclosed is therefore handed on, and enclosure sealing is
used to place a fresh seal around it in the new configuration.  Once a seal is
valid in the current configuration, further finite growth is again handled by
enlargement sealing.  The precise way in which these
operations are organised along the slide dynamics is developed in
Chapter~\ref{s.seal-sequence-layer}; the present chapter develops the local
sealing tools themselves.

We work for now in a single Bernoulli product field on a fixed
layer.  Suppose that a seal is already in place and that an exploration
reaches the layer at one or more vertices.  These
\hfff{seed}{\emph{seeds}} provide starting points for further growth.  A seed
is \emph{virgin} or \emph{sealed} according as it is or is not incident to
the unread exterior of the current seal.  A virgin seed has immediate access
to untouched product randomness.  From a sealed seed, by contrast, the
exploration first follows continuation through the current capture, querying
a captured edge when it is first needed.  Such continuation may die
out, may reach the side of the layer, may locate a virgin seed, or may meet the
woven fabric of a component of the seal.

The exploration is organised so that captured continuation is exhausted
before fresh exterior growth is pursued.  Encountering the woven fabric of
a seal component activates that component: the known open fabric is followed
as part of the captured scan, and may expose further virgin seeds.  Exterior
growth is launched from a virgin seed.  On open escape in the layer no further
sealing is needed; after a finite episode its newly revealed record is
absorbed by enlarging the seal, and the remaining continuation is reconsidered
relative to the enlarged capture.

We begin by constructing seals around prescribed finite regions.  In the
second section, we study enlargement of a seal already in force, including the
conditional estimates needed for successive episodes.  The third section gives
the initial-seal estimate and organises further growth into a \emph{paused
exploration}, alternating captured continuation with episodes launched from
virgin seeds.

We end with a freestanding control problem that puts the sealing techniques
into action.  Beginning from a given finite record of discovered open and
closed edges, an Encloser seeks to accumulate a large but finite volume of
open growth by repeatedly running percolation explorations in a layer.  She
plants seeds adaptively, guided by what earlier explorations have revealed,
and watches the ensuing growth.  She may stop and collect her winnings at
any stage, but an open escape wipes out her score.  The sealing estimates
show that, however she plays, her winnings have a stretched-exponential
tail beyond the isoperimetric scale of the initial finite record.

\section{Seal construction}

We gather the local ingredients for constructing seals.  We first describe
seeds and their exploration relative to an existing seal, then develop the
home-box and boundary-crossing geometry used to surround finite growth by an
enclosing quilt.  Finally, we complete such a quilt to a seal, controlling how
far the completion search must range.

\subsection{Seals and virgin seeds}
\label{ss.seals-virgin-seeds}

We first formulate the sealing mechanism on a Bernoulli-\(p\) product field
in the fixed layer.  The sealing constructions use only this local product
structure.  Their insertion into the fixed-shift explorations of the slide
dynamics is carried out in Chapter~\ref{s.seal-sequence-layer}.

Write \(N\) for the fixed partition layer, with box lattice \(\overline N\) and mesh \(\Mesh\) as in Section~\ref{s.slab-setting}.  Recall that
$\MeshInt(S)$ denotes the mesh interior of a finite box collection $S$, namely the
captured mesh boxes lying strictly behind~$S$.

\begin{definition}[Seal]
\label{d.seal}
A \hfff{seal}{\emph{seal}} is a random finite quilt collection
$\mathscr S\subseteq\overline N$ for which, almost surely, $\mathscr S$ is
welcoming and
\begin{equation}
\label{e.seal-mesh-geometry}
   \mathscr S=(\mathscr S\cap\Mesh)^{\mathrm{ov}},
   \qquad
   \mathscr S\cap\Mesh=\partial_{\rm ext}^{\Mesh}\MeshInt(\mathscr S);
\end{equation}
and, for every deterministic finite collection $S$,
\[
   \{\mathscr S=S\}\in\cF_{\ECap(S)}.
\]
\end{definition}
The first identity says that the whole quilt collection is determined by its
mesh boxes; the second says that these mesh boxes form the outward-facing
boundary of its mesh interior.  Thus a seal is a welcoming stopping-domain quilt
collection; it may have several quilt components, and we also permit the empty
seal.
Such collections exist: Section~\ref{ss.trap-cluster-to-seal} constructs them.
When an output $\mathscr S'$ is defined only on an event $F$, the assertion
that it is a seal on $F$ means that its geometry and welcoming
property hold on $F$ and, for every deterministic $S$,
$\{F,\mathscr S'=S\}\in\cF_{\ECap(S)}$.

For a seal $\mathscr S$, the random edge set $\ECap(\mathscr S)$ is a stopping
domain: for every finite edge set $A$, the event
$\{\ECap(\mathscr S)=A\}$ is the union of the events $\{\mathscr S=S\}$ over
deterministic $S$ with $\ECap(S)=A$.  This is a countable union, since there are countably many finite box collections.  Put
$\cF_{\mathscr S}=\cF_{\ECap(\mathscr S)}$.  Lemma~\ref{l.freshness} then
gives the Bernoulli-$p$ product law on $\mathcal U(\mathscr S)$ conditionally
on $\cF_{\mathscr S}$.  Thus every edge of $\mathcal U(\mathscr S)$ is
unread relative to this conditioning.

\begin{definition}[Virgin and sealed seeds]
\label{d.virgin-seed}
Let \(\mathscr S\) be a seal.  A \emph{seed} is an
$\cF_{\mathscr S}$-measurable random vertex $s\in V(N)$.  It is \emph{virgin}
if at least one edge of \(\mathcal U(\mathscr S)\) is incident to \(s\), and
\emph{sealed} otherwise.
\end{definition}

In the explorations to be described shortly, open connected sets will be
identified within the layer.  When such a set meets the woven fabric
$\mathcal W(Q)$ of a nearest-neighbour component $Q$ of the seal~$\mathscr S$,
we say that $Q$ is \emph{activated}.  Since $\mathcal W(Q)$ is connected
and its edge states are already part of the seal record, the whole fabric
may then be incorporated into the reached open set without a new product-field
query.  The captured scan follows its as-yet-unprocessed open connections;
any resulting exterior opportunities are handled from virgin seeds.

In use, a supplied seed is first followed through the current capture.  For
the ordinary seal just defined, the captured states belong to
$\cF_{\mathscr S}$; in the relative formulation below, an unread captured
edge is instead queried when it is first needed.  This search may die
out or locate one or several virgin seeds, recording any side vertices
reached along the way.  No exterior edge is queried merely in locating such
a seed.

When a virgin seed~$s$ is selected for growth, order its incident edges
deterministically and let~$e_s$, the \emph{launch edge}, be the first one in
$\mathcal U(\mathscr S)$.  Let~$L_s$ be the unique boxlet to which~$e_s$ is
assigned.  The pair~$(e_s,L_s)$ specifies where the ensuing exterior exploration is
launched.  Lemma~\ref{l.unread-edge-launch-cell}
gives
$E(L_s)\subseteq\mathcal U(\mathscr S)$, and every overlapping box having
\(L_s\) as a boxlet is exterior.  Thus the bounded connected boxlet \(L_s\) is
wholly fresh when the virgin seed is selected for growth, before any of its
remaining edges are queried.

For a vertex set or subgraph \(A\) in~\(N\), write
\[
   \LayerBoxsupp(A)=\Boxsupp(A)\cap\overline N.
\]
Cluster explorations use a fixed deterministic rule to select the next
eligible edge.  An edge whose state is already recorded is processed using
that state, while an unread edge is queried when it is first needed.  The
estimates for these explorations do not depend on the particular processing order.

To grow a virgin seed \(s\) is to reveal its open cluster \(O_s\) in \(N\),
following the fixed states in \(\ECap(\mathscr S)\) and querying only edges in
\(\mathcal U(\mathscr S)\).  We regard this search, including any returns
through the capture and later exits back into unread territory, as a single
\emph{exploration episode}.  If \(O_s\) is finite, the episode ends only after every edge of
\(\mathcal U(\mathscr S)\) incident to \(V(O_s)\) has been queried.  Let
\(\mathcal R_s\) be the set of edges queried in the episode, and let
\[
       \mathcal B_s=\LayerBoxsupp(V(O_s))
\]
be the layer-box support of this episode.  Every queried edge is either open
in \(O_s\) or closed and incident to a vertex of \(O_s\), so the queried
record introduces no additional layer box.  The event
\(\{|E(O_s)|=\infty\}\) is the open-escape event.

On a finite outcome, the current seal is enlarged to capture the newly
queried record while preserving a product exterior.  This has two stages.
First we surround the new record by an enclosing quilt, an auxiliary apt
barrier, found by a stopping exploration, whose capture contains the newly
revealed material.  This quilt need not be a seal, however, and, in the
second \emph{seal completion} step, the mesh boundary is pushed out so as to
form a welcoming collection; the new seal is obtained as the resulting
reduction.
We begin with the path geometry used to recognise a crossing of such a quilt.

\subsection{Home boxes and boundary crossing}
\label{ss.home-box-boundary-crossing}

Recall the layer geometry of Section~\ref{s.overlapping-boxes-boxlets}.
Axial \(K\)-columns cut the tilted layer into boxlets, and neighbouring boxlets
form the overlapping boxes used for aptness and quilts.  For the path arguments
below we define for each vertex a \emph{home box} from this overlapping collection
in which the vertex lies safely away from the transverse boundary.  Thus
leaving the home box produces a \(K/2\)-scale open path inside it which, when
the box is apt, is absorbed into its woven cluster.

For \(w\in V(N)\), define
\[
   \mu_i(w)=\left\lceil\frac{w_i}{K}+\frac12\right\rceil,
   \qquad i\in\intint{d-1}.
\]
The box
\[
       B_k(\mu(w))=\Obox{k,\mu(w)}
\]
is the \emph{home box} of \(w\).  Coordinatewise,
\begin{equation}
\label{e.home-half-offset}
   K\bigl(\mu_i(w)-\tfrac32\bigr)
     <w_i\le
   K\bigl(\mu_i(w)-\tfrac12\bigr).
\end{equation}
The parallelotope \(\Lambda_{k,\mu(w)}\) has transverse walls at
\(K(\mu_i(w)-2)\) and \(K\mu_i(w)\).  Hence \(w\) lies in its middle half,
with clearance at least \(K/2\) from either transverse wall.

For a vertex set \(A\subseteq V(N)\), define its \emph{home-box trace} by
\[
       \mu(A)=\{\mu(w):w\in A\}\subseteq\Z^{d-1};
\]
for a subgraph use its vertex set.
For a set \(C\) in the box-index lattice \(\Z^{d-1}\), or in the
identified mesh lattice, write \(C^{[r]}\) for the set of sites at
\(\ell^\infty\)-distance at most \(r\) from~\(C\).

\begin{lemma}[Properties of the home-box trace]
\label{l.middle-box-basic}
Let \(K\ge2\), and let \(A\) be a subgraph of the layer \(N=N_k\).
\begin{enumerate}
\item If \(A\) is connected, then \(\mu(A)\) is connected in \(\Z^{d-1}\).
\item Suppose \(A\) is connected and \(u\in\mu(A)\).  If \(A\) has a
      vertex outside \(B_k(u)\), then \(A\) contains a path
      in \(B_k(u)\) of Euclidean diameter at least \(K/2\).
\item Every home box belongs to the layer support:
      \[
          \{B_k(u):u\in\mu(A)\}\subseteq\LayerBoxsupp(A),
      \]
      while, after identifying boxes with transverse indices,
      \[
          \LayerBoxsupp(A)\subseteq\mu(A)^{[1]} .
      \]
\end{enumerate}
\end{lemma}

\noindent{\bf Proof: (1).}
One lattice step changes at most one of the first \(d-1\)
coordinates, by one.  The same is therefore true of the map \(\mu\), so the
image of a path is a path after consecutive repetitions are suppressed.

{\bf (2).}
Choose \(w\in V(A)\) with \(\mu(w)=u\) and follow a path in \(A\)
from \(w\) to a vertex outside \(B_k(u)\).  Let \(e\) be the first edge of
this path not belonging to \(B_k(u)\), and let \(z\) be the endpoint of \(e\)
reached just before \(e\) is traversed.  Since \(e\in E(N_k)\), its upper
endpoint has the correct layer height, so its failure to belong to the home
box is caused by a transverse inequality in \eqref{e.layer-box-region}.
Here \(z\) is \(e^\uparrow\) or its neighbour along the violating transverse
coordinate.  Since membership of an edge in \(B_k(u)\) is determined by
the position of its upper endpoint in \(\Lambda_{k,u}\), \(z\) lies on or
beyond the corresponding transverse wall.  By \eqref{e.home-half-offset}, \(w\) has
clearance at least \(K/2\) from that wall, and hence \(|z-w|\ge K/2\).  The path segment from \(w\) to \(z\)
consists entirely of edges of \(B_k(u)\), and gives the required path.

{\bf (3).}
For the first assertion, take any edge \(e\in E(N)\) incident to
\(w\); such an edge exists because \(w\in V(N)\).  Its upper endpoint
\(e^\uparrow\) belongs to \(\mathsf A_k\), and its transverse coordinates
differ from those of \(w\) by at most one.  By~\eqref{e.home-half-offset}
and \(K\ge2\), this unit displacement remains inside the transverse region of
\(B_k(\mu(w))\).  Hence \(e\) is assigned to that box, so the home box
meets \(A\).

Conversely, suppose \(B_k(v)\) meets a vertex \(w\in V(A)\).
By the definition of \(B_k(v)\), \(w\) is either in
\(\Lambda_{k,v}\) or is a lower endpoint of an edge whose upper endpoint
lies there.  Thus its transverse coordinates differ by at most one lattice
step from those of a point satisfying \eqref{e.layer-box-region}.
Comparing these bounds with~\eqref{e.home-half-offset} gives
\(|v_i-\mu_i(w)|\le1\) for every \(i<d\), proving the second inclusion. \qed

The home-box trace turns a crossing of a finite box-index enclosure into a
macroscopic open path on its boundary.

\begin{lemma}[Episode support at an exterior boundary]
\label{l.episode-support-crossing}
Let \(H\subseteq\Z^{d-1}\) be finite and \(*\)-connected, and put
\(I=\partial_{\rm ext}^{\Z^{d-1}}H\).  Suppose that \(O_s\) is not
contained in one overlapping box and that, for some \(w_-\in V(O_s)\),
\(\mu(w_-)\) lies in \(H\) or in a bounded component of
\(\Z^{d-1}\setminus H\).  If
\(\mathcal B_s\) contains a box of
\(\Ext_{\overline N}(I)\), then \(O_s\) contains an open path of
Euclidean diameter at least \(K/2\) in some box of \(I\).
\end{lemma}

\noindent{\bf Proof.} By Lemma~\ref{l.middle-box-basic}(1), the home-box trace \(\mu(O_s)\) is
nearest-neighbour connected.  If it meets~\(I\), choose
\(u\in\mu(O_s)\cap I\).  The home box \(B_k(u)\) belongs to
\(I\), and Lemma~\ref{l.middle-box-basic}(2) gives the required path.

Suppose that the trace does not meet \(I\).  Starting from
\(\mu(w_-)\), nearest-neighbour connectivity then keeps it in \(H\) together
with the bounded components of \(\Z^{d-1}\setminus H\).  Indeed, along any
nearest-neighbour path from this set to the unbounded component, the first site
in the unbounded component has a predecessor in \(H\), and hence belongs to
\(I\).

Let \(B_k(v)\) support a vertex of \(O_s\).  Lemma~\ref{l.middle-box-basic}(3)
places \(v\) within one \(*\)-step of some \(u\in\mu(O_s)\).  If \(u\in H\),
then a \(*\)-neighbour \(v\) in the unbounded component of
\(\Z^{d-1}\setminus H\) belongs to \(I\) by definition.  If instead
\(u\) lies in a bounded complementary component, resolve the \(*\)-step from
\(u\) to such a \(v\) into nearest-neighbour steps inside the corresponding unit
coordinate cube.  Since its endpoints lie in different components of
\(\Z^{d-1}\setminus H\), an intermediate site lies in~\(H\); hence \(v\) is
\(*\)-adjacent to~\(H\) and again belongs to~\(I\).  Thus every layer
box supporting \(O_s\) belongs to
\(H\), to a bounded complementary component of~\(H\), or to~\(I\).
Each of these sets is contained in \(\Capt_{\overline N}(I)\):
any path from \(H\) or one of its bounded complementary
components to infinity must meet~\(I\).  This contradicts the assumed
exterior support.  The trace must therefore
meet \(I\), and the first paragraph applies. \qed

\subsection{Enclosing quilts}
\label{ss.trap-formation}

In the product-field setting of Section~\ref{ss.seals-virgin-seeds}, let
\(\mathscr S\) be the seal already in place and let \(s\) be a virgin seed
whose episode \(O_s\) is finite, with queried record \(\mathcal R_s\) and
layer-box support \(\mathcal B_s\).  We now construct the auxiliary quilt
used in seal enlargement.  Conditionally on \(\cF_{\mathscr S}\), the
boundary-crossing criterion above produces an enclosing quilt around the part
of this support not already captured by \(\mathscr S\).  The search works
around finite inapt obstacles in the exterior, using their apt outer
boundaries to push the barrier outward until the new record is captured.
More precisely, we construct a quilt
$\mathcal T_s\subseteq\Ext_{\overline N}(\mathscr S)\cup\mathscr S$ such that
$\mathcal B_s\subseteq\Capt_{\overline N}(\mathscr S)\cup
\Capt_{\overline N}(\mathcal T_s)$.  Thus every box that meets an edge or
vertex of the explored cluster lies in the capture of \(\mathscr S\) or of
$\mathcal T_s$; the new quilt may use a portion of \(\mathscr S\) as a cap.
Let \(\mathsf{Loc}_s\) be the event that \(O_s\) is contained in one
overlapping box.  On this event \(\mathcal R_s\) has overlapping-box support
of cardinality bounded by a constant depending only on \(d\), and is absorbed
directly by a seal completion.  The enclosing quilt is needed only on
\(\mathsf{Loc}_s^{\rm c}\).

We identify the transverse box-index lattice with \(\Z^{d-1}\) by
\(u\leftrightarrow B_k(u)\).  For a finite \(A\subseteq\Z^{d-1}\), write
\(\partial_{\rm ext}^{\Z^{d-1}}A\) for its exterior \( * \)-boundary: the
sites in the unbounded nearest-neighbour component of
\(\Z^{d-1}\setminus A\) that are \( * \)-adjacent to \(A\).
Thus \( * \)-adjacency enters the definition of the boundary, whereas
Lemma~\ref{l.extbdy-slab} supplies its nearest-neighbour connectivity when
\(A\) is finite and \( * \)-connected.

\begin{lemma}[Subcritical comparison clusters]
\label{l.subcritical-star-components}
Let $\Delta_*=3^{d-1}-1$ be the degree of the $*$-adjacency graph on
$\Z^{d-1}$.  There is $p_\star>0$ with $\Delta_*p_\star<1$ and constants
$\lambda_*,c_*>0$ such that, for Bernoulli site percolation $\xi$ of density
$p_\star$ and every finite $R\subseteq\Z^{d-1}$,
\[
  \EE_{p_\star}\exp\{\lambda_*|\mathcal C_{p_\star}(R)|\}
     \le \exp\{c_*|R|\},
\]
where $\mathcal C_{p_\star}(R)$ is the union of the $*$-components of $\xi$
that meet $R$.
\end{lemma}

\noindent{\bf Proof.}
Explore these components generation by generation along $\xi$-$*$-paths from
$R$.  The exploration is dominated by a Galton--Watson process with $|R|$
initial ancestors and offspring distribution
$\operatorname{Bin}(\Delta_*,p_\star)$.  Choosing
$\Delta_*p_\star<1$, the total progeny has an exponential moment, giving the
claim. \qed

\begin{lemma}[Finite inapt components]
\label{l.finite-inapt-components}
For $K$ sufficiently large, conditionally on $\cF_{\mathscr S}$, every
inapt-box $*$-component contained in $\Ext_{\overline N}(\mathscr S)$ is
almost surely finite, uniformly in the realised seal.  Moreover, there are
$C,c>0$ such that, uniformly in the realised seal and in the choice of an
exterior box~$B$,
\[
 \PP\bigl(|\mathcal C_{\rm in}(B)|\ge m\mid\cF_{\mathscr S}\bigr)
 \le Ce^{-cm},\qquad m\ge1,
\]
where $\mathcal C_{\rm in}(B)$ is the inapt $*$-component of $B$, and is
empty when $B$ is apt.
\end{lemma}

\noindent{\bf Proof.} Decompose over the deterministic values $S$ of the seal.  On
$\{\mathscr S=S\}$, the exterior edge states have their product law by
Lemma~\ref{l.freshness}, and the deterministic exterior-box estimate applies.
Conditionally on $\cF_{\mathscr S}$, the inapt-box field is one-dependent:
aptness of a box is read from the edges in its boxlets, and boxes at
$\ell^\infty$-distance at least two have disjoint edge sets.  The one-box
inaptness probability is at most an arbitrarily small $\rho$ once $K$ is
large, by Lemma~\ref{l.low-density-unusable-exterior-boxes}.  By Liggett--Schonmann--Stacey~\cite{LSS}, the inapt-box field is
dominated, uniformly in $S$, by the Bernoulli field $\xi$ of
Lemma~\ref{l.subcritical-star-components} once $K$ is large enough.  The
singleton case of that lemma gives the displayed estimate and, in particular,
almost-sure finiteness. \qed

\begin{lemma}[The enclosing search stays outside the old capture]
\label{l.enclosing-search-exterior}
Let \(S\) be deterministic, satisfy \eqref{e.seal-mesh-geometry}, and be welcoming.
\begin{enumerate}
\item No box of $\Capt_{\overline N}(S)\setminus S$ is $*$-adjacent to a box
of $\Ext_{\overline N}(S)$.
\item Let $\mathcal V$ be a finite $*$-component of apt boxes in
$\Ext_{\overline N}(S)\cup S$.  Then
\[
   \partial_{\rm ext}^{\overline N}\mathcal V
       \subseteq\Ext_{\overline N}(S).
\]
If $H$ is the inapt-box $*$-component containing this boundary, then
$H\subseteq\Ext_{\overline N}(S)$, while
\[
   \partial_{\rm ext}^{\overline N}H
       \subseteq\Ext_{\overline N}(S)\cup S.
\]
\end{enumerate}
\end{lemma}

\noindent{\bf Proof: (1).}
Put $A=S\cap\Mesh$.  If
$B\in\Capt_{\overline N}(S)\setminus S$, then $B\notin S=A^{\mathrm{ov}}$,
so $\MeshSupp(B)$ has a member $M\notin A$; by
\eqref{e.mesh-transfer-full-exterior}, no such member lies in $U^*(A)$.
Thus $M$ belongs to a bounded $*$-component of $\Mesh\setminus A$.
If $C\in\Ext_{\overline N}(S)$, the same formula gives
$M'\in\MeshSupp(C)\cap U^*(A)$.  When $B\sim_*C$, the coordinate
form of mesh supports gives $M\sim_*M'$ (or $M=M'$), contradicting the
choice of the two $*$-components.

{\bf (2).}
Let $I=\partial_{\rm ext}^{\overline N}\mathcal V$.  Suppose
$B\in I\cap\Capt_{\overline N}(S)$.  The box $B$ is not in $S$: otherwise,
being apt and $*$-adjacent to $\mathcal V$, it would belong to the same apt
$*$-component.  Hence $B$ lies strictly inside the capture.  Choose
$C\in\mathcal V$ with $C\sim_*B$.  Part~(1) excludes
$C\in\Ext_{\overline N}(S)$, so $C\in S$.

Take a nearest-neighbour path from $B$ to infinity which avoids
$\mathcal V$, and let $C'$ be its first contact with $S$ before it reaches
the old exterior; such a contact exists by part~(1).  Resolve the single
$*$-step from $C$ to $B$ into nearest-neighbour steps and concatenate with
the path from $B$ to $C'$.  Part~(1) keeps every non-seal box of this finite
path within the capture.  Breaking at its visits to $S$ and applying
Lemma~\ref{l.mesh-boundary-screening} shows that $C$ and $C'$ lie in the
same nearest-neighbour component of $S$.  This component is apt and hence
contained in $\mathcal V$, contradicting that the chosen path avoids
$\mathcal V$.  Therefore $I\subseteq\Ext_{\overline N}(S)$.

Every box of $I$ is inapt by maximality of $\mathcal V$.  Its inapt
$*$-component $H$ cannot cross the apt set $S$, and part~(1) forbids a
$*$-step from the old exterior directly to a box strictly inside the capture.
Thus $H\subseteq\Ext_{\overline N}(S)$.  The same argument applied to a
box of $\partial_{\rm ext}^{\overline N}H$ gives
$\partial_{\rm ext}^{\overline N}H\subseteq
\Ext_{\overline N}(S)\cup S$. \qed

\begin{proposition}[Enclosing quilt on non-escape]
\label{p.finite-trap-on-nonescape}
Almost surely on the event
$\{s\notconn\infty\ \text{in}\ N\}\cap\mathsf{Loc}_s^{\rm c}$, there is a
quilt $\mathcal T_s\subseteq\Ext_{\overline N}(\mathscr S)\cup\mathscr S$
whose capture, together with that of \(\mathscr S\), contains
\(\mathcal B_s\).
\end{proposition}

\noindent{\bf Proof.} By Lemma~\ref{l.finite-inapt-components}, conditionally on $\cF_{\mathscr S}$, every inapt-box $*$-component is almost surely finite.  Hence this holds $\PP$-almost surely; work on this event.
The enclosing quilt is built in successive rounds $j=0,1,\ldots$, until the stopping condition below is met, adding at each round
an apt region $\cV_j$ and a finite inapt-box $*$-component $H_j$.  The quilt after round $j$ is \(\mathcal T_j=\partial_{\rm ext}^{\Z^{d-1}}H_j\), the exterior boundary of \(H_j\) in the box-index lattice.  Since $H_j$ is finite and
$*$-connected, Lemma~\ref{l.extbdy-slab} makes $\mathcal T_j$
nearest-neighbour connected.  Its boxes are apt, so it is a quilt.

The apt $*$-components used below are taken in
$\Ext_{\overline N}(\mathscr S)\cup\mathscr S$.  Lemma
\ref{l.enclosing-search-exterior} will place the inapt component $H_j$ in the
old exterior and its apt boundary $\mathcal T_j$ in
$\Ext_{\overline N}(\mathscr S)\cup\mathscr S$ at every round.

Round zero splits according to whether some apt
$B\in\Ext_{\overline N}(\mathscr S)$ contains an open subpath of $O_s$ of
diameter at least $K/2$.  To decide the case, test in a fixed order the
finitely many exterior boxes containing such a subpath, querying their
unread edges in a fixed order.  These candidates are identified from $O_s$
and belong to $\mathcal B_s$.

\textit{Round zero when such a box exists.}
Let $B_0$ be the first apt box among these candidates, and let $z_0$ be
its transverse index.  Fix one such subpath in $B_0$,
using a deterministic tie-break.  By aptness, this subpath is contained in
$\mathcal W(B_0)$.  Let $\cV_0$ be the apt-box $*$-component of $z_0$ in
$\Ext_{\overline N}(\mathscr S)\cup\mathscr S$.  Were $\cV_0$ infinite, choose an infinite self-avoiding $*$-path in it
starting at $B_0$.  Consecutive boxes on this path have intersecting
woven clusters by Lemma~\ref{l.woven-fabric-note}(1).  This open
subpath in $B_0$ is absorbed into its woven cluster, so chaining along
the $*$-path would connect $O_s$ to infinity, contrary to the event in force.
Thus $\cV_0$ is finite.  Set
$I_0 = \partial_{\rm ext}^{\Z^{d-1}}\cV_0$.  Lemma~\ref{l.extbdy-slab}
makes $I_0$ nearest-neighbour connected, and
Lemma~\ref{l.enclosing-search-exterior} places $I_0$ in the old exterior.
By maximality of $\cV_0$ as an apt $*$-component in
$\Ext_{\overline N}(\mathscr S)\cup\mathscr S$, every box of $I_0$ is inapt.
Let $H_0$ be the inapt-box $*$-component containing $I_0$; the same lemma
places $H_0$ in the old exterior.

\textit{Round zero when no such box exists.}
On \(\mathsf{Loc}_s^{\rm c}\),
$O_s$ leaves every overlapping box.  Thus, for each \(u\in\mu(O_s)\), it has a vertex \(w\) with
\(\mu(w)=u\), lying in the middle half of \(B_k(u)\), and it has a vertex
outside \(B_k(u)\).
Lemma~\ref{l.middle-box-basic}(2) yields an open path of diameter at least $K/2$
inside $B_k(u)$.  Take $z_0=\mu(s)$.  The edge $e_s$ is incident to $s$, so the proof of
Lemma~\ref{l.middle-box-basic}(3) gives $e_s\in E(B_k(z_0))$; equivalently,
$L_s$ is a boxlet of this home box.
Lemma~\ref{l.unread-edge-launch-cell} therefore makes $B_k(z_0)$ exterior to
$\mathscr S$.  The preceding paragraph supplies a diameter-$K/2$ open path
there, so the present case makes $B_k(z_0)$ inapt.  Let $H_0$ be its
inapt-box $*$-component.  Lemma~\ref{l.enclosing-search-exterior}(1), together
with aptness of the old seal, keeps this whole component in the old exterior.

\textit{Closing round zero.}
$H_0$ is finite on the event in force.  By Lemma~\ref{l.extbdy-slab} in
$\Z^{d-1}$, the boundary $\mathcal T_0 = \partial_{\rm ext}^{\Z^{d-1}} H_0$ is
nearest-neighbour connected.  No box of $\mathcal T_0$ belongs to the inapt
component $H_0$, so every such box is apt and $\mathcal T_0$ is a quilt.
Lemma~\ref{l.enclosing-search-exterior}(1) gives
$\mathcal T_0\subseteq\Ext_{\overline N}(\mathscr S)\cup\mathscr S$.
In the first round-zero case, $I_0\subseteq H_0$ separates the finite
apt component $\cV_0$ from infinity.  Let $w_0$ be the first vertex of the chosen
open subpath in $B_0$ under a fixed vertex order.  Lemma~\ref{l.middle-box-basic}(3) places its
home-box index within one $*$-step of the index of $B_0$.  It therefore
either lies on $\mathcal T_0$ or lies in $H_0$ or in a bounded component of
$\Z^{d-1}\setminus H_0$: a $*$-step from the bounded side to the unbounded
component, if resolved into nearest-neighbour steps, meets $H_0$, and its
unbounded endpoint then belongs to $\mathcal T_0$.  In the second case take
$w_0=s$, for which $\mu(w_0)=z_0\in H_0$.  We shall use the following
alternative at each round:
\begin{equation}
\label{e.enclosing-interior-anchor}
 \mu(w_j)\in\mathcal T_j
 \quad\hbox{or}\quad
 \mu(w_j)\text{ lies in }H_j\text{ or a bounded component of }
       \Z^{d-1}\setminus H_j.
\end{equation}

\textit{Iterating.}
Write $\Capt_{\overline N}(\mathcal T_j)$ for the running capture and test
whether
\[
   \mathcal B_s\,\subseteq\,
   \Capt_{\overline N}(\mathscr S)\cup\Capt_{\overline N}(\mathcal T_j)\, .
\tag{$\star_j$}
\]
If~($\star_j$) holds, let $J=j$ and end the construction.  Otherwise,
$\mathcal B_s$ contains a box in
$\Ext_{\overline N}(\mathcal T_j)$.  Use the
alternative~\eqref{e.enclosing-interior-anchor}.  If
$\mu(w_j)\in\mathcal T_j$, Lemma~\ref{l.middle-box-basic}(2) gives an open
path of diameter at least $K/2$ in the home box $B_k(\mu(w_j))$, which is a
box of~$\mathcal T_j$.  If $\mu(w_j)\notin\mathcal T_j$,
Lemma~\ref{l.episode-support-crossing}, applied with $H=H_j$, gives such a
path in a box of~$\mathcal T_j$.  In either case let $z'_j\in\mathcal T_j$
be the first box in the fixed order carrying one of these paths, and
fix one such path there by a deterministic tie-break.  Both choices
are determined by $O_s$ and the current boundary $\mathcal T_j$.
Repeating the round-zero construction from $z'_j$
gives the apt-box $*$-component
$\cV_{j+1}$ containing $z'_j$.  Were it infinite, an infinite self-avoiding $*$-path from $z'_j$ and
Lemma~\ref{l.woven-fabric-note}(1), applied successively along that path,
would connect the open path just identified in $O_s$ to infinity.  Thus
$\cV_{j+1}$ is finite.
Because the boxes of the quilt
$\mathcal T_j$ are apt and nearest-neighbour connected,
$\mathcal T_j\subseteq\cV_{j+1}$.  Put
$I_{j+1}=\partial_{\rm ext}^{\Z^{d-1}}\cV_{j+1}$.
Lemma~\ref{l.extbdy-slab} makes $I_{j+1}$ nearest-neighbour connected, and
Lemma~\ref{l.enclosing-search-exterior} places $I_{j+1}$ in the old exterior.
Maximality of $\cV_{j+1}$ as an apt $*$-component in
$\Ext_{\overline N}(\mathscr S)\cup\mathscr S$ then makes every box of
$I_{j+1}$ inapt.  Let $H_{j+1}$ be the inapt-box $*$-component containing
$I_{j+1}$; Lemma~\ref{l.enclosing-search-exterior} places it in the old
exterior.  Every path from $\Capt_{\overline N}(\mathcal T_j)$ to infinity meets
$\mathcal T_j\subseteq\cV_{j+1}$.  After its last visit to the finite set
$\cV_{j+1}$, the path enters its exterior boundary $I_{j+1}$.  Thus
$I_{j+1}\subseteq H_{j+1}$ separates the whole old capture from
infinity.  For
$\mathcal T_{j+1}=\partial_{\rm ext}^{\Z^{d-1}}H_{j+1}$ this gives the stronger
nesting
\[
 \Capt_{\overline N}(\mathcal T_j)
   \subseteq
 \Capt_{\overline N}(\mathcal T_{j+1})\setminus\mathcal T_{j+1}.
\]
In particular $z'_j\in\mathcal T_j\subseteq\cV_{j+1}$ lies on the
bounded side of $I_{j+1}\subseteq H_{j+1}$.  Let $w_{j+1}$ be the first vertex of
the chosen subpath in $B_k(z'_j)$ under the fixed vertex order.  As in the closing of round zero,
its home-box index either belongs to $\mathcal T_{j+1}$ or lies in $H_{j+1}$
or in a bounded component of $\Z^{d-1}\setminus H_{j+1}$.  Thus
\eqref{e.enclosing-interior-anchor} is available at the next round.

\textit{Termination.}
Whenever~($\star_j$) fails, the selected box $z'_j$ contains an open subpath
of $O_s$ and hence belongs to the finite set $\mathcal B_s$.  The next
boundary places the whole preceding capture strictly inside the next
one.  Hence $z'_j$ remains strictly captured on every later round and cannot
be selected again.  Thus~($\star_j$) can fail for at most
$|\mathcal B_s|$ indices $j$, and so holds at some finite least index $J$;
define $\mathcal T_s=\mathcal T_J$.  The components
$H_0,\ldots,H_J$ are successively contained in the captures of the
later boundaries.  With $H_s:=H_J$ we have
$\mathcal T_s=\partial_{\rm ext}^{\Z^{d-1}}H_s$, and~($\star_J$) is the asserted
inclusion.  Lemma~\ref{l.enclosing-search-exterior}(1) gives
$\mathcal T_s\subseteq\Ext_{\overline N}(\mathscr S)\cup\mathscr S$.
\qed

For a finite quilt collection $\mathscr Q$, put
\begin{equation}
\label{e.mesh-region-enclosing-quilt}
 R(\mathscr Q)
 =\{M\in\Mesh:E(M)\cap E(B)\ne\emptyset
       \text{ for some }B\in\Capt_{\overline N}(\mathscr Q)\}.
\end{equation}
Thus \(R(\mathscr Q)\) is the collection of mesh boxes which meet the edge set
of an overlapping box captured by \(\mathscr Q\).  This notation is used only
for the enclosing quilt, which need not itself have the disjoint-mesh form of
a seal.

\begin{lemma}[Stopped enclosing record]
\label{l.stopped-enclosing-record}
On \(\{O_s\text{ is finite}\}\cap\mathsf{Loc}_s^{\rm c}\), let \(\mathcal G_s\) be the
\(\sigma\)-field generated by \(\cF_{\mathscr S}\) and by the identities and
states of the edges queried in the cluster exploration and in the construction
of \(\mathcal T_s\).  Put \(\mathcal M_s=R(\mathcal T_s)\).  Then
$\mathcal M_s$, $\mathcal T_s$, and $\mathcal R_s$ together with the states
of its edges are $\mathcal G_s$-measurable, every
edge newly queried after conditioning on $\cF_{\mathscr S}$ belongs to
\[
   \bigcup_{M\in \mathcal M_s}E(M),
\]
and, conditionally on $\mathcal G_s$, the remaining edge states have their
product Bernoulli law.
\end{lemma}

\noindent{\bf Proof.} The finite-cluster exploration is a stopping
exploration, and the preliminary round-zero aptness tests concern only boxes
of $\mathcal B_s$.  Subsequent searches determine the apt $*$-components
$\cV_j$ and inapt $*$-components $H_j$ by examining their boxes and
$*$-neighbours in fixed order.  A $*$-neighbour outside the component
belongs either to its exterior boundary or to a bounded component of
its complement.  The exterior boundary of $\cV_j$ is
$I_j\subseteq H_j$, while that of $H_j$ is $\mathcal T_j$.
Every nearest-neighbour path to infinity from either component or any
of its bounded complementary components meets its exterior boundary.
Thus the boxes examined in these searches are captured by
$\mathcal T_j$, except for any already in the old capture.
The successive captures are nested, and the terminal inclusion
~($\star_J$) places $\mathcal B_s$ in the union of the old and terminal
captures.  Since $\mathcal B_s$ contains the round-zero candidate boxes
and every box meeting an edge of $\mathcal R_s$, every tested box and
every box meeting an edge of $\mathcal R_s$ belongs to
$\Capt_{\overline N}(\mathscr S)\cup\Capt_{\overline N}(\mathcal T_s)$.

Edges in the old capture are already fixed under $\cF_{\mathscr S}$ and are
not newly queried.  Every newly queried edge therefore lies in an overlapping
box captured by $\mathcal T_s$, so its unique mesh box belongs to $\mathcal M_s$ by
\eqref{e.mesh-region-enclosing-quilt}.  At each stage, the next box to test and the decision to continue or
terminate are determined by the record then available; the path and
anchor choices are fixed as above.  The cluster episode and subsequent
enclosing search therefore form a stopping exploration of the product
exterior.  Lemma~\ref{l.freshness} gives the conditional Bernoulli product
law for the unqueried edges. \qed

\begin{lemma}[Isoperimetry for a separating box set]
\label{l.column-isoperimetry}
Let $d\geq3$.  Let $H\subseteq\Z^{d-1}$ be finite and let $G$ be the union of the bounded
components of $\Z^{d-1}\setminus H$.  There is $c_d>0$ such that
\[
   |H|\, \geq \, c_d|G|^{(d-2)/(d-1)}\, .
\]
\end{lemma}

\noindent{\bf Proof.} For each coordinate direction \(i\in\intint{d-1}\), every
axis-parallel line that meets $G$ also meets $H$: starting from a point of a
bounded complementary component and moving to infinity along the line, one
must cross the separating set.  Thus the coordinate projection $\pi_i(G)$
has cardinality at most $|H|$.  The discrete Loomis--Whitney inequality
\cite{LoomisWhitney} gives
\[
 |G|^{d-2}
 \, \leq \,
 \prod_{i=1}^{d-1}|\pi_i(G)|
 \, \leq \,
 |H|^{d-1}\, .
\]
Rearrangement proves the claim. \qed

\begin{lemma}[Captured volume and the separating set]
\label{l.trap-charged-to-inapt-hulls}
Let $d\ge3$.  There is \(a=a(d)>0\) such that, if \(H\subseteq\Z^{d-1}\) is finite and
\(*\)-connected and
\(\mathcal T=\partial_{\rm ext}^{\Z^{d-1}}H\), then
\[
 |H|\, \geq \, a\,|\Capt_{\overline N}(\mathcal T)|^{(d-2)/(d-1)}\, .
\]
\end{lemma}

\noindent{\bf Proof.} Put $\alpha=(d-2)/(d-1)$.  Each site of \(\mathcal T\) is \(*\)-adjacent to \(H\), and a site of \(H\) has at most \(3^{d-1}-1\) such neighbours.  Hence \(|\mathcal T|\leq(3^{d-1}-1)|H|\).
By definition,
\(\Capt_{\overline N}(\mathcal T)\setminus\mathcal T\) is the union of the
bounded components of \(\Z^{d-1}\setminus\mathcal T\).  Applying Lemma~\ref{l.column-isoperimetry} to \(\mathcal T\) gives \(|\Capt_{\overline N}(\mathcal T)\setminus\mathcal T|\leq C_d|\mathcal T|^{(d-1)/(d-2)}\).  Consequently \(|\Capt_{\overline N}(\mathcal T)|\leq|\mathcal T|+C_d|\mathcal T|^{(d-1)/(d-2)}\leq C'_d(|H|+|H|^{(d-1)/(d-2)})\).
If \(|H|=0\), the assertion is trivial.  Otherwise
\(1/\alpha=(d-1)/(d-2)\ge1\), so the last display is at most
\(2C'_d|H|^{1/\alpha}\).  Rearranging gives the asserted lower bound. \qed

\begin{corollary}[Enclosing-quilt cardinality]
\label{c.trap-cardinality-tail}
There are \(K_0\in\N\) and \(C,c\in(0,\infty)\) such that, for all
\(K\ge K_0\), conditionally on \(\cF_{\mathscr S}\),
\[
 \PP\bigl(
    |\mathcal T_s|\, \ge \, m,\ O_s\text{ finite},\ \mathsf{Loc}_s^{\rm c}
    \mid\cF_{\mathscr S}
 \bigr)\, \le \, Ce^{-cm},\qquad m\, \ge \, 1\, .
\]
\end{corollary}

\noindent{\bf Proof.} On the event in force,
\(\mathcal T_s=\partial_{\rm ext}^{\Z^{d-1}}H_s\), where \(H_s\) is one finite
inapt \(*\)-component, and
\(|\mathcal T_s|\le(3^{d-1}-1)|H_s|\).
Condition additionally on any auxiliary randomness used to select the virgin
seed; the estimate below is uniform in this choice, so we may average over it
at the end.  Recall that \(L_s\) is the boxlet assigned to \(e_s\). Choose, by the same
fixed rule as in Proposition~\ref{p.uniform-open-escape-exterior-seed}, an
overlapping box \(B_0\) having \(L_s\) as a boxlet.  Let \(b_0\) be its
transverse index.  Lemma~\ref{l.unread-edge-launch-cell} makes \(B_0\)
exterior to \(\mathscr S\), while \(s\in V(L_s)\subseteq V(B_0)\) gives
\(B_0\in\mathcal B_s\).  The terminal inclusion in
Proposition~\ref{p.finite-trap-on-nonescape} therefore places \(B_0\) in
\(\Capt_{\overline N}(\mathcal T_s)\).

Since $b_0\in\Capt_{\overline N}(\mathcal T_s)$, either
$b_0\in\mathcal T_s$ or $b_0$ lies in a bounded component of
$\Z^{d-1}\setminus\mathcal T_s$.  In the first case
$\operatorname{dist}_\infty(b_0,H_s)\le1$.  In the second case each of two
opposite coordinate rays from $b_0$ meets the connected set
$\mathcal T_s$.  If their first-hit distances are $r_+$ and $r_-$, any path
in $\mathcal T_s$ joining the two first-hit sites has at least
$r_++r_-+1$ vertices.  Hence
$\min\{r_+,r_-\}\le |\mathcal T_s|$, and, since every site of
$\mathcal T_s$ is $*$-adjacent to~$H_s$,
\[
  \operatorname{dist}_\infty(b_0,H_s)
     \le |\mathcal T_s|+1\le C|H_s|.
\]
Thus, on $\{|H_s|=n\}$, the component $H_s$ contains a site among at most
$Cn^{d-1}$ indices at $\ell^\infty$-distance at most $Cn$ from $b_0$.
Lemma~\ref{l.finite-inapt-components} and a union bound therefore give,
uniformly in the auxiliary choice,
\[
 \PP\bigl(|H_s|\ge m,\ O_s\text{ finite},\ \mathsf{Loc}_s^{\rm c}
       \mid\cF_{\mathscr S}\bigr)
 \le \sum_{n\ge m} Cn^{d-1}e^{-cn}
 \le C'e^{-c'm}.
\]
Since $|\mathcal T_s|\le(3^{d-1}-1)|H_s|$, this proves the assertion after
changing constants.  Averaging over the auxiliary choice completes the proof.
\qed

\begin{corollary}[Enclosing-region tail]
\label{c.trap-tail}
There exist $K_0\in\N$ and $C,c\in(0,\infty)$ such that, for all
$K\geq K_0$, conditionally on $\cF_{\mathscr S}$,
\[
   \PP\!\left(
      \bigl|\Capt_{\overline N}(\mathcal T_s)\bigr|\, \geq \, r,\;
      s\notconn\infty\ \text{in}\ N,\quad \mathsf{Loc}_s^{\rm c}
      \,\middle|\, \cF_{\mathscr S}
   \right)
\, \leq \,
   C\,e^{-c r^{(d-2)/(d-1)}}\, ,
   \qquad r\, \geq \, 1\, .
\]
\end{corollary}

\noindent{\bf Proof.} Write \(\mathcal T_s=\partial_{\rm ext}^{\Z^{d-1}}H_s\).  By Lemma~\ref{l.trap-charged-to-inapt-hulls}, \(|\Capt_{\overline N}(\mathcal T_s)|\ge r\) implies \(|H_s|\ge a r^{(d-2)/(d-1)}\).  The anchored component estimate proved in
Corollary~\ref{c.trap-cardinality-tail} gives an exponential tail for
\(|H_s|\).  Substitution yields the result. \qed

\subsection{Seal completion}
\label{ss.trap-cluster-to-seal}

The enclosing quilt $\mathcal T_s$ has been constructed relative to the
seal $\mathscr S$ already in place.  Together, their captures contain the
newly revealed finite record, but we have not yet obtained a new seal.  To
enlarge $\mathscr S$ so that it also captures this record, we combine its
mesh interior with the mesh region $R(\mathcal T_s)$ supplied by the quilt.
The operation that turns this union into a seal is \emph{seal completion}.
More generally, starting from any finite mesh region $D$, we form its mesh boundary
$A(D)=\partial_{\rm ext}^{\Mesh}D$ and the associated overlapping-box
boundary $A(D)^{\mathrm{ov}}$, and check aptness at its boundary boxes and
the welcoming inequality~\eqref{e.welcoming-halo} at its halo boxes.
If either condition fails at a box~$B$, we absorb a fixed mesh
neighbourhood of $B$ into $D$ and recompute the boundary.  Thus a defect
is repaired locally, but the repair may reveal another nearby defect,
and the boundary moves outward only as far as
this chain of failures requires.  The process stops once the candidate is
welcoming.  Every edge inspected along the
way is swallowed by the terminal capture, leaving the exterior unread.  After
the final reduction, the resulting mesh-boundary geometry and welcoming
condition therefore give a seal.

The local defect event $\mathsf E_B$ was defined in
\eqref{e.absolute-seal-bad-box}.  By \eqref{e.failed-seal-test-density}, its
probability is at most $C_d\delta_K$, and the field
$(\mathbf 1_{\mathsf E_B})$ is one-dependent.  Define the corresponding
$3$-thickened defect mesh set by
\[
   \mathcal B_K
   =\bigcup_{\{B:\,\mathsf E_B\}}
       \MeshSupp(B)^{[3]}\subseteq\Mesh.
\]
The thickening makes
$(\mathbf 1_{\{M\in\mathcal B_K\}})_{M\in\Mesh}$ $7$-dependent, and
\[
   \sup_{M\in\Mesh}\PP(M\in\mathcal B_K)
      \le C'_d\delta_K.
\]
Let $\xi$ and $p_\star$ be as in
Lemma~\ref{l.subcritical-star-components}.  By
Liggett--Schonmann--Stacey~\cite{LSS}, there is $\delta_d>0$ such that
$\mathcal B_K$ is dominated by $\xi$ whenever
$C'_d\delta_K\le\delta_d$.  We henceforth take $K\ge K_0$ so that this
inequality holds.  The rooted exponential-moment bound of
Lemma~\ref{l.subcritical-star-components} then holds uniformly for
$\mathcal B_K$.

For a finite mesh region $D$, the overlapping-box boundary
$T=A(D)^{\mathrm{ov}}$ is a candidate seal boundary.  We call the aptness
condition at a box $C\in T$ a \emph{seal test}, and the welcoming
inequality~\eqref{e.welcoming-halo} at $C\in\Halo(T)$ a \emph{halo test}.
A seal test is specified by its box, and a halo test by its box and
captured-side conditioning set $\mathcal I_C(T)$.  When a seal test at
$C$ is evaluated, every edge of $E(C)$ is queried unless its state is
already known; a halo test queries only its specified captured-side edges.
The specification identifies a test independently of whether or when it is
evaluated.  The collection of required tests, however, depends on the
candidate boundary.

\begin{lemma}[Locality of a boundary update]
\label{l.frontier-addition-range}
Let $D\subseteq\Mesh$ be finite, put $A=A(D)$ and $S=A^{\mathrm{ov}}$, and suppose
$B\in S\cup\Halo(S)$.  Then
\[
        \MeshSupp(B)\subseteq A^{[1]}.
\]
Put $H_B=\MeshSupp(B)^{[1]}$ and $D'=D\cup H_B$.  Then
\[
   A(D')\setminus A(D)\subseteq \MeshSupp(B)^{[2]}.
\]
\end{lemma}

\noindent{\bf Proof.}
If $B\in S$, then $\MeshSupp(B)\subseteq A$.  If $B\in\Halo(S)$,
its mesh support contains a member $M_0$ of $\MCap(S)$.
Lemma~\ref{l.mesh-boundary-interface} gives $M_0\in A$.
Every other member of $\MeshSupp(B)$ is at mesh $\ell^\infty$-distance
at most one from $M_0$.  This proves the first inclusion.

Let $M\in A(D')\setminus A(D)$.  Any unbounded path witnessing
$M\in A(D')$ also avoids $D$, so $M$ lies in the unbounded complement of
$D$.  Since $M$ is $*$-adjacent to $D'$ but does not belong to $A(D)$, one
of its $*$-neighbours in $D'$ lies in $D'\setminus D\subseteq H_B$.
Hence $M\in H_B^{[1]}=\MeshSupp(B)^{[2]}$, as required.
\qed

\begin{lemma}[Monotonicity of the mesh exterior]
\label{l.mesh-exterior-monotone}
If $D\subseteq D'\subseteq\Mesh$ are finite, then
\[
   U^*(A(D'))\subseteq U^*(A(D)),
   \qquad
   \MCap(A(D)^{\mathrm{ov}})\subseteq\MCap(A(D')^{\mathrm{ov}}).
\]
\end{lemma}

\noindent{\bf Proof.}
Write $A=A(D)$, $A'=A(D')$, and $S=A^{\mathrm{ov}}$.  First note that, for every finite mesh set $E$, each site of
$U^*(A(E))$ lies in the unbounded nearest-neighbour component of
$\Mesh\setminus E$.  Indeed, resolve an infinite $*$-path in
$\Mesh\setminus A(E)$ into nearest-neighbour steps.  Each inserted site either belongs to $A(E)$, hence not to $E$, or is
$*$-adjacent to a site of the original path and therefore belongs to
$U^*(A(E))$.  In the latter case it lies outside $\MeshInt(A(E)^{\mathrm{ov}})$, while
Lemma~\ref{l.mesh-boundary-capture}(1) gives $E\subseteq \MeshInt(A(E)^{\mathrm{ov}})$.
Thus the resolved path avoids $E$.

Now let $x\in U^*(A')$.  The preceding observation puts $x$ in the
unbounded nearest-neighbour component of $\Mesh\setminus D'$, hence also in
that of $\Mesh\setminus D$.  If $x\in A$, then
$x$ does not belong to $D'$ and, being $*$-adjacent to $D\subseteq D'$, it belongs to
$A'$, a contradiction.  If $x\in\MeshInt(S)$, an infinite $*$-path from $x$ to
infinity in $\Mesh\setminus A'$ must meet $A$ before reaching $U^*(A)$.
Its first such site lies on the same infinite path and hence belongs to
$U^*(A')$, contradicting the case just excluded.  Thus $x\in U^*(A)$.
The capture inclusion follows from~\eqref{e.mesh-transfer-capture}. \qed

When a mesh region is enlarged, its candidate boundary and halo may
change.  A test required for the new boundary is called \emph{newly required}
relative to the old boundary if no test with the same specification was
required before the update.  Thus a seal test is newly required if its
box is new to the boundary, while a halo test is newly required if its box
is new to the halo or its captured-side conditioning set differs.  In
particular, a box may remain in the halo while further boxlets of it
enter the capture.  Its captured-side conditioning set then grows, and
the test with this enlarged conditioning set is newly required.

\begin{lemma}[Locality of newly required tests]
\label{l.changed-seal-tests-local}
Let $D\subseteq\Mesh$ be finite, put $S=A(D)^{\mathrm{ov}}$, and let
$B\in S\cup\Halo(S)$.  Set $D'=D\cup H_B$ and $S'=A(D')^{\mathrm{ov}}$.
If a test for $S'$ at a box~$C$ is newly required relative to~$S$, then
\[
   \MeshSupp(C)\cap\MeshSupp(B)^{[2]}\ne\varnothing.
\]
In particular, at most $C_d$ tests are newly required relative to~$S$.
\end{lemma}

\noindent{\bf Proof.}
Write $A=A(D)$ and $A'=A(D')$.  If $C\in S'\setminus S$, then
$\MeshSupp(C)\subseteq A'$ but $\MeshSupp(C)\not\subseteq A$; hence some
$M\in\MeshSupp(C)$ belongs to $A'\setminus A$.  Lemma
\ref{l.frontier-addition-range} gives
$M\in\MeshSupp(B)^{[2]}$.

Suppose now that the newly required test at $C$ is a halo test.
Lemma~\ref{l.mesh-exterior-monotone} shows that $C$, exterior for $S'$, was
exterior for $S$, and that its captured mesh support can only grow.  A newly required halo test therefore
supplies a mesh box $M\in\MeshSupp(C)$ which is captured after the update but not before it.
Lemma~\ref{l.mesh-boundary-interface}, applied to $S'$, gives $M\in A'$;
previously $M$ was not captured, so $M\notin A$.  Lemma
\ref{l.frontier-addition-range} again gives
$M\in\MeshSupp(B)^{[2]}$.  Bounded local multiplicity proves the last
assertion. \qed

\begin{lemma}[Boundary growth under an addition]
\label{l.boundary-growth-under-addition}
For finite $D,H\subseteq\Mesh$,
\[
 \partial_{\rm ext}^{\Mesh}(D\cup H)
 \subseteq
 \partial_{\rm ext}^{\Mesh}D\ \cup\ H^{[1]},
\]
and hence
\[
 |A(D\cup H)|\le |A(D)|+C_d|H|,
 \qquad
 |A(D\cup H)^{\mathrm{ov}}|\le |A(D)^{\mathrm{ov}}|+C_d|H|.
\]
The same assertions hold if bounded complementary components are filled
before the boundary is taken.
\end{lemma}

\noindent{\bf Proof.}
Let $M\in A(D\cup H)\setminus A(D)$.  Since $M$ lies in the
unbounded complement of~$D$ but not on its exterior boundary, it
has a $*$-neighbour in $D\cup H$ and none in $D$.  Thus a neighbour
lies in $H$, and $M\in H^{[1]}$.
Filling bounded complementary components does not change the exterior
boundary.  The first cardinality bound follows from
$|H^{[1]}|\le 3^{d-1}|H|$.  Each box of
$A(D\cup H)^{\mathrm{ov}}\setminus A(D)^{\mathrm{ov}}$ has a
mesh-support member in $H^{[1]}$, and a fixed mesh box belongs to the
support of at most $3^{d-1}$ overlapping boxes.  This gives the second
bound.
\qed

The lemmas show that adjoining $H_B$ can create mesh-boundary sites only
near~$B$, and that any newly required test is likewise based at a nearby box.
We now specify the completion search.  Start from a finite mesh set $D_0$.  At a
stage with current set $D$, form
\[
   A=A(D)=\partial_{\rm ext}^{\Mesh}D,
   \qquad S=A^{\mathrm{ov}}.
\]
Examine the boxes of $S$ in a fixed order, reading all the edge states in
$E(B)$ when testing a box $B$.  If $B$ is inapt, replace
$D$ by $D\cup H_B$ and restart.  Once every seal box is apt, examine the
halo boxes in a fixed order.  For $B\in\Halo(S)$ compute $r_B(S)$ and, if
\[
       q_B(S)>\eta_{K,r_B(S)},
\]
replace $D$ by $D\cup H_B$ and restart.  If a later candidate requires a test with the same specification,
retain the earlier outcome.  Stop precisely
when $S$ is welcoming.  By \eqref{e.absolute-seal-bad-box}, every update is
based at a box $B$ for which $\mathsf E_B$ occurs.

\begin{definition}[Seal completion]
\label{d.seal-completion}
For a finite mesh set $D$, let $S_*$ be the terminal overlapping-box
collection produced by the preceding search, and let
$\Seal(D)=S_*^\circ$ denote its reduction from
Lemma~\ref{l.mesh-boundary-capture}(2).
\end{definition}

\begin{proposition}[Completion gives a seal]
\label{p.weak-frontier-to-completed-seal}
There are $K_0\in\N$ and $C,c>0$ such that, for every $K\ge K_0$
and every deterministic finite mesh set $D$, the completion terminates
almost surely.  Let $D_*$ be the value of the mesh set when it terminates.
Then
\begin{equation}
\label{e.completion-rooted-tail}
 \PP\bigl(|D_*\setminus D|\ge C|A(D)|+v\bigr)\le Ce^{-cv},
 \qquad v\ge0.
\end{equation}
Its terminal value $\Seal(D)$ is a seal and every mesh box of $D$ belongs to
$\MeshInt(\Seal(D))$.  For every deterministic output $S$,
\[
   \{\Seal(D)=S\}\in\cF_{\ECap(S)}.
\]
In particular, conditionally on $\cF_{\Seal(D)}$, the edge states on
$\mathcal U(\Seal(D))$ have their Bernoulli product law.
\end{proposition}

\noindent{\bf Proof.}
Let $D_t$ be the mesh set after $t$ updates, with $D_0=D$, and put
$R=A(D)^{[1]}$.  Write
\[
  \mathcal C_K(R)
  =\bigcup\{C:C\text{ is a $*$-component of }\mathcal B_K,
                    \ C\cap R\ne\varnothing\}.
\]
If the update from $D_t$ to $D_{t+1}$ occurs, let $B_t$ be the box of the
failed test and let $\tau_t$ denote its specification.  We claim, by induction
on $t$, that $\MeshSupp(B_t)^{[3]}\subseteq\mathcal C_K(R)$ whenever this
update occurs.  Indeed, $\mathsf E_{B_t}$ occurs by
\eqref{e.absolute-seal-bad-box}.  If $\tau_t$ was required by the initial
candidate, then $\MeshSupp(B_t)$ meets~$R$, and its $3$-thickening belongs
to $\mathcal C_K(R)$.

Otherwise, choose the largest $j<t$ for which $\tau_t$ was not required
by the candidate based on $D_j$.  The update at $B_j$ made this test newly
required, so Lemma~\ref{l.changed-seal-tests-local} gives
\[
  \MeshSupp(B_t)\cap\MeshSupp(B_j)^{[2]}\ne\varnothing.
\]
Since both $\mathsf E_{B_j}$ and $\mathsf E_{B_t}$ occur, their thickened mesh
supports belong to the same $*$-component of $\mathcal B_K$.
The induction hypothesis therefore proves the claim.  Every update adds
$H_{B_t}=\MeshSupp(B_t)^{[1]}$, and hence, at every finite stage,
\[
       D_t\setminus D\subseteq\mathcal C_K(R).
\]

Every box $B_t$ carrying a failed test has mesh support meeting
$A(D_t)$: directly for a seal test, and by
Lemma~\ref{l.mesh-boundary-interface} for a halo test.
Since $A(D_t)$ is disjoint from $D_t$, every update strictly enlarges
$D_t$.  The components of $\mathcal B_K$ are almost surely finite, so
$\mathcal C_K(R)$ is almost surely finite.  There can therefore be only
finitely many updates.  At each stage there are finitely many tests, and
thus the search terminates almost surely.  We may now put $D_*$ equal to
its terminal mesh set, obtaining
\begin{equation}
\label{e.completion-bad-hull}
       D_*\setminus D\subseteq\mathcal C_K\bigl(A(D)^{[1]}\bigr).
\end{equation}

The rooted exponential-moment bound above gives
constants $\lambda_*,c_*>0$ from
Lemma~\ref{l.subcritical-star-components} such that, for every finite mesh set $R$,
\[
 \EE\exp\{\lambda_*|\mathcal C_K(R)|\}\le \exp\{c_*|R|\}.
\]
Since $|A(D)^{[1]}|\le C_d|A(D)|$, Markov's inequality applied to
\eqref{e.completion-bad-hull} proves~\eqref{e.completion-rooted-tail}, after
changing constants.

Put
\[
        A_*=A(D_*),\qquad S_*=A_*^{\mathrm{ov}},
        \qquad S=S_*^\circ=\Seal(D).
\]
At termination every box of $S_*$ is apt and every halo inequality holds, so
$\mathrm{Wel}(S_*)$ occurs.  Lemma~\ref{l.reduction-preserves-welcome} gives
$\mathrm{Wel}(S)$.  Lemma~\ref{l.mesh-boundary-capture}(2) gives
\[
   S=(S\cap\Mesh)^{\mathrm{ov}},
   \qquad
   S\cap\Mesh
   =\partial_{\rm ext}^{\Mesh}
       \bigl((\Mesh\cap\Capt_{\overline N}(S))\setminus(S\cap\Mesh)\bigr),
\]
which is \eqref{e.seal-mesh-geometry}.  The same lemma also shows
that every mesh box of the initial set $D$ belongs to $\MeshInt(S)$.

It remains to verify the stopping-domain condition.  At a stage with
mesh set $D_t$, put $T=A(D_t)^{\mathrm{ov}}$.  Every edge queried while
testing this candidate belongs to $E(T)$: a seal test at $B\in T$ uses
only $E(B)$, while a halo test at $B\in\Halo(T)$ uses only
$\mathcal I_B(T)\subseteq E(T)$ by~\eqref{e.halo-reads-seal}.
Since $T=A(D_t)^{\mathrm{ov}}$, the mesh support of each box of $T$
lies in $A(D_t)\subseteq\MCap(T)$, and hence
$E(T)\subseteq\ECap(T)$.  As $D_t\subseteq D_*$,
Lemma~\ref{l.mesh-exterior-monotone} gives
$\ECap(T)\subseteq\ECap(S_*)$.  Reduction preserves capture, so
$\ECap(S_*)=\ECap(S)$.  Thus every edge queried by the completion
search belongs to $\ECap(S)$.

If a second configuration agrees with the first on $\ECap(S)$, the fixed
orders reproduce the same tests, updates, $D_*$ and $S$.  Hence
$\{\Seal(D)=S\}\in\cF_{\ECap(S)}$.  Lemma~\ref{l.freshness} gives the product
law on $\mathcal U(S)$.
\qed

\section{Seal enlargement and layer episodes}

With a seal already in force, further finite exploration is absorbed by
enlarging its capture while preserving the product law in the unread exterior.
We establish the conditional seal-completion estimate that permits this
enlargement and control the accompanying boundary growth.

A small cap added to a seal may close the mouth of a large cavern.  The capture
then gains the cavern at once, though little new boundary has been built and
the exploration may have reached only a few vertices within the cavern.  Our estimates
keep track of what was reached and what had to be built, rather than the
volume brought into the capture.

\begin{figure}[htbp]
\centering
\includegraphics[width=0.75\textwidth]{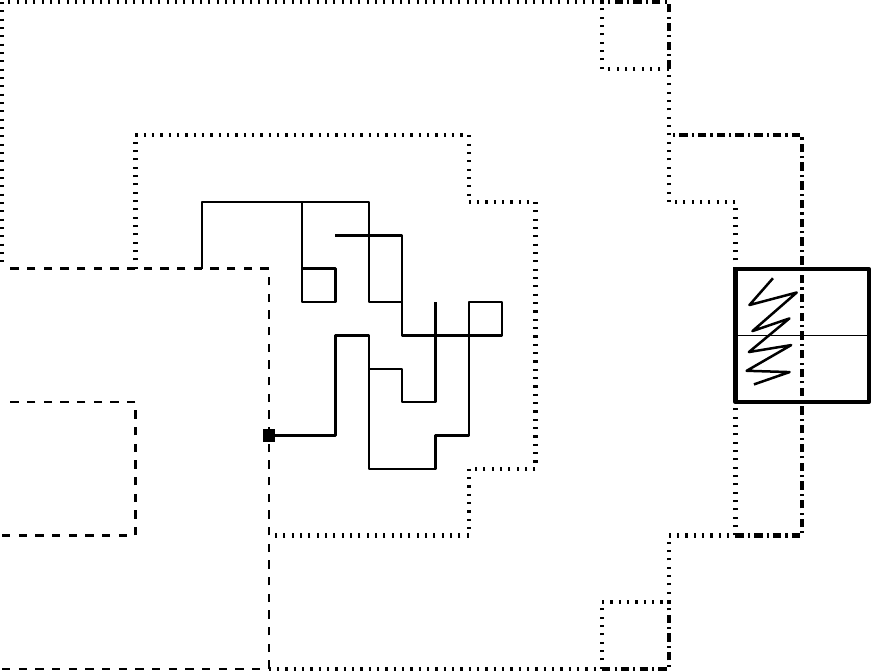}
\caption{Seal enlargement after a growth episode. The square-marked
virgin seed lies on the boundary of an existing seal, delimited by
dashed lines in the lower-left. From the seed proceeds an open
exploration that reenters the seal. A dotted-line enclosing quilt,
together with the old seal, captures this growth. 
The capture of the quilt is transferred to the mesh and combined with the old mesh interior. The resulting candidate boundary is represented schematically by the outermost contour on the right, following the dotted outline with the local adjustments marked by bold dash-dotted lines.
This boundary serves as the starting point for an
outward search to construct a new seal. The box on the right lies in
the halo of the candidate boundary, with its captured-side conditioning
set indicated by sketched lines. As the search continues, the box may
remain in the halo even as its captured side comes to include one or other of its right-hand constituent boxlets; or it may leave the halo
as the enclosed region expands.}\label{f.halotest}
\end{figure}

In this light, consider a finite exploration episode from a virgin seed, such as that depicted in Figure~\ref{f.halotest}.
The search may venture into the exterior of the existing seal, return through
its capture, and emerge again before coming to an end.  The old capture and a
new enclosing quilt together house the record of the episode; the quilt may
close against part of the old seal.  With the
excursion enclosed, the candidate boundary need not yet be welcoming.  Seal
completion builds outward from this candidate, enlarging the enclosed mesh
region whenever a required test fails.

As the boundary moves outward, a halo box may remain exterior while additional
boxlets enter its captured-side conditioning set.  A welcoming test passed
under the earlier conditioning need not pass under the new.  The successive
square-root allowances introduced in Section~\ref{ss.visible-welcoming-exteriors}
are designed for such changes.  Captured-side conditioning sets consist of
whole boxlets, of which an overlapping box has only a dimension-dependent
number.  The possible failures can therefore be controlled together, rather
than successively along the adaptive completion search.  A relative defect
field specified from the old seal covers these failures, along with
inaptitude.  Locality of newly required tests then confines the further
spread of repairs to subcritical defect components rooted near the initial
obligations.

The estimates accordingly distinguish the open growth reached outside the
old capture from the new boundary that had to be built.  The former is bounded
by the volume enclosed by the quilt and has a stretched-exponential tail.
The size of the quilt and the further mesh region built during completion
have exponential tails and control the increase in seal size.

We first adapt seal completion to an existing seal and a partially revealed
history.  We then control the boundary added when a finite exterior record is
absorbed.  Finally, we combine this enlargement bound with open escape in the
layer to obtain the estimates for an entire virgin-seed episode, first under
ordinary conditioning and then with the unread field preserved for successive
episodes.

\subsection{Relative seals and conditional completion}

The seal-completion procedure of Section~\ref{ss.trap-cluster-to-seal} absorbs a mesh neighbourhood of each failed-test box. We now use it with a seal already in force. Further exploration may require a larger capture, but the old seal has passed its tests, and much of that certification remains valid.
  We must control the repairs prompted
by the enlargement while preserving the product law of the still-unread
exterior.

\begin{definition}[Seals relative to a conditioning field]
\label{d.relative-seal}
Let $\mathcal G$ be a $\sigma$-field and let
$\mathscr S\subseteq\overline N$ be a finite $\mathcal G$-measurable quilt
collection.  We call $\mathscr S$ a \emph{$\mathcal G$-seal} if the following
hold: almost surely, $\mathscr S$ satisfies \eqref{e.seal-mesh-geometry} and is
welcoming; the states on $E(\mathscr S)$ are $\mathcal G$-measurable;
and, conditionally on $\mathcal G$, the edge states on
$\mathcal U(\mathscr S)$ have their Bernoulli-$p$ product law.  Thus
$\mathcal U(\mathscr S)$ is the unread exterior edge set relative to $\mathcal G$.  Equivalently, we
say that~$\mathscr S$ is a seal relative to~$\mathcal G$.
\end{definition}

Every ordinary seal is an $\cF_{\mathscr S}$-seal, by the stopping-domain
property and Lemma~\ref{l.freshness}.  Unlike $\cF_{\mathscr S}$, the
field~$\mathcal G$ need not reveal
all captured edges: cluster exploration and sealing may leave some of
them unread.  The seal and its boundary-edge states are known, while the
exterior remains conditionally product; no product assertion is made for
unread captured edges.  The product law on the random exterior is understood on
$\{\mathscr S=S\}$ for deterministic~$S$.  Virgin and sealed seeds are defined relative
to~$\mathcal G$.

We first identify the new tests caused by enlargement.  Given a seal
$\mathscr S$, ordinary or relative, let
$D_1=\MeshInt(\mathscr S)$ be the captured mesh boxes behind its boundary.
Further suppose given a finite mesh set  $D_2$ that is to be enclosed together with $D_1$.
Form the candidate boundary by taking the exterior mesh boundary of \(D_0=D_1\cup D_2\) and passing to the associated overlapping-box boundary:
\[
      S^{(0)}=A(D_0)^{\mathrm{ov}}.
\]
A test inherited from~$\mathscr S$ with the same specification has already
passed.  Recall that a seal test is specified by its box, while a halo
test is specified by its box and its captured-side conditioning set.
Consequently, a halo box may remain in place but require a new test if
additional boxlets enter the capture.  The boxes carrying tests newly
required by the first candidate form
\begin{align*}
 J&=J^{\rm seal}\cup J^{\rm halo},\\
 J^{\rm seal}&=S^{(0)}\setminus\mathscr S,\\
 J^{\rm halo}&=\bigl\{B\in\Halo(S^{(0)}):
       B\notin\Halo(\mathscr S)
       \text{ or }\mathcal I_B(S^{(0)})\ne\mathcal I_B(\mathscr S)\bigr\}.
\end{align*}
These are the initial tests needed to make the candidate boundary welcoming;
in the enlargement application, their number will be controlled by the
enclosing quilt.  A failed test prompts an outward update, which may in turn
require further tests.  While a box remains exterior, its captured-side
conditioning set can only gain whole boxlets.  The finite-boxlet reserve
controls possible failures as this conditioning set grows.  Relative to the
old seal, we specify a defect event covering failure under any strict
enlargement of the set, as well as inaptitude of the box.  Its specification
is fixed before the adaptive completion search, whether or not the search
eventually requires any of these tests.

For a deterministic boundary $S$ satisfying~\eqref{e.seal-mesh-geometry},
use $\mathscr L_B(S)$ from Section~\ref{ss.visible-welcoming-exteriors}
when $B\in\Halo(S)$, and put $\mathscr L_B(S)=\varnothing$
for $B\in\Ext_{\overline N}(S)\setminus\Halo(S)$.
Thus $I_B(\mathscr L_B(S))=\mathcal I_B(S)$ for exterior boxes.
Let $\mathscr A$ range over sets satisfying
$\mathscr L_B(S)\subsetneq\mathscr A\subsetneq\mathscr L(B)$.
Define the \emph{relative defect event}
\begin{equation}
\label{e.relative-seal-bad-box}
 \mathsf E_B(S)
 =\{B\text{ is not apt}\}
  \cup
  \bigcup_{\mathscr A}
  \{q_{B,\mathscr A}>\eta_{K,|\mathscr A|}\}.
\end{equation}

The completion search tests only its current candidate, not
$\mathsf E_B(S)$.  The latter event dominates failures that may arise
as the boundary changes.

Define the \(3\)-thickened relative defect mesh set by
\[
 \mathcal B_K(S)
 =\bigcup_{\{B\in\Ext_{\overline N}(S):\,\mathsf E_B(S)\}}
       \MeshSupp(B)^{[3]}.
\]

\begin{lemma}[Relative defect boxes]
\label{l.relative-bad-boxes}
For sufficiently large fixed $K$, let \(S\) be deterministic and satisfy
\eqref{e.seal-mesh-geometry}.  On $\mathrm{Wel}(S)$, uniformly in $S$,
\[
 \PP\bigl(\mathsf E_B(S)\mid\cF_{\ECap(S)}\bigr)
       \le C_d\delta_K,
 \qquad B\in\Ext_{\overline N}(S).
\]
Conditionally on $\cF_{\ECap(S)}$, the field
$(\mathbf1_{\mathsf E_B(S)})_{B\in\Ext_{\overline N}(S)}$ is
one-dependent.  Moreover, on $\mathrm{Wel}(S)$,
\[
  \sup_{M\in\Mesh}
  \PP\bigl(M\in\mathcal B_K(S)\mid\cF_{\ECap(S)}\bigr)
  \le C'_d\delta_K.
\]
Conditionally on $\cF_{\ECap(S)}$, the field
$(\mathbf1_{\{M\in\mathcal B_K(S)\}})_{M\in\Mesh}$ is
$7$-dependent.  On $\mathrm{Wel}(S)$, its $*$-components are
stochastically dominated by subcritical Bernoulli $*$-components and have a
uniform exponential size tail.  The same assertions hold, almost surely,
when $S$ is the realised value of a relative seal and the conditioning is on
its relative-seal field.
\end{lemma}

\noindent{\bf Proof.}
If $B$ is outside the halo, then
$q_{B,\varnothing}=\PP(B\text{ is not apt})\le\varepsilon_K=\eta_{K,0}$.
If $B$ is in the halo, welcoming gives
\[
 q_{B,\mathscr L_B(S)}\le\eta_{K,|\mathscr L_B(S)|}.
\]
For every strict enlargement $\mathscr A\supsetneq\mathscr L_B(S)$,
Lemma~\ref{l.sqrt-welcoming-note} therefore gives, conditionally on the old
captured information, probability at most $\delta_K$ that the corresponding
term in~\eqref{e.relative-seal-bad-box} occurs.  The probability that $B$
itself is inapt is at most $\delta_K$ as well.  There are at most $2^{m_d}$
possible enlargements, so the first assertion follows by a union bound.

After the captured edge states are fixed, $\mathsf E_B(S)$ depends only on
the remaining edge states on $E(B)$.  Boxes at transverse index distance at
least two use disjoint remaining edge sets, giving one-dependence.  The
$3$-thickening then gives the stated $7$-dependence of the mesh field, while
a union bound gives the displayed mesh-site estimate.  By
Liggett--Schonmann--Stacey~\cite{LSS}, for $K$ large enough that
$C'_d\delta_K$ is below a fixed subcritical threshold for the $*$-graph, the
component assertion follows.  For a relative $\mathcal G$-seal, condition instead on $\mathcal G$
and its realised value $S$.  In an exterior box $B$, the captured-side
edges $\mathcal I_B(S)$ lie in $E(S)$ by~\eqref{e.halo-reads-seal}
and thus have known states; the other edges of $E(B)$ lie in
$\mathcal U(S)$ and are conditionally product.  Consequently,
Lemma~\ref{l.sqrt-welcoming-note} applies to each strict enlargement of the
captured-side boxlet set, giving the same $C_d\delta_K$ bound by a union
bound.  Moreover, exterior boxes at transverse index distance at least two
depend on disjoint subsets of the unread product field, so their defect
indicators are conditionally independent.  The mesh-field bounds follow as
before.
\qed

The conditional completion result bounds the additional mesh region in terms
of the number $|J|$ of initially required tests.  The new history may contain
additional information about edges within the old capture as well as within
the finite region~$D_2$; the remaining old-exterior edges are required to
retain their conditional product law.

\begin{lemma}[Conditional seal completion]
\label{l.conditional-rooted-completion}
For all sufficiently large $K$, let $\mathscr S$ be a $\mathcal G$-seal,
let $\mathcal H\supseteq\mathcal G$ be a $\sigma$-field, and let
$D_2\subseteq\Mesh$ be finite and $\mathcal H$-measurable.
Put $D_1=\MeshInt(\mathscr S)$ and $D_0=D_1\cup D_2$, and let~$J$
be the set of boxes carrying tests newly required by the candidate
$A(D_0)^{\mathrm{ov}}$, as specified above.  Suppose that,
conditionally on $\mathcal H$, the edge states on
\[
 \mathcal U(\mathscr S)\setminus\bigcup_{M\in D_2}E(M)
\]
have their Bernoulli-$p$ product law.  Run seal completion from $D_0$, treating tests inherited unchanged
from the old seal as already passed.
The completion terminates almost surely.  Let $D_*$ be its terminal mesh
set and put $H=D_*\setminus D_0$.  Then, for constants
$c,C>0$ uniform in $K$, $\mathcal G$, $\mathcal H$, $\mathscr S$, $D_2$
and the conditioned history,
\begin{equation}
\label{e.relative-completion-tail}
  \PP\bigl(|H|\ge C|J|+v\mid\mathcal H\bigr)\le Ce^{-cv},
  \qquad v\ge0.
\end{equation}
If $\mathcal Q$ is the set of edges queried by this completion and
\[
 \mathcal H^+=\mathcal H\vee
     \sigma(\mathcal Q,\omega|_{\mathcal Q}),
\]
then the reduction $\mathscr S^+$ of the terminal collection is a
$\mathcal H^+$-seal.  Moreover
\[
        D_1\cup D_2\subseteq \MeshInt(\mathscr S^+).
\]
\end{lemma}

\noindent{\bf Proof.}
Condition on~$\mathcal H$.  The old seal, $D_2$, the initial candidate
$S^{(0)}$ and $J$ are then fixed.  We first check that the tests used in
completion depend only on states already available under~$\mathcal H$
and on the stipulated product field.  This also permits us to use the
relative defect estimate under the stronger conditioning.

Let $D\supseteq D_0$ be the mesh set at any subsequent stage.  A box
$C\in A(D)^{\mathrm{ov}}$ has $\MeshSupp(C)\subseteq A(D)$, disjoint
from~$D_2$.  The same separation holds for a halo box.  Indeed, if
$C\in\Halo(A(D)^{\mathrm{ov}})$ and
$M\in D_2\cap\MeshSupp(C)$, then
$M\in\MCap(A(D)^{\mathrm{ov}})$, and
Lemma~\ref{l.mesh-boundary-interface} puts~$M$ in $A(D)$, a
contradiction.  Therefore, throughout the search,
\begin{equation}
\label{e.relative-completion-interior-separation}
 E(C)\cap\bigcup_{M\in D_2}E(M)=\varnothing
 \qquad\text{for }C\in A(D)^{\mathrm{ov}}\cup\Halo(A(D)^{\mathrm{ov}}).
\end{equation}
We also have that every current test box belongs either to the old
seal or to its exterior.  Indeed, the seal identities
\eqref{e.seal-mesh-geometry} give
$\mathscr S=A(D_1)^{\mathrm{ov}}$.  By
Lemma~\ref{l.mesh-exterior-monotone} and
\eqref{e.mesh-transfer-full-exterior}, the old capture is contained in
that of the current boundary; thus no old captured box can be a
current halo box.  An old captured box outside~$\mathscr S$ cannot be a
current seal box either: by~\eqref{e.mesh-transfer-capture} and
\eqref{e.mesh-transfer-full-exterior}, its mesh support meets~$D_1$,
whereas the support of a current seal box lies in $A(D)$, disjoint
from~$D$.
Thus tests use only known old-seal states or edges of the old exterior,
and by the separation just proved they never consult the potentially
nonfresh states within~$D_2$.  In the enlargement application, this keeps
the partial revelations of the enclosed episode away from the completion
tests.  Form $\mathcal B_K^{D_2}(\mathscr S)$ from
$\mathcal B_K(\mathscr S)$ by omitting the defect boxes~$B$ for which
$E(B)$ meets $\bigcup_{M\in D_2}E(M)$.  The observation following
Lemma~\ref{l.sqrt-welcoming-note} applies under~$\mathcal H$ to
\eqref{e.welcome-martingale} and~\eqref{e.welcome-reserve-transition}.
The proof of Lemma~\ref{l.relative-bad-boxes} then gives conditional
$7$-dependence, subcritical domination and a uniform exponential
component tail for this truncated field.

Next we identify the tests that can actually fail.  Any test of
$S^{(0)}$ whose box is outside~$J$ is an inherited test with unchanged
specification, and has passed already.  Moreover, a test of the old seal
which is absent from the initial candidate cannot reappear later with
its old specification.  To see this for a seal test, note that, whenever
$D\supseteq D_0$,
\[
      A(D_1)\cap A(D)\subseteq A(D_0).
\]
Indeed, a site of the first boundary is $*$-adjacent to~$D_1$; if it
also lies in~$A(D)$, it is connected to infinity outside~$D$, hence
outside~$D_0$.  The corresponding old seal box must therefore already
belong to~$S^{(0)}$.  For an old halo test, suppose that its original
specification is required at a later candidate $T=A(D)^{\mathrm{ov}}$.
By mesh-exterior monotonicity, its box $B$ is exterior for~$S^{(0)}$,
and capture monotonicity gives
\[
  \mathcal I_B(\mathscr S)
     \subseteq\mathcal I_B(S^{(0)})
     \subseteq\mathcal I_B(T).
\]
The first and last sets agree by the supposition, so the middle set
agrees as well.  Since $B$ was an old halo box, the same halo test was
already required at~$S^{(0)}$.  Hence all inherited certifications that
can be used later are already accounted for at the outset.

Suppose a test fails at a later candidate $T=A(D)^{\mathrm{ov}}$.
By the preceding geometric observation, its box~$B$ belongs to the old
seal or its exterior.  Tests on the old seal are seal tests and have
already passed.  Hence a failed test occurs only at an exterior box
of~$\mathscr S$.
If it is a seal test, failure implies that $B$ is inapt, the first
event in~\eqref{e.relative-seal-bad-box}.  If it is a halo test,
capture monotonicity gives
\[
    \mathcal I_B(\mathscr S)\subseteq\mathcal I_B(T).
\]
Equality would give a test already passed by the old welcoming seal,
so the inclusion is strict.  Since~$B$ is still a halo box, the larger
set omits at least one boxlet.  Failure therefore gives one of the
remaining terms in~\eqref{e.relative-seal-bad-box}.  In either case,
every failed test at an old-exterior box is covered by
$\mathsf E_B(\mathscr S)$, whether or not its box belongs to~$J$.

The locality of newly required tests allows the outward repairs to be traced
to the initial obligations~$J$.  A failed test not among these obligations
must have become necessary after an earlier nearby update.  The relative
defect field controls the further propagation of repairs.  Put
\[
 R_J=\bigcup_{B\in J}\MeshSupp(B)^{[3]},
 \qquad
 \mathcal C_K^{\mathscr S}(R_J)
 =\bigcup\{C:C\text{ is a $*$-component of }
       \mathcal B_K^{D_2}(\mathscr S),\ C\cap R_J\ne\varnothing\}.
\]
Write $D_t$ for the mesh set after $t$ updates, starting at~$D_0$,
and let $B_t$ and $\tau_t$ be the box and specification of the test
that causes the update from $D_t$ to~$D_{t+1}$.  We claim inductively
that either $B_t\in J$ or
$\MeshSupp(B_t)^{[3]}\subseteq\mathcal C_K^{\mathscr S}(R_J)$.

For $B_t\notin J$, the preceding argument gives
$\mathsf E_{B_t}(\mathscr S)$; separation
\eqref{e.relative-completion-interior-separation} makes its thickened
mesh support part of $\mathcal B_K^{D_2}(\mathscr S)$.  Since a test of
the initial candidate outside~$J$ was already successful,
$\tau_t$ was not required initially.  Choose the largest $j<t$ for
which it was not required by the candidate based on~$D_j$.
The update at~$B_j$ made it newly required, and
Lemma~\ref{l.changed-seal-tests-local} gives
\[
 \MeshSupp(B_t)\cap\MeshSupp(B_j)^{[2]}\ne\varnothing.
\]
If $B_j\in J$, the thickened support of~$B_t$ lies in a defect
component meeting~$R_J$.  Otherwise the induction hypothesis places
the thickened support of~$B_j$ in such a component, and the two
supports belong to the same component.  This proves the claim, and
since each update adds $H_{B_t}=\MeshSupp(B_t)^{[1]}$, it gives the
pathwise bound
\[
  D_t\setminus D_0\subseteq
  R_J\cup\mathcal C_K^{\mathscr S}(R_J)
  \qquad\text{at every finite stage }t.
\]
Each failed test has mesh support meeting~$A(D_t)$, using
Lemma~\ref{l.mesh-boundary-interface} in the halo case, and so each
update strictly enlarges~$D_t$.  The right-hand side is almost surely
finite under~$\mathcal H$; hence only finitely many updates are possible.
There are finitely many tests at each stage, so completion terminates.
Its growth $H=D_*\setminus D_0$ satisfies the same containment.
The conditional subcritical domination of the truncated field
$\mathcal B_K^{D_2}(\mathscr S)$ established above, together with the
exponential-moment bound in Lemma~\ref{l.subcritical-star-components}, yields
\[
 \EE\!\left[
    e^{\lambda|\mathcal C_K^{\mathscr S}(R_J)|}
       \middle|\mathcal H\right]
    \le e^{c_0|R_J|}
\]
for uniform positive $\lambda,c_0$.  Since $|R_J|\le C_d|J|$,
exponential Markov proves~\eqref{e.relative-completion-tail}.

It remains to check that this finite repair returns a seal under the
updated history, so the operation can be repeated.  At termination all
required tests have passed, including inherited ones.  The terminal
candidate is therefore welcoming; its reduction~$\mathscr S^+$ is
welcoming and satisfies~\eqref{e.seal-mesh-geometry} by
Lemmas~\ref{l.reduction-preserves-welcome} and
\ref{l.mesh-boundary-capture}(2).  The latter lemma also gives
\[
     D_1\cup D_2\subseteq\MeshInt(\mathscr S^+).
\]
All completion queries lie on a current seal box or on the captured
part of a current halo box.  Monotonicity of capture places these
queries in $\ECap(\mathscr S^+)$, as in the proof of
Proposition~\ref{p.weak-frontier-to-completed-seal}.  Replaying the
fixed-order search from the information in~$\mathcal H^+$ determines
$\mathscr S^+$.  The edge states of every terminal seal box were either
known under~$\mathcal H$ or queried during completion, so they too
are $\mathcal H^+$-measurable.

Finally, mesh-exterior monotonicity and capture of~$D_2$ give
\[
   \mathcal U(\mathscr S^+)
     \subseteq \mathcal U(\mathscr S)
                    \setminus\bigcup_{M\in D_2}E(M).
\]
Conditionally on~$\mathcal H$, each query is selected from this product
field using the answers to preceding queries.  For every deterministic
finite edge set~$Q$, the event $\{\mathcal Q=Q\}$ is measurable with
respect to $\mathcal H\vee\cF_Q$.  Lemma~\ref{l.freshness}, applied
conditionally on~$\mathcal H$, leaves the states on
$\mathcal U(\mathscr S^+)$ with their Bernoulli product law given
$\mathcal H^+$.  Hence $\mathscr S^+$ is an $\mathcal H^+$-seal.
\qed

When no old seal is present, the same argument confines completion growth
to a neighbourhood of the initial mesh boundary and the defect components
meeting it.  We record this conditional form for later use.

\begin{corollary}[Conditional bound on completion growth]
\label{c.conditional-completion-hull}
For sufficiently large fixed $K$, let $D\subseteq\Mesh$ be finite
and $\mathcal H$-measurable.  Suppose that, conditionally on~$\mathcal H$,
the edges of~$E(N)$ outside
$\bigcup_{M\in D}E(M)$ have their Bernoulli-$p$ product law, and run seal
completion from~$D$.  Let $D_*$ be the terminal mesh set.  There is a fixed
$r_0<\infty$ and a $7$-dependent mesh field~$\mathcal B_K^D$, conditionally
dominated by a fixed subcritical Bernoulli field, such that
\begin{equation}
\label{e.conditional-completion-hull}
 D_*\setminus D
 \subseteq
 A(D)^{[r_0]}\cup
 \mathcal C_K^D\bigl(A(D)^{[r_0]}\bigr),
\end{equation}
where $\mathcal C_K^D(R)$ is the union of the $*$-components of
$\mathcal B_K^D$ meeting~$R$.  Moreover, for some $\lambda,C>0$ and every
mesh site~$z$, the component~$C^D(z)$ of~$\mathcal B_K^D$ containing~$z$
(with the empty set when $z\notin\mathcal B_K^D$) satisfies
\[
 \EE\left[e^{\lambda|C^D(z)|}\mid\mathcal H\right]\le C.
\]
The constants are uniform in~$D$ and in the conditioned history.
\end{corollary}

\noindent{\bf Proof.}
Apply Lemma~\ref{l.conditional-rooted-completion} directly with the empty
seal, trivial $\mathcal G$, and $D_2=D$.  Since
$\mathcal U(\varnothing)=E(N)$, its conditional product-law hypothesis is
precisely the assumption of the corollary.

Let $\mathcal B_K^D$ be the $3$-thickened defect mesh set formed from
$\mathcal B_K$ by omitting every defect box $B$ whose edge set meets
$\bigcup_{M\in D}E(M)$.  The proof of Lemma~\ref{l.conditional-rooted-completion} gives
conditional $7$-dependence and subcritical domination for this field.  With the empty
seal, $J$ consists of the boxes carrying the seal and halo tests of
$A(D)^{\mathrm{ov}}$.  Lemma~\ref{l.frontier-addition-range} places their
mesh supports in $A(D)^{[1]}$, so $R_J\subseteq A(D)^{[4]}$.
The pathwise inclusion
$H\subseteq R_J\cup\mathcal C_K^D(R_J)$ from that proof therefore gives
\eqref{e.conditional-completion-hull}, with $r_0=4$.
The rooted exponential moment is the same uniform subcritical component
estimate used in the lemma.
\qed

\subsection{Seal enlargement}
\label{ss.efficient-extension}

We now turn conditional completion into a boundary-growth estimate for
seal enlargement.  Two geometric facts about mesh boundaries allow us to
compare the enlarged seal with the enclosing quilt and the further mesh
region added during completion.

\begin{lemma}[Boundary of an enclosed region]
\label{l.boundary-enclosed-region}
There is $C_d<\infty$ such that, for every finite quilt collection
$\mathscr Q$, $|A(R(\mathscr Q))|\le C_d|\mathscr Q|$.
Every $M\in A(R(\mathscr Q))$, viewed as an overlapping box with even
transverse index, satisfies $d_\infty(M,\mathscr Q)\le3$.
\end{lemma}

\noindent{\bf Proof.}
Fix $M\in A(R(\mathscr Q))$.  There is a $*$-neighbour
$M'\in R(\mathscr Q)$, and by the definition of $R(\mathscr Q)$ there is
$B\in\Capt_{\overline N}(\mathscr Q)$ with
$M'\in\MeshSupp(B)$.  In the transverse indexing of overlapping boxes, we have
\[
       d_\infty(M,M')\le2,
       \qquad d_\infty(M',B)\le1.
\]
The mesh box $M$ is not captured by $\mathscr Q$, since otherwise
$M\in R(\mathscr Q)$.  Thus any coordinatewise nearest-neighbour interpolation
from $B$ to $M$ meets $\mathscr Q$, and every box on this interpolation is at
$\ell^\infty$-distance at most $3$ from $M$.  Hence
$d_\infty(M,\mathscr Q)\le3$.  Bounded local multiplicity now gives
$|A(R(\mathscr Q))|\le C_d|\mathscr Q|$.
\qed

Second, exterior boundaries behave well under unions.

\begin{lemma}[Exterior boundary of a union]
\label{l.union-of-captured-quilts}
For finite mesh sets $D_1,D_2$,
\[
 A(D_1\cup D_2)\subseteq A(D_1)\cup A(D_2).
\]
Consequently, if $S_i=A(D_i)^{\mathrm{ov}}$, then
\[
 |A(D_1\cup D_2)^{\mathrm{ov}}|\le C_d(|S_1|+|S_2|).
\]
Moreover
\[
 A(D_1)\setminus A(D_1\cup D_2)
 \subseteq \MeshInt(A(D_1\cup D_2)^{\mathrm{ov}}).
\]
\end{lemma}

\noindent{\bf Proof.}
If $M\in A(D_1\cup D_2)$, it has a $*$-neighbour in at least one of
$D_1,D_2$, while an unbounded path from $M$ avoids their union.  Hence
$M\in A(D_1)\cup A(D_2)$.  The cardinality bound follows because each mesh
box supports at most $3^{d-1}$ overlapping boxes and
$A(D_i)\subseteq S_i$.

For the final assertion, put $D_0=D_1\cup D_2$.  By
Lemma~\ref{l.mesh-exterior-monotone},
$\MCap(A(D_1)^{\mathrm{ov}})\subseteq
\MCap(A(D_0)^{\mathrm{ov}})$.  Thus every $M\in A(D_1)$ is
captured by $A(D_0)^{\mathrm{ov}}$.  If $M\notin A(D_0)$, the
mesh-transfer identity~\eqref{e.mesh-transfer-capture} places it in
$\MeshInt(A(D_0)^{\mathrm{ov}})$, as required.
\qed

\begin{proposition}[Boundary growth under seal enlargement]
\label{p.trap-cluster-to-seal}
For sufficiently large fixed $K$, let $\mathscr S$ be a seal, $s$ a
virgin seed, and $\mathcal T_s$ the quilt of
Proposition~\ref{p.finite-trap-on-nonescape}.  There is a seal $\mathscr S_s$ on
$\{O_s\text{ is finite}\}\cap\mathsf{Loc}_s^{\rm c}$ such that
\[
 \Capt_{\overline N}(\mathscr S)
    \subseteq\Capt_{\overline N}(\mathscr S_s),
 \qquad
 \MeshSupp(\mathcal R_s)\subseteq\MeshInt(\mathscr S_s).
\]
There is also a non-negative integer-valued random variable $Z_s$ such that
\begin{equation}
\label{e.efficient-boundary-growth}
   |\mathscr S_s|\le |\mathscr S|+CZ_s
\end{equation}
and, uniformly over the conditioning $\cF_{\mathscr S}$,
\begin{equation}
\label{e.nonlocal-efficient-boundary-tail}
 \PP\bigl(Z_s\ge u,\ O_s\text{ finite},\ \mathsf{Loc}_s^{\rm c}
       \mid\cF_{\mathscr S}\bigr)
 \le Ce^{-cu},\qquad u\ge0.
\end{equation}
\end{proposition}

\noindent{\bf Proof.}
Let $\mathcal G_s$ be the field of
Lemma~\ref{l.stopped-enclosing-record}, and put
\[
 D_1=\MeshInt(\mathscr S),
 \qquad
 D_2=R(\mathcal T_s),
 \qquad
 D_0=D_1\cup D_2,
 \qquad
 S^{(0)}=A(D_0)^{\mathrm{ov}}.
\]
Here $D_1$ retains the old captured mesh interior, while
$R(\mathcal T_s)$ is the mesh support of the entire overlapping-box capture
of the enclosing quilt.  Their union is the region to be enclosed, and the
proof bounds the additional boundary needed to seal it.
By that lemma, $D_2$ is $\mathcal G_s$-measurable and, conditionally on
$\mathcal G_s$, the exterior edge states outside
$\bigcup_{M\in D_2}E(M)$ retain their product law.  Thus the conditioning
hypotheses of Lemma~\ref{l.conditional-rooted-completion} hold with
$\mathcal G=\cF_{\mathscr S}$, $\mathcal H=\mathcal G_s$ and this~$D_2$.

It remains only to bound the number of boxes carrying tests required for the
initial candidate but not for the old seal.  Let
$J_s^{\rm seal}$ and $J_s^{\rm halo}$ be the two sets in that lemma.  If
$B\in J_s^{\rm seal}$, then some member of $\MeshSupp(B)$ belongs to
$A(D_0)\setminus A(D_1)$; Lemma~\ref{l.union-of-captured-quilts} places this
mesh box in $A(D_2)$.  If $B\in J_s^{\rm halo}$, its captured mesh support at
$S^{(0)}$ contains a mesh box not captured by $\mathscr S$; by
Lemma~\ref{l.mesh-boundary-interface} and the same union lemma, this mesh box
also belongs to $A(D_2)$.  A mesh box belongs to the support of only boundedly
many overlapping boxes.  Hence Lemma~\ref{l.boundary-enclosed-region} gives
\[
       |J_s^{\rm seal}|+|J_s^{\rm halo}|
       \le C_d|A(D_2)|
       \le C_d|\mathcal T_s|.
\]
With these initial testing obligations bounded, apply
Lemma~\ref{l.conditional-rooted-completion} with~$D_2$.  If
$H_s^{\rm comp}$ is the set of mesh boxes adjoined after $D_0$, then
\[
 \PP\bigl(|H_s^{\rm comp}|\ge C|\mathcal T_s|+v\mid\mathcal G_s\bigr)
 \le Ce^{-cv}.
\]
Let the resulting seal be $\mathscr S_s$.  The placement conclusion of the
lemma gives $D_1\cup D_2\subseteq\MeshInt(\mathscr S_s)$.
Lemma~\ref{l.stopped-enclosing-record} therefore yields
$\MeshSupp(\mathcal R_s)\subseteq\MeshInt(\mathscr S_s)$.
Since $D_1\subseteq D_*$, Lemma~\ref{l.mesh-exterior-monotone} and
the mesh-transfer formula~\eqref{e.mesh-transfer-full-exterior}, followed
by preservation of capture under reduction
(Lemma~\ref{l.mesh-boundary-capture}(2)), give
$\Capt_{\overline N}(\mathscr S)\subseteq
 \Capt_{\overline N}(\mathscr S_s)$.

The relative-seal conclusion of Lemma~\ref{l.conditional-rooted-completion}
gives the product exterior after the new queries are adjoined.  Moreover the
old stopping-domain record, the queries used to obtain $\mathcal T_s$, and
the completion queries all lie in $\ECap(\mathscr S_s)$.  Since $s$ is
$\cF_{\mathscr S}$-measurable and the old capture is contained in the new
one, the captured edge states determine the old seal and the seed.  Rerunning
the subsequent fixed-order explorations then shows that, for every
deterministic possible output $S$,
\[
 \{O_s\text{ is finite},\ \mathsf{Loc}_s^{\rm c},\ \mathscr S_s=S\}
 \in\cF_{\ECap(S)}.
\]
This verifies the stopping-domain condition for the output on the finite
nonlocal event, as specified after Definition~\ref{d.seal}.

Every box of $S^{(0)}$ which is not already in $\mathscr S$ has a mesh
support member in
$A(D_0)\setminus A(D_1)\subseteq A(D_2)$.  Bounded local multiplicity and
Lemma~\ref{l.boundary-enclosed-region} therefore give
\[
       |S^{(0)}|\le |\mathscr S|+C|\mathcal T_s|.
\]
Lemma~\ref{l.boundary-growth-under-addition} then yields for the terminal
collection before reduction, and hence also after reduction,
\[
 |\mathscr S_s|
 \le |\mathscr S|+C|\mathcal T_s|+C|H_s^{\rm comp}|.
\]
Set $Z_s=|\mathcal T_s|+|H_s^{\rm comp}|$, accounting for the enclosing
quilt and the further mesh boxes added during seal completion.  Corollary
\ref{c.trap-cardinality-tail}, the preceding conditional estimate and the
tower property prove~\eqref{e.nonlocal-efficient-boundary-tail}.
\qed

\subsection{Open escape and episode bounds}

We now combine open escape in the layer with the enlargement bounds to control
the open growth during an episode.  For a seed~$s$ and an incoming seal~$\mathscr S$,
write
\[
 G_s=\left|V(O_s)\setminus
      \bigcup_{M\in\MCap(\mathscr S)}V(M)\right|,
\]
the number of vertices of~$O_s$ outside the previously captured mesh boxes.
Thus $G_s$ measures open growth beyond the old capture, not the potentially
much larger volume newly captured by the enlarged seal.
We first assemble the one-episode estimates for an ordinary seal,
conditioning on the entire old capture.

\begin{proposition}[Open escape and seal enlargement]
\label{p.seed-escape-or-seal-extension}
Let $d\ge3$ and $\alpha=(d-2)/(d-1)$.  For sufficiently large fixed $K$
there are $c_{\rm esc}>0$ and $C,c\in(0,\infty)$, depending only on
$d,p,K$, such that, for every seal
$\mathscr S$ and every virgin seed $s$,
\[
 \PP\bigl(s\longleftrightarrow\infty\text{ in }N
       \mid\cF_{\mathscr S}\bigr)\ge c_{\rm esc}.
\]
There is a seal $\mathscr S_s$ on
$\{O_s\text{ is finite}\}$ such that
\[
 \Capt_{\overline N}(\mathscr S)
    \subseteq\Capt_{\overline N}(\mathscr S_s),
 \qquad
 \MeshSupp(\mathcal R_s)\subseteq \MeshInt(\mathscr S_s).
\]
Moreover,
\begin{equation}
\label{e.efficient-local-gain}
 \PP\bigl(G_s\ge r,\ O_s\text{ finite}\mid\cF_{\mathscr S}\bigr)
 \le Ce^{-cr^\alpha},\qquad r\ge1,
\end{equation}
and the enlargement has a boundary increment $Z_s$ satisfying
\eqref{e.efficient-boundary-growth} and
\begin{equation}
\label{e.efficient-boundary-tail}
 \PP\bigl(Z_s\ge u,\ O_s\text{ finite}\mid\cF_{\mathscr S}\bigr)
 \le Ce^{-cu},\qquad u\ge0.
\end{equation}
\end{proposition}

\noindent{\bf Proof.}
Condition on $\cF_{\mathscr S}$.  On each event
$\{\mathscr S=S,s=v\}$ with deterministic $S,v$, the exterior
product law follows from Lemma~\ref{l.freshness}, and
Proposition~\ref{p.uniform-open-escape-exterior-seed} gives the escape
bound.  Its constant is uniform in $S,v$, so the assertion holds under
the stated conditioning.

Assume $O_s$ is finite.  On $\mathsf{Loc}_s$, let $\mathcal H_s^{\rm loc}$
be the field generated by $\cF_{\mathscr S}$ and the finite cluster-query
record, and put
\[
       D_2^{\rm loc}=\MeshSupp(\mathcal R_s).
\]
The local event makes $|D_2^{\rm loc}|\le C(d)$, while the cluster exploration
is a stopping exploration and all newly queried edges lie in
$\bigcup_{M\in D_2^{\rm loc}}E(M)$.  Hence the hypotheses of
Lemma~\ref{l.conditional-rooted-completion} hold with
$\mathcal H=\mathcal H_s^{\rm loc}$ and $D_2=D_2^{\rm loc}$.
The comparison of newly required seal and halo tests in the proof of
Proposition~\ref{p.trap-cluster-to-seal}, together with
$|D_2^{\rm loc}|\le C(d)$, bounds their number by a constant depending
only on $d$.  If $H_s^{\rm loc}$ is the mesh
set adjoined by the completion, then
\[
  \PP(|H_s^{\rm loc}|\ge v,\,\mathsf{Loc}_s,\,O_s\text{ finite}
       \mid\cF_{\mathscr S})\le Ce^{-cv}.
\]
The output $\mathscr S_s$ satisfies
\[
       \MeshInt(\mathscr S)\cup D_2^{\rm loc}
       \subseteq \MeshInt(\mathscr S_s).
\]
Thus $\MeshSupp(\mathcal R_s)\subseteq \MeshInt(\mathscr S_s)$.
As in Proposition~\ref{p.trap-cluster-to-seal}, mesh-exterior
monotonicity and preservation of capture under reduction give
$\Capt_{\overline N}(\mathscr S)\subseteq
 \Capt_{\overline N}(\mathscr S_s)$.  The old stopped record, the cluster
queries and the completion queries are all captured by $\mathscr S_s$, so
the measurability argument of Proposition~\ref{p.trap-cluster-to-seal}
applies.  Thus $\mathscr S_s$ is a seal in the sense of
Definition~\ref{d.seal}.

Since $|D_2^{\rm loc}|\le C(d)$, Lemmas
\ref{l.union-of-captured-quilts} and~\ref{l.boundary-growth-under-addition}
give
\[
   |\mathscr S_s|\le|\mathscr S|+C(1+|H_s^{\rm loc}|).
\]
On this event set $Z_s=1+|H_s^{\rm loc}|$.  This gives
\eqref{e.efficient-boundary-growth} and the required exponential tail after
changing constants.  Also $G_s\le C(d,K)$.

On $\mathsf{Loc}_s^{\rm c}$, Proposition~\ref{p.finite-trap-on-nonescape}
gives $\mathcal T_s$.  Every vertex counted by $G_s$ belongs to the revealed
record outside the old capture and hence lies in an overlapping box captured
by $\mathcal T_s$.  Since a layer box has at most $C(d,K)$ vertices,
\[
   G_s\le C(d,K)|\Capt_{\overline N}(\mathcal T_s)|.
\]
Corollary~\ref{c.trap-tail} gives \eqref{e.efficient-local-gain}, while
Proposition~\ref{p.trap-cluster-to-seal} gives the stated seal placement and
boundary-increment estimate.  Combining the two cases proves the proposition.
\qed

The result just proved concerns an ordinary seal: it conditions on the
entire capture and produces another seal with the stopping-domain property.
For successive layer episodes, however, we wish to condition only on the
information actually revealed.  Edges within the enlarged capture may
remain unread, and their states must retain their joint product law with
those of the unread exterior.  The \emph{joint unread-field condition}
for a relative $\mathcal G$-seal $\mathscr S$ and a $\mathcal G$-measurable
set $\mathcal E^{\rm rd}\subseteq E(N)$ of revealed edges requires that
the states on $\mathcal E^{\rm rd}$ are $\mathcal G$-measurable,
$\mathcal U(\mathscr S)\cap\mathcal E^{\rm rd}=\varnothing$, and,
conditionally on $\mathcal G$,
\[
 \mathcal L\bigl((\omega(e))_{e\in E(N)\setminus\mathcal E^{\rm rd}}\mid\mathcal G\bigr)
 =\bigotimes_{e\in E(N)\setminus\mathcal E^{\rm rd}}{\rm Bernoulli}(p).
\]

\begin{proposition}[Relative open escape and seal enlargement]
\label{p.relative-seal-play}
Let a relative $\mathcal G$-seal $\mathscr S$ and a revealed-edge set
$\mathcal E^{\rm rd}$ satisfy the joint unread-field condition just
specified.  Put $\alpha=(d-2)/(d-1)$; the constants below may be chosen
as in Proposition~\ref{p.seed-escape-or-seal-extension}.  Let $s$ be
$\mathcal G$-measurable and suppose that an edge of
$\mathcal U(\mathscr S)$ is incident to $s$.  Then
\[
 \PP\bigl(s\longleftrightarrow\infty\text{ in }N\mid\mathcal G\bigr)
 \ge c_{\rm esc}.
\]
On $\{O_s\text{ is finite}\}$ the same construction gives an enlarged
conditioning field $\mathcal G_s$ and a $\mathcal G_s$-seal
$\mathscr S_s$ such that
\[
 \Capt_{\overline N}(\mathscr S)
    \subseteq\Capt_{\overline N}(\mathscr S_s),
 \qquad
 \MeshSupp(\mathcal R_s)\subseteq \MeshInt(\mathscr S_s).
\]
There is a non-negative integer $Z_s$ such that
\[
 |\mathscr S_s|\le |\mathscr S|+CZ_s,
\]
and
\[
 \PP(G_s\ge r,\ O_s\text{ finite}\mid\mathcal G)
    \le Ce^{-cr^\alpha},\qquad r\ge1,
\]
\[
 \PP(Z_s\ge u,\ O_s\text{ finite}\mid\mathcal G)
    \le Ce^{-cu},\qquad u\ge0.
\]
Let $\mathcal Q_s^{\rm all}$ be the set of previously unread edges queried
during the episode and seal construction.  The enlarged history and
revealed-edge set may be taken to be
\[
 \mathcal G_s=\mathcal G\vee
   \sigma(\mathcal Q_s^{\rm all},\omega|_{\mathcal Q_s^{\rm all}}),
 \qquad
 \mathcal E_s^{\rm rd}=\mathcal E^{\rm rd}\cup\mathcal Q_s^{\rm all}.
\]
On the finite outcome, $\mathcal E_s^{\rm rd}$ is $\mathcal G_s$-measurable,
its edge states are $\mathcal G_s$-measurable,
$\mathcal U(\mathscr S_s)\cap\mathcal E_s^{\rm rd}=\varnothing$, and the
states on $E(N)\setminus\mathcal E_s^{\rm rd}$ have their joint
Bernoulli-$p$ product law conditionally on $\mathcal G_s$.  Thus the
joint unread-field condition is preserved, including for edges left
unread inside the new capture.
\end{proposition}

\noindent{\bf Proof.}
Condition on $\mathcal G$.  The seal, the seed, its boundary-edge states,
and the revealed-edge set $\mathcal E^{\rm rd}$ are then fixed, while the
states on $E(N)\setminus\mathcal E^{\rm rd}$ are jointly Bernoulli-$p$.
The event $E_\infty^{\mathscr S}(L_s)$ in the proof of
Proposition~\ref{p.uniform-open-escape-exterior-seed} uses only seal edges
and edges of $\mathcal U(\mathscr S)$, all of the latter still unread.
It is increasing in these product edges, as is the event that all edges of
the launch boxlet are open.  The conditional Harris--FKG argument of that
proof therefore gives the same escape constant.

For a finite episode, follow the fixed-order cluster search, querying a
previously unread edge only when its state is needed, and on
$\mathsf{Loc}_s^{\rm c}$ construct the enclosing quilt as in
Proposition~\ref{p.finite-trap-on-nonescape}.  The latter construction
uses the old seal-edge states, which are $\mathcal G$-measurable, and the
joint product field on $\mathcal U(\mathscr S)$.  In particular, for an
exterior box the captured-side edges used in its tests lie in
$E(\mathscr S)$, while the remaining edges belong to
$\mathcal U(\mathscr S)$.  Consequently the inapt-component estimate
holds conditionally on $\mathcal G$.  The launch box $B_0$ used in the
proof of Corollary~\ref{c.trap-cardinality-tail} is determined by the
seal and seed.  The anchored-component counting argument of that
corollary, and then the enclosing-region bound of
Corollary~\ref{c.trap-tail}, carry over with the same constants.  The enclosing-quilt construction is unchanged.

Unlike the ordinary-seal search, the present exploration may query
an unread edge inside $\ECap(\mathscr S)$.  Such an edge need not
lie in $R(\mathcal T_s)$, the mesh region supplied by the new quilt.
The geometry in the proof of
Lemma~\ref{l.stopped-enclosing-record} nevertheless places
every tested box and every box meeting the cluster record in
$\Capt_{\overline N}(\mathscr S)\cup
 \Capt_{\overline N}(\mathcal T_s)$.
Thus every newly queried edge in the \emph{old exterior} belongs to
$\bigcup_{M\in R(\mathcal T_s)}E(M)$.  The remaining newly queried edges
lie in the old capture.  Since all seal-edge states are already known,
these latter queries lie strictly behind the old seal, in mesh boxes of
$D_1=\MeshInt(\mathscr S)$, by the mesh-edge partition and
\eqref{e.mesh-transfer-capture}.

Let $\mathcal H_s$ record $\mathcal G$ together with the identities and
states of the cluster and enclosing-quilt queries.  On the nonlocal finite
outcome, take $D_2=R(\mathcal T_s)$.  The searches are stopping
explorations of the conditional product field on
$E(N)\setminus\mathcal E^{\rm rd}$.
In particular, after conditioning on $\mathcal H_s$, the edges of
$\mathcal U(\mathscr S)$ outside $\bigcup_{M\in D_2}E(M)$ still
have their joint Bernoulli-$p$ product law: none of these exterior edges
was queried.  This is exactly the product-law hypothesis of
Lemma~\ref{l.conditional-rooted-completion}.  The comparison of newly
required tests in Proposition~\ref{p.trap-cluster-to-seal} is geometric
and gives $|J|\le C|\mathcal T_s|$.  Conditional completion therefore
produces a $\mathcal G_s$-seal, with completion excess having a uniform
exponential tail beyond $C|\mathcal T_s|$.  The boundary-growth bound
of that proposition and Corollary~\ref{c.trap-cardinality-tail} give the
stated exponential tail on taking
$Z_s=|\mathcal T_s|+|H_s^{\rm comp}|$, where $H_s^{\rm comp}$ is the
mesh set adjoined by completion.  Every vertex counted by $G_s$ lies in a
box captured by $\mathcal T_s$, so Corollary~\ref{c.trap-tail} gives
the asserted stretched-exponential tail for~$G_s$.

On $\mathsf{Loc}_s$, take instead
$D_2=\MeshSupp(\mathcal R_s)$ and let $\mathcal H_s$ record just
$\mathcal G$ and the cluster queries.  Here $|D_2|\le C(d)$,
and every exterior edge newly queried during the episode belongs to
$\bigcup_{M\in D_2}E(M)$.
The same stopped-query argument verifies the conditional product-law
hypothesis of completion.  As in the local case of
Proposition~\ref{p.seed-escape-or-seal-extension}, the number of newly
required tests is bounded, the completion excess $H_s^{\rm loc}$ has an
exponential tail, and $G_s\le C(d,K)$.  Here take
$Z_s=1+|H_s^{\rm loc}|$.  Combining the two cases proves the probability
estimates and the boundary-increment assertion.

In both cases completion starts from
$D_1\cup D_2$ and gives
$D_1\cup D_2\subseteq\MeshInt(\mathscr S_s)$.
Mesh-exterior monotonicity, the mesh-transfer formula and preservation of
capture under reduction give
$\Capt_{\overline N}(\mathscr S)\subseteq
 \Capt_{\overline N}(\mathscr S_s)$.
Newly queried edges in the old capture have mesh support in $D_1$,
while those in the old exterior have mesh support in $D_2$.  Hence
$\MeshSupp(\mathcal R_s)\subseteq\MeshInt(\mathscr S_s)$, and all
completion queries lie in the edge capture of $\mathscr S_s$.
The terminal tests also give the measurable boundary-edge states, as in
Lemma~\ref{l.conditional-rooted-completion}.

Finally, let $\mathcal Q_s^{\rm all}$ and $\mathcal E_s^{\rm rd}$ be as in
the statement.  The successive fixed-order queries, including those made
inside the old capture, form a stopping exploration of the joint product
field on $E(N)\setminus\mathcal E^{\rm rd}$.  On the finite outcome,
Lemma~\ref{l.freshness}, applied conditionally on $\mathcal G$, gives the
joint Bernoulli-$p$ product law on $E(N)\setminus\mathcal E_s^{\rm rd}$
after the history has been enlarged to $\mathcal G_s$.  The query set and
$\mathcal E_s^{\rm rd}$ are $\mathcal G_s$-measurable, and the states of
all edges in $\mathcal E_s^{\rm rd}$ are now known.  Every newly queried
edge lies in the new capture.  Also, the initial inclusion
$\mathcal U(\mathscr S)\cap\mathcal E^{\rm rd}=\varnothing$ places every
previously revealed edge in the old capture, and hence in the new one.
Thus $\mathcal U(\mathscr S_s)\cap\mathcal E_s^{\rm rd}=\varnothing$.
Together with the welcoming geometry and the measurable boundary-edge
states furnished by completion, this shows that $\mathscr S_s$ is a
$\mathcal G_s$-seal and proves the asserted preservation of the unread
product field. \qed

The relative episode already includes returns through the old capture and
encounters with its known woven fabric.  For iteration we also need to
control newly created mesh-boundary sites, not merely the positive change in
seal cardinality.  The same construction supplies this strengthening and a
joint bound on open growth and boundary cost.

\begin{lemma}[Episode growth and boundary costs]
\label{l.effective-layer-episode}
Assume the hypotheses of Proposition~\ref{p.relative-seal-play}, and run
its exhaustive virgin-seed episode from~$s$.  Put
$\alpha=(d-2)/(d-1)$.  Conditionally on $\mathcal G$, open escape in
$N$ has probability at least~$c_{\rm esc}$.  On the finite outcome, let
$\mathscr S_s$ and $\mathcal G_s$ be the enlarged relative seal and
history furnished by that proposition, and let $x=G_s$.
There are constants $C,\lambda>0$ and a non-negative integer-valued,
$\mathcal G_s$-measurable cost~$Z_s$ such that
\begin{equation}
\label{e.episode-new-boundary-sites}
 \bigl| (\mathscr S_s\cap\Mesh)\setminus(\mathscr S\cap\Mesh)\bigr|
 +\bigl(|\mathscr S_s|-|\mathscr S|\bigr)_+
 \le C Z_s,
\end{equation}
and
\[
  \EE\left[
    \exp\{\lambda(x^\alpha+Z_s)\}
       \one_{\{O_s\text{ finite}\}}
     \middle|\mathcal G\right]\le C.
\]
The joint unread-field condition is preserved for the enlarged history,
seal and revealed-edge record.  The constants are
uniform in the incoming history, seal and seed.
\end{lemma}

\noindent{\bf Proof.}
We retain the choice of $Z_s$ made in the proof of
Proposition~\ref{p.relative-seal-play}.  On the finite outcome with
$\mathsf{Loc}_s^{\rm c}$, write
\[
 D_1=\MeshInt(\mathscr S),\qquad D_2=R(\mathcal T_s),
 \qquad D_*=D_1\cup D_2\cup H_s^{\rm comp},
\]
where $H_s^{\rm comp}$ is the mesh set adjoined during completion.
Applying Lemma~\ref{l.union-of-captured-quilts} twice gives
\[
 A(D_*)\setminus A(D_1)
 \subseteq A(D_2)\cup A(H_s^{\rm comp}).
\]
By Lemma~\ref{l.boundary-enclosed-region},
$|A(D_2)|\le C|\mathcal T_s|$; bounded mesh degree gives
$|A(H_s^{\rm comp})|\le C|H_s^{\rm comp}|$.
Reduction only removes boundary sites, and
$\mathscr S\cap\Mesh=A(D_1)$.  Thus the number of new
mesh-boundary sites is at most
$C(|\mathcal T_s|+|H_s^{\rm comp}|)$.
The proof of Proposition~\ref{p.relative-seal-play} takes
$Z_s=|\mathcal T_s|+|H_s^{\rm comp}|$ and controls the positive
seal-cardinality increment by the same quantity.

On the finite outcome with $\mathsf{Loc}_s$, the set
$D_2=\MeshSupp(\mathcal R_s)$ has bounded size.  The same argument
gives both boundary bounds with $Z_s=1+|H_s^{\rm loc}|$, as in the
relative proposition's proof.
Its conditional exponential tail, together with the
stretched-exponential bound for $G_s$, yields the asserted joint
exponential moment by conditional Cauchy--Schwarz, after reducing
$\lambda$.  The escape and unread-field conclusions are those of
Proposition~\ref{p.relative-seal-play}. \qed

\section{Initial seals and paused exploration}

We now organise the estimates for individual episodes into a continuing
layer search.  An initial seal captures a prescribed finite record;
thereafter, captured scans alternate with episodes launched from virgin
seeds, and the seal is enlarged after each finite episode.

\subsection{The initial seal}
\label{ss.initial-seal}

Here we apply enclosure sealing to a finite open-and-closed record.  The mesh
support of this record is taken as the first finite mesh set in the same
completion algorithm.  Because the prescribed edges lie in this mesh set, the
boundary and halo tests are made in the remaining product field.

Let
$\mathfrak W=(D^{\rm op},D^{\rm cl})$, where $D^{\rm op}$ and
$D^{\rm cl}$ are disjoint finite edge sets in $N$ with
$|D^{\rm op}|+|D^{\rm cl}|\ge1$.  Write
$q=|D^{\rm op}|+|D^{\rm cl}|$, and introduce the
\emph{imposed-record event}
\[
 \operatorname{Imp}(\mathfrak W)
 =\{\omega(e)=1\ (e\in D^{\rm op}),\
      \omega(e)=0\ (e\in D^{\rm cl})\}.
\]

\phantomsection
\label{conv.empty-seal}
\noindent\emph{Empty-seal convention.}
For an empty prescribed record we take
$\Seal(\varnothing)=\varnothing$.  Then
$\MCap(\varnothing)=\ECap(\varnothing)=\varnothing$ and
$\mathcal U(\varnothing)=E(N)$; no sealing query is made.  The same
convention applies in Lemma~\ref{l.conditional-rooted-completion},
Proposition~\ref{p.relative-seal-play}, and Lemma~\ref{l.effective-layer-episode}
when the seal in force is empty and the full unread product field is available.

\begin{lemma}[Seal size controls captured volume]
\label{l.seal-cardinality-controls-volume}
With $\alpha=(d-2)/(d-1)$, for every finite seal~$\mathscr S$,
\[
 |\Capt_{\overline N}(\mathscr S)|
 \le C\bigl(1+|\mathscr S|\bigr)^{1/\alpha}.
\]
\end{lemma}

\noindent{\bf Proof.}
The region $\Capt_{\overline N}(\mathscr S)\setminus\mathscr S$ is the union of the bounded nearest-neighbour components of
$\Z^{d-1}\setminus\mathscr S$.  Apply
Lemma~\ref{l.column-isoperimetry} with $H=\mathscr S$ and then add
$|\mathscr S|$.
\qed

Because an initial seal may be placed after a layer history is already
known, we give the estimate in a conditional form.

\begin{proposition}[Initial-seal estimate]
\label{p.initial-seal-estimate}
Let $d\ge3$ and $\alpha=(d-2)/(d-1)$.  For all sufficiently large
$K$, there are $A,C,c>0$, depending only on $d$ and $p$, with the following
property.  Let $\mathcal H$ be a $\sigma$-field and let
$\mathfrak W=(D^{\rm op},D^{\rm cl})$ be an $\mathcal H$-measurable finite
record of disjoint edge sets in $N$, of size
$q=|D^{\rm op}|+|D^{\rm cl}|\ge1$.  Suppose that the prescribed states
hold almost surely and are recorded in~$\mathcal H$, and that, conditionally
on $\mathcal H$, the edge states
outside $D^{\rm op}\cup D^{\rm cl}$ have their joint Bernoulli-$p$
product law.

Put $Z=\MeshSupp(D^{\rm op}\cup D^{\rm cl})$, and let
$\Seal(\mathfrak W)=\Seal(Z)$ be the seal completion from~$Z$.  This completion terminates almost surely, reads no edge of
$\bigcup_{M\in Z}E(M)$, and gives
\[
 Z\subseteq\MeshInt(\Seal(\mathfrak W)).
\]
In particular, all prescribed edges are captured and none is a seal edge.
Let $\mathcal Q$ be the completion query set and put
\[
 \mathcal H^+=\mathcal H\vee\sigma(\mathcal Q,\omega|_{\mathcal Q}),
 \qquad
 \mathcal E_+^{\rm rd}=D^{\rm op}\cup D^{\rm cl}\cup\mathcal Q.
\]
Then $\Seal(\mathfrak W)$ is a seal relative to $\mathcal H^+$,
$\mathcal U(\Seal(\mathfrak W))\cap\mathcal E_+^{\rm rd}=\varnothing$,
and the states on $E(N)\setminus\mathcal E_+^{\rm rd}$ have their joint
Bernoulli-$p$ product law conditionally on $\mathcal H^+$.  All states
on $\mathcal E_+^{\rm rd}$ are $\mathcal H^+$-measurable.

If $V=|\Capt_{\overline N}(\Seal(\mathfrak W))|$, then
\begin{equation}
\label{e.initial-seal-tail}
 \PP\bigl(V\ge R\mid\mathcal H\bigr)
 \le Ce^{Cq}e^{-cR^\alpha},\qquad R\ge1,
\end{equation}
and the right-hand side is at most $1/2$ at
$R=Aq^{1/\alpha}$.  For a deterministic record~$\mathfrak W$, these conclusions hold
under $\PP(\cdot\mid\operatorname{Imp}(\mathfrak W))$ with trivial
initial $\mathcal H$; in this case $\Seal(\mathfrak W)$ is also
an ordinary seal.
\end{proposition}

\noindent{\bf Proof.}
Condition on~$\mathcal H$ and put
\[
 Z=\MeshSupp(D^{\rm op}\cup D^{\rm cl}),
 \qquad S^{(0)}=A(Z)^{\mathrm{ov}}.
\]
The mesh edge sets partition~$E(N)$
(Lemma~\ref{l.disjoint-twoK-mesh}), so $|Z|\le q$.
Every site of $A(Z)$ is $*$-adjacent to~$Z$, whence
$|A(Z)|\le Cq$.  By Lemma~\ref{l.frontier-addition-range}, every box in
$S^{(0)}\cup\Halo(S^{(0)})$ has mesh support contained in
$A(Z)^{[1]}$.  Since each mesh box belongs to the support of at most
$3^{d-1}$ overlapping boxes,
\[
 |S^{(0)}|+|\Halo(S^{(0)})|\le C|A(Z)|\le Cq.
\]
Take the empty seal as the old $\mathcal G$-seal, with $\mathcal G$
trivial, and apply Lemma~\ref{l.conditional-rooted-completion} with
$D_1=\varnothing$ and $D_2=D_0=Z$.  Its product-law hypothesis holds:
conditionally on~$\mathcal H$, the edge states outside
$\bigcup_{M\in Z}E(M)$ are independent Bernoulli-$p$ variables.
Since the old seal and its halo are empty, the initial testing obligations
are exactly $J=S^{(0)}\cup\Halo(S^{(0)})$.  With $H$ the mesh boxes
adjoined during completion, the conditional estimate gives
\begin{equation}
\label{e.initial-seal-weak-hull}
 \PP\bigl(|H|\ge Cq+u\mid\mathcal H\bigr)
 \le Ce^{-cu},\qquad u\ge0.
\end{equation}
The placement conclusion gives
$Z\subseteq\MeshInt(\Seal(\mathfrak W))$, so all prescribed edges
are captured.  Moreover, the separation property
\eqref{e.relative-completion-interior-separation}, with $D_2=Z$,
shows that neither boundary nor halo tests query any edge of
$\bigcup_{M\in Z}E(M)$.  In particular, the completion reads no
prescribed edge.

Lemma~\ref{l.conditional-rooted-completion} makes the reduction of
the terminal collection an $\mathcal H^+$-seal.  Every completion query is
in the terminal capture, as shown in the stopping-domain argument of
Proposition~\ref{p.weak-frontier-to-completed-seal}; the prescribed
edges also lie there.  Thus its exterior edge set is disjoint
from~$\mathcal E_+^{\rm rd}$.  Given~$\mathcal H$, the completion queries
are made adaptively in the product field outside the prescribed record.
The conditional form of Lemma~\ref{l.freshness} therefore gives the
joint Bernoulli-$p$ product law on every edge of
$E(N)\setminus\mathcal E_+^{\rm rd}$ given $\mathcal H^+$, including
those still unread inside the capture.  By construction, every state
on~$\mathcal E_+^{\rm rd}$ is $\mathcal H^+$-measurable.

For a deterministic record, $Z$ is fixed, and
Proposition~\ref{p.weak-frontier-to-completed-seal} gives that
$\Seal(Z)$ is an ordinary seal under~$\PP$.  Since completion
consults no edge of~$E(Z)$, the output is independent of
$\operatorname{Imp}(\mathfrak W)$.  The event
$\operatorname{Imp}(\mathfrak W)$ concerns only edges captured by the
output, so the ordinary-seal stopping-domain conclusion also holds
under the indicated conditional measure.

Let $D_*$ be the terminal mesh set, so that $D_*=Z\cup H$.
Lemma~\ref{l.boundary-growth-under-addition} gives
\[
 |A(D_*)^{\mathrm{ov}}|
 \le |A(Z)^{\mathrm{ov}}|+C|H|\le Cq+C|H|.
\]
Reduction only deletes boxes, and hence
\begin{equation}
\label{e.initial-seal-boundary-from-H}
 |\Seal(\mathfrak W)|\le Cq+C|H|.
\end{equation}
Lemma~\ref{l.seal-cardinality-controls-volume} implies
$V\le C(1+|\Seal(\mathfrak W)|)^{1/\alpha}$, so
$V\ge R$ entails
$|\Seal(\mathfrak W)|\ge cR^\alpha-C$.
When $R^\alpha$ exceeds a sufficiently large multiple of~$q$,
\eqref{e.initial-seal-weak-hull} and
\eqref{e.initial-seal-boundary-from-H} yield
$\PP(V\ge R\mid\mathcal H)\le Ce^{Cq}e^{-cR^\alpha}$.
For the remaining values of~$R$, this bound holds after enlarging
its prefactor, since $q\ge1$.  This proves
\eqref{e.initial-seal-tail}.  At $R=Aq^{1/\alpha}$ its right-hand
side is $C\exp\{(C-cA^\alpha)q\}$; choose~$A$ large enough.
\qed

\begin{corollary}[Boundary size of the initial seal]
\label{c.initial-seal-boundary}
In the setting of Proposition~\ref{p.initial-seal-estimate}, there are
$A_{\rm bd},C,c>0$, depending only on $d$ and $p$, such that, for
$u\ge0$,
\[
 \PP\bigl(|\Seal(\mathfrak W)|\ge A_{\rm bd}q+u
       \mid\mathcal H\bigr)
 \le Ce^{-cu}.
\]
For a deterministic record the same estimate holds under
$\PP(\cdot\mid\operatorname{Imp}(\mathfrak W))$.
\end{corollary}

\noindent{\bf Proof.}
Combine \eqref{e.initial-seal-weak-hull} with
\eqref{e.initial-seal-boundary-from-H}.
\qed

\subsection{Paused layer exploration}
\label{ss.paused-seed-exploration}

The continuing layer search takes as input a relative seal, its
conditioning history and revealed-edge record, and a finite set of supplied
seeds.  For a prescribed non-empty record~$\mathfrak W$,
Proposition~\ref{p.initial-seal-estimate} supplies the seal
$\Seal(\mathfrak W)$, history~$\mathcal H^+$ and revealed-edge
record~$\mathcal E_+^{\rm rd}$; the seeds are supplied separately.
For an empty record in a fresh layer, the seal and revealed-edge record
are empty, and no sealing query is made.

\begin{definition}[Paused layer exploration]
\label{d.paused-layer-exploration}
Let~$\mathscr S_0$ be a relative $\mathcal G$-seal and
$\mathcal E_0^{\rm rd}\subseteq E(N)$ a $\mathcal G$-measurable set of
revealed layer edges.  Suppose that these objects satisfy the joint
unread-field condition stated before
Proposition~\ref{p.relative-seal-play}.  Let~$P\subseteq V(N)$ be a
finite $\mathcal G$-measurable set of supplied seeds.  Fix deterministic orders
on vertices and edges.

The state consists of the reached vertex set~$V^{\rm rch}$, the set
$\mathcal Q$ of edges processed by the cluster search, the set
$\mathcal E^{\rm rd}$ of layer edges whose states are recorded, and
the current seal~$\mathscr S$.  Initially,
\[
 (V^{\rm rch},\mathcal Q,\mathcal E^{\rm rd},\mathscr S)
   =(P,\varnothing,\mathcal E_0^{\rm rd},\mathscr S_0).
\]
In particular, $\mathcal Q\subseteq\mathcal E^{\rm rd}$.
Queries made solely for sealing enlarge the revealed-edge record and
conditioning history, but not~$\mathcal Q$.  Whenever an open edge
is processed, its endpoints enter~$V^{\rm rch}$.

For the current state put
\[
\begin{aligned}
 \mathcal P_{\rm cap}
   &=\{e\in\ECap(\mathscr S)\setminus\mathcal Q:
                   e\cap V^{\rm rch}\ne\varnothing\},\\
 \mathcal P_{\rm ext}
   &=\{e\in\mathcal U(\mathscr S)\setminus\mathcal Q:
                   e\cap V^{\rm rch}\ne\varnothing\}.
\end{aligned}
\]
These are the pending captured and exterior edges.  Recompute them as the
state changes.  The seal update after a finite exterior episode is the
relative enlargement of Proposition~\ref{p.relative-seal-play}.
For such an episode, let~$s$ denote its seed, $e_s$ its launch edge,
and~$\mathcal R_s$ its query record.  The update uses these episode data together with the current conditioning
history, seal and revealed-edge record.
Put
\[
 D_1=\MeshInt(\mathscr S),\qquad
 D_2=\begin{cases}
       \MeshSupp(\mathcal R_s),&\mathsf{Loc}_s,\\
       R(\mathcal T_s),&\mathsf{Loc}_s^{\rm c},
     \end{cases}
 \qquad D_0=D_1\cup D_2.
\]
In the second case, $\mathcal T_s$ is the enclosing quilt constructed
from the episode record as in Section~\ref{ss.trap-formation}.  Run
conditional seal completion from~$D_0$ as in
Lemma~\ref{l.conditional-rooted-completion}, treating tests inherited
unchanged from~$\mathscr S$ as already passed.  Denote the resulting
relative seal, enlarged conditioning history and revealed-edge set by
$\mathscr S_s$, $\mathcal G_s$ and $\mathcal E_s^{\rm rd}$,
as in Proposition~\ref{p.relative-seal-play}.  The search follows
the cycle below.

\begin{enumerate}
\item If $\mathcal P_{\rm cap}$ is non-empty, process its first edge in
the fixed order, using its recorded state or querying it if unread.
Adjoin the edge to~$\mathcal Q$; if newly queried, adjoin it to
$\mathcal E^{\rm rd}$ and record its state in the conditioning
history.  If open, adjoin its other endpoint to~$V^{\rm rch}$.
Repeat until $\mathcal P_{\rm cap}=\varnothing$, recomputing the pending
sets after each update.  This iteration is the \emph{captured scan}.
Meeting a seal component activates its woven fabric; its known open
edges are processed by the same rule.

\item If $\mathcal P_{\rm cap}$ is empty and $\mathcal P_{\rm ext}$
is non-empty, choose the first reached vertex~$s$ incident to a pending
exterior edge, and then its first such edge~$e_s$.  Thus~$s$ is a
virgin seed.  Run the exhaustive episode of
Proposition~\ref{p.relative-seal-play}, processing~$e_s$ first and
including any returns through the capture.  Open escape ends the search.
On a finite outcome, perform the seal update just specified.  Adjoin
the episode's reached vertices and processed edges to~$V^{\rm rch}$
and~$\mathcal Q$, and set
$(\mathscr S,\mathcal E^{\rm rd})=(\mathscr S_s,\mathcal E_s^{\rm rd})$.
The current conditioning history becomes~$\mathcal G_s$.  Return to
step~(1).

\item If both pending sets are empty, return the reached side vertices
to the application, retaining the current seal, records and conditioning
history.  Later finite seed sets measurable from the retained history may be
adjoined to~$V^{\rm rch}$, with the search resumed from this state.
\end{enumerate}
\end{definition}

Enlargement changes the capture and hence the pending-edge sets: a pending
exterior edge may become captured and be processed in the next scan.

\begin{lemma}[Reduction to layer episodes]
\label{l.reduction-to-virgin-contacts}
The exterior searches in the paused layer exploration form an adapted
sequence of exhaustive virgin-seed episodes.  At the start of each episode,
the joint unread-field condition holds for the current history, seal and
revealed-edge record, and the conditional estimates of
Lemma~\ref{l.effective-layer-episode} apply.

Almost surely, there is no infinite succession of finite episodes.  If the
search finishes without open escape, all in-layer continuation from the
supplied seeds has been processed.
\end{lemma}

\noindent{\bf Proof.}
At each stage the displayed state determines the two pending-edge sets.
In a captured scan, a known edge state is used directly and an unread
captured edge is queried at a stopping time.  Conditional freshness
(Lemma~\ref{l.freshness}) preserves the joint product law on the
remaining unread edges.  The seal stays relative to the enlarged history:
its geometry and boundary states are known, and the captured scan makes
no exterior query.  Once the scan is exhausted, the seed and launch edge
selected from $\mathcal P_{\rm ext}$ are determined before the launch
edge is queried.  The hypotheses of
Proposition~\ref{p.relative-seal-play} therefore hold.

On a finite outcome, that proposition captures the entire episode query
record in an enlarged relative seal, preserving the joint unread-field law.
By exhaustiveness, every edge of the old exterior incident to the explored
cluster has been queried.  Since the new exterior is contained in the old
exterior and contains no queried edge, no vertex of that cluster is incident
to the new exterior.  With the enlarged seal in force, the captured scan resumes on the
pending edges it now captures.  The joint unread-field condition thus holds when the next exterior
episode begins, if one occurs.  When no pending edge remains, the in-layer
continuation is exhausted.

Finally, Lemma~\ref{l.effective-layer-episode} gives a conditional
probability of open escape of at least~$c_{\rm esc}>0$ at the start of
each exterior episode.  Thus the probability that at least
$n$ episodes occur and all of them have finite outcomes is at most
$(1-c_{\rm esc})^n$.  This rules out an infinite succession of finite
episodes; it does not constrain an application that continues supplying
new seed sets after completed searches. \qed

\section{The Encloser control problem}
\label{ss.one-life-percolation-control}

We close the chapter with a freestanding problem of adaptive control
that illustrates the use of seals to limit finite open growth.

The Encloser's terrain is an infinite layer $N=N_k$, equipped with
supercritical Bernoulli-$p$ bond percolation conditioned on the edges
in a given finite set $D^{\rm op}$ being open and those in a disjoint
finite set $D^{\rm cl}$ being closed.  The Encloser is informed only of
these edge states at the outset, so that to her the rest of the terrain
has the Bernoulli-$p$ law.

At any stage, including the outset, the Encloser may stop or plant
seeds at a non-empty finite set of vertices of~$N$.  The open components in the layer containing these seeds are then
explored.  If any is infinite, \emph{open escape} occurs and play ends
with no winnings.  Otherwise she learns these components and the closed
edges on their boundaries, and may stop or begin another round,
choosing further seeds in light of everything revealed so far.  She may seed anywhere in the
layer, including territory already explored, and may use auxiliary randomness independent of the percolation
configuration in making her choices.

If she stops, her winnings are the number of vertices in the union of
the finite open components discovered during play, each counted once.
If she never stops, her winnings are zero.  Write
$B_\pi$ for this outcome under an adaptive strategy~$\pi$.
How likely can she make it that her winnings exceed a given threshold?

Write $\mathfrak W=(D^{\rm op},D^{\rm cl})$ for the initial record, and put
\[
  q=|D^{\rm op}|+|D^{\rm cl}|,\qquad
  \PP_{\mathfrak W}(\,\cdot\,)
    =\PP(\,\cdot\mid\operatorname{Imp}(\mathfrak W)).
\]
For an empty record, $\PP_{\varnothing}=\PP$.

The initial information may offer the Encloser some favourable opening moves.
For example, take $D^{\rm op}=\varnothing$ and let $D^{\rm cl}$ consist
of the edges crossing the boundary of a large box in the layer.
By seeding every vertex in the box in a single round, she can secure
its entire volume without risking open escape.  For a box of transverse
side length~$L$, this gives winnings of order $L^{d-1}$ from a record
of size $q\asymp L^{d-2}$.
Thus a gain of order $q^{(d-1)/(d-2)}$ may be guaranteed by the
initial information.
We bound the probability of gains beyond an allowance on
this scale, uniformly over the record and the Encloser's strategy.

\begin{theorem}[Encloser control]
\label{t.one-life-percolation-control}
Let $d\ge3$, and recall that $\alpha=(d-2)/(d-1)$.  There are
$A,C,c>0$ such that, for every finite initial record of size~$q$ and
$r\ge0$,
\[
   \sup_{\pi}
   \PP_{\mathfrak W}\bigl(B_\pi\ge Aq^{1/\alpha}+r\bigr)
   \le Ce^{-cr^\alpha} \, ,
\]
where the supremum is taken over all adaptive strategies as just described.
The constants may depend on $d,p$ and the layer width, but are uniform
in the initial record, including the empty record, for which the offset
vanishes.
\end{theorem}

The box example exhibits complete shielding from open escape.
When gaps are left in the closed boundary, growth from prospective seeds
may be only partially shielded.  Bounding the tail of the Encloser's
winnings requires controlling this shielding effect as the record grows
under her adaptive choices.

We use seals for this purpose.  A seal encloses the discovered material
in a finite capture, whose entire volume we allow for in bounding her
winnings.  Outside the capture, the unread edges retain their product
law.  The screening argument of
Section~\ref{s.exterior-percolation-open-escape} replaces passages
through the capture in auxiliary routes to infinity by passages along
the seal's open fabric.  Consequently, a virgin-seed episode has a
uniformly positive chance of open escape, while an episode with finite
outcome can be absorbed by enlarging the seal at an exponentially
controlled boundary cost.

The next result controls the accumulation of these boundary costs
under adaptive play.

\begin{lemma}[Accumulated cost before escape]
\label{l.adaptive-sealed-chain}
Let $(\cF_n)_{n\geq0}$ be a filtration and let
$(Y_n)_{n\geq1}$ be an adapted sequence with values in $[0,\infty]$.  Put
\[
   \tau_{\rm esc}\, = \, \inf\{n\, \geq \, 1:Y_n\, = \, \infty\}\, ,
   \qquad
   X_n\, = \,
   \begin{cases}
      Y_n,&Y_n<\infty\, ,\\
      0,&Y_n\, = \, \infty\, .
   \end{cases}
\]
Fix $c_{\rm esc}\in(0,1]$, $C,c>0$ and $\alpha\in(0,1]$.
Suppose that, for every $n\ge1$, on $\{\tau_{\rm esc}>n-1\}$
and conditionally on $\cF_{n-1}$,
\[
\begin{split}
   \PP(Y_n=\infty\mid\cF_{n-1})&\, \geq \, c_{\rm esc}\, ,\\
   \PP(X_n\, \geq \, r,\ Y_n<\infty\mid\cF_{n-1})
      &\, \leq \, C e^{-cr^\alpha},\qquad r\, \geq \, 1\, ,
\end{split}
\]
There are $\lambda>0$ and $K_1<\infty$, depending
only on these constants, such that
\[
   \EE\!\left[
      \exp\!\left\{\lambda\sum_{1\le n<\tau_{\rm esc}}X_n^\alpha\right\}
      \middle|\cF_0\right]
\, \leq \, K_1\, .
\]
In particular, $\tau_{\rm esc}<\infty$ almost surely.
\end{lemma}

\noindent{\bf Proof.} Write $a=1-c_{\rm esc}<1$.  On
$\{\tau_{\rm esc}>n-1\}$, conditional integration of the finite-outcome
tail gives, for $0<\xi<c$,
\[
\begin{aligned}
\EE\!\left[\one_{\{Y_n<\infty\}}e^{\xi X_n^\alpha}
    \middle|\cF_{n-1}\right]
&\leq a+\xi\int_0^1 e^{\xi t}\,dt
  +C\xi\int_1^\infty e^{-(c-\xi)t}\,dt\\
&=a+(e^\xi-1)+\frac{C\xi e^{-(c-\xi)}}{c-\xi}
  =:\vartheta(\xi).
\end{aligned}
\]
Here the bound by one is used for $0\leq t<1$, and the assumed
tail bound for $t\geq1$.  Since $\vartheta(\xi)\to a<1$ as
$\xi\searrow0$, choose $\lambda>0$ such that
$\vartheta:=\vartheta(\lambda)<1$.  It follows that
\[
\EE\!\left[\one_{\{\tau_{\rm esc}>n\}}e^{\lambda\sum_{j=1}^{n}X_j^\alpha}
  \middle|\cF_{n-1}\right]
\leq\vartheta\,\one_{\{\tau_{\rm esc}>n-1\}}
  e^{\lambda\sum_{j=1}^{n-1}X_j^\alpha}.
\]
Iterating conditionally on $\cF_0$ gives
\[
 \EE\!\left[
   \one_{\{\tau_{\rm esc}>n\}}
   \exp\!\left\{\lambda\sum_{j=1}^{n}X_j^\alpha\right\}
   \middle|\cF_0\right]\leq\vartheta^n.
\]
In particular, $\PP(\tau_{\rm esc}>n\mid\cF_0)\leq\vartheta^n$,
so $\tau_{\rm esc}<\infty$ almost surely.  Partitioning according
to the first infinite outcome, and bounding the contribution from
$\{\tau_{\rm esc}=n+1\}$ by the corresponding contribution from
$\{\tau_{\rm esc}>n\}$, yields
\[
   \EE\!\left[
      \exp\!\left\{\lambda
        \sum_{1\le n<\tau_{\rm esc}}X_n^\alpha\right\}
      \middle|\cF_0\right]
\leq\sum_{n\geq0}\vartheta^n=(1-\vartheta)^{-1}.
\]
\qed

Under the reference dynamics to be introduced in
Chapter~\ref{s.model-problems-accrual}, cluster deposition in a given layer
proceeds along the slide schedule, until rejection of a proposed slide sends
the trajectory to the cemetery state.  Furnish this dynamics to the
Encloser, and she has an adaptive strategy, albeit a defective one: arrival
at the cemetery spells bankruptcy.  We do not formulate our analysis of
slide dynamics in terms of the Encloser problem, but the preceding lemma
will be used in Chapter~\ref{s.seal-sequence-layer} to control the
accumulated costs of fresh growth and seal enlargement.

\noindent{\bf Proof of Theorem~\ref{t.one-life-percolation-control}.}
Fix a strategy~$\pi$.  We construct seals to analyse its play.
Let $\mathscr S^0$ be the initial seal obtained from
Proposition~\ref{p.initial-seal-estimate} when~$q>0$, and set
$\mathscr S^0=\varnothing$ when~$q=0$.
Let $\mathcal H_0$ be the $\sigma$-field generated by the prescribed
record and the seal-construction transcript.  Corollary~\ref{c.initial-seal-boundary}
and the empty-seal convention give constants $A_0,C_0,c_0>0$ such that
\begin{equation}
\label{e.encloser-initial-boundary-tail}
 \PP_{\mathfrak W}\bigl(|\mathscr S^0|\ge A_0q+u\bigr)
 \le C_0e^{-c_0u},\qquad u\ge0.
\end{equation}
After this construction all still-unread edges have their joint
conditional Bernoulli-$p$ law.  Adjoin the Encloser's independent auxiliary randomness to the
exploration history at the outset.  The Encloser sees only the information
specified by the rules, which is contained in this enlarged history.
Her choices are thus adapted to it, and the conditional product law on
the unread edges is preserved.

Use the paused search of Section~\ref{ss.paused-seed-exploration}
to process the Encloser's seed sets.  By
Lemma~\ref{l.reduction-to-virgin-contacts}, this is an adapted succession
of exhaustive virgin-seed episodes, separated by captured scans.  Index
the enlarged exploration history, containing $\mathcal H_0$, by its
successive decision and edge-processing steps, writing $\mathcal G_t$
for the history at step~$t$ and $\mathcal G_\infty$ for its
limiting history.  Let $T_n$ be the
stopping time just after the seed and launch edge of the $n$\textsuperscript{th} episode
have been selected, but before that edge is queried; set $T_n=\infty$
if no such episode is launched.  Put $\cF_{n-1}=\mathcal G_{T_n}$.
These stopped fields are nested.  On $\{T_n<\infty\}$, the selected
seed and launch edge, the seal then in force, and the set of still-unread
layer edges are all $\cF_{n-1}$-measurable.  The event
$\{T_n=\infty\}$ is also $\cF_{n-1}$-measurable.  On a finite outcome let $Z_n$ be the boundary
cost furnished by
Lemma~\ref{l.effective-layer-episode}, and put $Y_n=Z_n$;
on open escape put $Y_n=\infty$.  If no $n$\textsuperscript{th} episode is undertaken,
put $Y_n=\infty$ and use this cemetery convention for all later
indices.  Define $X_n=Z_n$ on a finite episode and $X_n=0$ otherwise.
Each $Y_n$ is $\cF_n$-measurable, since the outcome is known
before the next launch, or from the limiting history if there is none.
The conditional escape bound and the exponential-moment
estimate for~$Z_n$ in Lemma~\ref{l.effective-layer-episode}
verify the hypotheses of Lemma~\ref{l.adaptive-sealed-chain}
with exponent one, uniformly in the past.  (In the cemetery case the
escape bound is automatic.)

Let $\tau_{\rm esc}=\inf\{n\ge1:Y_n=\infty\}$, and let~$\sigma$
be the number of finite episodes completed when the Encloser stops;
set $\sigma=\infty$ if she never stops.  If she stops without open escape, no further episode is launched.
The cemetery convention then gives $\tau_{\rm esc}=\sigma+1$ and
$\sum_{1\le n\le\sigma}Z_n=\sum_{1\le n<\tau_{\rm esc}}X_n$.
Lemma~\ref{l.adaptive-sealed-chain} and the tower property, since
$\mathcal H_0\subseteq\cF_0$, give $\lambda>0$ and $K_1<\infty$
such that
\begin{equation}
\label{e.encloser-accumulated-cost}
 \EE_{\mathfrak W}\!\left[
   \exp\!\left\{\lambda\sum_{1\le n<\tau_{\rm esc}}X_n\right\}
       \middle|\mathcal H_0\right]\le K_1.
\end{equation}

On the event that the Encloser stops without open escape, let
$\mathscr S^{\rm fin}$ denote the seal in force when she stops.
Every captured scan processes only edges in the current edge capture.
Each finite exterior episode has its full query record captured by the
enlarged seal, and captures are nested, by
Proposition~\ref{p.relative-seal-play}.  By the exhaustion assertion of
Lemma~\ref{l.reduction-to-virgin-contacts}, every layer edge incident
to a discovered component has been processed.  Thus all these edges
belong to the final captured edge set~$\ECap(\mathscr S^{\rm fin})$.
The mesh partition~\eqref{e.mesh-edge-partition} assigns each to a
captured mesh box containing both endpoints.  Since each such box has
at most~$C(d,K)$ vertices,
\[
 B_\pi\le C(d,K)|\MCap(\mathscr S^{\rm fin})|
       \le C(d,K)|\Capt_{\overline N}(\mathscr S^{\rm fin})|.
\]
Summing the positive seal-cardinality increments in
\eqref{e.episode-new-boundary-sites} over the finite episodes, with no
seal change during captured scans, gives
\[
 |\mathscr S^{\rm fin}|
   \le |\mathscr S^0|+C\sum_{1\le n\le\sigma}Z_n.
\]
Together with Lemma~\ref{l.seal-cardinality-controls-volume}, this
implies, on the same event,
\begin{equation}
\label{e.encloser-score-by-cost}
 B_\pi^\alpha
 \le C\left(1+|\mathscr S^0|+
                  \sum_{1\le n\le\sigma}Z_n\right).
\end{equation}

Put $U=(|\mathscr S^0|-A_0q)_+$.  The initial tail
\eqref{e.encloser-initial-boundary-tail} (and $U=0$ when $q=0$)
gives $\EE_{\mathfrak W}e^{\eta U}\le C$ for sufficiently small
$\eta>0$, uniformly in~$\mathfrak W$.  Taking $\eta$ smaller if
necessary, use \eqref{e.encloser-score-by-cost} and then condition on
$\mathcal H_0$ in \eqref{e.encloser-accumulated-cost} to obtain
\[
 \EE_{\mathfrak W}\!\left[
     \one_{\{\pi\text{ stops without escape}\}}
     e^{\eta B_\pi^\alpha}\right]
 \le Ce^{Cq}.
\]
For $R>0$, unsuccessful play contributes no probability to
$\{B_\pi\ge R\}$.  Exponential Markov gives
\[
 \PP_{\mathfrak W}(B_\pi\ge R)
   \le C\exp\{Cq-\eta R^\alpha\}.
\]
The same bound at $R=0$ is immediate after enlarging~$C$.
Finally, for $R=Aq^{1/\alpha}+r$,
\[
 R^\alpha\ge2^{\alpha-1}\bigl(A^\alpha q+r^\alpha\bigr).
\]
Choose $A$ large enough to absorb the term~$Cq$ and adjust the
remaining constants.  The bound is uniform in~$\pi$ and
$\mathfrak W$, proving the theorem. \qed

\chapter{Wavebreakers}
\label{s.wavebreakers}

\hfff{wavebreaker}{Wavebreakers} arrest the percolation explorations that form midpieces in slide
operations.  They are built from thin rows of test layers and are used against
explorations arriving from either side.  During the
IVF-to-LLF transition, a prospective fine height must withstand
northgoing waves launched by the succession of slides to its south.  When we
come to renewal survival in Chapter~\ref{s.renewal-survival}, renewals already
produced must instead be protected from southgoing waves launched by later
slides to their north.

The two breaker conditions reflect these uses.  At a very-fine height,
\({\bf B}^{+}_{\rm vf}\)
imposes a thin profile on the whole input cluster along a fixed
positive-offset breaker window north of the indexing height: in each of three rows, the cluster has bounded occupation in many separated short test layers.  The typical-band consequence of the global surface bound, Proposition~\ref{p.nobigcluster}, is used to show that, with overwhelming probability,  such thin whole-cluster windows occur with positive density.  As the schedule rises,
their thin profiles  are restored whenever they are tested, allowing the windows to be used repeatedly against
northgoing waves generated by active slides below them.  At a fine height,
\({\bf B}^{+}_{\rm f}\) imposes the weaker condition in which the height-indexed upcluster replaces the whole cluster.  The upcluster is retained, up to translation,
by a slide resample; when the slide at that fine height occurs, the thin profile
of its upcluster is used to prepare single-use breakers against later southgoing waves.  Very-fine breaker windows, which are reusable against northgoing waves, are employed in
the fine-delivery argument of
Chapters~\ref{s.seal-sequence-layer}--\ref{s.shadow-bounds-delivery-fine};
the single-use southgoing breakers prepared at fine heights are used to prove
renewal survival in Chapter~\ref{s.renewal-survival}.

In a single test layer with bounded occupation, finite passage can be
prevented in two ways.  An arrival with few vertices on the near side may be arrested when all the as-yet-unread layer edges incident to these vertices are closed.   For a larger arrival, place an initial seal
around the fixed record in the layer.  With fixed positive probability this
seal has bounded captured volume and hence covers only a bounded part of the
near side; one of the remaining arrival vertices is then a virgin seed, from
which open escape has uniformly positive probability.  Neither outcome is a
finite crossing of the layer.  Stacking many disjoint bounded-occupation layers turns this
fixed arrest probability into a strong wavebreaker estimate.

A breaker row used against northgoing waves may be compromised and yet restored
before its next test.  Suppose that a slide below the breaker deposits
midpiece material in its row, and that the schedule next performs an active
slide at a height still south of the row.  The new cut lies below the deposit
and therefore cuts through it.  By the cut-it-and-lose-it mechanism of
Chapter~\ref{s.dmu}, the deposited material then leaves the row.
Consequently, the fixed record present at the next test is contained in the
transported input whole-cluster record, so the input thin profile governs the
row anew when the next northgoing midpiece arrives.
  Proposition~\ref{p.breaker-persistence} is the
formal recovery statement.  Successive recovered breakers attenuate
northgoing waves launched far to the south and localise the advance-phase
accrual studied in
Chapters~\ref{s.seal-sequence-layer}--\ref{s.shadow-bounds-delivery-fine}.

The southgoing use begins at the slide which supplies the breaker.  Its new
midpiece may itself enter the breaker window, which consists of three
consecutive rows.  With suitably chosen row parameters, the lower row
typically arrests the midpiece.  This permits the middle and upper rows to be
reserved as separate single-use breakers against later waves from the north.
They protect renewals produced farther south as the schedule advances north.

Section~\ref{s.wavebreaker-local-estimates} develops the one-layer consequence of the TipToe dichotomy: the probability of growth that neither perishes on entry nor escapes openly is bounded away from one.  It then specialises the fixed-shift record to a test layer, isolating the fresh layer field to which the estimate will later be applied.
Section~\ref{s.wavebreaker-rows} assembles test layers into breaker rows and windows,
bounds the probability of finite passage through a row, and establishes
recovery of the very-fine profile before a northgoing test.  Section~\ref{ss.successive-northgoing-wavebreakers}
places, in the path of a northgoing exploration, first the breaker window
immediately north of the active slide and then the successive very-fine breaker
windows farther north, whose thin profiles are recovered when they are reached.
The exploration must pass these obstacles one by one, and their conditional passage bounds multiply.
Section~\ref{ss.breaker-stands-postpassage} shows that, in a slide resampling,
the lower row just north of the adjusted slide height usually arrests the
new midpiece before it reaches the middle and upper rows.  It then shows
that, when left untouched, these two rows furnish single-use breakers for
later southgoing tests.

Throughout we use the
fixed-shift single-slide law \(\Rlaw^v_{\uclust,\dclust}\) of
Section~\ref{s.onetime}.  The retained pair \((\uclust,\dclust)\), together
with the proposed head~\(v\), determines the discovered-open and
discovered-closed edge sets \(O_v\) and \(C_v\) of
Definition~\ref{d.prime}.  Their states are fixed, and every other edge is
supplied by the auxiliary Bernoulli field.

\section{Local layer estimates}
\label{s.wavebreaker-local-estimates}

When an exploration approaches a test layer, some of its edges are already
fixed by the retained clusters, while the remaining edges preserve their
Bernoulli law.  Stopping the exploration at the layer's boundary identifies
the vertices through which fresh growth may enter.  We first bound the
probability of finite growth from such an entrance set, and then identify
the fixed-shift records to which the estimate applies.
 
\subsection{Layer records and a one-layer estimate}
\label{ss.wavebreaker-crossing-test}

This subsection isolates the local estimate used at every breaker layer.  Its
input is a conditioning field which leaves part of the layer in product law,
together with a finite set of vertices through which an exploration may enter
that unread part.  The later sections obtain this data by stopping the
northgoing or southgoing fixed-shift exploration and verify the required
entrance property.

Fix a real number~\(a\), put \(b=a+2K\), and let \(N=\slab_{a,b}\), with the
layer-box structure \(\overline N=\overline N(a)\) of
Section~\ref{s.overlapping-boxes-boxlets}.  The southside and northside of~\(N\) are
\(\partial_-N=\vertbdry_a\) and \(\partial_+N=\vertbdry_b\).
Fix possible retained values \(\uclust,\dclust\) and a proposed head
\(v\in\Z^d\) with \(v\cvl\in[u,u+3r]\), as in
Subsection~\ref{ss.explicit-slide-law}, and work under
\(\Rlaw^v_{\uclust,\dclust}\).  For the slide height~\(h\) under
consideration, write \(\mathfrak u^v\) for the representative of~\(\uclust\)
whose head is~\(v\), and set
\begin{equation}
\label{e.fixed-shift-open-closed-data}
 \begin{split}
 O_v&=E(\mathfrak u^v)\cup E(\dclust),\\
 C_v
   &=\bigl(E(\downspace_h)\cap\partial_E(\dclust)\bigr)
     \cup
     \bigl(E(\upspace_{\hadj(v)})
             \cap\partial_E(\mathfrak u^v)\bigr).
 \end{split}
\end{equation}
Under \(\Rlaw^v_{\uclust,\dclust}\), the edges in \(O_v\) are open, those in
\(C_v\) are closed, and every edge outside \(O_v\cup C_v\) receives its state
from the auxiliary Bernoulli field~\(\omega'\), as in
Definition~\ref{d.prime}.
Fix deterministic orders on the vertices and edges of~\(N\), used only to
make measurable choices when several candidates are available.

The information already fixed inside~\(N\) is recorded at the level of both
vertices and edges.  A finite \emph{layer record} is a triple
\(\mathfrak W=(R,D^{\rm op},D^{\rm cl})\), where the \emph{recorded vertex
set} \(R\subseteq V(N)\) is finite and
\(D^{\rm op},D^{\rm cl}\subseteq E(N)\) are disjoint.  We require the
\emph{full-boundary property}
\begin{equation}
\label{e.layer-record-full-boundary}
 D^{\rm op}\cup D^{\rm cl}
   =\{e\in E(N):e\text{ is incident to }R\},
 \qquad
 D^{\rm op}\subseteq\{e\in E(N):e\text{ has both endpoints in }R\}.
\end{equation}
The sets \(D^{\rm op}\) and \(D^{\rm cl}\) record edges fixed open and closed.
Write \(\mathfrak W_{\rm E}=(D^{\rm op},D^{\rm cl})\) for the
edge part of the layer record.  Every layer edge incident to~\(R\) has
therefore been determined, and no edge in the complement
\(F(\mathfrak W)=E(N)\setminus(D^{\rm op}\cup D^{\rm cl})\) is incident
to~\(R\).  We also allow the empty record.

The edge record \(\mathfrak W_{\rm E}\) need not determine~\(R\): a vertex
outside~\(R\) whose layer neighbours all lie in~\(R\) can be added without
changing either edge set.  In each application, the
record is the restriction of the fixed-shift data
\eqref{e.fixed-shift-open-closed-data} to~\(N\):
\begin{equation}
\label{e.fixed-shift-layer-record}
 R=\bigl(V(\dclust)\cup V(\mathfrak u^v)\bigr)\cap V(N),\qquad
 D^{\rm op}=O_v\cap E(N),\qquad
 D^{\rm cl}=C_v\cap E(N).
\end{equation}
In applications, the layer's position ensures that only one of the two
clusters in \eqref{e.fixed-shift-layer-record} contributes vertices to~\(R\):
the shifted upcluster in northgoing ones and the downcluster in southgoing
ones.  Lemma~\ref{l.fixed-shift-layer-record} gives the corresponding
full-boundary criterion.  The entrance set is supplied later by the stopped
exploration.

Call \(\mathfrak W\) \emph{thin} when \(|R|\le C_1\), where \(C_1\) is the
local breaker threshold of Definition~\ref{d.fvf}.  From
\eqref{e.layer-record-full-boundary},
\begin{equation}
\label{e.thin-layer-record-size}
 |D^{\rm op}|+|D^{\rm cl}|\le 2d|R|\le 2dC_1 .
\end{equation}
The initial-seal estimate (Proposition~\ref{p.initial-seal-estimate})
applies to the edge record~\(\mathfrak W_{\rm E}\),
with constants depending only on the fixed local threshold.  Every vertex
of~\(R\) is incident to a layer edge, so a non-empty~\(R\) has a
non-empty edge record.

We next state the conditional experiment to which the estimate applies.  Let
\(\psi\in\{-,+\}\) name the side of entry: \(\psi=-\) for northgoing
travel and \(\psi=+\) for southgoing travel.  Let \(\mathcal H\) be a \(\sigma\)-field
containing the fixed-shift data; in the applications, it also records the
exploration stopped at the entry side of~\(N\).  Suppose that a layer
record~\(\mathfrak W\) and a finite set
\(A_\psi\subseteq\partial_\psi N\setminus R\) are
\(\mathcal H\)-measurable.  Assume that every edge in \(D^{\rm op}\) is
open and every edge in \(D^{\rm cl}\) is closed, with these states
\(\mathcal H\)-measurable, and that, conditionally on
\(\mathcal H\), the edge states on \(F(\mathfrak W)\) are independent
Bernoulli variables of parameter~\(p\).  We call \(A_\psi\) the
\emph{entrance set}.  At such a stop, the exploration has reached the
entry side while leaving~\(F(\mathfrak W)\) unread.  Its reached vertices
outside~\(R\) incident to these edges form~\(A_\psi\); layer edges
from~\(R\) into~\(V(N)\setminus R\) are already recorded closed.

Let \(\mathcal C_\psi=\mathcal C_N(A_\psi;\mathfrak W)\) be
the union of the components in
\[
 \bigl(V(N)\setminus R,
       \{e\in F(\mathfrak W):e\text{ is open}\}\bigr)
\]
which meet~\(A_\psi\).  For an empty entrance set, take
\(\mathcal C_\psi=\varnothing\).  Put
\begin{equation}
\label{e.local-layer-obstructions}
 \begin{split}
 \mathsf{Shut}_\psi
   &=\{\text{every edge of }F(\mathfrak W)
          \text{ incident to }A_\psi\text{ is closed}\},\\
 \mathsf{Esc}_\psi
   &=\{|V(\mathcal C_\psi)|=\infty\},\\
 \mathsf{FinGrow}_\psi
   &=\bigl(\mathsf{Shut}_\psi\cup\mathsf{Esc}_\psi\bigr)^{\rm c}.
 \end{split}
\end{equation}
The event \(\mathsf{Shut}_\psi\) blocks entry into the unread layer graph,
while \(\mathsf{Esc}_\psi\) gives an infinite open continuation wholly
inside~\(N\).  The complementary event \(\mathsf{FinGrow}_\psi\) says that
some unread entrance edge is open and that every component generated from
the entrance set is finite.  With the record and entrance set fixed,
these events depend only on the unread edge states.  We may reveal
those edges in an order convenient for the estimate.  After the stopped
fixed-shift search has been specified, Lemma~\ref{l.fresh-interface-verification}
identifies the entrance set at a finite pre-entry stop and places finite
passage through the layer in \(\mathsf{FinGrow}_\psi\).

Recall that \(\mathcal U(S)\) is the exterior edge set of a finite box
collection~\(S\).  For a finite mesh set~\(D\), consider its overlapping-box boundary
\(S=A(D)^{\mathrm{ov}}\) from Lemma~\ref{l.mesh-boundary-capture}.
For \(\psi\in\{-,+\}\), put
\begin{equation}
\label{e.covered-side-set}
 \Gamma_\psi(S)
\, = \, \bigl\{x\in\partial_\psi N:
       \text{no edge of }\mathcal U(S)\text{ is incident to }x\bigr\}\, .
\end{equation}
If \(x\notin\Gamma_\psi(S)\), then \(x\) is incident to an edge of
\(\mathcal U(S)\).  By Lemma~\ref{l.unread-edge-launch-cell}, the
boxlet assigned to the first such edge lies wholly in \(\mathcal U(S)\),
and every overlapping box containing it is exterior.

\begin{lemma}
\label{l.breaker-side-covered}
There is \(c_\Gamma=c_\Gamma(d,K)<\infty\) such that, for every finite
boundary collection \(S=A(D)^{\mathrm{ov}}\),
\[
 |\Gamma_-(S)|\vee|\Gamma_+(S)|
\, \le \, c_\Gamma\,|\Capt_{\overline N}(S)|\, .
\]
\end{lemma}

\noindent{\bf Proof.} Fix \(x\in\Gamma_\psi(S)\).  A side vertex is incident to a layer
edge, and every such edge lies in \(\ECap(S)\) by
\eqref{e.covered-side-set}.  Its unique mesh box therefore belongs to
\(\MCap(S)\) and contains \(x\) as an endpoint.  A mesh box has at most
\(c_\Gamma(d,K)\) vertices on either side of the layer.  Summing over the
captured mesh boxes gives \(|\Gamma_\psi(S)|\le c_\Gamma|\MCap(S)|\).
By \eqref{e.mesh-edge-partition},
\(\MCap(S)\subseteq\Capt_{\overline N}(S)\), proving the claim.
\qed

The initial-seal estimate enters the one-layer argument through a fixed
bounded-volume event.  For a layer record~\(\mathfrak W\), let
\(\Seal(\mathfrak W_{\rm E})\) be the initial seal of its edge record,
constructed as in Section~\ref{ss.initial-seal}.
For a deterministic record, Proposition~\ref{p.initial-seal-estimate}
gives the estimate under
\(\PP(\,\cdot\mid\operatorname{Imp}(\mathfrak W_{\rm E}))\),
the conditioning event that the prescribed edge states hold.  Under the
fixed-shift experiment these states are already known, and the proposition's
conditional-history formulation applies directly.  Write \(A_{\rm sc}\) for the scale constant in that proposition,
and choose a fixed integer
\(R_\star\ge A_{\rm sc}(2dC_1)^{1/\alpha}\).  The initial-seal
tail bound then gives, for every non-empty thin record,
\begin{equation}
\label{e.small-initial-seal-prob}
 \PP\bigl(
   |\Capt_{\overline N}(\Seal(\mathfrak W_{\rm E}))|\, \le \, R_\star
   \bigm|\operatorname{Imp}(\mathfrak W_{\rm E})
 \bigr)\, \ge \, \tfrac12\, .
\end{equation}
For the empty record take \(\Seal(\mathfrak W_{\rm E})=\emptyset\), so the display
holds trivially.  Under the conditional experiment just specified, the same
one-half lower bound holds after conditioning on~\(\mathcal H\): the
initial-seal construction
sees on \(F(\mathfrak W)\) exactly the product field used in
Proposition~\ref{p.initial-seal-estimate}.  Set
\(M_\star=c_\Gamma R_\star\) and
\(p_{\rm arr}=\min\{(1-p)^{2dM_\star},\tfrac12c_{\rm esc}\}>0\).

\begin{proposition}[One-layer finite-growth estimate]
\label{p.tiptoe-one-layer}
Let \(\mathcal H\) contain the fixed-shift data.  Let
\(\mathfrak W=(R,D^{\rm op},D^{\rm cl})\) be a thin
\(\mathcal H\)-measurable layer record, and let
\(A_\psi\subseteq\partial_\psi N\setminus R\) be finite and
\(\mathcal H\)-measurable.  Suppose that every edge in \(D^{\rm op}\) is open
and every edge in \(D^{\rm cl}\) is closed, with these states
\(\mathcal H\)-measurable, and that,
conditionally on~\(\mathcal H\), the edge states on~\(F(\mathfrak W)\)
have their Bernoulli product law.  Then, almost surely,
\[
 \Rlaw^v_{\uclust,\dclust}
   \bigl(\mathsf{FinGrow}_\psi\mid\mathcal H\bigr)
\, \le \, 1-p_{\rm arr}\, .
\]
Here \(p_{\rm arr}=p_{\rm arr}(d,p,K,C_1)>0\) is the constant fixed above.
\end{proposition}

The bound also holds on any \(\mathcal H\)-measurable event on which
its hypotheses hold: the proof applies conditionally on that event.

The proof uses the two sides of the TipToe dichotomy: \emph{thin,
instantly perish} in the small-entrance case and \emph{thick, openly escape}
in the large-entrance case.

\noindent{\bf Proof of Proposition~\ref{p.tiptoe-one-layer}.} Write
\(A=A_\psi\).  If \(A=\varnothing\), then
\(\mathsf{Shut}_\psi\) occurs and the assertion is immediate.

\emph{Small entrance set.}  Suppose \(|A|\le M_\star\), and let
\(F_A=\{e\in F(\mathfrak W):e\text{ is incident to }A\}\).
Conditionally on~\(\mathcal H\), the states of these edges are independent
Bernoulli variables of parameter~\(p\), and \(|F_A|\le2dM_\star\).  If every edge of~\(F_A\) is closed, then
\(\mathsf{FinGrow}_\psi\) does not occur.
Consequently
\[
 \Rlaw^v_{\uclust,\dclust}
   (\mathsf{FinGrow}_\psi^{\rm c}\mid\mathcal H)
 \ge (1-p)^{2dM_\star}.
\]

\emph{Large entrance set.}  Suppose \(|A|>M_\star\).  First reveal the
unread layer edges needed to construct the initial seal
\(\mathscr S=\Seal(\mathfrak W_{\rm E})\).
Let \(\mathcal Q\) be its completion query set, and put
\[
 \mathcal H^+=\mathcal H\vee\sigma(\mathcal Q,\omega|_{\mathcal Q}),
 \qquad
 \mathcal E_+^{\rm rd}=D^{\rm op}\cup D^{\rm cl}\cup\mathcal Q.
\]
Taking \(\mathscr S=\varnothing\) when \(\mathfrak W_{\rm E}\) is empty,
as prescribed by the \hyperref[conv.empty-seal]{empty-seal convention},
Proposition~\ref{p.initial-seal-estimate} gives an
\(\mathcal H^+\)-seal.  The states on \(\mathcal E_+^{\rm rd}\) are
\(\mathcal H^+\)-measurable, the seal exterior contains none of these
edges, and the edges outside \(\mathcal E_+^{\rm rd}\) retain their
joint Bernoulli product law conditionally on \(\mathcal H^+\).
These are the joint unread-field hypotheses of
Proposition~\ref{p.relative-seal-play}.

On the event \(\{|\Capt_{\overline N}(\mathscr S)|\le R_\star\}\),
the seal geometry \eqref{e.seal-mesh-geometry} permits the use of
Lemma~\ref{l.breaker-side-covered}, giving
\(|\Gamma_\psi(\mathscr S)|\le M_\star<|A|\).  Choose the first
\(s\in A\setminus\Gamma_\psi(\mathscr S)\).  This choice is
\(\mathcal H^+\)-measurable, and Lemma~\ref{l.unread-edge-launch-cell}
gives a fresh exterior boxlet at its first incident exterior edge.
Proposition~\ref{p.relative-seal-play} therefore gives, conditionally
on \(\mathcal H^+\), an open escape from~\(s\) in~\(N\) with probability
at least \(c_{\rm esc}\).

Since every layer edge from \(R\) to its complement is recorded closed
by \eqref{e.layer-record-full-boundary}, the open escape from
\(s\notin R\) stays outside \(R\) and uses only edges of
\(F(\mathfrak W)\).  It therefore witnesses \(\mathsf{Esc}_\psi\),
so \(\mathsf{FinGrow}_\psi\) does not occur.  The small-capture event
has conditional probability at least one half by
\eqref{e.small-initial-seal-prob}, so
\(\Rlaw^v_{\uclust,\dclust}
(\mathsf{FinGrow}_\psi^{\rm c}\mid\mathcal H)\ge c_{\rm esc}/2\)
in this case.

In either case the conditional probability of
\(\mathsf{FinGrow}_\psi^{\rm c}\) is at least \(p_{\rm arr}\), which is
the stated bound.
\qed

\begin{proposition}[Successive finite-growth events]
\label{p.tiptoe-successive}
Fix \(\psi\in\{-,+\}\), and let \(N_1,\ldots,N_m\) be layers
encountered in that direction.  Let \(\mathcal H_0\) be an initial
conditioning field, and put \(B_0=\Omega\).  For \(j\in\intint{m}\), let \(\mathcal H_j\supseteq\mathcal H_0\) be a
\(\sigma\)-field with \(B_{j-1}\in\mathcal H_j\).  Suppose that on
\(B_{j-1}\) there are a thin \(\mathcal H_j\)-measurable layer record
\(\mathfrak W_j=(R_j,D_j^{\rm op},D_j^{\rm cl})\) and a finite
\(\mathcal H_j\)-measurable entrance set
\(A_j\subseteq\partial_\psi N_j\setminus R_j\) satisfying the
conditional-law hypotheses of Proposition~\ref{p.tiptoe-one-layer}
for~\(N_j\).  Let \(\mathsf{FinGrow}_j\) be the event
\eqref{e.local-layer-obstructions} for this data, and let
\(B_j\subseteq B_{j-1}\cap\mathsf{FinGrow}_j\).  Then
\[
 \Rlaw^v_{\uclust,\dclust}(B_m\mid\mathcal H_0)
 \le(1-p_{\rm arr})^m
\]
almost surely.
\end{proposition}

\noindent{\bf Proof.} Proposition~\ref{p.tiptoe-one-layer} gives
\(\Rlaw^v_{\uclust,\dclust}(\mathsf{FinGrow}_j\mid\mathcal H_j)
\le 1-p_{\rm arr}\) on \(B_{j-1}\).  Since
\(B_{j-1}\in\mathcal H_j\), the tower property yields
\[
 \Rlaw^v_{\uclust,\dclust}(B_j\mid\mathcal H_0)
 \le (1-p_{\rm arr})
      \Rlaw^v_{\uclust,\dclust}(B_{j-1}\mid\mathcal H_0).
\]
Iteration proves the claim. \qed

\subsection{Fixed-shift interfaces}
\label{ss.breaker-fixed-shift-interfaces}

The next two lemmas record the elementary fixed-shift geometry used later.
First, an undiscovered edge can meet the shifted upcluster only along its bed.

\begin{lemma}[Bed attachment for the shifted upcluster]
\label{l.ssr-upcluster-bed-interface}
Fix admissible data \((\uclust,\dclust,v)\), and let
$e\notin E(\mathfrak u^v)$ be incident to a vertex of~$\mathfrak u^v$.  If $e$
is undiscovered in the fixed-shift partition of Definition~\ref{d.prime},
then its endpoint in~$\mathfrak u^v$ belongs to~$\Bed(\mathfrak u^v)$ and
\(e\notin E(\upspace_{\hadj(v)})\).
\end{lemma}

\noindent{\bf Proof.}
If $e\in E(\upspace_{\hadj(v)})$, then
\(e\in E(\upspace_{\hadj(v)})\cap\partial_E(\mathfrak u^v)\),
and Definition~\ref{d.prime} declares~$e$ discovered closed.  Thus an
undiscovered incident edge lies in \(E(\downspace_{\hadj(v)})\).
Every vertex of the shifted upcluster is incident to one of its edges
(the retained upcluster has positive height), and these edges belong to
\(E(\upspace_{\hadj(v)})\).  Thus the upcluster endpoint of~\(e\) lies in
\(\vertbdry_{\hadj(v)}\cap V(\mathfrak u^v)=\Bed(\mathfrak u^v)\),
as in Lemma~\ref{l.midpiece}.
\qed

Second, when a test layer lies on one side of the relevant slide cut and is clear of the opposite fixed graph, every layer edge incident to vertices of the retained cluster is recorded open or closed, while all other layer edges are undiscovered.

\begin{lemma}[Fixed-shift layer record]
\label{l.fixed-shift-layer-record}
Suppose that either
\(E(N)\subseteq E(\upspace_{\hadj(v)})\) and
\(V(N)\cap V(\dclust)=\emptyset\), or
\(E(N)\subseteq E(\downspace_h)\) and
\(V(N)\cap V(\mathfrak u^v)=\emptyset\).
Then the restriction \eqref{e.fixed-shift-layer-record} has the full-boundary
property \eqref{e.layer-record-full-boundary}, and every edge of
\(F(\mathfrak W)\) is undiscovered.
\end{lemma}

\noindent{\bf Proof.}
Since \(\hadj(v)\ge h\), the edge sets
\(E(\upspace_{\hadj(v)})\) and \(E(\downspace_h)\) are disjoint.
Thus in the first case the downcluster terms of
\eqref{e.fixed-shift-open-closed-data} contribute no layer edges, and the
recorded vertex set is \(R=V(\mathfrak u^v)\cap V(N)\).
Definition~\ref{d.prime} fixes its layer edges open and every other incident
layer edge closed.  In the second case the same argument applies to
\(\dclust\), with the upcluster terms absent.  In either case a fixed-open
layer edge has both endpoints in~\(R\), every layer edge incident to~\(R\)
is recorded, and every remaining layer edge lies outside \(O_v\cup C_v\),
hence is undiscovered.
\qed

\section{Breaker rows and their records}
\label{s.wavebreaker-rows}

For a finite exploration to cross a breaker row, it must make finite passage
through a succession of test layers.  The one-layer estimate can therefore
be applied repeatedly, provided sufficiently many of their fixed records
are thin.

\subsection{Rows of test layers}
\label{ss.row-of-thin-slabs}

The three rows of short slabs specified in Definition~\ref{d.fvf}
supply the test layers of a breaker window.
Recall from Section~\ref{ss.random-mesh} that
\(\overline H=(H_1,\ldots,H_L)\) is the random schedule of attempted slide
heights: the step at~\(H_i\) is active or dormant according as the fine
condition holds or fails.  Fix a realised schedule height~\(h=H_i\),
\(\rho\in\{0,1,2\}\) and \(j\in\intint{i_{\rm br}}\).

\begin{definition}[Breaker-window geometry]
\label{d.breaker-window-geometry}
The breaker profile geometry was fixed in
\eqref{e.breaker-profile-band}--\eqref{e.breaker-enlarged-interval}.
For the chosen \(h,\rho,j\), let
\[
 N_{h,\rho,j}
   =\slab_{a_{\rho,j}(h),\,a_{\rho,j}(h)+2K}
\]
be the physical test layer, equipped with the translated layer-box lattice
\(\overline N(a_{\rho,j}(h))\) of Section~\ref{s.overlapping-boxes-boxlets}.  The corresponding enlarged interval, recalled from
\eqref{e.breaker-enlarged-interval}, is
\[
 J^\sharp_{\rho,j}(h)
   =[a_{\rho,j}(h)-r,a_{\rho,j}(h)+2K].
\]
The \emph{breaker row} \(\mathscr B_h^\rho\) is the ordered list
\((N_{h,\rho,j}:j\in\intint{i_{\rm br}})\), and
\(\mathscr T_h^\rho=\bigcup_jN_{h,\rho,j}\) is its union graph.  The
\hfff{breaker-window}{breaker window \(\mathscr B_h\)} consists of the three
rows \(\mathscr B_h^0,\mathscr B_h^1,\mathscr B_h^2\), ordered from south
to north.  Consecutive test layers have southern sides separated by
\(g_{\rm br}\), and the whole window of height intervals
\(J_{\rho,j}(h)\) lies in \([\hadj+b_0,\hadj+b_0+C_0]\), where
\(C_0=3g_{\rm br}i_{\rm br}\).
\end{definition}

\begin{lemma}[Vertex range of a breaker test layer]
\label{l.breaker-layer-vertex-range}
For every \(h,\rho,j\),
\[
 V(N_{h,\rho,j})
 \subseteq
 \{x\in\Z^d:x\cvl\in J^\sharp_{\rho,j}(h)\}.
\]
Consequently, if the discovered-open graph
\(\dclust\cup\mathfrak u^v\) has at most \(C_1\) vertices
in \(J^\sharp_{\rho,j}(h)\), then its layer record in
\(N_{h,\rho,j}\) is thin.
\end{lemma}

\noindent{\bf Proof.} Write \(a=a_{\rho,j}(h)\).  If \(x\in V(N_{h,\rho,j})\), then \(x\) is an
endpoint of an edge~\(e\) whose greater endpoint height satisfies
\(a<h^+(e)\le a+2K\).  The height difference between the endpoints of a
nearest-neighbour edge is at most~\(r\).  Hence
\(a-r<x\cvl\le a+2K\), proving the inclusion.  For the second assertion,
the fixed vertex set of the layer record is
\[
 R=\bigl(V(\dclust)\cup V(\mathfrak u^v)\bigr)
      \cap V(N_{h,\rho,j}),
\]
so the inclusion gives
\(|R|\le C_1\). \qed

\begin{lemma}[Schedule separation of breaker windows]
\label{l.scheduled-breaker-separation}
Let \(H_i<H_j\) be distinct attempted heights.  Every breaker window
associated with~\(H_i\) lies strictly below~\(H_j\).  Consequently, the
breaker windows are pairwise disjoint and ordered from south to north.
\end{lemma}

\noindent{\bf Proof.} At an attempted height~\(H_i\), the adjusted height is at most~\(H_i+3r\).
Thus \eqref{e.breaker-profile-band} and the breaker-window geometry place the
whole breaker window below height~\(H_i+3r+b_0+C_0\).  On the other hand,
\(H_{i+1}-H_i\geq2K(g+1)\), while \eqref{e.random-mesh-gap} gives
\(g>3r+b_0+C_0\).  Hence the breaker window lies strictly below
\(H_{i+1}\), and therefore below~\(H_j\).  The window at~\(H_j\) lies
strictly above~\(H_j\), since its vertices have height greater than
\(\hadj+b_0-r>H_j\).  This proves the disjointness. \qed

The local occupation threshold~\(C_1\) was fixed in the very-fine
abundance argument through the requirement~\eqref{e.vf-abundance-choice}.
The constants
\(R_\star,M_\star,p_{\rm arr}\) chosen in Section~\ref{ss.wavebreaker-crossing-test}
therefore depend on this fixed local threshold but not on the row length
or layer separation.  The row length,
its admissible separation and the schedule gap are fixed compatibly, together
with the positive-depth allowance, in
Lemma~\ref{l.delivery-constant-compatibility}.  Lemma~\ref{l.nhigh} is uniform
in \(i_{\rm br}\), \(g_{\rm br}\), \(b_0\) and \(C_0\), so increasing
the row length does not reduce the proved density of very-fine heights.

For a breaker window in a fixed-shift experiment, use the restriction
\eqref{e.fixed-shift-layer-record} in every test layer~\(N_{h,\rho,j}\).
Thus its fixed vertex set is
\[
 R_{h,\rho,j}
 =\bigl(V(\dclust)\cup V(\mathfrak u^v)\bigr)
    \cap V(N_{h,\rho,j}),
\]
and its open and closed edge sets are \(O_v\cap E(N_{h,\rho,j})\) and
\(C_v\cap E(N_{h,\rho,j})\).
The alternatives of Lemma~\ref{l.fixed-shift-layer-record} ensure that
the layer meets at most one of these clusters.
Lemma~\ref{l.breaker-layer-vertex-range} shows that a
\(C_1\)-occupation bound in \(J^\sharp_{\rho,j}(h)\) makes this record thin.

\begin{figure}[htbp]
\centering
\includegraphics[width=0.75\textwidth]{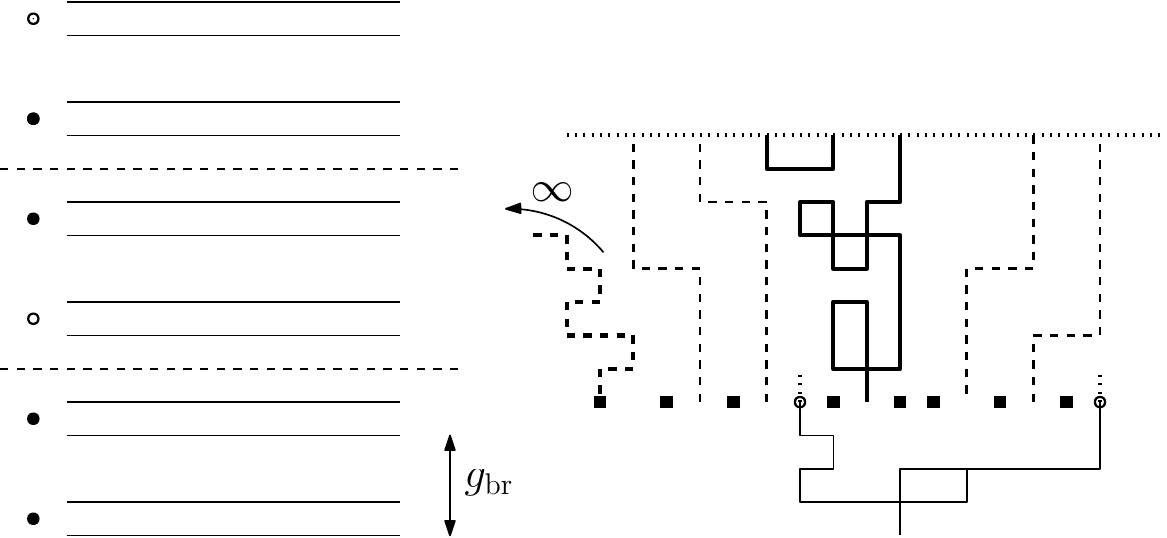}
\caption{The breaker window that lies directly to the north of an attempted slide height is schematically depicted in a planar sketch. {\em Left:} With $i_{\rm br} = 2$, three pairs of test layers form the lower, middle and upper rows. Filled discs mark initially thin layers. As the slide height rises towards the window, this initial thinness continues to govern the data fixed when a slide is attempted. {\rm Right}:
On such a slide, an exploration approaches the south side of a given thin test layer. The bold-marked discovered open material in the layer is enclosed in the dash-bordered bounded-cardinality seal. When the exploration approaches thinly, as the contact at the two circled locations depicts, the local closure of layer edges such as the vertical dotted pair bars entry into the layer. Should the exploration arrive thickly through the locations marked with squares, these points of entry will escape the bounded grasp of the southside seal boundary. The dashed path offers open escape in the layer from such a virgin seed.}\label{f.breaker}
\end{figure}

\begin{definition}[Admissible successive layer tests]
\label{d.admissible-breaker-tests}
Fix \(\psi\in\{-,+\}\), a conditioning field~\(\mathcal H_0\), and
\(\mathcal H_0\)-measurably selected layers \(N_1,\ldots,N_m\), ordered in
direction~\(\psi\).  Put \(E_0=\Omega\).  A family
\((\mathcal H_j,E_j)_{j=1}^m\) is \emph{admissible for these tests} if
\(\mathcal H_j\supseteq\mathcal H_0\), \(E_{j-1}\in\mathcal H_j\), and, on~\(E_{j-1}\), \(N_j\) has a thin
\(\mathcal H_j\)-measurable record and finite \(\mathcal H_j\)-measurable
entrance set satisfying Proposition~\ref{p.tiptoe-one-layer}, and
\(E_j\subseteq E_{j-1}\cap\mathsf{FinGrow}_j\).
\end{definition}

\begin{proposition}[Row finite-growth bound]
\label{p.breaker-passage}
For \(m=i_{\rm br}/2\) selected layers of a row, suppose that
\((\mathcal H_j,E_j)_{j\in\intint{m}}\) is admissible for these tests,
with initial conditioning field~\(\mathcal H_0\).  Then
\begin{equation}
\label{e.q-row-definition}
 \Rlaw^v_{\uclust,\dclust}(E_m\mid\mathcal H_0)
 \le q_{\rm row}:=(1-p_{\rm arr})^{i_{\rm br}/2}.
\end{equation}
The same bound holds for every event contained in~\(E_m\).
\end{proposition}

\noindent{\bf Proof.} This is Proposition~\ref{p.tiptoe-successive} with
\(B_j=E_j\). \qed

Since \(p_{\rm arr}>0\) is independent of the row length,
\(q_{\rm row}\) tends to zero exponentially as the even integer
\(i_{\rm br}\) tends to infinity.  The separation \(g_{\rm br}\ge2K\)
and offset \(b_0>r\) may be kept fixed, with
\(C_0=3g_{\rm br}i_{\rm br}\).

At a fine height, \({\bf B}^{+}_{\rm f}\) supplies the required
thin-layer majority in each of its three rows that Figure~\ref{f.breaker} illustrates.
Proposition~\ref{p.breaker-persistence} below supplies the same majority when a
very-fine window is recovered for a later northgoing test.

\subsection{Recovery before a northgoing test}
\label{ss.breaker-tested-persistence}

For a very-fine height~\(h\), Definition~\ref{d.fvf} bounds the whole input
cluster in the enlarged intervals of at least half the test layers in each
row.  As the schedule advances towards the window, earlier midpiece deposits
are cut away, while the recorded vertex sets in the window come from a
translate of the input cluster.  The next two propositions show that the
thin-layer majority is recovered when a later northgoing exploration tests
the window.

Since the test layers move with the adjusted height, we use the
state-explicit notation \(N_{h,\rho,j}(S)\) and
\(J^\sharp_{\rho,j}(S,h)\) for the physical test
layer and enlarged interval obtained by evaluating \(\hadj\) in the
state~\(S\).  For a proposed head~\(v\), write
\(N_{h,\rho,j}(v)\) and \(J^\sharp_{\rho,j}(v,h)\) for the same objects
with adjusted height \(h+(v\cvl-u)\).  We use the same convention for
\(\mathscr B_h(v)\) and \(\mathscr T_h^\rho(v)\).

Fix a scheduled height~\(h\), and suppose that the active heights preceding a
northgoing test of its breaker window are
\(\ell_1<\cdots<\ell_N<h\).  Let \(S_0\) be the state before the slide at
\(\ell_1\), and, for \(1\leq i<N\), let \(S_i\) be the state after the slide
at~\(\ell_i\).  At~\(\ell_N\), fix an admissible proposed head~\(v\), and
let \(\mathfrak u^v\) be the translate of
\(\upcluster_{\ell_N}(S_{N-1})\) with head~\(v\).  Put
\(\tau=v-\head(S_0)\), and let \(R^{\rm pre}_{h,\rho,j}\) denote the
recorded vertex set in \(N_{h,\rho,j}(v)\), with ``pre'' indicating
that the layer record is determined before exploration of that layer.

\begin{proposition}[Transported breaker content]
\label{p.transported-breaker-content}
For every \(\rho\in\{0,1,2\}\) and \(j\in\intint{i_{\rm br}}\),
\begin{equation}
\label{e.transported-breaker-trace}
 R^{\rm pre}_{h,\rho,j}
 =V(\mathfrak u^v)\cap V\bigl(N_{h,\rho,j}(v)\bigr)
 =\tau\cdot\Biggl(
    V\bigl(\upcluster_{\ell_N}(S_0)\bigr)
       \cap V\bigl(N_{h,\rho,j}(S_0)\bigr)
   \Biggr).
\end{equation}
\end{proposition}

\noindent{\bf Proof.} The fixed graphs are
\(\downcluster_{\ell_N}(S_{N-1})\) and \(\mathfrak u^v\).
The downcluster lies in~\(\downspace_{\ell_N}\), so all its vertices have
height at most~\(\ell_N\).  Lemma~\ref{l.breaker-layer-vertex-range},
\(v\cvl\geq u\), and \(b_0>r\) place every vertex of
\(N_{h,\rho,j}(v)\) strictly above \(h>\ell_N\).
Thus the recorded vertex set in \(N_{h,\rho,j}(v)\) consists exactly
of the vertices of \(\mathfrak u^v\) in that layer, proving the first
equality in
\eqref{e.transported-breaker-trace}.

At each active slide below~\(\ell_N\), Lemma~\ref{l.higher-ucs-shift}
translates \(\upcluster_{\ell_N}\) by the selected shift.  Composing these
translations with the proposed translation at~\(\ell_N\) carries
\(\head(S_0)\) to~\(v\), so the total translation is by~\(\tau\):
\[
 \mathfrak u^v
   =\tau\cdot\upcluster_{\ell_N}(S_0).
\]
The two adjusted heights differ by~\(\tau\cvl\), since
\(h+(v\cvl-u)=\hadj(S_0,h)+\tau\cvl\).
Thus translation by~\(\tau\) carries \(N_{h,\rho,j}(S_0)\) onto
\(N_{h,\rho,j}(v)\).  Intersecting the two translated vertex sets proves
the second equality. \qed

The result just proved also explains the whole-cluster strengthening in
\({\bf B}^{+}_{\rm vf}\).  When translated back to~\(S_0\), the recorded
vertex set of a test layer is contained in \(\upcluster_{\ell_N}(S_0)\),
but need not be contained in \(\upcluster_h(S_0)\); since \(\ell_N<h\),
Proposition~\ref{p.dmutwoheights}(2) gives
\(\upcluster_h(S_0)\subseteq\upcluster_{\ell_N}(S_0)\), generally strictly.
The very-fine condition controls the ambient input cluster, and hence these
recorded vertex sets in the required majority of layers.

\begin{proposition}[Recovery at a tested breaker]
\label{p.breaker-persistence}
In the setting described before Proposition~\ref{p.transported-breaker-content}, suppose
\(h\) satisfies \({\bf B}^{+}_{\rm vf}(h)\) in~\(S_0\).  At every
such northgoing test and in each row \(\rho\in\{0,1,2\}\),
\(|R^{\rm pre}_{h,\rho,j}|\le C_1\) for at least \(i_{\rm br}/2\) indices~\(j\);
thus the corresponding pre-test records are thin.
\end{proposition}

\noindent{\bf Proof.} Fix a row~\(\rho\).  By \({\bf B}^{+}_{\rm vf}(h)\), at least
\(i_{\rm br}/2\) indices~\(j\) satisfy
\[
 \bigl|V(\csetuthree(S_0))\cap
       \{x:x\cvl\in J^\sharp_{\rho,j}(S_0,h)\}\bigr|
 \leq C_1.
\]
For each such~\(j\), Proposition~\ref{p.transported-breaker-content} transports the
vertex set of the pre-test layer record back into
\(V(\upcluster_{\ell_N}(S_0))\cap V(N_{h,\rho,j}(S_0))\).  Since
\(\upcluster_{\ell_N}(S_0)\subseteq\csetuthree(S_0)\), Lemma
\ref{l.breaker-layer-vertex-range} gives
\(|R^{\rm pre}_{h,\rho,j}|\leq C_1\).  Thus these records are thin.
\qed

\section{Successive northgoing wavebreakers}
\label{ss.successive-northgoing-wavebreakers}

A northgoing exploration that makes finite passage through one breaker window
may endure far enough to meet the next.  This section places the breaker window
just north of the present slide and the later recovered breaker windows in one
exhaustive fixed-shift exploration, so that passage beyond the sequence requires
passage through its breakers one by one and the corresponding stopped bounds
multiply.

\subsection{Fixed-shift search and layer entry}
\label{ss.successive-stopped-tests}

Fix an attempted-schedule index~$q_0$, put \(h_0=H_{q_0}\), suppose that the
slide at~\(h_0\) is active, and fix admissible values
\((\uclust,\dclust,v)\) of its fixed-shift data.

Under~\(\Rlaw^v_{\uclust,\dclust}\), we first specify the cluster search
to which all the stopped tests belong.
The object being determined is precisely the cluster~\(\mc C(0)\) in the
reconstructed fixed-shift configuration.  The retained downcluster is already
known open and connected to the origin, so there is no reason to rediscover it
edge by edge: it is traversed at time zero, and its unknown departures issue
from~\(\Top(\dclust)\).  If the two retained graphs already meet, the shifted
upcluster is activated at time zero; otherwise it is traversed only after the
search first reaches one of its vertices.  From then on its discovered-open
edges are known, while Lemma~\ref{l.ssr-upcluster-bed-interface} places every
undiscovered departure from it at its bed.

Put
\[
 \mathcal E_v^\circ=E(\mathbb Z^d)\setminus(O_v\cup C_v)
\]
for the fixed-shift undiscovered edge set.  The shifted upcluster needs no edge-by-edge
exploration: whenever the reached set first meets~$V(\mathfrak u^v)$, all of
$V(\mathfrak u^v)$ is declared reached at once, since every edge of
$\mathfrak u^v$ is already known to be open.  This reveals no edge of~$\mathcal E_v^\circ$.

Fix deterministic orders on vertices and edges.  These orders are used only
to resolve choices among simultaneously available queries.

\begin{definition}[Fixed-shift cluster exploration]
\label{d.fixed-shift-origin-search}
After~\(n\) queries of edges in~\(\mathcal E_v^\circ\), let
\(\mathcal V_n\) be the reached vertex set and
\(\mathcal Q_n\subseteq\mathcal E_v^\circ\) the set of edges already
queried.  Put~$\mathcal Q_0=\emptyset$ and begin with
$\mathcal V_0=V(\dclust)$.  If
$V(\dclust)\cap V(\mathfrak u^v)\neq\emptyset$, also adjoin all of
$V(\mathfrak u^v)$ to~$\mathcal V_0$.  Let the current pending list be the
deterministically ordered set
\begin{equation}
\label{e.fixed-shift-pending-list}
 \mathcal L_n=
 \{e\in\mathcal E_v^\circ\setminus\mathcal Q_n:
          e\text{ is incident to }\mathcal V_n\}.
\end{equation}
When \(\mathcal L_n\neq\emptyset\), a scheduling rule chooses the next
edge~\(e_n\) from~\(\mathcal L_n\), measurably from the information already
revealed.  Set
\(\mathcal Q_{n+1}=\mathcal Q_n\cup\{e_n\}\).  If~\(e_n\) is closed, put
\(\mathcal V_{n+1}=\mathcal V_n\).  If~\(e_n\) is open, first put
\(\mathcal V_{n+1}=\mathcal V_n\cup V(e_n)\); if this set meets
$V(\mathfrak u^v)$, adjoin all of~$V(\mathfrak u^v)$ to it as well.
The list \(\mathcal L_{n+1}\) is then given again by
\eqref{e.fixed-shift-pending-list}.  Edges in~\(O_v\) are already known
open: the downcluster is traversed at initialisation and the shifted
upcluster when first reached.  Edges in~\(C_v\) are already known closed.
Thus neither class requires a query.

Let \(\cF_n\) be the \(\sigma\)-field generated by the fixed-shift data and
the responses to the first~\(n\) queries of edges in~\(\mathcal E_v^\circ\).  A scheduling is
\emph{exhaustive} if no pending edge is left unqueried forever.  More locally,
a stage with eligible edge set~\(D\) is exhaustive on~\(D\) when it chooses
queries from \(\mathcal L_n\cap D\) and can end finitely only when this set is
empty.
\end{definition}

By Definition~\ref{d.prime}, undiscovered departures from the
downcluster issue from~\(\Top(\dclust)\); by
Lemma~\ref{l.ssr-upcluster-bed-interface}, those from the shifted upcluster
issue from its bed.  Thus \(\mathcal L_0\) consists of the former together
with the latter if the two graphs already meet.  Otherwise its remaining unqueried bed departures enter the pending list
when the shifted upcluster is first reached.

\begin{lemma}[Completeness of the fixed-shift search]
\label{l.fixed-shift-search-completeness}
If \(\mathcal L_n=\emptyset\) at a finite query time~\(n\), then
\(\mathcal V_n\) is the vertex set of the fixed-shift open cluster of the
origin, and the fixed data together with the queried edges of~$\mathcal E_v^\circ$ determine
that cluster and its incident closed boundary.  Conversely, if this origin
cluster is finite, every exhaustive scheduling of the search reaches such a
time after finitely many queries.
\end{lemma}

\noindent{\bf Proof.}
Every vertex inserted into~\(\mathcal V_n\) is joined to the origin by fixed
open edges and queried-open edges of~$\mathcal E_v^\circ$, so \(\mathcal V_n\) is always
contained in the reconstructed cluster.  If \(\mathcal L_n=\emptyset\), every edge of~$\mathcal E_v^\circ$
incident to~\(\mathcal V_n\) has been queried and found either open, in
which case its endpoints have been inserted, or closed.  The fixed boundary
edges in~\(C_v\) are already closed.  Hence no open edge leaves
\(\mathcal V_n\), proving equality with the reconstructed cluster and determining
its tested boundary.  If the reconstructed cluster is finite, only finitely many edges of~$\mathcal E_v^\circ$ are incident to it.  An exhaustive scheduling processes all of
them after finitely many queries, when the pending list is empty. \qed

In Section~\ref{ss.multi-layer-paused}, we will develop a more
elaborate exploration that follows the cluster through several monitored
layers, pausing to construct or enlarge seals as needed.  This introduces
auxiliary seal queries and changes the scheduling, but retains the
reachability rule of the present search.  Lemma~\ref{l.multi-layer-regional-completeness}
shows that the completed cluster record agrees with that of this unpaused
search.

For the breaker tests we use height fronts to schedule this search.  For a
real~\(c\), write
\begin{equation}
\label{e.height-front-pending-lists}
 \mathcal L_n^-(c)=\mathcal L_n\cap E(\downspace_c),
 \qquad
 \mathcal L_n^+(c)=\mathcal L_n\cap E(\upspace_c).
\end{equation}
A northgoing front at height~\(c\) repeatedly queries the first edge of
\(\mathcal L_n^-(c)\) until this set is empty; a southgoing front does the
same with~\(\mathcal L_n^+(c)\).  If the stage begins at query count~\(s\),
its completion time is
\(\inf\{n\ge s: \mathcal L_n^-(c)=\emptyset\}\) northgoing, or the
same expression with~\(\mathcal L_n^+(c)\) southgoing; it is~\(\infty\)
if the stage is not reached or does not finish.  As the front advances, newly discovered
open connections may lead back into heights already passed.  Any pending
edges exposed there are processed before the front advances again, as part
of the same search.

For a test layer \(N=\slab_{a,b}\), a northgoing pre-entry stop is the finite
completion of the front at~\(a\), before the front is raised to~\(b\).
Every edge of~\(N\) has greater endpoint height in~\((a,b]\), so no
layer edge is eligible at this stage.  Raising the front to~\(b\) admits the
fresh layer edges; finite passage means that this exhaustion terminates and
a vertex of~\(\partial_+N\setminus R\) is reached.  The southgoing convention is
the reverse one: exhaust the front at~\(b\), then lower it to~\(a\), and use
\(\partial_-N\) as the far side.  Between successive test layers the front is
moved across the intervening gap and the newly eligible pending edges are
exhausted before the next layer is admitted.  For the southgoing pre-entry
stage, likewise~\(E(N)\cap E(\upspace_b)=\emptyset\).  As neither
protocol admits the layer earlier, Lemma~\ref{l.fixed-shift-layer-record}
ensures that~\(F(\mathfrak W)\) remains unread at either pre-entry stop.

We record a consequence of front exhaustion that is used below, and
again in Chapter~10 when we show that reaching a monitored layer forces
passage through the intervening breakers.
Suppose the northgoing front at~$c>\hadj(v)$ finishes and the reconstructed
origin cluster contains a vertex~$x\notin V(\mathfrak u^v)$ of height greater
than~$c$.  At that stop some reached vertex
$z\in\vertbdry_c\setminus V(\mathfrak u^v)$ is incident to a still-unqueried
open edge of~$\mathcal E_v^\circ$ with greater endpoint height above~$c$.
Every reached vertex outside the shifted upcluster has height at most~$c$,
so~$x$ is not yet reached.  Follow an open path from the downcluster to~$x$
until it first leaves the reached set at the stop.  Its departing edge is not fixed open,
since the retained clusters are traversed on contact, and is not eligible
below the front, since those pending edges have been exhausted.  It is
therefore an unread open edge with greater endpoint above~$c$.
Lemma~\ref{l.ssr-upcluster-bed-interface} rules out an unread departure
from the shifted upcluster above~$c$.  The departing edge consequently crosses
height~$c$ from a reached vertex~$z$ on~$\vertbdry_c$ outside the
shifted upcluster.  At the southern pre-entry front of a breaker test
layer, this edge belongs to its unread field (the layer has width~$2K\ge r$),
so it gives the required entrance contact.  At the northern completion
front it gives a reached vertex on the far side, outside the shifted
upcluster.  In particular, once the front has been exhausted at a layer's
entry side, any first fresh entry into its unread part must occur through
a vertex on that side already reached by the exploration.

\begin{lemma}[Layer entrance at a finite stop]
\label{l.fresh-interface-verification}
Let \(N=\slab_{a,b}\) satisfy one of the two alternatives of
Lemma~\ref{l.fixed-shift-layer-record}.  Take \(\psi=-\) in the first case
and \(\psi=+\) in the second, and put \(c=a\) for \(\psi=-\) and \(c=b\)
for \(\psi=+\).  Let~\(T\) be the completion time of the
corresponding pre-entry stage at~\(c\), and suppose
\(T<\infty\).  Let
\(V_{\rm pre}=\mathcal V_T\) and put
\[
 A_\psi=\{x\in V_{\rm pre}\cap(\partial_\psi N\setminus R):
 \exists e\in F(\mathfrak W)\text{ incident to }x\}.
\]
Then \(A_\psi\) is finite and \(\cF_T\)-measurable.  If the front is then
moved across~\(N\), finite passage of the fresh search through~\(N\) is
contained in \(\mathsf{FinGrow}_\psi\).
\end{lemma}

\noindent{\bf Proof.}
The finite stopped search determines \(\mathcal V_T\), and hence the finite
set~\(A_\psi\), without querying~\(F(\mathfrak W)\).  The fixed-shift
search adds no vertex of \(V(N)\setminus R\) without querying an open edge.
By time~\(T\), such queries have been confined to \(E(\downspace_a)\) when
\(\psi=-\), or \(E(\upspace_b)\) when \(\psi=+\).  Hence
\(V_{\rm pre}\cap(V(N)\setminus R)\subseteq\partial_\psi N\).
By Lemma~\ref{l.slab-entry-faces}, any first fresh entry into the unread
layer uses an edge incident to~\(A_\psi\), so \(\mathsf{Shut}_\psi\)
prevents passage.  If \(\mathsf{Esc}_\psi\) occurs, an infinite complementary open
component is reachable from the entrance set.  Once the front is moved across
the layer, every edge of this component is eligible for the same exhaustive
search, so the layer stage cannot terminate finitely.  By the full-boundary
property from Lemma~\ref{l.fixed-shift-layer-record}, every layer edge from
\(V(N)\setminus R\) into~\(R\) is closed.  Thus a fresh crossing cannot be
completed by entering the retained trace.  Hence finite passage is contained
in \(\mathsf{FinGrow}_\psi\). \qed

\subsection{Successive recovered breaker windows}
\label{ss.successive-recovered-breaker-windows}

Now let \(q_0<q_1<\cdots<q_R\), with \(h_s=H_{q_s}\), where
\(h_1,\ldots,h_R\) are initially very fine, and put
\(\mathscr B_s=\mathscr B_{h_s}(v)\), \(1\leq s\leq R\).  Each window is
positioned using the adjusted height \(h_s+(v\cvl-u)\) of the current proposed
head.  Let~\(F_s\) be the union of the complements
\(F(\mathfrak W)\) of the fixed-shift layer records in the test layers of
\(\mathscr B_s\).  Write \(\partial_-\mathscr B_s\) for the southern side
of its first test layer and \(\partial_+\mathscr B_s\) for the northern side
of its last test layer.  Lemma~\ref{l.scheduled-breaker-separation} orders
these windows from south to north and makes their edge sets pairwise disjoint.
Every test layer~\(N\) in~\(\mathscr B_s\) satisfies the first alternative
of Lemma~\ref{l.fixed-shift-layer-record}: its layer edges have greater
endpoint height above \(h_s+(v\cvl-u)+b_0>\hadj(v)\), while
Lemma~\ref{l.breaker-layer-vertex-range} and \(b_0>r\) place its vertices
strictly above \(h_s>h_0\), hence outside~\(V(\dclust)\).

We now apply the northgoing front scheduling to these recovered breaker
windows.  For \(1\leq s\leq R\), let
\(n_s^-\) be the query count at the finite completion of the front at
\(\partial_-\mathscr B_s\), provided that the reached set has a contact
incident to~\(F_s\).  Put \(n_s^-=\infty\) if this stage is never reached,
is infinite, or finishes without such a contact, and set
\[
 \mathcal R_s^-=\cF_{n_s^-}\quad\hbox{on }\{n_s^-<\infty\},
 \qquad
 \mathsf{Ent}_s=\{n_s^-<\infty\}.
\]
On~\(\mathsf{Ent}_s\), no edge of~\(F_s\) has been queried, and every reached
vertex incident to~\(F_s\) lies on~\(\partial_-\mathscr B_s\).
By Lemma~\ref{l.freshness}, the edge states on~\(F_s\) retain their product
law at~\(\mathcal R_s^-\); at subsequent finite layer stops the same
holds for the edges still unread.
Proposition~\ref{p.breaker-persistence} supplies at least \(i_{\rm br}/2\)
thin layers in each row.  The height and vertex-range observations for~\(\mathscr B_s\)
place these layers under the first alternative of
Lemma~\ref{l.fixed-shift-layer-record}, giving their full-boundary records;
Lemma~\ref{l.fresh-interface-verification}
supplies their entrance sets and finite-growth inclusions.  For the selected layers, let
\(E_j\) be finite passage through the first \(j\) layers.  This data forms an
admissible family in the sense of Definition~\ref{d.admissible-breaker-tests}.

Starting from~\(\mathcal R_s^-\), raise the front successively across the
test layers and intervening gaps of~\(\mathscr B_s\).  Let \(n_s^+\) be the
query count when this window test terminates finitely, and put
\(n_s^+=\infty\) when the test is not begun or the search is infinite.  Define
\[
 P_s=\{n_s^+<\infty,\ \mathcal V_{n_s^+}\cap
 (\partial_+\mathscr B_s\setminus V(\mathfrak u^v))\ne\varnothing\}.
\]
By Lemma~\ref{l.slab-entry-faces}, finite passage through the window
entails finite passage through each selected test layer in order.  Thus
\(P_s\subseteq E_{i_{\rm br}/2}\), and Proposition~\ref{p.breaker-passage} gives
\begin{equation}
\label{e.recovered-window-passage-bound}
 \Rlaw^v_{\uclust,\dclust}(P_s\mid\mathcal R_s^-)
 \le q_{\rm row}
 \qquad\hbox{on }\mathsf{Ent}_s.
\end{equation}

Define
\begin{equation}
\label{e.successive-northgoing-events}
 A_0=\Omega,\qquad
 A_s=A_{s-1}\cap\mathsf{Ent}_s\cap P_s,
 \qquad 1\leq s\leq R.
\end{equation}
On~$A_{s-1}\cap\mathsf{Ent}_s$, height ordering gives
\[
 n_1^-\leq n_1^+\leq n_2^-\leq n_2^+\leq\cdots
 \leq n_{s-1}^-\leq n_{s-1}^+\leq n_s^-.
\]
Consequently
$A_{s-1}\cap\mathsf{Ent}_s\in\cF_{n_s^-}=\mathcal R_s^-$.  The tests of~$\mathscr B_s$ belong to the fixed-shift experiment for
the slide at index~$q_0<q_s$, while~$h_s$ is still ahead in the
attempted-slide schedule.  Proposition~\ref{p.breaker-persistence}
therefore supplies the thin pre-test records.

\begin{proposition}[Successive northgoing wavebreakers]
\label{p.successive-northgoing-wavebreakers}
For $0\leq s\leq R$,
\[
 \Rlaw^v_{\uclust,\dclust}(A_s\mid\cF_0)
 \leq q_{\rm row}^s.
\]
\end{proposition}

\noindent{\bf Proof.}
For $s\geq1$, the tower property at~$\mathcal R_s^-$ and
\eqref{e.recovered-window-passage-bound} give
\[
 \Rlaw^v_{\uclust,\dclust}(A_s\mid\cF_0) =\EE_{\Rlaw^v_{\uclust,\dclust}}\!\left[ \one_{A_{s-1}\cap\mathsf{Ent}_s} \Rlaw^v_{\uclust,\dclust}(P_s\mid\mathcal R_s^-) \middle|\cF_0\right] \leq q_{\rm row} \Rlaw^v_{\uclust,\dclust}(A_{s-1}\mid\cF_0).
\]
Iteration proves the claim. \qed

\subsection{The current-slide breaker window and recovered breakers}

Fix an active slide at~\(h_0\), with admissible fixed-shift data
\((\uclust,\dclust,v)\).  Its breaker window
\(\mathscr B_0=\mathscr B_{h_0}(v)\) is the three-row window of
Definition~\ref{d.breaker-window-geometry}, positioned relative to the
adjusted cut \(h_0+(v\cvl-u)\).  It is the first breaker window in the
northgoing sequence.  Farther north, let \(h_1<\cdots<h_R\) be later
initially very-fine scheduled heights and put
\(\mathscr B_s=\mathscr B_{h_s}(v)\), \(1\le s\le R\).  
We call \(\mathscr B_0\) the \emph{supplying breaker window}, to distinguish it from the later recovered breaker windows: its thin profile is supplied directly by \({\bf B}^{+}_{\rm f}\) at the active slide, rather than recovered from an earlier very-fine profile.
We consider a fixed-shift experiment under
\(\Rlaw^v_{\uclust,\dclust}\), exploring the reconstructed cluster by the
procedure of Section~\ref{ss.successive-stopped-tests}, with
\(\mathscr B_0\) tested first.  The query filtration is \((\cF_n)\),
and the shifted upcluster is traversed only upon contact, as in that procedure.

Let \(n_0^-\) be the query count at the finite completion of the
northgoing front at \(\partial_-\mathscr B_0\), provided that the reached set
has a contact incident to an unread edge of~\(\mathscr B_0\).  Put
\(n_0^-=\infty\) if this front search is infinite or finishes without such a
contact, and write \(\mathcal R_0^-=\cF_{n_0^-}\) on
\(\{n_0^-<\infty\}\).  If \(n_0^-<\infty\), raise the front successively
through the layers and gaps of the window, using the protocol of
Section~\ref{ss.successive-stopped-tests}.  Any newly pending edges below the current front are exhausted
before it advances again.  Let
\(n_0^+\) be the query count at finite completion of the window test, or
infinity if the test is not begun or the search is infinite.  These are
stopping times for~\((\cF_n)\).  Put
\(P_0=\{n_0^+<\infty,\ \mathcal V_{n_0^+}\cap
(\partial_+\mathscr B_0\setminus V(\mathfrak u^v))\ne\varnothing\}\).  After passage, test
\(\mathscr B_1,\ldots,\mathscr B_R\) in order, with the same stopping
rules and passage events~\(P_s\) as in the preceding subsection.
For a window whose last test layer is \(\slab_{a,b}\), a vertex is
\emph{strictly north of the window} when its height is greater than~\(b\).
We now bound the probability of fresh open growth north of these windows,
beyond the mere presence of the shifted upcluster.

\begin{proposition}
\label{p.origin-and-successive-wavebreakers}
In this experiment, for \(0\le s\le R\), let \(A_s^\circ\) be the event that the reconstructed cluster is finite and contains an undiscovered open edge with an endpoint
strictly north of~\(\mathscr B_0\) and of the first \(s\) recovered windows.  Then
\[
 \Rlaw^v_{\uclust,\dclust}(A_s^\circ\mid\cF_0)
 \leq q_{\rm row}^{\,s+1}.
\]
The estimate is uniform in the fixed fine record and in~$v$.
\end{proposition}

\noindent{\bf Proof.}
Consider first~\(\mathscr B_0\).  Every fresh departure from
\(\dclust\) starts from its top, at height at most~\(h_0\).  If a fresh edge is instead
incident to the shifted upcluster, Lemma~\ref{l.ssr-upcluster-bed-interface}
places its shifted-upcluster endpoint in the bed, at height at most~\(\hadj(v)\).
The enlarged interval for the first test layer of~\(\mathscr B_0\) begins above
\(\hadj(v)+b_0-r>\hadj(v)\).  Thus any fresh continuation reaching a test layer of~\(\mathscr B_0\)
starts strictly to its south; Lemma~\ref{l.slab-entry-faces}
places its first entry on the southern side.

For a test layer~\(N=N_{h_0,\rho,j}(v)\) of~\(\mathscr B_0\), every
layer edge lies in \(E(\upspace_{\hadj(v)})\), while
Lemma~\ref{l.breaker-layer-vertex-range} and \(b_0>r\) place every
layer vertex strictly above~\(h_0\), hence outside~\(V(\dclust)\).
Thus the first alternative of Lemma~\ref{l.fixed-shift-layer-record} applies,
and the fixed vertex set is
\[
 R=V(\mathfrak u^v)\cap V(N).
\]
The condition \({\bf B}^{+}_{\rm f}\) in
the pre-slide state gives, in each row, at least \(i_{\rm br}/2\) indices
whose upcluster occupation in \(J^\sharp_{\rho,j}(h_0)\) is at most~\(C_1\).
Translation to~\(v\) preserves adjusted coordinates, and Lemma~
\ref{l.breaker-layer-vertex-range} therefore makes the corresponding
records thin.  Lemma~\ref{l.fixed-shift-layer-record} also identifies
the complementary layer edges as undiscovered.  At the pre-entry stop,
none of these edges has been queried.  The fixed-shift product law and
Lemma~\ref{l.freshness} give the required conditional product law, while
Lemma~\ref{l.fresh-interface-verification} places finite passage in the
corresponding \(\mathsf{FinGrow}_-\) events.  Thus successive selected-layer
passage is an admissible family whose terminal event contains window passage.
Proposition~\ref{p.breaker-passage} consequently bounds the
conditional probability of~\(P_0\) by \(q_{\rm row}\), given
\(\mathcal R_0^-\) on~\(\{n_0^-<\infty\}\).
Since \(P_0\subseteq\{n_0^-<\infty\}\), the tower property gives
\begin{equation}
\label{e.current-slide-window-passage}
 \Rlaw^v_{\uclust,\dclust}(P_0\mid\cF_0)\leq q_{\rm row}.
\end{equation}

The later breaker windows are edge-disjoint from~\(\mathscr B_0\).
On~\(P_0\), their fresh edge sets~\(F_s\) remain unread.  At each subsequent
pre-entry stop, Lemma~\ref{l.freshness} preserves the conditional product
law of the still-unread layer edges.
A fresh continuation that crosses a later recovered window from south to
north must first enter each of its test layers from the south, by
Lemma~\ref{l.slab-entry-faces}.  Proposition~\ref{p.breaker-persistence}
supplies the thin records in every recovered window.  The height and
vertex-range check for these test layers in the preceding subsection,
together with the fixed-shift record, entrance and freshness lemmas used
there, makes successive
selected-layer passage an admissible family whose terminal event contains
window passage.  Thus the conditional bound
\eqref{e.recovered-window-passage-bound} applies to each recovered
window in the present experiment.

With~\(A_t\) defined for these tests by~\eqref{e.successive-northgoing-events},
the event
\(P_0\cap A_{t-1}\cap\mathsf{Ent}_t\) is measurable at~\(\mathcal R_t^-\):
the supplying-window test and the preceding recovered-window tests finish
before that pre-entry stop.  The tower property therefore gives
\[
 \Rlaw^v_{\uclust,\dclust}(P_0\cap A_t\mid\cF_0)
 \leq q_{\rm row}
 \Rlaw^v_{\uclust,\dclust}(P_0\cap A_{t-1}\mid\cF_0),
 \qquad 1\leq t\leq R.
\]
Iteration yields
\[
 \Rlaw^v_{\uclust,\dclust}(P_0\cap A_s\mid\cF_0)
 \leq q_{\rm row}^{\,s}
 \Rlaw^v_{\uclust,\dclust}(P_0\mid\cF_0).
\]

Suppose \(A_s^\circ\) occurs, and choose an open path from the origin
whose final edge is an undiscovered edge witnessing this event.  Write~$x$
for its endpoint strictly north of the windows.  By
Lemma~\ref{l.ssr-upcluster-bed-interface}, $x$ lies outside the shifted
upcluster: an undiscovered edge cannot meet that cluster at these heights.
After its last visit to \(V(\dclust)\cup V(\mathfrak u^v)\), the path
uses only undiscovered edges.  This fresh continuation departs either
from the downcluster, at height at most~\(h_0\), or from the shifted
upcluster, through its bed at height at most~\(\hadj(v)\), by
Lemma~\ref{l.ssr-upcluster-bed-interface}.  It therefore crosses every
intervening breaker window from south to north, as described by
Lemma~\ref{l.slab-entry-faces}.

To check that the front-scheduled search registers these crossings, use
the front-exhaustion observation in the paragraph preceding
Lemma~\ref{l.fresh-interface-verification}.  At the southern pre-entry
front of each window, it gives a reached contact with the unread first
test layer; at the northern completion front, it gives a reached vertex
outside the shifted upcluster on the far side.  On the finite-cluster
event, the successive stages finish finitely.  Thus the supplying window
records~$P_0$ and each of the first~$s$ recovered windows records
$\mathsf{Ent}_t\cap P_t$.  Hence, using the iteration above and
\eqref{e.current-slide-window-passage}, we obtain
\[
 \Rlaw^v_{\uclust,\dclust}(A_s^\circ\mid\cF_0)
 \leq
 \EE_{\Rlaw^v_{\uclust,\dclust}}\!\left[
   \one_{P_0}q_{\rm row}^{\,s}\mid\cF_0\right]
 \leq q_{\rm row}^{\,s+1},
\]
and complete the proof of Proposition~\ref{p.origin-and-successive-wavebreakers}.
\qed

\section{Two future breakers after an active slide}
\label{ss.breaker-stands-postpassage}

When a slide is performed at a fine height, its breaker window has three
rows to the north of the cut.  The new midpiece may enter this window, but
to reach the middle or upper row it must first pass through the lower one.
The thin test layers of the lower row make such passage unlikely.
The new midpiece thus typically leaves the middle and upper rows untouched,
so that they may act as separate breakers against southgoing explorations
launched by later slides farther north, protecting the region to their
south.  We first bound the probability that the supplying midpiece reaches
either of these rows, and then establish persistence and passage bounds
when both remain untouched.

\subsection{Reserved rows at the supplying slide}

Let \(S\) and \(S'\) be the input and output states of the slide at~\(h\),
 and let \(v_h=\head(S')\) be its selected head.  The upcluster
translation is \(\tau_h=v_h-\head(S)\).  The middle and upper rows,
\(\mathscr B_h^1\) and \(\mathscr B_h^2\), are the \emph{reserved rows}.
Throughout this section the reserved layers and rows are fixed in the
positions determined by that slide:
\(N_{h,\rho,j}=N_{h,\rho,j}(v_h)\) and
\(\mathscr T_h^\rho=\mathscr T_h^\rho(v_h)\), with the same convention
for \(\mathscr B_h^\rho\).  Define
\[
 \mathsf{Clean}_h
 =\bigl\{V(\midpiece_h(S'))\cap
 V(\mathscr T_h^1\cup\mathscr T_h^2)=\emptyset\bigr\}.
\]
Thus the midpiece formed at~\(h\) has not entered either reserved test graph.  In
particular it has not queried a layer edge there: such a query can occur only
after the fresh exploration has reached one of the edge's endpoints, which
is a vertex of the test graph.

\begin{proposition}[Finite entry into the reserved rows]
\label{p.fixed-shift-clean-future-breakers}
Fix a fine input state at height~\(h\) and admissible fixed-shift data
\((\uclust,\dclust,v)\).  Let \(\mathsf{FinHit}_h^v\) be the event that the cluster~\(\mc C(0)\) in
the reconstructed configuration is finite and contains an undiscovered open
edge incident to \(V(\mathscr T_h^1(v)\cup\mathscr T_h^2(v))\).  Then
\[
 \Rlaw^v_{\uclust,\dclust}(\mathsf{FinHit}_h^v\mid\cF_0)
 \leq q_{\rm row},
\]
uniformly in the fine input state and proposed head.  On \(\Accept_v\), \(\mathsf{FinHit}_h^v\) occurs exactly when
\(\mathsf{Clean}_h\) fails.
\end{proposition}

\noindent{\bf Proof.}
Use the cluster search of
Section~\ref{ss.successive-stopped-tests}, testing the lower row
\(\mathscr B_h^0(v)\) first.  Departures from the downcluster and from the
activated upcluster bed lie strictly south of this row.  Thus, on
\(\mathsf{FinHit}_h^v\), Lemma~\ref{l.slab-entry-faces} forces passage
through every selected layer of the lower row; finiteness of the reconstructed
cluster makes each passage finite.  The fine profile provides at least
\(i_{\rm br}/2\) thin layers.  The record, freshness and entrance checks for these layers are the
same as those for~\(\mathscr B_0\) in the proof of
Proposition~\ref{p.origin-and-successive-wavebreakers}; hence their successive
passage events form an admissible family whose terminal event contains passage
through the lower row.  Proposition~\ref{p.breaker-passage} therefore gives the factor
\(q_{\rm row}\) at the finite pre-entry field; the tower property gives the
displayed bound.  On \(\Accept_v\), the acceptance definition
\eqref{e.fixed-shift-acceptance} identifies the retained downcluster and
upcluster with~\(\dclust\) and~\(\mathfrak u^v\).  The undiscovered open
edges of the reconstructed cluster are therefore exactly the midpiece edges, whose
endpoints form \(V(\midpiece_h(S'))\), by
Proposition~\ref{p.dmutwoheights}(4).
This identifies \(\mathsf{FinHit}_h^v\) with failure of cleanliness there.
\qed

The preceding estimate is uniform in the proposed head.  Together with
Corollary~\ref{c.tether-head-summability}, it gives the corresponding
estimate for the active slide itself.

\begin{corollary}[Cleanliness of the two reserved rows]
\label{c.accepted-clean-future-breakers}
There are \(C_{\rm cl}<\infty\) and \(\delta_{\rm cl}\in(0,1)\) such that,
for every fine input state~\(S\) at height~\(h\),
\begin{equation}
\label{e.clean-future-breakers}
 \slide_h(S,\mathsf{Clean}_h^c)
 \le C_{\rm cl}q_{\rm row}^{\delta_{\rm cl}}.
\end{equation}
For the allowance in~\eqref{e.allowance-design}, \(\delta_{\rm cl}\) may be
chosen independently of \(\Cfine\), while
\begin{equation}
\label{e.clean-future-Cfine-dependence}
 C_{\rm cl}\le C(2+\Cfine)^{D_{\rm cl}}
\end{equation}
for fixed \(C,D_{\rm cl}<\infty\), uniformly in the breaker-profile
parameters \(b_0,C_0,C_1,i_{\rm br}\).  Thus the probability that the newly
formed midpiece reaches the middle or upper row can be made arbitrarily small
by increasing the number~\(i_{\rm br}\) of test layers in each row.
\end{corollary}

\noindent{\bf Proof.}
Fix the retained pair of~\(S\), and write
\[
 Z=\sum_{v\in\fubuthree}
     \Rlaw^v_{\uclust,\dclust}(\Accept_v),
 \qquad
 \nu(v)=Z^{-1}\Rlaw^v_{\uclust,\dclust}(\Accept_v).
\]
By Lemma~\ref{l.qfine}, \(Z\ge c>0\), and by
Corollary~\ref{c.tether-head-summability},
\(\sum_{\nu(v)<t}\nu(v)\le C_{\rm T}t^{\gamma_{\rm T}}\).
Lemma~\ref{l.probform} and Proposition~\ref{p.fixed-shift-clean-future-breakers}
give
\[
 \slide_h(S,\mathsf{Clean}_h^c)
 \le
 \sum_v \min\{q_{\rm row}/Z,\nu(v)\}.
\]
Put \(a=q_{\rm row}/Z\).  For \(a\le1\), split the sum at
\(t=a^{1/2}\).  The small atoms contribute at most
\(C_{\rm T}a^{\gamma_{\rm T}/2}\), while there are at most
\(a^{-1/2}\) remaining heads and their total contribution is at most
\(a^{1/2}\).  Enlarging the constant covers \(a>1\).  Thus one may take
\(\delta_{\rm cl}=\tfrac12\min\{1,\gamma_{\rm T}\}\).
The stated dependence on~\(\Cfine\) follows from
\eqref{e.tether-Cfine-dependence} and the fixed lower bound on~\(Z\). \qed

For \(\rho\in\{1,2\}\) and \(j\in\intint{i_{\rm br}}\), put
\[
 R_{h,\rho,j}
   =V(\tau_h\cdot\upcluster_h(S))\cap V(N_{h,\rho,j}).
\]
We refer to these as the \emph{traces at height~\(h\)}.  Translation by~\(\tau_h\)
preserves adjusted coordinates, so the fine breaker profile at~\(h\) and
Lemma~\ref{l.breaker-layer-vertex-range} give, in each reserved row, at
least \(i_{\rm br}/2\) indices~\(j\) for which \(|R_{h,\rho,j}|\le C_1\).

\subsection{Persistence and passage at a later test}

Now fix \(\rho\in\{1,2\}\) and a later active slide at height~\(H>h\),
and suppose that no midpiece
formed after the slide at~\(h\) and before this slide has entered
\(\mathscr T_h^\rho\).  Condition on the retained upcluster--downcluster pair
\((\uclust,\dclust)\) for the slide at~\(H\) and on an admissible proposed
head~\(v\), whose adjusted height is \(\hadj(v)=H+(v\cvl-u)\).
The following records use the current fixed-shift data, with layer
positions determined by the earlier slide at~\(h\).  For each test layer put
\[
 R_{h,\rho,j}^{\rm pre}
   =V(\dclust)\cap V(N_{h,\rho,j}),
 \qquad
 D_{h,\rho,j}^{{\rm pre},{\rm op}}
   =E(\dclust)\cap E(N_{h,\rho,j}),
\]
and let \(D_{h,\rho,j}^{{\rm pre},{\rm cl}}\) be the remaining layer edges
incident to \(R_{h,\rho,j}^{\rm pre}\).  Set
\[
 \mathfrak W_{h,\rho,j}^{\rm pre}
 =\bigl(R_{h,\rho,j}^{\rm pre},
        D_{h,\rho,j}^{{\rm pre},{\rm op}},
        D_{h,\rho,j}^{{\rm pre},{\rm cl}}\bigr),
 \qquad
 F_{h,\rho,j}^{\rm pre}
 =E(N_{h,\rho,j})\setminus
   \bigl(D_{h,\rho,j}^{{\rm pre},{\rm op}}
        \cup D_{h,\rho,j}^{{\rm pre},{\rm cl}}\bigr).
\]
Thus the record and its complement are determined by the current fixed-shift
data \((\uclust,\dclust,v)\), with no additional edge-state conditions
inherited from the slide at~\(h\).

For the southgoing test of \(\mathscr B_h^\rho\), use the fixed-shift
cluster exploration of Definition~\ref{d.fixed-shift-origin-search}
with the southgoing front protocol of Section~\ref{ss.successive-stopped-tests}.
Exhaust the front at the northern side of the last test layer, leaving
\(\bigcup_jF_{h,\rho,j}^{\rm pre}\) unread.  If this pre-entry search is
finite and the reached set contains a vertex incident to an edge of
\(\bigcup_jF_{h,\rho,j}^{\rm pre}\), let \(T_{h,\rho}^+\) be its
terminal query count; otherwise put
\(T_{h,\rho}^+=\infty\).  On
\(\{T_{h,\rho}^+<\infty\}\), let \(\Xi_{h,\rho}^+\) be the finite
set of edges queried by the pre-entry search and
\(\mathcal R_{h,\rho}^+\) its query filtration stopped at
\(T_{h,\rho}^+\), including the fixed-shift data.  Write
\[
G_{h,\rho}^{\rm pre}=\dclust\cap\mathscr T_h^\rho
\]
for the discovered-open graph in the reserved row before the fresh
component is adjoined.  On \(\{T_{h,\rho}^+<\infty\}\), lower the front successively through the
test layers and gaps of the row as in
Section~\ref{ss.successive-stopped-tests}, and let
\(\mathsf{Pass}^{\downarrow}_{h,\rho}\) be the event that this test
terminates finitely and the reached set at its completion meets
\(\partial_-N_{h,\rho,1}\setminus V(\dclust)\).  Set this event to
be empty when \(T_{h,\rho}^+=\infty\).

It remains to show that a reserved row untouched until a later slide retains the thin records needed for this southgoing test.

\begin{proposition}[Persistence before a southgoing test]
\label{p.southgoing-breaker-persistence}
Fix \(\rho\in\{1,2\}\). On \(\mathsf{Clean}_h\), consider a later active slide at height \(H>h\) with admissible data as above. Suppose no midpiece formed by an intervening slide has entered \(\mathscr T_h^\rho\).
Then, for every \(j\in\intint{i_{\rm br}}\),
\begin{equation}
\label{e.southgoing-trace-containment}
 R_{h,\rho,j}^{\rm pre}\subseteq R_{h,\rho,j},
\end{equation}
and
\[
 G_{h,\rho}^{\rm pre}
 \subseteq
 \bigl(\tau_h\cdot\upcluster_h(S)\bigr)\cap\mathscr T_h^\rho.
\]
The triple \(\mathfrak W_{h,\rho,j}^{\rm pre}\) is a layer record in the
sense of \eqref{e.layer-record-full-boundary}, and every edge of
\(F_{h,\rho,j}^{\rm pre}\) is undiscovered in the current fixed-shift
partition.  Moreover, on \(\{T_{h,\rho}^+<\infty\}\),
\begin{equation}
\label{e.southgoing-stopping-domain}
 \Xi_{h,\rho}^+\cap\bigcup_jF_{h,\rho,j}^{\rm pre}=\emptyset,
\end{equation}
and, conditionally on \(\mathcal R_{h,\rho}^+\), the current percolation
states on \(\bigcup_jF_{h,\rho,j}^{\rm pre}\) have their Bernoulli product
law.  If the exploration of undiscovered edges reaches a test layer, its first entry is through its northern side.
\end{proposition}

\noindent{\bf Proof.} The breaker geometry gives the margin needed at the vertex level.  Since
\(v_h\cvl-u\leq3r\), every vertex of the reserved row has height at most
\[
 h+(v_h\cvl-u)+b_0+C_0
 \leq h+3r+b_0+C_0
 <h+g.
\]
If \(H>h\) is any later attempted height, then
\(H-h\geq2K(g+1)\).  Since \(K\geq1\) and
\(g>3r+b_0+C_0>r\),
\[
 H-r-(h+g)\geq2K(g+1)-g-r>0.
\]
Together with the preceding bound, this gives
\begin{equation}
\label{e.southgoing-row-margin}
 \sup\{x\cvl:x\in V(\mathscr T_h^\rho)\}
 < H-r
 \leq \hadj(v)-r,
\end{equation}
where the last inequality uses \(v\cvl\geq u\).  Thus no vertex of the
current upcluster can lie in the reserved row: every vertex of
\(\upspace_{\hadj(v)}\), and hence of that upcluster, has height strictly
greater than \(\hadj(v)-r\).

Since \(b_0>r\), Lemma~\ref{l.breaker-layer-vertex-range} also places
every vertex of either reserved row strictly above the adjusted cut used
in the slide at~\(h\).  Immediately after that slide,
\(\mathsf{Clean}_h\) therefore leaves in \(\mathscr T_h^\rho\) only
the retained translate \(\tau_h\cdot\upcluster_h(S)\).

Consider successively the active slides completed after height~\(h\) and
before the current slide at height~\(H\).  Applying
\eqref{e.southgoing-row-margin} at each such height, with its input head
and its proposed head in place of~\(v\), excludes both the input and the
translated upcluster from the reserved row.  Material already
present there is retained only through the downcluster; any remaining cluster
material in the row belongs to the midpiece of the pre-slide state and is removed by the
slide.  By hypothesis the new midpiece does not enter
\(\mathscr T_h^\rho\).  Hence the cluster trace in the reserved row can
only decrease from one such state to the next.  Just before the slide at height~\(H\)
it is therefore contained in the corresponding trace at height~\(h\), and its retained part
in the current fixed-shift data is precisely the downcluster part.  This
proves \eqref{e.southgoing-trace-containment} and the displayed graph
inclusion.

By \eqref{e.southgoing-row-margin}, every edge of a reserved test layer
lies in \(E(\downspace_H)\).  Definition~\ref{d.prime} therefore fixes
open the layer edges of~\(\dclust\), fixes
closed every other layer edge incident to
\(R_{h,\rho,j}^{\rm pre}\), and declares every remaining layer edge
undiscovered.  Each fixed-open layer edge has both endpoints in
\(R_{h,\rho,j}^{\rm pre}\).  This is exactly the record
\(\mathfrak W_{h,\rho,j}^{\rm pre}\) and its complement
\(F_{h,\rho,j}^{\rm pre}\).  The stopping rule gives
\eqref{e.southgoing-stopping-domain}.  Lemma~\ref{l.freshness}, applied to
the current auxiliary product field on these undiscovered edges and
localised on \(\{T_{h,\rho}^+<\infty\}\) as explained after that lemma, gives
the asserted conditional Bernoulli product law.

Finally, fresh departures are launched from \(\Top(\dclust)\) and, after
attachment, from \(\Bed(\mathfrak u^v)\).  By
Lemma~\ref{l.geometrybasics}(1), the former vertices have heights in
\((H-r,H]\); the latter have heights in
\((\hadj(v)-r,\hadj(v)]\), hence strictly above~\(H-r\).  Meanwhile
\eqref{e.southgoing-row-margin} places the whole reserved row strictly below
\(H-r\).  A continuation of the undiscovered-edge exploration from either seed set
that reaches a test layer therefore first enters it through its northern side,
by Lemma~\ref{l.slab-entry-faces}. \qed

\begin{proposition}[Later passage of a reserved row]
\label{p.breaker-stands-postpassage}
Under the hypotheses of
Proposition~\ref{p.southgoing-breaker-persistence}, for \(\rho\in\{1,2\}\), on
\(\{T_{h,\rho}^+<\infty\}\),
\[
 \Rlaw^v_{\uclust,\dclust}
 \bigl(\mathsf{Pass}^{\downarrow}_{h,\rho}
       \mid\mathcal R_{h,\rho}^+\bigr)
 \leq q_{\rm row}.
\]
\end{proposition}

\noindent{\bf Proof.} By \eqref{e.southgoing-trace-containment}, the current
fixed vertex set in each test layer is contained in the corresponding trace
at height~\(h\).  The fine breaker profile at~\(h\) therefore supplies
at least \(i_{\rm br}/2\) layers in the row with thin current records
for the successive tests.  Proposition~\ref{p.southgoing-breaker-persistence}
identifies their full incident
open--closed records and their current undiscovered complements, and gives
the stopping-domain separation \eqref{e.southgoing-stopping-domain}.
Lemma~\ref{l.freshness} leaves the edge states on these complements
under the product law at every subsequent finite layer stop, while
Lemma~\ref{l.fresh-interface-verification} places finite passage in the
corresponding \(\mathsf{FinGrow}_+\) events.
Thus successive selected-layer passage is an admissible family whose terminal
event contains row passage.  Proposition~\ref{p.breaker-passage} gives the
displayed bound. \qed

\begin{proposition}[Successive southgoing breaker rows]
\label{p.successive-southgoing-breakers}
Fix an active slide at height~$H$ and admissible fixed-shift data
$(\uclust,\dclust,v)$.  Let
$\mathscr B_{h_{a_1}}^{\rho_1},\ldots,
\mathscr B_{h_{a_m}}^{\rho_m}$ be reserved rows of earlier active slides,
listed from north to south, with their positions fixed at those earlier
slides.  Suppose that $\mathsf{Clean}_{h_{a_s}}$ holds for every
$1\le s\le m$, and that no midpiece formed after the slide at
$h_{a_s}$ and before the current slide has entered
$\mathscr T_{h_{a_s}}^{\rho_s}$.  Both reserved rows from a given earlier
slide may occur among these $m$ rows.

Run the fixed-shift origin-cluster search with a single southgoing height
front through these rows in order.  For row~$s$, let $\mathsf{Ent}_s$ be
the event of a finite pre-entry stop with an entrance contact, and let
$P_s$ be finite passage through the row, using the one-row conventions
above with stops taken in this common exploration.  Put $A_0=\Omega$ and
$A_s=A_{s-1}\cap\mathsf{Ent}_s\cap P_s$ for $1\le s\le m$.
Then
\begin{equation}
\label{e.successive-southgoing-row-bound}
 \Rlaw^v_{\uclust,\dclust}(A_m\mid\cF_0)
 \le q_{\rm row}^{\,m}.
\end{equation}
On $\Accept_v$, the event~$A_m$ agrees with passage of all $m$ rows
by the output midpiece.
\end{proposition}

\noindent{\bf Proof.} The two reserved rows of a slide are ordered and have
disjoint edge sets by the breaker-window geometry; rows from distinct
slides have the same properties by
Lemma~\ref{l.scheduled-breaker-separation}.  Write $T_s^+$ for the
pre-entry query count of row~$s$ in the common exploration, with
$T_s^+=\infty$ when the stage is not reached or does not finish with
an entrance contact, as in the one-row convention.  On
$A_{s-1}\cap\{T_s^+<\infty\}$, the passage
events for the earlier rows have been determined before this stop, and
this event belongs to the stopped field~$\cF_{T_s^+}$.  The
southgoing front has exhausted all eligible edges north of row~$s$.
Rows may share boundary vertices, but their edge sets are disjoint;
the front has not yet admitted the complementary edges of row~$s$.

Proposition~\ref{p.southgoing-breaker-persistence} supplies the current
thin records and their undiscovered complements in row~$s$.
At the finite pre-entry stop, Lemma~\ref{l.freshness} gives the
conditional Bernoulli product law on the still-unread edges; the search
uses only its own edge-query filtration.  At subsequent finite layer
stops the same freshness assertion applies, and
Lemma~\ref{l.fresh-interface-verification} places finite passage in the
successive finite-growth events.  Thus the argument of
Proposition~\ref{p.breaker-stands-postpassage}, now at this stop in the
common exploration, gives
\[
 \Rlaw^v_{\uclust,\dclust}(P_s\mid\cF_{T_s^+})
 \le q_{\rm row}
 \quad\hbox{on }A_{s-1}\cap\mathsf{Ent}_s.
\]
The tower property, using the stated stopped-field measurability, gives
$\Rlaw^v_{\uclust,\dclust}(A_s\mid\cF_0)
\le q_{\rm row}\Rlaw^v_{\uclust,\dclust}(A_{s-1}\mid\cF_0)$.
Iteration proves~\eqref{e.successive-southgoing-row-bound}.

On acceptance the reconstructed origin cluster is finite, so the
exhaustive search finishes by Lemma~\ref{l.fixed-shift-search-completeness}.
Its undiscovered open edges are exactly the output midpiece edges, by
Proposition~\ref{p.dmutwoheights}(4).  Every midpiece component meeting
the southern side of a reserved row is joined to a departure from the
retained downcluster top or shifted upcluster bed, both above the row.
To see that the search registers passage, take an open path from the
origin to a midpiece vertex on the southern side of the row.  If this
vertex is not yet reached when a southgoing front at height~$c$ has been
exhausted, consider the first edge of the path leaving the reached set.
It is an unread open edge (the retained open graphs are traversed on
contact) and has greater endpoint height at most~$c$, since all eligible
pending edges have been processed.  For either face of the row, its
departing vertex cannot belong to either retained cluster:
their fresh departures are from above the row.  This vertex was therefore
reached by an edge of greater endpoint height above~$c$, and so lies in
$\vertbdry_c\setminus V(\dclust)$.  At the northern pre-entry front the
chosen midpiece vertex lies too far south to have been reached, and the
departing edge supplies an entrance contact in the northernmost test layer.
At exhaustion of the front at the southern face, either the chosen vertex
or the departing vertex gives a reached vertex on the southern side outside
the downcluster.  Thus the
search records both entrance and finite passage through the row.
Conversely, the full closed boundary of the retained
downcluster in the row prevents the fresh passage event from being
witnessed solely by vertices of that downcluster.  It therefore reaches
the southern side through midpiece edges.  This proves the asserted
identification for each row and hence for their intersection. \qed

This chapter completes the general tools that Part~\ref{part.two} provides
for analysing slide dynamics.
Propositions~\ref{p.successive-northgoing-wavebreakers} and
\ref{p.origin-and-successive-wavebreakers} bound successive northgoing
passages for the fine-delivery argument.
Proposition~\ref{p.fixed-shift-clean-future-breakers} bounds
finite entry into the reserved rows; on \(\mathsf{Clean}_h\),
Propositions~\ref{p.southgoing-breaker-persistence},
\ref{p.breaker-stands-postpassage} and
\ref{p.successive-southgoing-breakers} give their persistence and later
southgoing passage bounds.  Part~\ref{part.three} transfers these estimates
to the accepted-slide dynamics for renewal production and survival.

\part{Renewal Production and Survival}\label{part.three}

\chapter{Monitoring cluster accrual under slides}
\label{s.seal-sequence-layer}

The tools developed in Part~\ref{part.two} are applied in Part~\ref{part.three}
first to prove that renewal points are produced and then to show, in
Chapter~\ref{s.renewal-survival}, that at least one such renewal survives to
the terminal cluster.

With overwhelming probability, the input contains a positive density of very
fine heights, and the random schedule selects a positive density of them.
Our aim in this and the next two chapters is to prove that a positive density
of these selected heights remain fine when the rising slide schedule brings
them to limelight.  The transport of higher upclusters preserves the retained
northside data in adjusted-height coordinates, and the substantive task is to establish
the southside allowance conditions at limelight, simultaneously across
all depths.

Fix a \hfff{target-layer}{target height~\(\ell\) and a depth~\(k\)}.
The condition \({\bf A}_{\rm f}^{-}(k)\) bounds the target-height downcluster
in the corresponding unit band.  Place this band in its partition
layer~\(\cS(k)\).  Since
\(\downcluster_\ell\cap\cS(k)\subseteq\csetuthree\cap\cS(k)\),
it is enough to control the whole conditioned cluster in this layer.  The
attempted slide heights pass from its south side to its north side without
meeting its interior, dividing the schedule into \emph{advance} and
\emph{post-passage} phases.

\begin{figure}[htbp]
\centering
\includegraphics[width=0.75\textwidth]{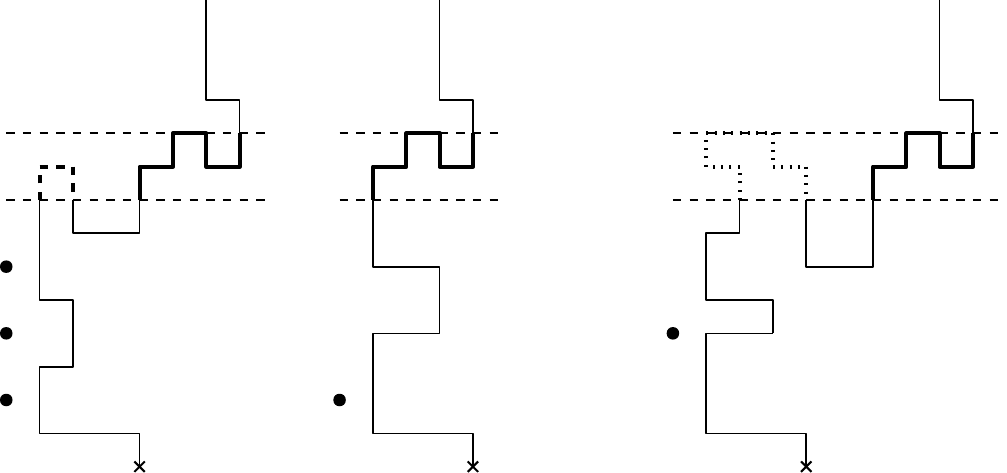}
\caption{The advance phase in the axial case $d = 2$ and $\vecell=e_2$. The three discs in the left panel mark the  attempted slide heights along this phase. {\em Left:} The input cluster emerges from the cross-marked origin, travelling through the dash-sided layer  towards the target height~$\ell$. The bold dashed fragment in the layer will exit at the first active slide. {\em  Middle:} 
After this slide, the surviving upcluster material in the layer is transported three units to the left.
 {\em Right:} A second active slide occurs when the schedule reaches the middle-disc height. The surviving host is transported five units to the right, while the newly sampled midpiece furnishes the dotted guest. Should the slide at the higher-disc height be dormant, this layer material will form the host--guest pair when the slide height passes over the layer at the upcoming moment of passage.}\label{f.advancephase}
\end{figure}

During the advance phase illustrated in Figure~\ref{f.advancephase}, the cluster material in the layer consists of a
\emph{host}, belonging to the transported upcluster, and a possibly empty
\emph{guest}, created when the current midpiece reaches the layer.  The host
may shrink as the rising cut removes components on midpiece exit, while
cut-it-and-lose-it removes a guest at the next active slide in the advance
phase.  Consequently, at passage the only possible guest is the one created
by the last active advance slide.  This guest, together with the surviving
host, determines the \hfff{retained-midpiece}{retained and midpiece
records~\(\RetRec\) and~\(\MidRec\)}, which form the passage record from
which later accrual in the layer is monitored.

The sealing constructions of Chapter~\ref{s.seed-exploration-seal-extension}
provide a way to enclose this passage record and follow subsequent growth in
the layer.  We seek a \emph{passage seal} that captures the surviving host
together with any guest left by the last active advance slide.  An initial
seal encloses the host, and enlargement absorbs any part of the guest outside
its capture.  The resulting passage seal encloses the complete passage
record, and its captured mesh region supplies the initial condition for the
subsequent sealing process.

\begin{figure}[htbp]
\centering
\includegraphics[width=0.75\textwidth]{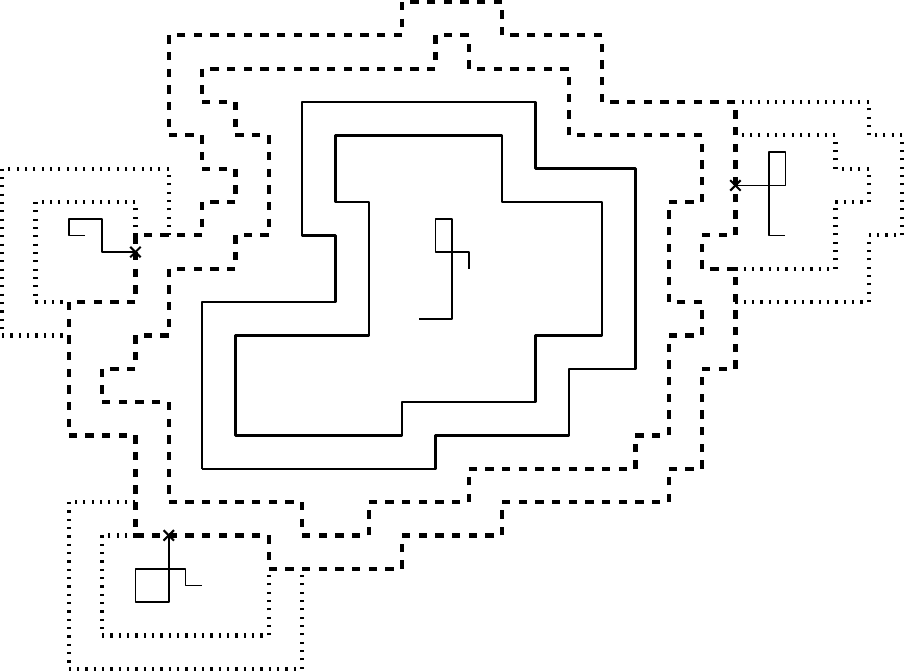}
\caption{With $d=3$ and $\vecell = e_3$, this planar depiction offers a schematic overhead view of a layer undergoing seal replacement and enlargement during a slide in the post-passage phase.
The discovered open structure lies in the centre, enclosed by the old seal, delimited by unbroken lines. Under replacement, a dash-bordered seal valid in the slide's resampled environment surrounds the old seal and traps the discovered cluster material. 
Exploration begins and finite growth episodes occur from three cross-marked virgin seeds at the replaced seal boundary. Seal enlargement captures this growth via the dotted seal caps. The outermost boundary borders the seal after the slide's sealing steps have been completed. The construction passes to the next active slide post-passage in an ongoing process of seal growth.}\label{f.seal}
\end{figure}

After passage, the slide schedule continues northwards towards~\(\ell\).
At each later active slide, fresh resampling may change edges certifying the
previous seal, so that seal cannot simply be carried into the new
configuration.  Its captured mesh region, however, remains a well-defined
finite region containing the inherited layer record.  Before further cluster
exploration in the layer, we place a fresh seal around this region in the new
fixed-shift field.  This is \emph{replacement}.  The ensuing exploration
may reveal additional cluster material as it continues through the captured
region, sometimes activating the known woven
fabric of the seal and finding virgin seeds at its exterior.  After a finite
exterior episode, \emph{enlargement} extends
the current seal to enclose it.  Thus the passage seal gives rise to
successive replacements and enlargements, of which Figure~\ref{f.seal} depicts one round, as the schedule advances towards
limelight.  The captured regions are nested, each fresh seal enclosing the
capture inherited from the preceding stage and any further finite growth.
For a layer at depth~\(k\), only \(O(1+k)\) attempted steps remain between
passage and limelight.  This auxiliary sealing construction accompanies the
accepted slide dynamics~\(\Qdyn\) without altering its cluster evolution.

This is the picture along the accepted slide trajectory.  Its quantitative
analysis, however, begins under a different probability law.  The seal and
wavebreaker estimates rely on the product law of unread edges, available
directly in the unconditioned fixed-shift experiment.  An active slide under
\(\Qdyn\), by contrast, conditions this experiment on acceptance.  Rather
than repeatedly impose this conditioning as estimates are accumulated along
the schedule, we establish the estimates of the present chapter one proposed
slide at a time under the fixed-shift law.  Their iteration in the next
chapter uses a \emph{reference dynamics}: the proposed head is chosen
according to the same cluster-dependent rule as under~\(\Qdyn\), but the
corresponding fixed-shift experiment is run without conditioning on
acceptance; rejection sends the reference trajectory to a cemetery state.
A trajectory-level comparison will then permit control obtained under the
reference law to be transferred to the accepted dynamics.

The construction of the passage seal presents a further difficulty.  The
enlargement that encloses the guest must be performed in the fixed-shift field
of the slide that creates it.  But the index of the last active slide in the
advance phase is not a stopping time: when the guest is formed, it is not yet
known whether this slide will be the last.  We therefore prepare, at each
active advance slide whose midpiece reaches the layer, the seal that would be
needed were that slide the last: an initial enclosure of the host, enlarged
to absorb the guest.  At passage, the appropriate candidate is selected; if
the last active advance slide did not reach the layer, the guest is absent and an initial seal is
constructed around the host alone.  These are the \emph{prospective passage seals}.
Although the advance phase may contain order~\(u\) attempted steps, the
preparations are affordable.  A midpiece formed far below the layer must
cross intervening wavebreaker windows to reach it, and the probability of
such a reach decays geometrically with the number of breakers.

So far we have followed one monitored layer.  The southside allowances,
however, must be controlled simultaneously across many depths via a
\emph{joint multi-layer accrual estimate}.  For a proposed slide south of a
monitored layer, the midpiece is unlikely to reach the layer if many breaker
windows intervene.  Among layers that can be reached after only a few
breakers, those lying close below their targets require the strongest
control, because comparatively little further accumulation can be tolerated
before the corresponding fine allowance is threatened.  A notion of
\emph{depth regularity}---later obtained on a positive-density subset of IVF
heights except with exponentially small probability---will rule out having
too many layers that are simultaneously close to the target height and separated from the slide
by only a few breakers, keeping the prospective seal preparation across the
monitored family under control.

A single proposed slide may encounter several monitored layers.  These
visits share the same shift and midpiece and unfold against a common revealed
history; their contributions cannot be treated independently.  We therefore
use one fixed-shift cluster exploration, pausing it whenever the outer search
supplies continuation to a monitored layer.  The layer exploration follows
the available cluster material, enclosing finite new growth before the outer
search resumes.  The auxiliary sealing operations change the scheduling of
queries but not the cluster determined by the underlying fixed-shift search.
The disjoint layer edge sets preserve the conditional product law on their
unread exteriors.  The contributions nevertheless remain dependent, and it
is the proposed slide, rather than the individual layer, that provides the
appropriate unit of probabilistic control.

For the joint estimate, the quantities ultimately at stake are the cluster
volumes in the monitored layers, since these are what the southside
allowances bound.  To control them, however, we must also control the regions
captured by the seals placed around the accumulating cluster material.
As in the Encloser control problem of
Section~\ref{ss.one-life-percolation-control}, the seal's capture is, in
effect, conceded territory: a subsequent cluster exploration may range
freely through it, following known open connections and querying still-unread
captured edges, without facing the uniform open-escape risk of fresh exterior
growth.  The volume of this captured territory is controlled through the
size of the seal.  Within one fixed-shift experiment we therefore control
both finite exterior gain and the boundary cost of enlarging the seal
around finite layer growth, including newly created mesh-boundary sites.
The cumulative outward
displacement caused by successive replacements is a distinct contribution,
treated in the next chapter.  Prospective passage sealing brings a further
cost, since possible initial seals must be prepared during the advance phase.

We collect these contributions one proposed slide at a time.  For a proposed
slide~\(b\), the \emph{one-slide cost}~\(\SlideCost_b\) combines the
finite-gain charge across the monitored layers, the boundary charge from
finite layer episodes and the excess cost of prospective passage
preparations.  For a suitably depth-regular family of monitored
layer--target pairs, Proposition~\ref{p.ssr-fixed-shift-input} establishes a
conditional exponential-moment bound for this cost, uniformly over the
preceding history and the proposed head, while preserving the geometric
attenuation supplied by intervening breakers.  At passage,
Corollary~\ref{c.ssr-last-advance-selection} bounds the weighted cost
of the selected prospective preparations by the accumulated seal costs,
while Proposition~\ref{p.ssr-passage-seal} gives the remaining
passage-seal excess a uniform exponential tail under the conditional
fixed-shift law and along the augmented accepted slide trajectory.  These
are the principal quantitative outputs of this chapter, providing the one-slide and passage
controls needed to bound accrual along the full schedule.

Chapter~\ref{s.model-problems-accrual} assembles these estimates under the
reference dynamics.  Its marked-kernel identity verifies that the auxiliary
sealing information does not alter the reference cluster evolution when
averaged out.  The chapter also controls the cumulative effect of
post-passage replacement.  Although the one-slide estimate bounds finite
cluster growth and the new seal boundary created by enlargement,
replacement may move the boundary of the captured region even without new
cluster material.  Following this \emph{replacement front} through
successive slides gives a linear-speed bound, up to an exponentially tailed
excess.  Combined with the one-slide and passage estimates, this yields a
joint upper tail for accrual across the monitored layers.  The chapter also
establishes the Tether comparisons used to return to the accepted slide
dynamics.

Chapter~\ref{s.shadow-bounds-delivery-fine} applies the joint accrual
estimate to the selected very-fine target heights.  A suitably spaced and
depth-regular family is monitored, and a tower-and-shadow argument bounds
the number of targets threatened by excessive southside accrual.  Together
with the northside data already carried to limelight, and the separate
depth-zero control, this proves the IVF-to-LLF transition and hence
positive-density renewal production.  Chapter~\ref{s.renewal-survival}
addresses the remaining problem of survival, using reserved breaker units
and a restart argument to show that, with overwhelming probability, at least
one produced renewal survives to the terminal cluster.

In the present chapter, the first section fixes the monitored layers and
constructs a paused fixed-shift exploration that follows several of them at
once.  The second treats prospective passage sealing, using breaker
attenuation to control the preparations and establishing the passage-seal
estimate.  The final section justifies post-passage replacement, verifies
that the augmented exploration projects to the fixed-shift cluster law, and
completes the conditional exponential-moment estimate for the one-slide cost.

\section{Monitored layers and the paused exploration}
\label{ss.monitored-layers-paused}

We first fix the monitored layers and distinguish their advance and
post-passage phases.  We then introduce the augmented histories needed for
conditioning and the paused exploration that follows these layers jointly.
This construction preserves the fixed-shift cluster law and the conditional
product law on unread edges.

\subsection{Target layers and the two phases}
\label{ss.accrual-and-sealing}

For this chapter, condition on a realisation
\(\overline H=(H_1,\ldots,H_L)\) of the attempted-height schedule in
Definition~\ref{d.qdyn}, and fix a scheduled target height
\(\ell=H_n\).  At step~\(i\leq n\), the fine test at height~\(H_i\) is made in
state~\(S_{i-1}\).  If it succeeds, the slide is \emph{active} and the
accepted resample is performed; otherwise the step is \emph{dormant} and
\(S_i=S_{i-1}\).  Thus, under the present conditioning, the attempted-height
schedule is fixed while the active subsequence is determined by the evolving
cluster states.

Fix \(k\in\Nzero\) and let
\(U_\ell(k)=\slab_{\ell-k-1,\,\ell-k}\) be the unit band addressed by
the condition \({\bf A}_{\rm f}^{-}(k)\).  Recall from
Section~\ref{ss.random-mesh} that the layer partition is delimited by the
integer-valued faces \(\beta_q=K(2q+1/2)\).  Since the endpoints of
\(U_\ell(k)\) are consecutive integers, this unit band meets no partition
face in its interior.  Equation~\eqref{e.layer-upper-endpoint} and
Lemma~\ref{l.layer-partition} therefore imply that \(E(U_\ell(k))\) is
contained in a unique partition layer.
\begin{definition}[The layer carrying the depth-\(k\) band]
\label{d.layer-at-k}
For \(k\in\Nzero\), let \(q(k)\in\Z\) be the unique integer such that
\[
   E(U_\ell(k))\subseteq E(N_{2q(k)}),
   \qquad
   \cS(k)=N_{2q(k)}=\slab_{s^-(k),s^+(k)}.
\]
We call \(\cS(k)\) the monitored layer associated with the
target--depth pair~\((\ell,k)\).
\end{definition}

Every positive attempted height is itself a partition face by the construction
of the random mesh in Section~\ref{ss.random-mesh}.  Thus the attempted slide
schedule crosses a monitored layer at most at its faces, never through its
interior.  We shall often work with the layer and target without recording
the particular unit band.  By a \emph{layer--target pair} we
mean a pair \((\cS,H)\), where \(H\) is a scheduled target height and
\(\cS=\slab_{s^-,s^+}\) is a partition layer with \(s^+\leq H\).  Its
\emph{geometric depth} is
\[
   \delta=H-s^+.
\]
If the unit band at depth~\(k\) below~\(H\) is assigned to~\(\cS\), then
\(|k-\delta|\leq2K\); in particular
\(\avf(k)\asymp\avf(\delta)\), with constants depending only on~\(K\).

\begin{definition}[Advance and post-passage phases]
\label{d.two-phases}
For a layer--target pair \((\cS,H)\), with
\(\cS=\slab_{s^-,s^+}\), the attempted steps at heights \(H_i\leq s^-\)
form the \hfff{two-phases}{\emph{advance phase}}, while those with
\(s^+\leq H_i<H\) form the \emph{post-passage phase}.  The \emph{passage}
of~\(\cS\) 
is its crossing by the attempted-height schedule, from the last attempted height at or below \(s^-\) to the first at or above \(s^+\).
 If no active step in the advance phase
occurs, the input state~\(S_0\) is the last pre-passage cluster state.
\end{definition}

No positive attempted height lies strictly between \(s^-\) and~\(s^+\), so
these two phases account for every attempted step strictly before the target step.
By the construction in Section~\ref{ss.random-mesh}, successive attempted
heights satisfy
\[
   H_{i+1}-H_i\geq 2K(g+1).
\]
For the layer carrying the depth-\(k\) band, one has
\(|\delta-k|\leq2K\), and hence \(\delta=O(1+k)\).  Together with the
spacing bound, this shows that its post-passage phase contains
\(O(1+k)\) attempted steps.

\begin{lemma}
\label{l.ssr-two-phase-cluster}
Let $S$ be the state after an active slide at height~$h$, and let
$\cS=\slab_{s^-,s^+}$ be a partition layer.  If $h\leq s^-$, then
\begin{equation}
\label{e.ssr-advance-cluster-decomposition}
 E(\csetuthree(S))\cap E(\cS)
\, = \, \bigl(E(\upcluster_h(S))\cap E(\cS)\bigr)
  \,\sqcup\,
  \bigl(E(\midpiece_h(S))\cap E(\cS)\bigr)\, .
\end{equation}
If $h\ge s^+$, then
\begin{equation}
\label{e.ssr-postpassage-cluster-decomposition}
 E(\csetuthree(S))\cap E(\cS)
\, = \, \bigl(E(\downcluster_h(S))\cap E(\cS)\bigr)
  \,\sqcup\,
  \bigl(E(\midpiece_h(S))\cap E(\cS)\bigr)\, .
\end{equation}
\end{lemma}

\noindent{\bf Proof.}
When $h\leq s^-$, the height-$h$ downcluster has no edge in~$E(\cS)$.
When $h\ge s^+$, one has $\hadj(h)\ge h\ge s^+$, so the height-$h$
upcluster has no edge in~$E(\cS)$.  The two identities now follow from
the DMU partition of the cluster into downcluster, upcluster and midpiece.
\qed

The upcluster and midpiece restrictions in~\eqref{e.ssr-advance-cluster-decomposition}
are the translated \emph{host} and the present \emph{guest},
as Section~\ref{ss.roadmap-iterated} describes them.  In the post-passage
decomposition~\eqref{e.ssr-postpassage-cluster-decomposition},
the retained cluster material in the layer belongs to the downcluster.
Section~\ref{ss.passage-preparation} constructs the passage record from
the final advance host and guest and their recorded closed boundary edges.

For the remainder of this chapter, fix a finite collection of
layer--target pairs \((\cS_j,\ell_j)\), indexed by~\(j\), whose partition
layers \(\cS_j=\slab_{s_j^-,s_j^+}\) are distinct; among the target heights, coincidences are permitted.  Write
\[
   E_j=E(\cS_j),
   \qquad \delta_j=\ell_j-s_j^+.
\]
The edge sets~$E_j$ are pairwise disjoint by
Lemma~\ref{l.layer-partition}, even when adjacent layer graphs share
boundary vertices. 

\subsection{Augmented histories and fixed-shift kernels}
\label{ss.ssr-augmented-histories}

Write
\[
   \mathcal K_t
   =\sigma(S_0,\overline H,S_1,\ldots,S_t).
\]
Under the slide dynamics~\(\Qdyn\), this is the cluster-state history up to
and including step~\(t\), together with the attempted-height schedule.  When
the reference dynamics is introduced in
Chapter~\ref{s.model-problems-accrual}, we use the same notation for its
cluster-state history, defined up to and including any cemetery transition.

To control accrual simultaneously across the monitored layers, we shall need
a more detailed history.  As the attempted schedule advances, monitored layers
may acquire seals; within an active step, the fixed-shift search moves between
the outer region and individual monitored layers, returning outside whenever
the available layer work is exhausted.  The edge queries and their outcomes,
together with the seals and captures produced, must be retained for
conditioning later one-step experiments.

Accordingly, for~\(0\leq t\leq L\), let~\(Y_t\) denote the history
through the first~\(t\) attempted steps, augmented by these query and sealing
records.  We use this notation first with the accepted slide
dynamics~\(\Qdyn\), and later with the reference dynamics constructed in
Chapter~\ref{s.model-problems-accrual}.  The index~\(t\) runs
on the attempted-step scale; within a single fixed-shift experiment a separate
paused filtration tracks successive queries.  At an active attempted
step~\(b\in\intint{L}\), the augmented history retains the proposed
head~\(V_b\) as a \emph{proposal mark}; at a dormant step, \(V_b\) is a fixed
null mark.  Put
\[
   \cH_t=\sigma(Y_t),
   \qquad
   \mathcal K_t\subseteq\cH_t.
\]
The next subsection specifies the paused exploration and the replacement of
post-passage seals; Section~\ref{ss.passage-preparation} treats the prospective
seals prepared before passage.

Fix an attempted step~\(b\in\intint{L}\), an incoming augmented history
\(Y_{b-1}=y\) at which the step is active, and a proposed head~\(v\in\Z^d\)
with \(v\cdot\vecell\in[u,u+3r]\), as in Definition~\ref{d.prime}.
The retained cluster data in~\(y\), together with~\(v\), determine the
fixed open and closed edge sets~$O_v$ and $C_v$.  Supply
\(\mathcal E_v^\circ=E(\Z^d)\setminus(O_v\cup C_v)\) with a new independent
Bernoulli-\(p\) product field.  On this one field, perform the fixed-shift
reconstruction together with the paused exploration and associated sealing
operations.  The law of this augmented experiment defines the \emph{unconditioned
augmented fixed-shift kernel}~\(\PP_b^v(y,\cdot)\); write~\(\EE_b^v\) for
its expectation kernel.  The reconstructed configuration is defined even
when the exploration does not complete; its possibly infinite query history
then has no finite terminal record.

For each monitored layer carrying an inherited post-passage seal in~\(y\),
write~\(\mathscr S_j^-\) for the seal produced by the preceding active
post-passage experiment, or for the passage seal if no such experiment has
occurred.  The \emph{incoming-capture condition} for such a layer is
\begin{equation}
\label{e.ssr-incoming-capture}
 \MeshSupp\Bigl(
   \bigl[E(\csetuthree(S_{b-1}))
     \cup\partial_E(\csetuthree(S_{b-1}))\bigr]\cap E_j
 \Bigr)
 \subseteq\MCap(\mathscr S_j^-).
\end{equation}
Here~\(S_{b-1}\) is the incoming cluster state in~\(y\); the left-hand
side is the mesh support of its complete open-and-closed cluster record
restricted to~\(E_j\).
We call~\((y,v)\) \emph{admissible} when~\(y\) is a live history consistent
with the preceding augmented update rules, \(v\) satisfies the head-height
condition above, and the incoming-capture condition holds for every layer
carrying such a seal.  Lemma~\ref{l.ssr-inherited-record-captured} establishes
this condition at passage and its preservation through successive active
post-passage steps.
For random admissible incoming histories we use
\(\PP_b^v(A\mid\cH_{b-1})=\PP_b^v(Y_{b-1},A)\), and similarly for
\(\EE_b^v\), as notation for evaluation of these kernels.  When no
auxiliary information has yet been created, the cluster projection reduces
to the ordinary fixed-shift law determined by the incoming cluster state.
The ensuing fixed-shift estimates are uniform over admissible~\((y,v)\).
The fresh field is unconditioned on~\(\Accept_v\); Chapter~\ref{s.model-problems-accrual}
specifies the proposal rule and the acceptance correction for the trajectory
laws.

\subsection{The multi-layer paused exploration}
\label{ss.multi-layer-paused}

We first prepare the post-passage replacement seals, then specify how the
outer cluster search is interleaved with exploration in the monitored layers.

Consider an active attempted step~$b\in\intint{L}$ of height~$H_b$
with an admissible incoming pair~$(y,v)$.  Only
monitored pairs whose target has not yet reached limelight are relevant to
this one-step experiment.  Split their indices as
\begin{equation}
\label{e.paused-current-layer-classes}
 \mathcal I_b^{\rm adv}=\{j:H_b\leq s_j^-\},\qquad
 \mathcal I_b^{\rm post}=\{j:s_j^+\leq H_b<\ell_j\},\qquad
 \mathcal I_b=\mathcal I_b^{\rm adv}\sqcup\mathcal I_b^{\rm post}.
\end{equation}
Since no attempted height lies strictly between~$s_j^-$ and~$s_j^+$, every
monitored pair with~$H_b<\ell_j$ belongs to exactly one of the two classes in
\eqref{e.paused-current-layer-classes}.

Under~$\PP_b^v$ the retained downcluster and shifted upcluster are fixed.
Recall from Definition~\ref{d.prime} that
\[
 \begin{split}
 O_v&=E(\mathfrak u^v)\cup E(\dclust),\\
 C_v&=\bigl(E(\downspace_{H_b})\cap\partial_E(\dclust)\bigr)
       \cup
       \bigl(E(\upspace_{\hadj(v)})\cap\partial_E(\mathfrak u^v)\bigr).
 \end{split}
\]
The edges in~$O_v$ are fixed open, those in~$C_v$ are fixed closed, and
\begin{equation}
\label{e.paused-undiscovered-set}
   \mathcal E_v^\circ=E(\mathbb Z^d)\setminus(O_v\cup C_v)
\end{equation}
is supplied by the fresh fixed-shift Bernoulli field.  At the start of the
present fixed-shift experiment, every edge of~$\mathcal E_v^\circ$ is unread.

\subsubsection{Fresh seals in the post-passage layers}
Before the cluster search begins, we refresh the seals in all post-passage
monitored layers, processing $j\in\mathcal I_b^{\rm post}$ in a fixed
deterministic order.  A seal constructed in an earlier active fixed-shift
experiment is certified only in that configuration: under the present
resampling, it need no longer be a seal.  Its captured mesh region, however,
is retained and used to build its replacement.

For~$j\in\mathcal I_b^{\rm post}$, the incoming-capture condition
\eqref{e.ssr-incoming-capture} says that~$\MCap(\mathscr S_j^-)$
contains the mesh support of the complete open-and-closed cluster record
of the incoming state in~$E_j$.  Since~$b$ is post-passage,
$s_j^+\le H_b\le\hadj(v)$.  By~\eqref{e.layer-upper-endpoint}, every
edge of~$E_j$ has upper endpoint at height at most~$s_j^+$.  Hence
\[
 E_j\cap E(\upspace_{\hadj(v)})=\varnothing,
 \qquad
 E_j\subseteq E(\downspace_{H_b}),
\]
and Definition~\ref{d.prime} identifies all fixed cluster data in this layer
with retained downcluster data:
\[
 O_v\cap E_j=E(\dclust)\cap E_j,
 \qquad
 C_v\cap E_j=\partial_E(\dclust)\cap E_j.
\]
The first set is open retained-downcluster data.  Every edge of the
second set is closed in the incoming configuration: if an edge of
$E(\downspace_{H_b})$ incident to~$\dclust$ were open, its other endpoint
would also belong to~$\dclust$.  Thus both sets are part of the incoming
cluster record in~$E_j$, and~\eqref{e.ssr-incoming-capture} gives
\begin{equation}
\label{e.postpassage-fixed-record-inherited}
 \MeshSupp\bigl((O_v\cup C_v)\cap E_j\bigr)
 \subseteq \MCap(\mathscr S_j^-).
\end{equation}

Thickening the inherited mesh capture by one mesh step and filling in the
resulting bounded holes, define the \emph{replacement carrier}
\begin{equation}
\label{e.postpassage-rebuild-carrier}
 D_{b,j}^{\rm repl}
 =\Capt_{\Mesh}\!\left(
     \MCap(\mathscr S_j^-)^{[1]}
   \right).
\end{equation}
The thickening gives the inherited capture a mesh collar before the
fresh seal is constructed.
Run the empty-seal completion of
Lemma~\ref{l.conditional-rooted-completion} from the carrier
$D_{b,j}^{\rm repl}$ in the fixed-shift field on~$E_j$, and denote its output
by~$\mathscr S_{b,j}^{\rm repl}$.  The layer edge sets are pairwise
disjoint by Lemma~\ref{l.layer-partition}.  With the fresh replacement seals
in place, the paused cluster exploration can begin.

\subsubsection{The paused cluster exploration}

Here we combine the fixed-shift cluster search of
Definition~\ref{d.fixed-shift-origin-search} with the within-layer search of
Definition~\ref{d.paused-layer-exploration}.  The outer phase processes pending
cluster edges outside the currently monitored layers, without processing edges
in those layers.  Once no such outer work remains, the exploration selects
a currently monitored layer with pending edges and runs the within-layer
search from their reached endpoints.  An advance-phase layer receives its initial seal only when first
reached; a
post-passage layer already carries~$\mathscr S_{b,j}^{\rm repl}$.  When no
further work remains in the layer, its reached side vertices return to the
outer search.  The two phases alternate until open escape occurs or no pending
cluster edge remains; in the latter case~$\Accept_v$ is determined from the
completed finite fixed-shift record.  We now specify the interleaving
precisely.

Put~$\mathcal H_-^{\rm pa}=\sigma(\cH_{b-1},v)$, where the superscript
$\mathrm{pa}$ indicates the \emph{paused} exploration.
For~$e\in\mathcal E_v^\circ$, write~$\varepsilon(e)$ for its fixed-shift
state.  The post-passage replacements above are successive stopping
constructions in this field.  Let~$\mathcal E_0^{\rm rd}$ be the set of edges
of~$\mathcal E_v^\circ$ queried by these replacements, and let
$\mathcal H_0^{\rm pa}$ be the history obtained by adjoining this queried set
and its revealed edge states to~$\mathcal H_-^{\rm pa}$.  Starting from the
fixed-shift product field at~$\mathcal H_-^{\rm pa}$, apply
Lemma~\ref{l.ssr-fresh-postpassage-rebuild} successively in the fixed
replacement order.  Pairwise disjointness of the layer edge sets leaves the
next layer unread at each step.  Thus, conditionally on~$\mathcal H_0^{\rm pa}$,
\begin{equation}
\label{e.paused-initial-product}
 \bigl(\varepsilon(e):e\in
   \mathcal E_v^\circ\setminus\mathcal E_0^{\rm rd}\bigr)
 \quad\hbox{has Bernoulli-$p$ product law.}
\end{equation}
Since later preliminary replacements use layer edge sets disjoint from that of
an earlier replacement, each~$\mathscr S_{b,j}^{\rm repl}$ remains a relative
seal for~$\mathcal H_0^{\rm pa}$.
We use the reached, processed and revealed records of
Definition~\ref{d.paused-layer-exploration} for the full fixed-shift
search.  Thus~$\mathcal E^{\rm rd}$ is the set of edges of
$\mathcal E_v^\circ$ whose fixed-shift states have been revealed,
$\mathcal Q\subseteq\mathcal E^{\rm rd}$ is the set already processed
by the cluster search, and~$V^{\rm rch}$ is the reached vertex set.
The pending cluster edges are those of
$\mathcal E_v^\circ\setminus\mathcal Q$ incident to~$V^{\rm rch}$.

Initialise
\[
 \mathcal J_{\rm seal}=\mathcal I_b^{\rm post},
 \qquad
 \mathscr S_j=\mathscr S_{b,j}^{\rm repl}
   \quad(j\in\mathcal I_b^{\rm post}).
\]
An advance-phase index enters~$\mathcal J_{\rm seal}$ when an initial seal is
placed in its layer.  For~$j\in\mathcal I_b$, let
\begin{equation}
\label{e.current-layer-regions}
 \mathcal U_j^{\rm cur}=
 \begin{cases}
   E_j\cap\mathcal E_v^\circ,&j\notin\mathcal J_{\rm seal},\\
   \mathcal U(\mathscr S_j),&j\in\mathcal J_{\rm seal}.
 \end{cases}
\end{equation}
In the sealed case, the exterior is already contained in~$\mathcal E_v^\circ$,
since the seal captures the fixed layer record.  A layer with a seal in
force may still contain unread edges in its capture; these are available
to the captured scan, while unread exterior edges are reserved for layer
episodes.

\begin{definition}[Multi-layer paused exploration]
\label{d.multi-layer-paused}
Use the fixed vertex and edge orders of
Definition~\ref{d.fixed-shift-origin-search}.  Between phases the state is
\[
 \bigl(V^{\rm rch},\mathcal Q,\mathcal E^{\rm rd},
       \mathcal J_{\rm seal},(\mathscr S_j)_{j\in\mathcal J_{\rm seal}}\bigr),
\]
where~$V^{\rm rch}$ is the reached vertex set.  Initially
\[
 V^{\rm rch}=V(\dclust),\qquad
 \mathcal Q=\varnothing,\qquad
 \mathcal E^{\rm rd}=\mathcal E_0^{\rm rd}.
\]
If~$V(\dclust)\cap V(\mathfrak u^v)\ne\varnothing$, adjoin
$V(\mathfrak u^v)$ to~$V^{\rm rch}$.  Thereafter, whenever processing an open
edge first makes~$V^{\rm rch}$ meet~$V(\mathfrak u^v)$, adjoin all of
$V(\mathfrak u^v)$ to~$V^{\rm rch}$ without a further query.

In an \emph{outer phase}, repeatedly process pending edges of~$\mathcal E_v^\circ$
that do not belong to~$E_j$ for any~$j\in\mathcal I_b$, updating the reached
set after each step.  Thus the outer search may process edges in any
partition layer not currently monitored, including layers whose targets have
already reached limelight.  Continue until no eligible pending edge remains.  The
reached endpoints of still-pending edges in~$\bigcup_{j\in\mathcal I_b}E_j$
then supply seeds to the monitored layers for this step.  Choose the least
such index~$j$, then the first reached vertex~$x$ incident to a pending
edge of~$E_j$, and its first such edge~$e$, using the fixed vertex and
edge orders.  Record the
seed--edge pair~$(x,e)$ before processing~$e$ or, when needed,
constructing an initial seal.

On the first visit to an advance-phase layer, this selection takes place
before any edge of~$E_j$ outside the fixed data has been read.  If the chosen
layer has no seal in force, construct its initial seal from the fixed
open-and-closed record~$(O_v\cap E_j,C_v\cap E_j)$ by
Proposition~\ref{p.initial-seal-estimate}, using the
\hyperref[conv.empty-seal]{empty-seal convention} when this record is empty;
the corresponding index then enters~$\mathcal J_{\rm seal}$.  Enlarge
$\mathcal E^{\rm rd}$ by the edges queried in this operation and recompute the
current regions.  Then begin the layer phase with all presently reached
vertices incident to pending edges of~$E_j$ as its seed set.

In a \emph{layer phase} for~$\cS_j$, first exhaust the captured scan:
process every pending edge of~$E_j\cap\mathcal E_v^\circ$ lying in
$\ECap(\mathscr S_j)$, adjoining it to~$\mathcal Q$.  Use its recorded
state if available; otherwise query it and adjoin it to~$\mathcal E^{\rm rd}$.
Continue through all open captured edges as they become incident to the
reached set.  Thus touching the woven fabric of a current seal component
merely adds its known open continuation to the captured scan.  Record
reached side vertices for the later outer phase.

When the captured scan is exhausted, if a pending edge of
$E_j\cap\mathcal E_v^\circ$ remains, choose the first reached vertex~$s$
incident to one in the fixed vertex order.  Its pending edges lie in
$\mathcal U(\mathscr S_j)$, so~$s$ is a virgin seed.  Select the first
incident exterior edge~$e_s$ in the fixed edge order, and run the exhaustive
episode of Proposition~\ref{p.relative-seal-play} from~$s$, beginning
with~$e_s$.  Open escape terminates the construction.  During
the episode every edge processed by the cluster search enters~$\mathcal Q$,
whether its state was already revealed or is first queried then.  On a finite
outcome replace~$\mathscr S_j$ by the output seal, enlarge~$\mathcal E^{\rm rd}$
by the edges first queried in the episode and seal update, and restart the
layer phase with the recomputed capture and~$V^{\rm rch}$.  Queries made
only for sealing do not enter~$\mathcal Q$.

If no pending edge of~$E_j\cap\mathcal E_v^\circ$ remains after the captured
scan, the layer phase is complete.  Return to the outer phase; any edge
leaving the layer from a reached side vertex is still pending there.
Iterate these phases until no pending cluster edge remains or open escape
occurs.
\end{definition}

\subsubsection{Properties of the paused exploration}

The paused-step filtration~$(\mathcal H_t^{\rm pa})$ is obtained by adjoining
one cluster query or one stopped seal query at a time.  An edge state already
revealed may be used without advancing the filtration.

\begin{lemma}[Adaptedness and first layer pause]
\label{l.ssr-paused-filtration}
At every finite paused time, the current state and the next move of the paused
exploration are measurable from the paused history, before any new product-field
query required by that move is made.  If~$H_b\leq s_j^-$ and~$\cS_j$ is first
paused while no seal is in force there, then before this pause no edge of~$E_j$
outside the fixed data $(O_v\cap E_j,C_v\cap E_j)$ has been read.
\end{lemma}

\noindent{\bf Proof.}
All choices are made by fixed deterministic orders from the displayed state.
The post-passage replacement seals have been constructed before the paused
cluster search starts.  The outer phase processes no edges in the currently
monitored layers.
Within a layer phase, membership of a pending edge
in the current capture is determined by the current seal; the captured scan
processes known open fabric as well as newly queried edges, and a virgin
seed is selected before its incident unread exterior edge is queried.
For a layer with~$H_b\leq s_j^-$ and no seal, the outer phase makes no
queries in~$E_j$.  The preliminary replacements and seal operations in other
layers query only their own, pairwise-disjoint layer edge sets.  No seal
operation has yet occurred in~$E_j$ itself.  Hence none of these edges has
been read before the first pause. \qed

The paused exploration sometimes learns an edge state while constructing or
enlarging a seal.  This does not count as processing that edge in the cluster
search: if the cluster search later reaches the edge, it uses the state already
revealed.  With this distinction, its cluster part can be compared directly
with the fixed-shift search of Definition~\ref{d.fixed-shift-origin-search}.

\begin{lemma}[Agreement with the fixed-shift search]
\label{l.multi-layer-regional-completeness}
On a finite completed outcome, suppressing the auxiliary seal-only queries
from the paused exploration leaves exactly the completed cluster record---its reached
cluster together with its incident closed boundary---of the fixed-shift search
of Definition~\ref{d.fixed-shift-origin-search}.  In particular,~$\Accept_v$
is measurable from the terminal paused history.  Conversely, if the fixed-shift
origin cluster is finite, the paused exploration completes after finitely many
steps, almost surely.
\end{lemma}

\noindent{\bf Proof.}
The paused exploration uses the same fixed data and the same fixed-shift edge
states, and it enlarges~$V^{\rm rch}$ only along open pending edges, apart from
the immediate inclusion of the shifted upcluster when it is first reached.
Thus~$V^{\rm rch}$ is always contained in the fixed-shift cluster.

To show that the terminal reached set is the vertex set of the fixed-shift
origin cluster, suppose that the paused exploration completes finitely while
some vertex of that cluster is absent from~$V^{\rm rch}$.  Follow a fixed
open route from the retained downcluster to such a vertex.  The retained
downcluster is present initially, and the whole shifted upcluster is adjoined
as soon as it is reached.  At the first departure of this route from
$V^{\rm rch}$, there is an open edge~$e$ incident to~$V^{\rm rch}$
whose other endpoint has not been reached.  This edge cannot have been
processed by the cluster search, since processing an open edge adds both
endpoints to~$V^{\rm rch}$.  Moreover, $e\notin C_v$ because it is open,
while $e\notin O_v$ because the downcluster is reached initially and the
shifted upcluster is adjoined in full when first reached.  Thus
$e\in\mathcal E_v^\circ\setminus\mathcal Q$.  This edge is pending,
contrary to finite completion.  Hence the reached vertex sets agree.

It remains to identify the incident closed boundary.  At finite completion,
every edge of~$\mathcal E_v^\circ$ incident to~$V^{\rm rch}$ has been
processed, since otherwise it would remain pending.  The processed closed
edges, together with those edges of~$C_v$ incident to the reached set,
therefore give the same incident closed boundary as the fixed-shift search, as in
Lemma~\ref{l.fixed-shift-search-completeness}.

For the converse, suppose that the fixed-shift origin cluster is finite.
Only finitely many edges of~$\mathcal E_v^\circ$ are incident to this cluster,
and every cluster-processing step processes a previously unprocessed one of
them.  Auxiliary seal completions terminate almost surely.  Open escape from
a layer episode is impossible, since its seed is connected to the origin;
on the finite-cluster event each such episode finishes finitely.  An outer
phase continuing indefinitely, or infinitely many successive visits to
monitored layers, would therefore require infinitely many cluster-processing
steps: each layer visit begins with pending cluster work.  Neither can occur,
so the paused exploration completes finitely almost surely. \qed

Together with the preceding identification of the cluster record, the later
sealing estimates require control of where the queried layer record lies
relative to the output seal.

\begin{lemma}[Layer record and output capture]
\label{l.ssr-cluster-completion}
For~$j\in\mathcal I_b$, on a finite completed outcome the fixed retained data
together with the restriction of the paused cluster record to~$E_j$ determine the
open-and-closed record of the fixed-shift origin cluster in~$E_j$.  On acceptance,
this record is precisely the union of the retained layer record and the layer
edges processed by the cluster search.  Whenever
the cluster search processes an edge of~$E_j$, a seal is already in force
there; every such edge has mesh support in the captured mesh set of the
output seal.
\end{lemma}

\noindent{\bf Proof.}
The assertions concerning the layer record follow from
Lemma~\ref{l.multi-layer-regional-completeness} and the fact that on acceptance
both retained clusters belong to the origin cluster.  A seal is in force before
the first edge of~$E_j$ is processed.  The captured scan processes only edges
in its capture.  Exterior edges are processed only in virgin-seed episodes.
At the start of each such episode, the joint unread-field condition gives
$\mathcal U(\mathscr S_j)\cap\mathcal E^{\rm rd}=\varnothing$,
so every exterior edge processed in the episode belongs to its newly queried
record.  Proposition~\ref{p.relative-seal-play} places the complete queried
record of every finite episode inside the output seal.  Seal enlargement
preserves the previous capture, so induction over the layer phases gives the
terminal assertion. \qed

We now turn from the cluster record to freshness of the unread fixed-shift
field.  After the preliminary post-passage replacements,
\eqref{e.paused-initial-product} gives one joint Bernoulli-$p$ product field on
the still-unread edges, across all monitored layers and including edges inside
current captures.  For either remaining seal operation---construction of an
initial seal in an advance-phase layer or finite enlargement after a
virgin-seed episode---let $\mathcal G$ be the incoming history and let~$\mathcal G'$ be obtained by
adjoining exactly the edge states newly queried by the operation.

\begin{lemma}[Freshness under seal operations]
\label{l.layer-internal-seal-subroutines}
Consider the construction of an initial seal in an advance-phase layer
or a finite seal enlargement following a virgin-seed episode.  If,
conditionally on the incoming history~$\mathcal G$, the unread edges of~$\mathcal E_v^\circ$ have
Bernoulli-$p$ product law, then the edges still unread after the operation
retain this product law, and the output is a seal relative to the enlarged
history~$\mathcal G'$.
\end{lemma}

\noindent{\bf Proof.}
If the retained record for initial sealing is empty, the empty seal is
chosen by \hyperref[conv.empty-seal]{convention} and no queries are made,
so that the assertion is immediate.  Otherwise the conditional completion
argument in the proof of Proposition~\ref{p.initial-seal-estimate}
applies to the current product field.  After a finite virgin-seed episode, the seal enlargement is
the stopping completion of Proposition~\ref{p.relative-seal-play}.  Its joint unread-field
conclusion leaves the still-unread complement product and gives a seal
relative to the enlarged history. \qed

For prospective sealing in the next section, we will need the unread-field
product law at stopping times of the paused exploration, not merely at the
completion of individual seal operations.

\begin{lemma}[Product law along the paused exploration]
\label{l.family-paused}
If~$\tau$ is a finite stopping time for the paused filtration, then,
conditionally on~$\mathcal H_\tau^{\rm pa}$, the still-unread edges of
$\mathcal E_v^\circ$ have their joint Bernoulli-$p$ product law.  At a finite
terminal time the same statement holds on~$\Accept_v$.
\end{lemma}

\noindent{\bf Proof.}
Equation~\eqref{e.paused-initial-product} is the base case.  At a cluster
step, an already revealed state leaves the unread field unchanged; otherwise
Lemma~\ref{l.ssr-paused-filtration} makes the queried edge measurable before
its state is read, and revealing that edge leaves the unread remainder product.
At a seal-query step, the next queried edge is selected from the preceding
responses of the stopping construction, so Lemma~\ref{l.freshness} gives the
same conclusion; Lemma~\ref{l.layer-internal-seal-subroutines} records the
relative-seal conclusion when the operation ends.  Induction gives the product
statement at every finite deterministic paused time; for a finite stopping
time~$\tau$, apply it on each event~$\{\tau=t\}\in\mathcal H_t^{\rm pa}$.
At the terminal time, Lemma~\ref{l.multi-layer-regional-completeness} makes
$\Accept_v$ measurable from the terminal paused history, so restricting to
$\Accept_v$ adds no further conditioning. \qed

\section{Prospective passage seals and their cost}
\label{ss.passage-and-breaker-attenuation}

The paused exploration supplies a fixed-shift description from which we can
prepare prospective passage seals before knowing which active slide in the
advance phase will be last.  Intervening breakers control the weighted cost
of these preparations, and at passage we select the appropriate seal.  The
stopping-time product law of Lemma~\ref{l.family-paused} gives the freshness
needed for prospective sealing.

\subsection{Prospective records and seals}
\label{ss.passage-preparation}

Whenever an active slide in the advance phase reaches a monitored layer, we construct the
candidate seal that would be used were that slide the last active one before passage.
If the last active slide in the advance phase reaches the layer, its candidate is
selected pathwise at passage; otherwise no earlier candidate is used and an
initial seal is constructed there.

For a finite open-and-closed edge record
\(\mathfrak R=(D^{\rm op},D^{\rm cl})\) in a layer, put
\[
 |\mathfrak R|_{\rm rec}=|D^{\rm op}|+|D^{\rm cl}|,
 \qquad
 \MeshSupp(\mathfrak R)
 =\{M\in\Mesh:E(M)\cap(D^{\rm op}\cup D^{\rm cl})\ne\emptyset\}.
\]
This record size is distinct from the coarse cardinality~\(|\cdot|_*\)
introduced in Section~\ref{s.overlapping-boxes-boxlets}.  An empty retained
record uses the \hyperref[conv.empty-seal]{empty-seal convention} of
Chapter~\ref{s.seed-exploration-seal-extension}.

For an active slide~$b$ in the advance phase, at height~$H_b$, let
$\RetRec_{b,j}=(D_{{\rm R},b,j}^{\rm op},D_{{\rm R},b,j}^{\rm cl})$ be the
fixed layer record $(O_v\cap E_j,C_v\cap E_j)$ for the proposed head~$v$.
Since $H_b\leq s_j^-$, the height-$H_b$ downcluster has no edge in~$E_j$.
On acceptance, it is the restriction to~$\cS_j$ of the retained record specifying
$\upcluster_{H_b}(S_b)$.  Run the multi-layer paused exploration of
Definition~\ref{d.multi-layer-paused}.  On an accepted outcome, put
\begin{align*}
 D_{{\rm M},b,j}^{\rm op}
 &\, = \, \bigl(E(\csetuthree(S_b))\cap E(\cS_j)\bigr)
      \setminus D_{{\rm R},b,j}^{\rm op}\, ,\\
 \MidRec_{b,j}
 &\, = \, \bigl(D_{{\rm M},b,j}^{\rm op},D_{{\rm M},b,j}^{\rm cl}\bigr)\, ,
\end{align*}
where $D_{{\rm M},b,j}^{\rm cl}$ is the incident closed boundary queried in
exposing these open components.  The open part of~$\MidRec_{b,j}$ is thus
precisely the layer part of the height-$H_b$ midpiece.
For a fixed layer we abbreviate these records to~\(\RetRec_b\) and
\(\MidRec_b\).

\subsubsection{Constructing the prospective seals}

The sealing constructions of Chapter~\ref{s.seed-exploration-seal-extension}
place a seal around a finite prescribed record and leave its unread exterior
conditionally product.  The construction is performed separately in each
fixed-shift experiment; a later active advance step has its own
prospective-seal variable.

Fix a fixed-shift slide~$b$ and a partition layer~$\cS_j$ for which
$H_b\leq s_j^-$.
Let~\(\Reach_{b,j}\) be the event that the paused exploration selects a
pending continuation into~$\cS_j$ while no seal is yet in force there.  On
this event let~\(\cH_{b,j}^{\rm init}\) be the stopped history at the first
such selection, after the continuation has been selected and before either
its edge is read or the initial seal is sampled.  By Lemma~\ref{l.ssr-paused-filtration}, before this
pause no edge of~\(E_j\) outside~\(\RetRec_{b,j}\) has been read.

On~\(\Reach_{b,j}\), if~\(|\RetRec_{b,j}|_{\rm rec}>0\), construct the
initial seal~\(\Seal_{b,j}^{\rm R}\) of
Proposition~\ref{p.initial-seal-estimate}; if the retained record is empty,
take~\(\Seal_{b,j}^{\rm R}=\varnothing\) under the \hyperref[conv.empty-seal]{empty-seal convention}.
Continue the paused exploration.  After every finite layer episode
in~$\cS_j$, use the relative seal update of
Proposition~\ref{p.relative-seal-play}.  The output seal arising on
\(\Reach_{b,j}\cap\Accept_v\) is the \emph{prospective passage seal}
\(\ProsSeal_{b,j}\); set \(\ProsSeal_{b,j}=\varnothing\) on the complementary event.

\begin{lemma}[Product law for advance-phase layer work]
\label{l.ssr-seal-first-product}
On~\(\Reach_{b,j}\), conditionally on~\(\cH_{b,j}^{\rm init}\), the edges
of~\(E_j\) outside~\(\RetRec_{b,j}\) have Bernoulli product law.  The initial
seal is empty under this convention when the retained record is empty, or is a
stopping construction on this field.  Immediately before every later
virgin-seed episode in~\(\cS_j\), its starting seed is measurable from
the stopped history and the still-unread edges of~\(E_j\) have their
joint conditional Bernoulli product law; in particular, so does the
exterior of the current seal.
\end{lemma}

\noindent{\bf Proof.}
At the initial pause, Lemma~\ref{l.family-paused} gives the conditional
product law for the unread layer edges, which by
Lemma~\ref{l.ssr-paused-filtration} are precisely those outside~\(\RetRec_{b,j}\).  Immediately before
each later virgin-seed episode, the starting seed is measurable from the
stopped history by the paused selection rule, and the same
product law follows from Lemma~\ref{l.family-paused} at that stopping
time.  The initial-seal construction and every later finite seal update
are stopping constructions in the same layer field, so the paused
product-law statement applies again when each operation ends. \qed

\begin{proposition}[Prospective passage seal]
\label{p.ssr-prospective-seal}
On~$\Reach_{b,j}\cap\Accept_v$, the output seal~$\ProsSeal_{b,j}$ is
a relative seal for the terminal paused history of the fixed-shift experiment.
Its captured mesh set contains the supports of~$\RetRec_{b,j}$ and of the
complete midpiece record~$\MidRec_{b,j}$.  Conditionally on the terminal
paused history, its exterior has Bernoulli product law.
\end{proposition}

\noindent{\bf Proof.}
On~$\Reach_{b,j}\cap\Accept_v$, the paused exploration completes finitely,
almost surely, by Lemma~\ref{l.multi-layer-regional-completeness}.  The
initial seal captures~$\MeshSupp(\RetRec_{b,j})$ by
Proposition~\ref{p.initial-seal-estimate}; when the retained record is empty,
there is nothing to capture.  Subsequent seal enlargements preserve this
capture.  Lemma~\ref{l.ssr-cluster-completion} identifies the completed
layer record and places the mesh support of every layer edge processed by
the cluster search inside the output capture.  In particular, this covers
the non-retained midpiece record, including its incident closed boundary.

The initial seal is relative to the history at its construction by
Lemma~\ref{l.layer-internal-seal-subroutines}, and each finite enlargement
preserves the relative-seal property.  Between enlargements, captured scans
process only edges inside the current capture, while operations in other
layers leave the exterior of~$\cS_j$ untouched.  Thus the output seal is
relative to the terminal paused history.  Finally,
Lemma~\ref{l.family-paused} gives the joint conditional product law of the
still-unread edges at the terminal stopping time.  Since all revealed
layer edges lie in the output capture, the exterior of~$\ProsSeal_{b,j}$
is still unread and has the asserted product law. \qed

Earlier prospective passage seals and the edge data revealed in constructing
them remain recorded in the augmented history.  At a later active advance
step the construction above uses the new records
\((\RetRec_{b',j},\MidRec_{b',j})\) and defines the new variable
\(\ProsSeal_{b',j}\).  Thus, if~\(b_j^\star\) denotes the last active step in
the advance phase for~\(\cS_j\), the candidate available from that step at
passage is~\(\ProsSeal_{b_j^\star,j}\); no earlier candidate is selected.
The earlier constructions remain auxiliary data in the augmented history,
whose compatibility with the cluster law is established later in
Lemma~\ref{l.ssr-dynamic-projection}.

\subsection{Prospective cost and breaker attenuation}
\label{ss.prospective-preparation-breaker}

The one-slide cost is a non-negative random variable associated with a
proposed slide, designed to control the finite cluster growth and seal growth
produced during its fixed-shift experiment, together with the prospective
sealing needed in the advance phase.  Here we specify it precisely.  The
prospective term is less immediate than the others to estimate: the initial
seals are sampled adaptively across the layers actually reached.  Breaker
attenuation and depth regularity control the total weighted mass of those
layers, giving a uniform moment for prospective preparation.  The same reach
attenuation will later localise the moment of the full one-slide cost.

\subsubsection{Components of the one-slide cost}

\begin{definition}[One-slide cost]
\label{d.ssr-one-slide-charges}
Fix a proposed slide~\(b\) and shift~\(v\).  Let~\(\mathsf{Finite}_b^v\) be
the event that the fixed-shift origin cluster is finite, equivalently that
the fixed-shift cluster exploration of
Section~\ref{ss.successive-northgoing-wavebreakers} terminates finitely.
In particular, \(\Accept_v\subseteq\mathsf{Finite}_b^v\).  We set
\begin{equation}
\label{e.ssr-slide-cost-definition}
 \SlideCost_b
 =
 \one_{\mathsf{Finite}_b^v}
 \bigl(\GainCost_b+\BdCost_b+\PrepCost_b\bigr),
\end{equation}
where the three terms record finite cluster gain, boundary growth and
prospective preparation, as follows.

\smallskip
\noindent\emph{Finite gain.}
Enumerate the virgin-seed episodes of
Definition~\ref{d.multi-layer-paused} in their global paused order.  At
episode \(n\), let~\(x_n\) be the number of vertices reached outside the
capture in force when that episode begins, provided the episode is finite,
and set~\(x_n=0\) on open escape.  If \(\mathcal N_{b,j}\) is the
set of finite episodes assigned to~\(\cS_j\), put
\[
 X_{b,j}=\sum_{n\in\mathcal N_{b,j}}x_n,
 \qquad
 \GainCost_b=\sum_jX_{b,j}^\alpha .
\]
Thus~\(X_{b,j}\) counts only fresh finite gain from the corresponding layer
episodes.  Material reached through the existing capture, including recorded
woven fabric of the current seal, is instead controlled through the seal
cardinality, by Lemma~\ref{l.seal-cardinality-controls-volume}.

\smallskip
\noindent\emph{Boundary growth.}
For the same episode~\(n\), let~\(Z_n\) be the cost~\(Z_s\) supplied by
Lemma~\ref{l.effective-layer-episode} on a finite outcome, and set~\(Z_n=0\)
on open escape.  Writing~\(\mathscr S\) and~\(\mathscr S^+\) for the seals
before and after the episode, fix~\(C_{\rm grow}\) so that
\eqref{e.episode-new-boundary-sites} gives, for every finite episode,
\[
  \bigl|(\mathscr S^+\cap\Mesh)\setminus(\mathscr S\cap\Mesh)\bigr|
  +\bigl(|\mathscr S^+|-|\mathscr S|\bigr)_+
  \le C_{\rm grow}Z_n.
\]
Put
\[
  \BdCost_{b,j}=C_{\rm grow}\sum_{n\in\mathcal N_{b,j}}Z_n,
  \qquad
  \BdCost_b=\sum_j\BdCost_{b,j}.
\]
Thus~\(\BdCost_{b,j}\) controls both the positive seal-cardinality change and
the new mesh-boundary sites created by finite episodes in the layer.
Preliminary replacement of an inherited post-passage seal is not included in
this term.

\smallskip
\noindent\emph{Prospective preparation.}
The initial-seal bound gives an exponential tail for the seal size beyond
a linear allowance for the retained record.  On~\(\Reach_{b,j}\), put
\[
 \InitExc_{b,j}
 =
 \bigl(|\Seal_{b,j}^{\rm R}|
       -A_{\rm bd}|\RetRec_{b,j}|_{\rm rec}\bigr)_+,
\]
where~\(A_{\rm bd}\) is the constant from
Corollary~\ref{c.initial-seal-boundary}, and set~\(\InitExc_{b,j}=0\)
off~\(\Reach_{b,j}\).  Give layer~\(j\) the depth weight
\begin{equation}
\label{e.ssr-depth-weight}
 \mathsf w_j=(1+\delta_j)^{-\sigma_{\rm wt}},
 \quad \text{for }\sigma_{\rm wt}=\frac32,
\end{equation}
and put
\[
 \PrepCost_b
 =\sum_{j:H_b\leq s_j^-}\mathsf w_j\InitExc_{b,j}.
\]

For later use, write
\[
 \SealCost_b
 =
 \one_{\mathsf{Finite}_b^v}
 \bigl(\BdCost_b+\PrepCost_b\bigr).
\]
Thus
\[
 \SlideCost_b
 =
 \one_{\mathsf{Finite}_b^v}\GainCost_b+\SealCost_b.
\]
The dependence on the proposed head~\(v\) is suppressed.  At a dormant
attempted step, set
\[
 \GainCost_b=\BdCost_b=\PrepCost_b=\SealCost_b=\SlideCost_b=0.
\]
\end{definition}

The choice~\(\sigma_{\rm wt}=3/2\) has \(\sigma_{\rm wt}>1\), and with
\(\Bfine=3(d+2)\) it also satisfies
\(\Bfine\alpha-(d-1)-\sigma_{\rm wt}>1\) in every dimension under
consideration.  The latter margin is used in the positive-depth shadow estimate
of Chapter~\ref{s.shadow-bounds-delivery-fine}.

\subsubsection{Breaker attenuation}

Let
\(\mathcal H_{\rm vf}=\{H_i:i\in\ivf\}\) and
\(D_{\rm wb}=3r+b_0+C_0\).  At a northgoing test, the edges of the
breaker window indexed by~\(h\in\mathcal H_{\rm vf}\) have greater
endpoint heights in~\((h+b_0,h+D_{\rm wb}]\), and its vertices have heights
in~\((h+b_0-r,h+D_{\rm wb}]\).  For a monitored pair with~\(H_b\leq s_j^-\), put
\begin{equation}
\label{e.ssr-breaker-count-definition}
 M_{b,j}
 =
 \#\{h\in\mathcal H_{\rm vf}:H_b<h,\
       h+D_{\rm wb}+r<s_j^-\}.
\end{equation}
Thus~\(M_{b,j}\) counts only later breaker windows lying between the present
slide and the south face of~\(\cS_j\).  The extra~\(r\) ensures that
every endpoint of an edge of~\(E_j\) lies north of these windows.

For the monitored pair~\((\cS_j,\ell_j)\), write
\[
 \mathcal B_j
 =\{b:H_b\leq s_j^-\}\cup\{b:s_j^+\leq H_b<\ell_j\}.
\]
For \(\kappa\in(0,1)\), define
\begin{equation}
\label{e.ssr-slide-weight}
 a_{b,j}=\kappa^{M_{b,j}}
 \qquad\text{when }H_b\leq s_j^-,
\end{equation}
set \(a_{b,j}=1\) when \(s_j^+\leq H_b<\ell_j\), and set \(a_{b,j}=0\)
when \(b\notin\mathcal B_j\).  The value of~\(\kappa\) will be chosen in
the final one-slide estimate, after the uniform cost moment is localised by
the probability of reaching an advance-phase layer.

As~\(j\) varies, the breaker windows counted by~\(M_{b,j}\) form nested
initial segments.  Put
\[
 R_b=\max_{j:H_b\leq s_j^-}M_{b,j},
\]
with~\(R_b=0\) if no~\(j\) satisfies~\(H_b\leq s_j^-\), and list their union in
increasing height order as
\(h_{b,1}<\cdots<h_{b,R_b}\), and put
\(\mathscr B_{b,q}=\mathscr B_{h_{b,q}}\).  For each~\(j\), the counted
windows are
\(\mathscr B_{b,1},\ldots,\mathscr B_{b,M_{b,j}}\), all lying wholly below
the south face of~\(\cS_j\).

Condition on an admissible value of~$\cH_{b-1}$ and on a proposed
shift~$v$.  By~\eqref{e.paused-undiscovered-set}, the whole fixed-shift
undiscovered field~$\mathcal E_v^\circ$ has its Bernoulli-$p$ product law
under~$\PP_b^v$.  Run the front-scheduled cluster search of
Section~\ref{ss.successive-stopped-tests} on this field, testing the recovered
windows~$\mathscr B_{b,1},\ldots,\mathscr B_{b,R_b}$ in height order, with
its own query filtration.  At each pre-entry stop,
Lemma~\ref{l.freshness} preserves the product law of the still-unread test
edges.  Every indexed height is tested before its limelight, so
Proposition~\ref{p.breaker-persistence} supplies the thin transported
records required by Proposition~\ref{p.successive-northgoing-wavebreakers}.
We may therefore apply that proposition conditionally on~$\cH_{b-1}$.

This auxiliary breaker search is not appended to the paused filtration.  It is a
parallel exploration of the same fixed-shift configuration, used only to
measure how many breaker windows a finite connection can pass.  In particular,
seal-only revelations from the paused construction do not enter its query
filtration; the two searches are compared pathwise through their common
fixed-shift open configuration.

Let~$A_q^{(b)}$ be the event~\eqref{e.successive-northgoing-events} for the
first~$q$ windows in this auxiliary breaker search.  With
$\mathsf{Finite}_b^v$ as defined above, put
\(T_b=\one_{\mathsf{Finite}_b^v}\max\{0\leq q\leq R_b:A_q^{(b)}\text{ occurs}\}\).
Then
\begin{equation}
\label{e.ssr-crossing-geometric}
 \PP_b^v(T_b\geq k\mid\cH_{b-1})
 \leq q_{\rm row}^k,
 \qquad k\geq0.
\end{equation}

\begin{lemma}[Breaker attenuation before layer entry]
\label{l.ssr-breaker-attenuation}
For every $j$ with $H_b\leq s_j^-$,
\[
 \Reach_{b,j}\cap\mathsf{Finite}_b^v
 \subseteq\{T_b\geq M_{b,j}\}.
\]
Consequently
\begin{equation}
\label{e.ssr-reach-probability}
 \PP_b^v(\Reach_{b,j}\cap\mathsf{Finite}_b^v\mid\cH_{b-1})
 \leq q_{\rm row}^{M_{b,j}}.
\end{equation}
\end{lemma}

\noindent{\bf Proof.}
Suppose that $\Reach_{b,j}\cap\mathsf{Finite}_b^v$ occurs.  We first
extract from the reach an open path to~$\cS_j$ whose terminal subpath in the
undiscovered edge set
starts south of every breaker window counted by~$M_{b,j}$ and ends north of
all of them.  Since the auxiliary breaker search is run in the same fixed
configuration, this will force finite passage through every counted window.

On~$\Reach_{b,j}$, let~$(x,e)$ be the seed--edge pair selected at the first
layer pause represented by the stopped history~$\cH_{b,j}^{\rm init}$.
Thus~$x\in V^{\rm rch}$, $e\in\mathcal U_j^{\rm cur}$ is
incident to~$x$, no prospective seal has yet been installed in~$\cS_j$, and
$e$ has not yet been processed by the cluster search.  By
Lemma~\ref{l.multi-layer-regional-completeness}, the reached set has the
same property as in the fixed-shift exploration: every vertex of~$V^{\rm rch}$ is
joined to the retained downcluster by fixed open retained edges and open edges of~$\mathcal E_v^\circ$ processed
by the cluster exploration.  Trimming such a path at its last vertex in
the downcluster, and using the full boundary recorded in
Definition~\ref{d.prime}, gives a finite open path~$\gamma$ in the fixed-shift
configuration from~$\Top(\downcluster_{H_b})$ to~$x$.  If the path first reaches the shifted upcluster, its fixed-open edges make
all of~$\mathfrak u^v$ reachable without further product queries; these edges
remain part of the same open connection.  No edge discovered solely by a seal query is used in~$\gamma$.

There is nothing to prove when $M_{b,j}$ is zero, so suppose that it is positive.
Before the prospective seal is installed,
\eqref{e.current-layer-regions} gives
$\mathcal U_j^{\rm cur}=E(\cS_j)\cap\mathcal E_v^\circ$, so~$e$ is
undiscovered in the fixed-shift partition.  Every endpoint of an edge of~$E_j$
has height at least~$s_j^--r$.  If $x$ belonged to~$V(\upwell^v)$,
Lemma~\ref{l.ssr-upcluster-bed-interface}, applied to the undiscovered
incident edge~$e$, would put~$x$ in~$\Bed(\upwell^v)$ and hence at height at
most~$\hadj(v)\leq H_b+3r$.  This is impossible because every counted
window has index~$h>H_b$ and
\[
 H_b+3r<h+b_0-r<h+b_0\leq h+D_{\rm wb}<s_j^--r\leq x\cvl .
\]
Here the first strict inequality follows from the schedule gap, and the
last from the extra~$r$ in~\eqref{e.ssr-breaker-count-definition}.  Thus
$x\notin V(\upwell^v)$.

Let~$y$ be the last vertex of~$\gamma$ in
\(\downwell\cup\upwell^v\), and let~$\gamma^\circ$ be the terminal subpath
from~$y$ to~$x$.  If
$y\in V(\downwell)\setminus V(\upwell^v)$, then the first edge of
$\gamma^\circ$ leaves the downcluster, so~$y\in\Top(\downwell)$ and
$y\cvl\leq H_b$.  If $y\in V(\upwell^v)$, the first edge of
$\gamma^\circ$ is an open undiscovered departure from the shifted upcluster;
Lemma~\ref{l.ssr-upcluster-bed-interface} then gives
$y\in\Bed(\upwell^v)$ and~$y\cvl\leq\hadj(v)\leq H_b+3r$.  Therefore, for
every breaker window indexed by~$h$ and counted by~$M_{b,j}$,
\begin{equation}
\label{e.ssr-breaker-height-chain}
 y\cvl\leq H_b+3r
 <h+b_0-r<h+b_0\leq h+D_{\rm wb}
 <s_j^--r\leq x\cvl .
\end{equation}
By the choice of~$y$, the interior of~$\gamma^\circ$ uses no edge of the fixed
upcluster or downcluster; all its open edges belong to the auxiliary Bernoulli
field.

The height chain~\eqref{e.ssr-breaker-height-chain} shows that
$\gamma^\circ$ makes a fresh open crossing of every counted window, from
below its first test layer to above its last.  We check that the
auxiliary front-scheduled search registers these crossings at the prescribed
stops.  Apply the front-exhaustion observation in the paragraph preceding
Lemma~\ref{l.fresh-interface-verification} to the vertex~$x$ outside the
shifted upcluster.  At the southern pre-entry front of each window it
gives a reached vertex incident to an unread edge of the first test layer,
which is~$\mathsf{Ent}_q$; at the northern completion front it gives a
reached vertex outside the shifted upcluster on the far side, which is~$P_q$.
On~$\mathsf{Finite}_b^v$ the successive front stages finish finitely.
Thus~$A_{M_{b,j}}^{(b)}$ occurs and~$T_b\geq M_{b,j}$.  The probability
bound follows from~\eqref{e.ssr-crossing-geometric}.
\qed

\subsubsection{Depth regularity and prospective preparation}

For a fixed proposed slide, the quantity~$M_{b,j}$ counts the breaker
windows that its midpiece must pass before it can reach~$\cS_j$.  The geometric
depth~$\delta_j$ records how far the layer lies below its target; larger depth
comes with more room in the corresponding fine allowance for later accumulation.
Simultaneous control is endangered if many monitored pairs have both small
breaker count and small depth.  The following condition rules out precisely
this crowding.

\begin{definition}[Depth regularity]
\label{d.ssr-depth-regular}
Fix \(C_{\rm reg}<\infty\).  A finite collection of layer--target pairs whose
partition layers are pairwise distinct is
\emph{\(C_{\rm reg}\)-depth regular} if, for every attempted step~$b$ and every
$t,\rho\geq0$,
\begin{equation}
 \#\{j:H_b\leq s_j^-,\ M_{b,j}\leq t,\ \delta_j\leq \rho\}
 \leq C_{\rm reg}(1+t+\rho).
\end{equation}
\end{definition}

Depth regularity is a property of the monitored family relative to the input
and attempted-height schedule, which determine the IVF set and the counts
\(M_{b,j}\).  Lemma~\ref{l.shadow-breaker-count} establishes it for the
selected target families used later.

On a finite fixed-shift outcome, every layer actually reached has breaker count
at most~$T_b$.  Depth regularity may therefore be summed against the decreasing
depth weights.

\begin{lemma}[Weighted mass of the reached layers]
\label{l.ssr-weighted-reached-mass}
For a \(C_{\rm reg}\)-depth-regular collection, with depth weights
\(\mathsf w_j=(1+\delta_j)^{-\sigma_{\rm wt}}\) as in
\eqref{e.ssr-depth-weight},
\begin{equation}
 \sum_{j:H_b\leq s_j^-}
       \mathsf w_j\one_{\Reach_{b,j}\cap\mathsf{Finite}_b^v}
 \leq C(1+C_{\rm reg})(1+T_b),
\end{equation}
where the constant \(C\) is independent of the number of pairs.
\end{lemma}

\noindent{\bf Proof.}
By Lemma~\ref{l.ssr-breaker-attenuation}, every index in the sum satisfies
$M_{b,j}\leq T_b$.  Divide the depths into the ranges
\(2^k\leq1+\delta_j<2^{k+1}\), \(k\geq0\).  Depth regularity bounds the number of
indices in the \(k\)\textsuperscript{th} range by
\(C_{\rm reg}(1+T_b+2^{k+1})\), while each such index has weight at most
\(2^{-\sigma_{\rm wt}k}\).  Hence the weighted sum is at most
\[
 C C_{\rm reg}\sum_{k\geq0}
 \bigl((1+T_b)2^{-\sigma_{\rm wt}k}
       +2^{(1-\sigma_{\rm wt})k}\bigr)
 \leq C(1+C_{\rm reg})(1+T_b),
\]
because~\(\sigma_{\rm wt}>1\).  This proves the claim.
\qed

\begin{lemma}[Exponential moment of prospective preparation]
\label{l.ssr-prospective-preparation-moment}
For a \(C_{\rm reg}\)-depth-regular collection, there are
$\lambda_{\rm prep}>0$ and $C_{\rm prep}<\infty$ such that, for every active
slide~$b$ and every admissible incoming pair~$(y,v)$,
\begin{equation}
\label{e.ssr-preparation-moment}
 \EE_b^v\left[
   \exp\{2\lambda_{\rm prep}
             \one_{\mathsf{Finite}_b^v}\PrepCost_b\}
   \middle|\cH_{b-1}\right]
 \leq C_{\rm prep}.
\end{equation}
The constants may depend on \(C_{\rm reg}\) and are otherwise independent of
the monitored collection.  Quantitatively, for any fixed \(q_0\in(0,1)\), there are
\(A_{\rm prep}(q_0)<\infty\) and a fixed exponent~\(D_{\rm prep}<\infty\)
such that
\begin{equation}
\label{e.ssr-preparation-Creg-dependence}
 \lambda_{\rm prep}^{-1}+C_{\rm prep}
 \leq A_{\rm prep}(q_0)(1+C_{\rm reg})^{D_{\rm prep}}
\end{equation}
holds uniformly over all breaker rows with \(q_{\rm row}\leq q_0\).
\end{lemma}

\noindent{\bf Proof.}
A single initial-seal excess already has a uniform conditional exponential
moment.  We must control how these excesses accumulate over the layers actually reached.
At each adaptive initial-seal time the layer, and hence its depth weight, is
known before any seal is sampled.  The initial-seal moment may therefore be
iterated with these predictable weights, while breaker attenuation and depth
regularity control their total mass.

There is at most one initial-seal operation for each monitored layer during a
fixed slide.  In the filtration of the paused exploration, let~$\tau_q$ be the
time immediately before the $q$\textsuperscript{th} such operation, with
$\tau_q=\infty$ if fewer than~$q$ operations occur.  At
$\tau_q<\infty$, Lemma~\ref{l.ssr-paused-filtration} shows that the selected
pending continuation which triggers the operation and its layer index~$j(q)$
are known, while no edge of~$E_{j(q)}$ outside the fixed retained record has
yet been read.  Thus, when that record is
non-empty, the field on which the initial seal is sampled still has its
conditional product law; when it is empty, the operation makes no seal query.
Put
\[
 \mathcal F_q^-=\mathcal H_{\tau_q}^{\rm pa},
 \qquad
 \mathcal H_\infty^{\rm pa}
 =\sigma\!\left(\bigcup_{t<\infty}\mathcal H_t^{\rm pa}\right).
\]
For each~$q$, set
\[
 (w_q,I_q)=
 \begin{cases}
  (\mathsf w_{j(q)},\InitExc_{b,j(q)}),&\tau_q<\infty,\\
  (0,0),&\tau_q=\infty.
 \end{cases}
\]
Then~$w_q$ and the occurrence of the operation are
$\mathcal F_q^-$-measurable.  Let~$\mathcal G_q$ be the field after the
operation has been completed, with
$\mathcal G_q=\mathcal F_q^-$ in the null case.  Let~$N_{\rm mon}$ be
the number of monitored pairs.  There are at most~$N_{\rm mon}$ operations, so these fields give the deterministic event-time
horizon~$1\leq q\leq N_{\rm mon}$.  Put
$\mathcal G_0=\cH_{b-1}$; then
\[
 \mathcal G_{q-1}\subseteq\mathcal F_q^-\subseteq\mathcal G_q.
\]

If the retained record at the~$q$\textsuperscript{th} operation is empty, then the \hyperref[conv.empty-seal]{convention} gives an empty initial seal and~$I_q=0$.  Otherwise Lemma~\ref{l.ssr-seal-first-product} lets us
apply Corollary~\ref{c.initial-seal-boundary} conditionally on
$\mathcal F_q^-$.  Thus, for fixed constants~$C_{\rm R},c_{\rm R}>0$,
\[
 \PP_b^v(I_q\geq u\mid\mathcal F_q^-)
 \leq C_{\rm R}e^{-c_{\rm R}u},
 \qquad u\geq0.
\]
For~$0\leq\lambda<c_{\rm R}$, integration of this tail gives
\[
 \EE_b^v[e^{\lambda I_q}\mid\mathcal F_q^-]
 \leq
 1+\frac{C_{\rm R}\lambda}{c_{\rm R}-\lambda}
 =:\Phi_{\rm R}(\lambda).
\]
In particular~$\Phi_{\rm R}(\lambda)\searrow1$ and
$\log\Phi_{\rm R}(\lambda)=O(\lambda)$ as~$\lambda\searrow0$.  The same
bound is immediate in the empty-record and null cases.  Choose
$4\lambda<c_{\rm R}$ and put
$c_\lambda=\log\Phi_{\rm R}(4\lambda)$.  Since~$w_q\in[0,1]$ is
$\mathcal F_q^-$-measurable, conditional Jensen for the concave map
$t\mapsto t^{w_q}$ gives
\[
 \EE_b^v\!\left[e^{4\lambda w_q I_q}\mid\mathcal F_q^-\right]
 \leq e^{w_qc_\lambda}.
\]
Consequently
\begin{equation}
\label{e.ssr-preparation-supermartingale}
 \mathcal M_n
 =\exp\left\{4\lambda\sum_{q=1}^n w_qI_q
              -c_\lambda\sum_{q=1}^n w_q\right\},
 \qquad 0\leq n\leq N_{\rm mon},
\end{equation}
is a non-negative supermartingale with respect to
$(\mathcal G_n)_{0\leq n\leq N_{\rm mon}}$, and hence
\(\EE_b^v[\mathcal M_{N_{\rm mon}}\mid\cH_{b-1}]\leq1\).

Put~$\Sigma_b=\sum_{q=1}^{N_{\rm mon}}w_q$.  On~$\mathsf{Finite}_b^v$, every operation
is attached to a reached layer and there is at most one operation in each
layer.  Therefore Lemma~\ref{l.ssr-weighted-reached-mass} gives
\begin{equation}
\label{e.ssr-preparation-weight-mass}
 \Sigma_b
 \leq
 \sum_{j:H_b\leq s_j^-}
      \mathsf w_j\one_{\Reach_{b,j}\cap\mathsf{Finite}_b^v}
 \leq C(1+C_{\rm reg})(1+T_b).
\end{equation}
On the same event,~$\PrepCost_b=\sum_{q=1}^{N_{\rm mon}} w_qI_q$, so
\eqref{e.ssr-preparation-supermartingale} yields
\[
 e^{2\lambda\PrepCost_b}
 =\mathcal M_{N_{\rm mon}}^{1/2}e^{c_\lambda\Sigma_b/2}.
\]
Using this identity only on~$\mathsf{Finite}_b^v$, conditional
Cauchy--Schwarz, \eqref{e.ssr-preparation-weight-mass} and the
supermartingale bound give
\begin{align}
 \EE_b^v\!\left[
   e^{2\lambda\one_{\mathsf{Finite}_b^v}\PrepCost_b}
   \middle|\cH_{b-1}\right]
 &\leq 1+
 \EE_b^v\!\left[
   \one_{\mathsf{Finite}_b^v}e^{c_\lambda\Sigma_b}
   \middle|\cH_{b-1}\right]^{1/2} \notag\\
 &\leq 1+e^{c_\lambda C(1+C_{\rm reg})/2}
 \EE_b^v\!\left[
   \one_{\mathsf{Finite}_b^v}
   e^{c_\lambda C(1+C_{\rm reg})T_b}
   \middle|\cH_{b-1}\right]^{1/2}.
\label{e.ssr-preparation-before-geometric}
\end{align}
Here the first~$1$ accounts for the complement of
$\mathsf{Finite}_b^v$.

It remains to choose~$\lambda$.  By
\eqref{e.ssr-crossing-geometric}, for every~$\eta>0$ such that
$q_{\rm row}e^\eta<1$,
\begin{align*}
 \EE_b^v[e^{\eta T_b}\mid\cH_{b-1}]
 &\leq
 1+(e^\eta-1)
   \sum_{k\geq1}e^{\eta(k-1)}
      \PP_b^v(T_b\geq k\mid\cH_{b-1})\\
 &\leq
 1+(e^\eta-1)
   \sum_{k\geq1}e^{\eta(k-1)}q_{\rm row}^{k}<\infty.
\end{align*}
Since~$c_\lambda\searrow0$ with~$\lambda$, choose~$\lambda>0$ so small
that
$c_\lambda C(1+C_{\rm reg})< -\tfrac12\log q_{\rm row}$.
Then~\eqref{e.ssr-preparation-before-geometric} is uniformly bounded, which
proves~\eqref{e.ssr-preparation-moment} after decreasing~$\lambda$ if
necessary and calling the resulting value~$\lambda_{\rm prep}$.

For the quantitative assertion, fix~$q_0<1$ and suppose
$q_{\rm row}\leq q_0$.  The constant in
Lemma~\ref{l.ssr-weighted-reached-mass} is linear in
$1+C_{\rm reg}$, while
$c_\lambda=\log\Phi_{\rm R}(4\lambda)=O(\lambda)$ as
$\lambda\searrow0$.  Taking
\(\lambda=c(1+C_{\rm reg})^{-1}\), with~$c>0$ sufficiently small depending
only on~$q_0$ and the fixed seal constants, makes the last geometric moment
uniformly bounded over all~$q_{\rm row}\leq q_0$.  Thus
$\lambda_{\rm prep}^{-1}$ and the resulting moment constant are bounded by a
fixed power of~$1+C_{\rm reg}$, proving
\eqref{e.ssr-preparation-Creg-dependence}.
\qed

\subsection{Selection of the passage seal}
\label{ss.passage-seal-selection}

The prospective-preparation estimate controls the candidates prepared in the advance phase one slide at a time.  At passage, we take the candidate prepared by the actual last active slide as the passage seal, when there is one; otherwise we construct the passage seal directly from the retained record.  In either case the resulting passage seal supplies the first finite captured mesh region for post-passage monitoring.

\subsubsection{The retained passage record}

Let $\mathfrak C_0^{\rm rec}$ be the complete cluster record of the initial state~$S_0$: the open
edges of $\csetuthree(S_0)$ together with their incident closed boundary.  For
the monitored layers, let $\mathcal U_0$ be the layer edges outside
this record, and write~$\omega_0$ for the input bond configuration
underlying~$S_0$.

\begin{lemma}[Input-record exterior]
\label{l.ssr-input-record-exterior}
Conditional on $\mathfrak C_0^{\rm rec}$ and on the attempted-height schedule,
\[
 \mathcal L\left(
   \omega_0|_{\mathcal U_0}
   \middle|\mathfrak C_0^{\rm rec},\overline H
 \right)
\, = \, \bigotimes_{e\in\mathcal U_0}{\rm Bernoulli}(p)\, .
\]
The same assertion holds after choosing any collection of targets by a rule measurable
with respect to the input-cluster record and the schedule.
\end{lemma}

\noindent{\bf Proof.}
The complete record prescribes the cluster edges to be open and their incident
boundary edges to be closed.  These closed edges already isolate the recorded
cluster from the rest of the graph, so specifying the complete record
imposes no further condition on the exterior.  Product disintegration gives the displayed product identity.  The schedule is independent of the input,
and the IVF target selection is measurable from the recorded cluster and the
schedule.
\qed

Accordingly, the histories used below reveal the input only through its
complete cluster record, leaving the ambient input exterior unread.  This
point is relevant only when the advance phase contains no active slide.

Fix one of the monitored layer--target pairs~$(\cS,\ell)$, with
$\cS=\slab_{s^-,s^+}$ and geometric depth \(\delta=\ell-s^+\).  Let
\(b_1<\cdots<b_{n_{\rm adv}}\) be the indices of the active steps in the advance phase for
this layer, put \(h_i^{\rm adv}=H_{b_i}\), and write
\(S_i^{\rm adv}=S_{b_i}\) for the state after the \(i\)\textsuperscript{th} such step, with
\(S_0^{\rm adv}=S_0\).  When \(n_{\rm adv}\geq1\), put
\(h_{\rm adv}=h_{n_{\rm adv}}^{\rm adv}\) and let \(\tau_i\) be the translation map
selected at step~\(b_i\); write
\(\tau_{\rm adv}=\tau_{n_{\rm adv}}\circ\cdots\circ\tau_1\) and
\(S_{\rm adv}=S_{n_{\rm adv}}^{\rm adv}\).  When
\(n_{\rm adv}=0\), put \(\tau_{\rm adv}=\operatorname{id}_{\Z^d}\) and
\(S_{\rm adv}=S_0\).

\begin{lemma}[Transport of the retained upcluster]
\label{l.ssr-retained-transport}
If $n_{\rm adv}\ge1$, then
\begin{equation}
\label{e.ssr-retained-transport}
 \upcluster_{h_{\rm adv}}(S_{\rm adv})\, = \,
 \tau_{\rm adv}\bigl(\upcluster_{h_{\rm adv}}(S_0)\bigr)\, .
\end{equation}
The fresh edges sampled at step~$b_{n_{\rm adv}}$ are disjoint from the upcluster in
\eqref{e.ssr-retained-transport}.
\end{lemma}

\noindent{\bf Proof.}
For $i<n_{\rm adv}$, Lemma~\ref{l.higher-ucs-shift}, applied with
$h_i^{\rm adv}<h_{\rm adv}$, gives
\[
 \upcluster_{h_{\rm adv}}(S_i^{\rm adv})
 =\tau_i\bigl(\upcluster_{h_{\rm adv}}(S_{i-1}^{\rm adv})\bigr).
\]
At $i=n_{\rm adv}$, the retained-data identity in the slide law gives the
same equality with $h_{\rm adv}=h_{n_{\rm adv}}^{\rm adv}$.  Composition
proves~\eqref{e.ssr-retained-transport}.  The last assertion
follows because the output law fixes the shift class of the height-$h_{\rm adv}$
upcluster.
\qed

If \(n_{\rm adv}\ge1\), put \(b^\star=b_{n_{\rm adv}}\) and set
\(\RetRec=\RetRec_{b^\star}\) and \(\MidRec=\MidRec_{b^\star}\).  Thus
\(\RetRec=(D_{\rm R}^{\rm op},D_{\rm R}^{\rm cl})\) is the restriction to
\(\cS\) of the retained record specifying
\(\upcluster_{h_{\rm adv}}(S_{\rm adv})\); its closed part consists of the
incident upper-space boundary edges used to specify that component.  If
\(n_{\rm adv}=0\), let~\(\RetRec\) be the restriction to~\(\cS\) of
\(\mathfrak C_0^{\rm rec}\) and take~\(\MidRec\) to be empty.  We also
record the retained vertex trace in the layer:
\[
 V_{\rm R}\, = \,
 \begin{cases}
 V(\upcluster_{h_{\rm adv}}(S_{\rm adv}))\cap V(\cS),&n_{\rm adv}\geq1,\\
 V(\csetuthree(S_0))\cap V(\cS),&n_{\rm adv}=0.
 \end{cases}
\]
Thus \(V_{\rm R}\) also includes retained vertices at which the cluster meets
the vertex set of the layer without using an open edge of~\(E(\cS)\).

\begin{lemma}[Retained-record bound]
\label{l.ssr-retained-band}
If $\ell$ is initially very fine, then
\[
 |V_{\rm R}|+|\RetRec|_{\rm rec}\, \leq \, C\avf(\delta)\, .
\]
\end{lemma}

\noindent{\bf Proof.}
Write \(\Delta_{\rm adv}=(\tau_{\rm adv}(0))\cvl\).  Since
\(\tau_{\rm adv}(\head(S_0))=\head(S_{\rm adv})\) and both head heights lie
in the window $[u,u+3r]$ fixed in Section~\ref{s.onetime},
\[
 |\Delta_{\rm adv}|=|\hmax(S_{\rm adv})-\hmax(S_0)|\leq3r.
\]  Every vertex of
\(\tau_{\rm adv}^{-1}(V_{\rm R})\) has height in
\[
 J=[s^- -\Delta_{\rm adv}-r,s^+-\Delta_{\rm adv}].
\]
Since $s^+=\ell-\delta$ and $s^+-s^-\leq\Dmesh$, this interval meets only a
bounded number, depending on the standing geometry, of the unit bands in
Remark~\ref{r.new-adjacent-band-form}.  Any southside index~$q$ of such a
band satisfies \(|q-\delta|\leq C\), while a northside index can occur only
when \(\delta,q\leq C\).

For $n_{\rm adv}\geq1$, Lemma~\ref{l.ssr-retained-transport} gives
\(\tau_{\rm adv}^{-1}(V_{\rm R})\subseteq V(\csetuthree(S_0))\); when
$n_{\rm adv}=0$, this follows directly from the definition of~\(V_{\rm R}\).
The strengthened very-fine allowances and
Remark~\ref{r.new-adjacent-band-form} therefore give
\[
 |V_{\rm R}|\leq C\avf(\delta).
\]
Every edge of \(D_{\rm R}^{\rm op}\cup D_{\rm R}^{\rm cl}\) is incident
to a vertex of~\(V_{\rm R}\): this is immediate for the retained open edges,
while \(D_{\rm R}^{\rm cl}\) is their retained component's incident boundary
inside the layer edge set.  The degree bound hence gives
\(|\RetRec|_{\rm rec}\leq2d\,|V_{\rm R}|\leq C\avf(\delta)\).
\qed

Between consecutive active steps in the advance phase for~$\cS$, the current state of every
edge represented in the midpiece record outside the retained record is
resampled.  The retained part at the later step is the restriction of the
nested higher upcluster, transported by the intervening shift.  Thus only the
retained record and the midpiece record of the last active step in the advance phase
contribute to the passage record; if there is no active step in the advance phase, the
midpiece record is empty.

\begin{lemma}[Cluster decomposition at passage]
\label{l.ssr-passage-decomposition}
In the last pre-passage cluster state~$S_{\rm adv}$,
\[
 E(\csetuthree(S_{\rm adv}))\cap E(\cS)
\, = \, D_{\rm R}^{\rm op}\,\sqcup\,D_{\rm M}^{\rm op}\, .
\]
Here $D_{\rm M}^{\rm op}$ is the open part of~$\MidRec$, and is empty
when $n_{\rm adv}=0$.
\end{lemma}

\noindent{\bf Proof.}
For $n_{\rm adv}\geq1$, identity~\eqref{e.ssr-advance-cluster-decomposition} of
Lemma~\ref{l.ssr-two-phase-cluster} partitions the layer cluster into its
height-$h_{\rm adv}$ upcluster and midpiece parts.  These are respectively
$D_{\rm R}^{\rm op}$ and $D_{\rm M}^{\rm op}$ by the definitions above and
Lemma~\ref{l.ssr-cluster-completion}.  For $n_{\rm adv}=0$, the identity is the
definition of the complete input-cluster record~$\RetRec$.
\qed

For the monitored family, write~$\RetRec_j$ and~$\MidRec_j$ for these
passage records in layer~$\cS_j$.

\subsubsection{Selection and passage sealing}

On~\(\Reach_{b,j}\cap\Accept_v\), set
\[
 \ProsExc_{b,j}
 =
 \bigl(|\ProsSeal_{b,j}|
       -A_{\rm bd}|\RetRec_{b,j}|_{\rm rec}\bigr)_+,
\]
and put \(\ProsExc_{b,j}=0\) otherwise.  Iterating the boundary increment in
\eqref{e.episode-new-boundary-sites} gives
\begin{equation}
\label{e.ssr-prospective-excess}
 \ProsExc_{b,j}
 \leq
 \InitExc_{b,j}+\BdCost_{b,j}.
\end{equation}

\begin{corollary}[Weighted excess of prospective seals selected at passage]
\label{c.ssr-last-advance-selection}
Let $b_j^\star$ be the last active step in the advance phase relative to~$\cS_j$.
If there is no such step, set $b_j^\star=0$, let $\Reach_{0,j}=\varnothing$,
and put $\ProsExc_{0,j}=0$ and $\MidRec_{0,j}=(\varnothing,\varnothing)$.  Put
\(\PassPrep_j=\ProsExc_{b_j^\star,j}\).
Then, on every trajectory for which the fixed-shift outcome at each
$b_j^\star\geq1$ is accepted, and in particular on the accepted slide
trajectory,
\begin{equation}
\label{e.ssr-selected-passage-energy}
 \sum_j\mathsf w_j\PassPrep_j
 \leq
 \sum_b\SealCost_b.
\end{equation}
\end{corollary}

\noindent{\bf Proof.}
If~\(b_j^\star=0\), then~\(\PassPrep_j=0\) by definition.  For
\(b_j^\star\geq1\), \eqref{e.ssr-prospective-excess} gives
\(\mathsf w_j\PassPrep_j\leq
\mathsf w_j\InitExc_{b_j^\star,j}+\BdCost_{b_j^\star,j}\),
because $\mathsf w_j\leq1$.  The first term is one of the non-negative terms of
$\PrepCost_{b_j^\star}$, and the second is one of the terms of~$\BdCost_{b_j^\star}$.
On such a trajectory, every slide contributing to the displayed costs is
finite, so
\(\SealCost_b=\BdCost_b+\PrepCost_b\).
Summing over the indices with~\(b_j^\star\geq1\) proves the claim.
\qed

At passage, the advance-stage preparation is converted into the passage seal whose captured mesh region is carried into post-passage monitoring.  Process passage operations chronologically,
using a fixed deterministic order for layers passed at the same attempted
height, and let \(\cH_j^{\rm pre}\) denote the full multi-layer history
immediately before the passage operation for~\(\cS_j\).  If
\(\Reach_{b_j^\star,j}\) occurred, the operation selects the prospective seal
already constructed at step~\(b_j^\star\); thus this history includes that
construction, its auxiliary choices, and all earlier query and sealing
records.  If \(\Reach_{b_j^\star,j}\) did not occur, the operation instead samples
an initial seal in the last pre-passage configuration~\(S_{\rm adv}\), before
the first post-passage fixed-shift experiment, and \(\cH_j^{\rm pre}\) is
taken before this new sampling.  When \(b_j^\star=0\), no active step in the advance phase has occurred for
\(\cS_j\), and the restriction of this history to~\(E_j\) contains no input
information beyond the complete input-cluster record~\(\mathfrak C_0^{\rm rec}\).

The unread layer field retains its product law until passage.  After the last
active step~\(b_j^\star\) in the advance phase, every remaining step
in that phase is dormant, leaving the cluster state unchanged; the fine tests determining
these dormant steps depend only on that state and the fixed attempted-height
schedule.  Thus, when~\(b_j^\star\geq1\), the conditional Bernoulli product
law at the terminal paused history of that step, given by
Lemma~\ref{l.family-paused}, persists until passage.  When~\(b_j^\star=0\),
the same law follows from Lemma~\ref{l.ssr-input-record-exterior}.
Earlier passage operations in other layers query disjoint edge sets, so in
both cases the still-unread edges of~\(E_j\) have conditional Bernoulli
product law given~\(\cH_j^{\rm pre}\).

\begin{proposition}[Passage seal]
\label{p.ssr-passage-seal}
Suppose that the target~$\ell_j$ is initially very fine.  The passage operation for~$\cS_j$
produces a finite seal~$\Seal_j^0$, possibly empty under the \hyperref[conv.empty-seal]{empty-seal convention}.  There are constants \(A_{\rm pass},C,c\in(0,\infty)\) such that, on defining
\(\PassExc_j=(|\Seal_j^0|-A_{\rm pass}\avf(\delta_j)-\PassPrep_j)_+\), one has
\begin{align}
 |\Seal_j^0|
 &\, \leq \, A_{\rm pass}\avf(\delta_j)+\PassPrep_j+\PassExc_j\, ,
 \\
 \PP_{\rm fs}(\PassExc_j\, \geq \, u\mid\cH_j^{\rm pre})
 &\, \leq \, Ce^{-cu},\qquad u\, \geq \, 0\, .
\end{align}
Here \(\PP_{\rm fs}\) denotes the conditional fixed-shift law given
\(\cH_j^{\rm pre}\), with the realised admissible shift at the last active
step in the advance phase held fixed when such a step exists.  A new seal, when needed, is
constructed from the still-unread layer field; if~\(b_j^\star=0\), this is
the input exterior of Lemma~\ref{l.ssr-input-record-exterior}.  The constants
are uniform in the realised shift.  The same conditional bound holds for
the augmented history accompanying the accepted slide dynamics~\(\Qdyn\).

Moreover,~\(\MCap(\Seal_j^0)\) contains the mesh support of the
complete open-and-closed layer record of the last pre-passage cluster state.
\end{proposition}

\noindent{\bf Proof.}
Lemma~\ref{l.ssr-retained-band} gives
\(|\RetRec_j|_{\rm rec}\leq C_{\rm ret}\avf(\delta_j)\).

Suppose first that~\(\Reach_{b_j^\star,j}\) occurred.  Take
\(\Seal_j^0=\ProsSeal_{b_j^\star,j}\).  Proposition~\ref{p.ssr-prospective-seal} gives the capture used in the
final assertion below.  The definitions of
\(\ProsExc_{b_j^\star,j}\) and~\(\PassPrep_j\) give
\[
 |\Seal_j^0|
 \le A_{\rm bd}|\RetRec_j|_{\rm rec}+\PassPrep_j.
\]
Choose~\(A_{\rm pass}\ge A_{\rm bd}C_{\rm ret}\).  Then~\(\PassExc_j=0\).

Suppose instead that~\(\Reach_{b_j^\star,j}\) did not occur.  No edge of~\(E_j\)
outside~\(\RetRec_j\) is then processed during the last active slide in the advance phase, so~\(\MidRec_{b_j^\star,j}\) is empty and
\(\PassPrep_j=0\).  If~\(|\RetRec_j|_{\rm rec}=0\), take
\(\Seal_j^0=\varnothing\) by the \hyperref[conv.empty-seal]{empty-seal convention}.  Otherwise construct the initial seal of
Proposition~\ref{p.initial-seal-estimate} from~\(\RetRec_j\).  The still-unread field used in this construction has Bernoulli product
law conditional on~\(\cH_j^{\rm pre}\), as established above.
Corollary~\ref{c.initial-seal-boundary} therefore gives
\[
 \PP_{\rm fs}\left(
   \bigl(|\Seal_j^0|-A_{\rm bd}|\RetRec_j|_{\rm rec}\bigr)_+\ge u
   \mid\cH_j^{\rm pre}
 \right)\le Ce^{-cu}.
\]
The same choice of~\(A_{\rm pass}\) bounds~\(\PassExc_j\) by the random
variable in the last display.  The product-law argument preceding the
proposition also applies along the accepted slide trajectory.  Since the seal-size estimate is uniform in the
realised shift, the same conditional bound holds for its augmented history.

It remains to verify the assertion about the complete layer record.  If
\(\Reach_{b_j^\star,j}\) occurred, then~\(\Seal_j^0\) is the output seal
of the completed paused exploration.  Its initial construction captures the
fixed retained layer record~\(\RetRec_j\), and
Lemma~\ref{l.ssr-cluster-completion} puts every additional layer edge
queried by the cluster inside the output capture.  The same lemma identifies
the fixed retained data together with these queried edges as the
complete open-and-closed layer record of the resulting cluster.  If
\(\Reach_{b_j^\star,j}\) did not occur, no edge of~\(E_j\) outside
\(\RetRec_j\) is processed, so that record is already complete
there and the initial seal captures it.  This includes the case
\(b_j^\star=0\), where~\(\RetRec_j\) is the complete input-cluster record.
Thus in every case the passage seal captures the complete layer record of the
last pre-passage state.  This proves the proposition.
\qed

\section{Post-passage seal replacement}
\label{ss.post-passage-resealing}

The paused construction of Section~\ref{ss.multi-layer-paused} begins each
active post-passage step by replacing the inherited seal in the fresh
fixed-shift field.  We now justify that replacement and show that the
incoming-capture condition~\eqref{e.ssr-incoming-capture} propagates from
the passage seal through every accepted post-passage step.  We then
project the augmented construction to the underlying fixed-shift cluster law
and combine the preceding estimates into the one-slide exponential moment.

\subsection{Fresh replacement and inherited capture}

\begin{lemma}[Fresh post-passage replacement]
\label{l.ssr-fresh-postpassage-rebuild}
Fix an active fixed-shift step~$b$, an admissible proposed head~$v$,
and $j\in\mathcal I_b^{\rm post}$.  Let
$\mathcal G\supseteq\sigma(\cH_{b-1},v)$ be an incoming paused history such
that, conditionally on~$\mathcal G$, the still-unread edges of
$\mathcal E_v^\circ$ have their joint Bernoulli-$p$ product law and no edge of
$E_j\cap\mathcal E_v^\circ$ has yet been read in the present fixed-shift
step.  Suppose also that the incoming-capture condition
\eqref{e.ssr-incoming-capture} holds for~$j$.

Starting from the carrier~$D_{b,j}^{\rm repl}$ of
\eqref{e.postpassage-rebuild-carrier}, run the empty-seal completion of
Lemma~\ref{l.conditional-rooted-completion} in the fixed-shift field on~$E_j$,
and let~$\mathcal G'$ be~$\mathcal G$ enlarged by the edges queried in this
replacement and their revealed states.  Denote the output by
$\mathscr S_{b,j}^{\rm repl}$.  Then
\begin{equation}
\label{e.postpassage-rebuild-capture}
 \MCap(\mathscr S_j^-)
 \subseteq D_{b,j}^{\rm repl}
 \subseteq \MeshInt(\mathscr S_{b,j}^{\rm repl}),
\end{equation}
and~$\mathscr S_{b,j}^{\rm repl}$ is a relative~$\mathcal G'$-seal.
Conditionally on~$\mathcal G'$, the edges of~$\mathcal E_v^\circ$ still unread
have their joint Bernoulli-$p$ product law.
\end{lemma}

\noindent{\bf Proof.}
By~\eqref{e.postpassage-fixed-record-inherited},
\[
 \MeshSupp\bigl((O_v\cup C_v)\cap E_j\bigr)
 \subseteq \MCap(\mathscr S_j^-)
 \subseteq D_{b,j}^{\rm repl}.
\]
The carrier~$D_{b,j}^{\rm repl}$ is~$\mathcal G$-measurable.  By hypothesis, every still-unread
edge of the fixed-shift field, in particular the edges of~$E_j$ available to
the completion outside the carrier, has the conditional product law given
$\mathcal G$.  The completion tests query no edge of
$\bigcup_{M\in D_{b,j}^{\rm repl}}E(M)$.  The empty-seal instance of
Lemma~\ref{l.conditional-rooted-completion}, with
$D_2=D_{b,j}^{\rm repl}$, gives the second inclusion in
\eqref{e.postpassage-rebuild-capture} and the relative-seal assertion.
Lemma~\ref{l.freshness} leaves the still-unread complement with its joint
conditional Bernoulli-$p$ product law. \qed

\begin{lemma}[Inherited layer record is captured]
\label{l.ssr-inherited-record-captured}
Fix a monitored layer~$\cS_j$ and an active post-passage step~$b$, on a
trajectory for which the preceding active fixed-shift outcomes have been
accepted.  Let~$\mathscr S_j^-$ be the seal produced for~$\cS_j$ by the
preceding active post-passage experiment, with the passage seal used if no
earlier active post-passage step has occurred.  Then the incoming-capture
condition~\eqref{e.ssr-incoming-capture} holds for~$j$.
\end{lemma}

\noindent{\bf Proof.}
Before the first active post-passage step, the incoming state is the last
pre-passage state, since intervening dormant steps change nothing.  The final
assertion of Proposition~\ref{p.ssr-passage-seal} puts the mesh support of its
complete layer record inside~$\MCap(\Seal_j^0)$, giving the base case.

For the induction step, suppose that the complete incoming layer record at one
active post-passage step lies in the capture of the preceding seal.  This
supplies the capture hypothesis of Lemma~\ref{l.ssr-fresh-postpassage-rebuild};
the product-field and unread-layer hypotheses hold for the fixed-order
preliminary replacements by the construction leading to
\eqref{e.paused-initial-product}.  The replacement seal therefore
contains the inherited capture in its mesh interior:
\[
 \MCap(\mathscr S_j^-)
 \subseteq \MeshInt(\mathscr S_{b,j}^{\rm repl}).
\]
During the ensuing paused exploration,
Lemma~\ref{l.ssr-cluster-completion} puts every additional layer edge queried
by the cluster inside the output capture, while seal enlargement preserves
the preceding capture.  The same lemma identifies the fixed retained data
together with these queried edges as the complete open-and-closed layer record
of the output cluster.  On the accepted outcome this is the state carried
forward to the next active step, so the output seal satisfies the induction
hypothesis there. \qed

The quantity~$\BdCost_b$ records seal extensions and new mesh-boundary sites
arising from finite episodes in the present cluster search.  It does not
include increases in seal size caused by replacing inherited seals before
that search begins.  The outward displacement arising from successive
replacements is controlled in Section~\ref{ss.rebuilt-seals-to-limelight}
through the geometry of the successive mesh boundaries.

\subsection{Projection to the fixed-shift law}

The augmented history retains seal and seal-query information that is absent
from the cluster history.  At a fixed proposed head, the auxiliary
constructions reveal states of the reconstructed configuration without
changing its law.  We state this compatibility before restricting to the
finite-cluster event required for an accepted slide.

\begin{lemma}[Fixed-shift cluster projection]
\label{l.ssr-dynamic-projection}
Fix an active step~$b$ and an admissible incoming pair~$(y,v)$.  Under the
augmented fixed-shift kernel~$\PP_b^v$, perform the preliminary post-passage
replacements and then the paused cluster exploration of
Definition~\ref{d.multi-layer-paused}.  If~$\omega''_v$ denotes the
reconstructed configuration, then, for every~$A\in\mc B$,
\begin{equation}
\label{e.ssr-fixed-shift-cluster-projection}
 \PP_b^v(y,\{\omega''_v\in A\})
 =\Rlaw^v_{\uclust,\dclust}(A),
\end{equation}
where $(\uclust,\dclust)$ is the retained pair determined by the incoming
cluster state in~$y$.
\end{lemma}

\noindent{\bf Proof.}
The retained open and closed edges are fixed by the incoming data, and
all edges in~$\mathcal E_v^\circ$ receive their states from the fresh
Bernoulli-$p$ product field.  The preliminary replacements and later seal
operations only reveal states of this same field; they do not change the
reconstructed configuration, even when the paused exploration does not
complete.  This proves~\eqref{e.ssr-fixed-shift-cluster-projection}.
On finite completed outcomes, Lemma~\ref{l.multi-layer-regional-completeness}
also identifies the paused cluster record with that of the ordinary fixed-shift
search. \qed

On~$\Accept_v$ the origin cluster is finite, and the paused exploration
completes almost surely by Lemma~\ref{l.multi-layer-regional-completeness}.
Restricting~\eqref{e.ssr-fixed-shift-cluster-projection} to this event and
projecting the finite cluster gives the ordinary accepted fixed-shift law.

\subsection{The one-slide exponential moment}

We now combine the bounds for finite gain, boundary growth and
prospective preparation to obtain a conditional exponential-moment
bound for their sum, the one-slide cost~\(\SlideCost_b\).

\begin{proposition}[One-slide exponential moment for layer accrual]
\label{p.ssr-fixed-shift-input}
For a finite \(C_{\rm reg}\)-depth-regular collection of layer--target pairs,
there are moment constants \(C,\lambda>0\) and an attenuation
base~\(\kappa\in(0,1)\) such that, for
every active slide~\(b\) and every admissible incoming pair~\((y,v)\),
\begin{equation}
\label{e.ssr-fixed-shift-moment}
 \EE_b^v[e^{2\lambda\SlideCost_b}\mid\cH_{b-1}]
 \leq
 \exp\left\{C\sum_{j:b\in\mathcal B_j}a_{b,j}\right\}.
\end{equation}
Once \(C_{\rm reg}\) is fixed, the moment constants are independent of the
monitored collection, the conditioning histories and the positive-depth
allowance coefficient~\(\Cfine\).

Quantitatively, fix~\(q_0\in(0,1)\).  Uniformly over breaker rows with
\(q_{\rm row}\leq q_0\), the moment constants may be chosen so that
\begin{equation}
\label{e.ssr-input-Creg-dependence}
 \lambda^{-1}+C
 \leq A(q_0)(1+C_{\rm reg})^D
\end{equation}
for some~\(A(q_0)<\infty\) and~\(D<\infty\) independent of~\(C_{\rm reg}\),
while the attenuation base may be taken as
\begin{equation}
\label{e.ssr-final-attenuation-base}
 \kappa=q_{\rm row}^{1/2}\leq q_0^{1/2}<1.
\end{equation}
Since~\(i_{\rm br}\geq2\) and
\(q_{\rm row}=(1-p_{\rm arr})^{i_{\rm br}/2}\), the choice
\(q_0=1-p_{\rm arr}<1\) shows in particular that the moment constants can be
chosen independently of~\(i_{\rm br}\); the attenuation base
decreases as the row is lengthened.
\end{proposition}

\noindent{\bf Proof.}
Lemma~\ref{l.ssr-prospective-preparation-moment} gives a conditional
exponential moment for the preparation contribution.  We first obtain the
corresponding uniform moment for finite gain and boundary growth.

The proof of the Encloser theorem in
Section~\ref{ss.one-life-percolation-control} uses
Lemma~\ref{l.adaptive-sealed-chain} to control costs accumulated over an
adaptive succession of virgin-seed episodes.  We use the same mechanism here,
with episodes from several monitored layers and with both finite open gain
and seal enlargement contributing to their costs.

Enumerate the episodes in their global paused order, and, as in the Encloser
proof, index the full paused history by successive decisions and queries,
writing~$\mathcal G_t$ for the history at step~$t$ and~$\mathcal G_\infty$
for its limit.  Let~$T_n$ be the stopping time just after the seed and
launch edge of the $n$\textsuperscript{th} episode have been selected, but
before the edge is queried; put~$T_n=\infty$ if no such episode is launched.
Set~$\cF_{n-1}=\mathcal G_{T_n}$.  These fields are nested, and each episode
outcome is~$\cF_n$-measurable, being known before the next launch or from
the limiting history if there is none.  At an actual launch,
Lemma~\ref{l.ssr-paused-filtration} makes the seed and launch edge measurable
from~$\cF_{n-1}$, and Lemma~\ref{l.family-paused} gives the joint
conditional product law on the still-unread layer edges.  Thus
Lemma~\ref{l.effective-layer-episode} applies uniformly in the incoming
seal, monitored layer and seed.

For the $n$\textsuperscript{th} actual episode, set~$W_n^{\rm ep}=\infty$
on open escape and, on a finite outcome, put
\[
 W_n^{\rm ep}=\bigl(x_n^\alpha+C_{\rm grow}Z_n\bigr)^{1/\alpha}.
\]
If no $n$\textsuperscript{th} episode is launched, set
$W_n^{\rm ep}=\infty$ and~$x_n=Z_n=0$, with the same cemetery convention
for all later indices.  The escape-probability bound of
Lemma~\ref{l.effective-layer-episode} applies at every actual launch;
in the cemetery case the escape probability is one.  After decreasing
the episode lemma's exponential parameter by a fixed factor to absorb
$C_{\rm grow}$, its exponential-moment estimate and conditional Markov's
inequality give the finite-value stretched-exponential tail required by
Lemma~\ref{l.adaptive-sealed-chain}.  Applied to these stopped histories,
that lemma gives a conditional exponential-moment bound at~\(\cF_0\).  Since
\(\cH_{b-1}\subseteq\cF_0\), the tower property yields constants
\(\lambda_{\rm gb},C_{\rm gb}>0\) such that
\[
 \EE_b^v\!\left[
   \exp\left\{2\lambda_{\rm gb}
      \sum_n\bigl(x_n^\alpha+C_{\rm grow}Z_n\bigr)\right\}
   \middle|\cH_{b-1}\right]
 \leq C_{\rm gb}.
\]
Here the sum ranges over the layer episodes actually run, with the escape
episode contributing zero as in Definition~\ref{d.ssr-one-slide-charges}.

For each layer~\(j\), concavity gives
\[
 X_{b,j}^\alpha
 =\left(\sum_{n\in\mathcal N_{b,j}}x_n\right)^\alpha
 \leq\sum_{n\in\mathcal N_{b,j}}x_n^\alpha.
\]
Summing over the layers and using the definition of~\(\BdCost_b\) therefore
shows pathwise that
\[
 \GainCost_b+\BdCost_b
 \leq\sum_n\bigl(x_n^\alpha+C_{\rm grow}Z_n\bigr).
\]
Consequently
\begin{equation}
\label{e.ssr-gain-boundary-uniform-moment}
 \EE_b^v\!\left[
   e^{2\lambda_{\rm gb}(\GainCost_b+\BdCost_b)}
   \middle|\cH_{b-1}\right]
 \leq C_{\rm gb}
\end{equation}
uniformly in~\(b\) and~\(v\).  These constants are independent of
\(C_{\rm reg}\) and of the breaker row.

Set
\[
 \mu=\frac12\bigl(\lambda_{\rm prep}\wedge\lambda_{\rm gb}\bigr).
\]
Conditional Cauchy--Schwarz, Lemma~\ref{l.ssr-prospective-preparation-moment},
\eqref{e.ssr-gain-boundary-uniform-moment} and
\eqref{e.ssr-slide-cost-definition} give
\begin{equation}
\label{e.ssr-slide-cost-uniform-moment}
 \EE_b^v\!\left[e^{2\mu\SlideCost_b}\mid\cH_{b-1}\right]
 \leq
 \bigl(C_{\rm prep}C_{\rm gb}\bigr)^{1/2}
 =:C_{\rm sl}.
\end{equation}
Indeed, after Cauchy--Schwarz the two factors have exponents~\(4\mu\),
which are at most~\(2\lambda_{\rm prep}\) and~\(2\lambda_{\rm gb}\),
respectively.  This is the only place where the preparation and gain--boundary
moments are combined.

We now recover the breaker attenuation.  If some monitored layer is
post-passage at step~\(b\), then
\(\sum_{j:b\in\mathcal B_j}a_{b,j}\geq1\).  Hence
\eqref{e.ssr-slide-cost-uniform-moment} already implies the desired bound at
exponent~\(\mu\), after choosing the constant in the right-hand exponential at
least~\(\log C_{\rm sl}\).

Suppose instead that every relevant monitored layer is in the advance phase,
and put
\[
 F_b^{\rm sl}
 =\mathsf{Finite}_b^v\cap\{\SlideCost_b>0\}.
\]
Every positive contribution to~\(\GainCost_b\),~\(\BdCost_b\) or~\(\PrepCost_b\)
comes from a reached monitored layer.  Hence
\begin{equation}
\label{e.ssr-slide-positive-requires-reach}
 F_b^{\rm sl}
 \subseteq
 \bigcup_{j:H_b\leq s_j^-}
   \Reach_{b,j}\cap\mathsf{Finite}_b^v.
\end{equation}
Since~\(\SlideCost_b\) vanishes off~\(F_b^{\rm sl}\), conditional
Cauchy--Schwarz and~\eqref{e.ssr-slide-cost-uniform-moment} give
\begin{align*}
 \EE_b^v\!\left[e^{\mu\SlideCost_b}-1\mid\cH_{b-1}\right]
 &\leq
 \PP_b^v(F_b^{\rm sl}\mid\cH_{b-1})^{1/2}
 \EE_b^v\!\left[(e^{\mu\SlideCost_b}-1)^2\mid\cH_{b-1}\right]^{1/2}\\
 &\leq C_{\rm sl}^{1/2}\,
 \PP_b^v(F_b^{\rm sl}\mid\cH_{b-1})^{1/2}.
\end{align*}
By~\eqref{e.ssr-slide-positive-requires-reach}, the union bound and
\eqref{e.ssr-reach-probability},
\[
 \PP_b^v(F_b^{\rm sl}\mid\cH_{b-1})
 \leq\sum_{j:H_b\leq s_j^-}q_{\rm row}^{M_{b,j}}.
\]
Since~\((\sum_j u_j)^{1/2}\leq\sum_j u_j^{1/2}\) for~$u_j\geq0$, and
\eqref{e.ssr-final-attenuation-base}, we obtain
\[
 \EE_b^v\!\left[e^{\mu\SlideCost_b}\mid\cH_{b-1}\right]
 \leq
 1+C_{\rm sl}^{1/2}
   \sum_{j:H_b\leq s_j^-}\kappa^{M_{b,j}}
 \leq
 \exp\left\{C\sum_{j:b\in\mathcal B_j}a_{b,j}\right\},
\]
where~\(C\) may be chosen to dominate both~\(C_{\rm sl}^{1/2}\) and
\(\log C_{\rm sl}\).  Taking~\(2\lambda=\mu\) proves
\eqref{e.ssr-fixed-shift-moment} in both phases.

For the quantitative assertion, fix~\(q_0<1\) and assume
\(q_{\rm row}\leq q_0\).  By
\eqref{e.ssr-preparation-Creg-dependence},
\(\lambda_{\rm prep}^{-1}+C_{\rm prep}\) is bounded by a fixed power of
\(1+C_{\rm reg}\), uniformly over these rows, while
\(\lambda_{\rm gb}^{-1}+C_{\rm gb}\) is bounded by an absolute constant.  Since
\[
 \lambda=\frac14
 \bigl(\lambda_{\rm prep}\wedge\lambda_{\rm gb}\bigr),
 \qquad
 C_{\rm sl}=(C_{\rm prep}C_{\rm gb})^{1/2},
\]
the chosen~\(\lambda^{-1}\) and the constant~\(C\) in the right-hand
exponential are bounded by a fixed power of~\(1+C_{\rm reg}\).  This gives
\eqref{e.ssr-input-Creg-dependence}.  The bound on~\(\kappa\) is exactly
\eqref{e.ssr-final-attenuation-base}.
\qed

\chapter{Reference-dynamics multi-layer accrual and the~accepted-slide~comparison}
\label{s.model-problems-accrual}

The conditional exponential-moment bound for the one-slide cost of a finite
family of monitored layers, Proposition~\ref{p.ssr-fixed-shift-input}, is now
iterated along the attempted schedule.  We use for this purpose a reference
dynamics, formally specified in the first section.  At each active step it
draws the proposed head from the same conditional head-placement law as the
accepted slide, but runs the corresponding fixed-shift experiment without
conditioning on acceptance; rejection sends the trajectory to a cemetery
state.  Its augmented form carries the sealing and paused-exploration data of
the preceding chapter.  The \emph{marked-kernel identity} verifies that this
auxiliary data leaves the reference cluster evolution unchanged when averaged
out.

Starting from the passage seal controlled in the preceding chapter,
post-passage seal growth proceeds by successive replacement and enlargement;
the task here is to control its cumulative effect along the schedule.  The
one-slide estimate controls finite cluster growth in the monitored layers and
the new seal boundary created when enlargement absorbs that growth.  It does
not by itself control the cumulative outward displacement caused by
replacement.  At each later active slide, the inherited captured region is
enclosed by a fresh seal in the new fixed-shift field, and the outer mesh
boundary may thereby move even when no new cluster material is added.  We
control this displacement by following the outer mesh boundary through
successive replacements.  A single replacement moves the inherited boundary
only through a bounded collar and the subcritical defect components encountered
by seal completion.  Iterating this control gives a replacement front that
advances at linear speed, up to an exponentially tailed excess.  Combined
with the preceding chapter's bounds on finite layer growth and on boundary
newly created by enlargement, this controls the terminal seal and, by
isoperimetry, the captured layer volume.  The resulting joint estimate is the
reference-law multi-layer accrual bound used for the positive-depth
southside conditions.

To transfer these estimates to the accepted dynamics, we record the proposed
head as a mark while keeping the auxiliary sealing information separate.
After this auxiliary information is averaged out, the accepted trajectory is
recovered from the reference trajectory by the appropriate acceptance
correction factors.  The Tether is used in two ways.  A
trajectory-level moment bound obtained in this chapter will control the
acceptance correction when the positive-depth shadow estimate is transferred
to~\(\Qdyn\) in the next.  A one-step comparison transfers uniform
fixed-shift rarity estimates to the accepted-slide law.  The latter is used
for the depth-zero southside condition and later for the southgoing
breaker units.

The chapter first establishes the marked-kernel identity.  It then follows
the outer mesh boundary through successive post-passage seal replacements,
proves the linear-speed estimate for the resulting replacement front and the
consequent terminal seal bound.  The third section iterates the one-slide
moment under the reference law and combines it with the replacement-front
estimate to obtain the multi-layer accrual bound.  The final section records
the two Tether comparisons used to pass from reference-law estimates to the
accepted dynamics.

\section{The marked-kernel identity}

The reference dynamics uses the same attempted-height schedule and initial cluster
law as the dynamic law of Definition~\ref{d.qdyn}, and, as long as the
trajectory has not entered the cemetery state, uses the same rule to decide
whether the slide is active or dormant.  A dormant step is therefore unchanged.  The difference lies in the fixed-shift
representation of an active slide.  If~\(C\) is the current cluster state,
the retained pair~\((\uclust,\dclust)\) is determined by~\(C\).
For a proposed head~\(v\), define the unconditioned fixed-shift cluster
kernel on measurable~\(A\subseteq\mathfrak C_u\) by
\[
 R_C^v(A)=\Rlaw^v_{\uclust,\dclust}
       (\clusteruthree\cap\pi^{-1}(A)).
\]
This is a subprobability measure on~\(\mathfrak C_u\); assign its missing
mass to the cemetery state, and continue to write~\(R_C^v\) for the
resulting probability kernel.  The configuration-level compatibility of this
kernel with the augmented fixed-shift experiment is
Lemma~\ref{l.ssr-dynamic-projection}.  Since acceptance is determined by
the reconstructed cluster, we also view the configuration-defined
event~\(\Accept_v\) as a subset of the cluster-state space~\(\mathfrak C_u\).
Put
\[
 w_C(v)=R_C^v(\Accept_v),\qquad
 \Theta_C=\sum_vw_C(v),\qquad
 \pi_C(v)=\frac{w_C(v)}{\Theta_C}.
\]
Lemma~\ref{l.probform}, in the form~\eqref{e.fixed-shift-head-law}, gives
the fixed-shift representation of the active step in~\(\Qdyn\): draw the
proposed head~\(V_b\) with law~\(\pi_C\), then run the corresponding
fixed-shift experiment conditioned on~\(\Accept_{V_b}\).  Lemma~\ref{l.qfine}
and Proposition~\ref{p.finesum} give \(0<\Theta_C<\infty\).  The
\hfff{reference-dynamics}{\emph{reference dynamics}} changes only the
conditioning in this representation: draw~\(V_b\) from the same law and run
the fixed-shift experiment unconditioned.  Record~\(V_b\) as a proposal mark
and retain the projected fixed-shift output even on rejection (the cemetery
state itself if the reconstruction lies outside~\(\mathfrak C_u\)).
On rejection the reference trajectory enters the absorbing cemetery state.

Section~\ref{ss.ssr-augmented-histories} specifies the augmented
histories~\(Y_b\) and, for each proposed head~\(v\), the augmented
fixed-shift kernel~\(\PP_b^v\).  The
\hfff{augmented-reference-dynamics}{\emph{augmented reference dynamics}} is
the corresponding lift of the reference dynamics: draw~\(V_b\) from
\(\pi_C\), conditionally on the marked cluster history and independently of
the remaining auxiliary information, and then use~\(\PP_b^{V_b}\).  Dormant
steps and transitions to the cemetery state occur just as above.
Probability and expectation under the augmented reference dynamics are
denoted by~\(\PP_{\rm ref}\) and~\(\EE_{\rm ref}\).
For the augmented accepted dynamics, which starts from the same initial
law and uses the augmented fixed-shift kernel conditioned on acceptance
at each active step, we write~\(\PP_{\rm acc}\).

For the marked-kernel calculation, let~\(G_b\) be the projection of~\(Y_b\)
that retains the attempted-height schedule, cluster history and proposal
marks, together with any rejected fixed-shift output, and let~\(z_b\)
comprise the remaining sealing, replacement and paused-search information.
Thus~\(Y_b\) is represented here by~\((G_b,z_b)\).  These state spaces
are standard Borel, so the conditional law
of~\(z_b\) given~\(G_b\) admits a measurable probability-kernel version; let
\(G\mapsto\mu_G\) denote one such version under~\(\PP_{\rm ref}\).  Statements
below involving~\(\mu_G\) are understood for~\(\PP_{\rm ref}\)-almost every
live marked history~\(G\).  Write~\(C=C(G)\) for its current cluster state.
For a proposed head~\(v\) with~\(\pi_C(v)>0\), let
\(\widehat R_G^v(z,dC',dz')\) be the cluster-and-auxiliary output kernel
induced by~\(\PP_b^v\) in the~\((G,z)\) representation.
Put~\(G'=(G,v,C')\).
As we verify shortly in Proposition~\ref{p.seal-build-up-noninterference},
the following fixed-head compatibility identity holds on this full-measure
set:
\begin{equation}
\label{e.marked-kernel-disintegration}
 \int\mu_G(dz)\,\widehat R_G^v(z,dC',dz')
 =R_C^v(dC')\,\mu_{G'}(dz').
\end{equation}
The factor~\(\mu_{G'}\) is the conditional auxiliary law described above;
the content of the identity is that averaging the auxiliary input leaves the
ordinary fixed-shift cluster marginal.

Define the one-step acceptance correction by
\[
 \ell_G(v,C')=
 \begin{cases}
  \one_{\{C'\in\Accept_v\}}/w_C(v),&w_C(v)>0,\\
  0,&w_C(v)=0
 \end{cases},
\]
and abbreviate~\(\ell=\ell_G(V,C')\) for the random proposal~\(V\)
and fixed-shift output~\(C'\).  For each attempted step~\(b\), write
\(\ell_b\) for its correction factor, taking~\(\ell_b=1\) at dormant steps.
A rejecting active step has~\(\ell_b=0\) and sends the reference trajectory
to the absorbing cemetery state.  The attempted-height schedule continues
formally to its end, without further slide operations; its remaining steps
also carry~\(\ell_b=1\).  We call a reference history \emph{live} if the
cemetery state has not been entered.

\begin{lemma}[Marked trajectory comparison]
\label{l.marked-whole-block-change}
Suppose that the compatibility identity~\eqref{e.marked-kernel-disintegration}
holds for
\(\PP_{\rm ref}\)-almost every live marked history~\(G\) and every proposed
head~\(v\) with~\(\pi_C(v)>0\).  Suppose also that, for each such~\(G,v\),
\(\mu_G\)-almost every auxiliary input is admissible, so that
Lemma~\ref{l.ssr-dynamic-projection} applies before averaging.
For such~\(G\), the averaged marked reference kernel is
\[
 K_G^{\rm ref}(v,dC',dz')
 =\pi_C(v)\int\mu_G(dz)\,\widehat R_G^v(z,dC',dz').
\]
The one-step correction satisfies
\(\EE_{\rm ref}[\ell\mid G]=1\), and the averaged marked transition
under~\(\PP_{\rm acc}\) is
\begin{equation}
\label{e.accepted-marked-transition}
 K_G^{\rm acc}(v,dC',dz')
 =\ell_G(v,C')K_G^{\rm ref}(v,dC',dz')
 =\frac1{\Theta_C}
   \one_{\{C'\in\Accept_v\}}R_C^v(dC')\,\mu_{G'}(dz').
\end{equation}
In particular, the cluster marginal of~\(K_G^{\rm acc}\) is the accepted
cluster transition and, conditionally on the marked output history
\(G'=(G,v,C')\), the auxiliary output has law~\(\mu_{G'}\).  The cluster
projection of~\(\PP_{\rm acc}\) is the accepted slide dynamics~\(\Qdyn\).

For successive steps, put~\(L_n=\prod_{b=1}^n\ell_b\).  The process~\((L_n)\) is a mean-one \(\PP_{\rm ref}\)-martingale with respect to the filtration
\((\sigma(G_n))_{n\ge0}\).  Moreover, for every non-negative
\(\sigma(G_n,z_n)\)-measurable random variable~\(F\),
\begin{equation}
\label{e.whole-block-change-of-measure}
 \EE_{\rm acc}[F]=\EE_{\rm ref}[L_nF].
\end{equation}
\end{lemma}

\noindent{\bf Proof.}
By~\eqref{e.marked-kernel-disintegration},
\[
 K_G^{\rm ref}(v,dC',dz')
 =\pi_C(v)R_C^v(dC')\,\mu_{G'}(dz').
\]
Multiplication by~\(\ell_G(v,C')\), together with
\(\pi_C(v)/w_C(v)=1/\Theta_C\) whenever~\(w_C(v)>0\), gives
\eqref{e.accepted-marked-transition}.  Summing its total mass over~\(v\)
gives
\(\Theta_C^{-1}\sum_vR_C^v(\Accept_v)
=\Theta_C^{-1}\sum_vw_C(v)=1\).  Thus
\(\EE_{\rm ref}[\ell\mid G]\) equals one, and the factor~\(\mu_{G'}(dz')\) identifies
the stated conditional law of the auxiliary output.

To identify the accepted augmented law, fix a live incoming augmented
history~\((G,z)\) and a proposed head~\(v\) with~\(\pi_C(v)>0\)
for which the incoming pair is admissible.
Lemma~\ref{l.ssr-dynamic-projection} gives that the acceptance probability
under~\(\widehat R_G^v(z,\cdot)\) is~\(w_C(v)\), independently of~\(z\).
Conditioning the complete augmented fixed-shift step on acceptance
therefore multiplies its unconditioned kernel by
\(\ell_G(v,C')\).  With the same proposal rule~\(\pi_C\), this identifies
the accepted augmented step of the~\(Y_b\)-construction as the
likelihood-ratio transform of the reference step.  The two laws also agree on dormant steps and on passage-seal operations
between slides, and after entry into the cemetery, the correction factors are
neutral.  Iterating these full augmented kernels gives
\eqref{e.whole-block-change-of-measure} for bounded
\(\sigma(G_n,z_n)\)-measurable functions, and monotone convergence gives
the stated non-negative case.

Since~\(L_n\) is measurable with respect to~\(G_n\), the conditional law of
\(z_n\) given~\(G_n\) under~\(\PP_{\rm acc}\) is~\(\mu_{G_n}\) on histories
of positive accepted probability.  The averaged accepted transition is
therefore~\(K_G^{\rm acc}\), as displayed above.  Finally, the conditional
mean-one identity, together with the neutral factors at dormant and
post-cemetery steps, makes~\((L_n)\) a mean-one~\(\PP_{\rm ref}\)-martingale
with respect to~\((\sigma(G_n))_{n\ge0}\). \qed

\begin{corollary}[Moments of the correction factors]
\label{c.marked-likelihood-moment}
Under the hypotheses of Lemma~\ref{l.marked-whole-block-change}, at every
live active step and for~\(0<\theta<1\),
\begin{equation}
\label{e.marked-step-renyi-moment}
 \EE_{\rm ref}[\ell^{1+\theta}\mid G]
 =\frac1{\Theta_C}\sum_{v:\,w_C(v)>0}w_C(v)^{1-\theta}
 =\Theta_C^{-\theta}\sum_v\pi_C(v)^{1-\theta}.
\end{equation}
If deterministic constants~\(M_b\geq1\) bound the right-hand side
at every live active step~\(b\), then
\[
 \EE_{\rm ref}[L_n^{1+\theta}]
 \leq\prod_{b=1}^nM_b.
\]
\end{corollary}

\noindent{\bf Proof.}
Condition on~\(G\) and on a chosen mark~\(V=v\) with
\(\pi_C(v)>0\).  The acceptance indicator has mean~\(w_C(v)\), so
\[
 \EE_{\rm ref}[\ell^{1+\theta}\mid G]
 =\sum_{v:\,w_C(v)>0}\pi_C(v)w_C(v)^{-\theta}.
\]
Substitution of~\(\pi_C(v)=w_C(v)/\Theta_C\) gives both expressions in
\eqref{e.marked-step-renyi-moment}.  At dormant and post-cemetery steps,
\(\ell_b=1\), so the conditional moment is one.  Applying these bounds
successively to~\(L_b^{1+\theta}=L_{b-1}^{1+\theta}\ell_b^{1+\theta}\)
proves the displayed block bound.
\qed

The proposal mark at a rejecting reference step may be projected out after
\eqref{e.whole-block-change-of-measure} is used.  In all applications below,
a dormant attempted step has~\(\SlideCost_b=0\).

It remains to verify~\eqref{e.marked-kernel-disintegration} for the
augmented reference dynamics.  For this we check that the conditional
auxiliary inputs are admissible, allowing the fixed-shift projection lemma
to be applied before disintegration.

\begin{proposition}[Auxiliary sealing preserves the fixed-shift kernel]
\label{p.seal-build-up-noninterference}
For \(\PP_{\rm ref}\)-almost every live marked reference history~$G$, with current
cluster state~$C$, and every proposed head~$v$ with
\(\pi_C(v)>0\), run the augmented unconditioned fixed-shift step using the
paused exploration of Definition~\ref{d.multi-layer-paused}, including
its post-passage seal replacements.
Averaging the input auxiliary coordinate against its conditional law~$\mu_G$
leaves the cluster marginal~$R_C^v$.  Moreover, with~$G'=(G,v,C')$,
disintegration gives
\[
 \int\mu_G(dz)\,\widehat R_G^v(z,dC',dz')
 =R_C^v(dC')\,\mu_{G'}(dz').
\]
\end{proposition}

\noindent{\bf Proof.}
On a live reference history, every preceding active fixed-shift outcome
has been accepted.  For~\(\PP_{\rm ref}\)-almost every such marked
history~$G$, Lemma~\ref{l.ssr-inherited-record-captured} shows pathwise that the
incoming-capture condition~\eqref{e.ssr-incoming-capture} is satisfied;
the remaining admissibility requirements are part of the augmented
construction.  Thus~$\mu_G$ assigns
mass one to admissible auxiliary inputs.

Condition first on an input auxiliary coordinate~$z$ in this support.
Lemma~\ref{l.ssr-dynamic-projection} identifies the cluster projection of the
resulting augmented fixed-shift step with the ordinary kernel~$R_C^v$.  This
marginal does not depend on~$z$.  Averaging against~$\mu_G$ therefore leaves
the same cluster marginal~$R_C^v$.

Since~$\mu_{G'}$ is a version of the conditional auxiliary law given the
marked output~$G'=(G,v,C')$, disintegration gives the displayed
factorisation. \qed

\begin{corollary}
\label{c.marked-kernel-seal-verification}
For \(\PP_{\rm ref}\)-almost every trajectory, the fixed-head compatibility
identity \eqref{e.marked-kernel-disintegration} holds at every active step
and for every head~\(v\) with~\(\pi_C(v)>0\).  Consequently, the marked
trajectory comparison of Lemma~\ref{l.marked-whole-block-change} applies
to the augmented reference and accepted laws.
\end{corollary}

\noindent{\bf Proof.}
The reference proposal rule depends only on the marked reference history and
therefore does not inspect the auxiliary coordinate.
Proposition~\ref{p.seal-build-up-noninterference} supplies the compatibility
identity~\eqref{e.marked-kernel-disintegration} for almost every live
marked history, and its proof establishes the required admissibility of
almost every auxiliary input. \qed

Corollary~\ref{c.marked-kernel-seal-verification} identifies the marked
cluster projection of~\(\PP_{\rm ref}\) with the ordinary reference law.
Fixed-shift estimates that are uniform in the proposed head may therefore be
averaged over the reference proposals.  Write
\(\mathsf{IVF}_j=\{S_0\in\veryfine_{\ell_j}\}\) for the event that the target
\(\ell_j\) is initially very fine; this event is determined by the input record
and attempted-height schedule contained in the pre-passage history.  On
\(\mathsf{IVF}_j\), Proposition~\ref{p.ssr-passage-seal} gives a uniform
fixed-shift exponential tail for the passage-seal excess~\(\PassExc_j\).
For reference histories that have entered the cemetery state before passage,
set~\(\PassExc_j=0\).  Averaging over the reference proposals gives the next estimate.

\begin{lemma}[Passage seal under the reference law]
\label{l.ssr-reference-passage-tail}
Under the augmented reference dynamics,
\[
 \one_{\mathsf{IVF}_j}\,
 \PP_{\rm ref}(\PassExc_j\geq u\mid\cH_j^{\rm pre})
 \leq Ce^{-cu},\qquad u\geq0.
\]
\end{lemma}

\noindent{\bf Proof.}
Since~$\PassExc_j\geq0$, the case~$u=0$ is immediate with~$C\geq1$.
For~$u>0$, it suffices, by the cemetery convention, to consider live
reference histories on~$\mathsf{IVF}_j$.  Condition on~$\cH_j^{\rm pre}$,
which contains the realised
proposal marks and the shift at the last active step in the advance phase.
If such a step exists, Lemma~\ref{l.family-paused} gives the product law on
the unread edges of~$E_j$ at its terminal paused history, including on the
accepted outcome.  As established before Proposition~\ref{p.ssr-passage-seal},
this law persists until passage.  If no such step exists,
Lemma~\ref{l.ssr-input-record-exterior} gives the same input product law.
Proposition~\ref{p.ssr-passage-seal}, uniformly in the realised shift, then
gives the displayed tail.
\qed

\section{Post-passage seals and the replacement front}
\label{ss.rebuilt-seals-to-limelight}

Under the augmented reference dynamics, the seal in a monitored layer is
replaced at each active post-passage step before that layer is explored.  The
new reconstruction starts from the mesh region captured by the preceding
seal.  By Lemma~\ref{l.ssr-inherited-record-captured}, this region already
contains the mesh support of the fixed layer data, and a fresh seal is
completed around it in the new fixed-shift field.

To control the terminal seal, we follow its outer mesh boundary through
successive replacements and enlargements.  This allows us to account
separately for the boundary supplied at passage or created by subsequent
layer growth, and for the displacement of its descendants under replacement.
The passage seal starts a collection of boundary lineages.  Later finite
layer growth may force an enlargement, starting further lineages at sites
of the output boundary that were absent from the replacement boundary.  We call
the sites at which these lineages begin \emph{boundary roots}.  Successive
replacements propagate the existing lineages: each site of the replacement
boundary lies within
a bounded deterministic collar of the preceding boundary, apart from
controlled excursions around subcritical defect components encountered by seal
completion.  Assigning such sites as descendants of sites on the preceding
boundary gives the root--descendant structure that will define the
replacement front.

\subsection{Boundary roots and one-step replacement}

For a seal~$\mathscr S$, write
\[
       A_{\rm mesh}(\mathscr S)=\mathscr S\cap\Mesh
\]
for its mesh boundary.  For a monitored pair~$(\cS_j,\ell_j)$, list the active
post-passage steps strictly before the limelight attempt at~$\ell_j$ as
\[
       b_{j,1}<\cdots<b_{j,N_j}.
\]
By the spacing bound following Definition~\ref{d.two-phases}, the number
of attempted heights in any height interval is bounded by a constant times
one plus its length.  Hence
\begin{equation}
\label{e.postpassage-generation-count}
       N_j\le C(1+\delta_j),
       \qquad \delta_j=\ell_j-s_j^+.
\end{equation}
Put~$\mathscr S_{j,0}^{\rm out}=\Seal_j^0$, the passage seal.  We number
these replacements by generations~$1,\ldots,N_j$, with the passage seal at
generation zero.  The~$r$\textsuperscript{th} generation occurs at step~$b_{j,r}$,
and we write~$\mathscr S_{j,r}^{\rm repl}$ for the fresh replacement seal
of Lemma~\ref{l.ssr-fresh-postpassage-rebuild}; when the step is accepted,
let~$\mathscr S_{j,r}^{\rm out}$ be the output seal after all layer episodes
have been exhausted.  On histories live at the target's limelight time,
write~$\Seal_j^{\rm LL}=\mathscr S_{j,N_j}^{\rm out}$.

We now record the root sets formally.  The passage boundary is the initial
root set, and later finite layer episodes create further roots only where they
create new mesh-boundary sites.  Set
\[
 \mathcal R_{j,0}=A_{\rm mesh}(\mathscr S_{j,0}^{\rm out}),
 \qquad
 \mathcal R_{j,r}
 =A_{\rm mesh}(\mathscr S_{j,r}^{\rm out})
     \setminus A_{\rm mesh}(\mathscr S_{j,r}^{\rm repl}),
 \quad 1\le r\le N_j,
\]
and put
\[
       R_j^{\partial}=\sum_{r=0}^{N_j}|\mathcal R_{j,r}|.
\]
An enlargement may create new mesh-boundary sites while absorbing others
into the interior, with little net change in boundary cardinality.  By
summing the new sites at each stage, $R_j^{\partial}$ records the creation
of boundary lineages.

Recall the depth weight
$\mathsf w_j=(1+\delta_j)^{-\sigma_{\rm wt}}$ from
\eqref{e.ssr-depth-weight}, and put
\[
 W_j=\PassPrep_j+
     \sum_{b:\,s_j^+\le H_b<\ell_j}\BdCost_{b,j}.
\]

\begin{lemma}[Boundary roots and accumulated slide variables]
\label{l.boundary-roots-slide-domination}
If the reference history is live at the limelight time of~$\ell_j$, then
\begin{equation}
\label{e.boundary-root-pathwise}
 R_j^{\partial}\le
 C\bigl(A_{\rm pass}\avf(\delta_j)+\PassExc_j+W_j\bigr).
\end{equation}
For every finite monitored collection, if the reference history is live at
the limelight time of each target in the collection, then
\begin{equation}
\label{e.prospective-domination}
 \sum_j\mathsf w_jW_j\le2\sum_b\SlideCost_b.
\end{equation}
\end{lemma}

\noindent{\bf Proof.}
By the definition of~\(\PassExc_j\),
\[
 |\mathcal R_{j,0}|\le |\Seal_j^0|
 \le A_{\rm pass}\avf(\delta_j)+\PassPrep_j+\PassExc_j.
\]
At a post-passage active step, write the successive finite layer episodes as
$\mathscr S^0\to\mathscr S^1\to\cdots\to\mathscr S^m$, with
$\mathscr S^0=\mathscr S_{j,r}^{\rm repl}$ and
$\mathscr S^m=\mathscr S_{j,r}^{\rm out}$.  Every site of
$A_{\rm mesh}(\mathscr S^m)\setminus A_{\rm mesh}(\mathscr S^0)$ first
appears at one of these extensions.  Lemma~\ref{l.effective-layer-episode}
and the choice of~$C_{\rm grow}$ in Definition~\ref{d.ssr-one-slide-charges}
therefore give
\[
 |\mathcal R_{j,r}|\le \BdCost_{b_{j,r},j}.
\]
Summing proves~\eqref{e.boundary-root-pathwise}, after changing the fixed
constant.

For~\eqref{e.prospective-domination}, every contributing step is finite,
and hence
\[
 \BdCost_b\le\SealCost_b\le\SlideCost_b.
\]
The selected-passage estimate~\eqref{e.ssr-selected-passage-energy}
and~$\mathsf w_j\le1$ give
\begin{align*}
 \sum_j\mathsf w_jW_j
 &=\sum_j\mathsf w_j\PassPrep_j
   +\sum_j\sum_{b:s_j^+\le H_b<\ell_j}
       \mathsf w_j\BdCost_{b,j}\\
 &\le\sum_b\SealCost_b+\sum_b\BdCost_b
 \le2\sum_b\SlideCost_b.\tag*{\qedsymbol}
\end{align*}

The preceding lemma controls how many boundary lineages are started by
passage and enlargement.  We next control how one seal replacement propagates
an existing lineage.  Distances below are mesh $\ell^\infty$ distances.

\begin{lemma}[One replacement of an outer boundary]
\label{l.one-rebuild-front}
There are constants~$\theta_{\rm rep}>0$ and~$M_{\rm rep}<\infty$ with the following property.
Fix a monitored layer~$\cS_j$ at an active step in its post-passage phase.
Let~$\mathcal H^-$ be the paused history immediately before this
replacement, including the proposed head and any earlier
replacements in the same step, and condition on~$\mathcal H^-$.  Let
$\mathscr S^-$ be the output seal of the preceding generation, and let
$\mathscr S^{\rm repl}$ be the fresh replacement seal.  Put
\[
       A^-=A_{\rm mesh}(\mathscr S^-),
       \qquad
       A^{\rm repl}=A_{\rm mesh}(\mathscr S^{\rm repl})
\]
for their mesh boundaries.  If~$A^-=\varnothing$, then~$A^{\rm repl}=\varnothing$.
Otherwise, assign each~$y\in A^{\rm repl}$ to a nearest site
$a(y)\in A^-$ in mesh distance, breaking ties by a fixed order on the
mesh, and put
\[
       T(x)=\{y\in A^{\rm repl}:a(y)=x\},\qquad x\in A^-.
\]
Then, uniformly over the conditioned history, for every~$x\in A^-$,
\begin{equation}
\label{e.one-rebuild-weighted-offspring}
 \EE_{\rm ref}\left[
   \sum_{y\in T(x)}e^{\theta_{\rm rep} d_\infty(x,y)}
   \middle|\mathcal H^-\right]\le M_{\rm rep}.
\end{equation}
\end{lemma}

\noindent{\bf Proof.}
If~$A^-=\varnothing$, then the preceding seal and its captured mesh set
are empty.  The empty-seal convention makes the fresh replacement empty
as well.  Assume henceforth that~$A^-\ne\varnothing$.

By Lemma~\ref{l.ssr-inherited-record-captured}, the mesh support of the fixed
layer data of the new reconstruction is already contained in~$\MCap(\mathscr S^-)$.  Thus
the carrier~\eqref{e.postpassage-rebuild-carrier} is obtained from that
inherited capture by one mesh thickening followed by filling bounded
complementary components.

Write~$D$ for this carrier, and recall that
$A(D)=\partial_{\rm ext}^{\Mesh}D$ is its exterior mesh boundary.  Since
$\mathscr S^-$ is a seal, its defining mesh geometry
\eqref{e.seal-mesh-geometry} gives
$\mathscr S^-=(A^-)^{\mathrm{ov}}$.  The second identity in~\eqref{e.seal-mesh-geometry} identifies~$A^-$ as
the exterior mesh boundary of~$\MeshInt(\mathscr S^-)$.  By
\eqref{e.mesh-transfer-capture}, the inherited capture is
$A^-\sqcup\MeshInt(\mathscr S^-)$.  Its exterior boundary therefore lies
in a fixed neighbourhood of~$A^-$.  One mesh thickening moves the boundary
only a bounded distance, while filling bounded complementary components
cannot create a remote exterior boundary.  Thus, for a fixed~$c_0$,
\begin{equation}
\label{e.rebuild-initial-boundary-local}
       A(D)\subseteq (A^-)^{[c_0]}.
\end{equation}

Run seal completion from~$D$ in the fresh fixed-shift field.  The mesh
support of the retained record lies in~$D$, and no completion test reads an edge in~$\bigcup_{M\in D}E(M)$.  Conditionally on~$\mathcal H^-$, every
edge of~$E_j\setminus\bigcup_{M\in D}E(M)$ is unread in the current fresh
field and hence has its Bernoulli product law.  Corollary~\ref{c.conditional-completion-hull}, applied conditionally on
$\mathcal H^-$, therefore gives the truncated defect mesh set
$\mathcal B_K^D$ and a fixed~$r_0$ such that
\[
 D_*\setminus D
 \subseteq
 A(D)^{[r_0]}\cup
 \mathcal C_K^D\bigl(A(D)^{[r_0]}\bigr).
\]
Lemma~\ref{l.boundary-growth-under-addition} then places every site of the
replacement boundary in a fixed neighbourhood of~$A(D)$ or of a component
of~$\mathcal B_K^D$.  Reduction can only remove boundary components.

For each~$x\in A^-$ let~$\Phi(x)$ consist of a sufficiently large fixed
neighbourhood of~$x$, together with a fixed neighbourhood of every
$*$-component of~$\mathcal B_K^D$ meeting it.  The neighbourhood sizes
are chosen in terms of~$c_0$ and~$r_0$, allowing one additional mesh step
for boundary growth.  The preceding paragraph and
\eqref{e.rebuild-initial-boundary-local} give the covering relation
\begin{equation}
\label{e.one-rebuild-cover}
       A^{\rm repl}\subseteq\bigcup_{z\in A^-}\Phi(z).
\end{equation}

Corollary~\ref{c.conditional-completion-hull} gives a uniform conditional
exponential moment for every rooted component of~$\mathcal B_K^D$.
Component diameter is bounded by component cardinality, and only a bounded
number of possible roots occurs in the fixed neighbourhood of a given site.
Choosing~$\theta_0>0$ below that exponential-moment parameter gives
$M_0<\infty$, uniform in the incoming history, such that for all~$z\in A^-$,
\[
 \EE_{\rm ref}\left[\sum_{y\in\Phi(z)}e^{\theta_0 d_\infty(z,y)}
       \middle|\mathcal H^-\right]\le M_0.
\]
The sets~$\Phi(z)$ are used here only to bound the offspring of the
nearest-site assignment in the lemma.  If~$y\in T(x)\cap\Phi(z)$,
then the choice of~$a(y)=x$ and the triangle inequality give
\[
 d_\infty(x,y)\le d_\infty(z,y),\qquad
 d_\infty(x,z)\le 2d_\infty(z,y).
\]
Set~$\theta_{\rm rep}=\theta_0/2$ and~$\eta=\theta_0/4$.  It follows that
\[
 e^{\theta_{\rm rep}d_\infty(x,y)}
 \le e^{-\eta d_\infty(x,z)}e^{\theta_0d_\infty(z,y)}.
\]
By~\eqref{e.one-rebuild-cover}, we may overcount the terms indexed by
$T(x)$ using the sets~$\Phi(z)$.  Taking conditional expectations yields
\begin{align*}
 \EE_{\rm ref}\left[
   \sum_{y\in T(x)}e^{\theta_{\rm rep}d_\infty(x,y)}
       \middle|\mathcal H^-\right]
 &\le \sum_{z\in A^-}e^{-\eta d_\infty(x,z)}
       \EE_{\rm ref}\left[
       \sum_{y\in\Phi(z)}e^{\theta_0d_\infty(z,y)}
       \middle|\mathcal H^-\right]\\
 &\le M_0\sum_{z\in\Mesh}e^{-\eta d_\infty(x,z)}
 \le M_{\rm rep}.
\end{align*}
The last sum is uniformly finite on the $(d-1)$-dimensional mesh;
we choose~$M_{\rm rep}\ge1$.  This proves~\eqref{e.one-rebuild-weighted-offspring}.
\qed

\subsection{The replacement front and terminal geometry}

We now control the cumulative spread of boundary lineages, first for one
monitored layer and then jointly across several.  At an accepted post-passage
step, sealing starts with the fresh replacement
seal and continues with any enlargement required by layer exploration,
resulting in the output seal.  Every site of the replacement boundary
inherits the ancestry of its assigned nearest site on the preceding boundary.
A site that remains on the output boundary keeps this ancestry; any output
boundary site absent from the replacement boundary belongs to~$\mathcal R_{j,r}$
and begins a new lineage.
By~\eqref{e.postpassage-rebuild-capture}, each replacement places the
preceding captured mesh region inside the new seal's mesh interior, while
subsequent enlargements preserve that interior.  Thus a mesh site that leaves
the outer boundary cannot return to it, and each site can be born as a
boundary root at most once.  Every site of the output boundary of generation~$r$ thus belongs to a
lineage rooted at some generation~$i\le r$.

The nearest-site assignment in Lemma~\ref{l.one-rebuild-front}
depends only on the preceding output seal and the fresh replacement seal,
so it is measurable from the augmented sealing history.  We use the augmented
histories~$\cH_b$ between attempted steps and the paused exploration
history of Section~\ref{ss.multi-layer-paused} within each step.  Before
each replacement, the proposed head and carrier are fixed by
the incoming paused history, while the current layer's fresh edges remain
unread.  The replacements take place in a fixed order on disjoint layer
edge sets, as in Lemma~\ref{l.ssr-fresh-postpassage-rebuild}.  Thus the
conditional estimate of Lemma~\ref{l.one-rebuild-front} applies at each
replacement, including after earlier replacements in the same slide.

For a root~$x\in\mathcal R_{j,i}$, let~$\widehat{\mathcal H}(x)$ be the
augmented history at the end of generation~$i$, when~$x$ is born, and
let $P_n(x)$ be the set of its descendants on the output boundary of
generation~$n\ge i$.  With~$\theta_{\rm rep}$ and~$M_{\rm rep}$
from Lemma~\ref{l.one-rebuild-front}, put
\[
 \Psi_n(x)=M_{\rm rep}^{-(n-i)}
       \sum_{y\in P_n(x)}e^{\theta_{\rm rep} d_\infty(x,y)}.
\]
On rejection, the boundary lineages are discontinued, with their
descendant sets thereafter taken to be empty.  Since~$P_i(x)=\{x\}$, we
have~$\Psi_i(x)=1$.  At a replacement, the triangle inequality and
\eqref{e.one-rebuild-weighted-offspring}, conditioned on the paused
history immediately before the replacement, bound the expected weight of
each existing lineage.  Subsequent same-slide enlargements can only remove
its descendants, and the zero assignment on rejection cannot increase the
weight.  Thus, by the tower property, $(\Psi_n(x))_{n\ge i}$ is a
non-negative supermartingale in the augmented history from its birth
history~$\widehat{\mathcal H}(x)$ onward.
Seal enlargement may introduce new boundary roots; successive seal replacements
propagate the lineages issuing from them.  We call the resulting collection
of descendants the \emph{replacement front}.

Choose~$v_{\rm rep}>0$ so that
$\theta_{\rm rep}v_{\rm rep}>\log M_{\rm rep}$.  A root born at
generation~$i$ has a linear-speed cone of mesh radius~$v_{\rm rep}(n-i)$
at generation~$n\ge i$.  To measure the greatest excess of descendants
beyond these cones, put
\[
 F_j=\max_{\substack{0\le i\le n\le N_j\\
                    x\in\mathcal R_{j,i},\ y\in P_n(x)}}
       \bigl(d_\infty(x,y)-v_{\rm rep}(n-i)\bigr)_+,
\]
with the empty maximum taken to be zero.  If the reference trajectory
enters the cemetery before the target's limelight time, set
$R_j^{\partial}=F_j=0$.

\begin{proposition}[Linear speed of the replacement front]
\label{p.rebuild-linear-speed}
There is~$c_{\rm rep}>0$ such that, for every
root~$x\in\mathcal R_{j,i}$ and every~$t\ge0$,
\begin{equation}
\label{e.single-root-repair-cone}
 \PP_{\rm ref}\left(
   \exists n\ge i,\ \exists y\in P_n(x):
   d_\infty(x,y)>v_{\rm rep}(n-i)+t
   \middle|\widehat{\mathcal H}(x)\right)
 \le e^{-c_{\rm rep}t}.
\end{equation}
Consequently, for every~$q\ge1$,
\begin{equation}
\label{e.one-layer-repair-front-tail}
   \PP_{\rm ref}(F_j>t,\ R_j^{\partial}\le q)\le q e^{-c_{\rm rep}t}.
\end{equation}
\end{proposition}

\noindent{\bf Proof.}
If a descendant at generation~$n$ lies beyond the displayed cone, then
\[
 \Psi_n(x)\ge
 \exp\{(\theta_{\rm rep} v_{\rm rep}-\log M_{\rm rep})(n-i)+\theta_{\rm rep} t\}
 \ge e^{\theta_{\rm rep} t}.
\]
The maximal inequality for non-negative supermartingales gives
\eqref{e.single-root-repair-cone}, with
$c_{\rm rep}=\theta_{\rm rep}$.  For the second assertion, enumerate the
roots in their birth order, with ties broken by a fixed order on mesh sites,
and consider the root indices~$k\le q$.  For a fixed index, introduce
a factor equal to one until a root with that index is created; if this
happens, the factor thereafter follows that root's stopped lineage
supermartingale.  The same maximal estimate bounds by~$e^{-c_{\rm rep}t}$
the probability
that such a root is created and has a descendant with excess greater
than~$t$.  On~$\{R_j^{\partial}\le q\}$ any root with such a descendant
has one of these indices.  Summing over them gives
\eqref{e.one-layer-repair-front-tail}. \qed

For several distinct monitored layers, the fixed-order replacements give
the product form needed for simultaneous control.  At each replacement,
the next layer's fresh field is still unread, even after
conditioning on the preceding replacements in the same slide.

\begin{proposition}[Joint replacement-front bound]
\label{p.joint-rebuild-front}
Fix a finite set~$I$ of monitored pairs whose partition layers are distinct.
For deterministic
$q_j\ge1$ and~$t_j\ge0$,
\begin{equation}
\label{e.joint-rebuild-front}
 \PP_{\rm ref}\left(
    \bigcap_{j\in I}\{F_j>t_j,\ R_j^{\partial}\le q_j\}
    \middle|\mathfrak C_0^{\rm rec},\overline H
 \right)
 \le \prod_{j\in I} q_j e^{-c_{\rm rep}t_j}.
\end{equation}
\end{proposition}

\noindent{\bf Proof.}
Condition on~$\mathfrak C_0^{\rm rec}$ and~$\overline H$ throughout.
Number the roots in each layer as they are created, breaking ties by a
fixed order on mesh sites.  On the event in~\eqref{e.joint-rebuild-front},
each layer~$j\in I$ has a root~$x\in\mathcal R_{j,i}$ with index
at most~$q_j$, born at some generation~$i$, and a descendant
$y\in P_n(x)$ at some generation~$n\ge i$ which exceeds the root's
linear-speed cone by more than~$t_j$:
\[
 d_\infty(x,y)>v_{\rm rep}(n-i)+t_j.
\]
For each fixed choice of positive integers~$k_j\le q_j$, $j\in I$,
we bound the probability that, in every layer~$j$, a root with
index~$k_j$ is created and has a descendant satisfying the displayed
inequality, with the trajectory live at the relevant limelight times.
We then sum over these choices.

For each~$j\in I$, associate a factor, initially equal to one, with the
index~$k_j$.  It remains one while fewer than~$k_j$ roots have been
created in that layer.  When the~$k_j$\textsuperscript{th} root is created,
the factor follows its lineage supermartingale~$\Psi$, starting at one
and stopped at the first output generation at which a descendant satisfies
the displayed excess condition.

The product contains one factor for each selected index throughout the
common attempted-step schedule.  We subdivide each active step into
successive stages corresponding to its layer replacements, in their
prescribed order, updating the relevant factor at each stage.  At a
replacement, condition on the paused history accumulated so far, including
the proposed head and the results of preceding replacements.  The carrier
is determined and the current layer's fresh field remains unread, because
the earlier replacements used disjoint layer edge sets.  By
Lemma~\ref{l.one-rebuild-front}, the conditional expected update of the
factor belonging to this layer is at most its present value; the other
factors remain fixed during that replacement.  Same-slide layer
exploration may remove descendants, but cannot increase the factors of
existing lineages; the creation of a root with a selected index leaves
its factor at one.  On rejection, all factors are set to zero, a downward
update.  Hence the product of the stopped factors is a non-negative
supermartingale in this refined augmented history.

From the limelight time of a given layer onwards, its factor is held fixed
unless the trajectory enters the cemetery.  Stop the product immediately
before the attempted slide at the latest limelight time among the selected
targets.  If a root is created at each selected index and the corresponding
descendants satisfy the displayed excess condition, while the trajectory
is live at all those limelight times, the product is at
least~$\exp\{c_{\rm rep}\sum_{j\in I}t_j\}$.
Earlier cemetery entry gives~$F_j=0$ for a target not yet at limelight and
cannot contribute to the event in~\eqref{e.joint-rebuild-front}.
Conditional Markov's inequality bounds this fixed-index event by
$\exp\{-c_{\rm rep}\sum_{j\in I}t_j\}$.  Summing over at most
$\prod_{j\in I}q_j$ choices of indices gives~\eqref{e.joint-rebuild-front}.
\qed

The replacement-front analysis has now isolated the two quantities governing
the terminal geometry.  The count~$R_j^{\partial}$ records how many boundary
lineages are started, while~$F_j$ controls the excess displacement of their
descendants beyond the linear-speed cone.  The next lemma combines these
controls to bound the terminal seal and, using the seal-volume estimate,
obtains from it the layer-volume bound needed in the next section.

\begin{lemma}[Terminal seal from the replacement front]
\label{l.terminal-seal-captured-volume}
Let~$i(j)$ be the attempted-step index for which the scheduled height equals
the target~$\ell_j$, and put~$S_j^{\rm LL}=S_{i(j)-1}$, the state at its
limelight time immediately before the attempted slide.  If the reference history is live at
this time, then
\begin{equation}
\label{e.terminal-seal-repair-front}
 |\Seal_j^{\rm LL}|
 \le C R_j^{\partial}(1+\delta_j+F_j)^{d-1},
\end{equation}
and
\begin{equation}
\label{e.terminal-cluster-repair-front}
 |\csetuthree(S_j^{\rm LL})\cap\cS_j|_*^\alpha
 \le C R_j^{\partial}(1+\delta_j+F_j)^{d-1}.
\end{equation}
The latter estimate also holds with the target downcluster in place of the
full cluster.
\end{lemma}

\noindent{\bf Proof.}
By the ancestry construction, every terminal mesh-boundary site descends
from a root born at some generation~$i\le N_j$.  The definition of~$F_j$
places that site within mesh distance~$v_{\rm rep}(N_j-i)+F_j$ of its root.  Thus the terminal mesh boundary is
contained in the union of~$R_j^{\partial}$ mesh balls of radius
$v_{\rm rep}N_j+F_j$.  Since the mesh dimension is~$d-1$,
\eqref{e.postpassage-generation-count} and bounded local multiplicity give
\[
 |\Seal_j^{\rm LL}|\le C R_j^{\partial}(1+\delta_j+F_j)^{d-1},
\]
which is~\eqref{e.terminal-seal-repair-front}.  Mergers and filled bounded
components can only delete outer-boundary sites and so require no additional
term.

The passage seal captures the complete layer record by
Proposition~\ref{p.ssr-passage-seal}.  This capture of the complete record
is preserved at each subsequent active step by the inherited-record
induction in Lemma~\ref{l.ssr-inherited-record-captured}, whose proof uses
Lemma~\ref{l.ssr-cluster-completion} to account for newly queried edges.
Dormant steps change nothing, so the terminal seal captures the complete
layer record.  Work in the overlapping-box lattice~$\overline N$ associated with the fixed
partition layer~$\cS_j$.  Using an incident cluster or tested closed-boundary
edge when a box meets the cluster only through a vertex shows that every
overlapping box meeting~$\csetuthree(S_j^{\rm LL})\cap\cS_j$ meets a
captured mesh box.  Hence
\[
 |\csetuthree(S_j^{\rm LL})\cap\cS_j|_*
 \le C|\Capt_{\overline N}(\Seal_j^{\rm LL})|
 \le C(1+|\Seal_j^{\rm LL}|)^{1/\alpha}.
\]
The second inequality is Lemma~\ref{l.seal-cardinality-controls-volume}.
The additive constant requires a separate check when~$R_j^{\partial}=0$.
In that case~\eqref{e.terminal-seal-repair-front} gives
$\Seal_j^{\rm LL}=\varnothing$, so its capture is empty and the first
inequality in the preceding display shows that the terminal layer
cluster is also empty.  For~$R_j^{\partial}\ge1$, the constant may be
absorbed into the right-hand side of~\eqref{e.terminal-seal-repair-front},
giving~\eqref{e.terminal-cluster-repair-front}.
The downcluster assertion follows by inclusion. \qed

\section{The reference-law multi-layer accrual bound}
\label{s.reference-law-multilayer-accrual}

The cluster material in a monitored layer at limelight can be large either
because many boundary-root lineages are born or because later replacements
carry existing lineages unusually far.  The replacement-front estimates above
control the second mechanism.  The passage seal supplies the initial roots,
while subsequent enlargements create further roots.
Proposition~\ref{p.ssr-fixed-shift-input} controls prospective passage
preparation, finite layer growth and the associated boundary creation, one
proposed slide at a time.  We now accumulate these costs along the reference
schedule and combine them with the passage-seal excess bound of
Lemma~\ref{l.ssr-reference-passage-tail} to control the number of roots.
Together with the replacement-front estimate, this will control layer accrual.

Throughout this section we condition on the complete input-cluster
record~\(\mathfrak C_0^{\rm rec}\) and attempted-height
schedule~\(\overline H\), with probability and expectation still taken
under the chapter-wide reference law.  The constants are uniform in the
realised record and schedule.

Depth regularity controls the monitored family from the viewpoint of one
proposed slide, preventing too many layers from having both small depth and
little breaker attenuation.  Viewed instead from a fixed layer, the cumulative
quantity is its reference exposure.  Recall from
Section~\ref{ss.prospective-preparation-breaker} that \(\mathcal B_j\) is the
set of attempted steps whose one-slide construction may contribute to
layer~\(j\).  The weight~\(a_{b,j}\) equals one at post-passage steps and,
in the advance phase, is attenuated by the number of intervening breaker
windows.  Its sum
\begin{equation}\label{e.reference-exposure}
 A_j^{\rm ref}=\sum_{b\in\mathcal B_j}a_{b,j}
\end{equation}
is the total reference exposure of layer~\(j\) to the surrounding slides.

To bound the accumulated one-slide costs, fix a finite stretch of the
attempted schedule and index its~\(n\) steps by
\(b=1,\ldots,n\).  Dormant steps are included, with~\(\SlideCost_b=0\) and
\(\ell_b=1\).  For a subcollection~$I$ of the monitored pairs, let
$\SlideCost_{b,I}$ denote the one-slide cost of
Definition~\ref{d.ssr-one-slide-charges} with every layer sum restricted to
$j\in I$.  Thus finite gain, boundary growth and prospective preparation are
all restricted to the same subcollection.

\begin{lemma}[Reference-law moment for accumulated slide variables]
\label{l.reference-block-cost}
For a finite \(C_{\rm reg}\)-depth-regular monitored collection, let
\(C,\lambda>0\) be the moment constants from
Proposition~\ref{p.ssr-fixed-shift-input}.  Then, for every
subcollection~$I$,
\[
 \EE_{\rm ref}\left[
   \exp\left\{2\lambda\sum_{b=1}^n\SlideCost_{b,I}\right\}
   \middle|\mathfrak C_0^{\rm rec},\overline H
 \right]
 \leq
 \exp\left\{C\sum_{j\in I}A_j^{\rm ref}\right\}.
\]
\end{lemma}

\noindent{\bf Proof.}
Put
\(A_{b,I}^{\rm step}=\sum_{j\in I:\,b\in\mathcal B_j}a_{b,j}\).
The proof of Proposition~\ref{p.ssr-fixed-shift-input} applies to the
restricted cost with this restricted exposure.  Indeed, if some~$j\in I$ is
post-passage at step~$b$, then~$A_{b,I}^{\rm step}\ge1$; otherwise a
positive restricted cost forces the reach of an advance-phase layer indexed
by~$I$, so the breaker-attenuation union bound is restricted to the same
indices.  Hence, conditionally on the full augmented pre-step history, every
admissible proposed head~\(v\) satisfies
\[
 \EE_b^v[e^{2\lambda\SlideCost_{b,I}}\mid\cH_{b-1}]
 \le e^{C A_{b,I}^{\rm step}}.
\]
The restricted cost is zero at dormant steps and after entry into the
cemetery.  Averaging over the reference proposal law and iterating over the
attempted schedule therefore gives
\[
 \EE_{\rm ref}\left[
   \exp\left\{2\lambda\sum_{b=1}^n\SlideCost_{b,I}\right\}
   \middle|\mathfrak C_0^{\rm rec},\overline H
 \right]
 \leq
 \exp\left\{C\sum_{b=1}^nA_{b,I}^{\rm step}\right\}
 \leq\exp\left\{C\sum_{j\in I}A_j^{\rm ref}\right\}.
\]
\qed

The accumulated-slide moment controls the~\(W_j\) term in the pathwise
boundary-root estimate~\eqref{e.boundary-root-pathwise}, while
Lemma~\ref{l.ssr-reference-passage-tail} controls the passage-seal
excess~\(\PassExc_j\).  Combining these two inputs gives a joint tail for the
numbers of boundary roots born in the monitored layers.  For a
subcollection~$I$, write
$\mathsf{IVF}_I=\bigcap_{j\in I}\mathsf{IVF}_j$ for the event that every
target indexed by~$I$ is initially very fine.  This event is determined by
the input record and attempted-height schedule, so its indicator is fixed
under the conditioning used throughout the section.  Let~$\mathcal L_I$
be the event that the reference history is live
at the limelight time of every~$\ell_j$, $j\in I$.

\begin{lemma}[Reference-law tail for boundary-root counts]
\label{l.reference-root-count}
For a finite $C_{\rm reg}$-depth-regular collection with distinct
partition layers, there are~$C,c>0$ such that, for every subcollection
$I$ and every deterministic~$q_j\ge C\avf(\delta_j)$,
\[
 \one_{\mathsf{IVF}_I}\,
 \PP_{\rm ref}\left(
   \mathcal L_I\cap
   \bigcap_{j\in I}\{R_j^{\partial}>q_j\}
   \middle|\mathfrak C_0^{\rm rec},\overline H
 \right)
 \le
 \exp\left\{
    C|I|+C\sum_{j\in I}A_j^{\rm ref}
    -c\sum_{j\in I}\mathsf w_j
          \bigl(q_j-C\avf(\delta_j)\bigr)
 \right\}.
\]
\end{lemma}

\noindent{\bf Proof.}
It suffices to work on~$\mathsf{IVF}_I$.  We use the above-defined
restricted one-slide costs~$\SlideCost_{b,I}$, evaluated on the fixed
monitored construction.  The termwise proof of~\eqref{e.ssr-selected-passage-energy}
remains valid when every layer sum is restricted to~$I$.
Applying Lemma~\ref{l.ssr-reference-passage-tail} successively in the
order of the deterministic passage positions gives an exponential moment
bound for $\sum_{j\in I}\mathsf w_j\PassExc_j$.
Restricting the termwise proof of~\eqref{e.prospective-domination} to~$I$
gives
$\sum_{j\in I}\mathsf w_jW_j\le2\sum_b\SlideCost_{b,I}$.
Lemma~\ref{l.reference-block-cost} therefore gives an exponential moment
for~$\sum_{j\in I}\mathsf w_jW_j$.  Conditional Cauchy--Schwarz then
yields constants~$\eta,C>0$ such that
\[
 \EE_{\rm ref}\left[
   \one_{\mathcal L_I}
   \exp\left\{\eta\sum_{j\in I}\mathsf w_j
          (\PassExc_j+W_j)\right\}
   \middle|\mathfrak C_0^{\rm rec},\overline H\right]
 \le
 \exp\left\{C|I|+C\sum_{j\in I}A_j^{\rm ref}\right\}.
\]
When~$R_j^{\partial}$ strictly exceeds~$q_j$, the pathwise
estimate~\eqref{e.boundary-root-pathwise} implies, after increasing the
deterministic multiple of~$\avf(\delta_j)$, that
$\PassExc_j+W_j\ge c(q_j-C\avf(\delta_j))$.  Applying the exponential
Markov inequality to~$\sum_{j\in I}\mathsf w_j(\PassExc_j+W_j)$
gives the stated bound. \qed

The root-count bound just proved complements the replacement-front estimate
of Proposition~\ref{p.joint-rebuild-front}.  For a prescribed volume
threshold, the terminal-seal estimate limits how much boundary-root creation
and front displacement can be tolerated.  The next result combines these
controls simultaneously across the monitored layers.

\begin{corollary}[Reference-law multi-layer accrual]
\label{c.reference-multilayer-accrual}
Let~$(\cS_j,\ell_j)$, $j\in I$, be a finite
$C_{\rm reg}$-depth-regular collection with distinct partition layers; target
heights may repeat.  There are positive constants~$C_{\rm acc}$, $C$ and~$c$
with the following property.  Fix deterministic volume thresholds~$L_j>0$.
For each~$j$, choose a root-count threshold~$q_j$ and a front-excess
threshold~$t_j$ such that
\[
 q_j\ge C_{\rm acc}\avf(\delta_j),\qquad t_j\ge0,
\]
and
\begin{equation}
\label{e.reference-accrual-threshold-condition}
       C_{\rm acc}q_j(1+\delta_j+t_j)^{d-1}
       \le L_j^\alpha.
\end{equation}
Recording the boundary-count contribution above~$q_j$ and the front-excess
contribution above~$t_j$ in the form
\begin{align}
 \varepsilon_j^{\rm root}&=\min\left\{1,
   \exp\left[C+CA_j^{\rm ref}
      -c\mathsf w_j\bigl(q_j-C_{\rm acc}\avf(\delta_j)\bigr)\right]\right\},
 \label{e.reference-source-envelope}\\
 \varepsilon_j^{\rm front}&=\min\{1,\ q_j e^{-c_{\rm rep}t_j}\},
 \label{e.reference-front-envelope}
\end{align}
we have
\[
 \one_{\mathsf{IVF}_I}\,
 \PP_{\rm ref}\left(
   \mathcal L_I\cap
   \bigcap_{j\in I}
   \left\{|\csetuthree(S_j^{\rm LL})\cap\cS_j|_*\ge L_j\right\}
   \middle|\mathfrak C_0^{\rm rec},\overline H
 \right)
 \le \prod_{j\in I}\bigl((\varepsilon_j^{\rm root})^{1/2}+(\varepsilon_j^{\rm front})^{1/2}\bigr).
\]
The same estimate holds with the target downcluster in place of the full
cluster.  The constants are uniform over the monitored collections, input-cluster
records and attempted-height schedules under consideration.  For any fixed
$q_0\in(0,1)$, if the breaker-row finite-passage bound $q_{\rm row}$
from~\eqref{e.q-row-definition} satisfies $q_{\rm row}\le q_0$,
the constants may moreover be chosen uniformly in $q_{\rm row}$,
independently of the positive-depth allowance coefficient~$\Cfine$,
and with at most polynomial dependence on~$C_{\rm reg}$.
The attempted-height schedule enters this dependence only through
$C_{\rm reg}$.
\end{corollary}

\noindent{\bf Proof.}
It suffices to work on~$\mathsf{IVF}_I$.  Put
$V_j=|\csetuthree(S_j^{\rm LL})\cap\cS_j|_*$.  On~$\mathcal L_I$,
Lemma~\ref{l.terminal-seal-captured-volume} gives, whenever
$R_j^{\partial}\le q_j$ and~$F_j\le t_j$,
\[
 V_j^\alpha\le Cq_j(1+\delta_j+t_j)^{d-1}.
\]
After increasing~$C_{\rm acc}$ in
\eqref{e.reference-accrual-threshold-condition}, the event~$\mathcal L_I\cap\{V_j\ge L_j\}$
is therefore contained in
\[
 E_j^{\rm root}\cup E_j^{\rm front},
 \qquad
 E_j^{\rm root}=\{R_j^{\partial}>q_j\},
 \quad
 E_j^{\rm front}=\{F_j>t_j,\ R_j^{\partial}\le q_j\}.
\]

For every~$J\subseteq I$, apply Lemma~\ref{l.reference-root-count}
only to those indices in~$J$ for which~$\varepsilon_j^{\rm root}<1$; the remaining factors equal
one.  This gives
\[
 \PP_{\rm ref}\left(
   \mathcal L_J\cap\bigcap_{j\in J}E_j^{\rm root}
   \middle|\mathfrak C_0^{\rm rec},\overline H\right)
 \le\prod_{j\in J}\varepsilon_j^{\rm root}.
\]
Likewise, apply Proposition~\ref{p.joint-rebuild-front} to the indices in
$I\setminus J$ for which~$\varepsilon_j^{\rm front}<1$.  Then
\[
 \PP_{\rm ref}\left(
   \bigcap_{j\in I\setminus J}E_j^{\rm front}
   \middle|\mathfrak C_0^{\rm rec},\overline H\right)
 \le\prod_{j\in I\setminus J}\varepsilon_j^{\rm front}.
\]
Cauchy--Schwarz and summation over the choices of~$J$ yield
\begin{align*}
 \PP_{\rm ref}\left(
     \mathcal L_I\cap
     \bigcap_{j\in I}(E_j^{\rm root}\cup E_j^{\rm front})
     \middle|\mathfrak C_0^{\rm rec},\overline H\right)
 &\le
 \sum_{J\subseteq I}
   \left(\prod_{j\in J}\varepsilon_j^{\rm root}\right)^{1/2}
   \left(\prod_{j\in I\setminus J}\varepsilon_j^{\rm front}\right)^{1/2}\\
 &=\prod_{j\in I}((\varepsilon_j^{\rm root})^{1/2}+(\varepsilon_j^{\rm front})^{1/2}),
\end{align*}
which proves the asserted bound.  The downcluster assertion is
by inclusion.  The one-slide moment constants are independent of~$\Cfine$ by
Proposition~\ref{p.ssr-fixed-shift-input}.  The passage-seal and
replacement-front constants, and the constants in the local rooted
completion bound, are likewise independent of this allowance coefficient.
These bounds are uniform over admissible schedules and histories; their
quantitative dependence on~$C_{\rm reg}$ comes from
\eqref{e.ssr-input-Creg-dependence}.
\qed

\section{Accepted shifts and the Tether}

The multi-layer estimate has been proved under the reference dynamics.  To
return to the accepted slide dynamics, it remains to restore the conditioning
on acceptance.  A fixed-shift estimate could in principle be badly distorted
at a head whose acceptance chance is very small.

Together with the single-slide SimpleJoin lower bound from Chapter~\ref{s.single-slide-positive-density}, the Tether makes this comparison effective.  It implies that the total accepted mass
carried by head positions whose individual accepted mass is below a threshold
$t$ decays as a power of~$t$; above that threshold there can be at most
$t^{-1}$ heads.  This lets us sum the uniformly rare fixed-shift events over the accepted
heads, obtaining an accepted-event probability bound that is a positive
power of the fixed-shift bound.  The same small-atom estimate also gives
a positive moment for the reciprocal mass of the selected head, which controls
the acceptance correction accumulated along a trajectory.
Lemma~\ref{l.accepted-rare-event-transfer} supplies the local
fixed-shift-to-accepted comparison used for depth zero in
Section~\ref{ss.zero-depth-delivery}.  The marked trajectory comparison of
Lemma~\ref{l.marked-whole-block-change}, together with the moment bound
of Corollary~\ref{c.tether-likelihood-moment}, then permits the completed
reference-law shadow estimate to be transferred to the accepted dynamics
through a single trajectory-level change of measure in
Section~\ref{ss.ivf-llf-transition}.

At an active reference step, let~\(C\) be the pre-slide cluster state and
recall the proposal notation
\[
 w_C(v)=R_C^v(\Accept_v),\qquad
 \Theta_C=\sum_v w_C(v),\qquad
 \pi_C(v)=\frac{w_C(v)}{\Theta_C}.
\]
Thus~\(\pi_C\) is the accepted head-placement law and~\(\Theta_C\) is
the total acceptance mass.  Lemma~\ref{l.ssr-dynamic-projection} and
Corollary~\ref{c.marked-kernel-seal-verification} show that, after the
auxiliary coordinate is averaged out, each fixed head has the ordinary
fixed-shift cluster law.  Corollary~\ref{c.tether-head-summability} therefore
applies to~\(\pi_C\): uniformly in the current fine state,
\[
   \sum_{v:\,\pi_C(v)<t}\pi_C(v)
   \leq C_{\rm T}t^{\gamma_{\rm T}},
   \qquad 0<t\leq1.
\]
The quantitative dependence of~\(C_{\rm T}\) is recorded in
\eqref{e.tether-Cfine-dependence}.  This is the small-atom estimate used in
both comparisons below.

First, this estimate gives a one-step comparison by which uniform rarity
under the fixed-shift laws transfers to positive-power rarity after
conditioning on acceptance.

\begin{lemma}[Accepted transfer of a uniformly rare fixed-shift event]
\label{l.accepted-rare-event-transfer}
At an active reference step with current cluster state~\(C\), let
\(E^v\) be a marked event in the augmented fixed-shift experiment for head
\(v\).  Suppose that, for~\(\PP_{\rm ref}\)-almost every live marked
pre-step history~\(G\),
\[
 \int\mu_G(dz)\,\widehat R_G^v(z,E^v)\le q,
 \qquad 0<q\le1,
\]
for every proposed head~\(v\) with~\(\pi_C(v)>0\).  Then there are
\(\delta_{\rm A}\in(0,1)\) and \(1\le C_{\rm A}<\infty\) such that,
for~\(\PP_{\rm acc}\)-almost every pre-step history~\(G\),
\begin{equation}
\label{e.accepted-rare-transfer}
 \PP_{\rm acc}(E^{V}\mid G)
 \le C_{\rm A}q^{\delta_{\rm A}}.
\end{equation}
For the allowance in~\eqref{e.allowance-design}, \(\delta_{\rm A}\) may be fixed independently of
\(\Cfine\), while
\begin{equation}
\label{e.accepted-rare-Cfine}
 C_{\rm A}\le C_{\rm tr}(2+\Cfine)^{D_{\rm A}}
\end{equation}
for some constants \(C_{\rm tr},D_{\rm A}<\infty\) independent of~\(\Cfine\).
\end{lemma}

\noindent{\bf Proof.}
Average the augmented fixed-shift kernel over the input auxiliary coordinate
with law~\(\mu_G\), and denote the resulting probability by
\(\overline R_G^v\).  Since
\(\overline R_G^v(E^v)\le q\) by hypothesis,
\(\overline R_G^v(E^v\cap\Accept_v)\le q\).  Corollary~\ref{c.marked-kernel-seal-verification} supplies the compatibility
identity~\eqref{e.marked-kernel-disintegration}, so that
\(\overline R_G^v(E^v\cap\Accept_v)\)
is also bounded above by \(\overline R_G^v(\Accept_v)=w_C(v)\).  Hence the accepted marked transition of
Lemma~\ref{l.marked-whole-block-change} gives
\[
 \PP_{\rm acc}(E^V\mid G)
 \le \sum_v\min\{q/\Theta_C,\pi_C(v)\}.
\]
Put \(a=q/\Theta_C\).  Lemma~\ref{l.qfine} gives
\(\Theta_C\ge c>0\).  For \(a\le1\), split the sum at
\(t=a^{1/2}\).  Corollary~\ref{c.tether-head-summability} bounds the
mass of \(\{\pi_C<t\}\) by \(C_{\rm T}a^{\gamma_{\rm T}/2}\), while there
are at most \(t^{-1}\) heads with \(\pi_C\ge t\), whose contribution is at
most \(a^{1/2}\).  Enlarging the constant covers \(a>1\).  Thus one may take
\(\delta_{\rm A}=\tfrac12\min\{1,\gamma_{\rm T}\}\).
Equation~\eqref{e.accepted-rare-Cfine} follows from
\eqref{e.tether-Cfine-dependence}.
\qed

In its second application, the small-atom estimate controls the accumulation
of acceptance corrections along a trajectory.  If~$V$ has law~$\pi_C$,
then~$\pi_C(V)^{-1}$ measures the reciprocal accepted mass of the
selected head.  Corollary~\ref{c.tether-head-summability} gives this
quantity a positive moment.  Combined with the uniform lower bound on
$\Theta_C$, Corollary~\ref{c.marked-likelihood-moment} turns this into a
moment bound for the one-step acceptance correction; these bounds can then
be iterated.

Fix once and for all a likelihood-moment exponent~$\theta$ satisfying
\begin{equation}\label{e.theta-choice}
 0<\theta<\min\{1,\gamma_{\rm T}/2\}.
\end{equation}

\begin{corollary}[Moments of the accepted-slide correction]
\label{c.tether-likelihood-moment}
There is~\(C_{\rm L}<\infty\), independent of the marked reference history,
such that at every active step
\[
 \EE_{\rm ref}[\ell_b^{1+\theta}\mid G_{b-1}]
 \leq e^{C_{\rm L}\theta}.
\]
If a finite stretch of the schedule contains~\(n\) attempted steps, then
\begin{equation}
\label{e.block-likelihood-renyi}
 \EE_{\rm ref}[L_n^{1+\theta}]
 \leq e^{C_{\rm L}\theta n}.
\end{equation}
Moreover, for the allowance in~\eqref{e.allowance-design}, with \(A_0,A_1\) fixed,
\begin{equation}
\label{e.likelihood-Cfine-log}
 C_{\rm L}\, \leq \, C\bigl(1+\log(2+\Cfine)\bigr).
\end{equation}
The constant is uniform in the breaker-row parameters \(b_0,C_0,C_1,i_{\rm br}\).
\end{corollary}

\noindent{\bf Proof.}
Let~\(V\) have law~\(\pi_C\), write~\(\EE_{\pi_C}\) for expectation
on this head-placement space, and put~\(\rho=\pi_C(V)^{-1}\).
Set~\(\theta_*=\min\{1,\gamma_{\rm T}/2\}\), so that~\(0<\theta<\theta_*\).
Corollary~\ref{c.tether-head-summability} and integration by parts give
\[
 \EE_{\pi_C}[\rho^{\theta_*}]
 \leq 1+\frac{C_{\rm T}\theta_*}{\gamma_{\rm T}-\theta_*}.
\]
Jensen's inequality and~\eqref{e.tether-Cfine-dependence} then yield
\[
 \log\EE_{\pi_C}[\rho^\theta]
 \leq \frac{\theta}{\theta_*}\log\EE_{\pi_C}[\rho^{\theta_*}]
 \leq C\theta\bigl(1+\log(2+\Cfine)\bigr),
\]
with~\(C\) independent of~\(\theta\).  The lower bound on~\(\Theta_C\) from Lemma~\ref{l.qfine}
depends on the fine allowance only through~\(\af(0)=A_1\), which was
fixed before~\(\Cfine\).  Combining this fixed lower bound
with the preceding estimate gives
\[
 \log\!\left(
 \Theta_C^{-\theta}\EE_{\pi_C}[\rho^\theta]\right)
 \leq C\theta\bigl(1+\log(2+\Cfine)\bigr).
\]
By Corollary~\ref{c.marked-likelihood-moment}, specifically
\eqref{e.marked-step-renyi-moment}, this proves the one-step displayed bound and
\eqref{e.likelihood-Cfine-log}.  As~\(\ell_b=1\) at dormant steps, iteration in the marked reference
filtration proves~\eqref{e.block-likelihood-renyi}.
\qed

\chapter{Shadow bounds, the fine transition, and renewal production}
\label{s.shadow-bounds-delivery-fine}

This chapter proves the IVF-to-LLF transition.  We show that, with overwhelming
probability, a positive fraction of the scheduled heights that are very fine
in the input are fine when they reach limelight.  Thus the
target heights themselves pass or fail the transition.  A target succeeds when
all of its fine tests hold at that moment.  The northside allowances and local
breaker profile are already carried from the input.  The substantial work of
the chapter is on the southside: at limelight, the target-height downcluster
must remain within all of the volume allowances imposed in the bands below the
target.

Our principal technical input is
Corollary~\ref{c.reference-multilayer-accrual}, the
\emph{reference-law multi-layer accrual} estimate.  It controls excessive
cluster volume simultaneously in a finite depth-regular family of monitored
layer--target pairs.  We obtain the required depth regularity by thinning the
positive-density supply of IVF heights with a one-sided Hardy--Littlewood
maximal argument.  This selects a positive-density set
$\mathcal W=\{h_1<\cdots<h_N\}$ of \emph{well-spaced} heights, none of
which is unduly isolated from other IVF heights to its south: every interval
extending southward from such a target contains linearly many such heights,
up to a fixed additive error.  All IVF heights, whether selected or not, remain available as
breakers.  Their breaker windows make every finite family of distinct partition
layers paired with selected targets depth regular.  For a layer--target pair at
depth $z\geq1$, the same breaker supply also gives the reference exposure bound
$A^{\rm ref}\le C(1+z)$.  Thus the target selection gives the depth regularity
required by Corollary~\ref{c.reference-multilayer-accrual}, together with a
reference exposure that grows only linearly in depth.

Fix now one selected target and one partition layer to its south, at depth
$z\geq1$, and let $V$ denote the volume of the target-height downcluster there
at limelight, taking $V=0$ if the reference trajectory was discontinued before
that time.  The corresponding fine condition allows a volume of order
$\af(z)$ in the bands assigned to this layer.  The preceding chapter traps the
layer cluster inside the capture of the post-passage seal and controls how this
capture grows.  Its terminal cluster estimate,
\eqref{e.terminal-cluster-repair-front} in
Lemma~\ref{l.terminal-seal-captured-volume}, is
$V^\alpha\le C R^{\partial}(1+z+F)^{d-1}$.  The two quantities here record
the two ways the seal capture may become large: $R^{\partial}$ counts the
boundary lineages present at passage together with those created by later seal
enlargement, while $F$ measures exceptional later spread of those lineages
under replacement.

Let $\Lambda$ be the volume allowed by the fine test.  We choose a \emph{root
threshold} $q$, limiting the number of boundary lineages, and a \emph{front
threshold} $\tau$, limiting their excess spread, subject to the volume budget
$Cq(1+z+\tau)^{d-1}\le \Lambda^\alpha$.  Corollary~\ref{c.reference-multilayer-accrual}
records the corresponding root and front errors in the schematic form
\[
 \varepsilon^{\rm root}
 \lesssim
 \exp\!\left\{CA^{\rm ref}
   -c(1+z)^{-\sigma_{\rm wt}}
      \bigl(q-C\avf(z)\bigr)\right\},
 \qquad
 \varepsilon^{\rm front}
 \lesssim q e^{-c\tau}.
\]
For one layer--target pair, the corollary bounds the large-accrual
probability, up to fixed constants, by
\mbox{$\sqrt{\varepsilon^{\rm root}}+\sqrt{\varepsilon^{\rm front}}$};
the multi-layer statement gives the corresponding product bound over a
depth-regular family.

For the fine test, $\Lambda$ is of order $\af(z)$ and hence grows as a sufficiently
high power of the depth.  Taking~$\tau$ proportional to
$(1+\log(2+\Cfine))(1+z)$ makes the front error exponentially small.  The
remaining slack allows~$q$ to be taken well above the inherited~$\avf(z)$
scale, and even after the polynomial depth factor~
$(1+z)^{-\sigma_{\rm wt}}$ is applied, the resulting excess dominates the
linear reference exposure $A^{\rm ref}\le C(1+z)$.  Thus the root error is
exponentially small as well.  Consequently, for suitable
constants $c,C,c_{\rm sh}>0$,
\begin{equation}
\label{e.one-pair-depth-decay}
   \PP_{\rm ref}\!\left(V\ge c_{\rm sh}\af(z)\right)
      \le C e^{-cz}, \qquad z\ge1.
\end{equation}
Thus a fixed layer becomes exponentially unlikely to spoil a selected target
as the target recedes northward.

How a layer may spoil targets to its north is encoded by a tower-and-shadow
description: the layer is a \emph{tower}, and it casts \emph{shade} on a
selected target when, at the target's limelight time, the target-height
downcluster accumulated in that layer violates one of the conditions
$\mathbf A^-_{\rm f}(k)$, $k\ge1$, assigned to that layer.  The shadow interval of a tower extends from
the first selected target at or north of it to the last target on which it casts
shade, and a target lying in no such interval is \emph{sunny}.  For a
collection of distinct towers, a lower bound on each shadow length is witnessed
by a sufficiently distant selected target on which that tower casts shade.
Applying Corollary~\ref{c.reference-multilayer-accrual} to the resulting
distinct layer--target pairs, and summing over the possible witnesses, gives
product upper bounds for the joint shadow-length tails under the reference
dynamics.  The root and front thresholds can moreover be chosen with enough
slack that these bounds yield exponential-moment control strong enough to
absorb the likelihood-ratio cost of reweighting from reference trajectories
to accepted slide trajectories.  This provides an exponential-moment comparison
for the total amount of shade; no independence of the actual shadows is asserted.

We therefore pursue the shadow argument under the reference dynamics until
the probability that too many selected targets lie in shadow has been bounded.  The
marked trajectory comparison of Chapter~\ref{s.model-problems-accrual} then
reweights this event once to the accepted slide dynamics.  Except with
exponentially small probability, at least $3N/4$ selected targets are sunny.

Depth zero is handled somewhat differently.  For the condition
$\mathbf A^-_{\rm f}(0)$, passage through the layer containing the target band
coincides with the target's arrival at limelight, so there is no post-passage
phase.  The retained host in the band is bounded deterministically.  A
depth-zero failure must therefore be caused by the fresh guest supplied by the
last active slide below the target.  For each such slide, we record how far
along the selected target list its midpiece may reach.  The breaker window attached to that slide and the recovered IVF breakers
between the slide and a later target give this reach a geometric tail under
every fixed shift.  The local accepted-shift comparison of
Chapter~\ref{s.model-problems-accrual} then reweights this fixed-shift
estimate over the selected head.  Summing the resulting reach variables
shows that, except with exponentially small probability, at most $N/4$
selected targets fail their depth-zero condition.

The parameter scales can be balanced so that reweighting the shadow estimate
to the accepted dynamics still leaves an exponentially small probability of
more than $N/4$ selected targets lying in shadow, while midpieces are
sufficiently unlikely to reach or pass the breaker rows used in the depth-zero
and later renewal arguments.  The two conclusions therefore give at least $N/2$ selected targets that
are sunny and satisfy their depth-zero condition.  At each such target every
southside fine condition holds; together with the northside data carried from
the input, this makes the target fine at limelight.  This proves the IVF-to-LLF transition,
Theorem~\ref{t.vf-to-fine-delivery}.

It remains to convert this supply of LLF heights into renewal production.  At
a scheduled height that is fine at limelight, the acting slide is active, and
the single-slide SimpleJoin estimate gives a uniformly positive conditional
probability that it produces a renewal point.  This converts the linear LLF
supply into positive-density renewal production, except with exponentially
small probability.  Whether one of these points survives the later slides is
the subject of Chapter~\ref{s.renewal-survival}.

The first section selects the well-spaced targets and proves the linear
exposure estimate.  We then treat the $k=0$ condition directly under the
accepted dynamics and organise the positive-depth conditions
$\mathbf A^-_{\rm f}(k)$, $k\ge1$, by partition layer.  After balancing the
parameters, we prove the joint shadow bound under the reference dynamics and
reweight it once to obtain the IVF-to-LLF transition.  The final section
proves positive-density renewal production.

{\bf Note on A.I. assistance.}
The random-geometric mechanisms feeding the accrual argument were developed by the authors.  These include the quilt, seal and wavebreaker constructions, and, in Chapter~\ref{s.seal-sequence-layer}, the host--guest description of a monitored layer, the division into advance and post-passage phases, the passage initial condition, and the replacement--enlargement sealing scheme used after passage.  The authors also developed the basic marginal mechanism for controlling accrual in a single monitored layer, of the form summarised above by~\eqref{e.one-pair-depth-decay}.

The more substantial conceptual contribution of GPT concerned the assembly of these ingredients into the joint control required for the fine transition.  In extended exchanges with the authors, GPT helped shape the reference-dynamics organisation of Chapter~\ref{s.model-problems-accrual}, the conditional exponential-moment and weighting apparatus used to gather one-slide costs across the schedule and across monitored layers, and the strategy of postponing the accumulated acceptance correction until a trajectory-level comparison at the end.  These ideas supplied the analytic framework in which the local and marginal controls are organised into the reference-law multi-layer accrual estimate of Corollary~\ref{c.reference-multilayer-accrual} and then used in the shadow argument of the present chapter to prove the IVF-to-LLF transition.

\section{Well-spaced targets and linear exposure}
\label{ss.well-spaced-exposure}

We now select a positive-density subfamily of the IVF heights whose elements
are well spaced.  The required property is that, for each selected target
\(h_j\), every interval \([x,h_j]\) contains at least
\(c(h_j-x)-C\) IVF heights.  These heights will supply the breaker windows
protecting an approach from the south.
The schedule parameters $m,g$ and breaker parameters
$i_{\rm br},g_{\rm br},b_0$ are treated as fixed in this section;
in Section~\ref{ss.delivery-constant-choice}, they are chosen together with
the fine-allowance scale to balance the accrual and acceptance estimates.
Recall that \(\mathcal H_{\rm vf}=\{H_i:i\in\ivf\}\) is the set of IVF
heights; enumerate its elements as \(e_1<\cdots<e_M\).  Fix \(\gamma\) as in
Lemma~\ref{l.sampling-ivf-points}, and set
\(\Omega_{\rm sat}=\{|\ivf|\ge\gamma L\}\).  For all sufficiently large~\(u\), the deterministic density of the
schedule permits us to fix \(c_E>0\) so that \(\gamma L\ge c_Eu\).
Thus \(M\ge c_Eu\) on \(\Omega_{\rm sat}\); bounded~\(u\) will be
absorbed by enlarging the final constants.
Put \(d_1=e_1\) and
\(d_r=e_r-e_{r-1}\) for \(2\le r\le M\), so that
\(\sum_{r=1}^M d_r\le u\).  Let \(B\) be a positive parameter to be chosen shortly, and call \(e_r\) well-spaced when
\(\max_{1\le q\le r}q^{-1}\sum_{t=r-q+1}^r d_t\le B\).
The one-sided Hardy--Littlewood maximal inequality on the finite gap sequence
gives
\[
 \#\{r:e_r\text{ is not well-spaced}\}
\, \le \, \frac{C_{\rm HL}}B\sum_{r=1}^M d_r
\, \le \, \frac{C_{\rm HL}}B u\, .
\]
Choose \(B\ge2C_{\rm HL}/c_E\), and let \(\mathcal W\) be the set of
well-spaced IVF heights.  Writing \(N=|\mathcal W|\), we have
\(N\ge M/2\ge c_Eu/2>0\) for the sufficiently large~\(u\) under consideration.
The input cluster record and attempted-height schedule determine
\(\mathcal W\).  Throughout this chapter, \(\PP_{\rm ref}\) and
\(\PP_{\rm acc}\) denote the corresponding reference and accepted laws
conditional on these data.

If \(h=e_r\in\mathcal W\) and \(0\le x\le h\), put
\(q=|\mathcal H_{\rm vf}\cap[x,h]|\).  If \(q<r\), then \(e_{r-q}<x\);
if \(q=r\), then \(x\ge0\).  In either case the well-spaced property gives
\(h-x\le Bq\), and hence
\begin{equation}
\label{e.well-spaced-backward-density}
   |\mathcal H_{\rm vf}\cap[x,h]|\, \ge \, B^{-1}(h-x)\, .
\end{equation}
Thus every selected target has a positive density of earlier IVF breaker
heights below it.  In the remainder of the shadow argument
\(h_1<\cdots<h_N\) are the ordered points of \(\mathcal W\); every point of
\(\mathcal H_{\rm vf}\), selected or not, remains available as a breaker.
On~$\Omega_{\rm sat}$ one has~$M\geq\gamma L$, while the whole schedule
has~$n=L$ attempted steps.  Hence the construction above gives the
parameter-uniform bound
\begin{equation}
\label{e.attempted-block-linear}
 n\leq C_{\rm att}N,
 \qquad C_{\rm att}=2/\gamma.
\end{equation}
\begin{lemma}[Breaker count and depth regularity]
\label{l.shadow-breaker-count}
There are positive constants $c_{\rm br},C_{\rm br}$ such that, on
$\Omega_{\rm sat}$, for every $h_j\in\mathcal W$ and every partition layer
$\cS=\slab_{s^-,s^+}$ with $s^+\leq h_j$, at geometric layer depth
$\delta=h_j-s^+$ and every attempted height $H_b\leq s^-$,
\begin{equation}
\label{e.shadow-depth-sensitive-count}
 M_{b,\cS}
 \geq\bigl(c_{\rm br}(s^--H_b)-C_{\rm br}\delta-C_{\rm br}\bigr)_+,
\end{equation}
where \(M_{b,\cS}\) is the quantity in
\eqref{e.ssr-breaker-count-definition}.  Consequently every finite collection of
distinct partition layers paired with heights in $\mathcal W$ is
\(C_{\rm reg}\)-depth regular in the sense of
Definition~\ref{d.ssr-depth-regular}, for one deterministic
\(C_{\rm reg}\) independent of the collection.
\end{lemma}

\noindent{\bf Proof.}
Fix $b,\cS,h_j$ and apply the backward-density bound
\eqref{e.well-spaced-backward-density} with $x=H_b$.  A height
$h\in\mathcal H_{\rm vf}\cap[H_b,h_j]$ is not counted by $M_{b,\cS}$ only if
$h=H_b$ or $h+D_{\rm wb}+r\geq s^-$.  The latter heights lie in
$[s^--D_{\rm wb}-r,h_j]$, whose length is at most
\[
 h_j-s^-+D_{\rm wb}+r
 =\delta+(s^+-s^-)+D_{\rm wb}+r.
\]
The layer width is the fixed value~$2K$, and the attempted heights have uniformly
bounded local density, so at most $C(1+\delta)$ IVF heights are excluded.  This proves
\eqref{e.shadow-depth-sensitive-count}.

If $M_{b,\cS}\leq t$ and $\delta\leq \rho$, the displayed estimate gives
$s^--H_b\leq C(1+t+\rho)$.  Distinct partition layers have distinct lower
faces, spaced by~$2K$, and hence there are at most
$C_{\rm reg}(1+t+\rho)$ such faces.
The constant depends only on the well-spacing bound~$B$, the mesh and schedule
parameters~$K,g$, and the breaker clearance~$D_{\rm wb}+r$; it is uniform over
the chosen collection.
\qed

\begin{lemma}[Linear reference exposure]
\label{l.linear-reference-exposure}
On~\(\Omega_{\rm sat}\), consider~\(h_j\in\mathcal W\)
and a partition layer~\(\cS=\slab_{s^-,s^+}\) with~\(s^+\leq h_j\), and put
\(\delta=h_j-s^+\).  For the pair \((\cS,h_j)\), write
\(a_{b,\cS}\) and \(A_{\cS,h_j}^{\rm ref}\) for the slide weights and
reference exposure of \eqref{e.ssr-slide-weight} and
\eqref{e.reference-exposure}, respectively.  Then there is a finite
constant~\(C\) such that
\begin{equation}
\label{e.exposure-linear}
   A_{\cS,h_j}^{\rm ref}\, \le \, C(1+\delta)\, .
\end{equation}
\end{lemma}

\noindent{\bf Proof.}
Successive attempted heights are separated by at least \(2K(g+1)\), so
for every interval \(I\) of integer heights,
\[
 \#\{b:H_b\in I\}\le 1+\frac{|I|}{2K(g+1)}.
\]
For a slide at height~\(H_b\leq s^-\),
\eqref{e.ssr-final-attenuation-base} gives
\(\kappa=q_{\rm row}^{1/2}\).  Since \(i_{\rm br}\geq2\),
\(q_{\rm row}\leq1-p_{\rm arr}\), and hence
\(\kappa\leq(1-p_{\rm arr})^{1/2}<1\).  The weight definition
\eqref{e.ssr-slide-weight} and Lemma~\ref{l.shadow-breaker-count} therefore give
\[
 a_{b,\cS}
\, \leq \,
 \exp\bigl\{-c\bigl(s^--H_b-C\delta-C\bigr)_+\bigr\}\, .
\]
The preceding bound on the number of attempted heights in an interval makes
the sum of these terms at most \(C(1+\delta)\): the slides within distance
\(O(\delta)\) of the layer contribute at most \(C(1+\delta)\), while the
contributions from slides farther south are geometrically summable.  The
post-passage interval \([s^+,h_j)\) contains \(O(\delta)\) attempted slides,
each with weight one.  This proves~\eqref{e.exposure-linear}.
\qed

The depth-regularity and linear-exposure conclusions above supply the geometric
input for the positive-depth multi-layer argument.  Before turning to it, we
treat the \(k=0\) southside condition separately.

\section{Depth-zero delivery}
\label{ss.zero-depth-delivery}

The depth-zero task is to show that a positive fraction of the selected IVF
heights satisfy their depth-zero fine condition at limelight.  The allowance
\(\af(0)=A_1\), fixed in Section~\ref{ss.positive-renewal-density} as a
sufficiently large multiple of \(\avf(0)=A_0\), provides room for the
retained host after its bounded vertical translation.  We count any fresh
midpiece entering the target band as a possible failure; breaker attenuation
will show that few targets are exposed to such intrusions.

Fix \(H=h_j\in\mathcal W\).  Lemma~\ref{l.ssr-retained-transport}
identifies the retained host at passage as a translate of an input upcluster
when an active slide has occurred; otherwise it is part of the input cluster.
The initial and current head heights both lie in \([u,u+3r]\), so the
translation changes vertex heights by at most~\(3r\).  The inverse image of
the target band \([H-1,H]\) therefore meets only a bounded number of initial
unit bands.  Initial very-fineness at~\(H\), together with
Remark~\ref{r.new-adjacent-band-form}, bounds the host's occupation of the
band, including retained vertices incident to no open edge of the layer, by
\(C_{\rm host}A_0\).  This constant is independent of the coarse scale~\(K\)
and of the later breaker and schedule choices.  Since \(A_1\) was chosen so
that \(C_{\rm host}A_0\leq A_1/2\), the transported host lies comfortably
within the depth-zero fine allowance.

The band \([H-1,H]\) lies in the partition layer whose north face is~\(H\).
No post-passage slide occurs in this layer before limelight.  Thus, for the
depth-zero condition to fail, the guest from the last active slide must reach
the target band.  We measure the reach of these prospective guests along the
ordered list \(h_1<\cdots<h_N\).

For an attempted step~\(b\) at height~\(H_b\), define its
\emph{depth-zero reach}~\(Z_b^0\) as follows.  Set \(Z_b^0=0\) when the
step is dormant.  Otherwise let \(j_b=\min\{j:H_b<h_j\}\), with \(Z_b^0=0\)
if this set is empty.  For \(j\geq j_b\), let \(\mathsf H^0_{b,j}\) be the
event that the midpiece created at step~\(b\) has a vertex in the target band
\([h_j-1,h_j]\), and put
\[
 Z_b^0=
 \max\Bigl(
   \{j-j_b+1:\mathsf H^0_{b,j}\text{ occurs}\}\cup\{0\}
 \Bigr).
\]
When \(Z_b^0>0\), the set of indices
\(\llbracket j_b,j_b+Z_b^0-1\rrbracket\) is the \emph{reach interval} of
step~\(b\).  It contains the index of every well-spaced target whose band
is reached by this midpiece, and may include intervening indices whose bands
are not reached.

\begin{lemma}[Pathwise bound on depth-zero failures]
\label{l.zero-depth-pathwise-domination}
Let \(D_N^0\subseteq\llbracket1,N\rrbracket\) consist of those indices
\(j\) for which the depth-zero condition \({\bf A}_{\rm f}^-(0)\) fails
at the limelight time of \(h_j\in\mathcal W\).  On~\(\Omega_{\rm sat}\), along
the accepted trajectory,
\begin{equation}
\label{e.zero-depth-pathwise-domination}
 |D_N^0|\leq\sum_{b=1}^n Z_b^0.
\end{equation}
\end{lemma}

\noindent{\bf Proof.}
If no active slide occurs before a target reaches limelight, its initially
very-fine input bound and \(A_1\geq A_0\) give the depth-zero condition.
Otherwise, for a failed target \(h_j\), let~\(b\) be the last active slide
before limelight.  By Lemma~\ref{l.ssr-passage-decomposition} and the
vertex-level DMU partition, the cluster material in the target band at
passage comprises retained host vertices and vertices of the midpiece
created at~\(b\).  Subsequent dormant attempts leave the cluster unchanged,
and there is no post-passage slide before limelight.  The retained host
satisfies the condition on its own, so failure forces
\(\mathsf H^0_{b,j}\) to occur.  Thus \(j\) belongs to
\(\llbracket j_b,j_b+Z_b^0-1\rrbracket\).  For each failed target we have identified
an active slide whose reach interval contains its index; summing the
lengths of these intervals proves the claim, even if they overlap.
\qed

It remains to control the sum of these depth-zero reaches under the
accepted law.  Recall the marked history~\(G_b\) of
Chapter~\ref{s.model-problems-accrual}, and put
\(\mathcal G_b^0=\sigma(G_b)\) under the accepted law.  Since no rejected
fixed-shift output occurs on an accepted trajectory, this filtration
records the attempted-height schedule, cluster history and proposal marks
up to and including step~\(b\), with a fixed null mark at dormant steps.  Given the
pre-step cluster, the step height, the shift mark and the output cluster,
the DMU decomposition identifies the new midpiece as the part left after
removing the retained downcluster and shifted retained upcluster.  Hence
\(Z_b^0\) is measurable with respect to \(\mathcal G_b^0\), the step-\(b\)
sigma-field of the marked filtration used by the accepted one-step comparison
in Lemma~\ref{l.accepted-rare-event-transfer}.

The current slide's breaker window and the recovered windows at intervening
IVF targets limit the reach of its midpiece.  The next lemma bounds the
depth-zero reach of each attempted step under the accepted law and gives the
exponential-moment estimate needed to sum these contributions.

\begin{lemma}[Accepted depth-zero reach tail]
\label{l.zero-reach-tail}
Let \(C_{\rm A}\) and \(\delta_{\rm A}\) be the constants of
Lemma~\ref{l.accepted-rare-event-transfer}, with
\(\delta_{\rm A}\) chosen independently of~\(\Cfine\), and set
\(C_{\rm z}=2C_{\rm A}\), \(\delta_{\rm z}=\delta_{\rm A}\), and
\(q_{\rm z}=q_{\rm row}^{\delta_{\rm z}}\).  Whenever the row has been
chosen so that \(q_{\rm z}<1/2\), every attempted step~\(b\) satisfies
\[
 \PP_{\rm acc}(Z_b^0\geq t\mid\mathcal G_{b-1}^0)
 \leq C_{\rm z}q_{\rm z}^t,
 \qquad t\geq1.
\]
Moreover
\begin{equation}
\label{e.zero-Cfine-dependence}
 C_{\rm z}\leq C(2+\Cfine)^{D_{\rm z}}
\end{equation}
for fixed \(C,D_{\rm z}<\infty\).  Consequently, for every \(\zeta>0\) with
\(e^\zeta q_{\rm z}<1\),
\begin{equation}
\label{e.zero-reach-mgf}
 \EE_{\rm acc}[e^{\zeta Z_b^0}\mid\mathcal G_{b-1}^0]
 \leq
 1+\frac{C_{\rm z}(e^\zeta-1)q_{\rm z}}
          {1-e^\zeta q_{\rm z}}.
\end{equation}
\end{lemma}

\noindent{\bf Proof.}
There is nothing to prove at a dormant step.  Condition on an active step,
fix \(j\geq j_b\), and fix an admissible proposed head~\(v\).  If
\(B_{\rm top}\) is the northern height of the current-slide breaker window
for this fixed shift, the window geometry~\eqref{e.breaker-profile-band}
gives \(B_{\rm top}\leq H_b+3r+b_0+C_0=H_b+D_{\rm wb}\).  The schedule
gap~\eqref{e.random-mesh-gap} then gives
\(H_b+D_{\rm wb}<H_{b+1}-2K\), so this window lies wholly below the
zero-depth layer of every later selected target.

For \(j_b\leq\ell<j\), the selected height~\(h_\ell\) is initially very
fine, and its breaker window has top at most~\(h_\ell+D_{\rm wb}\).
The clearance and attempted-height spacing give
\[
 h_\ell+D_{\rm wb}+r\leq h_\ell+g<h_j-2K.
\]
Thus these \(j-j_b\) recovered windows also lie below the zero-depth layer
at~\(h_j\).

Let \(E_{b,j}^v\) be the fixed-shift event that acceptance occurs and the
output midpiece reaches the height interval~\([h_j-1,h_j]\).  On this event a midpiece
vertex in the target band is incident to an undiscovered open edge, by the
fixed-shift record and the DMU edge decomposition.  Since the band lies
strictly north of the supplying window and the recovered windows at
\(h_{j_b},\ldots,h_{j-1}\), we have
\(E_{b,j}^v\subseteq A_{j-j_b}^\circ\) for the wavebreaker experiment
of Proposition~\ref{p.origin-and-successive-wavebreakers}.  That proposition
gives, conditionally on the stopped fixed record used by that search,
\(\Rlaw^v_{\uclust,\dclust}(E_{b,j}^v\mid\cF_0)
\leq q_{\rm row}^{\,j-j_b+1}\), uniformly in the retained record and in~\(v\).  Since
\(E_{b,j}^v\) is determined by the reconstructed cluster and the fixed
retained record, the marked-kernel identity
\eqref{e.marked-kernel-disintegration} gives the same bound for the
auxiliary-averaged kernel in the hypothesis of
Lemma~\ref{l.accepted-rare-event-transfer}.  Under the accepted law,
\(\mathsf H^0_{b,j}\) is the marked event \(E_{b,j}^{V_b}\).  Applying that lemma with the marked pre-step
history~\(G_{b-1}\), whose sigma-field is~\(\mathcal G_{b-1}^0\), gives
\(\PP_{\rm acc}(\mathsf H^0_{b,j}\mid\mathcal G_{b-1}^0)
\leq C_{\rm A}q_{\rm z}^{\,j-j_b+1}\).  If \(Z_b^0\geq t\), one of these
events occurs for \(j\geq j_b+t-1\).  Under the stated row condition,
summing the geometric series gives the asserted tail bound.
Equation~\eqref{e.accepted-rare-Cfine} gives
\eqref{e.zero-Cfine-dependence}.  The tail-sum identity gives
\eqref{e.zero-reach-mgf}.
\qed

The conditional exponential-moment bound is suited to the sum in
\eqref{e.zero-depth-pathwise-domination}: it can be iterated over the
\(n\leq C_{\rm att}N\) attempts without assuming that their reaches are
independent.  To make the cumulative moment cost small, the breaker rows
must make the one-step moment factor sufficiently close to one.  Take the exponent
in~\eqref{e.zero-reach-mgf} to be~\(\zeta_0=1\).  The final row choice in
Lemma~\ref{l.delivery-constant-compatibility} will ensure the condition in
the following proposition.

\begin{proposition}[Depth-zero delivery]
\label{p.zero-depth-delivery}
Suppose that the row parameters satisfy
\begin{equation}
\label{e.zero-row-standing-choice}
 e^{\zeta_0} q_{\rm z}<1,
 \qquad
 \log\left(
 1+\frac{C_{\rm z}(e^{\zeta_0}-1)q_{\rm z}}
          {1-e^{\zeta_0} q_{\rm z}}
 \right)
 <\frac{\zeta_0}{32C_{\rm att}}.
\end{equation}
Then, on \(\Omega_{\rm sat}\), conditional on the input record and schedule,
\[
 \PP_{\rm acc}(|D_N^0|\geq N/4)\leq e^{-N/8}.
\]
\end{proposition}

\noindent{\bf Proof.}
The first inequality in~\eqref{e.zero-row-standing-choice} gives
\(q_{\rm z}<e^{-\zeta_0}<1/2\), so Lemma~\ref{l.zero-reach-tail} applies.
Combine it with~\eqref{e.zero-depth-pathwise-domination}.  Since \(Z_b^0\) is
\(\mathcal G_b^0\)-measurable and satisfies
\eqref{e.zero-reach-mgf}, that bound may be iterated in the filtration
\((\mathcal G_b^0)_{b=0}^n\).  Using \(n\leq C_{\rm att}N\) and Markov's
inequality, we obtain
\[
 \PP_{\rm acc}(|D_N^0|\geq N/4)
 \leq
 \exp\{-\zeta_0 N/4+\zeta_0 N/32\}
 \leq e^{-N/8}.
\]
\qed

\section{Towers and shadows}
\label{ss.physical-shadows}

In our \hfff{towers-shadows}{tower-and-shadow model} for the positive-index southside conditions,
each partition layer is a tower that casts shade on a selected target when
some \({\bf A}^{-}_{\rm f}(k)\), \(k\ge1\), assigned to that layer fails.

Take \(\mathcal W=\{h_1<\cdots<h_N\}\).  List as
\(\mathcal Q_1,\ldots,\mathcal Q_{q_N}\), in south-to-north order, the
partition layers \(\mathcal Q_q=\slab_{a_q,b_q}\) whose north face satisfies
\(0\le b_q\le h_N\).  These are the only layers that need be considered when checking the southside
fine conditions at the selected targets.  Indeed every selected height is a face of the layer
partition.  If the band \([h_j-k-1,h_j-k]\) is not wholly below~\(0\), the
layer containing it lies south of the face~\(h_j\) and has north face in
\([0,h_j]\).  A unit band lying wholly below~\(0\) contains no cluster
vertex, since the origin is the cluster base.
Lemma~\ref{l.layer-partition} gives pairwise-disjoint layer edge sets,
common vertex faces for adjacent layers, and width \(b_q-a_q=2K\).  On
\(\Omega_{\rm sat}\), one has
\(q_N\le h_N/(2K)+2\le u/(2K)+2\) and \(N\ge c_Eu/2\),
and hence \(q_N\le C_QN\) for a fixed \(C_Q<\infty\), after enlarging it
to cover bounded~\(u\).

For every \(q\le q_N\), the set \(\{j:b_q\le h_j\}\) is non-empty by
construction.  Put \(j_0(q)=\min\{j:b_q\le h_j\}\).  For
\(j\ge j_0(q)\), write \(\delta_{q,j}=h_j-b_q\), and let
\(\mathcal D_{q,j}=\{k\in\Nzero:\cS(k)=\mathcal Q_q
\text{ when the target is }h_j\}\).  Here \(\delta_{q,j}\) is the \emph{geometric layer depth}; the elements of
\(\mathcal D_{q,j}\) are the fine-condition displacements.  The distinction is
important when \(b_q=h_j\): then \(\delta_{q,j}=0\), but the same layer also contains
bands at positive displacement.  Because the layer has width~\(2K\),
\(|\mathcal D_{q,j}|\le 2K+1\).  Every positive
\(k\in\mathcal D_{q,j}\) satisfies
\(k\ge \delta_{q,j}\) and \(k=\delta_{q,j}+O(K)\); when \(\delta_{q,j}=0\), positivity
also gives \(k\ge1\).  Indeed the top of the unit band is at most the north
face~\(b_q\), so \(h_j-k\le b_q\).
At the limelight time of \(h_j\), set
\(G_{q,j}:=|\downcluster_{h_j}\cap\mathcal Q_q|_*\).
Because \( |\cdot|_* \) is defined through vertex support, a cluster
vertex on a common face is counted in either adjacent layer containing that
face, regardless of which layer contains its incident cluster edges.

We say that \(\mathcal Q_q\) \emph{casts shade on} \(h_j\)
when \({\bf A}^{-}_{\rm f}(k)\) fails at the limelight time of \(h_j\)
for some positive \(k\in\mathcal D_{q,j}\).
Since the coarse-box scale \(K\) is fixed, each coarse box contains at most
\(C<\infty\) lattice vertices.  Hence the number of cluster vertices
in a set is at most \(C\) times its coarse cardinality.  Put
\(c_{\rm sh}=C^{-1}\).  If \(k\) is responsible for shade on
\(h_j\), then \(G_{q,j}\ge c_{\rm sh}\af(k)\).
Since \(\af\) is non-decreasing on the positive integers, the bound is at least
\(c_{\rm sh}\af(\delta_{q,j})\) when \(\delta_{q,j}\ge1\), and at least
\(c_{\rm sh}\af(1)\) when \(\delta_{q,j}=0\).

The shadow length of \(\mathcal Q_q\), measured in target indices, is
\begin{equation}
\label{e.shadow-interval-definition}
 X_q\, = \,
 \max\Bigl(
 \{j-j_0(q)+1:\mathcal Q_q\text{ casts shade on }h_j,\ j\, \ge \, j_0(q)\}
 \cup\{0\}\Bigr)\, .
\end{equation}
When \(X_q>0\), the \emph{shadow interval} of \(\mathcal Q_q\) is
\(\llbracket j_0(q),j_0(q)+X_q-1\rrbracket\).
It therefore indexes every target between the first target at or north of
the layer and the last target on which it casts shade.  A target
\(h_j\) is \emph{sunny} when its index \(j\) belongs to none
of these shadow intervals.  Since the selected heights are distinct integer
heights, \(\delta_{q,j}\ge j-j_0(q)\); for fixed~\(q\), the map
\(j\mapsto \delta_{q,j}=h_j-b_q\) is injective.

\section{Balancing the delivery estimates}
\label{ss.delivery-constant-choice}

The positive-depth fine allowance gives us one source of slack.  Increasing
\(\Cfine\) enlarges the volume budget available to the root and front
thresholds in Corollary~\ref{c.reference-multilayer-accrual}.  The same choice
has a cost when the reference estimates are reweighted to the accepted slide
dynamics.  We compensate by increasing the breaker-row length \(i_{\rm br}\), which makes
finite passage through a row exponentially less likely.  Longer rows in turn
require a wider gap between attempted slide heights.  Since all these scales
are fixed independently of~\(u\), the resulting loss in attempted-height
density is only by a fixed factor.

For a positive-index fine condition, Corollary~\ref{c.reference-multilayer-accrual}
controls large layer accrual through two alternatives: too many boundary roots,
or excessive spread of the replacement front.  In this section we prove a
compatibility lemma showing that the corresponding thresholds can be set at
the scale of the positive-index fine allowance, while choosing the breaker-row
and attempted-height schedule parameters so that the depth-zero and
shadow-moment requirements are also met.

Let \(C_{\rm acc}\ge1\), \(C\) and \(c\) be constants valid in
Corollary~\ref{c.reference-multilayer-accrual}, with the upper-bound
constants enlarged as needed to absorb the constant
in~\eqref{e.exposure-linear}.  For the positive-depth fine allowance
coefficient~\(\Cfine\), put
\[
 \rho_{\rm f}=1+\log(2+\Cfine).
\]
Let~\(D_{\rm rep}>0\) be a constant.  For~\(z\ge1\), set
\[
\begin{aligned}
 \tau_{\rm f}(z)&=D_{\rm rep}\rho_{\rm f}(1+z),\\
 q_{\rm r}(z)&=C_{\rm acc}\avf(z)
   +\frac{c_{\rm sh}^\alpha}{4C_{\rm acc}}
      \frac{\af(z)^\alpha}
       {(1+z+\tau_{\rm f}(z))^{d-1}},\\
 \eta_{\rm r}(z)&=\min\left\{1,
   \exp\left[C(1+z)
     -c(1+z)^{-\sigma_{\rm wt}}
          \bigl(q_{\rm r}(z)-C_{\rm acc}\avf(z)\bigr)\right]\right\},\\
 \eta_{\rm f}(z)&=\min\{1,\ q_{\rm r}(z)e^{-c_{\rm rep}\tau_{\rm f}(z)}\},\\
 \eta(z)&=\eta_{\rm r}(z)^{1/2}+\eta_{\rm f}(z)^{1/2},\\
 \overline G(r)&=2\sum_{z\geq\max\{1,r-1\}}\eta(z),\qquad r\geq1.
\end{aligned}
\]
Here \(\eta_{\rm r}\) and \(\eta_{\rm f}\) are the root and front error
bounds, and \(\eta\) combines them.  The sum \(\overline G(r)\) runs over
the possible witnessing depths for a shadow of length at least~\(r\).  The
factor~\(2\) accounts for the only multiplicity in this indexing: geometric
depths~\(0\) and~\(1\) both enter the estimate at~\(z=1\).

\begin{lemma}[Compatible parameter choice]
\label{l.delivery-constant-compatibility}
With \(\zeta_0=1\), let the even coarse scale~\(K\) be fixed
sufficiently large for the coarse-grained, apt-box, seal-completion and
layer-episode estimates, and let \(A_0,A_1,C_1\) be fixed as in
Section~\ref{ss.positive-renewal-density}.  The remaining allowance,
breaker-row and schedule parameters, together with the multiplier
\(D_{\rm rep}\) and a moment exponent \(t>0\), may be chosen so that
\(A_1\leq\af(1)\),
\begin{gather}
 \frac{t}{8}>C_{\rm L}C_{\rm att}+2,
 \label{e.compat-positive-margin}\\
 \log\left(
  1+(e^t-1)\sum_{r\geq1}e^{t(r-1)}\overline G(r)
 \right)<\frac{t}{8C_Q},
 \label{e.compat-positive-mgf}
\end{gather}
and \(\overline G(1)<1/2\).  For every~\(z\ge1\),
\begin{equation}
\label{e.compat-root-front-threshold}
 c_{\rm sh}^{\alpha}\af(z)^{\alpha}
 \ge C_{\rm acc}q_{\rm r}(z)(1+z+\tau_{\rm f}(z))^{d-1},
\end{equation}
so the root-count and front-excess thresholds satisfy
\eqref{e.reference-accrual-threshold-condition}.  The depth-zero condition
\eqref{e.zero-row-standing-choice} also holds.  Finally,
\begin{equation}
\label{e.compat-row-errors}
 C_{\rm cl}q_{\rm row}^{\delta_{\rm cl}}<\frac18,
 \qquad
 \frac{C_{\rm A}q_{\rm row}^{\delta_{\rm A}}}
      {1-q_{\rm row}^{\delta_{\rm A}}}<\frac14,
\end{equation}
with the constants of \eqref{e.clean-future-breakers}
and \eqref{e.accepted-rare-Cfine}.
\end{lemma}

The root/front threshold bound \eqref{e.compat-root-front-threshold} and
\(\overline G(1)<1/2\) are used in
Section~\ref{ss.joint-shadow-bounds}; both
\eqref{e.compat-positive-margin} and \eqref{e.compat-positive-mgf} are used in
Section~\ref{ss.ivf-llf-transition}.  We record the bounds in \eqref{e.compat-row-errors} for use in
Chapter~\ref{s.renewal-survival}.

\noindent{\bf Proof of Lemma~\ref{l.delivery-constant-compatibility}.}
We treat \(R\) as a parameter to be taken sufficiently large and set
\(\Cfine=e^R\).  The likelihood-moment exponent~\(\theta\) was fixed
in~\eqref{e.theta-choice}.  This choice remains admissible as~\(R\)
varies, since Corollary~\ref{c.tether-head-summability} allows
\(\gamma_{\rm T}\) to be chosen independently of~\(\Cfine\).
The fixed numerical multipliers \(B_{\rm row}\), \(D_2\) and
\(D_{\rm rep}\) will be chosen, in that order, independently of~\(R\);
only after these choices are made do we take \(R\) sufficiently large.  The resulting
breaker count, window length and schedule parameters are then determined from
\(R\).  The logarithmic dependence in~\eqref{e.likelihood-Cfine-log} makes the
acceptance cost of order~\(R\), while \(\af(z)^\alpha\) contains the factor
\(e^{\alpha R}\).  The remaining choices ensure that this gain dominates
the polynomial losses.

We strengthen the breaker rows by taking their number of test layers
proportional to~\(R\).  Choose a fixed admissible test-layer separation \(g_{\rm br}\) (for example
the least admissible multiple of \(2K\)).  Take the even breaker count
\(i_{\rm br}=2\lceil B_{\rm row}R\rceil\), with
\(C_0=3g_{\rm br}i_{\rm br}\), and take a fixed admissible \(b_0>r\).  Since
\(q_{\rm row}=(1-p_{\rm arr})^{i_{\rm br}/2}\), putting
$c_{\rm row}=-\log(1-p_{\rm arr})>0$ gives
\(q_{\rm row}\le e^{-c_{\rm row}B_{\rm row}R}\).
Thus a row becomes exponentially harder to pass as \(R\) grows.  The same
choice also drives down the probability that the midpiece formed at a slide
reaches the reserved rows.  Lengthening the rows does not alter the
uniformity of the comparison estimates: Corollaries~\ref{c.tether-head-summability}
and~\ref{c.tether-likelihood-moment} are uniform in the breaker-profile
parameters, and \(c_{\rm rep}\) is the fixed replacement-front rate of
Proposition~\ref{p.rebuild-linear-speed}.  Hence the constants in
\eqref{e.accepted-rare-Cfine}, \eqref{e.zero-Cfine-dependence} and
\eqref{e.clean-future-Cfine-dependence} are at most \(e^{DR}\), with fixed \(D\),
while \(\delta_{\rm A}\) and \(\delta_{\rm cl}\) are independent
of~\(R\).
Choose \(B_{\rm row}\) so large that
\[
  B_{\rm row}>\frac{D+1}
  {c_{\rm row}\min\{\delta_{\rm A},\delta_{\rm cl}\}}.
\]
With this choice, $q_{\rm row}$ decays quickly enough to offset
the growth of the accepted-law prefactors.  For all sufficiently large~\(R\), this gives
\eqref{e.zero-row-standing-choice} and \eqref{e.compat-row-errors}.

Take \(g\) to be the least integer satisfying the direction-independent
version of~\eqref{e.random-mesh-gap}, with \(r\) replaced by~\(1\), and
set \(m=2g\).  Since \(C_0=O(R)\), we have \(g+m=O(R)\).  Thus
the attempted-height density is of order \((1+R)^{-1}\), and remains positive
once~\(R\) is fixed.  The failure-count estimate of
Lemma~\ref{l.nhigh} is uniform in the breaker count and separation.
With \(m=2g\), Lemma~\ref{l.sampling-ivf-points} permits
\(\gamma=1/24\), independently of~\(R\).  On~\(\Omega_{\rm sat}\),
\(M=|\ivf|\ge\gamma L\), with \(L\asymp u/(m+g)\).  The constant
\(c_E\) in the well-spaced extraction may therefore be chosen with
\(c_E^{-1}\le C(m+g)\le C(1+R)\).  Since
\(N\ge M/2\ge c_Eu/2\) and there are \(O(u)\) relevant partition layers,
the layer-count constant \(C_Q\) of Section~\ref{ss.physical-shadows}
may likewise be taken with \(C_Q\le C(1+R)\).

We next track the depth-regularity constants.  The Hardy--Littlewood
cutoff in the well-spaced extraction may be taken as
\(B=2C_{\rm HL}/c_E=O(1+R)\).  Also
\(D_{\rm wb}=3r+b_0+C_0=O(R)\), while the upper and lower mesh spacings are
both of order \(m+g=O(R)\).  Inspecting the proof of
Lemma~\ref{l.shadow-breaker-count} now shows that the constants
\(c_{\rm br}^{-1}\) and \(C_{\rm br}\) are bounded by fixed powers of
\(1+R\).  Solving \eqref{e.shadow-depth-sensitive-count} for
\(s^- - H_b\), and then counting the distinct lower faces spaced by \(2K\),
gives, for some fixed \(D_1<\infty\),
\(C_{\rm reg}\le C(1+R)^{D_1}\).  For all sufficiently large~\(R\),
\(q_{\rm row}\le1/2\).  The moment bound~\eqref{e.ssr-input-Creg-dependence}
with \(q_0=1/2\) has polynomial dependence on \(1+R\), while
\eqref{e.ssr-final-attenuation-base} gives \(\kappa\le1/\sqrt2\).
The bounds on \(c_{\rm br}^{-1}\), \(C_{\rm br}\) and the attempted-height
spacing are also polynomial in \(1+R\).  The geometric summation of the
slide weights in Lemma~\ref{l.linear-reference-exposure} therefore gives
the same polynomial bound for the constant in~\eqref{e.exposure-linear}.
The final uniformity assertion of
Corollary~\ref{c.reference-multilayer-accrual}, applied with \(q_0=1/2\),
then allows \(C_{\rm acc}\), \(C\) and \(c^{-1}\) to be bounded by fixed
powers of \(1+R\).

Finally, \eqref{e.likelihood-Cfine-log} gives
$C_{\rm L}\le C(1+R)$.  Choose the shadow moment exponent
$t=D_2(1+R)$, with~$D_2$ sufficiently large for
\eqref{e.compat-positive-margin}.  The definition of~$q_{\rm r}$ and the
polynomial bounds on~$C_{\rm acc}$ give, for~$z\ge1$,
\(\log(2+q_{\rm r}(z))\le C_*(1+R)+C_*\log(1+z)\),
where~$C_*$ is fixed independently of~$D_{\rm rep}$: increasing
$D_{\rm rep}$ only increases the denominator in~$q_{\rm r}$.
Choose~$D_{\rm rep}$ once so that
\[
 \frac{c_{\rm rep}D_{\rm rep}}{2}\ge C_*+3D_2+1.
\]
Since~$\rho_{\rm f}\ge1+R$, this choice gives enough front decay to
absorb the moment weight~$e^{tz}$.

We record the root term separately.  For~$z\ge1$,
\[
 \af(z)^\alpha
 =e^{\alpha R}A_0^{2\alpha}(1+z)^{\Bfine\alpha},
\]
and
\[
 1+z+\tau_{\rm f}(z)\le C(1+R)(1+z).
\]
Hence, after the polynomial dependence of the reference-law constants on
$C_{\rm reg}$ is inserted, there is a fixed polynomial~$P$ such that
\[
 (1+z)^{-\sigma_{\rm wt}}
   \bigl(q_{\rm r}(z)-C_{\rm acc}\avf(z)\bigr)
 \ge
 P(R)^{-1}e^{\alpha R}(1+z)^{p_*},
 \qquad
 p_*=\Bfine\alpha-(d-1)-\sigma_{\rm wt}>1.
\]
The strict inequality is the margin recorded after
\eqref{e.ssr-depth-weight}.  Thus
\begin{equation}
\label{e.shadow-root-envelope-large-R}
 \eta_{\rm r}(z)^{1/2}
 \le
 \exp\left\{P(R)(1+z)
      -P(R)^{-1}e^{\alpha R}(1+z)^{p_*}\right\}.
\end{equation}

For the front term, the bound on~$\log(2+q_{\rm r}(z))$ and the choice
of~$D_{\rm rep}$ give
\begin{equation}
\label{e.shadow-front-envelope-large-R}
 e^{tz}\eta_{\rm f}(z)^{1/2}
 \le \exp\{-(2D_2+1)(1+R)(1+z)\}
\end{equation}
for~$z\ge1$ and all sufficiently large~$R$.
By the definition of~$\overline G$, reversing the sums gives
\[
 \sum_{r\ge1}e^{t(r-1)}\overline G(r)
 =2\sum_{z\ge1}\eta(z)\sum_{r=1}^{z+1}e^{t(r-1)}
 \le \frac{2}{1-e^{-t}}\sum_{z\ge1}e^{tz}\eta(z).
\]
For sufficiently large~$R$, the negative exponent in
\eqref{e.shadow-root-envelope-large-R} dominates both its positive term
and the moment weight~$e^{tz}$, uniformly for~$z\ge1$.
Together with~\eqref{e.shadow-front-envelope-large-R} and~$p_*>1$,
this gives, for another fixed polynomial~$P_1$,
\[
 \sum_{r\ge1}e^{t(r-1)}\overline G(r)
 \le
 C e^{-c e^{\alpha R}/P_1(R)}
 + C e^{-(2D_2+1)(1+R)}.
\]
Multiplication by~$e^t-1\le e^{D_2(1+R)}$ leaves a front contribution
at most~$C e^{-(D_2+1)(1+R)}$, while the root contribution tends to zero
still faster.  Consequently
\[
 (e^t-1)\sum_{r\ge1}e^{t(r-1)}\overline G(r)\longrightarrow0.
\]
Since~$C_Q\le C(1+R)$ and~$t=D_2(1+R)$, the right-hand side of
\eqref{e.compat-positive-mgf} is bounded away from zero uniformly in~$R$.
This proves~\eqref{e.compat-positive-mgf}; the same estimate without the
moment weight gives~$\overline G(1)<1/2$ for all sufficiently large~$R$.

It remains to verify that the thresholds just chosen are admissible in
Corollary~\ref{c.reference-multilayer-accrual}.  After multiplication by
$C_{\rm acc}(1+z+\tau_{\rm f}(z))^{d-1}$, the second term in the definition
of~$q_{\rm r}(z)$ uses exactly one quarter of the available quantity
$c_{\rm sh}^\alpha\af(z)^\alpha$.  Thus it suffices to verify
\[
 C_{\rm acc}^2\avf(z)(1+z+\tau_{\rm f}(z))^{d-1}
 \le \frac34c_{\rm sh}^\alpha\af(z)^\alpha,
 \qquad z\ge1.
\]
Since $C_{\rm acc}$ grows at most polynomially in~$1+R$, the ratio of
the left-hand side to the right-hand side is, after enlarging the fixed
polynomial~$P$ if necessary, at most
\[
 P(R)e^{-\alpha R}(1+z)^{d+2-\Bfine\alpha}.
\]
This tends uniformly to zero because~$\Bfine\alpha>d+2$.  Thus the
threshold condition holds for every~$z\ge1$ when~$R$ is sufficiently large.
Taking~$R$ larger if necessary also ensures
\(A_1\leq\af(1)\), completing the choice of parameters.

The one-pair estimate \eqref{e.one-pair-depth-decay} from the chapter
opening now follows for a selected IVF target, with the convention that
$V=0$ on reference histories discontinued before limelight.  For a single
layer--target pair at depth~$z$, take the volume, root and front thresholds
to be $c_{\rm sh}\af(z)$, $q_{\rm r}(z)$ and $\tau_{\rm f}(z)$.
These choices satisfy \eqref{e.reference-accrual-threshold-condition}, and
the bound follows from \eqref{e.shadow-root-envelope-large-R} and
\eqref{e.shadow-front-envelope-large-R}. \qed

\section{Joint shadow bounds from the multi-layer estimate}
\label{ss.joint-shadow-bounds}

The shadow variables are indexed by distinct partition layers.  A lower bound
on several shadow lengths may be witnessed by choosing, for each layer, a
target at the required distance on which it still casts shade.  Each choice
gives a fixed layer--target family, depth regular by
Lemma~\ref{l.shadow-breaker-count}, to which
Corollary~\ref{c.reference-multilayer-accrual} applies.  Summing over
the possible witnessing indices gives product-form upper bounds on the
joint tails of the shadow lengths, without assuming independence of the
shadows themselves.  These bounds in turn yield exponential-moment
control of the total shadow length.

The bounds are uniform in the fixed input record and attempted-height
schedule.  On~\(\Omega_{\rm sat}\), every monitored target below belongs to
\(\mathcal H_{\rm vf}\), so the factor \(\one_{\mathsf{IVF}_I}\) in
Corollary~\ref{c.reference-multilayer-accrual} is identically one for the
collections considered here.  On the reference space, a target is counted as
shaded only if the reference history is live at its limelight time; the
qualification is automatic under~\(\PP_{\rm acc}\).

Take \(\Cfine,i_{\rm br},g,m,D_{\rm rep}\) and the shadow moment
exponent \(t>0\) as chosen in
Lemma~\ref{l.delivery-constant-compatibility}.  Then \(\overline G(1)<1/2\), and both
\eqref{e.compat-positive-margin} and \eqref{e.compat-positive-mgf} hold.

\begin{proposition}[Joint shadow bound]
\label{p.shadow-tail}
On~\(\Omega_{\rm sat}\), for every
\(I\subseteq\intint{q_N}\) and integers \(r_q\ge1\),
\begin{equation}
\label{e.joint-shadow-tail}
 \PP_{\rm ref}\left(
   X_q\ge r_q\text{ for every }q\in I\right)
 \le \prod_{q\in I}\overline G(r_q).
\end{equation}
\end{proposition}

\noindent{\bf Proof.}
By the definition of~\(X_q\), on the joint event in the proposition there is,
for each~\(q\in I\), an index~\(j_q\) satisfying
\[
 j_0(q)+r_q-1\le j_q\le N
\]
such that \(\mathcal Q_q\) casts shade on~\(h_{j_q}\).  We take a union bound
over the possible choices of these witnessing indices~\((j_q)_{q\in I}\).
Fix one such choice.  For each~\(q\), put
\[
 \delta_q=\delta_{q,j_q}=h_{j_q}-b_q,
 \qquad z_q=\max\{\delta_q,1\}.
\]
Every positive displacement responsible for shade on~$h_{j_q}$ is at
least~$z_q$.  Since the fine allowance is non-decreasing on the positive
integers, shade therefore implies
\[
 G_{q,j_q}\ge c_{\rm sh}\af(z_q).
\]

The layers \(\mathcal Q_q\), \(q\in I\), are distinct, and
Lemma~\ref{l.shadow-breaker-count} makes the corresponding layer--target
pairs a \(C_{\rm reg}\)-depth-regular family.  With the input record and
attempted-height schedule fixed, these layer--target pairs, their
depths~\(z_q\) and the corresponding thresholds are deterministic.
By the shading convention, the corresponding joint shading event is
contained in the event that the reference history is live at the limelight time of every~\(h_{j_q}\).
Apply Corollary~\ref{c.reference-multilayer-accrual} to the pairs
\((\mathcal Q_q,h_{j_q})\), taking \(L_q=c_{\rm sh}\af(z_q)\),
\(q_{\rm r}(z_q)\) for the root count and \(\tau_{\rm f}(z_q)\) for the
front excess.  Since \(\delta_q\le z_q\), monotonicity of~\(\avf\)
gives \(q_{\rm r}(z_q)\ge C_{\rm acc}\avf(\delta_q)\), and
\eqref{e.compat-root-front-threshold} gives
\eqref{e.reference-accrual-threshold-condition} at depth~\(\delta_q\).
Moreover
\[
 \mathsf w_q=(1+\delta_q)^{-\sigma_{\rm wt}}
   \ge (1+z_q)^{-\sigma_{\rm wt}},
 \qquad
 A_{\mathcal Q_q,h_{j_q}}^{\rm ref}
   \le C(1+\delta_q)\le C(1+z_q).
\]
Also, by monotonicity of the very-fine allowance and the definition
of~\(q_{\rm r}\),
\[
 q_{\rm r}(z_q)-C_{\rm acc}\avf(\delta_q)
 \ge q_{\rm r}(z_q)-C_{\rm acc}\avf(z_q)\ge0.
\]
Thus the boundary-count contribution in
\eqref{e.reference-source-envelope} is bounded by~\(\eta_{\rm r}(z_q)\),
while the front-excess contribution in
\eqref{e.reference-front-envelope} is bounded by~\(\eta_{\rm f}(z_q)\).
Corollary~\ref{c.reference-multilayer-accrual} therefore gives, for these
fixed indices,
\[
 \PP_{\rm ref}\left(
   \mathcal Q_q\text{ casts shade on }h_{j_q}
       \text{ for every }q\in I\right)
 \le \prod_{q\in I}\eta(z_q).
\]

For fixed~\(q\), the map \(j\mapsto \delta_{q,j}\) is injective, while the passage
from \(\delta_q\) to~\(z_q\) has multiplicity at most two: only geometric
depths~\(0\) and~\(1\) both give \(z_q=1\).  Since
\(j_q\ge j_0(q)+r_q-1\), we have
\(z_q\ge\max\{1,r_q-1\}\).  The sum over possible choices of
\((j_q)_{q\in I}\) factorises because each index ranges over its own
interval.  The union bound therefore gives at most
\[
 \prod_{q\in I}
  \sum_{z\ge\max\{1,r_q-1\}}2\eta(z)
 =\prod_{q\in I}\overline G(r_q),
\]
which is~\eqref{e.joint-shadow-tail}.
\qed

Introduce an auxiliary non-negative integer-valued random variable~\(Y\),
independent of the reference trajectory, with tail probabilities
\[
 \PP_{\rm ref}(Y\ge r)=\overline G(r),\qquad r\ge1.
\]
The sequence \(\overline G(r)\) is non-increasing and tends to zero, with
\(\overline G(1)<1/2\), so these tails specify a probability law and
\(\PP_{\rm ref}(Y=0)=1-\overline G(1)>0\).  At the chosen value of~\(t\),
the tail-sum identity and \eqref{e.compat-positive-mgf} give
\[
 \log\EE_{\rm ref}e^{tY}<\frac{t}{8C_Q}.
\]
The denominator \(C_Q\) accounts for the bound \(q_N\le C_QN\) on the
number of partition layers.

\begin{lemma}[Shadow moment comparison]
\label{l.shadow-moment-comparison}
If \(Y_1,\ldots,Y_{q_N}\) are independent copies of \(Y\), then, for every
\(t\ge0\) for which \(\EE_{\rm ref}e^{tY}\) is finite,
\[
 \EE_{\rm ref}\exp\left\{t\sum_{q=1}^{q_N}X_q\right\}
 \le \EE_{\rm ref}\exp\left\{t\sum_{q=1}^{q_N}Y_q\right\}
 =\bigl(\EE_{\rm ref}e^{tY}\bigr)^{q_N}.
\]
\end{lemma}

\noindent{\bf Proof.}
For every non-negative integer \(x\),
\[
 e^{tx}\, = \, 1+(e^t-1)\sum_{r\ge1}e^{t(r-1)}\one_{\{x\ge r\}}\, .
\]
Multiply this identity over \(q=1,\ldots,q_N\) and expand.  All terms are
non-negative, so the resulting sums may be compared term by term.  Each
non-constant term is a positive multiple of the indicator of a joint event
\(\{X_q\ge r_q:q\in I\}\).  Apply~\eqref{e.joint-shadow-tail}; the constant
term is unchanged.  The resulting expansion is exactly the corresponding one
for independent \(Y_q\), giving the desired comparison.
\qed

\section{The IVF-to-LLF transition}
\label{ss.ivf-llf-transition}

Sunshine encodes the positive-depth fine conditions: a target outside every
shadow interval is sunny, and the sum of the shadow lengths bounds the number
of targets that are not.  Up to this point, the positive-depth argument has
been carried out under the reference dynamics, avoiding acceptance
conditioning at each successive slide.  The marked trajectory comparison was
designed to postpone this cost: we now apply it once, to the event that too
many selected targets lie in shadow, and thereby transfer the estimate to the
accepted slide dynamics.

\begin{proposition}[Many sunny targets under the accepted law]
\label{p.accepted-sunny-very-fine-heights}
On~\(\Omega_{\rm sat}\),
\[
 \PP_{\rm acc}\left(
   \text{fewer than }3N/4\text{ selected targets are sunny}\right)
 \le e^{-cN}
\]
for some~\(c>0\).
\end{proposition}

\noindent{\bf Proof.}
Let \(D_N\) be the set of target indices covered by the shadow intervals.
Since the layer \(\mathcal Q_q\) covers at most \(X_q\) targets,
\(|D_N|\leq\sum_{q=1}^{q_N}X_q\).  Using \(q_N\le C_QN\), the moment bound
for~\(Y\) and Lemma~\ref{l.shadow-moment-comparison}, Markov's inequality gives
\[
 \PP_{\rm ref}\left(|D_N|\geq\frac N4\right)
 \leq
 e^{-tN/4}\bigl(\EE_{\rm ref}e^{tY}\bigr)^{q_N}
 \leq e^{-tN/8}.
\]
Let~\(\mathsf E_N^+=\{|D_N|\ge N/4\}\).  Apply the marked trajectory
comparison of Lemma~\ref{l.marked-whole-block-change} to~\(\mathsf E_N^+\):
\[
 \PP_{\rm acc}(\mathsf E_N^+)
 =\EE_{\rm ref}[L_n\one_{\mathsf E_N^+}].
\]
The iterated one-step estimate of Corollary~\ref{c.tether-likelihood-moment}
gives \(\EE_{\rm ref}[L_n^{1+\theta}]\le e^{C_{\rm L}\theta n}\).
H\"older's inequality then gives
\[
 \PP_{\rm acc}(\mathsf E_N^+)
 \leq
 \EE_{\rm ref}[L_n^{1+\theta}]^{1/(1+\theta)}
 \PP_{\rm ref}(\mathsf E_N^+)^{\theta/(1+\theta)}
 \leq
 \exp\left\{\frac{\theta}{1+\theta}
   \left(C_{\rm L}n-\frac t8N\right)\right\}.
\]
By~\eqref{e.attempted-block-linear}, \(n\le C_{\rm att}N\).  The margin
\eqref{e.compat-positive-margin} therefore makes the reference decay
\(tN/8\) dominate the likelihood-ratio cost \(C_{\rm L}n\), and the last
display is at most~\(e^{-2\theta N/(1+\theta)}\).
Outside~\(\mathsf E_N^+\), at least \(3N/4\) targets are sunny.
\qed

\begin{proposition}[Sunny targets and depth-zero delivery give LLF]
\label{p.vf-delivered-as-fine}
On~\(\Omega_{\rm sat}\),
\[
 \PP_{\rm acc}\left(
   \text{fewer than }N/2\text{ selected IVF heights are fine at limelight}
 \right)\le Ce^{-cN}.
\]
\end{proposition}

\noindent{\bf Proof.}
If the conclusions of Propositions~\ref{p.zero-depth-delivery} and
\ref{p.accepted-sunny-very-fine-heights} both hold, at most \(N/4\) selected
targets fail \({\bf A}_{\rm f}^-(0)\), and at most \(N/4\) are not sunny.
Thus at least \(N/2\) targets satisfy the depth-zero southside condition and
are sunny.  Fix one such target~\(h_j\).  For \(k\geq1\), consider the unit
band in~\({\bf A}_{\rm f}^-(k)\).  If it lies wholly below height~\(0\), it
contains no cluster vertex, because the origin is the base, and the condition
holds.  Otherwise, by Definition~\ref{d.layer-at-k}, it belongs to a unique
partition layer~\(\mathcal Q_q\).  Since the band lies south of~\(h_j\), this
layer is among those listed, so \(q\le q_N\), and
\(k\in\mathcal D_{q,j}\).  Because~\(h_j\) is sunny,
\(\mathcal Q_q\) does not cast shade on it, and hence
\({\bf A}_{\rm f}^-(k)\) holds.  Thus every southside condition holds.

The northside conditions were already valid at time zero by
Lemma~\ref{l.vf-implies-fine}.  Lemma~\ref{l.higher-ucs-shift} transports the
northside allowance profile rigidly: a lower-slide translation~\(\tau\)
changes ordinary vertex height and the head height by the same amount, so it
preserves the adjusted height \(x\cvl-(\hmax(S)-u)\).  Hence every vertex
count in a fixed adjusted-height interval, and thus every northside allowance
\({\bf A}_{\rm f}^+(k)\), is invariant.  Corollary~\ref{c.breaker-points-persist}
likewise transports the local breaker occupation profile.  Thus the target is
fine at limelight.  A union bound with
Propositions~\ref{p.zero-depth-delivery} and
\ref{p.accepted-sunny-very-fine-heights} gives the stated estimate.
\qed

By Lemma~\ref{l.marked-whole-block-change}, the cluster projection of
the accepted dynamics~\(\PP_{\rm acc}\) has law~\(\Qdyn\).  As such, the
next result is obtained by averaging the conditional estimate of
Proposition~\ref{p.vf-delivered-as-fine} over the input cluster record and
attempted-height schedule.

\begin{theorem}[IVF-to-LLF transition]
\label{t.vf-to-fine-delivery}
Let
\[
 N_{\llf}(u)\, = \, \#\{h\in[0,u]:h\text{ is scheduled and fine at limelight}\}\, .
\]
There are \(\delta_{\rm f},c,C>0\) such that
\[
   \Qdyn\bigl(N_{\llf}(u)<\delta_{\rm f}u\bigr)\, \le \, Ce^{-cu}\, .
\]
\end{theorem}

\noindent{\bf Proof.}
Lemma~\ref{l.sampling-ivf-points}, together with the deterministic
relation \(L\asymp u\) for the now fixed random mesh, gives
\(\Qdyn(\Omega_{\rm sat}^c)\le Ce^{-cu}\).  On~\(\Omega_{\rm sat}\),
\(M=|\mathcal H_{\rm vf}|=|\ivf|\ge\gamma L\ge c_Eu\), and the deterministic well-spaced extraction gives
\(N=|\mathcal W|\ge M/2\ge c_Eu/2\).
Every selected target is a scheduled height, so \(N_{\llf}(u)\) is at
least the number of selected targets that are fine at limelight.
Averaging the estimate of Proposition~\ref{p.vf-delivered-as-fine}
as described above therefore yields
\[
 \Qdyn\left(\Omega_{\rm sat}\cap
   \{N_{\llf}(u)<N/2\}\right)\le Ce^{-c'u},
\]
for some~\(c'>0\), since~\(N\ge c_Eu/2\) on~\(\Omega_{\rm sat}\).
Moreover, \(N/2\ge c_Eu/4\).  Adding the probability of
\(\Omega_{\rm sat}^c\) proves the theorem with
\(\delta_{\rm f}=c_E/4\), after changing~\(c\) and~\(C\).
\qed

\section{Renewal production}
\label{ss.renewal-production-from-fine}

A slide attempted at a scheduled height is active precisely when that height
is fine at limelight.  At each such active slide, SimpleJoin occurs with
conditional probability at least \(p_{\rm ren}>0\), uniformly in the revealed
past; its midvertex is then a renewal point.  This is
Corollary~\ref{c.uniform-renewal-chance-active-slide}.  Thus the linear
fine-height supply yields a positive density of renewal-producing slides.
This inference concludes the chapter, even as the proof of renewal survival
in the next depends on the underlying estimates rather than on the
conclusion itself.

\begin{proposition}[Positive-density renewal production]
\label{p.many-produced-renewals-from-fine}
Let \(N_{\rm prod}(u)\) be the number of scheduled heights in \([0,u]\) at
which an active slide produces a renewal point.  There are
\(a_{\rm prod},c,C>0\) such that
\[
   \Qdyn\bigl(N_{\rm prod}(u)<a_{\rm prod}u\bigr)\leq Ce^{-cu}.
\]
\end{proposition}

\noindent{\bf Proof.}
By Theorem~\ref{t.vf-to-fine-delivery}, outside an event of probability
\(Ce^{-cu}\) there are at least \(\delta_{\rm f}u\) LLF heights.  For
\(u\ge K\), every scheduled height lies in the range \(2r<h<u-2r\), so
Corollary~\ref{c.uniform-renewal-chance-active-slide} applies to each of them.
For the attempted slide at schedule index~\(i\), put
\[
 A_i=\one\{S_{i-1}\in\fine_{H_i}\}.
\]
Set \(\xi_i=0\) on \(\{A_i=0\}\); on \(\{A_i=1\}\), let \(\xi_i\) indicate that
the output state~\(S_i\) realises \(\success\) at height~\(H_i\).  The
height-\(H_i\) DMU decomposition is determined by~\(S_i\), so \(\xi_i\) is
\(\mathcal K_i\)-measurable.  Moreover, Lemma~\ref{l.success} shows that every
step with \(\xi_i=1\) produces a renewal point.  Hence
\[
 N_{\llf}(u)=\sum_{i=1}^L A_i,
 \qquad
 N_{\rm prod}(u)\ge\sum_{i=1}^L \xi_i.
\]
The variable \(A_i\) is \(\mathcal K_{i-1}\)-measurable.  On \(\{A_i=1\}\),
the transition law~\eqref{e.cluster-trajectory-markov} and the uniformity in
\((\downwell,[\upwell])\) in
Corollary~\ref{c.uniform-renewal-chance-active-slide} give
\[
 \Qdyn(\xi_i=1\mid\mathcal K_{i-1})\ge p_{\rm ren}.
\]
Choose \(\lambda>0\) so small that
\[
 e^{\lambda p_{\rm ren}/2}
 \bigl(1-p_{\rm ren}+p_{\rm ren}e^{-\lambda}\bigr)<1.
\]
Conditionally on \(\mathcal K_{i-1}\), the factor
\(\exp\{\lambda(p_{\rm ren}A_i/2-\xi_i)\}\) has mean at most one.  Indeed, it
is identically one when \(A_i=0\); when \(A_i=1\), its conditional mean is
\[
 e^{\lambda p_{\rm ren}/2}
 \bigl(1-p_i+p_i e^{-\lambda}\bigr),
 \qquad
 p_i=\Qdyn(\xi_i=1\mid\mathcal K_{i-1})\ge p_{\rm ren},
\]
which is at most the preceding display because the bracketed expression is
decreasing in~\(p_i\).  Iterating over the attempted schedule therefore yields
\[
 \EE_{\Qdyn}\exp\left\{\lambda\left(
   \frac{p_{\rm ren}}2N_{\llf}(u)-\sum_{i=1}^L \xi_i
 \right)\right\}\le1.
\]
Since \(N_{\rm prod}(u)\ge\sum_i\xi_i\), Markov's inequality now gives,
with \(a_{\rm prod}=p_{\rm ren}\delta_{\rm f}/4\),
\[
 \Qdyn\left(
  N_{\rm prod}(u)<a_{\rm prod}u,\
  N_{\llf}(u)\ge\delta_{\rm f}u
 \right)
 \le \exp\{-\lambda p_{\rm ren}\delta_{\rm f}u/4\}.
\]
Combining this estimate with Theorem~\ref{t.vf-to-fine-delivery}, and
enlarging \(C\) to cover bounded~\(u\), proves the proposition.
\qed

\chapter{Renewal survival and breaker queues}
\label{s.renewal-survival}

Proposition~\ref{p.many-produced-renewals-from-fine} shows that linear
renewal production occurs with overwhelming probability along the trajectory
governed by the slide dynamics~$\Qdyn$.  The terminal cluster of that trajectory has the law
of~$\mc C(0)$ under the static measure~$\Q$ by
Proposition~\ref{p.qdyn-output}.  It is in this terminal cluster that
Theorem~\ref{t.intermediate-renewal}(2) requires a renewal point with
overwhelming probability.  We cannot therefore observe a renewal during the
dynamics, stop at that moment and declare victory.  We must show that under~$\Qdyn$ at least
one point produced in this way retains its renewal status when the schedule
ends.  This is the renewal-survival problem treated in the present chapter.

Suppose, then, that as the schedule advances north, an active slide produces a
renewal point.  Is this point ever undone?
A dormant slide changes nothing.  Each active slide, however, presents a
threat: its newly formed midpiece may travel south and add another cluster
vertex to the renewal band, destroying the renewal property.  Yet the same
active slide also supplies protection against later threats from farther
north.  Indeed, the fine condition which permits the
slide includes a triple-row breaker window immediately to its north.  The midpiece formed
in the slide typically reaches only the lower row, leaving the middle and upper
rows untouched and standing ready to resist southgoing attacks generated at
still later slide heights.

This suggests treating the middle and upper rows of each breaker window as
separate units and monitoring how many such units remain untested north of the
renewal point currently being followed and south of the advancing slide
schedule.  Active slide by active slide, a southgoing midpiece may remove some
of the standing units, while the new slide may place further units beyond
them.  The count evolves as a one-dimensional breaker queue.
If this queue is exhausted, we abandon the renewal then being followed: its
survival has not been established.  At the next renewal produced together
with two unused breaker units
farther north, we begin a new queue and follow that renewal instead.  Thus the
tracked renewal changes when the queue is restarted.  If the construction
reaches the terminal slide with a current queue still positive, the renewal
attached to its last restart has survived, and so the terminal cluster has a
renewal point.

In this broad sense survival is easier than production.  Renewal production
required the simultaneous control, up to limelight, of the several southside
cardinality conditions which define a fine height, and hence the monitoring of
downcluster accrual in many layers.  The survival argument keeps track of one
renewal and is essentially one-dimensional in nature, following one queue at
a time.

There is nevertheless one respect in which survival is less favourable.  The
wavebreakers used in fine delivery faced midpieces advancing from the south.
Before passage, the rising schedule could cut through an old midpiece, and the
cut-it-and-lose-it mechanism of Section~\ref{ss.protection-height} restored
the thin breaker record for the next test.  Here the attack comes from the
north.  If a southgoing midpiece passes a breaker row below it, no later cut
is compelled to remove that midpiece.  The breaker units in the survival
argument must therefore be treated as single-use: once tested, a unit is
removed from the queue whether it is passed or arrests the attack.  The task is to arrange that
the resulting queue nonetheless has positive step-by-step drift.

The three-row breaker construction was given in
Definition~\ref{d.breaker-window-geometry}.  The lower row confronts the
midpiece formed at the active slide itself.  Corollary~\ref{c.accepted-clean-future-breakers}
shows that the midpiece does not reach the middle or upper row with conditional
probability at least \(1-\varepsilon_{\rm row}\), where
\(\varepsilon_{\rm row}=C_{\rm cl}q_{\rm row}^{\delta_{\rm cl}}<1/8\) by
\eqref{e.compat-row-errors}.  On this event the middle and upper rows remain untested.  Both become
breaker units; once up and running, these two units are equally valid
defences against an attack from the north.  Thus the conditional mean number
of units made available at an active slide is at least~$2(1-\varepsilon_{\rm row})$.

For the loss, consider the next southgoing midpiece.  The standing
units are ordered from north to south.  The northernmost unit is removed in
every case, even if the midpiece stops before reaching it; its supporting row
is the first barrier protecting the units farther south.  A second unit is
removed only if that row is passed, and each further removal requires passage
through another row.  Proposition~\ref{p.successive-southgoing-breakers} bounds the
fixed-shift probability of passing several such rows in one southgoing
exploration.  With \(\delta_{\rm A}\) and \(C_{\rm A}\) from
\eqref{e.accepted-rare-transfer}, put
\(\qbarrow=q_{\rm row}^{\delta_{\rm A}}\).  Lemma~\ref{l.geometric-breaker-loss}
transfers the joint event that several successive units are passed to the
accepted law and gives the tail bound \(C_{\rm A}\qbarrow^{t-1}\) for the
queue-loss variable to be at least~\(t\).  By \eqref{e.compat-row-errors}, its
mean is at most
\[
 m_{\rm q}=1+\frac{C_{\rm A}\qbarrow}{1-\qbarrow}<\frac54.
\]
The mean number of units supplied is at least
\(2(1-\varepsilon_{\rm row})>7/4\), so the queue has mean drift greater than
\(1/2\).

The proof implements the restarted queue by a greedy rule.  If the current
queue is abandoned, we continue along the active slide schedule until a later
active slide realises the single-slide SimpleJoin construction of
Chapter~\ref{s.single-slide-positive-density}.  That slide produces the next
renewal to be followed, with two untouched breaker units already available to
start its protection queue.  Positive drift gives a uniform non-abandonment
estimate for each queue.  The waiting time for this next
renewal-production event and the duration of a failed attempt have uniform
exponential tails.  Hence, if there are linearly many active slides,
either a queue survives or, except with exponentially small probability,
linearly many successive attempts occur.  The conditional failure bounds
for these attempts then multiply by the tower property.  With overwhelming
probability, the procedure therefore reaches the terminal slide with a
positive current queue, which is precisely the asserted survival of some
produced renewal.

As the discussion indicates, the natural clock is active-slide time, which we
first introduce formally.  We then state precisely when the middle and upper
rows remain untouched and may enter the queue.  Next, we define the queue
and prove its geometric loss estimate.  We then establish queue persistence
and use a deterministic lemma based on the DMU decomposition of
Chapter~\ref{s.dmu} to verify that a queue which survives to the end leaves
its renewal as the unique cluster vertex in its renewal band.  Finally, the
greedy restart procedure shows that some such queue survives with overwhelming
probability.  This completes the passage from renewal production along~$\Qdyn$
to the terminal~$\Q$-statement of Theorem~\ref{t.intermediate-renewal}(2).

In this chapter, unmarked $\PP$ and $\E$ denote probability and
expectation under~$\Qdyn$.

\section{Active-slide time}
\label{ss.slide-index-time-queue}

Recall from Section~\ref{ss.candidate-heights} that
\(\mathcal K_i=\sigma(S_0,\overline H,S_1,\ldots,S_i)\) is the information
available after the first~$i$ attempted slides.  The $i$\textsuperscript{th}
attempt is active exactly when $S_{i-1}\in\fine_{H_i}$.  Let
\[
 N=\sum_{i=1}^{L}\one\{S_{i-1}\in\fine_{H_i}\}
\]
be the number of active slides, denoted \(N_{\llf}(u)\) in
Theorem~\ref{t.vf-to-fine-delivery}.
For \(1\le j\le L\), let
\[
 \kappa_j
 =\min\left(\left\{1\le i\le L:
       \sum_{a=1}^{i}\one\{S_{a-1}\in\fine_{H_a}\}=j\right\}
       \cup\{L+1\}\right).
\]
Thus \(\kappa_j\le L\) exactly when \(j\le N\), and, for \(i\le L\),
\(\{\kappa_j=i\}\in\mathcal K_{i-1}\).  Set
\[
 \mathcal G_0=\mathcal K_0,
 \qquad
 \mathcal G_j=\mathcal K_{\kappa_j\wedge L},\quad 1\le j\le L,
 \qquad
 \mathcal G_{L+1}=\mathcal K_L.
\]
The stopped times \(\kappa_j\wedge L\) are non-decreasing, so
\((\mathcal G_j)_{j=0}^{L+1}\) is a filtration.  For \(j>N\), we have
\(\mathcal G_j=\mathcal K_L\).  For \(1\le j\le L\), let
\(\mathcal G_{j-}\) record the information available before the \(j\)\textsuperscript{th} active
slide is sampled, and let \(\widehat S_{j-}\) and \(\widehat S_j\) denote
the cluster immediately before and after that slide:
\[
 \mathcal G_{j-}=\mathcal K_{\kappa_j-1},\qquad
 \widehat S_{j-}=S_{\kappa_j-1},\qquad
 \widehat S_j=S_{\kappa_j\wedge L}.
\]
The random time \(\kappa_j-1\) is a stopping time and equals \(L\) when
\(j>N\).  Thus the definitions remain valid even when the
$j$\textsuperscript{th} active slide does not occur.  The schedule is known
at time zero, and dormant steps leave the cluster unchanged.  Thus
\(\mathcal G_{j-}=\mathcal G_{j-1}\): intervening dormant attempts reveal no
new cluster information, and the existence of the next active slide is
determined by the preceding active state.

On \(\{j\le N\}\), write \(h_j=H_{\kappa_j}\).  Then
\(\widehat S_{j-}\in\fine_{h_j}\).  In particular, the input cluster and
attempted height are known at this pre-slide time, while the outcome of
the slide remains random.  Its conditional law is
\begin{equation}
\label{e.active-slide-stopped-kernel}
 \operatorname{Law}(\widehat S_j\mid\mathcal G_{j-})
 =\overline\slide_{h_j}(\widehat S_{j-},\cdot)
 \qquad\text{on }\{j\le N\}.
\end{equation}
Indeed, on \(\{\kappa_j=i\}\in\mathcal K_{i-1}\), this is the
one-step identity~\eqref{e.cluster-trajectory-markov} at attempted index~\(i\).
Since the events \(\{\kappa_j=i\}\), \(i\le L\), partition \(\{j\le N\}\),
these one-step identities combine to give the display.  It permits the
one-step estimates to be iterated along the active-slide schedule up to its
random terminal index~\(N\).

\section{Two unused breaker units at an active slide}
\label{ss.usable-breaker-birth}

At an active index~$j\le N$, the breaker window at~$h_j$ has the three
rows of Definition~\ref{d.breaker-window-geometry}.  Its lower row
$\mathscr B_{h_j}^0$ confronts the midpiece formed at the slide itself.  The
middle and upper rows are the reserved rows of
Section~\ref{ss.breaker-stands-postpassage}, intended for later attacks.

Because the adjusted cut moves with the translated upcluster, the fine
input bounds on its occupation of the breaker rows remain unchanged
after the slide.

Recall that, for $\rho\in\{0,1,2\}$, $\mathscr T_h^\rho$ denotes the union graph
of the test layers which constitute $\mathscr B_h^\rho$.  For
$\rho\in\{1,2\}$, write $\mathfrak b_{j,\rho}$ for the \emph{breaker unit}
labelled by the active index~$j$ and the row index~$\rho$.  Its supporting
row is $\mathscr B_{h_j}^\rho$, and its test graph is
$\mathscr T_{h_j}^\rho$.  Thus a
breaker unit is a labelled, single-use candidate for a protection
queue.  For \(1\le j\le L\), define
\[
 \mathsf U_j\, = \,
 \Bigl\{j\le N,\ V\bigl(\midpiece_{h_j}(\widehat S_j)\bigr)
       \cap V\bigl(\mathscr T_{h_j}^1\cup\mathscr T_{h_j}^2\bigr)
\, = \, \emptyset\Bigr\}\, ,
 \qquad \beta_j\, = \, \one_{\mathsf U_j}\, .
\]
Thus $\beta_j=0$ for $j>N$, and $\beta_j=1$ precisely when the midpiece
formed at the active slide avoids the test graphs $\mathscr T_{h_j}^1$ and $\mathscr T_{h_j}^2$
of the middle and upper breaker rows.  The event $\mathsf U_j$
is $\mathcal G_j$-measurable.  Equation~\eqref{e.active-slide-stopped-kernel},
Corollary~\ref{c.accepted-clean-future-breakers} and
\eqref{e.compat-row-errors} give, on \(\{j\le N\}\),
\begin{equation}
\label{e.usable-birth}
   \PP(\beta_j\, = \, 0\mid\mathcal G_{j-})
\, \le \, \varepsilon_{\rm row}
:=C_{\rm cl}q_{\rm row}^{\delta_{\rm cl}}<\frac18\, .
\end{equation}

On $\mathsf U_j$, the current midpiece leaves both reserved rows
untouched.  The retained downcluster lies below their test graphs, so the
transferred upcluster bounds and Lemma~\ref{l.breaker-layer-vertex-range}
give at least $i_{\rm br}/2$ thin test layers in each row.

Let $\mathsf R_j$ be the event that $j\le N$ and $\success$ occurs
at the active slide at~$h_j$.  On this event the midpiece is the two-edge
bridge from $b=\head(\downcluster)$ to $b+2e_d$, the latter vertex lying
in the bed of the output upcluster.  Since this vertex has height at most
$\hadj(\widehat S_j,h_j)$, the whole bridge lies at or below the adjusted
cut.  By \eqref{e.breaker-profile-band} and
Lemma~\ref{l.breaker-layer-vertex-range}, every vertex of either reserved
test graph has height strictly greater than
$\hadj(\widehat S_j,h_j)+b_0-r>\hadj(\widehat S_j,h_j)$.
Thus $\mathsf R_j\subseteq\mathsf U_j$.  Put
$J_j=\one_{\mathsf R_j}$ for $1\le j\le L$; in particular, $J_j=0$ when
$j>N$.  This is $\mathcal G_j$-measurable.  For all sufficiently large~$u$ we also
have $u\ge K$.  By the endpoint choice in Section~\ref{ss.random-mesh},
every scheduled height satisfies $2r<H_i<u-2r$, so
\eqref{e.active-slide-stopped-kernel} and
Corollary~\ref{c.uniform-renewal-chance-active-slide} give,
on $\{j\le N\}$,
\begin{equation}
\label{e.queue-start-probability}
 \PP(J_j=1\mid\mathcal G_{j-})
 =\PP(\mathsf R_j\mid\mathcal G_{j-})\ge p_{\rm ren}>0.
\end{equation}
On $\{J_j=1\}$, let $\varrho_j$ denote the midvertex of the
SimpleJoin bridge, which is a renewal point by Lemma~\ref{l.success}.
The pair $\mathfrak b_{j,2},\mathfrak b_{j,1}$ is then available to begin a
protection queue.

\section{The protection queue and its one-step estimate}
\label{ss.breaker-loss}

Fix $1\le\sigma\le L$ and work on $\{J_\sigma=1\}$.
Consider the renewal $\varrho_\sigma$ produced at~$h_\sigma$.
A standing breaker unit will be recorded by its label
$\mathfrak b_{a,\rho}$.  We say that it is \emph{passed} at a later
active slide when the output midpiece has a vertex on the southern side
of the southernmost test layer of its supporting row
$\mathscr B_{h_a}^\rho$.  On acceptance, this agrees with finite passage
in the southgoing fixed-shift exploration, by
Proposition~\ref{p.successive-southgoing-breakers}.
The breaker windows belonging to distinct attempted heights are disjoint by
Lemma~\ref{l.scheduled-breaker-separation}; the labels therefore inherit an unambiguous
north-to-south order.

\begin{definition}[The protection queue]
\label{d.active-slide-queue}
For a birth index~$\sigma$ with $J_\sigma=1$, begin the
\hfff{protection-queue}{protection queue} with the ordered list
\[
 \mathbf U_0^{(\sigma)}
\, = \, \bigl(\mathfrak b_{\sigma,2},\mathfrak b_{\sigma,1}\bigr)\, ,
 \qquad Q_0^{(\sigma)}\, = \, \bigl|\mathbf U_0^{(\sigma)}\bigr|\, = \, 2\, .
\]
We define the subsequent lists and the abandonment time $\tau_\sigma$
recursively.  At step $1\le n\le N-\sigma$, on $\{\tau_\sigma\ge n\}$,
immediately before the active slide at~$h_{\sigma+n}$, write
\(\mathbf U_{n-1}^{(\sigma)}=(\mathfrak b_1,\ldots,\mathfrak b_q)\), where
\(q=Q_{n-1}^{(\sigma)}\) and each $\mathfrak b_i$ is the label of the
$i$\textsuperscript{th} standing breaker unit.  The list is ordered from
north to south.  Let $K_n^{(\sigma)}$ be the
number of consecutive supporting rows passed by
$\midpiece_{h_{\sigma+n}}(\widehat S_{\sigma+n})$, beginning with that of
$\mathfrak b_1$, and put $D_n^{(\sigma)}=K_n^{(\sigma)}+1$.

If $D_n^{(\sigma)}>q$, abandon the queue.  Otherwise remove
$\mathfrak b_1,\ldots,\mathfrak b_{D_n^{(\sigma)}}$ and then, when $\beta_{\sigma+n}=1$, insert
the newly available pair at the north end.  Abandon the queue at this step
also if the updated list is empty.  On $\{D_n^{(\sigma)}\le q\}$, the updated list is
\[
 \mathbf U_n^{(\sigma)}\, = \,
 \begin{cases}
  \bigl(\mathfrak b_{\sigma+n,2},\mathfrak b_{\sigma+n,1}\, ,
        \mathfrak b_{D_n^{(\sigma)}+1},\ldots,\mathfrak b_q\bigr)\, ,
       &\beta_{\sigma+n}\, = \, 1\, ,\\[1mm]
  \bigl(\mathfrak b_{D_n^{(\sigma)}+1},\ldots,\mathfrak b_q\bigr)\, ,
       &\beta_{\sigma+n}\, = \, 0\, ,
 \end{cases}
\]
If $D_n^{(\sigma)}=q$, no units from the previous list remain.
Set $Q_n^{(\sigma)}=|\mathbf U_n^{(\sigma)}|$.  Thus
\[
 Q_n^{(\sigma)}
\, = \, Q_{n-1}^{(\sigma)}-D_n^{(\sigma)}+2\beta_{\sigma+n}\, .
\]
The abandonment time is
\[
 \tau_\sigma\, = \, \inf\Bigl\{1\, \le \, n\, \le \, N-\sigma\, :
 D_n^{(\sigma)}>Q_{n-1}^{(\sigma)}
 \ \hbox{or}\
 Q_{n-1}^{(\sigma)}-D_n^{(\sigma)}+2\beta_{\sigma+n}\, = \, 0
 \Bigr\}\, ,
\]
where $\inf\emptyset=\infty$.  Under the first condition, every
standing row has been passed, and no update is performed.  When the first
condition fails, the second means that the last standing row remains
unpassed, but conservative removal leaves no unit for a later defence.
We count this as abandonment even if the schedule ends at that slide.

The first unpassed unit, when present, occupies position
$D_n^{(\sigma)}$.  It is removed along with any passed units, since the
midpiece may have entered its test graph without passing through the row,
compromising the thin record needed for a later test.  In particular,
every attack removes at least one unit, even when its midpiece dies before
reaching the first supporting row.

The queue survives to the terminal state exactly when $\tau_\sigma=\infty$.
When a birth index is fixed, the superscript $(\sigma)$ is omitted.
\end{definition}

\begin{lemma}[The standing breaker list]
\label{l.standing-breaker-list}
Fix $1\le\sigma\le L$ and work on $\{J_\sigma=1\}$.
For $0\le n\le N-\sigma$, on $\{\tau_\sigma>n\}$, write
\(\mathbf U_n^{(\sigma)}=(\mathfrak b_{a_1,\rho_1},\ldots,
\mathfrak b_{a_q,\rho_q})\).  Then:
\begin{enumerate}
\item\label{it.standing-breaker-order}
We have $q=Q_n^{(\sigma)}$.  The supporting rows have pairwise disjoint
edge sets and are listed in north-to-south order.  No midpiece
formed at an active index in $\{a_s,\ldots,\sigma+n\}$ has entered
$\mathscr T_{h_{a_s}}^{\rho_s}$, for any $1\le s\le q$.

\item\label{it.standing-breaker-separation}
Every supporting row in the list lies strictly north of
$\vertbdry_{\varrho_\sigma\cvl}$ and strictly south of every later attempted
height.  Consequently a path descending from a later cut
to $\vertbdry_{\varrho_\sigma\cvl}$ meets the test layers supporting the
listed units in their listed order.

\item\label{it.standing-breaker-persistence}
If $n<N-\sigma$, then at the next active slide $h_{\sigma+n+1}$ the
hypotheses of Proposition~\ref{p.southgoing-breaker-persistence} hold in the
current fixed-shift experiment for each listed unit
$\mathfrak b_{a_s,\rho_s}$, with $h=h_{a_s}$ and $\rho=\rho_s$.
Thus its conclusions, including the conditional Bernoulli product law at a
finite pre-entry stop, apply to each such unit.
\end{enumerate}
\end{lemma}

\noindent{\bf Proof.} At $n=0$, the event $J_\sigma=1$ includes $\mathsf U_\sigma$, so the
midpiece formed at~$h_\sigma$ enters neither supporting row.  Their
north-to-south order is $\mathfrak b_{\sigma,2}$ followed by
$\mathfrak b_{\sigma,1}$.  In the SimpleJoin construction,
$\varrho_\sigma=b+e_d$ is the midvertex of a bridge whose upper endpoint
$b+2e_d$ lies in the bed of the output upcluster.  Thus
$\varrho_\sigma\cvl\le\hadj(\widehat S_\sigma,h_\sigma)$, while
\eqref{e.breaker-profile-band} and
Lemma~\ref{l.breaker-layer-vertex-range} place every vertex of either
reserved row strictly above that adjusted height.  Hence both supporting
rows lie strictly north of the renewal level.  Lemma~\ref{l.scheduled-breaker-separation}
places every later attempted
height strictly north of them.

Suppose that (1) and (2) hold for
$\mathbf U_{n-1}^{(\sigma)}=(\mathfrak b_1,\ldots,\mathfrak b_q)$.
The current midpiece passes the rows of
$\mathfrak b_1,\ldots,\mathfrak b_{K_n}$ in order.  The next unit is removed by the conservative
rule, whether the midpiece stops before reaching it or enters it without
passing it.  Every unit in the retained tail
$\mathfrak b_{D_n+1},\ldots,\mathfrak b_q$ lies farther south.
If the current midpiece entered one of these lower rows, a component
meeting it would have an attachment to the retained clusters above all
the standing rows, by Lemma~\ref{l.midpiece}(3) and
\eqref{e.southgoing-row-margin}.  It would therefore cross the southern
side of the first unpassed row, contradicting the definition of~$K_n$.
Thus no retained row has been entered, and (1) and (2) persist.

When $\beta_{\sigma+n}=1$, the current midpiece avoids the two new test
graphs, and the pair
$\mathfrak b_{\sigma+n,2},\mathfrak b_{\sigma+n,1}$ is inserted.
Lemma~\ref{l.scheduled-breaker-separation} places this pair north of
every older unit and south of every later attempted height.  It also lies
north of the renewal level: the random-mesh
gap gives $h_{\sigma+n}>h_\sigma+4r$, whereas
$\varrho_\sigma\cvl\le h_\sigma+r$.  A path descending from north of a
test layer to south of it must enter its vertex set through the northern side,
by Lemma~\ref{l.slab-entry-faces}.  Applying this observation to the
height-ordered layers gives the path statement in (2).  This proves
(1) and (2) by induction.

For each unit in $\mathbf U_n^{(\sigma)}$, its birth event $\beta_{a_s}=1$
is the event $\mathsf{Clean}_{h_{a_s}}$ of Chapter~\ref{s.wavebreakers}.
The no-entry assertion in (1) covers all intervening midpieces.
Proposition~\ref{p.southgoing-breaker-persistence} therefore applies at the
next active slide, whenever it exists, proving (3).
\qed

For the queue estimates, recall $\qbarrow=q_{\rm row}^{\delta_{\rm A}}$
and put
\[
 m_{\rm q}:=1+\frac{C_{\rm A}\qbarrow}{1-\qbarrow}<\frac54,
 \qquad
 \delta_{\rm q}:=2(1-\varepsilon_{\rm row})-m_{\rm q}>\frac12\, .
\]
The inequalities follow from~\eqref{e.compat-row-errors}.
The next lemma uses the successive southgoing estimate to control
loss under the accepted slide law, and derives the drift and moment bounds.
Its joint exponential-moment estimate allows for dependence between supply
and loss.

\begin{lemma}[Geometric queue loss and positive drift]
\label{l.geometric-breaker-loss}
Uniformly in the birth index~$\sigma$, the step~$n$, the number of standing
units and the pre-slide history, for every $1\leq n\leq N-\sigma$, on
$\{\tau_\sigma\ge n\}$,
\[
   \PP(D_n^{(\sigma)}\, \geq \, t \mid \mathcal G_{(\sigma+n)-})
\, \leq \, C_{\rm A}\qbarrow^{\,t-1},\qquad t\, \ge \, 1\, .
\]
With $X_n^{(\sigma)}=2\beta_{\sigma+n}-D_n^{(\sigma)}$,
\[
 \E(X_n^{(\sigma)}\mid\mathcal G_{(\sigma+n)-})
\, \ge \, \delta_{\rm q}>0\, .
\]
Moreover, for $0<\theta_{\rm q}<-\log \qbarrow$, put
\[
 T_{\rm q}(\theta_{\rm q})\, = \,
  (e^{\theta_{\rm q}}-1)\frac{C_{\rm A}\qbarrow e^{\theta_{\rm q}}}
                         {1-\qbarrow e^{\theta_{\rm q}}}\, ,
 \qquad M_{\rm q}(\theta_{\rm q})\, = \, e^{\theta_{\rm q}}+T_{\rm q}(\theta_{\rm q})\, .
\]
Then
\begin{equation}
\label{e.queue-laplace-bound}
 \E\bigl(e^{-\theta_{\rm q} X_n^{(\sigma)}}
       \mid\mathcal G_{(\sigma+n)-}\bigr)
\, \le \, \Phi(\theta_{\rm q}):=
 e^{-2\theta_{\rm q}}M_{\rm q}(\theta_{\rm q})
 +(1-e^{-2\theta_{\rm q}})
   \bigl(\varepsilon_{\rm row}e^{\theta_{\rm q}}+T_{\rm q}(\theta_{\rm q})\bigr)\, .
\end{equation}
In particular $\Phi(0)=1$ and $\Phi'(0)=-\delta_{\rm q}<0$, so there is
$\theta_{{\rm q},0}>0$ for which $\Phi(\theta_{{\rm q},0})<1$.
\end{lemma}

\noindent{\bf Proof.} By
Lemma~\ref{l.standing-breaker-list}(\ref{it.standing-breaker-separation}),
passage of a standing row entails passage of every standing row farther north.
Thus, for $2\le t\le Q_{n-1}^{(\sigma)}+1$,
$\{D_n^{(\sigma)}\ge t\}$ is the event that the current midpiece
passes the first $t-1$ rows of $\mathbf U_{n-1}^{(\sigma)}$.
By Lemma~\ref{l.standing-breaker-list}, these rows have remained unentered
since they were reserved and are admissible for the successive southgoing
test of Proposition~\ref{p.successive-southgoing-breakers}.
Fix the pre-slide history and write $E_{t-1}$ for the event that the
accepted output midpiece passes these rows.  For each proposed head~$v$,
let $A_{t-1}^v$ be their joint finite-passage event under the ordinary
fixed-shift law.  Proposition~\ref{p.successive-southgoing-breakers} gives
\[
 \Rlaw^v_{\uclust,\dclust}(A_{t-1}^v\mid\cF_0)
 \le q_{\rm row}^{\,t-1}.
\]
On $\Accept_v$, this is precisely the passage event~$E_{t-1}$,
by the final assertion of that proposition.

Write
\[
 Z=\sum_v\Rlaw^v_{\uclust,\dclust}(\Accept_v),
 \qquad
 \nu(v)=Z^{-1}\Rlaw^v_{\uclust,\dclust}(\Accept_v).
\]
The accepted cluster kernel of Lemma~\ref{l.probform} therefore gives
\[
 \begin{aligned}
 \overline\slide_H(S,E_{t-1})
 &=\frac1Z\sum_v
   \Rlaw^v_{\uclust,\dclust}(A_{t-1}^v\cap\Accept_v)\\
 &\le\sum_v\min\{q_{\rm row}^{\,t-1}/Z,\nu(v)\}
 \le C_{\rm A}q_{\rm row}^{\delta_{\rm A}(t-1)}
 =C_{\rm A}\qbarrow^{\,t-1}.
 \end{aligned}
\]
The final inequality is the same small-atom summation estimate used in
the proof of Lemma~\ref{l.accepted-rare-event-transfer}, with its
constants $C_{\rm A}$ and $\delta_{\rm A}$.  In particular its prefactor
is used only once for the joint passage event.  The active-slide
pre-history specifies the rows used to define $E_{t-1}$.  Thus
\eqref{e.active-slide-stopped-kernel} gives
\[
 \PP(D_n^{(\sigma)}\ge t\mid\mathcal G_{(\sigma+n)-})
 \le C_{\rm A}\qbarrow^{\,t-1}.
\]
For $t=1$, the bound follows from $C_{\rm A}\ge1$; for
$t>Q_{n-1}^{(\sigma)}+1$ the event is empty.  This proves the tail
estimate in all cases.

For an $\N$-valued random variable~$D$, a $\sigma$-field~$\mathcal H$, and
$z>1$, the conditional tail-sum identity is
\[
 \E(z^D\mid\mathcal H)
\, = \, z+(z-1)\sum_{t\ge2}z^{t-1}\PP(D\ge t\mid\mathcal H)\, .
\]
For $D=D_n^{(\sigma)}$ and
$\mathcal H=\mathcal G_{(\sigma+n)-}$, applying the tail bound to the
ordinary tail sum and to the displayed identity gives, respectively,
\[
 \E(D_n^{(\sigma)}\mid\mathcal G_{(\sigma+n)-})\, \le \, m_{\rm q}
\quad\text{and}\quad
 \E(e^{\theta_{\rm q} D_n^{(\sigma)}}\mid\mathcal G_{(\sigma+n)-})
\, \le \, M_{\rm q}(\theta_{\rm q})\, .
\]
Together with~\eqref{e.usable-birth}, the first inequality yields
$\E(X_n^{(\sigma)}\mid\mathcal G_{(\sigma+n)-})\ge\delta_{\rm q}$.

For the joint exponential moment, use
\[
 e^{\theta_{\rm q} D_n^{(\sigma)}-2\theta_{\rm q} \beta_{\sigma+n}}
\, = \, e^{-2\theta_{\rm q}}e^{\theta_{\rm q} D_n^{(\sigma)}}
 +(1-e^{-2\theta_{\rm q}})e^{\theta_{\rm q} D_n^{(\sigma)}}
       \one_{\{\beta_{\sigma+n}\, = \, 0\}}\, .
\]
The same tail sum, now intersected with $\{\beta_{\sigma+n}=0\}$, and
$\PP(\beta_{\sigma+n}=0\mid\mathcal G_{(\sigma+n)-})
\le\varepsilon_{\rm row}$ give
\[
 \E(e^{\theta_{\rm q} D_n^{(\sigma)}}\one_{\{\beta_{\sigma+n}\, = \, 0\}}\mid
       \mathcal G_{(\sigma+n)-})
\, \le \, \varepsilon_{\rm row}e^{\theta_{\rm q}}+T_{\rm q}(\theta_{\rm q})\, .
\]
This proves~\eqref{e.queue-laplace-bound}.  Finally,
$M_{\rm q}'(0)=m_{\rm q}$, so
$\Phi'(0)=m_{\rm q}-2(1-\varepsilon_{\rm row})=-\delta_{\rm q}$.
\qed

\section{Queue persistence and renewal preservation}
\label{ss.queue-persistence}

Fix $1\le\sigma\le L$ and work on $\{J_\sigma=1\}$, suppressing the
superscript $(\sigma)$.  On $\{\tau_\sigma>n,\ \sigma+n\le N\}$,
\[
 Q_n=2+\sum_{i=1}^n X_i,\qquad X_i=2\beta_{\sigma+i}-D_i.
\]
Loss precedes replenishment, so a positive updated count does not by itself
preclude abandonment.  The geometric-loss and Laplace bounds of
Lemma~\ref{l.geometric-breaker-loss} control abandonment from a standing
reserve; a short run of favourable slides then gives a uniform chance of
survival.

\begin{lemma}[Uniform queue-survival probability]
\label{l.non-ruin-probability}
There is $s>0$, depending only on $C_{\rm A}$, $\qbarrow$ and
$\varepsilon_{\rm row}$, such that, for every renewal born at an active
slide~$\sigma$ with $J_\sigma=1$,
\[
   \PP\bigl(\tau_\sigma=\infty\mid\mathcal G_\sigma\bigr)\, \ge \, s\, .
\]
\end{lemma}

\noindent{\bf Proof.} Choose $\theta_{\rm q}\in(0,-\log\qbarrow)$ so that
$\chi_{\rm q}:=\Phi(\theta_{\rm q})<1$, and put
$C_{\rm q}:=C_{\rm A}/[\qbarrow(1-\chi_{\rm q})]$.
We first bound the probability of abandonment from a standing reserve.
For $0\le m<n$, on $\{\tau_\sigma>m,\ \sigma+m\le N\}$, we claim that
\begin{equation}
\label{e.queue-first-abandonment}
 \PP(\tau_\sigma=n\mid\mathcal G_{\sigma+m})
 \, \le \, \frac{C_{\rm A}}{\qbarrow}\,
 \chi_{\rm q}^{\,n-m-1}e^{-\theta_{\rm q}Q_m}\, .
\end{equation}
Indeed, if the queue is still alive before an existing step~$n$
with $q\ge1$ standing units, abandonment at that step requires
$D_n\ge q$.  The geometric-loss bound and
$\theta_{\rm q}<-\log\qbarrow$ therefore give
\begin{equation}
\label{e.queue-loss-at-reserve}
 \PP(D_n\ge q\mid\mathcal G_{(\sigma+n)-})
 \, \le \, C_{\rm A}\qbarrow^{q-1}
 \, \le \, \frac{C_{\rm A}}{\qbarrow}e^{-\theta_{\rm q}q}\, .
\end{equation}
At a subsequent step~$j\ge1$, on $\{\tau_\sigma>j-1,\ \sigma+j\le N\}$, a non-abandoned update has
$Q_j=Q_{j-1}+X_j$.  Thus \eqref{e.queue-laplace-bound} gives
\[
 \E\!\left[e^{-\theta_{\rm q}(Q_{j-1}+X_j)}
       \one_{\{\tau_\sigma>j\}}
       \mid\mathcal G_{(\sigma+j)-}\right]
 \le \chi_{\rm q}e^{-\theta_{\rm q}Q_{j-1}}\, .
\]
As explained in Section~\ref{ss.slide-index-time-queue}, the information
after the current active slide determines whether another will occur;
if none does, no further abandonment is possible.  Consequently, for $j\ge m$,
\[
 \E\bigl[e^{-\theta_{\rm q}Q_j}
       \one_{\{\tau_\sigma>j,\ \sigma+j\le N\}}
       \mid\mathcal G_{\sigma+m}\bigr]
 \, \le \, \chi_{\rm q}^{\,j-m}e^{-\theta_{\rm q}Q_m}
 \qquad\text{on }\{\tau_\sigma>m,\ \sigma+m\le N\}.
\]
Combining~\eqref{e.queue-loss-at-reserve} with the preceding display
at $j=n-1$ proves~\eqref{e.queue-first-abandonment}.  Summing over $n>m$
gives
\begin{equation}
\label{e.queue-ruin-from-height}
 \PP(m<\tau_\sigma<\infty\mid\mathcal G_{\sigma+m})
 \, \le \, C_{\rm q}e^{-\theta_{\rm q}Q_m}
 \qquad\text{on }\{\tau_\sigma>m,\ \sigma+m\le N\}.
\end{equation}

Choose an integer~$n_{\rm fav}$ such that
$C_{\rm q}e^{-\theta_{\rm q}(2+n_{\rm fav})}\le\tfrac12$.
At each subsequent active slide while the queue is alive,
\[
 \PP(\beta_{\sigma+n}\, = \, 1,D_n\, = \, 1\mid\mathcal G_{(\sigma+n)-})
\, \ge \, 1-\varepsilon_{\rm row}-C_{\rm A}\qbarrow\, =: \, a_0>0\, .
\]
Indeed, the complementary event is contained in
$\{\beta_{\sigma+n}=0\}\cup\{D_n\ge2\}$.  By
\eqref{e.compat-row-errors}, \(\varepsilon_{\rm row}<1/8\) and
\(C_{\rm A}\qbarrow<1/4\), so \(a_0>5/8\).  On the displayed event
the reserve increases by one.
Starting at the birth history, expose up to $n_{\rm fav}$ subsequent active
slides, stopping if the schedule ends first.  At each slide that occurs,
the conditional probability of the favourable event is at least~$a_0$;
if the schedule ends early, no further attack occurs.  Since the existence
of the next active slide is determined by the current history, successive
conditioning gives probability at least~$a_0^{n_{\rm fav}}$ that every
slide encountered in this trial is favourable.  When this happens, either
fewer than $n_{\rm fav}$ slides occur and the queue survives to the terminal
state, or $n_{\rm fav}$ occur and the reserve reaches~$2+n_{\rm fav}$.
In the latter case, \eqref{e.queue-ruin-from-height}, applied with
$m=n_{\rm fav}$, yields an upper bound of one half on the conditional
probability of later abandonment.
Hence $s=\tfrac12a_0^{n_{\rm fav}}$ is valid in both cases.
\qed

\begin{lemma}[Exponential tail of a failed attempt]
\label{l.failed-attempt-tail}
For the abandonment time $\tau_\sigma$ of
Definition~\ref{d.active-slide-queue}, conditional on the birth history
$\mathcal G_\sigma$,
\[
   \PP\bigl(k\, \le \, \tau_\sigma<\infty\mid\mathcal G_\sigma\bigr)
\, \le \, Ce^{-ck}\, ,
\]
for constants $C,c>0$ depending only on $C_{\rm A}$, $\qbarrow$
and $\varepsilon_{\rm row}$.
\end{lemma}

\noindent{\bf Proof.} Use $\theta_{\rm q}$, $\chi_{\rm q}$ and $C_{\rm q}$
from the proof of Lemma~\ref{l.non-ruin-probability}.  Since $Q_0=2$,
summing~\eqref{e.queue-first-abandonment} with $m=0$ over $n\ge k$
gives, for $k\ge1$,
\[
 \PP(k\le\tau_\sigma<\infty\mid\mathcal G_\sigma)
 \le \frac{C_{\rm A}}{\qbarrow}e^{-2\theta_{\rm q}}
     \sum_{n\ge k}\chi_{\rm q}^{\,n-1}
 =C_{\rm q}e^{-2\theta_{\rm q}}\chi_{\rm q}^{\,k-1}.
\]
The assertion follows with $c=-\log\chi_{\rm q}>0$ and
$C=C_{\rm q}e^{-2\theta_{\rm q}}\chi_{\rm q}^{-1}$.
\qed

\begin{lemma}[Queue survival preserves the renewal level]
\label{l.queue-survival-preserves-renewal}
Let $J_\sigma=1$, and let $\varrho_\sigma$ be the renewal born at the active
slide of index~$\sigma$.  If $\tau_\sigma=\infty$, then $\varrho_\sigma$ is the only vertex of
the terminal cluster in its renewal band, and hence is a renewal point of
$S_L$ in the sense of Definition~\ref{d.renewal}.
\end{lemma}

\noindent{\bf Proof.} For brevity write $\varrho=\varrho_\sigma$.  The SimpleJoin event
$\mathsf R_\sigma$ gives
\begin{equation}
\label{e.queue-renewal-band-invariant}
 V\bigl(\csetuthree(\widehat S_\sigma)\bigr)
 \cap\vertbdry_{\varrho\cvl}\, = \, \{\varrho\}\, .
\end{equation}
We prove that this equality persists after each later active slide.  Suppose
it holds immediately before the active slide at~$h_{\sigma+n}$, and write
$\mathbf U_{n-1}^{(\sigma)}=(\mathfrak b_1,\ldots,\mathfrak b_q)$.  Since
$\tau_\sigma=\infty$, we have $D_n\le q$.  Thus $K_n=D_n-1<q$, so the
current midpiece does not pass the row of $\mathfrak b_{K_n+1}$.

The vertex~$\varrho$ belongs to the retained downcluster at the
current slide.  The SimpleJoin geometry gives
$h_\sigma<\varrho\cvl\le h_\sigma+r$; since $h_\sigma>2r$,
the origin lies below $\vertbdry_{\varrho\cvl}$.  Stop an open path
from~$0$ to~$\varrho$ at its first visit to this band.  Each edge
changes height by at most~$r$, and the induction hypothesis makes
$\varrho$ the only cluster vertex in the band.  Thus the stopped path
has maximal height~$\varrho\cvl$, giving
$\dbt(\varrho)\le\varrho\cvl$.  The schedule gap gives
$h_{\sigma+n}>h_\sigma+4r>\varrho\cvl$, so
Proposition~\ref{p.dmutwoheights}(1) places~$\varrho$ in the retained
downcluster.

By the slide reconstruction of Definition~\ref{d.prime}, the output
cluster consists of the retained downcluster, the newly sampled midpiece
and the transported upcluster.  The retained downcluster is a subgraph
of the pre-slide cluster and therefore contributes no new vertex in
$\vertbdry_{\varrho\cvl}$.  By~\eqref{e.southgoing-row-margin}, the transported
upcluster lies entirely north of the row of~$\mathfrak b_{K_n+1}$, while
Lemma~\ref{l.standing-breaker-list}(\ref{it.standing-breaker-separation})
places that row strictly north of the renewal band.  Hence any new vertex
in the renewal band would have to be reached by
$\midpiece_{h_{\sigma+n}}(\widehat S_{\sigma+n})$.  In particular,
a fresh edge cannot leave the retained downcluster at the renewal band:
such an edge belongs to $E(\downspace_{h_{\sigma+n}})$ and is fixed
closed by $C_v$ in Definition~\ref{d.prime}.  Fresh departures instead
start at the top of the retained downcluster or the bed of the shifted
upcluster, both north of that row.  A midpiece path to a new
vertex in the renewal band would therefore reach the southern side of
that row, passing $\mathfrak b_{K_n+1}$ in the sense specified before
Definition~\ref{d.active-slide-queue}.  This contradicts the definition
of $K_n$, so the identity~\eqref{e.queue-renewal-band-invariant} persists
after the slide.

Dormant attempted slides change no edge.  Induction over the remaining
active slides gives
\[
 V\bigl(\csetuthree(S_L)\bigr)
 \cap\vertbdry_{\varrho\cvl}\, = \, \{\varrho\}\, .
\]
Thus $\varrho$ remains a renewal point.
\qed

\section{The greedy restart argument}
\label{ss.greedy-restart}

When a queue is exhausted, the prospective renewal point is abandoned and
the search continues to the north.  The construction waits for another
renewal-producing slide and starts a new queue.  Such slides occur with
uniformly positive conditional probability, yielding a bound on the waiting
time.

\begin{lemma}[Waiting time for a renewal-producing slide]
\label{l.simplejoin-wait-tail}
Let $T$ be a stopping time of $(\mathcal G_j)_{j=0}^{L+1}$.
Put
\[
 \sigma(T)=\min\bigl(\{j:T<j\le L,\ J_j=1\}\cup\{L+1\}\bigr),
\]
and set $W_T=(\sigma(T)\wedge(N+1))-T$ on $\{T\le N\}$,
with $W_T=0$ otherwise.  Then
\[
 \PP(W_T\, \ge \, k\mid\mathcal G_T)
\, \le \, (1-p_{\rm ren})^{k-1},\qquad k\, \ge \, 1\, .
\]
\end{lemma}

\noindent{\bf Proof.} On $\{T\le N\}$, the event $\{W_T\ge k\}$ requires the next $k-1$
active slides to occur, with $J_j=0$ at each.  By
\eqref{e.queue-start-probability}, each has conditional failure probability
at most $1-p_{\rm ren}$ given $\mathcal G_{j-}=\mathcal G_{j-1}$,
which records earlier failures and whether this slide exists.
Conditioning on $\{T=t\}\in\mathcal G_t$ and iterating gives the bound.
The event $\{T\le N\}$ is $\mathcal G_T$-measurable, while $W_T=0$ on
its complement, so the bound holds as stated.
\qed

A produced renewal~$\varrho_j$ is said to \emph{survive} if
$\varrho_j\in R(S_L)$.

\begin{proposition}[Renewal survival with many active slides]
\label{p.renewal-survives}
There are $C,c>0$ such that
\[
 \Qdyn\bigl(\text{no renewal }\varrho_j\text{ with }J_j\, = \, 1
       \text{ survives in }S_L,\ N\, \ge \, \delta_{\rm f}u\bigr)
 \, \le \, \ Ce^{-cu}\, .
\]
\end{proposition}

\noindent{\bf Proof.}
Since the event $\{N\ge\delta_{\rm f}u\}$ depends on the entire slide trajectory,
conditioning on it need not preserve the one-step queue estimates.  We
therefore apply these estimates under $\Qdyn$ at successive stopping times,
bounding directly the event that $N\ge\delta_{\rm f}u$ but no produced
renewal survives.  Let
\[
 \sigma_1=\min\bigl(\{1\le j\le L:J_j=1\}\cup\{L+1\}\bigr).
\]
If $\sigma_1=L+1$, the construction stops unsuccessfully.  Otherwise follow
the queue born there.  When $\sigma_i\le L$, let
\[
 F_i=\{\tau_{\sigma_i}<\infty\};
\]
when $\sigma_i=L+1$, put $F_i=\emptyset$.  On~$F_i$, set
$T_i=\sigma_i+\tau_{\sigma_i}$ and
\[
 \sigma_{i+1}
 =\min\bigl(\{T_i<j\le L:J_j=1\}\cup\{L+1\}\bigr).
\]
A SimpleJoin slide cannot abandon a live queue: its two-edge bridge
lies above the standing rows, so $D=1$ and $\beta=1$, increasing the
reserve by one.  If \(\sigma_i\le L\) but \(F_i\) does not occur, the queue survives and
Lemma~\ref{l.queue-survival-preserves-renewal} shows that
$\varrho_{\sigma_i}$ remains a renewal point of~$S_L$.  If
\(\sigma_i=L+1\), no later birth occurs.  In either case the construction
stops; set $T_i=\sigma_{i+1}=L+1$ after stopping.

The restart times at which conditioning is used are stopping times of the
filtration $(\mathcal G_j)_{j=0}^{L+1}$.  Indeed $J_j$ is
$\mathcal G_j$-measurable, so the first later index with $J_j=1$ is a stopping
time, and the abandonment time of an adapted queue has the same property.
Induction therefore makes every $\sigma_i$ and $T_i$ stopping times.  The
value~$L+1$ records only that the construction has stopped.  The conditional
queue estimates established for a fixed birth index also apply at each
$\sigma_i$: restrict to $\{\sigma_i=j\}\in\mathcal G_j$ and apply the
corresponding bound at the deterministic index~$j$.

Put $T_0=0$.  On $F_1\cap\cdots\cap F_{i-1}$, let
\[
 W_i=(\sigma_i\wedge(N+1))-T_{i-1}.
\]
Thus $W_i$ measures the wait to the next queue birth, or to the terminal
index $N+1$ if no further birth occurs.  Set $W_i=0$ outside
$F_1\cap\cdots\cap F_{i-1}$.  Put
$A_i=\tau_{\sigma_i}=T_i-\sigma_i$ on~$F_i$, and set $A_i=0$ off~$F_i$.
Choose $\lambda>0$ so that $e^\lambda(1-p_{\rm ren})<1$ and
$e^\lambda\chi_{\rm q}<1$, where $\chi_{\rm q}<1$ is as in the proof of
Lemma~\ref{l.non-ruin-probability}.  Then
Lemmas~\ref{l.failed-attempt-tail} and
\ref{l.simplejoin-wait-tail} give finite constants
$K_{\rm w},K_{\rm a}$ such that, uniformly in the histories at the relevant stopping times,
\[
   \E\bigl[e^{\lambda W_i}\mid\mathcal G_{T_{i-1}}\bigr]\, \le \, K_{\rm w}\, ,
   \qquad
   \E\bigl[e^{\lambda A_i}\one_{F_i}\mid\mathcal G_{\sigma_i}\bigr]\, \le \, K_{\rm a}\, .
\]
For the $i$\textsuperscript{th} cycle, conditioning first at $\sigma_i$
and then at $T_{i-1}$ gives
\[
 \E\bigl[e^{\lambda(W_i+A_i)}\one_{F_i}\mid\mathcal G_{T_{i-1}}\bigr]
 \,\le\, K_{\rm a}\E\bigl[e^{\lambda W_i}\mid\mathcal G_{T_{i-1}}\bigr]
 \,\le\, K_{\rm w}K_{\rm a}\, .
\]
Put $M_{\rm cyc}=\max\{K_{\rm w},K_{\rm w}K_{\rm a},1\}$.  Iterating this
conditional calculation and appending the final wait gives, for every
$k\ge0$,
\begin{equation}
\label{e.greedy-cycle-overflow-moment}
 \E\left[
   \exp\left\{\lambda\left(
      \sum_{i\, = \, 1}^{k}(W_i+A_i)+W_{k+1}\right)\right\}
   \one_{F_1\cap\cdots\cap F_k}\right]\, \le \, M_{\rm cyc}^{k+1}\, .
\end{equation}

On $F_1\cap\cdots\cap F_k$ the $k$ failed cycles satisfy
$W_i+A_i=T_i-T_{i-1}$, so their durations telescope:
\[
   \sum_{i=1}^{k}(W_i+A_i)+W_{k+1}
   \, = \, \sigma_{k+1}\wedge(N+1).
\]
This identity also holds for $k=0$, with $T_0=0$.

Set $n_0=\delta_{\rm f}u$ and choose $\eta>0$ so small that
$\lambda-\eta\log M_{\rm cyc}>0$; put
$n_{\rm att}=\lfloor\eta n_0\rfloor$.  We may assume
that~$u$ is large enough that $n_{\rm att}\ge1$, since bounded~$u$ is covered by
enlarging the final constant.  Suppose first that,
on $\{N\ge n_0\}$, no renewal survives to the endpoint and fewer than $n_{\rm att}$ attempts
begin by active-slide time $N$.  If exactly $k<n_{\rm att}$ attempts begin,
then all of them
are abandoned by $N$ and $\sigma_{k+1}=L+1$.  The telescoping identity
therefore makes the total duration $N+1>n_0$.  By
\eqref{e.greedy-cycle-overflow-moment}, Markov's inequality and a union over
$0\le k<n_{\rm att}$,
\[
 \PP\bigl(\text{no surviving renewal},\,N\, \ge \, n_0,\ \sigma_{n_{\rm att}}>N\bigr)
\, \le \, \sum_{k=0}^{n_{\rm att}-1}e^{-\lambda n_0}M_{\rm cyc}^{k+1}
\, \le \, n_{\rm att}\,e^{-\lambda n_0}M_{\rm cyc}^{n_{\rm att}}
\, \le \, Ce^{-cu}\, .
\]

Suppose instead that at least $n_{\rm att}$ attempts begin by time $N$.
On $\{\sigma_i\le N\}$, Lemma~\ref{l.non-ruin-probability} gives
$\PP(F_i\mid\mathcal G_{\sigma_i})\le1-s$.  If no renewal survives,
all $n_{\rm att}$ attempts must fail.  At each birth~$\sigma_i$, the preceding
failures and the occurrence of that birth are known.  Repeated conditioning
at the birth times gives
\[
 \PP\bigl(\sigma_{n_{\rm att}}\le N,\,F_1\cap\cdots\cap F_{n_{\rm att}}\bigr)
\, \le \, (1-s)^{n_{\rm att}}\, \le \, Ce^{-cu}\, .
\]
Together, the two bounds prove the proposition.
\qed

The greedy rule begins a new queue only after the preceding one has been
abandoned.  At each new birth time $\sigma_i$, the failures of earlier queues
are therefore recorded in $\mathcal G_{\sigma_i}$, allowing the conditional
failure bounds to be iterated.

It remains to combine the queue estimate with the supply of active slides and
to invoke the invariance of the cluster law under the slide dynamics.

\begin{proposition}[Survival of a created renewal]
\label{p.local-survival}
There are $C,c>0$ such that
\[
   \Qdyn\bigl(\text{no renewal }\varrho_j\text{ with }J_j=1
       \text{ survives in }S_L\bigr)
   \, \le \, \ Ce^{-cu}\, .
\]
\end{proposition}

\noindent{\bf Proof.} We decompose on the supply event $\{N\ge \delta_{\rm f}u\}$.

Its complement is rare.  Theorem~\ref{t.vf-to-fine-delivery} gives at least
$\delta_{\rm f}u$ heights fine at limelight, except with probability $Ce^{-cu}$.
These are exactly the active slides, so
$\Qdyn(N<\delta_{\rm f}u)\le Ce^{-cu}$.

Proposition~\ref{p.renewal-survives} bounds the remaining joint event
by $Ce^{-cu}$.  Adding the bounds proves the claim. \qed

\noindent{\bf Proof of Theorem~\ref{t.intermediate-renewal}(2).}
By Proposition~\ref{p.qdyn-output}, the terminal cluster $S_L$ under the
slide dynamics~$\Qdyn$ has law $\overline\Q=\Q\circ\pi^{-1}$, the cluster
law induced by~$\Q$.  On the conditioning event~$\clusteruthree$,
$\csetuthree$ is precisely the origin cluster~$\mc C(0)$; hence
an event determined by $S_L$ has the same probability under $\Qdyn$ as
the corresponding event for~$\mc C(0)$ under~$\Q$.
On $\{J_j=1\}$, the SimpleJoin midvertex~$\varrho_j$ has height
$h_j<\varrho_j\cvl\le h_j+r<u$, where the final inequality uses the
scheduled-height restriction $h_j<u-2r$.  Proposition~\ref{p.local-survival}
therefore gives, for all sufficiently large~$u$,
\[
   \Q\bigl(\rlevel(\mc C(0))\cap(0,u)\neq\emptyset\bigr)
\, \ge \, 1-Ce^{-cu}\, .
\]
Enlarging $C$ covers the remaining bounded range~$u\ge r$, as required
for Theorem~\ref{t.intermediate-renewal}(2). \qed

With the exponential rarity of renewal-free clusters established, the slide
argument is complete.  Part~\ref{part.four} will harness this estimate in a
renewal-theoretic derivation of the mass gap and the Ornstein--Zernike
asymptotics.

\part{Ornstein--Zernike Consequences}\label{part.four}

\chapter{Asymptotic estimates via renewal theory}
\label{s.asp}

By analysing the slide dynamics, we have obtained
Theorem~\ref{t.intermediate-renewal}(2): under the point-to-plane
conditioning, the probability that a long finite cluster has no interior
renewal point is exponentially small in its height.  To derive the paper's
strongest renewal conclusion, Theorem~\ref{t.renewal}, according to which
renewal points populate the conditioned cluster with uniformly positive
density except with exponentially small probability, we now set aside the
iterated slide scheme and turn to classical renewal theory.  In the present
chapter this change of method yields the mass gap of
Theorem~\ref{t.mass-gap} and the sharp string asymptotic,
Theorem~\ref{t.stingwelltrap.sting}.  The latter leads to the sharp
point-to-plane cluster asymptotic stated in the introduction,
Theorem~\ref{t.stingwelltrap.cluster}.  The following chapter derives
renewal abundance and carries the resulting vector renewal description through to the
Ornstein--Zernike asymptotic and the analytic and strictly convex Wulff
geometry of Theorem~\ref{t.oz-intro}.

Renewal theory is naturally formulated for strings.  As recalled in
Section~\ref{ss.renewal-intro}, the renewal points of a string cut it uniquely
into sausages, the irreducible pieces that may be concatenated.  After the
occurrence weight of one sausage is tilted at the string decay rate~$\zeta$,
it will define the probability law~$\mu$; independent concatenation under
$\mu^{\otimes\N}$ then gives the associated one-dimensional renewal process.
A full cluster with at least two renewal points is recovered by adjoining a
lower and an upper cul-de-sac to the resulting middle string.
Theorem~\ref{t.intermediate-renewal}(2) first enters this construction in the
control of a single sausage.  A long sausage has no interior renewal point,
and closing a bounded collection of edges at its two endpoints turns it into
a full-lattice finite cluster of the same height without creating one.
Together with the definition of~$\zeta$, the theorem therefore gives the
normalisation of~$\mu$ and exponential control of the sausage height.  The
one-dimensional renewal theorem may then be applied to obtain the sharp
string asymptotic, and the endpoint culs-de-sac transfer this asymptotic to
full clusters.

The chapter has four sections. The first constructs the probability law~$\mu$
and proves its mass gap. The second derives the sharp string asymptotic by
renewal theory, then verifies the string decay input from~\cite{FriberghHammond}.
The third introduces the endpoint culs-de-sac and derives the sharp
point-to-plane cluster asymptotic. The fourth controls transverse displacement
for clusters and sausages. Its cluster estimate will also keep the endpoint
pieces localised when they enter the vector renewal representation.
Together, these results provide the needed elements for the
renewal-abundance and Ornstein--Zernike arguments of
Chapter~\ref{s.renewal-abundance-oz}.

\section{The sausage law and its mass gap}
\label{s.sausage-law-mass-gap}

We specify strings and sausages, show that the tilted weights define a
probability law~$\mu$, and prove the exponential moment bound of
Theorem~\ref{t.mass-gap}.
\begin{lemma}\label{l.stingbasics}
A string $S$ has a unique vertex of lowest height $\base(S)$ and a
unique vertex of greatest height $\head(S)$, both renewal points.
The renewal point set $R(S)$ is strictly and totally ordered by height.
In particular, every string has height at least~$r$.
\end{lemma}
\noindent\textbf{Proof.} Let $x,y$ be the endpoints in the definition of~$S$.
Since $S$ lies in $\slab_{x\cvl,y\cvl}$ and $x$ is its only vertex
in $\vertbdry_{x\cvl}$, every other vertex has height strictly greater
than~$x\cvl$.  Likewise, $y$ is the unique vertex of
greatest height.  Thus $\base(S)=x$ and $\head(S)=y$ are renewal points.
Lemma~\ref{l.geometrybasics}(2) gives the strict ordering of all renewal points.
Since $x\ne y$, the same lemma gives $\hgt(S)\ge r$. \qed

\begin{definition}\label{d.sausage}
A \hfff{sausage}{\emph{sausage}} is a string whose only renewal points are its base
and head.
\end{definition}

\begin{definition}\label{d.os}
The \emph{slab} of a string $S$ is $\slab_{\base(S)\cvl,\head(S)\cvl}$.
In specifying occurrence events, we regard $S$ as a fixed finite connected
subgraph of its slab, with its base and head as renewal points.
The string $S$ \hfff{occurrence-event}{\emph{occurs}} if every edge of $S$ is open and every edge
in its slab that is incident to $S$ but does not belong to $S$ is closed.
Write $\mathsf{O}(S)$ for the occurrence event.
Equivalently, the string occurs when it is the open component of its base
in its slab.
\end{definition}

By the boundary convention of Section~\ref{s.basics}, the occurrence
events of consecutive strings depend on disjoint edge sets.

For strings positioned end-to-end, $S_1\circ\cdots\circ S_m$
denotes their union as subgraphs.  When all strings are based at~$0$,
concatenation first translates each successive string to the head of its
predecessor.  Let $x_0,\ldots,x_m$ be the renewal points of a string~$S$,
listed in increasing height order, and set
$S_i=S\cap\slab_{x_{i-1}\cvl,x_i\cvl}$ for $1\le i\le m$.
The next result identifies these pieces and the converse construction.

\begin{lemma}[Cutting and concatenating sausages]
\label{l.sausage-concatenation}
Cutting a string $S$ at its renewal points yields sausages
$S_1,\ldots,S_m$, with $S=S_1\circ\cdots\circ S_m$.
Conversely, concatenation of any finite non-empty sequence of sausages based at~$0$
yields a string based at~$0$, whose renewal points are precisely the
endpoints of its constituent sausages.  If $T_i$ denotes the origin-based
translate of $S_i$, then
\[
  \mathsf O(S)=\bigcap_{i=1}^m\mathsf O(S_i),
  \qquad \hgt(S)=\sum_{i=1}^m\hgt(T_i).
\]
\end{lemma}
\noindent\textbf{Proof.} Every vertex of~$S$ in
$V(\slab_{x_{i-1}\cvl,x_i\cvl})$ can be joined to~$x_{i-1}$ by a path
in that slab.  By Lemma~\ref{l.slab-entry-faces}, a self-avoiding path
in~$S$ cannot leave and return: an excursion must re-enter through the
side it left, revisiting its unique renewal vertex.  Moreover, an edge
of~$S$ with both endpoints in the smaller slab's vertex set belongs to
that slab; otherwise both endpoints would lie in the lower renewal band.
Thus $S_i$ is connected.  Its endpoints remain renewal points.
An additional renewal point of~$S_i$ would also be one of~$S$:
the smaller slab contains every vertex of~$S$ in its renewal band.
This contradicts the consecutiveness of $x_{i-1},x_i$.
The half-open edge allocation gives $S=S_1\circ\cdots\circ S_m$.

Conversely, at a junction of two translated sausages, the lower one has
no vertex other than its head in the band of width~$r$ ending there,
while the upper one has no vertex other than its base at or below the
junction height.  No other sausage meets this band, since each has height
at least~$r$.  Consequently the pieces meet only at their junctions,
which are renewal points of the union.  A non-endpoint vertex cannot
become a renewal point on
concatenation, since it is not a renewal point of its own sausage.
The union is therefore a string with exactly the asserted renewal points,
and cutting it recovers its constituent sausages.

Finally, the renewal bands ensure that the vertices of the union
lying in the slab of~$S_i$ are exactly those of~$S_i$.  These slabs
partition the edges of the full slab, so the open and closed edge
prescriptions in Definition~\ref{d.os} agree on both sides of the
displayed occurrence-event identity.  Height additivity is immediate. \qed

As we next state,
the string probability $\PP \big( \sting_{[u,u+r)} \big)$ decays in the form $\exp \{ - (\zeta + o(1))u \}$.
In essence, the assertion is due to~\cite[Lemma~$8.5$]{FriberghHammond}, and we derive it from this input in the next section.
\begin{proposition}\label{p.zeta}
There exists $\zeta = \zeta(p,d,\vecell) \in (0,\infty)$ such that
\[
  \lim_{u\to\infty} u^{-1}\log\PP\bigl(\sting_{[u,u+r)}\bigr)\, = \, -\zeta\, .
\]
\end{proposition}

The string decay rate first gives a corresponding upper bound for finite
clusters.  Combined with the exponential rarity of a conditioned cluster
without an interior renewal, this yields an extra exponential factor for
renewal-free clusters.  An endpoint comparison will carry this gain back to
sausages and give the required control of their tilted weights.

\begin{lemma}[First-order point-to-plane upper bound]
\label{l.first-order-cluster-upper}
For every $\varepsilon>0$ there is $C_\varepsilon<\infty$ such that
\[
 \PP\bigl(\cluster_{[u,u+3r]}\bigr)
\, \le \, C_\varepsilon e^{-(\zeta-\varepsilon)u}\, ,
 \qquad u\, \ge \, 0\, .
\]
\end{lemma}

\noindent{\bf Proof.} Given a cluster counted on the left, with maximum
height~$H$, attach one open $e_d$-edge below its base and one above its
head, and close the other edges incident to the two new extreme vertices.
These vertices, of heights $-r$ and $H+r$, are the unique vertices of the
enlarged cluster in $\vertbdry_{-r}$ and $\vertbdry_{H+r}$.
Since every edge of the enlarged cluster lies in the slab between these
heights, the enlarged cluster is its own slab component and hence a string.  Translating by~$e_d$ gives a string based at the origin
whose height lies in $[u+2r,u+5r]$, and hence an occurrence of
$\sting_{[u+2r,u+5r]}$.  The modification changes the states of at most
$4d$ edges.  After undoing the translation, deleting the two added edges
and their outer endpoints recovers the original cluster.  The modified edges
are incident to the two new endpoints, so their locations are determined by
the resulting string and their original states admit only a bounded number
of choices.  Finite energy therefore gives
\(\PP(\cluster_{[u,u+3r]})\le C\PP(\sting_{[u+2r,u+5r]})\).
Cover the interval on the right by four half-open intervals of length~$r$.
For each $\varepsilon>0$, Proposition~\ref{p.zeta} bounds each of the four
probabilities by $e^{-(\zeta-\varepsilon)u}$ for all sufficiently large~$u$.
Enlarging the constant covers bounded~$u$. \qed

\begin{lemma}[Cluster upper bound without an intermediate renewal]
\label{l.no-renewal-cluster-upper}
There exist $C,c>0$ such that, for every $u\geq0$,
\[
 \PP\bigl(\cluster_{[u,u+3r]},\,
       \rlevel(\mc C)\cap(0,u)\, = \, \varnothing\bigr)
\, \leq \, Ce^{-\zeta u-cu}\, .
\]
\end{lemma}

\noindent\textbf{Proof.} For $u\geq r$, Theorem~\ref{t.intermediate-renewal}(2) bounds the conditional
probability of no renewal in $(0,u)$ by $Ce^{-c_0u}$.  Multiplying by
Lemma~\ref{l.first-order-cluster-upper} with
$\varepsilon<c_0/2$ gives the assertion, after decreasing the remaining
exponential rate.  Enlarging~$C$ covers $u<r$. \qed

\begin{definition}\label{d.sausagemu}
Let $Y$ denote the set of strings with base $0 \in \Z^d$.  Define the
\hfff{sausage-law}{tilted weight \(\mu : Y \to (0,\infty)\)} by
\[
  \mu(S)\, = \, \exp\{\zeta\cdot\hgt(S)\}\,\PP(\mathsf{O}(S))\, ,
\]
where $\zeta > 0$ is given by Proposition~\ref{p.zeta}.
\end{definition}

\begin{lemma}\label{l.percoprod}
Let $S$ be a string based at~$0$, write $S=S_1\circ\cdots\circ S_m$
for its successive sausage pieces, and let $T_i$ be the translate of~$S_i$
with base at~$0$.  Then
\[
  \mu(S)\, = \, \prod_{i=1}^m \mu(T_i), \qquad
  \PP(\mathsf{O}(S))\, = \, \prod_{i=1}^m \PP(\mathsf{O}(T_i))\, .
\]
\end{lemma}
\noindent\textbf{Proof.} Lemma~\ref{l.sausage-concatenation} identifies
$\mathsf O(S)$ with $\bigcap_i\mathsf O(S_i)$.  These events depend on
disjoint edge sets and are independent.  Translation invariance gives
$\PP(\mathsf O(S_i))=\PP(\mathsf O(T_i))$, proving the probability
identity; height additivity gives the identity for~$\mu$. \qed

\begin{definition}\label{d.sausageset}
Let $\sausageset$ denote the set of sausages with base $0 \in \Z^d$.
\end{definition}

For the endpoint edge counts below, record the number of coordinate directions
tangent to the height hyperplanes:
\[
 k_0=\#\{i\in\intint{d}:e_i\cvl=0\}.
\]

\begin{lemma}[Endpoint-closure comparison for strings]
\label{l.string-endpoint-closure}
Let $\mathcal A$ be any collection of strings based at the origin.  If
$\mc C$ denotes the full-lattice open cluster of the origin, then
\[
 \sum_{S\in\mathcal A}\PP\bigl(\mathsf O(S)\bigr)
\, = \, (1-p)^{-2d}\PP\bigl(\mc C\in\mathcal A\bigr)\, .
\]
Here the event on the right requires that $\mc C$ be finite and equal, as a
subgraph, to one of the strings in~$\mathcal A$.
\end{lemma}
\noindent\textbf{Proof.} Write $v\in e$ when the vertex~$v$ is an
endpoint of the edge~$e$.  For a string~$S$, set
\[
 B_-(S)=\{e\in E(\downspace_{\base(S)\cvl}):\base(S)\in e\},
 \qquad
 B_+(S)=\{e\in E(\upspace_{\head(S)\cvl}):\head(S)\in e\}.
\]
By the half-open edge allocation, the lower set has \(d+k_0\) edges and
the upper set has \(d-k_0\) edges: a coordinate direction with \(e_i\cvl=0\) contributes both horizontal incident
edges at the lower endpoint and none at the upper endpoint, while every
positive-height coordinate direction contributes one edge at each end.
Thus
\[
 |B_-(S)|+|B_+(S)|=2d.
\]
By the half-open edge allocation, every edge incident to~$V(S)$ but
outside the edge set of its slab belongs either to
$E(\downspace_{\base(S)\cvl})$ or to
$E(\upspace_{\head(S)\cvl})$.  The renewal conditions ensure that
only the base and head of~$S$ are incident to such edges on the respective
sides.  Thus these edges are exactly $B_-(S)\cup B_+(S)$, and their states
are not prescribed by $\mathsf O(S)$.  The event $\{\mc C=S\}$ occurs
if and only if $\mathsf O(S)$ occurs and all these endpoint edges are
closed.  Since they lie outside the slab and number~$2d$, independence gives
\(\PP(\mc C=S)=(1-p)^{2d}\PP(\mathsf O(S))\).
The events $\{\mc C=S\}$ are disjoint as~$S$ varies.  Summing proves the
identity. \qed

\begin{lemma}[Exponential tilted weight of a sausage]
\label{l.sausage-tilted-tail}
There exist $C,c>0$ such that, for every $u\geq0$,
\[
 \sum_{\substack{S\in\sausageset\\ \hgt(S)\, \geq \, u}}
     e^{\zeta\hgt(S)}\PP\bigl(\mathsf O(S)\bigr)
\, \leq \, Ce^{-cu}\, .
\]
Consequently, for some $\theta>0$,
\[
 \sum_{S\in\sausageset}
     e^{(\zeta+\theta)\hgt(S)}\PP\bigl(\mathsf O(S)\bigr)<\infty\, .
\]
\end{lemma}

\noindent{\bf Proof.} Put $I_n=[u+nr,u+(n+1)r)$ and $v_n=u+nr$.  Apply
Lemma~\ref{l.string-endpoint-closure} to the sausages whose heights lie in
$I_n$.  Since such a sausage has no interior renewal point, the resulting
full-lattice cluster has no renewal point at a height in $(0,v_n)$.
For some $c_0>0$, Lemma~\ref{l.no-renewal-cluster-upper} therefore gives
\[
 \sum_{\substack{S\in\sausageset\\ \hgt(S)\in I_n}}
   \PP\bigl(\mathsf O(S)\bigr)
 \, \le \, C\PP\bigl(\cluster_{[v_n,v_n+3r]},\,
                 \rlevel(\mc C)\cap(0,v_n)=\varnothing\bigr)
 \, \le \, C e^{-\zeta v_n-c_0v_n}\, .
\]
Since
$e^{\zeta\hgt(S)}\le e^{\zeta(v_n+r)}$ on $I_n$, the total tilted weight in this
window is at most $Ce^{-c_0v_n}$.  Summing over $n\geq0$
proves the first assertion.  Summation by height windows after multiplication
by $e^{\theta\hgt(S)}$, with $\theta<c_0$, proves the second. \qed

Next is the key algebraic identity.  Let $C(u)$ denote
the open component of $0$ in $\slab_{0,u}$, and
$\mc{S}_u$ the set of strings with base $0$ and $\head(S)\cvl = u$.

\begin{lemma}\label{l.essn}
For any $A \subseteq \mc{S}_u$,
\[
  \PP(C(u) \in A)\, = \, e^{-\zeta u}\sum_{S \in A}\mu(S)\, .
\]
\end{lemma}
\noindent\textbf{Proof.} For $S\in\mc S_u$,
$\mathsf O(S)=\{C(u)=S\}$ and $\hgt(S)=u$.  Sum over~$A$ and apply
Definition~\ref{d.sausagemu}. \qed

Within short height windows, string probabilities may be summed without double
counting.
\begin{lemma}[Disjoint string-height windows]
\label{l.string-window-disjointness}
If $I$ is a half-open interval of length at most~$r$, the positive-height
string events $\{\sting_h:h\in I,\ h>0\}$ are pairwise disjoint.
\end{lemma}

\noindent\textbf{Proof.} Suppose that strings with terminal heights $h<h'$ occur in the same
configuration.  The smaller-slab component is contained in the larger one,
so its head, of height~$h$, is a second vertex of the larger string in the
terminal band $(h'-r,h']$ whenever $h'-h<r$.  This contradicts renewal of
the larger head. \qed

\begin{proposition}\label{p.muprob}
The tilted string weight~$\mu$ on $\sausageset$ is a probability measure:
\[
 \sum_{S\in\sausageset}\mu(S)=1\, .
\]
\end{proposition}

\textbf{Proof of Proposition~\ref{p.muprob}.}
Write
\[
  \Lambda(t)\, = \, \sum_{S\in\sausageset}\PP(\mathsf O(S))\,e^{t\,\hgt(S)}
  \;\in\;[0,\infty]
\]
for the untilted sausage generating function; it is non-decreasing in $t$, and
$\Lambda(\zeta)=\sum_{S\in\sausageset}\mu(S)=\mu(\sausageset)$ by
Definition~\ref{d.sausagemu}.  By Lemmas~\ref{l.sausage-concatenation}
and~\ref{l.percoprod}, each string based at~$0$ corresponds uniquely to a
finite non-empty sequence of sausages based at~$0$, with additive heights and
multiplicative occurrence probabilities.  Hence, writing
$\Gamma(t)=\sum_{S}\PP(\mathsf O(S))\,e^{t\,\hgt(S)}$ for the
generating function summed over \emph{all} strings with base $0$, Tonelli gives
\[
  \Gamma(t)\, = \, \sum_{m\ge1}\Lambda(t)^m\, ,
\]
so that $\Gamma(t)<\infty$ if and only if $\Lambda(t)<1$.  By
Proposition~\ref{p.zeta} the string probability decays at exponential rate
$\zeta$, that is $u^{-1}\log\PP(\sting_{[u,u+r)})\to-\zeta$.  By Lemmas~\ref{l.essn} and~\ref{l.string-window-disjointness},
grouping string heights into the windows $[nr,(n+1)r)$, $n\ge1$, gives
$\Gamma(t)\le\sum_{n\ge1}e^{t(n+1)r}\,\PP\bigl(\sting_{[nr,(n+1)r)}\bigr)<\infty$ for
every $t\in[0,\zeta)$.

\emph{Upper bound $\mu(\sausageset)\le1$.}  For each $t\in[0,\zeta)$ the finiteness of
$\Gamma(t)$ forces $\Lambda(t)<1$.  Letting $t\nearrow\zeta$, monotone convergence
gives $\mu(\sausageset)=\Lambda(\zeta)=\lim_{t\nearrow\zeta}\Lambda(t)\le1$.

\emph{Lower bound $\mu(\sausageset)\ge1$.}  By
Lemma~\ref{l.sausage-tilted-tail}, $\Lambda$ is finite on an interval
$[\zeta,\zeta+\theta]$ and continuous there by dominated convergence.  If $\Lambda(\zeta)<1$, continuity would give
$\Lambda(t)<1$ for some $t>\zeta$.  The identity
$\Gamma(t)=\sum_{m\ge1}\Lambda(t)^m$ would then make $\Gamma(t)$ finite.  Grouping its
terms by the height windows $[nr,(n+1)r)$, $n\ge1$, would force
\(\PP(\sting_{[nr,(n+1)r)})\le\Gamma(t)e^{-tnr}\), contradicting Proposition~\ref{p.zeta}, which places the exponential decay
rate of this probability at~$\zeta$.  Thus $\Lambda(\zeta)\ge1$.

Combining the two bounds, $\sum_{S\in\sausageset}\mu(S)=\mu(\sausageset)=1$. \qed

\textbf{Proof of Theorem~\ref{t.mass-gap}.}
By Proposition~\ref{p.muprob}, $\mu$ is a probability measure, and the second
assertion of Lemma~\ref{l.sausage-tilted-tail} is exactly
$\E_\mu[e^{\theta\hgt(S)}]<\infty$ for some $\theta>0$. \qed

\section{Renewal theory and the sharp string asymptotic}
\label{s.provingsting}

With the sausage law normalised and its exponential moment established, we
turn to the renewal argument for the sharp string asymptotic.  We identify
the subgroup generated by sausage heights before applying the
appropriate renewal theorem to independent concatenations.  The requisite
finite mean follows from Theorem~\ref{t.mass-gap}.  The section concludes
by verifying the string decay input through comparison with the sausage
connections of~\cite{FriberghHammond}.

\begin{lemma}[Strings to high endpoints]
\label{l.high-endpoint-string}
For $x\in\Z^d$ with $x\cvl\geq3r$, a string from~$0$ to~$x$ exists.
\end{lemma}
\noindent\textbf{Proof.} Put $h=x\cvl$.  By Lemma~\ref{l.heightordering},
there is a path from~$0$ to~$x-e_d$ whose non-initial vertices have
heights in $(0,h-r]$.  Append the $e_d$-edge ending at~$x$.
The resulting path lies in $\slab_{0,h}$ and, by
Lemma~\ref{l.geometrybasics}(1), meets $\vertbdry_0$ only at~$0$
and $\vertbdry_h$ only at~$x$.
Prescribe its edges to be open and every other edge of the slab incident
to the path to be closed.  This finite prescription makes the open
component of~$0$ in the slab exactly the path, with $0$ and $x$ as
renewal points. \qed

\begin{lemma}[The sausage-height group]
\label{l.sausage-height-group}
The additive subgroup of~$\R$ generated by the heights of sausages
in $\sausageset$ equals the $\Z$-linear span of
$\{e_i\cvl:1\leq i\leq d\}$.  Consequently,
the sausage-height law is non-arithmetic when~$\vecell$ is non-lattice and
has span~$s$ when~$\vecell$ is lattice with span~$s$.
\end{lemma}

\noindent\textbf{Proof.} Every sausage endpoint lies in~$\Z^d$, so every
sausage height belongs to the $\Z$-linear span of the coordinate heights.  For the reverse
inclusion, fix~$i$ and choose $x\in\Z^d$ so that both $x\cvl$ and
$(x+e_i)\cvl$ exceed~$3r$.  By Lemma~\ref{l.high-endpoint-string},
strings from~$0$ to these two endpoints exist.  By
Lemma~\ref{l.sausage-concatenation}, cutting these strings into sausages
and translating the pieces to the origin expresses both endpoint heights
as sums of sausage heights.  Their difference is $e_i\cvl$, which
therefore belongs to the generated subgroup.

The group just identified is $s\Z$ when $\vecell$ is lattice and is
dense in $\R$ when $\vecell$ is non-lattice.  Thus the sausage-height law has
span~$s$ in the former case and is non-arithmetic in the latter. \qed

Let \hfff{renewal-chain}{$\mu^{\otimes\N}$} denote the law of a
sequence $S_1,S_2,\ldots$ of independent sausages drawn from $\mu$ and
concatenated end-to-end.  Write $Z_i=\hgt(S_i)$ and let
$\rlevel=\{0,Z_1,Z_1+Z_2,\ldots\}$ be the renewal levels.  Their renewal
measure is
\[
 U(A)\, = \, \E_{\mu^{\otimes\N}}\bigl[|\rlevel\cap A|\bigr]
\]
for Borel $A\subseteq[0,\infty)$.

For an interval $I\subseteq[0,\infty)$, write
\[
 \mathsf M_{\rm str}(I)
\, = \, \sum_{h\in I}\PP(\sting_h)\, .
\]
The sum runs over zero and the countable set of attainable positive heights,
with $\sting_0$ understood to hold identically and $\stingset_0$ to be
the one-vertex graph at the origin.
When $I\subset(0,\infty)$ is half-open and has length at most~$r$,
Lemma~\ref{l.string-window-disjointness} gives
$\mathsf M_{\rm str}(I)=\PP(\sting_I)$.  This identity does not hold
for a semi-infinite interval: events at terminal heights separated by at
least~$r$ may overlap.

\begin{lemma}[String--renewal identities and asymptotics]\label{l.zetarenewal}
\leavevmode
\begin{enumerate}
\item[\rm(1)] For $u\geq0$, $e^{\zeta u}\PP(\sting_u) = \mu^{\otimes\N}(u\in\rlevel)$.
\item[\rm(2)] For $u\geq0$ and any set $A$ of finite subsets of $(0,u)$,
\[
 e^{\zeta u}\PP(\sting_u,\,\rlevel(\stingset_u)\cap(0,u)\in A)
 \, = \, \mu^{\otimes\N}(u\in\rlevel,\,\rlevel\cap(0,u)\in A)\, .
\]
\item[\rm(3)] For $u>0$ and $k\geq2$,
\[
 e^{\zeta u}\PP(\sting_u,\,|R(\stingset_u)|=k)
 \, = \, \mu^{\otimes\N}(u\in\rlevel,\,|\rlevel\cap(0,u)|=k-2)\, .
\]
\end{enumerate}
\smallskip
\noindent If $\vecell$ is non-lattice, the renewal asymptotics are:
\begin{enumerate}
\item[\rm(4)] $\lim_{a\to\infty}U([a,a+\phi)) = \phi\,\mu_\infty$ for every $\phi>0$,
where $\mu_\infty = 1/\E_\mu[\hgt(S)]$ is the asymptotic renewal density; in
particular $U([u,u+r))\to r\mu_\infty$.
\item[\rm(5)] The local string probability and the semi-infinite
terminal-height mass satisfy
\[
\begin{split}
 e^{\zeta u}\PP(\sting_{[u,u+r)})
   &\longrightarrow \mu_\infty\zeta^{-1}(1-e^{-\zeta r})\, ,\\
 e^{\zeta u}\mathsf M_{\rm str}([u,\infty))
   &\longrightarrow \mu_\infty\zeta^{-1}\, .
\end{split}
\]
\end{enumerate}
\end{lemma}
\noindent{\bf Proof: (1,2,3).} For $u>0$, Lemmas~\ref{l.sausage-concatenation}
and~\ref{l.percoprod} identify strings of height~$u$ with finite initial
sequences of sausages whose cumulative height is~$u$, the weight of each
string being the product of the sausage weights.  The renewal points of the
string are exactly the initial point and the successive junctions; positive
sausage heights ensure that any given level is hit at most once.  Summing
Lemma~\ref{l.essn} over these sequences gives~(1), and recording the
interior junction heights gives~(2).  A string with $k$ renewal points
has $k-1$ sausages, giving~(3).  The cases $u=0$ in~(1) and~(2)
follow from the one-vertex convention.

\noindent{\bf (4).} The increment law has finite mean by
Theorem~\ref{t.mass-gap} and is non-arithmetic by
Lemma~\ref{l.sausage-height-group}.  Blackwell's renewal theorem
\cite[Chapter~XI]{Feller} gives
$U([a,a+\phi))\to\phi/\E_\mu\hgt(S)$ for every~$\phi>0$, proving~(4).
Equivalently, the translated measures $U(u+\,\cdot\,)$ converge vaguely on
$[0,\infty)$ to the measure $\mu_\infty\,dt$.

\noindent{\bf (5).} Summing the exact identity in~(1) over terminal heights gives
\[
 e^{\zeta u}\mathsf M_{\rm str}([u,u+r))
\, = \, \int_{[0,r)}e^{-\zeta t}\,U(u+dt)\, .
\]
Since the limiting measure does not charge the endpoints of~$[0,r)$,
the vague convergence just obtained makes this tend to
\(\mu_\infty\int_0^r e^{-\zeta t}\,dt
=\mu_\infty\zeta^{-1}(1-e^{-\zeta r})\).  The interval has length~$r$, so
$\mathsf M_{\rm str}([u,u+r))=\PP(\sting_{[u,u+r)})$.

The same identity with $[0,r)$ replaced by $[0,\infty)$ gives the
semi-infinite string mass.  By Lemma~\ref{l.stingbasics}, successive
renewal levels are separated by at least~$r$; hence
$U([a,a+r))\leq1$ for every~$a>0$.  Splitting $[0,\infty)$ into intervals
$[jr,(j+1)r)$ therefore bounds the contribution from $[mr,\infty)$,
uniformly in~$u$, by \(\sum_{j\geq m}e^{-\zeta jr}\).
On the finite interval $[0,mr)$ vague convergence applies.  Letting first
$u\to\infty$ and then $m\to\infty$ proves
\(e^{\zeta u}\mathsf M_{\rm str}([u,\infty))\to\mu_\infty/\zeta\). \qed

In the lattice case, the classical renewal theorem gives convergence along
the cover of~$\vecell$ to a positive limit of the probabilities
$\mu^{\otimes\N}(u\in\rlevel)$.  The exponential
moment of the sausage-height law yields an exponential rate of convergence.

\begin{lemma}\label{l.renewalrate}
Let $\vecell$ be lattice with span~$s$, and set
$g=s/\E_\mu[\hgt(S)]>0$.  There exist $C,c>0$ such that
for $u\in s\N$,
\[
  |\mu^{\otimes\N}(u\in\rlevel) - g|\, \leq \, Ce^{-cu}\, .
\]
\end{lemma}
\noindent\textbf{Proof.} Put
$f_k=\mu(\hgt(S)=ks)$ and $F(z)=\sum_{k\ge1}f_kz^k$.
By Lemma~\ref{l.sausage-height-group} and Theorem~\ref{t.mass-gap},
$\gcd\{k:f_k>0\}=1$ and $F$ is analytic on $|z|<R$ for some $R>1$.
Writing $v_n=\mu^{\otimes\N}(ns\in\rlevel)$, the renewal identity gives
$\sum_{n\ge0}v_nz^n=(1-F(z))^{-1}$ for $|z|<1$.

For $|z|<1$ we have $|F(z)|<1$.  On $|z|=1$, equality $F(z)=1$
forces $z^k=1$ whenever $f_k>0$, and hence $z=1$ by aperiodicity.
Moreover, $F'(1)=\E_\mu\hgt(S)/s=1/g$.  By isolation of zeros
and compactness of the unit circle, for some $\rho>1$,
\[
   \frac{1}{1-F(z)}-\frac{g}{1-z}
\]
is analytic on $|z|<\rho$: its singularity at~$1$ is removable,
and $1-F$ has no other zero in this disk.
For any $1<\rho_0<\rho$, Cauchy's coefficient estimate on
$|z|=\rho_0$ bounds the $n$\textsuperscript{th} Taylor coefficient of this function,
namely~$v_n-g$, by $C\rho_0^{-n}$.  This proves the lemma. \qed

Recall that $\simvecell$ and $\approxvecell$ denote asymptotic equivalence
along the cover of $\vecell$, the latter with exponentially small relative error.

\begin{theorem}[Sharp string asymptotics]
\label{t.stingwelltrap.sting}
There exists $\kappasting = \kappasting(d,p,\vecell) \in (0,\infty)$
such that, if $\vecell$ is lattice, $\PP(\sting_u) \approxvecell
\kappasting e^{-\zeta u}$; if $\vecell$ is non-lattice,
$\PP(\sting_{[u,u+r)}) \simvecell \kappasting e^{-\zeta u}$.
\end{theorem}
\noindent\textbf{Proof.} The non-lattice case is Lemma~\ref{l.zetarenewal}(5) with
\(\kappasting=\mu_\infty\zeta^{-1}(1-e^{-\zeta r})\).  In the lattice
case, Lemma~\ref{l.zetarenewal}(1) and Lemma~\ref{l.renewalrate}
give $\kappasting=g$ and the exponential error required by~$\approxvecell$. \qed

We finish the section by deriving the string decay input,
Proposition~\ref{p.zeta}, from~\cite[Lemma~8.5]{FriberghHammond}.
The sausage connections of that paper use a different slab-edge allocation
and different endpoint conditions.  We recall their reflected form before
comparing them with our strings.

Relabel the coordinate directions so that $e_1$ in
\cite{FriberghHammond} is our~$e_d$, with $e_d\cvl=r$.
For an edge~$e$, let $h^-(e)\le h^+(e)$ be its endpoint heights.  For vertices $b,t$ with $b\cvl<t\cvl$, set
\[
 E_{\rm sc}(b,t)=\{e:h^-(e)\ge b\cvl,\ h^-(e)<t\cvl\},
\]
and let $\mathcal K_{b,t}$ be the open component of~$b$ in this
restricted graph.  Write $\mc H^-(a)=\{z\in\Z^d:z\cvl\leq a\}$.
For a positive integer~$n$, let $\mathsf{SC}^+(n)$ be the event that
there is a unique vertex~$t$ with $t\cvl\in[n-r,n)$ for which
$\mathcal K=\mathcal K_{0,t}$ is finite, contains~$t$, and satisfies
\[
\begin{aligned}
 \mathcal K\cap\{v:v\cvl\leq r\}&=\{0,e_d\},\\
 \mathcal K\cap\{v:v\cvl\ge t\cvl-r\}&=\{t-e_d,t\}.
\end{aligned}
\]
These are conditions on vertex half-spaces.  For the large values of~$n$
used below, the two endpoint half-spaces are disjoint.  The event is the
image under $z\mapsto-z$ of the sausage-connection event from~$0$ to
$\mc H^-(-n)$ in~\cite[Section~8.2.1]{FriberghHammond}.

\begin{lemma}[Comparison with sausage connections]
\label{l.fh-sausage-comparison}
There are $c,C>0$ such that, for every sufficiently large real~$u$,
\begin{equation}
\label{e.fh-sausage-comparison}
 c\,\PP\bigl(\mathsf{SC}^+(\lfloor u\rfloor)\bigr)
\,\leq\, \PP\bigl(\sting_{[u,u+r)}\bigr)
\,\leq\, C\,\PP\bigl(\mathsf{SC}^+(\lceil u+3r\rceil)\bigr)\, .
\end{equation}
\end{lemma}

\noindent{\bf Proof.} The string slab with endpoints $b,t$ has edge set
\[
 E_{\rm str}(b,t)=\{e:h^+(e)>b\cvl,\ h^+(e)\le t\cvl\}.
\]

First suppose that a reflected sausage connection runs from~$0$ to~$t$.
The two endpoint conditions show that its restricted component~$\mathcal K$
becomes a string with the same endpoints once the edges in
$E_{\rm str}(0,t)\setminus E_{\rm sc}(0,t)$ incident to~$\mathcal K$
are closed.  Indeed, any such edge must be horizontal at height~$t\cvl$
and incident to~$t$: an edge crossing the lower cut from below cannot
meet~$\mathcal K$, by the lower endpoint condition.  Conversely, an edge of
$E_{\rm sc}(0,t)\setminus E_{\rm str}(0,t)$ crossing above the
upper cut cannot meet~$\mathcal K$: its lower endpoint would have
height strictly between~$t\cvl-r$ and $t\cvl$, where
$\mathcal K$ has no vertex.  The edges of this
difference incident to~$\mathcal K$ are therefore horizontal at the
base~$0$, and are already closed by the lower endpoint condition.
Thus no other edge prescriptions change.
We may append a bounded number of open $e_d$-edges above~$t$,
closing all other newly admitted edges incident to the old head and the
added vertices.  The resulting component is a string from~$0$ to its new head.

Conversely, given a string~$S$ from~$0$ to~$y$, adjoin the edge
$\{-e_d,0\}$ and a bounded positive number of consecutive $e_d$-edges
above~$y$.  Close the other edges admitted by the resulting
restricted graph and incident to the augmented graph
that are not already prescribed by~$\mathsf O(S)$.  Every such edge
is incident to the old base, the old head, or one of the added vertices:
for a vertex strictly between the old renewal bands, all incident
edges admitted by the new restriction were already inspected by the
string slab.  The new restricted component has exactly the two
required vertices in each endpoint half-space.  Translation by~$e_d$
then puts its lower endpoint at~$0$.

To verify the uniqueness required in $\mathsf{SC}^+(n)$ for the second
surgery, suppose two distinct terminal vertices satisfy the endpoint
conditions in the same window $[n-r,n)$.  The smaller restricted
component would be contained in the larger one, placing an extra vertex
in the upper endpoint half-space of the latter.  The case of equal
terminal heights is ruled out in the same way.

These surgeries have a uniform finite-energy cost and uniformly bounded
multiplicity.  In the first surgery, Lemma~\ref{l.string-window-disjointness}
identifies the output string uniquely in its height window; in the
second, the terminal vertex is unique.  To see the
multiplicity bound, fix the number of added steps.  The original terminal
vertex is recovered from the output head
by subtracting those steps in the $e_d$-direction (and undoing the fixed
translation in the second surgery).  Only the states of the bounded
set of edges at the affected endpoints and along the added chain may
have changed, so an output has at most $2^M$ possible preimages for
some fixed~$M$.  There are uniformly finitely many possible numbers
of steps.  Finite energy and bounded multiplicity therefore give the
claimed comparisons of event probabilities, with constants depending
only on $p,d,r$.

For the left inequality take $n=\lfloor u\rfloor$.  On
$\mathsf{SC}^+(n)$ the head height $h\in[n-r,n)$.  Append $k$ $e_d$-edges, where $k$ is the unique positive integer
such that $h+kr\in[u,u+r)$.
The number $k$ is uniformly bounded because $0\le u-n<1$.
The first surgery and the finite-energy comparison yield the claim.

For the right inequality take $n=\lceil u+3r\rceil$.
A string with head height $h\in[u,u+r)$ is extended downwards by one
$e_d$-edge and upwards by $k+1$ such edges, where $k\ge0$ is least
such that $h+(k+2)r\ge n-r$.  Then
$h+(k+2)r\in[n-r,n)$, and $k$ is uniformly bounded because
$0\le n-(u+3r)<1$.  The second surgery gives the required
comparison. \qed

\noindent\textbf{Proof of Proposition~\ref{p.zeta}.}
Lemma~8.5 of~\cite{FriberghHammond} gives
\[
 \PP\bigl(\mathsf{SC}^+(n)\bigr)
\,=\,\exp\{-\zeta n(1+o(1))\}, \qquad n\to\infty\text{ through integers},
\]
for some $\zeta\in(0,\infty)$; $\PP$ is invariant under $z\mapsto-z$.
In~\eqref{e.fh-sausage-comparison} both integer thresholds differ
from~$u$ by a bounded amount.  The two bounds thus give
$\PP(\sting_{[u,u+r)})=\exp\{-\zeta u(1+o(1))\}$
for all real $u\to\infty$, as required. \qed

\section{Endpoint culs-de-sac}
\label{ss.endpoint-culdesac}

The decomposition of a finite cluster at its extremal renewal points has
two endpoint pieces, a lower cul-de-sac and an inverted cul-de-sac, joined
by a middle string.  We introduce these endpoint pieces before establishing
the factorisation that transfers the sharp string asymptotic to clusters.

\begin{definition}[Cul-de-sac]
\label{d.culdesac-early}
A \hfff{culdesac}{\emph{lower cul-de-sac} with \emph{entry} $x\in\Z^d$} is a finite open
connected component $D$ of $\downspace_{x\cvl}$, with base~$0$, such that
\[
  V(D)\cap\vertbdry_{x\cvl}\, = \, \{x\}\, .
\]
Its head is then~$x$.  Write $\culdesac[x]$ for the event that the open
component of~$0$ in $\downspace_{x\cvl}$ is such a cul-de-sac with entry~$x$,
and $\culdesacset[x]$ for that component.  For $u\geq0$, put
$\culdesac_u=\bigcup_{x\cvl=u}\culdesac[x]$ and use
$\culdesacset_u$ for the corresponding component.

An \emph{upper}, or \emph{inverted}, \emph{cul-de-sac} with entry~$x$ is a
finite open connected component~$P$ of $\upspace_{x\cvl}$ with head~$0$ and
with~$x$ its only vertex on $\vertbdry_{x\cvl}$.  The band condition makes
$x$ the base of~$P$.  Write $\icds[x]$ for the event that the open
component of~$0$ in $\upspace_{x\cvl}$ is such an inverted cul-de-sac
with entry~$x$, and $\icdsset[x]$ for that component.  For $u\geq0$, put
$\icds_u=\bigcup_{x\cvl=-u}\icds[x]$ and use $\icdsset_u$ for the
corresponding component.  Thus the height index~$u$ records the positive
extent from the entry at height~$-u$ to the head at height~$0$.

For each height, the events in either union are pairwise disjoint, since
the corresponding component has a unique entry.

A cul-de-sac is \emph{irreducible} when its entry is its only renewal point;
for a lower cul-de-sac this says
\[
   R(\culdesacset[x])\, = \, \{x\}\, ,
\]
and analogously for an inverted cul-de-sac.

At height zero, the definitions force the entry to be~$0$ and the endpoint
component to be the one-vertex graph~$\{0\}$, whose renewal set is
$\{0\}$.  With $k_0$ as in the endpoint-edge count above, the half-open
edge allocation gives
\begin{equation}
\label{e.zero-culdesac-weights}
 \PP(\culdesac_0)=(1-p)^{d+k_0},
 \qquad
 \PP(\icds_0)=(1-p)^{d-k_0}.
\end{equation}
\end{definition}

Let $A_x$ be the set of edges incident to~$x$ in
$E(\upspace_{x\cvl})$.  By the half-open edge allocation,
$|A_x|=d-k_0$.  Let \(\mathsf G_x\) be the event that
\(\cluster_{x\cvl}\) occurs, \(\rlevel(\mc C)\cap(0,x\cvl)=\varnothing\),
and \(\head(\mc C)=x\).

\begin{lemma}[Cul-de-sac to cluster comparison]
\label{l.cds-cluster-comparison}
For every $x\in\Z^d$,
\[
  \PP\bigl(\culdesac[x],\ R(\culdesacset[x])\, = \, \{x\}\bigr)
\, \le \, (1-p)^{-(d-k_0)}\PP(\mathsf G_x)\, .
\]
\end{lemma}

\noindent{\bf Proof.} On the event on the left, $D=\culdesacset[x]$ is a finite open component
of $\downspace_{x\cvl}$, with base~$0$ and with~$x$ its unique vertex on
$\vertbdry_{x\cvl}$ and with no renewal point at a height in
$(0,x\cvl)$.  If every edge of~$A_x$ is closed, then the full-lattice
cluster of~$0$ equals~$D$ and realises~$\mathsf G_x$.  The edges of~$A_x$
are outside $E(\downspace_{x\cvl})$ and hence independent of the
cul-de-sac event.  These edges are all closed with probability
$(1-p)^{d-k_0}$, giving the result. \qed

The comparison just made has a counterpart for the upper endpoint, stated
separately because reflection through the origin does not preserve the half-open
upspace--downspace edge allocation.

\begin{lemma}[Inverted cul-de-sac to cluster comparison]
\label{l.icds-cluster-comparison}
For every $x\in\Z^d$ with $x\cvl<0$, put $u=-x\cvl$.  Then
\[
\begin{split}
 &\PP\bigl(\icds[x],R(\icdsset[x])=\{x\}\bigr)\\
 &\qquad\le (1-p)^{-(d+k_0)}
 \PP\bigl(\cluster_u,\rlevel(\mc C)\cap(0,u)=\varnothing,
           \head(\mc C)=-x\bigr)\, .
\end{split}
\]
\end{lemma}

\noindent{\bf Proof.} Let $P=\icdsset[x]$.  Since $P\subseteq\upspace_{x\cvl}$ and
$V(P)\cap\vertbdry_{x\cvl}=\{x\}$, every vertex of~$P$ other than~$x$
has height greater than $x\cvl$.  Thus $x$ is its base.  Moreover, only~$x$
is incident to an edge of $E(\downspace_{x\cvl})$ that can leave~$P$.
The $d+k_0$ edges of this edge set incident to~$x$ all leave~$P$;
close them.  They lie outside the edge set defining the inverted
cul-de-sac event, and after their closure the
full-lattice cluster of~$x$ equals~$P$.  Translation by~$-x$ makes this a
cluster based at~$0$, of height~$u$ and with head~$-x$.  Since
$R(P)=\{x\}$, the translated cluster has no renewal level in~$(0,u)$.
The closure of these edges is independent of the inverted cul-de-sac event
and has probability $(1-p)^{d+k_0}$, giving the asserted bound. \qed

We now make the three-piece decomposition precise.

On $\cluster_u\cap\{R(\mc C(0))\ne\varnothing\}$, write $\rmin$ and $\rmax$
for the minimum and maximum elements of $\rlevel(\mc C(0))$.

For $z\in\Z^d$ with $z\cvl>0$, write $0\sconn z$ for the event that the
open component of~$0$ in $\slab_{0,z\cvl}$ is a string with base~$0$ and
head~$z$.

\begin{lemma}[Pinned cul-de-sac--string decomposition]
\label{l.jell.cds}
For $z\in\Z^d$ with $z\cvl>0$, put
\[
\begin{split}
 W^{\base}(z)&\, = \, \PP\bigl(\culdesac[z],R(\culdesacset[z])=\{z\}\bigr)\, ,\\
 G(z)&\, = \, \PP(0\sconn z)\, ,\\
 W^{\head}(z)&\, = \, \PP\bigl(\icds[-z],R(\icdsset[-z])=\{-z\}\bigr)\, .
\end{split}
\]
At zero displacement put
\[
 W^{\base}(0)=(1-p)^{d+k_0},\qquad G(0)=1,\qquad
 W^{\head}(0)=(1-p)^{d-k_0},
\]
in accordance with~\eqref{e.zero-culdesac-weights}; here $G(0)=1$ represents
the one-vertex string at the origin.  At every other displacement of
non-positive height, all three weights are zero.
Let $a,b,x\in\Z^d$ satisfy
$0\le a\cvl\le b\cvl\le x\cvl$.  The probability that the finite cluster
based at~$0$ has head~$x$, lowest renewal point~$a$ and highest renewal
point~$b$ is
\begin{equation}
\label{e.pinned-cds-string}
   W^{\base}(a)\,G(b-a)\,W^{\head}(x-b)\, .
\end{equation}
Consequently, for $u>0$ and $j,k\in[0,u]$ with $j+k\le u$,
\begin{multline*}
  \PP(\cluster_u,\,\rmin\, = \, j,\,\rmax\, = \, u-k) \\
\, = \, \PP(\cds_j,\,\rlevel(\cdsset_j)=\{j\})
  \cdot\PP(\sting_{u-j-k})
  \cdot\PP(\icds_k,\,\rlevel(\icdsset_k)\, = \, \{-k\})\, .
\end{multline*}
The renewal points of the cluster coincide with those of the middle
string, translated to its position between the endpoint pieces.
The identity remains valid when any event determined by the middle string
in a translation-invariant way is imposed on both the left-hand side and
the middle factor on the right.
\end{lemma}

\noindent\textbf{Proof.} Let $C$ be the full cluster on the event in the first
assertion.  The renewal-band condition at~$a$ shows that any self-avoiding
path in~$C$ from~$0$ to a vertex of height at most~$a\cvl$ stays in
$\downspace_{a\cvl}$: an excursion above that height would require two
visits to~$a$.  Thus the component of~$0$ in this downspace is precisely
the part of~$C$ at or below the cut, and is a lower cul-de-sac with
entry~$a$.  Any renewal point of this component other than~$a$ has height
at most $a\cvl-r$; its renewal band contains no vertex of~$C$ outside
the component.  Such a point would therefore be a renewal point of~$C$
below~$a$, contrary to the choice of~$a$.  The lower piece is
irreducible.  The analogous argument from the head~$x$, cutting at~$b$,
gives an irreducible inverted cul-de-sac.  Lemma~\ref{l.cluster-renewal-string}
identifies the intervening piece as a string.  Translate the middle piece
by~$-a$ and the upper piece by~$-x$ to obtain the three origin-rooted
factors in~\eqref{e.pinned-cds-string}.

Conversely, place any three pieces counted by these factors at their
prescribed positions and take their union, using a one-vertex middle
piece when $a=b$.  Their renewal-band conditions ensure that they meet
only at the prescribed junctions and that, within each of the three edge
regions, the union has just the vertices of the corresponding piece.
Each piece is the complete open component in its region, so no open
edge leaves the union.  It is therefore the full-lattice cluster of~$0$,
with base~$0$ and head~$x$.  The junctions $a,b$ are renewal points and, by irreducibility of the
endpoint pieces, are extremal.
Every renewal point of the cluster belongs to the middle string and
remains a renewal point there.  Conversely, by
Lemma~\ref{l.geometrybasics}(2), the renewal band of any middle-string
renewal point other than~$a$ lies above the lower cut and no higher than
the upper cut, where the cluster and middle string have the same vertices.
Thus all middle-string renewal points are renewal points of the cluster.
This also covers the one-vertex middle string when $a=b$.

The occurrence events of the three pieces use disjoint edge sets.  An edge
extending strictly above a renewal cut is assigned to the region above that
cut, while an edge contained in or ending at the cut is assigned below:
$E(\downspace_{a\cvl})$ contains the lower piece,
$E(\slab_{a\cvl,b\cvl})$ the middle piece, and
$E(\upspace_{b\cvl})$ the upper piece.  Their occurrence probabilities
therefore multiply under product percolation, proving the pinned identity.
The same argument is unchanged after restricting the middle factor by a
middle-string event.

To obtain the scalar formula, sum~\eqref{e.pinned-cds-string} over
$a\cvl=j$, $(b-a)\cvl=u-j-k$ and $(x-b)\cvl=k$.  The change of variables
$(a,y,z)=(a,b-a,x-b)$ makes the three sums independent; translation invariance
identifies them with the three displayed height-indexed probabilities. \qed

\begin{proposition}[Exponential decay of irreducible culs-de-sac]
\label{p.estimates.cds}
There exist $C,c>0$ such that, for every real $u\geq0$,
\[
 \sum_{h\in[u,u+r)}
 \PP\bigl(\culdesac_h,\rlevel(\culdesacset_h)=\{h\}\bigr)
 \leq Ce^{-\zeta u-cu}\, .
\]
The same bound holds with summand
$\PP\bigl(\icds_h,\rlevel(\icdsset_h)=\{-h\}\bigr)$.
\end{proposition}
\noindent\textbf{Proof.} For an entry~$x$ of positive height
$h=x\cvl\in[u,u+r)$, Lemma~\ref{l.cds-cluster-comparison} bounds
its irreducible lower-cul-de-sac probability by a fixed multiple of
$\PP(\mathsf G_x)$.  The events $\mathsf G_x$ are disjoint as~$x$
varies, since the full cluster has a unique head.  Each is contained in
\(\cluster_{[u,u+3r]}\cap
\{\rlevel(\mc C)\cap(0,u)=\varnothing\}\).
Summing over entries and applying
Lemma~\ref{l.no-renewal-cluster-upper} proves the lower estimate, apart from a possible
zero-height term, which is absorbed by increasing~$C$.
For the upper estimate, use Lemma~\ref{l.icds-cluster-comparison} and sum over entries
$x$ with $-x\cvl\in[u,u+r)$.  The events on the right of that lemma are disjoint as~$x$ varies,
since the cluster based at~$0$ has the uniquely specified head~$-x$.
Each lies in the same cluster event used for the lower estimate.
The zero-height term is again absorbed into~$C$. \qed

\textbf{Proof of Theorem~\ref{t.stingwelltrap.cluster}.}
First suppose that $\vecell$ is lattice, with span~$s$.  For
$u\in\{0\}\cup s\N$, put
\[
\begin{aligned}
  \alpha_u&\, = \, e^{\zeta u}
    \PP\bigl(\culdesac_u,\rlevel(\culdesacset_u)=\{u\}\bigr)\, ,\\
  \beta_u&\, = \, e^{\zeta u}
    \PP\bigl(\icds_u,\rlevel(\icdsset_u)=\{-u\}\bigr)\, .
\end{aligned}
\]
Proposition~\ref{p.estimates.cds}, together with the zero-height
probabilities in~\eqref{e.zero-culdesac-weights}, gives
$\max\{\alpha_u,\beta_u\}\leq Ce^{-cu}$ for these~$u$.
Set $\gamma_i=e^{\zeta i}\PP(\sting_i)$ for $i\in s\N$, with
$\gamma_0=1$.  Theorem~\ref{t.stingwelltrap.sting} gives
$\gamma_i=\kappasting+O(e^{-ci})$.
The convolution
$\phi_u=\sum_{j+k\leq u}\alpha_j\beta_k\gamma_{u-j-k}$
satisfies $\PP(\cluster_u,R(\mc C(0))\neq\emptyset)
=e^{-\zeta u}\phi_u$ by Lemma~\ref{l.jell.cds}.  Put
$\phi_\infty=\kappasting(\sum_j\alpha_j)\sum_k\beta_k$.
The exponential bounds on the endpoint weights and the convergence of
$\gamma_i$ give
\[
 |\phi_u-\phi_\infty|
 \leq \sum_{j+k\leq u}\alpha_j\beta_k
       |\gamma_{u-j-k}-\kappasting|
       +\kappasting\sum_{j+k>u}\alpha_j\beta_k
 \leq Ce^{-cu}\, .
\]
Since $\cluster_u\cap\{R(\mc C(0))=\emptyset\}$ is contained in the event of
Lemma~\ref{l.no-renewal-cluster-upper}, we obtain
$\PP(\cluster_u)\approxvecell\kappacluster e^{-\zeta u}$ with
$\kappacluster=\phi_\infty$.

Now suppose that $\vecell$ is non-lattice.  Let $A^-$ and $A^+$ be the
atomic measures on $[0,\infty)$ obtained by assigning to an attainable
positive endpoint extent $h$ the respective masses
\[
 e^{\zeta h}\PP\bigl(\cds_h,\rlevel(\cdsset_h)\, = \, \{h\}\bigr)\, ,
 \qquad
 e^{\zeta h}\PP\bigl(\icds_h,\rlevel(\icdsset_h)\, = \, \{-h\}\bigr)\, ,
\]
and by assigning masses $(1-p)^{d+k_0}$ and $(1-p)^{d-k_0}$ at zero to
the lower and upper measures, respectively.
Proposition~\ref{p.estimates.cds} gives both measures exponential moments.
Summing the pinned decomposition of Lemma~\ref{l.jell.cds} over heads
and heights at least~$u$ gives
\[
\begin{aligned}
 &e^{\zeta u}\PP\bigl(\cluster_{[u,\infty)},
     R(\mc C(0))\ne\varnothing\bigr)\\
 &\qquad = \int A^-(dj)A^+(dk)\,
   e^{\zeta(u-j-k)}
   \mathsf M_{\rm str}([\max\{u-j-k,0\},\infty)).
\end{aligned}
\]
Here $\mathsf M_{\rm str}(\{0\})=1$, accounting for clusters whose lowest
and highest renewal points coincide.  For each fixed $j,k$, the last two
factors converge to $\mu_\infty/\zeta$ by Lemma~\ref{l.zetarenewal}(5).
For $v\geq0$, the string--renewal identity and the separation of successive
renewal levels by at least~$r$ give
\[
 e^{\zeta v}\mathsf M_{\rm str}([v,\infty))
 = \int_{[0,\infty)}e^{-\zeta t}\,U(v+dt)
 \leq \sum_{m\geq0}e^{-\zeta mr}
 = (1-e^{-\zeta r})^{-1}.
\]
For $v<0$, the expression
$e^{\zeta v}\mathsf M_{\rm str}([0,\infty))$ is bounded by the same
constant, by taking $v=0$ in the preceding bound.  Thus the integrand is
uniformly bounded, and the finiteness of $A^-$ and $A^+$ permits dominated
convergence, yielding the limit
\[
  \mu_\infty\zeta^{-1}
  A^-([0,\infty))A^+([0,\infty))\in(0,\infty)\, .
\]
Finally, summing Lemma~\ref{l.no-renewal-cluster-upper} over successive
height windows makes clusters with no renewal contribute
$O(e^{-\zeta u-cu})$.  This proves the non-lattice assertion and identifies
its constant $\kappacluster$. \qed

\section{Transverse rarity for finite clusters and sausages}
\label{ss.transverse-rarity}

Lemma~\ref{l.bigwell} supplies an exponential cost for a large external
boundary of the cluster's $K$-box family.  We show that a large diameter
forces such a boundary, obtaining transverse confinement bounds for finite
clusters and, by closing endpoint edges, for sausages.  The cluster bound
localises endpoint culs-de-sac in the vector renewal analysis; the sausage
bound controls transverse fluctuations and supplies the exponential moment
needed there.

Let $\pi$ denote orthogonal projection onto $\vecell^\perp$, so that
$\pi(x)=x-(x\cvl)\vecell$.  For a finite connected subgraph~$G$ of
$(\Z^d,\sim)$, put
\[
  \operatorname{tsp}(G)
  \, = \, \max_{x\in V(G)}|\pi(x-\base(G))|\, .
\]
This is the \emph{transverse span} of $G$ relative to its base.

\begin{lemma}[Diameter forces large coarse exterior boundary]
\label{l.tsp-to-boundary}
For every finite connected subgraph~$G$ of $(\Z^d,\sim)$, the quantity
$2K\sqrt d\,|\partial_{\rm ext}\overline G|$ is at least the diameter
of~$G$, and hence also its transverse span.
\end{lemma}

\noindent\textbf{Proof.} Recall from Definition~\ref{d.exterior} that
$\overline G$ is the family of $K$-boxes containing edges of~$G$, and
$\partial_{\rm ext}\overline G$ is its coarse exterior boundary.
For a finite non-empty $*$-connected family~$\mathcal F$ of $K$-boxes,
let $v_j(B)$ denote the $j$\textsuperscript{th} component of the index
of~$B$ and put
$\diam_j(\mathcal F)=\max_{B\in\mathcal F}v_j(B)
-\min_{B\in\mathcal F}v_j(B)$.
We claim that, for each coordinate~$j$,
$|\partial_{\rm ext}\mathcal F|\geq\diam_j(\mathcal F)+1$.
Indeed, choose $j'\ne j$.  By $*$-connectivity, every integer~$k$
between the extreme $j$-indices $k_{\min}$ and $k_{\max}$ is attained.
For each such~$k$, choose a box~$B_k^*\in\mathcal F$ with
$j$-index~$k$ and maximal $j'$-index among boxes at that level.
The neighbouring box $\widehat B_k=B_k^*+2Ke_{j'}$ lies in the coarse
exterior boundary: the ray continuing from it in the positive
$e_{j'}$-direction avoids~$\mathcal F$.
These boxes have distinct $j$-indices, proving the claim.

If $G$ has no edges, its diameter is zero and the assertion is immediate.
Otherwise, $\overline G$ is non-empty and $*$-connected, since boxes
containing consecutive edges of a path in~$G$ are $*$-neighbours.
Choose $x,y\in V(G)$ with $|x-y|=\diam(G)$ and a component~$j$ for which
$|x_j-y_j|\geq\diam(G)/\sqrt d$.
Both vertices are incident to edges of~$G$, so there are boxes in
$\overline G$ with indices~$v,w$ that contain~$x,y$, respectively.
Their $j$-indices satisfy
\[
  |v_j-w_j|\,\geq\,\frac{|x_j-y_j|}{2K}-1
  \,\geq\,\frac{\diam(G)}{2K\sqrt d}-1\, .
\]
Combining this with the claimed bound for $\mathcal F=\overline G$
proves the lemma. \qed

Recall that $r=e_d\cvl$ under the axes-ordering convention.

We obtain an exponential penalty for transverse span exceeding a 
large multiple of cluster height.

\begin{proposition}[Transverse rarity for finite clusters]
\label{p.cluster-transverse-rarity}
There exist constants $A, C, c > 0$ such that, for every $u \ge 0$ and
$t \ge 0$,
\[
  \PP\bigl(\cluster_{[u,u+3r]},\;
           \operatorname{tsp}(\mc{C})\, \ge \, A u + t\bigr)
\, \le \,
  C\,e^{-(\zeta+c)u}\,e^{-ct}\, .
\]
\end{proposition}

\noindent{\bf Proof.} Let $c_*>0$ be the exponential-rate constant in
the first assertion of Lemma~\ref{l.bigwell}, and put
$a=c_*/(2K\sqrt d)>0$.  Choose $A$ so large that $aA>\zeta$.

On~$\cluster_{[u,u+3r]}\cap\{\operatorname{tsp}(\mc C)\ge Au+t\}$,
Lemma~\ref{l.tsp-to-boundary} gives
$|\partial_{\rm ext}\overline{\mc C}|\ge m$, where
$m=(Au+t)/(2K\sqrt d)$.  When $m\geq1$,
\[
  \diam(\mc C)\geq\operatorname{tsp}(\mc C)
  \geq Au+t=2K\sqrt d\,m\geq K\, .
\]
The first assertion of Lemma~\ref{l.bigwell} therefore yields
\[
  \PP\bigl(\cluster_{[u,u+3r]},\;
           \operatorname{tsp}(\mc{C})\geq Au+t\bigr)
  \leq C\exp\{-a(Au+t)\}\, .
\]
For $m<1$, the threshold $Au+t$ is bounded, so the estimate follows after increasing~$C$.

Take $c=\tfrac12\min\{a,aA-\zeta\}>0$.  Then
$a(Au+t)\ge\zeta u+c(u+t)$, yielding the sought bound for all $u,t\geq 0$. \qed

The corresponding sausage estimate follows by closing endpoint edges.

\begin{corollary}[Transverse fluctuations of sausages]
\label{c.ordinary-sausage-transverse-fluctuation}
There exist constants $A,C,c,\theta>0$ such that, for every $u,t\geq0$,
\[
  \mu\bigl(\hgt(S)\in[u,u+r),\;
            \operatorname{tsp}(S)\geq Au+t\bigr)
  \,\leq\,Ce^{-c(u+t)}\, .
\]
Moreover, the sausage increment has an exponential moment:
\[
  \E_\mu\exp\bigl\{\theta|\head(S)-\base(S)|\bigr\}<\infty\, .
\]
\end{corollary}

\noindent\textbf{Proof.} Let $\mathcal A$ be the family of sausages in
 the event in the first display.  Lemma~\ref{l.string-endpoint-closure},
 together with $e^{\zeta\hgt(S)}\leq e^{\zeta(u+r)}$ on~$\mathcal A$,
 and Proposition~\ref{p.cluster-transverse-rarity} give
\[
\begin{aligned}
  \mu\bigl(\hgt(S)\in[u,u+r),\;\operatorname{tsp}(S)\geq Au+t\bigr)
  &\leq (1-p)^{-2d}e^{\zeta(u+r)}\\
  &\quad\times\PP\bigl(\cluster_{[u,u+3r]},\;
          \operatorname{tsp}(\mc C)\geq Au+t\bigr)
  \leq Ce^{-c(u+t)}\, .
\end{aligned}
\]
To obtain an unconditional transverse tail, fix $v\geq0$ and
partition sausage heights into windows $[kr,(k+1)r)$, $k\in\N_0$.
For $kr\leq v/(2A)$, the event $\operatorname{tsp}(S)\geq v$
implies $\operatorname{tsp}(S)\geq Akr+v/2$, so summing the first
bound over these windows gives $Ce^{-cv/2}$.  For the remaining windows,
$\hgt(S)>v/(2A)$; Lemma~\ref{l.sausage-tilted-tail} gives an
exponentially small contribution in~$v$.  Hence
$\mu(\operatorname{tsp}(S)\geq v)\leq Ce^{-c'v}$ for some $c'>0$.
The same lemma gives an exponential height tail, and
$|\head(S)-\base(S)|\leq\hgt(S)+\operatorname{tsp}(S)$.
The exponential moment follows by Cauchy--Schwarz, after choosing
$\theta>0$ sufficiently small. \qed

\chapter{Renewal abundance and Ornstein--Zernike consequences}
\label{s.renewal-abundance-oz}
Harnessing the renewal description developed in Chapter~\ref{s.asp}, we now
establish renewal abundance and derive the vector Ornstein--Zernike
consequences of the theory.  We begin by deriving our strongest renewal
assertion, Theorem~\ref{t.renewal}: under the conditioned cluster law,
linearly many renewal points are present with overwhelming probability.
The mass gap for the tilted sausage law~$\mu$ makes sparse renewal
exponentially unlikely for strings; the cul-de-sac--string decomposition
and sharp point-to-plane cluster asymptotic carry this conclusion over to
full clusters.

After the proof of renewal abundance, the next two sections pass from
scalar point-to-plane information to the full vector Ornstein--Zernike
picture.  In Section~\ref{ss.radial-identity}, the radial-identity argument
identifies $\zeta(\vecell)\vecell$ as the radial boundary point of the
finite-connectivity Wulff body.  In Section~\ref{ss.oz-consequences}, we pass
to the full vector increments of sausages and endpoint pieces, drawing on the
transverse estimates of Chapter~\ref{s.asp}.  For each fixed height direction
$\vecell$, the vector renewal law has a typical direction of travel
$\drift(\vecell)$, which generally differs from~$\vecell$.  The analysis identifies $\drift(\vecell)$ as the outward
normal at $\zeta(\vecell)\vecell$ and shows that the drift map parametrises
all directions.  It also yields the multivariate local-renewal asymptotic
together with the global analytic and strictly convex Wulff and equidecay
geometry.

The fourth and final section turns to point-to-point finite connectivity.
For a fixed height direction, the extreme-pair event gives the simplest
endpoint realisation of the renewal decomposition: two irreducible endpoint
pieces joined by the middle string, and hence a companion sharp asymptotic.
For unrestricted finite connectivity, the prescribed vertices may instead
lie inside their cluster.  A two-mark re-rooting then replaces the two bare
endpoint factors by marked endpoint factors, while leaving the middle renewal
kernel unchanged.  This yields the asymptotic for each fixed height direction.  Since the drift
map reaches every direction of travel, these local asymptotics combine to give
the global Ornstein--Zernike asymptotic of Theorem~\ref{t.oz-intro},
with a positive real-analytic amplitude.

\section{Renewal abundance from the mass gap}
\label{ss.renewal-abundance-proof}

Recall that $r=\max\{|e_i\cvl|:i\in\intint{d}\}$, equal to $e_d\cvl$ under the
axes-ordering convention; by Lemma~\ref{l.geometrybasics}(2), distinct renewal
levels in a finite connected subgraph differ by at least~$r$.

Let \(\mu\) be the tilted sausage law on \(\sausageset\),
\(\mu(S)=e^{\zeta\hgt(S)}\PP(\mathsf O(S))\).  Proposition~\ref{p.muprob}
says that \(\mu\) is a probability measure.  Under \(\mu^{\otimes\N}\),
write \hfff{renewal-measure}{\(T_n=Z_1+\cdots+Z_n\)}, \(n\geq1\), for the
successive positive renewal levels, with $T_0=0$.  Recall that the renewal measure is
\[
 U(A)=\E_{\mu^{\otimes\N}}\bigl[|\{T_n:n\geq0\}\cap A|\bigr]
\]
for Borel $A\subseteq[0,\infty)$.  For every attainable height~$h$,
Lemma~\ref{l.zetarenewal}(1) gives
\(e^{\zeta h}\PP(\sting_h)=U(\{h\})\); in the non-lattice case we use this
identity after summing over windows of width~$r$.

We introduce two measures $\eta^-$ and $\eta^+$ on $[0,\infty)$, recording
the exponentially tilted height weights of the lower and upper endpoint
pieces in the pinned cluster decomposition.  They are atomic measures, with
masses at attainable positive heights given by
\[
\begin{split}
 \eta^-(\{h\})
   &\, = \, e^{\zeta h}\PP\bigl(\cds_h,\rlevel(\cdsset_h)=\{h\}\bigr)\, ,\\
 \eta^+(\{h\})
   &\, = \, e^{\zeta h}\PP\bigl(\icds_h,\rlevel(\icdsset_h)=\{-h\}\bigr)\, .
\end{split}
\]
Here the inverted object is normalised with head at~$0$ and entry at
height~$-h$; translating its entry to~$0$ gives the equivalent upper
endpoint piece of positive extent~$h$ used in the cluster decomposition.
Their respective masses at zero are $(1-p)^{d+k_0}$ and $(1-p)^{d-k_0}$.
Thus the measures $\eta^\pm$ unify the endpoint weights denoted by
$\alpha_u,\beta_u$ in the lattice case and by $A^\pm$ in the non-lattice
case in the preceding chapter.  The cul-de-sac estimates obtained there
show that both measures are finite and have exponential moments.

\begin{lemma}[Clusters with fewer than two renewal points]
\label{l.fewer-two-renewals}
There exist $C,c>0$ such that, for every $u\geq0$,
\[
 \PP\bigl(\cluster_{[u,u+3r]},\,|R(\mc C)|<2\bigr)
\, \leq \, Ce^{-\zeta u-cu}\, .
\]
\end{lemma}

\noindent{\bf Proof.} If the cluster has no renewal point, Lemma~\ref{l.no-renewal-cluster-upper}
gives the required estimate.

If it has exactly one renewal point, the pinned decomposition
Lemma~\ref{l.jell.cds} has a degenerate middle string.  Writing $j$ and $k$ for the two endpoint extents and summing over
cluster heights in $[u,u+3r]$, we obtain from the product factorisation
\[
 \PP\bigl(\cluster_{[u,u+3r]},\,|R(\mc C)|=1\bigr)
 \, \leq \, e^{-\zeta u}(\eta^-\otimes\eta^+)
 \bigl\{(j,k):j+k\geq u\bigr\}\, .
\]
The exponential moments of the two endpoint measures make the last factor at
most $Ce^{-cu}$. \qed

\begin{lemma}[Sparse-renewal strings]
\label{l.sparse-renewal-strings}
There exist $a,C,c>0$ such that, for all large~$v$,
\[
 \sum_{h\in[v,v+r)}e^{\zeta h}
 \PP\bigl(\sting_h,\,|R(\stingset_h)|<av\bigr)
 \, \leq \, Ce^{-cv}\, ,
\]
where the sum is over attainable heights.
\end{lemma}

\noindent{\bf Proof.} Put $\beta=\E_\mu\hgt(S)$ and choose $a>0$ so small that $2a\beta<1$.  For
$n=\lceil av\rceil$, the exponential moment in
Theorem~\ref{t.mass-gap} gives, for sufficiently small $\lambda>0$,
$\log\E_\mu e^{\lambda\hgt(S)}\leq\beta\lambda+K\lambda^2$
for some $K<\infty$.  The exponential Markov inequality then yields
\[
 \begin{aligned}
 \mu^{\otimes\N}(T_n>v)
 &\leq \exp\{-\lambda v+n(\beta\lambda+K\lambda^2)\}\\
 &\leq C\exp\{-\lambda v(1-a\beta-aK\lambda)\}
 \leq Ce^{-cv}\, ,
 \end{aligned}
\]
where $\lambda$ is chosen small enough that $1-a\beta-aK\lambda>0$.

If a string of height $h\geq v$ has fewer than $av$ renewal points,
then fewer than $n$ positive renewal levels have occurred by~$h$, and hence
$T_n>h\geq v$.  By Lemma~\ref{l.zetarenewal}(2), summing the tilted
string weights over $h\in[v,v+r)$ counts at most one renewal level from
any chain, since distinct levels are separated by at least~$r$.  Thus
\[
 \sum_{h\in[v,v+r)}e^{\zeta h}\PP\bigl(\sting_h,\,|R(\stingset_h)|<av\bigr)
 \, \leq \, \mu^{\otimes\N}(T_n>v)
 \, \leq \, Ce^{-cv}\, .
\]
\qed

\noindent{\bf Proof of Theorem~\ref{t.renewal}.}
Choose $c_0>0$ so that
$c_0<a/2$, where $a$ is supplied by
Lemma~\ref{l.sparse-renewal-strings}.  We estimate directly the joint
probability of $\clusteruthree$ and fewer than $c_0u$ renewal points; this
avoids conditioning on an exact cluster height in the non-lattice case.

On a cluster with at least two renewal points, let $j$ and $k$ be the
extents of the lower and upper endpoint culs-de-sac, let $H\in[u,u+3r]$ be
the cluster height, and put $m=H-j-k$ for the height of the middle string.
By Lemma~\ref{l.jell.cds}, the cluster renewal points are exactly
those of this middle string.

First consider clusters for which $j>u/4$ or
$k>u/4$.  By Lemma~\ref{l.zetarenewal}(1), the tilted middle-string
mass is the renewal measure~$U$.  Since successive renewal levels are at
least~$r$ apart, the supremum of $U([x,x+3r])$ over $x\in\R$ is at most four.  The pinned decomposition
and $H\geq u$ therefore give
\[
 \PP\bigl(\clusteruthree,\,|R(\mc C)|\geq2,\,j>u/4\bigr)
 \, \leq \, Ce^{-\zeta u}\eta^-([u/4,\infty))\eta^+([0,\infty))
 \, \leq \, Ce^{-\zeta u-cu}\, ,
\]
and the same bound holds with $j$ and $k$ interchanged.  Together with
Lemma~\ref{l.fewer-two-renewals}, these estimates reduce the proof to the case
$j,k\leq u/4$.

In this case, $m\geq u/2$.  Define the tilted sparse-renewal measure on
attainable heights by
\[
 B_u(\{m\})=e^{\zeta m}\PP(\sting_m,|R(\stingset_m)|<c_0u).
\]
Partition $[u/2,\infty)$ into successive height windows of width~$r$.
Because $c_0<a/2$, the threshold $c_0u$ is at most $av$ at the left
endpoint~$v$ of every such window.  For all large~$u$,
Lemma~\ref{l.sparse-renewal-strings}, followed by summation over the windows,
gives \(B_u([u/2,\infty))\leq Ce^{-cu}\).
Applying Lemma~\ref{l.jell.cds} to the middle event, and using
$H\geq u$ again, we obtain
\[
\begin{aligned}
 &\PP\bigl(\clusteruthree,\,j,k\leq u/4,\,2\leq|R(\mc C)|<c_0u\bigr)\\
 &\qquad\leq e^{-\zeta u}\eta^-([0,\infty))\eta^+([0,\infty))
   B_u([u/2,\infty))\leq Ce^{-\zeta u-cu}\, .
\end{aligned}
\]
Combining the preceding estimates gives
\[
 \PP\bigl(\clusteruthree,\,|R(\mc C)|<c_0u\bigr)
\, \leq \, Ce^{-\zeta u-cu}\, .
\]
Theorem~\ref{t.stingwelltrap.cluster} implies
$\PP(\clusteruthree)\geq ce^{-\zeta u}$ for all large~$u$.  In the lattice
case, the height span~$s$ satisfies $s\leq r$, so $[u,u+3r]$ contains
an attainable height at which its exact-height asymptotic applies.
In the non-lattice case, subtract the cumulative tail above $u+3r$
from that at~$u$.  Bounding the open upper tail by the closed one,
the asymptotics give
$\liminf_{u\to\infty}e^{\zeta u}\PP(\clusteruthree)
\geq\kappacluster(1-e^{-3\zeta r})>0$.
Division by the resulting lower bound proves
\[
   \Q\bigl(|R(\csetuthree)|<c_0u\bigr)\, \leq \, Ce^{-cu}\, ,
\]
which is Theorem~\ref{t.renewal}. \qed

\section{The Wulff body and the radial identity}
\label{ss.radial-identity}

\begin{figure}[htbp]
\centering
\includegraphics[width=0.90\textwidth]{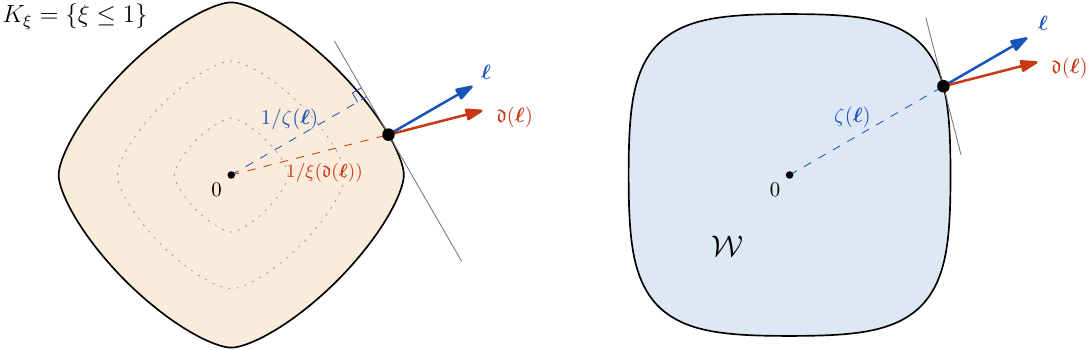}
\caption{The exponential rates of decay for forging a finite point-to-point or
point-to-plane connection in direction $\vecell\in S^{d-1}$ are
$\xi(\vecell)$ and $\zeta(\vecell)$.  The point-to-plane connection drifts in
the direction $\drift(\vecell)$, with $\vecell$ the normal to the
tangent plane to the equidecay surface $\partial K_\xi=\{\xi=1\}$ at the
point of intersection with the ray along $\drift(\vecell)$.
In our string analysis, this drift arises as the direction of the mean
increment under the associated exponentially tilted sausage law. In the Wulff
body~$\mc{W}$, the distance to the surface in direction $\vecell$ is
$\zeta(\vecell)$, the reciprocal of the distance from the origin to this
tangent plane; as we will see in the next section, the roles of $\vecell$
and $\drift(\vecell)$ at the point of tangency are reversed.
(Figure prepared with Claude.)}
\label{f.Wulffduality}
\end{figure}

To the finite-connectivity function
$\tau^{\rm f}(x)=\PP(0\fconn x)$ of Section~\ref{ss.oz-intro},
we associate the exponential-moment sum
\[
  \chi^{\rm f}(t)\, = \, \sum_{x\in\Z^d}e^{t\cdot x}\tau^{\rm f}(x)\, .
\]
The \hfff{moment-body}{moment domain}
$\mathcal D$, on which this sum is finite, and its closure, the Wulff
body~$\mathcal W$, connect the geometry of finite connectivity to the
point-to-plane decay rate~$\zeta(\vecell)$.  We write
\[
  \mathcal D\, = \, \{t\in\R^d:\chi^{\rm f}(t)<\infty\}\, ,
  \qquad \mathcal W\, = \, \overline{\mathcal D}\, .
\]
The set $\mathcal D$ is convex and symmetric.  It is bounded: for each
coordinate direction, prescribing a straight open path of length~$n$ and
closing its other incident edges gives
$\tau^{\rm f}(\pm ne_i)\ge e^{-Cn}$, so every $t\in\mathcal D$ satisfies
$|t\cdot e_i|\le C$.  The geometry of $\mathcal W$ is also
described by its support function and the associated unit ball:
\[
  \xi(x)\, = \, \sup_{t\in\mathcal W}t\cdot x\, ,
  \qquad K_\xi\, = \, \{x:\xi(x)\, \le \, 1\}\, .
\]
In this section, we prove that the boundary of $\mathcal W$ meets the ray
in direction~$\vecell$ at $\zeta(\vecell)\vecell$.  This radial identity is
depicted in the right-hand sketch in Figure~\ref{f.Wulffduality}. It establishes
that $0\in\operatorname{int}\mathcal D$, so that $\mathcal W$ is a compact
convex body and $\xi$ is finite and positive away from the origin.

We prepare for the radial identity with a volume-moment estimate that allows
us to pass from the point-to-plane theory to the moment sum.

\begin{lemma}[Volume moments under point-to-plane conditioning]
\label{l.point-plane-volume-moments}
Here $|\mc C|$ denotes the number of vertices in the cluster.
For every $q\ge1$ there is $C_q<\infty$ such that, for $u\ge0$,
\[
 \E\left[|\mc C|^q\, ;
          \cluster_{[u,u+3r]}\right]
\, \le \, C_q(1+u)^{qd/(d-1)}
       \PP\bigl(\cluster_{[u,u+3r]}\bigr)\, .
\]
\end{lemma}

\noindent{\bf Proof.} We combine the global surface bound of
Proposition~\ref{p.nobigcluster} with lattice isoperimetry to control
the cluster volume.  The lattice isoperimetric inequality, obtained from the discrete
Loomis--Whitney bound~\cite{LoomisWhitney},
applied in the $K$-box index lattice after filling finite complementary
components, gives
\(|C|\le C_K(1+|\partial_{\rm ext}\overline C|)^{d/(d-1)}\)
for every finite connected subgraph~$C$.  Indeed, the exterior
nearest-neighbour boundary of the filled box family is contained in
$\partial_{\rm ext}\overline C$, and each $K$-box contains at most $C_K$
vertices.

Put \(B=|\partial_{\rm ext}\overline{\mc C}|\) and
\(F_u=\cluster_{[u,u+3r]}\).
For $u\ge K$, Proposition~\ref{p.nobigcluster} gives a conditional
exponential tail for $B/u$, above a fixed threshold.  For $0\le u<K$,
a cluster of diameter less than~$K$ has $B\le M_K$, for a fixed
constant $M_K$.  Thus, for $m>M_K$,
the event $\{B\ge m\}$ forces diameter at least~$K$, and the first
assertion of Lemma~\ref{l.bigwell} gives an unconditional exponential
tail for $B$ on finite clusters.  Moreover,
$\inf_{0\le u<K}\PP(F_u)>0$: a bounded-length path in a preferred
coordinate direction can reach the window $[u,u+3r]$, and closing
its other incident edges witnesses $F_u$.  Dividing the unconditional
tail by this lower bound gives a conditional exponential tail in this
range as well.  Thus, for all $u,t\ge0$,
\[
 \PP\bigl(B\ge t(1+u)\mid F_u\bigr)\le Ce^{-ct}\, .
\]
Integrating this bound and applying the isoperimetric inequality proves
the assertion. \qed

\begin{proposition}[Radial identity]
\label{p.radial-identity}
For every $\vecell\in S^{d-1}$,
\[
  \sup\{s\, \ge \, 0:s\vecell\in\mathcal W\}\, = \, \zeta(\vecell)\, .
\]
In particular, $\zeta(\vecell)\vecell\in\partial\mathcal W$.
\end{proposition}

\noindent{\bf Proof.} Let $\mathscr C_0$ be the family of finite connected subgraphs of
$\Z^d$ whose base is the origin, and write $\mathsf O_{\rm full}(\Gamma)$
for the event that the full-lattice cluster of the origin is exactly~$\Gamma$.
Translation invariance and uniqueness of the base give the exact
re-rooting identity
\begin{equation}
\label{e.susceptibility-two-marks}
 \chi^{\rm f}(t)
\, = \, \sum_{\Gamma\in\mathscr C_0}\PP\bigl(\mathsf O_{\rm full}(\Gamma)\bigr)
   \sum_{a,b\in V(\Gamma)}e^{t\cdot(b-a)}\, .
\end{equation}
To see this, start on the left with the finite cluster of~$0$ and a marked
vertex $x$ in it.  Translate its base to~$0$.  The old origin and~$x$ become
an ordered pair $(a,b)$ in a unique member of $\mathscr C_0$, with
$x=b-a$.  The construction is reversible.

For every $a,b\in V(\Gamma)$, the height difference $(b-a)\cvl$ is at most
the height of~$\Gamma$.  Thus, for $s\ge0$, the two-mark identity gives
\[
 \chi^{\rm f}(s\vecell)
 \le \sum_{\Gamma\in\mathscr C_0}\PP\bigl(\mathsf O_{\rm full}(\Gamma)\bigr)
 |V(\Gamma)|^2 e^{s\,\hgt(\Gamma)}\, .
\]
Covering the possible heights by the intervals $[nr,(n+3)r]$ and
applying Lemma~\ref{l.point-plane-volume-moments} with $q=2$, we obtain
\[
 \chi^{\rm f}(s\vecell)
 \le C\sum_{n\ge0}e^{s(n+3)r}
 \E\bigl[|\mc C|^2;\cluster_{[nr,(n+3)r]}\bigr].
\]
By Lemma~\ref{l.first-order-cluster-upper}, for every $\varepsilon>0$
the right-hand side is at most
\[
 C_{s,\varepsilon}\sum_{n\ge0}(1+n)^{2d/(d-1)}
 e^{-(\zeta(\vecell)-\varepsilon-s)nr}.
\]
For $0\le s<\zeta(\vecell)$, choose
$0<\varepsilon<\zeta(\vecell)-s$.  The series converges, so
$s\vecell\in\mathcal D$.

To prove divergence above $\zeta(\vecell)$, we bound the susceptibility
below by the contribution from the heads of clusters that are strings.
By Lemma~\ref{l.string-endpoint-closure}, for $s>\zeta(\vecell)$,
\[
\begin{aligned}
 \chi^{\rm f}(s\vecell)
 &\ge e^{su}\PP\bigl(\mc C\text{ is a string of height in }[u,u+r)\bigr)\\
 &= (1-p)^{2d}e^{su}\PP\bigl(\sting_{[u,u+r)}\bigr).
\end{aligned}
\]
Proposition~\ref{p.zeta} shows that the last expression is
$\exp\{(s-\zeta(\vecell))u+o(u)\}$ as $u\to\infty$.
Thus $\chi^{\rm f}(s\vecell)=\infty$ whenever $s>\zeta(\vecell)$.

For each coordinate direction $v=\pm e_i$, choose
$0<s_v<\zeta(v)$.  The finiteness result proved above gives
$s_vv\in\mathcal D$.  By convexity, $\mathcal D$ contains the convex
hull of these $2d$ points, which has the origin in its interior.

Now suppose that $s\vecell\in\mathcal W$ for some
$s>\zeta(\vecell)$.  Since $\mathcal D$ is convex and contains the
origin in its interior, $s'\vecell\in\mathcal D$ whenever
$0\le s'<s$.  Choosing $\zeta(\vecell)<s'<s$ contradicts the
divergence established above.  Conversely, finiteness of the susceptibility
for every $0\le s<\zeta(\vecell)$ implies that
$\zeta(\vecell)\vecell\in\mathcal W$ by closure.  This proves the
radial identity. \qed

Proposition~\ref{p.radial-identity} identifies the boundary of the Wulff
body with the points $\zeta(\vecell)\vecell$, as $\vecell$ varies over the unit sphere.
Since this body is compact, convex and contains the origin in its interior,
$\zeta$ is continuous and bounded above and away from zero.  Convex duality also identifies $\zeta(\vecell)$ as the minimum of
$\xi(x)$ over vectors $x$ satisfying $x\cdot\vecell=1$.

\section{Vector renewal and Wulff geometry}
\label{ss.oz-consequences}

The radial identity describes the Wulff body through point-to-plane
decay.  We now turn to directional finite connectivity, using the full
displacements of sausages.  First, we study the vector increment law and
show how exponential tilting aligns its mean displacement with a prescribed
direction near the original drift.  In the next two subsections, estimates
for endpoint and marked pieces allow us to compare the convergence of the
sausage renewal series with that of the finite-connectivity susceptibility.
A local renewal theorem then gives the asymptotic for string connectivity,
and the comparison identifies its exponential decay rate with that of
finite connectivity and with the support function $\xi$ of the Wulff body.
Finally, varying the height direction extends these conclusions to all
directions of travel and establishes real analyticity and positive curvature
of the Wulff and equidecay surfaces.

\subsection{Vector increments and conjugate tilts}

Fix a height direction $\vecell_0$ and put
$\zeta=\zeta(\vecell_0)$ and $t_0=\zeta\vecell_0$.
We write $\mu$ for the associated sausage law, leaving its dependence
on $\vecell_0$ implicit.

For a sausage \(S\), define its \hfff{vector-increment}{vector increment}
by \(X(S)=\head(S)-\base(S)\in\Z^d\).
For this height direction, let $\pi$ be the orthogonal projection onto
$\vecell_0^\perp$, so that $\pi(z)=z-(z\cvlzero)\vecell_0$.
Write $\pi(S)=\pi(X(S))$, giving
\(X(S)=\hgt(S)\vecell_0+\pi(S)\), and use
\(\operatorname{tsp}(S)\) for the transverse span introduced in
Section~\ref{ss.transverse-rarity}.
The mass gap controls the longitudinal coordinate, while
Corollary~\ref{c.ordinary-sausage-transverse-fluctuation} supplies
exponential control of transverse fluctuations and, in particular, an
exponential moment for the full vector increment.
Thus \(\E_\mu X(S)\) is finite, with strictly positive component
\(\E_\mu\hgt(S)\) in the height direction.  The associated travel direction,
or \emph{drift}, is
\[
   \drift(\vecell_0)
   \, = \, \frac{\E_\mu X(S)}{|\E_\mu X(S)|}\, .
\]
For independent \(\mu\)-sausages concatenated end-to-end, the direction
of the cumulative displacement converges almost surely to this drift.

We first record the regularity facts about the increment laws needed for the
Ornstein--Zernike analysis.

The increment lattice \(\mathcal G_\mu\) is the subgroup of \(\Z^d\)
generated by displacements of finite strings with positive occurrence
probability, or equivalently by increments of sausages of positive
\(\mu\)-mass.
We also set \(Z(t)=\E_\mu e^{t\cdot X(S)}\), and write
\[
  \operatorname{Dom}Z=\{t\in\R^d: Z(t)<\infty\}
\]
for its domain of finiteness.

\begin{lemma}[Non-degeneracy and regularity of the increment law]
\label{l.kernel-regularity}
The following properties hold.
\begin{enumerate}
\item \textup{(Full lattice.)}
The increment lattice is \(\mathcal G_\mu=\Z^d\).

\item \textup{(Non-degenerate tilted covariance.)}
For every \(t\in\operatorname{int}\operatorname{Dom}Z\), the
exponentially tilted weights
\[
  \frac{e^{t\cdot X(S)}\mu(S)}{Z(t)},\qquad S\in\sausageset,
\]
define a probability measure under which the increment \(X(S)\) has
positive-definite covariance.

\item \textup{(Analytic increment transform.)}
There is a complex neighbourhood \(U\) of the origin such that \(Z(t)\) is
finite and analytic for \(t\in U\).
\end{enumerate}
\end{lemma}

\noindent{\bf Proof: (1).} Lemma~\ref{l.high-endpoint-string} shows that every
lattice point $x$ of sufficiently large height is the head of a
positive-probability string based at~$0$.  By
Lemma~\ref{l.sausage-concatenation}, its displacement~$x$ is a sum of
sausage increments.  Given any $a\in\Z^d$, choose $x$ so
that both $x\cvlzero$ and $(x+a)\cvlzero$ exceed~$3r$.  Then $x$ and $x+a$
belong to $\mathcal G_\mu$, as does their difference~$a$.  Thus
$\mathcal G_\mu=\Z^d$.  The renewal sum allows any number of
sausages, so the lattice generated by their increments is the relevant
one; individual convolution powers may still be periodic.

{\bf (2).} Fix \(v\in\R^d\setminus\{0\}\) and a
coordinate direction \(e_i\) with \(v\cdot e_i\ne0\).  It suffices to produce
two sausages \(S\) and \(S'\) of positive \(\mu\)-mass with
\(v\cdot(X(S')-X(S))\ne0\), so \(v\cdot X(S)\) is not \(\mu\)-almost
surely constant.  Any finite sausage graph has positive \(\mu\)-mass: prescribing
its edges open and all other edges of its slab incident to it closed
realises its occurrence event, so
\(\mu(S)=e^{\zeta\hgt(S)}\PP(\mathsf O(S))>0\).  It remains to exhibit
two such graphs whose increments have distinct scalar products with~\(v\).

When $i=d$, compare the one-edge sausage from $0$ to $e_d$ with the
graph formed by the path $0,e_d,2e_d$ and one side edge in direction
$e_j$, for some $j<d$.  Attach this edge at $0$ if $e_j\cvlzero>0$, and at
$e_d$ otherwise.  Its extra vertex shares the band $\vertbdry_r$ with
$e_d$ and, at its own height, shares a band with the base or with $e_d$.
Thus neither intermediate vertex is a renewal point, while the base and
head remain renewal points.  The second graph is a sausage of increment
$2e_d$, and the two increments differ by $e_d$.  For $i<d$, write
\(r_i=e_i\cvlzero\in[0,r]\), and recall from
Lemma~\ref{l.geometrybasics}(1)
that \(\vertbdry_t=\{w:w\cvlzero\in(t-r,t]\}\).  When \(0<r_i<r\), take
\(S=\{0,e_d\}\) and let \(S'\) be the path \(0,e_i,e_i+e_d\): the intermediate
vertex \(e_i\) sits at height \(r_i\), and \(0\) lies in \(\vertbdry_{r_i}\) as
well, because \(0\in(r_i-r,r_i]\), so \(e_i\) is not a renewal point; the
increments are \(X(S)=e_d\) and \(X(S')=e_i+e_d\).  When \(r_i=0\), take
\(S=\{0,e_d,2e_d,e_d+e_i\}\), with head \(2e_d\), and let \(S'\) be the path
\(0,e_d,e_d+e_i,2e_d+e_i\), with head \(2e_d+e_i\); the vertices \(e_d\) and
\(e_d+e_i\) both have height \(r\) and so share \(\vertbdry_r\), leaving no
internal renewal point.  When \(r_i=r\), take \(S=\{0,e_d\}\) and let \(S'\) be
the unit square \(\{0,e_d,e_i,e_d+e_i\}\), with head \(e_d+e_i\); here \(e_d\)
and \(e_i\) both have height \(r\) and share \(\vertbdry_r\).  In each case
\(X(S')=X(S)+e_i\), whose scalar product with~\(v\) is non-zero.
For any \(t\in\operatorname{int}\operatorname{Dom}Z\), the weights
in part~(2) sum to one, have finite second moments and give positive mass
to exactly the same sausages as~\(\mu\).  Thus the covariance remains
positive definite under each such tilt.

{\bf (3).} Corollary~\ref{c.ordinary-sausage-transverse-fluctuation} gives
\(\E_\mu e^{\theta |X(S)|}<\infty\) for some \(\theta>0\).  Therefore \(Z(t)=\E_\mu e^{t\cdot X(S)}\) converges
absolutely and locally uniformly for \(t\) in a complex neighbourhood of the
origin, and is analytic there. \qed

We use the vector generating function associated with the sausage law:
\begin{equation}
\label{e.lambda-vec}
   \Lambda_{\mathrm{vec}}(t)\, = \, \sum_{S\in\sausageset}e^{t\cdot\head(S)}\,
   \PP\bigl(\mathsf O(S)\bigr)\, .
\end{equation}
This is the vector version of the scalar height series $\Lambda$
appearing in the proof of Proposition~\ref{p.muprob}.  The height
tilt defining~$\mu$ gives
$\Lambda_{\mathrm{vec}}(t_0+s)=Z(s)$, and in particular
$\Lambda_{\mathrm{vec}}(t_0)=\Lambda(\zeta)=1$.
The string sum $\sum_{m\ge1}\Lambda_{\mathrm{vec}}(t)^m$ converges
exactly when $\Lambda_{\mathrm{vec}}(t)<1$.

At a point $t$ of the critical surface
$\{\Lambda_{\mathrm{vec}}=1\}$ near~$t_0$, the weights
$e^{t\cdot X(S)}\PP(\mathsf O(S))$ sum to one and define a probability law
on sausages.  At $t_0$ this law is $\mu$.  Its mean increment is
$\nabla\Lambda_{\mathrm{vec}}(t)$, which points in the outward normal
direction to the surface.  We call the
locally determined tilt for which this mean points along~$\hat n$ the
\emph{conjugate tilt} for~$\hat n$.

\begin{lemma}[Conjugate tilts near the drift]
\label{l.conjugate-tilt-interior}
For the fixed height direction~$\vecell_0$, there are an open
neighbourhood $\mathcal N=\mathcal N(\vecell_0)\subseteq S^{d-1}$ of
\(\drift(\vecell_0)\) and a real-analytic map $\hat n\mapsto t_{\mathrm{str}}(\hat n)$ on
$\mathcal N$ such that
\[
  \Lambda_{\mathrm{vec}}\bigl(t_{\mathrm{str}}(\hat n)\bigr)\, = \, 1\, ,
  \qquad
  t_{\mathrm{str}}(\hat n)-t_0
       \in\operatorname{int}\operatorname{Dom}Z\, .
\]
The outer normal to the critical surface
$\{\Lambda_{\mathrm{vec}}=1\}$ at
$t_{\mathrm{str}}(\hat n)$ is~$\hat n$.  In particular
$t_{\mathrm{str}}(\drift(\vecell_0))=t_0$.
\end{lemma}

\noindent{\bf Proof.} By Lemma~\ref{l.kernel-regularity}(3), $Z$ is analytic
near the origin, with $Z(0)=1$.  Its gradient at the origin is
$\E_\mu X(S)$, which is non-zero because
$X(S)\cdot\vecell_0=\hgt(S)>0$.  The level set $\{Z=1\}$ is
therefore an analytic hypersurface near the origin.

The Hessian of $\log Z$ is the covariance matrix of the corresponding
exponentially tilted increment law.  By Lemma~\ref{l.kernel-regularity}(2),
it is positive definite throughout $\operatorname{int}\operatorname{Dom}Z$.
The second fundamental form of $\{Z=1\}$ is therefore positive definite
on its tangent space.  Its Gauss map is locally an analytic diffeomorphism,
with value $\drift(\vecell_0)$ at the origin.  Inverting this map on a
neighbourhood of that direction and translating by~$t_0$ gives
$t_{\mathrm{str}}$ and the neighbourhood~$\mathcal N$. \qed

For the directions in $\mathcal N$ supplied by the lemma, call
\[
   \bigl\{t_{\mathrm{str}}(\hat n):\hat n\in\mathcal N\bigr\}
\]
the \emph{critical patch}.  As a sum of non-negative exponentials,
$\Lambda_{\mathrm{vec}}$ is convex and lower semicontinuous.  Its unit
sublevel set
\(\mathcal W_{\mathrm{str}}=\{t:\Lambda_{\mathrm{vec}}(t)\le1\}\)
is therefore closed and convex.  Unlike the compact Wulff body~$\mathcal W$,
this auxiliary set is unbounded: since every sausage has positive height,
$t\in\mathcal W_{\mathrm{str}}$ implies
$t-a\vecell_0\in\mathcal W_{\mathrm{str}}$ for every $a\ge0$.
Lemma~\ref{l.local-domain-identification} will identify the two sets
locally near~$t_0$.

For $\hat n\in\mathcal N$, put $t^*=t_{\mathrm{str}}(\hat n)$.
Convexity of $\Lambda_{\mathrm{vec}}$ gives, for every
$t\in\mathcal W_{\mathrm{str}}$,
\[
  1\, \ge\, \Lambda_{\mathrm{vec}}(t)
  \, \ge\, \Lambda_{\mathrm{vec}}(t^*)
   +\nabla\Lambda_{\mathrm{vec}}(t^*)\cdot(t-t^*)
  \, =\, 1+\nabla\Lambda_{\mathrm{vec}}(t^*)\cdot(t-t^*).
\]
Thus the hyperplane through~$t^*$ with normal~$\hat n$ supports the
whole set.  For non-zero~$x$ with $\hat x\in\mathcal N$, define the
\emph{string inverse correlation length} by
\[
   \xi_{\mathrm{str}}(x)=t_{\mathrm{str}}(\hat x)\cdot x.
\]
Equivalently, for these $x$,
\[
   \xi_{\mathrm{str}}(x)
   =\sup_{t\in\mathcal W_{\mathrm{str}}}t\cdot x<\infty.
\]
For a string of displacement~$x$ in this cone, undoing the conjugate tilt
contributes the factor $e^{-\xi_{\mathrm{str}}(x)}$.  The IID
local-renewal theorem, Theorem~\ref{t.iid-local-renewal}, will give the
accompanying power-law prefactor.

\subsection{Endpoint and marked-piece localisation}

In the renewal decomposition of a finite cluster, endpoint pieces flank
the string.  We first control their occurrence weights and transforms,
estimates that will also serve in the extreme-pair asymptotic of
Section~\ref{s.point-to-point-oz}.

\begin{lemma}[Endpoint cul-de-sac deficit and transform convergence]
\label{l.cds-deficit}
Let $A$ be the transverse slope of Proposition~\ref{p.cluster-transverse-rarity}.
\begin{enumerate}
\item \textup{(Scalar deficit.)}  There exist $C,c>0$ such that, for every
$z\in\Z^d$ with $z\cvlzero>0$,
\begin{equation}
\label{e.cds-deficit-cone}
   \PP\bigl(\culdesac[z],R(\culdesacset[z])\, = \, \{z\}\bigr)
\, \le \,
   C\,e^{-\zeta z\cvlzero}\,e^{-c\bigl(z\cvlzero+(|\pi(z)|-Az\cvlzero)_+\bigr)}\, .
\end{equation}
The same bound holds for \(W^{\head}(z)\).  In particular, both
endpoint weights are at most \(Ce^{-(\zeta+c)z\cvlzero}\).
\item \textup{(Transform convergence.)}  There exists $\delta>0$ such
that the endpoint transforms
\[
   \widehat W^{\base}(t)=\sum_z e^{t\cdot z}W^{\base}(z),
   \qquad
   \widehat W^{\head}(t)=\sum_z e^{t\cdot z}W^{\head}(z)
\]
converge absolutely and are analytic for
$t\in\mathbb C^d$ satisfying $|\operatorname{Re}t-t_0|<\delta$,
with $\operatorname{Re}t$ taken componentwise.
\end{enumerate}
\end{lemma}

\noindent{\bf Proof: (1).} By Lemma~\ref{l.cds-cluster-comparison},
\(W^{\base}(z)\le(1-p)^{-(d-k_0)}\PP(\mathsf G_z)\).
Proposition~\ref{p.estimates.cds}, applied at height $z\cvlzero$, gives
\(W^{\base}(z)\le Ce^{-(\zeta+c)z\cvlzero}\).  If $|\pi(z)|>Az\cvlzero$,
the condition $\head=z$ forces
$\operatorname{tsp}(\mc C)\ge |\pi(z)|$, and
Proposition~\ref{p.cluster-transverse-rarity} supplies the transverse factor
in~\eqref{e.cds-deficit-cone}.  After decreasing~$c$, these two estimates give
the asserted bound.  For the inverted event, apply
Lemma~\ref{l.icds-cluster-comparison} with entry~$-z$.  The resulting cluster,
after translation, has head~$z$, so the same scalar and transverse estimates
give the asserted bound.

{\bf (2).} By Lemma~\ref{l.jell.cds}, the zero-displacement endpoint
weights are finite and all other weights at non-positive height vanish.
For $z\cvlzero>0$, put $h=z\cvlzero$, so that
$z=h\vecell_0+\pi(z)$.  For $t=t_0+s\in\mathbb C^d$, part~(1) gives
\[
 \bigl|e^{t\cdot z}W^{\base}(z)\bigr|
\, \le \, C\exp\!\left\{
   |\operatorname{Re}s|\bigl(h+|\pi(z)|\bigr)
   -c\bigl(h+(|\pi(z)|-Ah)_+\bigr)\right\}\, .
\]
Choose $\delta>0$ so small that
$\delta(2+A)<c/4$.  For $|\operatorname{Re}s|<\delta$,
split the sum over~$z$ according as
$|\pi(z)|\le2Ah$ or $|\pi(z)|>2Ah$.  In the first region the displayed
weight is at most $Ce^{-c'h}$; in the second it is at most
$Ce^{-c'(h+|\pi(z)|)}$, for some $c'>0$.  The number of lattice points in
each unit height and transverse-radius range is polynomial, and is absorbed by
the exponential decay.  Thus the series for $\widehat W^{\base}$ converges
absolutely and locally uniformly on
$\{t\in\mathbb C^d:|\operatorname{Re}t-t_0|<\delta\}$ and is analytic there.
The same argument with the inverted estimate in part~(1) gives the conclusion
for~$\widehat W^{\head}$. \qed

Shrink $\mathcal N$ from Lemma~\ref{l.conjugate-tilt-interior}, if
necessary, so that $|t_{\mathrm{str}}(\hat n)-t_0|<\delta/2$ and
$\hat n\cvlzero>0$ for every $\hat n\in\mathcal N$.  Both endpoint
transforms are then analytic on a neighbourhood of the critical patch.
In the remaining local analysis for this fixed height direction~$\vecell_0$,
we continue to write $\mathcal N$ for such a neighbourhood of
$\drift(\vecell_0)$, reducing it further when a later argument requires
this and always keeping $\drift(\vecell_0)$ in its interior.

The re-rooted representation of the susceptibility
in~\eqref{e.susceptibility-two-marks} involves an ordered pair of cluster
vertices.  Cutting at renewal points separates the exponential weight of
the marks into contributions from their pieces and the sausages between
them.  We introduce the corresponding marked-piece weights.  For a sausage
or an irreducible lower or upper endpoint piece~$Q$, write $q^-$ and $q^+$
for its southern and northern anchors: for a sausage, its base and head;
for a lower endpoint piece, its base and entry; and for an upper endpoint
piece, its entry and head.  Put
\[
   V^\circ(Q)=V(Q)\setminus\{q^+\}.
\]
For an irreducible upper endpoint piece put instead $V^\circ(Q)=V(Q)$.
Equivalently, in a concatenated record a renewal vertex common to two
adjacent pieces is assigned to the northern piece, while an endpoint belonging
to only one piece remains with that piece.  Thus the lower endpoint and every
sausage are half-open at their northern anchor, whereas the upper endpoint
retains the cluster head.  For $t\in\R^d$, put
\[
 \begin{split}
  L_t(Q)&\, = \, \sum_{a\in V^\circ(Q)}e^{t\cdot(q^+-a)},\qquad
  R_t(Q)\, = \, \sum_{b\in V^\circ(Q)}e^{t\cdot(b-q^-)}\, ,\\
  D_t(Q)&\, = \, \sum_{a,b\in V^\circ(Q)}e^{t\cdot(b-a)}\, .
 \end{split}
\]
When the piece containing the first mark precedes that containing the
second, $L_t$ and $R_t$ sum the contributions from the two marked pieces,
while the intervening sausages contribute their full increment weights.
When both marks lie in one piece, their joint contribution is~$D_t$.
Reversing the order of the marked pieces gives the same factorisation
with $t$ replaced by~$-t$.  Clusters with fewer than two renewal points
have no string between distinct extremal renewal points and will be
controlled separately.

\begin{lemma}[Summability of marked pieces]
\label{l.marked-piece-summability}
Fix the height direction $\vecell_0$ and put
$t_0=\zeta(\vecell_0)\vecell_0$.  There are $\eta>0$ and a neighbourhood
$U_0$ of~$t_0$ with the following properties.
\begin{enumerate}
\item Under the sausage law,
\[
 \E_\mu\left[(1+|V(S)|)^2e^{\eta\operatorname{diam}(S)}\right]<\infty\, .
\]
Consequently, for each $F\in\{L,R,D\}$, the series
$\sum_S\PP(\mathsf O(S))F_t(S)$ is bounded uniformly for
$t\in U_0$; the same holds with $t$ replaced by~$-t$.
\item If $w(Q)$ denotes the occurrence probability of an irreducible lower
or upper endpoint piece and $X(Q)=q^+-q^-$, then
\[
 \sum_Q w(Q)e^{t_0\cdot X(Q)}
       (1+|V(Q)|)^2e^{\eta\operatorname{diam}(Q)}<\infty\, .
\]
For each $F\in\{L,R,D\}$, the sum $\sum_Q w(Q)F_t(Q)$, over
irreducible endpoint pieces of either type, is bounded uniformly for
$t\in U_0$; the same holds with $t$ replaced by~$-t$.
\item After decreasing $\eta$ and $U_0$ if necessary, the two-mark
contribution from finite clusters having fewer than two renewal points obeys,
locally uniformly for real $t\in U_0$,
\[
 \sum_{\substack{\Gamma\in\mathscr C_0\\ |R(\Gamma)|<2}}
   \PP\bigl(\mathsf O_{\rm full}(\Gamma)\bigr)
   e^{\eta\operatorname{diam}(\Gamma)}
   \sum_{a,b\in V(\Gamma)}e^{t\cdot(b-a)}<\infty .
\]
In particular these clusters make a finite real-analytic contribution to
$\chi^{\rm f}(t)$ on~$U_0$.
\end{enumerate}
Each of these series also defines an analytic function on a complex
neighbourhood of~$t_0$.
\end{lemma}

\noindent{\bf Proof.} For a sausage, its Euclidean diameter is at most a
constant multiple of $1+\hgt(S)+\operatorname{tsp}(S)$.  Since a connected subset of~$\Z^d$ of
diameter~$R$ has at most $C(1+R)^d$ vertices, Theorem~\ref{t.mass-gap} and
Corollary~\ref{c.ordinary-sausage-transverse-fluctuation} give the first
display for sufficiently small~$\eta$.

Every vertex of a piece lies between the heights of its anchors.  Hence, at
$t=t_0$, each of \(t_0\cdot(q^+-a)\), \(t_0\cdot(b-q^-)\), and
\(t_0\cdot(b-a)\) is at most \(t_0\cdot X(Q)\).  For a sausage,
$\PP(\mathsf O(S))e^{t_0\cdot X(S)}=\mu(S)$.  Perturbing $t_0$ to~$t$
multiplies any of the three marked weights by at most
$e^{|t-t_0|\operatorname{diam}(S)}$.  The first display therefore proves
the asserted local boundedness.  The same argument applies with $t$ replaced
by $-t$, since the absolute height difference between any two vertices of a
sausage is bounded by its height.

For part~\textup{(2)}, partition endpoint heights into windows
$[u,u+r)$, where $u=kr$ and $k\in\N_0$.  For each window, separate the
pieces with $\operatorname{tsp}(Q)<Au+1$ from those with
$Au+j\le\operatorname{tsp}(Q)<Au+j+1$, $j\ge1$.
Proposition~\ref{p.estimates.cds} bounds the total occurrence weight in
the first class by $Ce^{-\zeta u-cu}$.  For $j\ge1$,
Lemmas~\ref{l.cds-cluster-comparison} and
\ref{l.icds-cluster-comparison} compare irreducible endpoint pieces
with full finite clusters by closing at most $2d$ edges incident to their
entries, outside the defining restricted graphs.  Summing over entries
and applying Proposition~\ref{p.cluster-transverse-rarity} gives total
occurrence weight at most $Ce^{-\zeta u-c(u+j)}$ in the $j$\textsuperscript{th}
class, after reducing~$c$ if necessary.  In each class,
$\operatorname{diam}(Q)\le C(1+u+j)$ and $|V(Q)|$ is bounded by a
polynomial in $1+u+j$.  Since
$e^{t_0\cdot X(Q)}=e^{\zeta\hgt(Q)}\le e^{\zeta(u+r)}$,
the exponential deficit makes the resulting double sum over $k,j$
finite for sufficiently small~$\eta$.  The anchor inequalities and
the same perturbation argument for both signs of $t$ give the marked
assertions in part~\textup{(2)}.

For part~\textup{(3)}, partition cluster heights into intervals
$[u,u+r)$ with $u=kr$, $k\in\N_0$, and put
\[
 \mathcal T_{u,0}=\{\operatorname{tsp}(\mc C)<Au+1\},
 \qquad
 \mathcal T_{u,j}
 =\{Au+j\le \operatorname{tsp}(\mc C)<Au+j+1\},
 \quad j\ge1.
\]
Theorem~\ref{t.stingwelltrap.cluster} gives
$\PP(\cluster_{[u,u+3r]})\le Ce^{-\zeta u}$: in the non-lattice
case by the tail asymptotic, and in the lattice case by summing over the
bounded number of attainable heights in the window.  For $j\ge1$, two applications of Cauchy--Schwarz give
\begin{align*}
 &\E\bigl[|\mc C|^2;
      \cluster_{[u,u+3r]},\ |R(\mc C)|<2,\ \mathcal T_{u,j}\bigr]\\
 &\quad\le
 \E\bigl[|\mc C|^4;\cluster_{[u,u+3r]}\bigr]^{1/2}
 \Bigl[
   \PP\bigl(\cluster_{[u,u+3r]},|R(\mc C)|<2\bigr)
   \PP\bigl(\cluster_{[u,u+3r]},\mathcal T_{u,j}\bigr)
 \Bigr]^{1/4}\\
 &\quad\le C(1+u)^M e^{-\zeta u}e^{-c(u+j)} .
\end{align*}
Here the last line uses Lemmas~\ref{l.point-plane-volume-moments}
and \ref{l.fewer-two-renewals} and
Proposition~\ref{p.cluster-transverse-rarity}, together with the
point-to-plane bound just recalled.  For $j=0$ the same bound follows
by omitting the transverse-rarity factor and combining the first two
inputs with the point-to-plane bound.

The height-window and transverse-span conditions give
$\operatorname{diam}(\mc C)\le C_0(u+j+1)$.  Moreover, if the cluster height
lies in $[u,u+r)$, then for any two marks
$t_0\cdot(b-a)\le\zeta(u+r)$.  Consequently, for sufficiently small
$\eta>0$,
\[
 \sum_{k,j\ge0} e^{\zeta kr}
   e^{\eta C_0(kr+j+1)}
   C(1+kr)^M e^{-\zeta kr}e^{-c(kr+j)}<\infty .
\]
This proves the displayed exponential-diameter two-mark bound in
part~\textup{(3)} at $t=t_0$.  A perturbation from $t_0$ to~$t$ changes each
marked weight by at most
$e^{|t-t_0|\operatorname{diam}(G)}$; after shrinking $U_0$ the same summable
majorant therefore gives locally uniform absolute convergence for complex
$t$ and hence analyticity. \qed

\subsection{Convergence domains and local renewal asymptotics}

The comparison of convergence domains rests on the preceding marked-piece
estimates and the following decomposition.

\begin{lemma}[Marked renewal decomposition]
\label{l.marked-renewal-decomposition}
Consider triples $(C,a,b)$ in which $C$ is a finite cluster based at the
origin, $|R(C)|\ge2$, and $a,b\in V(C)$ are ordered marks.  Cut~$C$ at its
lowest and highest renewal points and at every renewal point between them.
Assign a vertex common to two adjacent pieces to the northern piece.  Then
\[
 V(C)=V^\circ(Q_-)\,\sqcup\,
      V^\circ(S_1)\,\sqcup\cdots\sqcup\,
      V^\circ(S_m)\,\sqcup\,V^\circ(Q_+),
\]
and each triple has a unique record consisting of the irreducible
endpoint pieces $Q_-,Q_+$, the sausages $S_1,\ldots,S_m$, and, for each
mark, its assigned piece and position within that piece.  Conversely,
concatenating such pieces at their prescribed anchors and restoring the marks
at their recorded positions determines one triple.  Under the edge-allocation
convention of Definition~\ref{d.os}, the occurrence probability of the
record is the product of the occurrence probabilities of its pieces.
\end{lemma}

\noindent{\bf Proof.} The ordered renewal points of~$C$ are intrinsic.
Lemma~\ref{l.jell.cds} determines the endpoint cuts and intervening string;
Lemma~\ref{l.sausage-concatenation} then determines the sausage cuts.
The explicit convention for
$V^\circ$ removes the northern anchor from the lower endpoint and each
sausage, while the upper endpoint retains its head.  Hence the displayed
union is disjoint and exhausts $V(C)$: each shared renewal vertex belongs to
the piece above it, and the cluster head belongs to the upper endpoint's
assigned vertex set.  This convention also covers a lower endpoint of zero
extent, whose assigned vertex set is empty.  In particular every mark has a
unique assigned piece and position within it.  This proves injectivity of the representation.  Concatenating the
pieces at their prescribed anchors and restoring the marked vertices
reverses the construction and proves surjectivity.  Edges extending strictly above a cut are assigned to the piece
above it, while edges contained in or ending at the cut are assigned below.
The piece occurrence events thus depend on disjoint edge sets, so their
probabilities multiply under product percolation. \qed

\begin{lemma}[Local identification of convergence domains]
\label{l.local-domain-identification}
Fix the height direction $\vecell_0$ and put
$t_0=\zeta(\vecell_0)\vecell_0$.  There is a neighbourhood $U$ of~$t_0$
such that
\[
   \mathcal D\cap U
\, = \, \{t\in U:\Lambda_{\mathrm{vec}}(t)<1\}\, ,
   \qquad
   \mathcal W\cap U
\, = \, \{t\in U:\Lambda_{\mathrm{vec}}(t)\, \le \, 1\}\, .
\]
Thus the finite-connectivity susceptibility and the $\vecell_0$-string
renewal series have the same boundary of convergence near~$t_0$.
\end{lemma}

\noindent{\bf Proof.} Choose fixed one-edge lower and upper extensions obtained by adjoining,
respectively, a pendant edge below the string base and above its head and
closing every other edge incident to the two new vertices.  These extensions
also close the possible exit edges at the original string base and head
identified in Lemma~\ref{l.string-endpoint-closure}, except for the two
added edges; there are at most $2d$ such exits altogether.
All prescribed edges lie outside the edge set defining the string,
so the extensions have a fixed positive occurrence probability
$c_{\rm ext}$, independent of the string.  Adjoining them to a string
from $0$ to $x$ produces a finite cluster containing both vertices.
The original string is recovered between the prescribed extensions,
so these full-cluster events are disjoint as the string varies.
The definition of $\chi^{\rm f}$ and the sausage concatenation identity
therefore give, for every real $t$,
\[
 \chi^{\rm f}(t)
\, \ge \, c_{\rm ext}\sum_{m\ge1}\Lambda_{\mathrm{vec}}(t)^m\, ,
 \qquad c_{\rm ext}>0\, ,
\]
with the inequality understood in $[0,\infty]$.  In particular,
$\mathcal D\subseteq\{t:\Lambda_{\mathrm{vec}}(t)<1\}$ for this
height direction.

For the reverse comparison use the exact two-mark identity
\eqref{e.susceptibility-two-marks} and the bijection of
Lemma~\ref{l.marked-renewal-decomposition}.  A cluster with at least two
renewal points has a unique ordered sequence
$Q_-,S_1,\ldots,S_m,Q_+$, its occurrence probability is the product of the
piece probabilities, and each mark is assigned to one piece.

Suppose first that the piece containing~$a$ precedes the piece containing~$b$.
Write $q_a^+$ for the northern anchor of the piece containing~$a$, and
$q_b^-$ for the southern anchor of the piece containing~$b$.  Then
\[
 e^{t\cdot(b-a)}
\, = \, e^{t\cdot(q_a^+-a)}
   \prod_{Q\,\text{strictly between the marks}}e^{t\cdot X(Q)}
   e^{t\cdot(b-q_b^-)}\, .
\]
Every piece strictly between the two marked pieces is a sausage.  Sausages
strictly before the first mark or after the second carry no tilt.  If the
piece containing~$b$ precedes that containing~$a$, the same formula holds
with $t$ replaced by~$-t$; if both marks lie in one piece, its contribution
is $D_t(Q)$.  Lemma~\ref{l.marked-piece-summability} controls the finitely many
possible types of the two marked pieces and of the endpoint pieces.

Now
\[
 \Lambda_{\mathrm{vec}}(0)
\, = \, \E_\mu e^{-t_0\cdot X(S)}<1\, ,
 \qquad
 \Lambda_{\mathrm{vec}}(-t_0)
\, = \, \E_\mu e^{-2t_0\cdot X(S)}<1\, ,
\]
because $t_0\cdot X(S)=\zeta\hgt(S)>0$.  Near $-t_0$ the longitudinal
factor is $e^{-2\zeta\hgt(S)}$, and
Corollary~\ref{c.ordinary-sausage-transverse-fluctuation} makes a small
transverse perturbation integrable.  Hence, after shrinking the
neighbourhood~$U$, Lemma~\ref{l.marked-piece-summability} applies throughout
$U$ and
\(\sup_{t\in U}\Lambda_{\mathrm{vec}}(-t)<1\).
By Lemma~\ref{l.marked-piece-summability}(3), the contribution from
clusters with fewer than two renewal points is uniformly bounded and may
be absorbed into $C_U$.  Summing the preceding factorisation over the
unmarked outer segments and then over the segment between the marks gives,
for $t\in U$ with $\Lambda_{\mathrm{vec}}(t)<1$,
\[
 \chi^{\rm f}(t)
\, \le \, \frac{C_U}{(1-\Lambda_{\mathrm{vec}}(0))^2}
 \left(
   1+\frac1{1-\Lambda_{\mathrm{vec}}(t)}
    +\frac1{1-\Lambda_{\mathrm{vec}}(-t)}
 \right)\, .
\]
In particular $\chi^{\rm f}(t)<\infty$ throughout the strict
sublevel set in~$U$.  Together with the preceding lower bound this proves
the first identity.  After replacing $U$ by a smaller neighbourhood with
closure in the one just used, continuity shows that
$\mathcal W\cap U\subseteq\{t\in U:\Lambda_{\mathrm{vec}}(t)\le1\}$:
indeed, points of $\mathcal W\cap U$ are limits of points of $\mathcal D$
lying in the larger neighbourhood.  Conversely,
$\nabla\Lambda_{\mathrm{vec}}$ does not vanish near~$t_0$, so every point
of the local level surface is a limit of points in the strict sublevel set,
which lies in~$\mathcal D$.  Hence
$\{t\in U:\Lambda_{\mathrm{vec}}(t)\le1\}\subseteq\mathcal W$,
and the second identity follows. \qed

The global lower bound also identifies the susceptibility convergence domain
with the interior of $\mathcal W$.  Since $\mathcal D$ is convex with
non-empty interior and $\mathcal W=\overline{\mathcal D}$, we have
$\operatorname{int}\mathcal W=\operatorname{int}\mathcal D\subseteq
\mathcal D$.  If $t\in\mathcal D\cap\partial\mathcal W$, the radial
identity gives $t=\zeta(\vecell)\vecell$, where
$\vecell=t/|t|$.  For this height direction the sausage generating
function equals one at~$t$, by Proposition~\ref{p.muprob}; the global
lower bound then forces $\chi^{\rm f}(t)=\infty$, a contradiction.
Consequently $\mathcal D=\operatorname{int}\mathcal W$.

We now turn from the convergence domain to the sharp string asymptotic.
The next theorem gives a multivariate local-renewal estimate for a general
lattice increment law.  Its proof uses a Fourier local central limit theorem
for a fixed number of increments and then sums over that number, taking
account of lattice periodicity.  The summation of fixed-increment local limit
estimates also appears in \cite[Section~4.4]{CampaninoIoffe};
an earlier multidimensional renewal local-limit analysis is given in \cite[Section~2]{Ioffe}.

\begin{theorem}[IID multivariate local-renewal theorem]
\label{t.iid-local-renewal}
Let $\nu$ be a probability law on $\Z^d$ such that
$w\cdot z>0$ for every $z\in\operatorname{supp}\nu$, for some
non-zero $w\in\mathbb R^d$, and put $Z_\nu(s)=\E_\nu e^{s\cdot X}$.
Suppose that:
\begin{enumerate}
\item $\E_\nu e^{\theta|X|}<\infty$ for some $\theta>0$;
\item the subgroup generated by $\operatorname{supp}\nu$ is $\Z^d$; and
\item $\operatorname{Cov}_\nu(X)$ is positive definite.
\end{enumerate}
There are a neighbourhood $U_\nu$ of the origin and a neighbourhood
$\mathcal N_\nu$ of the direction $\E_\nu X/|\E_\nu X|$ such that, for
$\hat n\in\mathcal N_\nu$, there is a unique $s(\hat n)\in U_\nu$ satisfying
\[
 Z_\nu(s(\hat n))\, = \, 1\, ,
 \qquad
 \frac{\E_\nu[Xe^{s(\hat n)\cdot X}]}
      {|\E_\nu[Xe^{s(\hat n)\cdot X}]|}\, = \, \hat n\, .
\]
The map $\hat n\mapsto s(\hat n)$ is real analytic.  There is a
positive real-analytic function $\Psi_\nu$ such that
\begin{equation}
\label{e.iid-local-renewal}
 \sum_{m\, \ge \, 1}\nu^{*m}(x)
\, = \,
 \Psi_\nu(\hat x)|x|^{-(d-1)/2}
 e^{-s(\hat x)\cdot x}\bigl(1+o_K(1)\bigr)\, ,
 \qquad \hat x\, = \, x/|x|\, ,
\end{equation}
where, for every compact $K\subset\mathcal N_\nu$, the error $o_K(1)$
tends to zero uniformly over $x$ with $\hat x\in K$ as $|x|\to\infty$.
\end{theorem}

\noindent{\bf Proof.} Write $v_s=\nabla\log Z_\nu(s)$ and
$\Sigma_s=\nabla^2\log Z_\nu(s)$.  The latter is the covariance matrix
under the exponentially tilted law.  Since tilting preserves the support,
item~(3) makes it positive definite throughout the interior of the
moment domain.  Thus $\{Z_\nu=1\}$ is a real-analytic, locally strictly
convex hypersurface near the origin, with outward normal $v_s/|v_s|$.
Its Gauss map is locally an analytic diffeomorphism, giving the unique
analytic choice $s=s(\hat n)$ asserted in the theorem.

Fix a compact $K\subset\mathcal N_\nu$, let $\hat n\in K$, and set
$s=s(\hat n)$.  Put $\nu_s(z)=e^{s\cdot z}\nu(z)$, a probability law
because $Z_\nu(s)=1$.  Then
\[
  \nu^{*m}(x)\, = \, e^{-s\cdot x}\nu_s^{*m}(x)\, .
\]
The sum over~$m$ is restricted by the coset structure of
$\operatorname{supp}\nu^{*m}$, which we now determine.  Let
\(D=\langle z-z':z,z'\in\operatorname{supp}\nu\rangle\) be the
difference lattice.  Positive-definiteness of the covariance makes
$D$ full rank.  If $a\in\operatorname{supp}\nu$, item~(2) says that
$\Z^d/D$ is generated by the class of~$a$; it is therefore a finite cyclic
group.  Denote its order, equivalently the index $[\Z^d:D]$, by~$q_{\rm lat}$.
The support of $\nu_s^{*m}$ is contained in $ma+D$.  Thus, for a fixed~$x$,
a non-zero term $\nu_s^{*m}(x)$ can occur only when~$m$ lies in one residue
class modulo~$q_{\rm lat}$; below we call such orders coset-compatible.

\smallskip
\noindent\emph{Uniform Fourier estimates.}
Choose a compact set $T$ in the interior of the moment domain whose
interior contains $\{s(\hat n):\hat n\in K\}$.  For $h\in T$ write
\[
 \nu_h(z)=\frac{e^{h\cdot z}}{Z_\nu(h)}\nu(z),\qquad
 v_h=\E_{\nu_h}X,\qquad
 \Sigma_h=\operatorname{Cov}_{\nu_h}(X).
\]
The support, and hence the difference lattice~$D$, is independent of~$h$.
Exponential moments hold uniformly in a common neighbourhood of zero,
and the least eigenvalue of $\Sigma_h$ is bounded away from zero on~$T$.
We claim that, uniformly in $h\in T$ and $y\in ma+D$,
\begin{equation}
\label{e.iid-uniform-fourier}
 \nu_h^{*m}(y)
 =\frac{q_{\rm lat}}{(2\pi m)^{d/2}\sqrt{\det\Sigma_h}}
 \exp\!\left\{-\frac{(y-mv_h)^{\mathsf T}\Sigma_h^{-1}(y-mv_h)}{2m}\right\}
 +o_T(m^{-d/2})\, .
\end{equation}
Off the coset $ma+D$ the probability vanishes.  The same argument gives
\begin{equation}
\label{e.iid-uniform-local-upper}
 \sup_{h\in T}\sup_{y\in\Z^d}\nu_h^{*m}(y)
 \le C_Tm^{-d/2}\, .
\end{equation}

To see both estimates, let $\phi_h$ be the characteristic function
of~$\nu_h$.  On the Fourier torus, $|\phi_h|=1$ exactly at the
$q_{\rm lat}$ characters annihilating~$D$, independently of~$h$.
Near the origin, uniformly for $h\in T$,
\[
 \log\phi_h(\theta)
 =i v_h\cdot\theta-\tfrac12\theta^{\mathsf T}\Sigma_h\theta
   +O_T(|\theta|^3).
\]
At each of the other exceptional frequencies the characteristic function
differs from its translate near zero by a fixed phase.  On the complement
of small neighbourhoods of these frequencies, compactness gives
$|\phi_h|\le\rho<1$ uniformly in~$h$; within them, positive-definite
covariance gives a uniform Gaussian bound.  Fourier inversion and the
rescaling $\theta\mapsto m^{-1/2}\theta$ give
\eqref{e.iid-uniform-fourier}, with an additive error uniform over all lattice endpoints~$y$.
The $q_{\rm lat}$ contributions add for $y\in ma+D$ and cancel otherwise.
Estimating the same Fourier integral in absolute value gives
\eqref{e.iid-uniform-local-upper}.  This is the characteristic-function
argument of \cite[Sections~2.2--2.3]{LawlerLimic}, with uniformity here
supplied by the common support lattice and compactness of~$T$.

Now let $v=v_s$, $\Sigma=\Sigma_s$, and, for $x=R\hat n$, put
\[
   m_0=\frac{R}{|v|},
   \qquad
   A=v^{\mathsf T}\Sigma^{-1}v.
\]
Fix $M>0$.  If $|m-m_0|\le M\sqrt{m_0}$, then
$|x-mv|=O_{K,M}(\sqrt m)$.  Thus the additive error in
\eqref{e.iid-uniform-fourier} is a relative error, and for
coset-compatible~$m$ it gives
\begin{equation}
\label{e.iid-fixed-m-llt}
 \nu_s^{*m}(x)
 \, = \, \frac{q_{\rm lat}}{(2\pi m)^{d/2}\sqrt{\det\Sigma}}
   \exp\!\left\{
      -\frac{(x-mv)^{\mathsf T}\Sigma^{-1}(x-mv)}{2m}
        \right\}
   \bigl(1+o_{K,M}(1)\bigr)\, .
\end{equation}
The error is uniform over $\hat n\in K$.  For any other~$m$ the
coset constraint still applies.

Because $x=m_0v$, the exponent in~\eqref{e.iid-fixed-m-llt} is exactly
$-A(m-m_0)^2/(2m)$.  Thus, for fixed~$M$, summing over the compatible
orders in this central window gives, as $R\to\infty$,
\[
 \begin{split}
 &\sum_{\substack{m\ge1:\,|m-m_0|\le M\sqrt{m_0}}}
       \nu_s^{*m}(x)  \\
 &\qquad =\,
 \frac{(2\pi)^{-d/2}|v|^{(d-1)/2}}
      {\sqrt{\det\Sigma}}\,
 R^{-(d-1)/2}
 \int_{-M}^{M} e^{-A u^2/2}\,du
 \,\bigl(1+o_{K,M}(1)\bigr).
 \end{split}
\]
Indeed, the spacing~$q_{\rm lat}$ of coset-compatible orders cancels the
factor~$q_{\rm lat}$ in~\eqref{e.iid-fixed-m-llt}.

We now control the orders outside the central window, beginning with
moderate deviations.  The derivative of the mean map under exponential
tilting is the positive-definite covariance matrix.
Thus we may choose $\delta>0$ so that, uniformly for $\hat n\in K$ and
$|m-m_0|\le\delta m_0$, a further tilt of~$\nu_s$ has mean~$x/m$.
Choose $\delta$ small enough that these further tilts belong to~$T$.
Applying~\eqref{e.iid-uniform-local-upper} under these tilts, and using
the quadratic lower bound on the Cram\'er rate near~$v$, gives
constants $C_K,c_K>0$ such that
\[
   \nu_s^{*m}(x)
   \le C_K m^{-d/2}
      \exp\!\left\{-c_K\frac{(m-m_0)^2}{m_0}\right\}.
\]
Consequently the contribution from
$M\sqrt{m_0}<|m-m_0|\le\delta m_0$, after multiplication by
$m_0^{(d-1)/2}$, tends to zero as $M\to\infty$, uniformly over~$K$.

We fix the cutoffs for the two extreme ranges first.  Choose $\beta>0$
so that
$M_\beta=\sup_{\hat n\in K}\E_{\nu_{s(\hat n)}}
e^{\beta\hat n\cdot X}$ is finite.  The exponential Markov inequality gives
\[
  \nu_s^{*m}(x)\le e^{-\beta R}M_\beta^m.
\]
Since $R=m_0|v_s|$ and $|v_s|$ is bounded away from zero on~$K$,
we may fix $c_-\in(0,1-\delta)$ sufficiently small that the sum over
$m<c_-m_0$ is exponentially small in~$m_0$.

For the other extreme, fix any $b>0$.  Since every increment lies in
$\{z:w\cdot z>0\}$, compactness gives
\[
 q_b=\sup_{\hat n\in K}\E_{\nu_{s(\hat n)}}e^{-bw\cdot X}<1.
\]
Thus $\nu_s^{*m}(x)\le e^{bw\cdot x}q_b^m$.  As
$w\cdot x=m_0|v_s|\,w\cdot\hat n$ and the latter coefficient is
uniformly bounded on~$K$, we may fix $c_+>1+\delta$ sufficiently large
that the sum over $m>c_+m_0$ is exponentially small in~$m_0$.

It remains to consider $c_-m_0\le m\le c_+m_0$ with
$|m-m_0|>\delta m_0$.  Here the required average $x/m$ stays a
uniformly positive distance from~$v_s$, since
$x/m=(m_0/m)v_s$.  The common exponential-moment neighbourhood then
gives a uniform Chernoff bound, making this remaining contribution
exponentially small in~$m_0$.  These estimates control all orders
outside the central window, uniformly over~$K$.

Letting first $R\to\infty$ and then $M\to\infty$, and using
\(\int_{\mathbb R}e^{-Au^2/2}\,du=\sqrt{2\pi/A}\), we obtain
\[
 \sum_{m\,\ge\,1}\nu_s^{*m}(x)
 \, = \,
 \frac{(2\pi)^{-(d-1)/2}|v|^{(d-1)/2}}
      {\sqrt{\det\Sigma}\sqrt A}
 R^{-(d-1)/2}\bigl(1+o_K(1)\bigr)\, .
\]
Multiplying the last estimate by the exponential factor from the tilt
identity proves~\eqref{e.iid-local-renewal}, with
\[
 \Psi_\nu(\hat n)
\, = \, \frac{(2\pi)^{-(d-1)/2}|v_{s(\hat n)}|^{(d-1)/2}}
        {\sqrt{\det\Sigma_{s(\hat n)}}
         \sqrt{v_{s(\hat n)}^{\mathsf T}
                    \Sigma_{s(\hat n)}^{-1}v_{s(\hat n)}}}\, .
\]
This function is positive and real analytic.  All estimates used above are
uniform when $\hat n$ ranges over a compact subset of $\mathcal N_\nu$, which
gives the stated uniform error. \qed

\begin{proposition}[String Ornstein--Zernike asymptotic on a cone]
\label{p.string-oz-cone}
Fix the height direction~$\vecell_0$.  After reducing the neighbourhood
$\mathcal N=\mathcal N(\vecell_0)$ if necessary, there is a positive
real-analytic function~$\Psi^{\mathrm{str}}$ on~$\mathcal N$ such that,
for every compact $K\subset\mathcal N$,
\begin{equation}
\label{e.oz-string-llt}
   \PP\bigl(0\sconn x\bigr)
\, = \,
   \Psi^{\mathrm{str}}(\hat x)\,|x|^{-(d-1)/2}
   e^{-\xi_{\mathrm{str}}(x)}\bigl(1+o_K(1)\bigr)\, ,
   \qquad \hat x=x/|x|\, ,
\end{equation}
where $o_K(1)\to0$ as $|x|\to\infty$, uniformly over
$x\in\Z^d$ with $\hat x\in K$.
\end{proposition}

\noindent{\bf Proof.} Lemma~\ref{l.kernel-regularity}(1)--(2)
gives the lattice and covariance properties of the increment law
under~$\mu$; Corollary~\ref{c.ordinary-sausage-transverse-fluctuation}
supplies its exponential moment.  Thus
Theorem~\ref{t.iid-local-renewal} applies.  Reduce
$\mathcal N$ so that $\mathcal N\subseteq\mathcal N_\mu$ and
$t_{\mathrm{str}}(\hat n)-t_0\in U_\mu$, the two neighbourhoods in that
theorem.  By local uniqueness, its tilt is
$s(\hat n)=t_{\mathrm{str}}(\hat n)-t_0$.

The exact string product identity gives
\(\PP(0\sconn x)=e^{-t_0\cdot x}\sum_{m\ge1}\mu^{*m}(x)\).
Theorem~\ref{t.iid-local-renewal} therefore yields
\eqref{e.oz-string-llt}, with $\Psi^{\mathrm{str}}=\Psi_\mu$ and
$\xi_{\mathrm{str}}(x)=t_{\mathrm{str}}(\hat x)\cdot x$.
In particular, $\xi_{\mathrm{str}}$ is the exponential decay rate of
string connectivity along every fixed lattice direction in this cone.
\qed

The asymptotic~\eqref{e.oz-string-llt} concerns the bare string.  Section~\ref{s.point-to-point-oz}
attaches the endpoint pieces and performs the marked re-rooting needed for
unrestricted finite connectivity.  We first identify the string rate with the
finite-connectivity rate and then globalise the resulting local description
through the drift map.

\begin{lemma}[Rate identification on a cone]
\label{l.fc-string-rate}
Fix $\vecell_0\in S^{d-1}$.  After reducing the neighbourhood
$\mathcal N=\mathcal N(\vecell_0)$ if necessary, the string rate agrees
with the support function~$\xi$ of~\eqref{e.xi-def} for every non-zero
$x$ with $x/|x|\in\mathcal N$.  Moreover, if $x_n\in\Z^d$ satisfies
$|x_n|\to\infty$ and $x_n/|x_n|\to\hat n\in\mathcal N$, then
\[
 -|x_n|^{-1}\log\PP(0\fconn x_n)\longrightarrow\xi(\hat n)\, .
\]
\end{lemma}

\noindent{\bf Proof.} Since $t_{\mathrm{str}}(\drift(\vecell_0))=t_0$
and the conjugate tilt depends continuously on direction, we may
reduce $\mathcal N$ so that all $t_{\mathrm{str}}(\hat n)$ lie in the
neighbourhood where Lemma~\ref{l.local-domain-identification} identifies
the boundaries of $\mathcal W$ and $\mathcal W_{\mathrm{str}}$.
Their common tangent hyperplane at $t_{\mathrm{str}}(\hat n)$ has outer
normal~$\hat n$.  Convexity makes it supporting for both full sets,
so $\xi_{\mathrm{str}}(\hat n)=\xi(\hat n)$; homogeneity gives the
first assertion.

The one-edge extensions in the proof of
Lemma~\ref{l.local-domain-identification} close a string into a finite
cluster with probability at least $c_{\rm ext}>0$.  Both original string
endpoints remain connected, and distinct strings give disjoint cluster
events.  Thus, without shifting the displacement,
\[
 \PP(0\fconn x)\ge c_{\rm ext}\PP(0\sconn x)\, .
\]
For $x=x_n$ as in the statement, the uniform string asymptotic of
Proposition~\ref{p.string-oz-cone} gives
\[
 \limsup_{n\to\infty}-|x_n|^{-1}\log\PP(0\fconn x_n)
 \le\xi_{\mathrm{str}}(\hat n)=\xi(\hat n).
\]
Conversely, for every $t\in\mathcal D$,
\(\PP(0\fconn x_n)\le e^{-t\cdot x_n}\chi^{\rm f}(t)\).
Taking the lower limit and then the supremum over $t\in\mathcal D$
gives the matching bound~$\xi(\hat n)$ by~\eqref{e.xi-def}. \qed

\subsection{Global Wulff geometry and the drift map}

We now vary the height direction.  By the radial identity, every boundary
point of the Wulff body has the form $\zeta(\vecell)\vecell$.  Near this
point, the renewal analysis for height direction~$\vecell$ identifies the
Wulff boundary with the critical surface of the associated sausage
transform.  Thus the local curvature and normal information can be assembled
over the whole boundary.  The drift map assigns to each height direction
the normalised mean sausage increment,
$\drift(\vecell)=\E_{\mu^{(\vecell)}}X/|\E_{\mu^{(\vecell)}}X|$.

\begin{theorem}[Global regularity, polarity and the drift map]
\label{l.drift-onto}
The point-to-plane rate $\zeta$ is real analytic on $S^{d-1}$.  The
boundaries of the Wulff body $\mathcal W$ and of its polar
$K_\xi=\{x:\xi(x)\le1\}$ are real analytic and have positive-definite second
fundamental form.  In particular both bodies are strictly convex.  The radial
map
\[
 \rho(\vecell)\, = \, \zeta(\vecell)\vecell
       \colon S^{d-1}\longrightarrow\partial\mathcal W
\]
and the Gauss map $G\colon\partial\mathcal W\to S^{d-1}$ are real-analytic
diffeomorphisms.  Moreover
\[
 \drift\, = \, G\circ\rho\, .
\]
Consequently the drift map is a real-analytic diffeomorphism of the sphere.
For each direction $\hat n$, its inverse
$\vecell(\hat n)=\drift^{-1}(\hat n)$ is the outward unit normal to
$\partial K_\xi$ at $\hat n/\xi(\hat n)$.  Also
$\bigcup_{\vecell}\mathcal N(\vecell)=S^{d-1}$, and
$\xi(\hat n)=\xi_{\mathrm{str}}^{(\vecell(\hat n))}(\hat n)$ in every
direction, where the superscript records the height direction defining the
string kernel.
\end{theorem}

\noindent{\bf Proof.} Fix $\vecell$, put $t_0=\zeta(\vecell)\vecell$, and let
$\nu_t$ be the probability law on sausages defined, for
$\Lambda_{\mathrm{vec}}(t)<\infty$, by
\(\nu_t(S)=\Lambda_{\mathrm{vec}}(t)^{-1}e^{t\cdot X(S)}\PP(\mathsf O(S))\).
Lemma~\ref{l.local-domain-identification} identifies $\partial\mathcal W$ near
$t_0$ with the level surface $\{\Lambda_{\mathrm{vec}}=1\}$.  On this surface
\[
 \nabla\log\Lambda_{\mathrm{vec}}(t)\, = \, \E_{\nu_t}X\, ,
 \qquad
 \nabla^2\log\Lambda_{\mathrm{vec}}(t)
\, = \, \operatorname{Cov}_{\nu_t}(X)\, .
\]
The gradient is non-zero, since $X\cdot\vecell>0$.  The support of~$\nu_t$ is
the support of the original sausage law, so
Lemma~\ref{l.kernel-regularity}(2) makes its covariance positive definite.
The analytic implicit-function theorem therefore gives an analytic boundary
chart at~$t_0$.  If $v$ is tangent to the level surface at~$t$, then
\begin{equation}
\label{e.second-form-covariance}
 \frac{v^{\mathsf T}\nabla^2\Lambda_{\mathrm{vec}}(t)v}
      {|\nabla\Lambda_{\mathrm{vec}}(t)|}
\, = \, \frac{\operatorname{Var}_{\nu_t}(v\cdot X)}
       {|\nabla\Lambda_{\mathrm{vec}}(t)|}>0\, .
\end{equation}
Indeed, on the tangent space the rank-one term relating
$\nabla^2\Lambda_{\mathrm{vec}}$ to
$\nabla^2\log\Lambda_{\mathrm{vec}}$ vanishes.  The left-hand side of
\eqref{e.second-form-covariance} is the second fundamental form with the
outward orientation.  Thus the chart is positively curved.  At~$t_0$,
$\nu_{t_0}=\mu^{(\vecell)}$, and its outer unit normal is
\begin{equation}
\label{e.drift-as-normal}
 \frac{\nabla\Lambda_{\mathrm{vec}}(t_0)}
      {|\nabla\Lambda_{\mathrm{vec}}(t_0)|}
\, = \, \frac{\E_{\mu^{(\vecell)}}X}
       {|\E_{\mu^{(\vecell)}}X|}
\, = \, \drift(\vecell)\, .
\end{equation}

By Proposition~\ref{p.radial-identity}, every point of
$\partial\mathcal W$ has the form $t_0=\zeta(\vecell)\vecell$.  The preceding
radial parametrisation therefore supplies one of the preceding charts at
every boundary point.  Hence $\partial\mathcal W$ is
real analytic and has positive-definite second fundamental form everywhere.
It cannot contain a non-trivial line segment, since at an interior point of
such a segment the normal curvature in its direction would be zero.  Thus
$\mathcal W$ is strictly convex.  The principal curvatures are continuous and
positive on the compact analytic boundary.  They, and hence the
Gauss--Kronecker curvature, are therefore bounded away from zero.

The radial identity also makes $\rho$ a homeomorphism, its inverse being
$t\mapsto t/|t|$.  To see that $\rho$ is analytic, take an analytic
defining chart $F=0$ at $t=\rho(\vecell)$.  The radial derivative
\(\partial_rF(r\vecell)|_{r=\zeta(\vecell)}=\nabla F(t)\cdot\vecell\)
is non-zero: if $n$ is the outer normal at~$t$, then $n\cdot t>0$
because $0\in\operatorname{int}\mathcal W$, and hence
$n\cdot\vecell>0$.  The analytic implicit-function theorem now makes
$\zeta$ and~$\rho$ analytic locally on the sphere, and the explicit
inverse makes $\rho$ an analytic diffeomorphism.

The Gauss map $G\colon\partial\mathcal W\to S^{d-1}$ is analytic.  Its
differential is the shape operator, which is invertible by
\eqref{e.second-form-covariance}.  Strict convexity gives uniqueness of the
support point with a prescribed outer normal, while compactness gives
existence; thus $G$ is bijective.  The analytic inverse-function theorem
makes it an analytic diffeomorphism.  Equation~\eqref{e.drift-as-normal}
now gives $\drift=G\circ\rho$, and proves all assertions concerning the drift
map and the covering of the sphere by the neighbourhoods $\mathcal N(\vecell)$.
Lemma~\ref{l.fc-string-rate} on these neighbourhoods gives
$\xi(\hat n)=\xi_{\mathrm{str}}^{(\vecell(\hat n))}(\hat n)$ in every
direction.

It remains to record explicitly the polar conclusion.  Since~$\xi$ is the
support function of~$\mathcal W$, $K_\xi$ is its polar body.  For
$n\in S^{d-1}$ let $t(n)=G^{-1}(n)$ and
$h(n)=n\cdot t(n)>0$.  The map \(P(n)=n/h(n)\) is an analytic parametrisation
of $\partial K_\xi$: for every
$s\in\mathcal W$, $s\cdot n\le h(n)$, with equality only at~$t(n)$, and
hence $s\cdot P(n)\le1$ with equality at~$t(n)$.  Every radial direction
$n$ occurs once, so $P$ is an analytic diffeomorphism from the sphere onto
$\partial K_\xi$.  The outer normal to $K_\xi$ at $P(n)=n/\xi(n)$ is
$t(n)/|t(n)|=\vecell(n)=\drift^{-1}(n)$, the reciprocal normal
correspondence illustrated in Figure~\ref{f.Wulffduality}.
Since $n\mapsto\vecell(n)$ is an analytic diffeomorphism, the Gauss map
of $\partial K_\xi$ is analytic with everywhere invertible differential.  Convexity
makes its second fundamental form non-negative; invertibility makes it
positive definite.  This proves the asserted analyticity and strict convexity
of the equidecay body. \qed

Writing $\vecell(\hat n)=\drift^{-1}(\hat n)$, the support point of
$\mathcal W$ with outer normal~$\hat n$ is
$\zeta(\vecell(\hat n))\vecell(\hat n)$.  Hence
\[
  \xi(\hat n)
  =\zeta(\vecell(\hat n))\,\vecell(\hat n)\cdot\hat n,
  \qquad
  \zeta(\hat n)\le\xi(\hat n)\le\zeta(\vecell(\hat n)).
\]
In particular, taking $\hat n=\drift(\vecell)$ gives
$\xi(\drift(\vecell))=\zeta(\vecell)\,\vecell\cdot\drift(\vecell)$.
For the first inequality, equality means that the radial boundary point
$\zeta(\hat n)\hat n$ is also the support point with outer normal~$\hat n$;
by strict convexity this holds exactly when $\drift(\hat n)=\hat n$.
For the second, equality is the Cauchy--Schwarz condition
$\vecell(\hat n)=\hat n$, which gives the same criterion.

\section{The point-to-point Ornstein--Zernike asymptotic}
\label{s.point-to-point-oz}

Here, we prove Theorem~\ref{t.oz-intro}.  The string
asymptotic of the preceding section has a rather direct application when
the prescribed vertices $0$ and $x$ are the base and head of their finite
cluster.  The renewal decomposition then consists of a long middle string
flanked by two endpoint pieces.  We begin with the resulting extreme-pair
asymptotic.

For unrestricted finite connectivity, the prescribed vertices need not be
extremal.  Treating them as marks in the renewal decomposition leads to
modified endpoint kernels, while the middle-string calculation remains
unchanged.  The two asymptotics also admit a comparison that identifies
the contribution of the additional endpoint configurations.

Fix a height direction $\vecell_0$, and let
$\mathcal N=\mathcal N(\vecell_0)$ be a neighbourhood of
$\drift(\vecell_0)$ in $S^{d-1}$ on which the string asymptotic and
rate identification of Section~\ref{ss.oz-consequences} hold.
Throughout the section, base and head are taken in the fixed
$\vecell_0$-height order, with its transverse lexicographic tie-breaker,
while target directions $\hat n$ range over~$\mathcal N$.
For a prescribed target direction~$\hat n$, the canonical choice of
height direction is $\drift^{-1}(\hat n)$, the outward unit normal to
$\partial K_\xi$ at $\hat n/\xi(\hat n)$ by
Theorem~\ref{l.drift-onto}.  These height frames will give the global
asymptotic of Theorem~\ref{t.oz-intro}.

For $x\in\Z^d$ with $x\cvlzero\neq0$, write
$0\fconnext_{\vecell_0}x$ for the event that $\mc C(0)$ is finite
and has $0$ and $x$ as its base and head, in either order.  This event
entails $0\fconn x$.  If $x/|x|\in\mathcal N$, then $x\cvlzero>0$,
so $0$ is the base and $x$ the head.

Both asymptotics will follow from a common convolution estimate for
general boundary kernels.

Recall from Lemma~\ref{l.jell.cds} the base- and head-rooted irreducible
cul-de-sac weights $W^{\base}$ and $W^{\head}$.  Their endpoint transforms are
\[
 \widehat W^{\base}(t)\, = \, \sum_z e^{t\cdot z}W^{\base}(z)\, ,
 \qquad
 \widehat W^{\head}(t)\, = \, \sum_y e^{t\cdot y}W^{\head}(y)\, .
\]

Write
\[
 G(w)=\PP(0\sconn w),\qquad G(0)=1,
 \qquad
 \widetilde G(w)=G(w)\mathbf 1\{w\ne0\}.
\]
Thus $\widetilde G$ is the occurrence kernel of a non-empty string.  The
renewal structure gives a useful bound on this kernel.  If $t$ lies on
the string critical surface, then
\[
 \nu_t(S)=e^{t\cdot X(S)}\PP(\mathsf O(S))
\]
is a probability law on sausages.  For $w\ne0$ the string decomposition gives
\begin{equation}
\label{e.tilted-string-bound}
 e^{t\cdot w}G(w)
 =\sum_{m\ge1}\nu_t^{*m}(w)
 \le \sum_{m\ge1}\PP_{\nu_t}\bigl(T_m=w\cvlzero\bigr)
 \le1,
\end{equation}
where $T_m$ are the strictly increasing renewal levels.  The same bound holds
at $w=0$ by the convention $G(0)=1$, and of course for
$\widetilde G\le G$.

Call a pair $(K^-,K^+)$ of non-negative kernels on $\Z^d$, neither
identically zero, \emph{admissible} on~$\mathcal N$ if their transforms are
analytic on a neighbourhood of the critical patch and, for every compact
$\mathcal K\Subset\mathcal N$, there are $c,C>0$ such that
\begin{equation}
\label{e.admissible-boundary-pair}
 \sum_u e^{t_{\mathrm{str}}(\hat n)\cdot u+c|u|}K^-(u)
 +\sum_v e^{t_{\mathrm{str}}(\hat n)\cdot v+c|v|}K^+(v)
 \le C,
 \qquad \hat n\in\mathcal K.
\end{equation}

\begin{lemma}[Boundary-kernel convolution]
\label{l.boundary-kernel-convolution}
Let $(K^-,K^+)$ be admissible and let $\mathcal G$ be either $G$ or~$\widetilde G$.  Then, uniformly for $\hat x=x/|x|$ in every compact
$\mathcal K\Subset\mathcal N$,
\begin{align*}
 &\sum_{u,v}K^-(u)\mathcal G(x-u-v)K^+(v)\\
 &\qquad =
 \widehat K^-(t_x)\,\Psi^{\mathrm{str}}(\hat x)\,
 \widehat K^+(t_x)\,
 |x|^{-(d-1)/2}e^{-\xi_{\mathrm{str}}(x)}\bigl(1+o_{\mathcal K}(1)\bigr),
\end{align*}
where $t_x=t_{\mathrm{str}}(\hat x)$.
\end{lemma}

\noindent{\bf Proof.}
Fix $\mathcal K\Subset\mathcal N$ and set $L=B\log|x|$.  From
\eqref{e.tilted-string-bound}, for either choice of~$\mathcal G$,
\[
 e^{t_x\cdot x}K^-(u)\mathcal G(x-u-v)K^+(v)
 \le e^{t_x\cdot u}K^-(u)e^{t_x\cdot v}K^+(v).
\]
Hence~\eqref{e.admissible-boundary-pair} makes the total contribution from
$|u|>L$ or $|v|>L$ at most
$Ce^{-t_x\cdot x-cL}$.  Taking $B$ large makes this
$o(|x|^{-(d-1)/2}e^{-\xi_{\mathrm{str}}(x)})$, uniformly on~$\mathcal K$.
Since neither kernel vanishes, both transforms are strictly positive and
bounded away from zero on the compact set of tilts corresponding to~$\mathcal K$.

For $|u|,|v|\le L$, put $\delta=u+v$.  Then
$|\delta|=O(\log|x|)$, and for all sufficiently large~$x$ the vector
$x-\delta$ is non-zero and has direction in a fixed compact subset of
$\mathcal N$.  Thus $\widetilde G(x-\delta)=G(x-\delta)$ and the uniform
string local limit~\eqref{e.oz-string-llt} applies.  Analyticity and
homogeneity of $\xi_{\mathrm{str}}$ give, uniformly for these $u,v$,
\[
 \xi_{\mathrm{str}}(x-\delta)
 =\xi_{\mathrm{str}}(x)-t_x\cdot\delta
   +O\!\left(\frac{(\log|x|)^2}{|x|}\right).
\]
The power-law factor and $\Psi^{\mathrm{str}}$ also change by $1+o(1)$.
Consequently
\[
 \mathcal G(x-u-v)
 =\Psi^{\mathrm{str}}(\hat x)|x|^{-(d-1)/2}
   e^{-\xi_{\mathrm{str}}(x)}e^{t_x\cdot u}e^{t_x\cdot v}
   \bigl(1+o_{\mathcal K}(1)\bigr)
\]
uniformly for $|u|,|v|\le L$.  Summing and then removing the logarithmic
truncation by~\eqref{e.admissible-boundary-pair} gives the asserted formula.
\qed

By Lemma~\ref{l.cds-deficit}(2), the endpoint transforms converge absolutely
in a complex tube around~$t_0$.  For every compact
$\mathcal K\Subset\mathcal N$, this tube also contains
$t_{\mathrm{str}}(\hat n)+c\sigma$ for all $\hat n\in\mathcal K$ and
$\sigma\in\{-1,1\}^d$, provided $c>0$ is small enough.  Since
$|u|\le |u|_1=\max_\sigma\sigma\cdot u$, summing the corresponding
transform bounds gives~\eqref{e.admissible-boundary-pair}.  The positive
zero-displacement weights show that neither kernel vanishes.  Thus
$(W^{\base},W^{\head})$ is admissible, and we obtain the extreme-pair
asymptotic with no further local-limit calculation.

\begin{theorem}[Extreme-pair Ornstein--Zernike asymptotic]
\label{t.oz-rooted}
Let $d\ge3$ and $p\in(p_c,1)$.  For every height direction
$\vecell_0\in S^{d-1}$, there is an open neighbourhood
$\mathcal N=\mathcal N(\vecell_0)$ of $\drift(\vecell_0)$ and a positive
real-analytic function $A_{\mathrm{ext}}^{\vecell_0}$ on~$\mathcal N$
such that, for every compact $\mathcal K\Subset\mathcal N$, uniformly as
$|x|\to\infty$ with $\hat x=x/|x|\in\mathcal K$,
\begin{align*}
   \PP\bigl(0\fconnext_{\vecell_0}x\bigr)
   &\, = \, A_{\mathrm{ext}}^{\vecell_0}(\hat x)\,|x|^{-(d-1)/2}\,
      e^{-\xi(x)}\,\bigl(1+o_{\mathcal K}(1)\bigr)\, , \\
   A_{\mathrm{ext}}^{\vecell_0}(\hat n)
   &\, = \, \widehat W^{\base}\bigl(t_{\mathrm{str}}(\hat n)\bigr)\,
      \Psi^{\mathrm{str}}(\hat n)\,
      \widehat W^{\head}\bigl(t_{\mathrm{str}}(\hat n)\bigr)\, .
\end{align*}
Here $t_{\mathrm{str}}$ and $\Psi^{\mathrm{str}}$ are the conjugate tilt and
string amplitude of Proposition~\ref{p.string-oz-cone} for this height
direction.
\end{theorem}

\noindent{\bf Proof.}
On the part of the extreme-pair event for which the cluster has at least two
renewal points, the pinned decomposition of Lemma~\ref{l.jell.cds} is exact.
The middle string then contains at least one sausage, so
\begin{equation}
\label{e.extreme-pair-convolution}
 \PP\bigl(0\fconnext_{\vecell_0}x\bigr)
 =\bigl[W^{\base}*\widetilde G*W^{\head}\bigr](x)
   +E^{\rm ext}_{<2}(x),
\end{equation}
where $E^{\rm ext}_{<2}$ is the contribution from clusters with fewer than
two renewal points and the prescribed extreme head~$x$.
By Lemma~\ref{l.marked-piece-summability}(3), the non-negative two-mark
series for clusters with fewer than two renewal points converges locally
uniformly on a neighbourhood $U_0$ of~$t_0$.  Reduce $\mathcal N$ if
necessary so that its critical patch lies in~$U_0$.  For each compact
$\mathcal K\Subset\mathcal N$, we may choose $\epsilon>0$ such that
$t_x+\epsilon\hat x\in U_0$ for all $\hat x\in\mathcal K$, with the
series uniformly bounded at these tilts.  The extreme-pair remainder is
bounded by the corresponding unrestricted two-mark coefficient, so
coefficient comparison gives
\[
 E^{\rm ext}_{<2}(x)
 \le C e^{-(t_x+\epsilon\hat x)\cdot x}
 =C e^{-\xi(x)-\epsilon|x|},
\]
where $t_x=t_{\mathrm{str}}(\hat x)$ and
Lemma~\ref{l.fc-string-rate} gives $t_x\cdot x=\xi(x)$.
Thus the remainder is exponentially smaller than the main term.
Lemma~\ref{l.boundary-kernel-convolution}, applied to the admissible pair
$(W^{\base},W^{\head})$, evaluates the first term in
\eqref{e.extreme-pair-convolution}; Lemma~\ref{l.fc-string-rate} identifies
its exponent with~$\xi$.  Positivity follows from the zero-extent
endpoint pieces, while analyticity follows from the two endpoint
transforms, the conjugate tilt and~$\Psi^{\mathrm{str}}$.
\qed

We now turn to unrestricted finite connectivity, the subject of
Theorem~\ref{t.oz-intro}.  The principal remaining work is to establish
its asymptotic in a fixed height frame~$\vecell_0$; the global theorem
will then follow by patching the fixed-frame results.  Since the
prescribed vertices $0$ and $x$
need not be extrema of their finite cluster, re-rooting it at its base
turns them into ordered marks.  The middle-string kernel is unchanged;
we must sum over their possible locations through effective boundary kernels
and control the cases in which their assigned pieces coincide or are
adjacent.  In
coefficient form, the exact identity~\eqref{e.susceptibility-two-marks} is
\begin{equation}
\label{e.two-mark-coefficient}
 \tau^{\rm f}(x)
 =\sum_{\Gamma\in\mathscr C_0}\PP\bigl(\mathsf O_{\rm full}(\Gamma)\bigr)
   \sum_{a,b\in V(\Gamma)}\mathbf 1\{b-a=x\}.
\end{equation}
Use the anchors $q^-,q^+$, the assigned vertex sets~$V^\circ$, and the
one-mark transforms $L_t,R_t$ from Lemma~\ref{l.marked-piece-summability}.

Let $\mathcal Q_-$ and $\mathcal Q_+$ denote the irreducible lower and upper
endpoint pieces, including the zero-extent conventions of
Definition~\ref{d.culdesac-early}, and let $w(Q)$ be the occurrence weight of
such a piece.  Then
\[
 \sum_{Q\in\mathcal Q_-}w(Q)=\widehat W^{\base}(0),\qquad
 \sum_{Q\in\mathcal Q_+}w(Q)=\widehat W^{\head}(0).
\]
For the lower mark~$a$, define the one-mark kernels
\begin{align*}
 J^{\base}(u)
 &=\sum_{Q\in\mathcal Q_-}w(Q)
     \sum_{a\in V^\circ(Q)}\mathbf 1\{q^+-a=u\},\\
 J^{\rm L}(u)
 &=\sum_{S\in\sausageset}\PP(\mathsf O(S))
     \sum_{a\in V^\circ(S)}\mathbf 1\{q^+-a=u\}.
\end{align*}
For the upper mark~$b$, define similarly
\begin{align*}
 J^{\rm R}(v)
 &=\sum_{S\in\sausageset}\PP(\mathsf O(S))
     \sum_{b\in V^\circ(S)}\mathbf 1\{b-q^-=v\},\\
 J^{\head}(v)
 &=\sum_{Q\in\mathcal Q_+}w(Q)
     \sum_{b\in V^\circ(Q)}\mathbf 1\{b-q^-=v\}.
\end{align*}
Put
\[
 \Lambda_0=\Lambda_{\mathrm{vec}}(0)
 =\E_\mu e^{-t_0\cdot X(S)}<1,
 \qquad t_0=\zeta(\vecell_0)\vecell_0,
\]
and define the marked boundary kernels
\begin{equation}
\label{e.marked-boundary-kernels}
 \begin{split}
 B^{\base}(u)&=J^{\base}(u)
   +\frac{\widehat W^{\base}(0)}{1-\Lambda_0}\,J^{\rm L}(u),\\
 B^{\head}(v)&=J^{\head}(v)
   +\frac{\widehat W^{\head}(0)}{1-\Lambda_0}\,J^{\rm R}(v).
 \end{split}
\end{equation}
If the lower mark lies in a sausage, the pieces below it consist of one lower
endpoint piece followed by an arbitrary number of unmarked sausages; their
total raw weight is
$\widehat W^{\base}(0)\sum_{j\ge0}\Lambda_0^j$.  Similarly,
when the upper mark lies in a sausage, the unmarked pieces above it contribute
$\widehat W^{\head}(0)/(1-\Lambda_0)$.

\begin{lemma}[Marked boundary localisation]
\label{l.marked-boundary-localization}
The transforms
\[
 \widehat B^{\base}(s)=\sum_u e^{s\cdot u}B^{\base}(u),\qquad
 \widehat B^{\head}(s)=\sum_v e^{s\cdot v}B^{\head}(v)
\]
are finite and analytic on a neighbourhood of
$\{t_{\mathrm{str}}(\hat n):\hat n\in\mathcal N\}$.  For every compact
$\mathcal K\Subset\mathcal N$ there are $c,C>0$ such that, uniformly for
$\hat n\in\mathcal K$,
\begin{equation}
\label{e.marked-boundary-exp}
 \sum_u e^{t_{\mathrm{str}}(\hat n)\cdot u+c|u|}B^{\base}(u)
 +\sum_v e^{t_{\mathrm{str}}(\hat n)\cdot v+c|v|}B^{\head}(v)
 \le C.
\end{equation}
\end{lemma}

\noindent{\bf Proof.}
The transforms of $J^{\base},J^{\rm L},J^{\rm R},J^{\head}$ are,
respectively, the lower-endpoint $L_s$, sausage $L_s$, sausage $R_s$ and
upper-endpoint $R_s$ transforms of
Lemma~\ref{l.marked-piece-summability}.  That lemma also gives an exponential
moment of the diameter of the marked piece after the critical tilt.  Since
$|q^+-a|$ and $|b-q^-|$ are at most this diameter, a small further factor
$e^{c|u|}$ or $e^{c|v|}$ remains integrable.  For $\mathcal K\Subset\mathcal N$, the tilts
$t_{\mathrm{str}}(\hat n)$ range over a compact subset of the common
convergence neighbourhood, giving uniform constants.
\qed

The marked-sausage terms in~\eqref{e.marked-boundary-kernels} have
positive weight, so neither boundary kernel vanishes.  Thus
$(B^{\base},B^{\head})$ is an admissible boundary pair.  The
longitudinal and lateral parts of~\eqref{e.marked-boundary-exp} come from the
same input: the mass gap controls the height of a sausage, while
Corollary~\ref{c.ordinary-sausage-transverse-fluctuation} charges an
exponential cost for transverse span beyond a fixed linear cone.  For endpoint
pieces the corresponding control is supplied by Lemma~\ref{l.cds-deficit}.

\begin{lemma}[Marked convolution]
\label{l.unrooted-marked-convolution}
For $x\cvlzero>0$,
\begin{equation}
\label{e.unrooted-marked-convolution}
 \tau^{\rm f}(x)
 =\bigl[B^{\base}*\widetilde G*B^{\head}\bigr](x)
  +E_{\rm near}(x)+E_{<2}(x),
\end{equation}
where $E_{\rm near}$ is the contribution from clusters with at least two
renewal points for which the two marks lie in the same or in adjacent assigned
pieces, and $E_{<2}$ is the contribution from clusters with fewer than two
renewal points.
\end{lemma}

\noindent{\bf Proof.}
Use~\eqref{e.two-mark-coefficient} and first suppose that the translated
cluster has at least two renewal points.  Since $b\cvlzero>a\cvlzero$, the assigned
piece containing~$a$ cannot lie above the assigned piece containing~$b$.
If the two marked pieces are separated by at least one complete sausage,
write them as $P_a,P_b$ and put
\[
 u=q^+(P_a)-a,\qquad v=b-q^-(P_b).
\]
The pieces strictly between them form a non-empty string.  If its displacement
is~$z$, then $x=b-a=u+z+v$.
The pieces strictly below~$P_a$ and strictly above~$P_b$ do not enter this
displacement.  By Lemma~\ref{l.marked-renewal-decomposition}, their occurrence
weights factor and may be summed first.  If a marked piece is an endpoint
piece this gives $J^{\base}$ or $J^{\head}$; if it is a sausage, the unmarked
outer pieces give the geometric factors in
\eqref{e.marked-boundary-kernels}.  The non-empty string strictly between the
marked pieces gives $\widetilde G$.  This is the convolution term in
\eqref{e.unrooted-marked-convolution}.

If the marked pieces coincide or are adjacent, the contribution is
$E_{\rm near}$.  Together with $E_{<2}$ these alternatives are disjoint and
exhaust the remaining terms in~\eqref{e.two-mark-coefficient}.
\qed

\begin{lemma}[Local marked remainder]
\label{l.unrooted-remainder}
After reducing $\mathcal N$ if necessary, for every compact
$\mathcal K\Subset\mathcal N$ there are $C,c>0$ such that,
uniformly for $x/|x|\in\mathcal K$,
\begin{equation}
\label{e.unrooted-remainder}
 E_{\rm near}(x)+E_{<2}(x)\le Ce^{-\xi(x)-c|x|}.
\end{equation}
\end{lemma}

\noindent{\bf Proof.}
Put $t_x=t_{\mathrm{str}}(\hat x)$.  We first record convergence of the two
remainder series in a neighbourhood of the critical patch.  If both marks lie
in one assigned piece, its two-mark weight is controlled by
Lemma~\ref{l.marked-piece-summability}; any unmarked sausages before or after
that piece sum geometrically because $\Lambda_0<1$.  If the marked pieces are
adjacent, their common anchor gives
\[
 e^{s\cdot(b-a)}
 =e^{s\cdot(q^+(P_a)-a)}e^{s\cdot(b-q^-(P_b))}.
\]
For a cluster with at least two renewals the possible adjacent types are
lower endpoint--sausage, sausage--sausage and sausage--upper endpoint.  Their
generating series are therefore bounded by products of the corresponding
one-mark transforms from Lemma~\ref{l.marked-piece-summability}, again times
geometric raw outer factors.  Finally, part~\textup{(3)} of that lemma gives a locally uniformly
convergent two-mark series for $E_{<2}$.  The neighbourhood $\mathcal N$ may be chosen small enough that
all these series converge on a common open neighbourhood of
$\{t_{\mathrm{str}}(\hat n):\hat n\in\mathcal N\}$.

Fix a compact set $\mathcal K\Subset\mathcal N$.  There is $\epsilon>0$
such that $t_x+\epsilon\hat x$ remains in that neighbourhood for every
$\hat x\in\mathcal K$, and the two generating series are uniformly bounded
there.  Since their coefficients are non-negative, coefficient comparison
gives
\[
 E_{\rm near}(x)+E_{<2}(x)
 \le C e^{-(t_x+\epsilon\hat x)\cdot x}
 =C e^{-\xi(x)-\epsilon|x|},
\]
where Lemma~\ref{l.fc-string-rate} identifies $t_x\cdot x=\xi(x)$.  This
proves the claim.
\qed

\begin{theorem}[Finite-connectivity asymptotic in a fixed frame]
\label{t.oz-unrooted}
For the height direction $\vecell_0$ and the neighbourhood
$\mathcal N=\mathcal N(\vecell_0)$ of Theorem~\ref{t.oz-rooted},
reduced if necessary, there is a positive real-analytic function
$A_{\rm f}^{\vecell_0}$ on $\mathcal N$ such that, uniformly for
$\hat x=x/|x|$ in every compact
$\mathcal K\Subset\mathcal N$,
\begin{equation}
\label{e.oz-unrooted}
 \PP(0\fconn x)
 =A_{\rm f}^{\vecell_0}(\hat x)|x|^{-(d-1)/2}e^{-\xi(x)}
   \bigl(1+o_{\mathcal K}(1)\bigr),
\end{equation}
and
\begin{equation}
\label{e.oz-unrooted-amplitude}
 A_{\rm f}^{\vecell_0}(\hat n)
 =\widehat B^{\base}\bigl(t_{\mathrm{str}}(\hat n)\bigr)
  \Psi^{\mathrm{str}}(\hat n)
  \widehat B^{\head}\bigl(t_{\mathrm{str}}(\hat n)\bigr).
\end{equation}
\end{theorem}

\noindent{\bf Proof.}
By Lemmas~\ref{l.unrooted-marked-convolution} and
\ref{l.unrooted-remainder}, the remainder is exponentially smaller than the
claimed main term.  Lemma~\ref{l.boundary-kernel-convolution}, applied to the
admissible pair $(B^{\base},B^{\head})$ and to~$\widetilde G$, evaluates the
convolution in~\eqref{e.unrooted-marked-convolution}.  Finally
Lemma~\ref{l.fc-string-rate} replaces $\xi_{\mathrm{str}}$ by~$\xi$ on the
cone.  Analyticity follows from the conjugate tilt, the string factor and the
two boundary transforms.  The two boundary kernels are non-zero,
as observed after Lemma~\ref{l.marked-boundary-localization}, so the
amplitude is positive.
\qed

\noindent{\bf Proof of Theorem~\ref{t.oz-intro}.}
For each $\hat n\in S^{d-1}$ take
$\vecell_{\hat n}=\drift^{-1}(\hat n)$.  The direction~$\hat n$ lies in the
interior of $\mathcal N(\vecell_{\hat n})$.  Choose a neighbourhood
$U_{\hat n}$ of~$\hat n$ whose closure is a compact subset of
$\mathcal N(\vecell_{\hat n})$.  By
compactness of the sphere, finitely many of these neighbourhoods cover
$S^{d-1}$.

Theorem~\ref{t.oz-unrooted} gives the asserted asymptotic uniformly on the
closure of each member of this finite cover.  On an overlap, fix a direction
$\hat n$ and choose lattice points $x_k$ with
$x_k/|x_k|\to\hat n$ inside the overlap.  Comparing the two local
asymptotics for the same quantity $\tau^{\rm f}(x_k)$ shows that their
amplitudes agree at~$\hat n$.  They therefore patch to a positive
real-analytic function $A_{\rm f}$ on the sphere.  Taking the maximum of the finitely many local
error terms shows that the resulting $o(1)$ is uniform in~$\hat x$.
Lemma~\ref{l.fc-string-rate} identifies the exponent with~$\xi$, while
Theorem~\ref{l.drift-onto} supplies the real-analytic and strictly convex Wulff
and equidecay surfaces.  In particular, the amplitude at~$\hat n$ may be
evaluated by~\eqref{e.oz-unrooted-amplitude} with the canonical choice
$\vecell_0=\drift^{-1}(\hat n)$.  Then
$t_{\mathrm{str}}(\hat n)=\zeta(\vecell_0)\vecell_0$, and the string factor
$\Psi^{\mathrm{str}}(\hat n)$ is determined by the mean and covariance of
the original tilted sausage law~$\mu^{(\vecell_0)}$.
\qed

We conclude with two observations about the geometry of long finite
connections that emerge from the preceding asymptotics.  The first
interprets the difference between the unrestricted and extreme-pair
amplitudes through one-sided endpoint occupations.  The second identifies
the drift of the tilted sausage law with the typical direction of the
conditioned cluster head.

For the fixed height direction~$\vecell_0$, the amplitudes
$A_{\rm f}^{\vecell_0}$ and $A_{\mathrm{ext}}^{\vecell_0}$
share the middle renewal factor~$\Psi^{\mathrm{str}}$.  Their ratio has a one-sided occupation
interpretation.  Fix $t=t_{\mathrm{str}}(\hat n)$ and recall the tilted
sausage law $\nu_t(S)=e^{t\cdot X(S)}\PP(\mathsf O(S))$.
Construct $\mathcal E_t^{\base}$ by sampling a lower endpoint piece~$Q$
with probability $e^{t\cdot X(Q)}w(Q)/\widehat W^{\base}(t)$ and
appending independent IID $\nu_t$-sausages.  Root this
one-sided environment at its actual base.  Construct
$\mathcal E_t^{\head}$ analogously from a tilted upper endpoint piece,
prepending independent IID $\nu_t$-sausages and rooting at its
actual head.

\Needspace{7\baselineskip}
\noindent The expected discounted occupations of the two environments are
\[
 H^{\base}(t)=\E_t^{\base}\!\left[
   \sum_{a\in V(\mathcal E_t^{\base})}e^{-t\cdot a}\right],
 \qquad
 H^{\head}(t)=\E_t^{\head}\!\left[
   \sum_{b\in V(\mathcal E_t^{\head})}e^{t\cdot b}\right].
\]
In each sum the actual extremum is at the origin.  If the marked vertex
lies in the endpoint piece, its discount cancels the endpoint tilt and
leaves the corresponding one-mark weight.  If it lies in a sausage, the
discount also untilts the endpoint piece and any intervening complete
sausages.  The former contributes $\widehat W^{\base}(0)$ or
$\widehat W^{\head}(0)$, and the latter give the geometric sum
$\sum_{j\ge0}\Lambda_0^j=(1-\Lambda_0)^{-1}$.
Thus~\eqref{e.marked-boundary-kernels} gives
\begin{equation}
\label{e.environment-amplitude}
 \frac{A_{\rm f}^{\vecell_0}(\hat n)}
      {A_{\mathrm{ext}}^{\vecell_0}(\hat n)}
 =\frac{\widehat B^{\base}(t)}{\widehat W^{\base}(t)}
  \frac{\widehat B^{\head}(t)}{\widehat W^{\head}(t)}
 =H^{\base}(t)H^{\head}(t).
\end{equation}
Lemma~\ref{l.marked-boundary-localization} and positivity of
$\widehat W^{\base}(t)$ and $\widehat W^{\head}(t)$ show that both
occupations are finite.  Each actual extremum has weight one, while
other vertices contribute positively, so both factors strictly exceed one.
Dividing the two
Ornstein--Zernike asymptotics also gives, uniformly as $|x|\to\infty$
with $\hat x\in\mathcal K\Subset\mathcal N$,
\[
 \PP\bigl(0\fconnext_{\vecell_0}x\mid 0\fconn x\bigr)
 =\frac{1}{H^{\base}(t_x)H^{\head}(t_x)}+o_{\mathcal K}(1),
 \qquad t_x=t_{\mathrm{str}}(\hat x).
\]
The occupation factors thus quantify the additional positions of the
prescribed vertices relative to the actual extrema.

\medskip
The vector $\drift(\vecell)$ entered our analysis as the normalised mean
increment under the tilted sausage law.  We substantiate that it is the
asymptotic direction of the cluster head under point-to-hyperplane conditioning.

\begin{corollary}[Cluster head direction]
\label{c.conditioned-head-drift}
Fix $\vecell$ and, on $\clusteruthree$, put
$X_u=\head(\mc C(0))$.  For every $\varepsilon>0$, there are
$C_\varepsilon,c_\varepsilon>0$ such that, for all sufficiently large~$u$,
\[
 \PP\left(\left|\frac{X_u}{|X_u|}-\drift(\vecell)\right|>\varepsilon
   \;\biggm|\;\clusteruthree\right)
 \le C_\varepsilon e^{-c_\varepsilon u}.
\]
\end{corollary}

\noindent{\bf Proof.} Set $t_0=\zeta(\vecell)\vecell$ and
$n_0=\drift(\vecell)$.  By Theorem~\ref{l.drift-onto}, $n_0$ is the
outward normal to $\mathcal W$ at~$t_0$.  Since $\xi$ is the support
function of~$\mathcal W$, the uniqueness of this normal and compactness
of the sphere give $\delta_\varepsilon>0$ such that
$\xi(x)\ge t_0\cdot x+\delta_\varepsilon|x|$ whenever $x\ne0$ and
$|x/|x|-n_0|>\varepsilon$.  Theorem~\ref{t.oz-intro} gives
$\PP(0\fconn x)\le C e^{-\xi(x)}$ uniformly in~$x$.  On
$\clusteruthree$, the head has height in $[u,u+3r]$ and norm at least~$u$.
Summing the finite-connectivity bound over its possible positions therefore
gives
\[
 \PP\left(\clusteruthree,\,
    \left|\frac{X_u}{|X_u|}-n_0\right|>\varepsilon\right)
 \le C e^{-\zeta u}\sum_{\substack{x\in\Z^d\\|x|\ge u}}
       e^{-\delta_\varepsilon|x|}
 \le C_\varepsilon e^{-(\zeta+c_\varepsilon)u}.
\]
As noted in the proof of Theorem~\ref{t.renewal},
$\PP(\clusteruthree)\ge c e^{-\zeta u}$ for all large~$u$.
Thus we obtain the corollary. \qed

\appendix

\chapter{The layer-box Antal--Pisztora estimate}
\label{s.nonaxial-ap-proof}

This appendix proves Proposition~\ref{p.nonaxial-ap}.  The bound asserted
by the proposition is that, for
$K\ge K_0$ and every overlapping layer box $B$,
\begin{equation}
\label{e.ap-box-estimate}
   \PP(B\text{ is not apt})\, \le \, Ce^{-cK}\, ,
\end{equation}
where $C,c>0$ and $K_0\in\N$ may be chosen independently of
$\vecell\in S^{d-1}$.
Aptness of $B$ is measurable with respect to the edge configuration on
$E(B)$, and the inapt-box field on $\Z\times\Z^{d-1}$ is one-dependent.  For
every height offset $a\in\R$, the same estimate and measurability hold,
with the same constants, for the boxes of $\overline N(a)$, and their inapt-box field on
$\Z^{d-1}$ is one-dependent.  This estimate supplies the apt-box input for the
quilt and seal construction beginning in Chapter~\ref{s.overlapping-quilts-seals}.

Aptness of $B$ asks for one open cluster that contains, in every boxlet,
an open path of diameter at least~$3K/4$ and also contains every open path
in $B$ of diameter at least $K/2$.  The classical axial Antal--Pisztora
input does not apply directly to these layer boxes: their height faces are
generally tilted, and their edges are allocated by upper endpoint.  We
therefore use that input at a smaller axial scale $\kappa=\kappa(d,p)$,
chosen sufficiently large and then fixed while $K$ tends to infinity.  The
apt $\kappa$-boxes of the overlapping array form a high-density site field.

We work with the \emph{jewel}, the collection of axial $\kappa$-boxes
lying wholly inside the layer-box region.  It reaches close to the sides
of this region, though its boundary is irregular where the axial boxes
meet the tilted faces.  In suitable thin \emph{working layers}, we can
select \emph{section arrays} of jewel boxes.  These may have a staircase
shape in physical space, but are parametrised by rectangles in
$\Z^{d-1}$.  The high density of apt $\kappa$-boxes supplies, with high
probability, a large connected woven fabric in each array.

Open connections through the jewel, made well away from the boundary of
the layer box, allow us to construct the cluster that will witness
aptness.  Within each boxlet of $B$, a large connected family of apt
$\kappa$-boxes supplies an open path running upwards through a substantial
part of the boxlet.  The union of these paths need not be connected.
Chains of overlapping rectangular arrays join the boxlet pieces to a
common component near the centre of $B$.  The open cluster containing
the resulting woven fabric is our candidate; it contains a long open
path in every boxlet.

We then show that every other macroscopic open path in $B$ must join this
cluster.  Such a path has substantial extent in some transverse
coordinate or in height, and therefore encounters order $K/\kappa$
successive working layers.  Here the jewel's reach towards the sides
becomes important.  Even when the path arrives near the irregular
boundary, a short connection within each fresh layer gives a fixed
positive chance of joining the large woven fabric of its section array.
With high probability, these layer fabrics are themselves connected
to the candidate cluster through portions of the jewel well away from
its boundary.
A macroscopic path remaining outside that cluster must therefore
repeatedly evade capture, an event of exponentially small probability.

The four sections that follow develop this outline into a proof.  They set up the geometry and high-density input, construct the
boxlet-reaching cluster, prove containment of macroscopic paths using
fresh layers, and assemble the estimate.

\medskip
\noindent\textbf{Constants.}
Throughout, $c_d$ and $C_d$ denote positive finite constants depending only on
the dimension, whose values may change between occurrences.  The scale
$\kappa$ will be chosen in terms of $d$ and $p$ alone; subsequent constants
may depend on it.  No constant in this
appendix depends on the height direction $\vecell$, on the box index $(k,v)$,
on a height offset $a$ used in place of $a_k$, or on $K$.  Uniformity in
$\vecell$ is secured by one design decision: the auxiliary array is axial, so
the tilt enters only through the geometry of the layer box, and there only
through the bounds $0\le\ell_i\le\ell_d$ and $\ell_d\ge d^{-1/2}$ of the
axes-ordering convention~\eqref{e.aoc}.

\section{Auxiliary geometry and the \texorpdfstring{$\kappa$}{kappa}-array}

We first record the uniform geometry of a layer box, then define the axial
$\kappa$-array and the jewel, and finally isolate the high-density site
input used on rectangular subarrays of that array.

\subsection{The geometry of a layer box}
\label{ss.nonaxial-geometry}

Fix a layer box $B=\Obox{k,v}$, and let
$\Lambda=\Lambda_{k,v}$ be the half-open parallelotope given
by~\eqref{e.layer-box-region}. The only vertices
of the edge graph $B$ outside $\Lambda\cap\Z^d$ are lower endpoints of edges
assigned to the box, and each is one lattice step from that region.
The same geometry estimates apply to offset layer boxes in
$\overline N(a)$: only the defining height interval changes.

The geometry becomes transparent in the \emph{defining coordinates}
\[
   \rho_i(x)=x_i\quad(i<d),
   \qquad
   \rho_d(x)=x\cvl .
\]
The map $\rho=(\rho_1,\ldots,\rho_d)$ sends $\Lambda$ to a half-open
cube of side $2K$, one interval per defining coordinate. We write
$\Xi=\{\pm\rho_1,\ldots,\pm\rho_d\}$ for the set of
\emph{signed} defining coordinates; the signs allow coordinate sweeps in
either direction.  For $\varphi\in\Xi$ and non-empty $A\subseteq\R^d$ we
call $\sup_A\varphi-\inf_A\varphi$ the \emph{$\varphi$-extent} of $A$.
The projection of $\Lambda$ onto the first $d-1$ coordinates is its
\emph{transverse footprint}.  For $x'\in\Z^{d-1}$ in that footprint, the
lattice points $(x',m)\in\Lambda$ form the $e_d$-\emph{fibre} above $x'$.
In expressions involving $x'$, $e_i$ ($i<d$) denotes the corresponding
basis vector of $\R^{d-1}$.

\begin{lemma}[Geometry of a layer box]
\label{l.layer-box-geometry}
There are $c_d,C_d\in(0,\infty)$ with the following properties.
\begin{enumerate}
\item\label{l.layer-box-geometry.bilip}
For every $y\in\R^d$,
\[
 c_d\,|y|
 \, \le \,
 \max\bigl\{|y_1|,\ldots,|y_{d-1}|,|y\cvl|\bigr\}
 \, \le \, |y| .
\]
\item\label{l.layer-box-geometry.inball}
The region $\Lambda$ is convex and has a unique largest inscribed open ball,
of radius $K$.  Its Euclidean diameter is at most $C_dK$.
\item\label{l.layer-box-geometry.fibres}
For $x'\in\Z^{d-1}$ in the transverse footprint of $\Lambda$, its
$e_d$-fibre is
$\{(x',m):m\in J_k(x')\cap\Z\}$, where $J_k(x')\cap\Z$ contains at least
$2K$ consecutive integers.  If $x'$ and $x'\pm e_i$, $i<d$, both lie in the
transverse footprint, the corresponding intervals $J_k$ are translates by
at most one.  Consequently any nearest-neighbour transverse path of length
$m$ in the footprint lifts, from any chosen point of its first fibre, to a
nearest-neighbour path in $\Lambda\cap\Z^d$ of length at most~$2m$.
\item\label{l.layer-box-geometry.lip}
Every signed defining coordinate is $1$-Lipschitz along lattice edges: if
$x\sim y$ and $\varphi\in\Xi$ then $|\varphi(y)-\varphi(x)|\le1$.
\end{enumerate}
\end{lemma}

\noindent{\bf Proof: (1).} The upper bound in part~\eqref{l.layer-box-geometry.bilip} is immediate. For the
lower, write $y=(y',y_d)$ and put
$M=\max\{|y_1|,\ldots,|y_{d-1}|,|y\cvl|\}$. From
$y_d\ell_d=y\cvl-\sum_{i<d}\ell_iy_i$ and $0\le\ell_i\le\ell_d$,
\[
 |y_d|
 \, \le \, \frac{|y\cvl|+\sum_{i<d}\ell_i|y_i|}{\ell_d}
 \, \le \, \bigl(d^{1/2}+d-1\bigr)M ,
\]
using $\ell_d\ge d^{-1/2}$. Hence $|y|\le C_dM$, which is the claim.

{\bf (2).} The region is an intersection of
$2d$ half-spaces, so it is convex. Its $d$ pairs of opposite faces have unit
normals $e_1,\ldots,e_{d-1},\vecell$ and Euclidean separation $2K$ apiece.
Those normals are linearly independent, the determinant of the matrix they
form being $\ell_d\ge d^{-1/2}$; so the $d$ mid-hyperplanes meet in a single
point, at distance exactly $K$ from each of the $2d$ faces. The open ball of
radius $K$ about it therefore lies in $\Lambda$. No ball of radius exceeding
$K$ fits between a pair of opposite supporting hyperplanes at separation
$2K$, so this is a largest inscribed ball. Its centre is unique: a ball of
radius $K$ must be centred on each of the $d$ mid-hyperplanes. Every
defining coordinate varies by at most $2K$ on $\Lambda$, so
part~\eqref{l.layer-box-geometry.bilip} bounds the diameter by $C_dK$.

{\bf (3).} By the fibre calculation of
Section~\ref{s.overlapping-boxes-boxlets}, $J_k(x')$ has length
$2K/\ell_d\ge2K$, while $J_k(x'\pm e_i)$ is a translate of $J_k(x')$ by at
most one.  If an integer level in one of two adjacent fibres is not in the
other, a single vertical step reaches a level common to both, since each
interval contains at least $2K\ge2$ consecutive integers.  Thus every
transverse step lifts using at most one additional vertical step, giving the
bound $2m$.

{\bf (4).} Apply the upper bound of part~\eqref{l.layer-box-geometry.bilip}
to the difference of the endpoints of a lattice edge. \qed

Fix henceforth a dimensional constant $c_{\rm span}>0$ such that, by
part~\eqref{l.layer-box-geometry.bilip}, every connected set of Euclidean
diameter at least $K/2$ has $\varphi$-extent at least $c_{\rm span}K$
for some $\varphi\in\Xi$.  We use this fixed constant in the absorption
argument.

\subsection{The axial \texorpdfstring{$\kappa$}{kappa}-array}
\label{ss.nonaxial-kappa-array}

The auxiliary array is the usual axial overlapping array at the smaller
scale~$\kappa$, aligned with the lattice axes.

\begin{definition}[The $\kappa$-array and aptness at scale $\kappa$]
\label{d.kappa-array}
Let $\kappa\in\N$. For $u\in\Z^d$ the \emph{$\kappa$-boxlet} $Q^\circ_u$ is the
subgraph induced on
\[
   \bigl\{x\in\Z^d:\ \kappa u_i\le x_i<\kappa(u_i+1),\ i\in\intint d\bigr\},
\]
and the \emph{$\kappa$-box} $Q_u$ is the subgraph induced on the union of the
$2^d$ $\kappa$-boxlets $Q^\circ_{u+\sigma}$, $\sigma\in\{0,1\}^d$: an
axis-parallel box of side $2\kappa$. Its edge set $E(Q_u)$ consists of the
edges with both endpoints in $Q_u$. Two $\kappa$-boxes share a
$\kappa$-boxlet exactly when their indices are at $\ell^\infty$-distance at
most one, and have disjoint edge sets when that distance is at least two.

The $\kappa$-box $Q_u$ is \emph{apt} if there is an open cluster $\mc W(Q_u)$
in $Q_u$ such that
\begin{enumerate}
\item for every $\kappa$-boxlet $Q^\circ$ of $Q_u$, the set
      $\mc W(Q_u)\cap Q^\circ$ contains an open path of Euclidean diameter at
      least $3\kappa/4$; and
\item every open path in $Q_u$ of Euclidean diameter at least $\kappa/2$ is
      contained in $\mc W(Q_u)$.
\end{enumerate}
These are the two aptness requirements of Definition~\ref{d.apt-box},
now imposed on the induced axial box~$Q_u$; we use the same terminology at
this auxiliary scale, calling $\mc W(Q_u)$ its \emph{woven cluster}.
\end{definition}

For the fixed layer box $B$, the \emph{jewel} is the collection of
$\kappa$-boxes $Q_u$ with $V(Q_u)\subseteq\Lambda$, where $\Lambda$ is the
parallelotope associated with $B$ above.  If
$Q_u$ is a jewel box, then every edge of $E(Q_u)$ has both endpoints in
$\Lambda$, hence has its upper endpoint there, and so belongs to $E(B)$
by~\eqref{e.cell-edge-partition}.  Thus the aptness test for every jewel box
uses only the edge configuration of $B$.

\begin{lemma}[Neighbouring apt $\kappa$-boxes]
\label{l.kappa-chaining}
\leavevmode
\begin{enumerate}
\item If $Q_u$ and $Q_{u'}$ are apt and $\|u-u'\|_\infty=1$, then
      $\mc W(Q_u)\cap\mc W(Q_{u'})\ne\emptyset$.
\item If $\mathscr F$ is a finite $*$-connected family of apt $\kappa$-boxes, then
      \[
         \mc W(\mathscr F):=\bigcup_{Q_u\in\mathscr F}\mc W(Q_u)
      \]
      is a connected open subgraph, which we call the \emph{woven fabric}
      of~$\mathscr F$. Any open path with a subpath of diameter at least
      $\kappa/2$ inside a box of $\mathscr F$ has that subpath contained in
      $\mc W(\mathscr F)$, and hence lies in the same open cluster as the
      woven fabric.
\end{enumerate}
\end{lemma}

\noindent{\bf Proof.} For part~(1), the boxes share a $\kappa$-boxlet.
The woven cluster of $Q_u$ contains there an open path of diameter at least
$3\kappa/4$, which lies in $Q_{u'}$ and hence in $\mc W(Q_{u'})$ by
its absorption property.  The two woven clusters therefore intersect.
Chaining these intersections along $*$-paths in $\mathscr F$ shows that
$\mc W(\mathscr F)$ is connected. If an open path has a subpath of
diameter at least $\kappa/2$ inside $Q_u\in\mathscr F$, aptness places that
subpath in $\mc W(Q_u)$ and hence in $\mc W(\mathscr F)$. \qed

\begin{lemma}[Apt $\kappa$-box properties]
\label{l.kappa-apt}
Apt $\kappa$-boxes have the following properties.
\begin{enumerate}
\item\label{l.kappa-apt.meas}
\textup{(Locality.)} $\{Q_u\text{ is apt}\}\in\sigma\bigl(E(Q_u)\bigr)$, so the
apt field $(\mathbf 1\{Q_u\text{ is apt}\})_{u\in\Z^d}$ is one-dependent.
\item\label{l.kappa-apt.rate}
\textup{(Rarity of failure.)}
$\PP(Q_u\text{ is not apt})\le\varepsilon_\kappa\to0$
as $\kappa\to\infty$, for every $u$.
\item\label{l.kappa-apt.dom}
\textup{(Domination.)} For every $\delta>0$ there is $\kappa_0(d,p,\delta)$
such that, for $\kappa\ge\kappa_0$, the apt field stochastically dominates
Bernoulli site percolation of parameter $1-\delta$ on $\Z^d$.
\end{enumerate}
\end{lemma}

\noindent{\bf Proof: (1).} Aptness of $Q_u$ is determined by the states of the edges with both endpoints
in $Q_u$, which establishes part~\eqref{l.kappa-apt.meas}; $\kappa$-boxes at index
distance at least two have disjoint edge sets.

{\bf (2).} We use the finite-box crossing and local-uniqueness estimate
of Pisztora~\cite[Theorem~3.1]{Pisztora96}; see also the good-box events
in~\cite[Section~2, (2.9) and (2.24)]{AntalPisztora}.  For $p>p_c$ and
any fixed positive fraction of the side length, this estimate gives, with
probability tending to one, a crossing open cluster in an axial cube and
uniqueness among open clusters of diameter at least that fraction.

Apply the crossing assertion to each of the $2^d$ constituent
$\kappa$-boxlets of $Q_u$, and local uniqueness in the full box $Q_u$
with diameter threshold $\kappa/2$.  The latter condition puts every open
path in $Q_u$ of Euclidean diameter at least $\kappa/2$ in one cluster.
Indeed, Euclidean diameter $\kappa/2$ implies $\ell^\infty$-diameter
at least $\kappa/(2\sqrt d)$, a fixed fraction of the side of $Q_u$.
Each boxlet crossing has diameter at least $\kappa-1\ge3\kappa/4$ for
large $\kappa$, so local uniqueness places all these crossings in the
same cluster.  This is aptness.  A union bound over the $2^d$ boxlets
and $Q_u$, followed by translation invariance, gives
$\varepsilon_\kappa\to0$ uniformly in~$u$.

{\bf (3).} This follows from locality, rarity and the
Liggett--Schonmann--Stacey theorem~\cite{LSS}. \qed

\subsection{High-density site input}
\label{ss.nonaxial-high-density}

Lemma~\ref{l.kappa-apt} says that, once $\kappa$ is large, the aptness
indicators of the $\kappa$-boxes dominate an independent and identically
distributed site field of density as close to one as we require.  Every use
of this fact below is on a rectangular
set of $\kappa$-indices with side lengths of order $K/\kappa$: the weaving
argument uses $d$-dimensional rectangles in the jewel, while a fresh layer is
pulled back to a $(d-1)$-dimensional rectangle.  We therefore isolate one
high-density statement that serves both purposes.  In the lemma, $X$ is an
abstract $\{0,1\}$-valued site field; in our applications, $X_v=1$ means that
the corresponding $\kappa$-box is apt.

Besides a dense crossing component, the statement records the two
consequences needed later: local uniqueness at linear scale and a uniform
rooted chance of belonging to the dense component.

\begin{lemma}[High-density rectangles]
\label{l.high-density-rectangles}
For each $r\ge2$ there are $\delta_\star(r),\eta(r)>0$ such that the
following holds. Fix $1\le A<\infty$ and $\alpha\in(0,1/4)$. There are
$c,C>0$ and $n_0<\infty$, possibly depending on $r,A,\alpha$, with the
following property. Let
\[
   R=I_1\times\cdots\times I_r\subseteq\Z^r,
   \qquad n:=\min_j|I_j|,\qquad |I_j|\le An,
\]
with $n\ge n_0$, and let $X=(X_v)_{v\in R}$ be a $\{0,1\}$-valued site
field, calling $v$ $X$-open when $X_v=1$. Suppose that $X$ stochastically
dominates Bernoulli site percolation of parameter $1-\delta_\star(r)$.

Let $\mathscr D_R$ be the unique nearest-neighbour connected $X$-open
component containing more than $7|R|/8$ sites, when such a
component exists, and put $\mathscr D_R=\emptyset$ otherwise.  When
$\mathscr D_R\ne\emptyset$, let $\mathscr G_R$ be the $*$-connected $X$-open
component containing it, and otherwise put $\mathscr G_R=\emptyset$.
The components $\mathscr D_R$ and $\mathscr G_R$ do not depend on
$\alpha$.  Let $\mathsf H_R(\alpha)$ be the $X$-measurable event on which
the following conditions hold:
\begin{enumerate}
\item\label{l.high-density-rectangles.cross}
$\mathscr D_R$ meets every face of $R$;
\item\label{l.high-density-rectangles.nn}
every nearest-neighbour connected $X$-open component of $R$ with
$\ell^\infty$-diameter at least $\alpha n$ is $\mathscr D_R$;
\item\label{l.high-density-rectangles.star}
every $*$-connected $X$-open component of $R$ with
$\ell^\infty$-diameter at least $\alpha n$ is $\mathscr G_R$.
\end{enumerate}
Then
\[
   \PP(\mathsf H_R(\alpha)^{\rm c})\le Ce^{-cn}.
\]
Moreover, uniformly in $v\in R$,
\begin{equation}
\label{e.high-density-rooted}
   \PP(v\in\mathscr D_R)\ge\eta .
\end{equation}
\end{lemma}

\noindent{\bf Proof.}
Couple $X$ above an independent and identically distributed Bernoulli
field $Y$ of parameter $1-\delta_\star(r)$.  Theorem~1.1 of
Deuschel--Pisztora~\cite{DeuschelPisztora1996} gives, at sufficiently high
density, a nearest-neighbour open cluster meeting every face and occupying
at least a prescribed proportion of a rectangle whose side lengths lie
between $n$ and $2n$.  Fix the density error in this input to be
$\varepsilon_r=2^{-r-5}$.  The required Bernoulli threshold then depends
only on~$r$.

For a general rectangle $R$, partition each coordinate interval into
$M_j=\lceil |I_j|/n\rceil$ consecutive blocks of nearly equal lengths.
If $M_j\ge2$, cover that interval by the unions of consecutive pairs of
blocks; for $M_j=1$, use the interval itself.  For sufficiently large $n$,
each covering interval has length between $n$ and $2n$, and consecutive
ones overlap in at least a third of either interval.  Taking products gives
a connected covering of $R$ by near-cubic rectangles, with adjacent
members overlapping in at least a third of their volume.  Each site of
$R$ lies in at most $2^r$ members.  The number of members may depend on
$A$, but these overlap and multiplicity bounds do not.

On the simultaneous Deuschel--Pisztora good events, the dense crossing
clusters of adjacent covering members intersect, since the missing volume
in each member is at most $\varepsilon_r$ of its volume.  Their union is
therefore a single $Y$-open component meeting every face of $R$, and its
complement has at most $2^r\varepsilon_r|R|=|R|/32$ sites.  In particular
this component occupies more than $7|R|/8$ sites.  A union bound over the
fixed number of members gives failure probability $Ce^{-cn}$ for $r\ge2$.

The same event establishes part~\eqref{l.high-density-rectangles.cross} for $X$:
the $Y$-cluster lies inside an $X$-open component of at least the same size.
Although $\mathsf H_R(\alpha)$ itself is not increasing, the two uniqueness
assertions will follow from a separate bound on $X$-closed obstructions.
The standard lattice-animal counting argument (see
Kesten~\cite[Lemma~5.1, Eq.~(5.22)]{Kesten1982}), applied to the
$*$-adjacency graph, gives $a_r<\infty$ such that the number of
$*$-connected sets of $m$ sites containing a prescribed site is at
most $a_r^m$.  Since an
$X$-closed set is $Y$-closed in the coupling, choose $\delta_\star(r)$
still smaller so that $a_r\delta_\star(r)<1$.  Then, for every fixed $\beta>0$,
there are $C_\beta,c_\beta>0$ such that
\begin{equation}
\label{e.high-density-animal}
 \PP\bigl(\text{there is a $*$-connected $X$-closed set in $R$ of
 size at least $\beta n$}\bigr)\le C_\beta e^{-c_\beta n}.
\end{equation}
Indeed, root the set at one of the $O(n^r)$ sites of $R$ and sum
$\sum_{m\ge\beta n}(a_r\delta_\star(r))^m$.  The same
dimension-dependent choice works for every fixed~$\beta$.

We use two elementary finite-box facts from
\cite[Lemmas~2.1 and~2.2]{GLA}.  If a set of sites separates two
nearest-neighbour connected subsets of a lattice box, it contains a
$*$-connected subset that still separates them.  Also, if $C$ is
$*$-connected and $x\notin C$, then the sites $*$-adjacent to $C$
that lie in the nearest-neighbour component of $x$ in the complement form a
nearest-neighbour connected set.  These facts are proved there for cubic lattice boxes; the arguments apply
equally to rectangular lattice boxes such as $R$.

Work now on the crossing event already obtained for $X$, on which
$\mathscr D_R$ meets every face.  A nearest-neighbour or $*$-connected set
meets every intermediate coordinate level between its extrema, since each
step changes any coordinate by at most one.

Suppose that $C$ is a nearest-neighbour $X$-open component, distinct
from $\mathscr D_R$, of $\ell^\infty$-diameter at least $\alpha n$.  The
closed sites of $R$ separate $C$ from $\mathscr D_R$, so by the first fact
they contain a $*$-connected subset $F$ which still separates them.
Choose a coordinate in which $C$ spans at least $\alpha n$ levels.  On every
such level choose a site of $C$ and a site of $\mathscr D_R$; the latter
exists because $\mathscr D_R$ is connected and crosses the two opposite
faces in that coordinate.  A nearest-neighbour path joining the two sites
inside that coordinate slice must meet $F$.  Different levels give different
sites of $F$, so $|F|\ge\alpha n$.  Estimate~\eqref{e.high-density-animal}
proves part~\eqref{l.high-density-rectangles.nn}.

For part~\eqref{l.high-density-rectangles.star}, let $C$ instead be a
$*$-connected $X$-open component distinct from $\mathscr G_R$ and of
$\ell^\infty$-diameter at least $\alpha n$.  Fix
$x\in\mathscr D_R$, and let $\Gamma$ be the sites of $R$ that are
$*$-adjacent to $C$ and lie in the nearest-neighbour component of $x$ in
$R\setminus C$.  By the second fact above, $\Gamma$ is nearest-neighbour
connected, hence $*$-connected.  Every site of $\Gamma$ is $X$-closed:
an $X$-open site $*$-adjacent to $C$ would belong to $C$ by maximality.  Choose a coordinate in which $C$ spans at least
$\alpha n$ levels.  For each such level choose
$c_j\in C$ and $d_j\in\mathscr D_R$ on that level.  Join $x$ to $d_j$
inside $\mathscr D_R$, and from $d_j$ follow a nearest-neighbour path in the
coordinate slice towards $c_j$.  The first site on this latter path nearest-neighbour adjacent
to $C$ belongs to $\Gamma$: the preceding path remains outside $C$ and
joins it to $x$, while nearest-neighbour adjacency implies $*$-adjacency.
It also lies on level~$j$.  Thus different levels
supply different sites of $\Gamma$, so $|\Gamma|\ge\alpha n$, and
\eqref{e.high-density-animal} again applies.  Together with the crossing
estimate, these two obstruction bounds give
$\PP(\mathsf H_R(\alpha)^{\rm c})\le Ce^{-cn}$.

It remains to prove~\eqref{e.high-density-rooted}.  For each $v\in R$ choose,
in every coordinate, the direction in which at least $(n-1)/2$ sites of $R$
remain.  A standard Peierls argument in the resulting orthant, using the first
finite-box fact above, gives $\eta_0=\eta_0(r)>0$, once
$\delta_\star(r)$ is sufficiently small, such that uniformly in $n$,
$R$, $v$ and $\alpha\in(0,1/4)$, the Bernoulli field $Y$ has with
probability at least $\eta_0$ a nearest-neighbour open path from $v$ to
$\ell^\infty$-distance $\alpha n$.  The same path is $X$-open.  On
$\mathsf H_R(\alpha)$ its $X$-open component is $\mathscr D_R$ by part~\eqref{l.high-density-rectangles.nn}.  Enlarging $n_0$ makes
$Ce^{-cn}\le\eta_0/2$, and~\eqref{e.high-density-rooted} follows with
$\eta(r)=\eta_0/2$. \qed

Choose $\delta_{\rm hd}=\min\{\delta_\star(d),\delta_\star(d-1)\}$.
By part~\eqref{l.kappa-apt.dom} of Lemma~\ref{l.kappa-apt}, we may henceforth
fix $\kappa=\kappa(d,p)$ large enough that the apt $\kappa$-box field
dominates Bernoulli site percolation of parameter $1-\delta_{\rm hd}$.

\section{A woven fabric that reaches every boxlet}
\label{ss.nonaxial-woven}

We now establish the weaving aspect of aptness: except with probability
exponentially small in~$K$, the jewel contains a single woven fabric that
reaches every boxlet of~$B$ by an open path of diameter at least~$3K/4$.
Throughout this section $K$ is sufficiently large relative to the fixed
auxiliary scale~$\kappa$.

\begin{lemma}[Rectangular core of a boxlet]
\label{l.cell-rectangular-core}
There is $c_d>0$ such that every boxlet $L$ of $B$ contains the lattice
points and nearest-neighbour edges of an axis-parallel rectangular prism
$\mathrm{Pr}(L)$.  Its first $d-1$ side lengths are at least $c_dK$, its
$e_d$-side length is at least $K$, and it lies at Euclidean distance at least
$c_dK$ from every face of $\Lambda$.
\end{lemma}

\noindent{\bf Proof.} Let $L=\Bbox{k,u}$ and choose in each transverse coordinate the centred
interval of $\lfloor K/(4d)\rfloor$ consecutive integers inside the range
$(K(u_i-1),Ku_i]$. As $x'$ ranges over the resulting transverse rectangle,
the endpoints of its $e_d$-fibre move by at most
$\lfloor K/(4d)\rfloor\sum_{i<d}\ell_i/\ell_d\le(d-1)K/(4d)\le K/4$.
Every fibre has length $2K/\ell_d\ge2K$, so the fibres over the rectangle
have a common interval of length at least $7K/4$. Take its centred
$e_d$-subinterval of length $K$: it lies at $e_d$-distance at least $3K/8$
from either endpoint of every fibre concerned, hence at Euclidean distance at
least $3K\ell_d/8\ge3K/(8d^{1/2})$ from the two height faces. The transverse
centring leaves at least $(K-\lfloor K/(4d)\rfloor)/2\ge K/4$ from each
transverse face.  Every lattice point of the prism belongs to the boxlet's
column set, and every edge between such points is a boxlet edge.
Reducing the constant proves the assertion. \qed

For each boxlet~$L$, fix the prism $\mathrm{Pr}(L)$ by the centred
construction in Lemma~\ref{l.cell-rectangular-core}, and set
\[
 R_L:=\{u\in\Z^d:V(Q_u)\subseteq\mathrm{Pr}(L)\}.
\]
Since the prism and the $\kappa$-boxes are axial, $R_L$ is a deterministic
index rectangle.  Its side lengths are of order $K/\kappa$, with dimensional
constants, and all its $\kappa$-boxes are subgraphs of~$L$ and lie a fixed
positive fraction of~$K$ inside~$\Lambda$.

\begin{lemma}[Crossing a boxlet with the $\kappa$-array]
\label{l.kappa-crossing}
There are $C,c>0$ such that, uniformly over the boxlets $L$ of~$B$,
$\PP(\mathsf H_{R_L}(1/8)^{\rm c})\le Ce^{-cK/\kappa}$.  On
$\mathsf H_{R_L}(1/8)$, the woven fabric $\mc W(\mathscr D_{R_L})$ contains
an open path $P_L\subseteq L$ of Euclidean diameter at least~$3K/4$.
\end{lemma}

\noindent{\bf Proof.} The probability estimate follows by applying
Lemma~\ref{l.high-density-rectangles} with $\alpha=1/8$, since
the minimum index side length of $R_L$ is of order $K/\kappa$
and its aspect ratio is bounded in terms of $d$.  On $\mathsf H_{R_L}(1/8)$,
$\mathscr D_{R_L}$ meets both opposite
$e_d$-faces of~$R_L$.  Choose boxes $Q_-$ and $Q_+$ of this component on
those faces.  As the prism has $e_d$-side length~$K$, their indices are
separated in the $e_d$-direction by $K/\kappa-O_d(1)$.

The family $\mathscr D_{R_L}$ is nearest-neighbour, hence $*$-connected,
so Lemma~\ref{l.kappa-chaining} makes its woven fabric connected.
Choosing $x_\pm\in\mc W(Q_\pm)$, we have
$|(x_+-x_-)\cdot e_d|\ge K-C_d\kappa$.
An open path in the woven fabric joining $x_-$ to $x_+$ therefore has
Euclidean diameter at least $K-C_d\kappa\ge3K/4$.
Every box in this fabric is a subgraph of~$L$, so the path lies in~$L$.
\qed

Thus, for each boxlet, we obtain a long open path carried by a
large apt component in an interior rectangular subarray.  A priori these
components may differ from one boxlet to the next.  We now anchor a
common apt component near the centre of $\Lambda$ and use interior joining
to connect the boxlet rectangles to it through the jewel.

Let $z_\ast$ be the centre of the largest inscribed ball in $\Lambda$,
and let $R_{\rm cen}$ be the $\kappa$-index rectangle of all boxes $Q_u$ whose vertices lie in the
axis-parallel cube of side $K/\sqrt d$ centred at $z_\ast$.  By
part~\eqref{l.layer-box-geometry.inball} of Lemma~\ref{l.layer-box-geometry},
this cube lies in the $K/2$-erosion of $\Lambda$, namely the set of points
of $\Lambda$ at distance at least $K/2$ from $\partial\Lambda$.  For $K/\kappa$ large,
the side lengths of $R_{\rm cen}$ are comparable to~$K/\kappa$,
with constants depending only on $d$.

Let $\mathscr G_{R_{\rm cen}}$ be the $*$-component supplied by
Lemma~\ref{l.high-density-rectangles}.  If it is non-empty, define
$\mathscr J^\ast$ to be the $*$-connected component of apt $\kappa$-boxes
in the jewel that contains it; otherwise put
$\mathscr J^\ast=\emptyset$.  Put
$\mc W=\mc W(\mathscr J^\ast)$ when $\mathscr J^\ast\ne\emptyset$ and
$\mc W=\emptyset$ otherwise.  Finally, when $\mc W\ne\emptyset$, let
$\mc C$ be the open cluster of $B$ containing $\mc W$, and put
$\mc C=\emptyset$ otherwise.  Thus the global component, its woven fabric and
the open cluster of $B$ containing it are defined on the whole configuration
space.

The next lemma establishes this joining through a bounded chain of
overlapping high-density rectangles in the jewel.

\begin{lemma}[Interior joining]
\label{l.interior-joining}
Fix $0<b<B_0<\infty$.  Let $R$ be a deterministic axis-parallel
$\kappa$-index rectangle whose side lengths lie between $bK/\kappa$
and $B_0K/\kappa$, and suppose that every $\kappa$-box indexed by
$R$ lies at Euclidean distance at least $bK$ from $\partial\Lambda$.
There are $C,c>0$, depending on the fixed bounds $b,B_0$ but uniform
in $K$, $\vecell$, $B$ and $R$, and an event $\mathsf{Join}(R)$ such that
\[
   \PP\bigl(\mathsf{Join}(R)^c\bigr)\le Ce^{-cK/\kappa},
\]
and, on $\mathsf{Join}(R)$,
\[
   \mathsf H_R(1/8)\cap\mathsf H_{R_{\rm cen}}(1/8)\text{ occurs}
   \qquad\hbox{and}\qquad
   \mathscr G_R\subseteq\mathscr J^\ast .
\]
\end{lemma}

\noindent{\bf Proof.}
Write $R$ and $R_{\rm cen}$ as products of index intervals and
interpolate linearly between their corresponding lower and upper
endpoints, rounding to integer endpoints.  Put
$b_0=\min\{b,1/(2\sqrt d)\}$.  Both rectangles have sides bounded above and
below by positive multiples of $K/\kappa$, depending only on $d,b,B_0$,
and their endpoints differ by at most $C_dK/\kappa$.  We may choose a
fixed number $m=m(d,b,B_0)$ of interpolation steps so that each endpoint
moves by at most $b_0K/(20\kappa)+O(1)$ per step.  For sufficiently
large $K/\kappa$, the successive rectangles
$R=R_0,\ldots,R_m=R_{\rm cen}$ have comparable side lengths, and each
intersection $O_j=R_{j-1}\cap R_j$ is a rectangle whose sides are at
least half the corresponding sides of either adjoining member.  Thus
their aspect ratios are bounded in terms of~$d$, $b$ and $B_0$.

Before rounding, the coordinate bounds of each intermediate rectangle
are obtained by linearly interpolating those of $R$ and $R_{\rm cen}$.
Its physical support is therefore contained in the convex hull of the two
endpoint supports.  Both lie in the $b_0K$-erosion of $\Lambda$, which
is convex.  Rounding the index endpoints changes the supports by at most
$O_d(\kappa)$.  Hence, for large $K/\kappa$, every intermediate
$\kappa$-box lies in the $(b_0/2)K$-erosion and belongs to the jewel.

Apply Lemma~\ref{l.high-density-rectangles} with $\alpha=1/8$ to
every $R_j$ and every $O_j$, taking $n$ to be the minimum index side
length in each case.  On the intersection of these good events,
$\mathscr D_{O_j}$ crosses $O_j$, and thus has index diameter at least
half the minimum side length of either adjoining rectangle, up to an
additive constant.  For sufficiently large $K/\kappa$, its $*$-component
in either adjoining rectangle is the
distinguished component supplied by
part~\eqref{l.high-density-rectangles.star} of Lemma~\ref{l.high-density-rectangles}.
Hence $\mathscr G_{R_{j-1}}$ and $\mathscr G_{R_j}$ belong to the same full
$*$-connected apt component of the jewel.  Chaining over $j$ shows that
$\mathscr G_R$ belongs to the component containing
$\mathscr G_{R_{\rm cen}}$, namely $\mathscr J^\ast$.  A union bound over the
bounded number of rectangle events gives the asserted probability estimate.
\qed

\begin{proposition}[The woven fabric reaches every boxlet]
\label{p.nonaxial-woven}
There is an event $\mathcal A\in\sigma(E(B))$ with
$\PP(\mathcal A^c)\le Ce^{-cK/\kappa}$, on which
$\mathscr J^\ast\ne\emptyset$ and, for every boxlet $L$ of $B$, the woven
fabric $\mc W$ contains an open path $P_L\subseteq L$ of Euclidean diameter
at least $3K/4$.
\end{proposition}

\noindent{\bf Proof.}
For each boxlet $L$, use the deterministic rectangle $R_L$ defined
before Lemma~\ref{l.kappa-crossing}.  Its side-length and interior-margin
bounds allow us to apply Lemma~\ref{l.interior-joining} with fixed
dimensional bounds $b,B_0$.
On $\mathsf{Join}(R_L)$, the event $\mathsf H_{R_L}(1/8)$ occurs, so
Lemma~\ref{l.kappa-crossing} supplies the path $P_L$ inside the woven
fabric of $\mathscr D_{R_L}$.  Since
$\mathscr D_{R_L}\subseteq\mathscr G_{R_L}\subseteq\mathscr J^\ast$,
the definition of the woven fabric places $P_L$ in $\mc W$.
There are $2^{d-1}$ boxlets, so intersecting these joining events
gives an event $\mathcal A$ satisfying the stated probability bound.  Every
box involved lies in the jewel, and hence $\mathcal A\in\sigma(E(B))$.
Thus, on $\mathcal A$, the cluster $\mc C$ satisfies the first aptness
requirement of Definition~\ref{d.apt-box}. \qed

\section{The layer sweep and absorption}

The global woven fabric now reaches every boxlet.  To establish absorption,
we follow the progress of a possible large cluster through successive edge
layers.  At each new layer, its explored history determines where it arrives,
while the edges of that layer offer a fresh chance to connect to a local woven
fabric.  We estimate these chances conditionally on the preceding exploration,
and establish separately that the local fabrics join the global one.

\subsection{Working layers and capture}
\label{ss.nonaxial-layers}

For $\varphi\in\Xi$ and an edge $e=\{x,y\}$ put
$\varphi^+(e)=\max\{\varphi(x),\varphi(y)\}$. For $s<t$ let
$\Theta_\varphi(s,t)$ be the half-open edge layer, namely the subgraph
of $B$ with edge set
\begin{equation}
\label{e.defining-edge-layer}
 E\bigl(\Theta_\varphi(s,t)\bigr)
   \, = \, \{e\in E(B):\ s<\varphi^+(e)\le t\},
\end{equation}
and vertex set consisting of the endpoints of these edges.
These are deterministic, and for a partition of the range of $\varphi^+$ into
half-open intervals their edge sets partition $E(B)$.

Fix an integer $W_0>2(\sqrt d+2)$, and for each $\varphi\in\Xi$ define
the deterministic $\varphi$-stack by
\[
   \Theta_{\varphi,j}
   =\Theta_\varphi\bigl(jW_0\kappa,(j+1)W_0\kappa\bigr),
   \qquad j\in\Z .
\]
Its members are the \emph{working layers} in direction $\varphi$; their edge
sets partition $E(B)$.  Fix the dimensional margin
$a_\star=c_{\rm span}/8$.  A working layer
$\Theta_\varphi(s,t)$ is \emph{usable} if
\[
   \inf_\Lambda\varphi+a_\star K\le s<t
   \le\sup_\Lambda\varphi-a_\star K .
\]
Thus the $2d$ stacks, and the collection of usable layers in each, are fixed
before the percolation configuration is examined.  For $u\in\Z^d$, designate
$z_u=(\kappa(u_i+1))_{i=1}^d$ as the central lattice vertex of~$Q_u$.

\begin{lemma}[Section arrays]
\label{l.section-arrays}
For every usable layer $\Theta$ there is a rectangular parameter box
$R_\Theta\subseteq\Z^{d-1}$, all of whose side lengths are of order
$K/\kappa$, and an injective map
$F_\Theta:R_\Theta\longrightarrow\Z^d$.  We call
$S(\Theta):=F_\Theta(R_\Theta)$ the \emph{section array} of $\Theta$.
The following properties hold:
\begin{enumerate}
\item\label{l.section-arrays.contained}
$S(\Theta)$ consists of indices $u$ for which
$Q_u$ lies in $\Theta$ and in the jewel;
\item\label{l.section-arrays.lipschitz}
nearest-neighbour indices of $R_\Theta$ are sent to
$*$-neighbouring indices of the $\kappa$-array; moreover, for
$r,r'\in R_\Theta$,
\[
 C_d^{-1}\|r-r'\|_\infty
 \;\le\;
 \|F_\Theta(r)-F_\Theta(r')\|_\infty
 \;\le\;
 C_d\|r-r'\|_\infty ;
\]
\item\label{l.section-arrays.core}
there is a central subrectangle $R_\Theta^\circ\subseteq R_\Theta$, with all
side lengths at least $b_dK/\kappa$, such that the boxes indexed by
$S^\circ(\Theta):=F_\Theta(R_\Theta^\circ)$ lie in the $c_dK$-erosion of
$\Lambda$ and their union has Euclidean diameter at most
$c_dK/(10\sqrt d)$, for fixed $b_d,c_d>0$ depending only on $d$;
\item\label{l.section-arrays.corebox}
there is an axis-parallel $d$-dimensional $\kappa$-index rectangle
$\widehat R_\Theta$ which contains $S^\circ(\Theta)$, has all side lengths of
order $K/\kappa$, and all of whose $\kappa$-boxes lie in a fixed
positive-fraction erosion of $\Lambda$.
\end{enumerate}
\end{lemma}

\noindent{\bf Proof.} For every vertex $v$ of $Q_u$, the designated
central vertex $z_u$ satisfies $|v_i-z_{u,i}|\le\kappa$ for each~$i$.
Consequently
\[
  |\rho_d(v)-\rho_d(z_u)|
     \le \kappa\sum_{i=1}^d\ell_i\le\sqrt d\,\kappa,
  \qquad
  |\rho_i(v)-\rho_i(z_u)|\le\kappa\quad(i<d).
\]
We use these bounds to keep every vertex of each selected box inside the
layer, and hence every edge of the box in $E(\Theta)$.

First take $\varphi=\varepsilon\rho_d$, $\varepsilon\in\{-1,1\}$, and write
$m=(s+t)/2$ for the midpoint of the usable layer
$\Theta=\Theta_\varphi(s,t)$.  Let
$U'_\Theta\subseteq\Z^{d-1}$ be the rectangular box of transverse
$\kappa$-indices for which the full transverse span of a $\kappa$-box lies in
the transverse footprint of $\Lambda$; each side has length
$2K/\kappa+O(1)$.  For $u'\in U'_\Theta$, choose the integer $u_d=u_d(u')$
by the floor rule
\[
   0\le\varepsilon m-\rho_d(z_{(u',u_d)})<\ell_d\kappa.
\]
Thus $|\varphi(z_{(u',u_d)})-m|<\kappa$, and every vertex $v$ of this
$Q_{(u',u_d)}$ satisfies
$|\varphi(v)-m|<(1+\sqrt d)\kappa<W_0\kappa/2$.
The box therefore lies in $\Theta$.  Its transverse span is inside
$\Lambda$ by construction, and its height span is inside $\Lambda$
because the layer is usable.  Hence the box belongs to the jewel.
Changing one component of $u'$ by one changes $u_d(u')$ by at most one,
because $\ell_i/\ell_d\le1$.  Thus nearest neighbours are sent to
$*$-neighbours.  Moreover the first $d-1$ coordinates of
$F_\Theta(u')=(u',u_d(u'))$ are exactly~$u'$, while
\[
 |u_d(u')-u_d(v')|
 \le \|u'-v'\|_1
 \le (d-1)\|u'-v'\|_\infty .
\]
This proves part~\eqref{l.section-arrays.lipschitz} in this case.

Suppose instead that $\varphi=\varepsilon\rho_i$, $i<d$.  Fix the
$i$\textsuperscript{th} $\kappa$-index by the floor rule
\[
   0\le\varepsilon m-z_{u,i}<\kappa .
\]
This choice is independent of the other indices.  Every vertex of the
resulting boxes has $|\varphi(v)-m|<2\kappa<W_0\kappa/2$, so their
edges lie in $\Theta$.  Since the layer is usable, the full
$i$\textsuperscript{th} transverse span of each box lies in $\Lambda$
for sufficiently large $K/\kappa$.  Let $z$ denote the remaining $d-2$
transverse indices, restricted so that their full $\kappa$-box spans lie in the transverse
footprint.  They form a rectangular box $Z_\Theta$, with each side of
length $2K/\kappa+O(1)$.

For a fixed $z$, consider integers $u_d$ satisfying
\[
 \inf_\Lambda\rho_d+\sqrt d\,\kappa
 <\rho_d(z_{(u',u_d)})
 \le\sup_\Lambda\rho_d-\sqrt d\,\kappa,
\]
where $u'$ combines the fixed $i$\textsuperscript{th} index and $z$.
The same inequalities determine a real interval $I(z)$; every integral
$u_d\in I(z)$ yields a $\kappa$-box inside the height interval of
$\Lambda$ and hence, by the transverse restrictions, in the jewel.
The interval $I(z)$ has length
$(2K-2\sqrt d\,\kappa)/(\ell_d\kappa)$, independent of $z$ and of order
$K/\kappa$.  Its endpoints have the same $z$-dependent translation, so
the number $N(z)=|I(z)\cap\Z|$ takes at most two consecutive values.
Changing $z$ by one coordinate step translates $I(z)$ by at most one,
because $\ell_j/\ell_d\le1$; hence its smallest integer $b(z)$ also
changes by at most one.  Put $M=\min_{z\in Z_\Theta}N(z)$ and retain the
levels $b(z),b(z)+1,\ldots,b(z)+M-1$.  At most one admissible level is
discarded for any fixed $z$, and $M\asymp K/\kappa$.  With
$R_\Theta=Z_\Theta\times\llbracket0,M-1\rrbracket$, define
$F_\Theta(z,h)$ by keeping the coordinates $z$, keeping the fixed
$i$\textsuperscript{th} index, and setting $u_d=b(z)+h$.
All these boxes belong to the jewel and lie in $\Theta$;
nearest neighbours again map to $*$-neighbours.

To check the metric estimate, let $r=(z,h)$, $r'=(z',h')$ and put
$D=\|F_\Theta(r)-F_\Theta(r')\|_\infty$.  Since the $z$-coordinates are
retained, $\|z-z'\|_\infty\le D$, while
$|b(z)-b(z')|\le\|z-z'\|_1$.  Hence
\[
 |h-h'|
 \le |u_d-u_d'|+|b(z)-b(z')|
 \le (d-1)D .
\]
Conversely,
\[
 |u_d-u_d'|
 \le |b(z)-b(z')|+|h-h'|
 \le (d-1)\|r-r'\|_\infty .
\]
This proves part~\eqref{l.section-arrays.lipschitz}, after enlarging $C_d$ if
necessary.

Finally take a sufficiently small fixed central fraction of each coordinate
range of $R_\Theta$.  Since the layer is usable, its defining-coordinate
range has a margin of at least $a_\star K-O_d(\kappa)$ from the
corresponding faces of $\Lambda$.  In each remaining defining coordinate,
central restriction gives a margin bounded below by a positive dimensional
multiple of $K$.  Hence the boxes indexed by the central patch lie in the
$c_dK$-erosion for some fixed $c_d>0$.  By the metric estimate in
part~\eqref{l.section-arrays.lipschitz}, the fraction can be chosen small enough
that the union of their physical supports has Euclidean diameter at most
$c_dK/(10\sqrt d)$.  The central parameter rectangle still has every
side of length at least $b_dK/\kappa$ for some fixed $b_d>0$, proving
part~\eqref{l.section-arrays.core}.

An axial bounding cube of this union has side at most
$c_dK/(10\sqrt d)$, up to $O_d(\kappa)$ rounding.  Enclose it, with
room for the full $\kappa$-boxes, in an axial cube of side
$c_dK/(5\sqrt d)$.  Every point of the larger cube is within
$c_dK/5+O_d(\kappa)$ of a point of the central patch, so the cube
lies in the $(c_d/2)K$-erosion for large $K/\kappa$.  The indices of
$\kappa$-boxes contained in this cube form a rectangle
$\widehat R_\Theta$, containing $S^\circ(\Theta)$ and having all side
lengths of order $K/\kappa$.  This proves
part~\eqref{l.section-arrays.corebox} and gives uniform dimensional bounds on
all the aspect ratios used above. \qed

Use $F_\Theta$ to regard the aptness indicators on $S(\Theta)$ as a
site field on the rectangular index set $R_\Theta$.  Restriction and
relabelling preserve stochastic domination, so by
part~\eqref{l.kappa-apt.dom} of Lemma~\ref{l.kappa-apt}, this field
dominates Bernoulli site percolation of parameter $1-\delta_{\rm hd}$ on $R_\Theta$.
In applying Lemma~\ref{l.high-density-rectangles}, take an aspect-ratio
bound from Lemma~\ref{l.section-arrays}, and choose
$\alpha_{\rm sec}\in(0,1/4)$ in a manner to be specified shortly.  Let
$\mathscr D_\Theta$ be the component $\mathscr D_{R_\Theta}$ of that
lemma for this pulled-back field, and define
$\mathscr J_\Theta=F_\Theta(\mathscr D_\Theta)$ and
$\mc W_\Theta=\mc W(\mathscr J_\Theta)$, with both objects empty when
$\mathscr D_\Theta=\emptyset$.  By
part~\eqref{l.section-arrays.lipschitz}, $\mathscr J_\Theta$ is $*$-connected.
Lemma~\ref{l.high-density-rectangles} gives, uniformly in the layer,
$\PP(\mathscr J_\Theta=\emptyset)\le Ce^{-cK/\kappa}$.  It also gives a
constant $q_{\rm sec}>0$, independent of $K$, of the layer
and of the prescribed index, such that
\begin{equation}
\label{e.site-connection}
   \PP\bigl(v\in\mathscr J_\Theta\bigr)\ge q_{\rm sec},
   \qquad v\in S(\Theta).
\end{equation}
Indeed this is~\eqref{e.high-density-rooted} applied to
$F_\Theta^{-1}(v)$.  This fixed positive lower bound will yield a uniform
chance of connection at each fresh layer of the sweep.

It remains to show that, wherever a cluster arrives at a usable layer, it
can be steered within that layer to the fresh section array just
constructed, using $O(\kappa)$ edges.

\begin{lemma}[Access and usable-layer count]
\label{l.useful-layers}
There are $C_0,c_0>0$, depending only on $d$, such that for
$K\ge K_0(d,\kappa)$:
\begin{enumerate}
\item\label{l.useful-layers.geom}
\textup{(Access.)} For every usable layer $\Theta$ and every vertex $w$
incident to an edge of $\Theta$, there is a path
$\gamma(w)$, determined by $w$ and $\Theta$ alone, of at most $C_0\kappa$
edges of $\Theta$, running from $w$ to the central vertex
$z_{u(w)}$ of a $\kappa$-box $Q_{u(w)}$ with $u(w)\in S(\Theta)$.
Every $\kappa$-box of $S(\Theta)$
whose edge set meets $\gamma(w)$ has index within
$C_1=C_0+2$ of $u(w)$.
\item\label{l.useful-layers.cross}
\textup{(Usable-layer count.)} Every connected $H\subseteq B$ with
$\operatorname{diam}H\ge K/2$ has $\varphi$-extent at least
$c_{\rm span}K$ for some $\varphi\in\Xi$.
For every $\varphi\in\Xi$ and $x\in V(B)$ satisfying
\[
   \varphi(x)\le\sup_\Lambda\varphi+1-c_{\rm span}K,
\]
the interval
$[\varphi(x)+W_0\kappa,\,\varphi(x)+c_{\rm span}K/2]$
contains the lower $\varphi$-levels of a consecutive run of at least
$c_0K/\kappa$ usable working layers.  In particular this applies when
$x$ is a $\varphi$-minimal vertex of a connected $H\subseteq B$ whose
$\varphi$-extent is at least $c_{\rm span}K$.
\end{enumerate}
\end{lemma}

\noindent{\bf Proof: (1).} Write $\Theta=\Theta_\varphi(s,t)$, so
$t-s=W_0\kappa$.  We first reach, in at most two edges of $\Theta$, a vertex
$y\in\Lambda$ with $s<\varphi(y)\le t$.  If $w\in\Lambda$ and
$\varphi(w)>s$, take $y=w$.  Otherwise choose, by a fixed rule, an edge of
$\Theta$ incident to~$w$.  If it has an endpoint in $\Lambda$ with
$\varphi>s$, take that endpoint, traversing the edge if necessary.

In the remaining cases the chosen edge is $\{z,o\}$ with $z\in\Lambda$,
$o=z-e_j\notin\Lambda$, $\varphi(z)\le s<\varphi(o)\le t$.
Because its lower endpoint has the larger $\varphi$-value, necessarily
$\varphi=-\rho_d$ or $\varphi=-\rho_j$ with $j<d$.
For $\varphi=-\rho_d$, the usability of the layer keeps $o$ away from
the lower height face, so $o$ lies outside through a transverse face and
$j<d$.  The step $z\mapsto z-e_d$ stays in $\Lambda$, raises
$\varphi$ by $\ell_d\ge\ell_j$, and crosses level~$s$.
For $\varphi=-\rho_j$, the layer is away from the transverse faces, so
$o$ can only lie below the lower height face.  Now
$o\cvl>a_k-\ell_j$ and $o\cvl\le a_k$, whence
$(o+e_d)\cvl>a_k-\ell_j+\ell_d\ge a_k$ and
$(o+e_d)\cvl\le a_k+\ell_d<b_k$.  Thus $o+e_d\in\Lambda$ and
$\varphi(o+e_d)=\varphi(o)\in(s,t]$.  The route $z,o,o+e_d$
uses two edges of $\Theta$.  This supplies the required $y$ from either
endpoint of the chosen edge; every edge used has
$s<\varphi^+\le t$.

Suppose first that $\varphi=\pm\rho_i$ with $i<d$.  Choose a transverse
index tuple $u'$ as in the construction of $S(\Theta)$: its
$i$\textsuperscript{th} index is the fixed one determined by the layer,
and each other transverse index is chosen so that its box centre is
nearest to the corresponding coordinate of $y$, among permitted indices.
Let $x'_*$ be the transverse coordinates of $z_{(u',u_d)}$;
they do not depend on $u_d$.
Then $\|x'_* - y'\|_1=O_d(\kappa)$, and both $\varphi(y)$ and the signed
$i$\textsuperscript{th} coordinate of $x'_*$ lie in $(s,t]$.

Move from $y'$ to $x'_*$ by a shortest transverse lattice path, changing
coordinates in a fixed order and moving the $i$\textsuperscript{th}
coordinate monotonically.  Lift it from $y$ using
part~\eqref{l.layer-box-geometry.fibres} of Lemma~\ref{l.layer-box-geometry}.
The lift stays in $\Lambda$ and uses at most twice as many edges;
its vertical corrections do not change $\varphi$, while each transverse
step keeps the signed $i$\textsuperscript{th} coordinate between its
initial and final values.  Hence every edge of the lift lies in $\Theta$.
At its final transverse position, the central heights of the selected
boxes of $S(\Theta)$ are spaced by $\kappa$ and reach within
$O_d(\kappa)$ of both ends of the vertical fibre: the admissible levels
were trimmed only by the fixed height margins and at most one index.
Choose $u_d$ whose central vertex is within $C_d\kappa$ vertical steps
of the lifted endpoint.  A final vertical path inside that fibre reaches
$z_{(u',u_d)}$, and all its edges lie in $\Theta$, since $\varphi$ is
unchanged by vertical moves.

Now suppose that $\varphi=\pm\rho_d$, and put $m=(s+t)/2$.  Set
$I=(m-\ell_d/2,\,m+\ell_d/2]$.
For every transverse lattice point $x'$ in the footprint of $\Lambda$,
there is exactly one integer level $h_I(x')$ with
$\varphi(x',h_I(x'))\in I$.  The resulting vertex lies in $\Lambda$,
since the usable layer is separated from both height faces.
First move vertically from $y$ to $(y',h_I(y'))$; this takes
$O_d(\kappa)$ steps, and all intermediate vertices are in $\Lambda$
and have $\varphi\in(s,t]$.

Choose transverse coordinates $x'_*$ of a box centre in the transverse
index range of $S(\Theta)$ by choosing the nearest permitted box
centres, as above; thus $\|x'_* -y'\|_1=O_d(\kappa)$.  Follow a shortest transverse path to
$x'_*$, after each transverse step making at most one vertical step to
return to the unique level prescribed by~$I$.  Indeed, a transverse step
changes $\varphi$ by at most $\ell_d$, while a vertical step changes it
by exactly $\ell_d$ in absolute value.  The interval $I$ lies well within
$(s,t]$, so both the transverse steps and their vertical corrections lie
in $\Theta$ and $\Lambda$.  At the final transverse position the floor
rule for $S(\Theta)$ supplies an index $u_d$ with
$|\varphi(z_{(u',u_d)})-m|<\kappa$.  The present height belongs to~$I$,
so the vertical separation from $z_{(u',u_d)}$ is at most
$(\kappa+\ell_d)/\ell_d=O_d(\kappa)$.  Move vertically within the fibre
to this central vertex; the path remains in $\Theta$ and $\Lambda$.

Make all choices of incident edges, rounding, lifts and nearest admissible
indices by fixed rules.  The two constructions then give a path
$\gamma(w)$ of at most $C_0\kappa$ edges.  Since a $\kappa$-box
has index diameter two,
any $Q_v$, $v\in S(\Theta)$, whose edge set meets $\gamma(w)$ has
index within $C_0+2=C_1$ of $u(w)$.

{\bf (2).} Choose two vertices of $H$ at Euclidean distance at least
$K/2$.  Part~\eqref{l.layer-box-geometry.bilip} of Lemma~\ref{l.layer-box-geometry}
gives $\varphi\in\Xi$ for which their $\varphi$-separation, and hence
the $\varphi$-extent of $H$, is at least $c_{\rm span}K$.

For the counting statement, fix $\varphi$ and $x$ as in the lemma,
and put $a=\varphi(x)$, $L=\inf_\Lambda\varphi$ and
$U=\sup_\Lambda\varphi$.  Since $x\in V(B)$,
$L-1\le a\le U+1-c_{\rm span}K$.  The allowed lower levels of usable
working layers lie in
$[L+a_\star K,\,U-a_\star K-W_0\kappa]$.
For large $K/\kappa$, the interval
\[
 [a+W_0\kappa,\,a+c_{\rm span}K/2]
 \ \cap\ [L+a_\star K,\,U-a_\star K-W_0\kappa]
\]
has length at least
$(c_{\rm span}/2-a_\star)K-O_d(\kappa)$.
Here the term $W_0\kappa$ in the upper bound accounts for the full
width of a working layer.  Since $a_\star=c_{\rm span}/8$ and consecutive
lower levels in the stack are spaced by $W_0\kappa$, this
intersection contains the lower levels of at least $c_0K/\kappa$
consecutive usable working layers, for a constant $c_0>0$ depending
only on~$d$.  Finally, if $x$ is $\varphi$-minimal in $H$ and
$H$ has $\varphi$-extent at least $c_{\rm span}K$, then
$a+c_{\rm span}K\le\sup_H\varphi\le U+1$;
so the counting condition holds. \qed

The next elementary observation is the point at which finite energy can be
used without asking a non-monotone aptness event to survive arbitrary
conditioning.

\begin{lemma}[Safe local opening]
\label{l.safe-kappa-opening}
Let $\mathscr H$ be a finite family of $\kappa$-boxes and obtain
$\omega^{\mathscr H}$ from $\omega$ by declaring every edge of every box in
$\mathscr H$ open. Then every $\kappa$-box that is apt in $\omega$ remains
apt in $\omega^{\mathscr H}$; every box of $\mathscr H$ is apt in
$\omega^{\mathscr H}$.
\end{lemma}

\noindent{\bf Proof.} A box of $\mathscr H$ is fully open and hence apt.  Let $Q$ be apt before the
opening. If its edge set is untouched there is nothing to prove. Otherwise a
changed edge of $Q$ lies in $Q\cap Q'$ for some $Q'\in\mathscr H$ with
$Q'$ a $*$-neighbour of $Q$.  The overlap is a connected union of
constituent $\kappa$-boxlets of $Q$, and $\mc W(Q)$ contains an open path
of Euclidean diameter at least $3\kappa/4$ in every such boxlet. Once the
overlap is made open, every changed
edge of $Q$ is therefore connected inside $Q$ to $\mc W(Q)$. The old
boxlet crossings remain, and every new open path of diameter at least
$\kappa/2$ that uses a changed edge belongs to the enlarged cluster
containing $\mc W(Q)$; the old long paths belonged to it already. Thus $Q$
remains apt. \qed

There are now two distinct objects with complementary roles.  The section
component $\mathscr J_\Theta$ is defined entirely from the fresh field in one
working layer and is therefore the correct target for the sequential capture
estimate.  The global component $\mathscr J^\ast$ was anchored in the central
jewel rectangle.  The two parts of the next lemma first use the local object
for capture and then identify it with the global one.

\begin{lemma}[Layer fabrics]
\label{l.layer-fabrics}
\leavevmode
\begin{enumerate}
\item\label{l.layer-fabrics.capture}
\textup{(Capture.)} Let $\Theta$ be a usable layer with lower interface at
$\varphi$-level $t$, and let $\mathcal F$ be the $\sigma$-field generated by
the edges $e\in E(B)$ with $\varphi^+(e)\le t$. Let $w$ be an
$\mathcal F$-measurable vertex incident to an edge of $\Theta$. Then there is
$c_{\rm cap}=c_{\rm cap}(d,p,\kappa)>0$ such that
\[
   \PP\bigl(w\ \text{is joined to}\ \mc W_\Theta\ \text{using edges of}\
       \Theta\ \bigm|\ \mathcal F\bigr)\,\ge\,c_{\rm cap}.
\]
\item\label{l.layer-fabrics.weave}
\textup{(Weave.)} There is an event $\mathsf{Weave}$ with
$\PP(\mathsf{Weave}^c)\le Ce^{-cK/\kappa}$ on which the event $\mathcal A$ of
Proposition~\ref{p.nonaxial-woven} occurs and, simultaneously for every usable
layer $\Theta$ in every one of the $2d$ stacks, $\mathscr J_\Theta$ is non-empty and
$\mathscr J_\Theta\subseteq\mathscr J^\ast$.
\end{enumerate}
\end{lemma}

\noindent{\bf Proof: (1).} Every edge of $\Theta$ is fresh given
$\mathcal F$. Let $\gamma=\gamma(w)$ and $u=u(w)$ be supplied by
part~\eqref{l.useful-layers.geom} of Lemma~\ref{l.useful-layers}.  Every
$v\in S(\Theta)$ for which the edge set of $Q_v$ meets~$\gamma$ lies within
index distance~$C_1$ of~$u$.  By
part~\eqref{l.section-arrays.lipschitz} of Lemma~\ref{l.section-arrays}, their
preimages under $F_\Theta$ lie within a fixed distance, depending only on
$d$, of $F_\Theta^{-1}(u)$.  Choose by a fixed rule a nearest-neighbour
connected family of parameter indices, of cardinality bounded in terms of
$d$ alone, which contains these preimages, and let $\mathscr H(w)$ be its
image under $F_\Theta$.  Then $\mathscr H(w)\subseteq S(\Theta)$ is
$*$-connected and contains~$u$ and every such~$v$.  All these choices are
$\mathcal F$-measurable.

Let $\mathcal E=\{u\in\mathscr J_\Theta\}$.  Since every $\kappa$-box
indexed by $S(\Theta)$ has its edges in $\Theta$, the apt field on this
array retains its unconditioned law given $\mathcal F$;
hence~\eqref{e.site-connection} gives
$\PP(\mathcal E\mid\mathcal F)\ge q_{\rm sec}$.

Now perform the insertion in two stages.  First open every edge of every box
$Q_v$ with $v\in\mathscr H(w)$.  By
Lemma~\ref{l.safe-kappa-opening}, every previously apt $\kappa$-box
remains apt and all boxes indexed by $\mathscr H(w)$ become apt.  On
$\mathcal E$, therefore, the old dense nearest-neighbour component
$\mathscr D_\Theta$, which contains $F_\Theta^{-1}(u)$, remains open and
connected.  Its size remains greater than $7|R_\Theta|/8$, so the
nearest-neighbour component containing it is the distinguished component
$\mathscr D_\Theta$ for the modified configuration as well.

Now also open every edge of~$\gamma$.  By the choice of $\mathscr H(w)$, if
$v\in S(\Theta)\setminus\mathscr H(w)$, then $Q_v$ has no edge in~$\gamma$;
hence this second stage changes the aptness of no such box, while the boxes
indexed by $\mathscr H(w)$ are already fully open.  In particular, no section box loses aptness, and the same density argument
shows that the modified
$\mathscr J_\Theta$ still contains~$u$.  The number of edges prescribed is at
most a constant $M_\kappa<\infty$ depending only on $d$ and $\kappa$.  The
box $Q_u$ is fully open and the opened path $\gamma$ ends in it, so
$w$ is joined within $\Theta$ to $\mc W(Q_u)$.  Since the modified
$\mathscr J_\Theta$ still contains $u$, we have
$\mc W(Q_u)\subseteq\mc W_\Theta$.  Thus $w$ is joined to
$\mc W_\Theta$ using only edges of $\Theta$.

Conditional on $\mathcal F$, the set of at most $M_\kappa$ edges to be
opened is fixed, and these edges are independent Bernoulli-$p$ variables.
Condition further on all other layer edges.  Whenever some assignment of
the prescribed edges yields $\mathcal E$, opening all of them yields the
capture event, as shown above.  Thus, with $a_\kappa=p^{M_\kappa}>0$,
\[
   \PP\bigl(w\leftrightarrow\mc W_\Theta\text{ in }\Theta
       \bigm|\mathcal F\bigr)
   \,\ge\,a_\kappa\PP(\mathcal E\mid\mathcal F)
   \,\ge\,a_\kappa q_{\rm sec}.
\]
Take $c_{\rm cap}=a_\kappa q_{\rm sec}$.

{\bf (2).} We use the central patch
$R_\Theta^\circ$ supplied in part~\eqref{l.section-arrays.core} of
Lemma~\ref{l.section-arrays}.  Its minimum side is at least
$b_dK/\kappa$, whereas all sides of $R_\Theta$ are at most
$C_dK/\kappa$.  We may therefore choose $\alpha_{\rm sec}$ sufficiently
small in terms of~$d$ that a crossing of $R_\Theta^\circ$ has index
diameter greater than $\alpha_{\rm sec}$ times the minimum side of
$R_\Theta$.  Apply Lemma~\ref{l.high-density-rectangles} with
$\alpha=\alpha_{\rm sec}$ to both $R_\Theta$ and $R_\Theta^\circ$.
On $\mathsf H_{R_\Theta}(\alpha_{\rm sec})\cap
\mathsf H_{R_\Theta^\circ}(\alpha_{\rm sec})$, let
$\mathscr D_\Theta^\circ$ be the dense nearest-neighbour component in the
central patch.  It crosses $R_\Theta^\circ$, so its index diameter exceeds
$\alpha_{\rm sec}$ times the minimum side of $R_\Theta$.  Hence its
nearest-neighbour component in $R_\Theta$ is $\mathscr D_\Theta$ by
part~\eqref{l.high-density-rectangles.nn} of Lemma~\ref{l.high-density-rectangles}.
Consequently
$F_\Theta(\mathscr D_\Theta^\circ)\subseteq\mathscr J_\Theta$.

Let $\widehat R_\Theta$ be the interior $d$-dimensional rectangle supplied by
part~\eqref{l.section-arrays.corebox} of Lemma~\ref{l.section-arrays}.  Apply
Lemma~\ref{l.interior-joining} to $\widehat R_\Theta$,
using the fixed dimensional side-length and margin bounds from
Lemma~\ref{l.section-arrays}.  The density just noted
gives $\mathscr D_\Theta^\circ$ parameter-index diameter of order
$K/\kappa$.  By the lower bi-Lipschitz bound of
part~\eqref{l.section-arrays.lipschitz} of Lemma~\ref{l.section-arrays},
$F_\Theta(\mathscr D_\Theta^\circ)$ is a $*$-connected apt set in
$\widehat R_\Theta$ of index diameter at least
$c_dK/\kappa$.  The sides of $\widehat R_\Theta$ are at most
$C_dK/\kappa$.  Choose a fixed dimensional
$\alpha_{\rm int}\in(0,1/8)$ sufficiently small that this diameter is
at least $\alpha_{\rm int}$ times the minimum side of
$\widehat R_\Theta$.  In addition to the joining event, require
$\mathsf H_{\widehat R_\Theta}(\alpha_{\rm int})$.  On this event the
$*$-component of $F_\Theta(\mathscr D_\Theta^\circ)$ in
$\widehat R_\Theta$ is $\mathscr G_{\widehat R_\Theta}$ by
part~\eqref{l.high-density-rectangles.star} of Lemma~\ref{l.high-density-rectangles}.
The joining lemma, which uses $\mathsf H_{\widehat R_\Theta}(1/8)$,
then gives
\[
   F_\Theta(\mathscr D_\Theta^\circ)
   \subseteq \mathscr G_{\widehat R_\Theta}
   \subseteq \mathscr J^\ast .
\]
Thus $\mathscr J_\Theta$ and $\mathscr J^\ast$ meet.  Since
$\mathscr J_\Theta$ is $*$-connected and $\mathscr J^\ast$ is a full
$*$-connected component of apt $\kappa$-boxes in the jewel, we have
$\mathscr J_\Theta\subseteq\mathscr J^\ast$.

There are $O(K/\kappa)$ usable layers in the $2d$ stacks.  For each layer we
have used only a fixed number of high-density rectangle events, including
the additional event at $\alpha_{\rm int}$, together with one application
of Lemma~\ref{l.interior-joining}.  Each such event has
failure probability $Ce^{-cK/\kappa}$.  The linear union factor is absorbed by the
exponential.  Adjoining the event $\mathcal A$ of
Proposition~\ref{p.nonaxial-woven} and taking $\mathsf{Weave}$ to be the
resulting intersection proves the claim. \qed

\subsection{Absorption}
\label{ss.nonaxial-absorption}

It remains to show that every open path in $B$ of diameter at least $K/2$
is contained in $\mc C$.  Could the TipToe dichotomy of
Chapter~\ref{s.heuristics} operate in the box~$B$?  Indeed it does, but in a
simpler, degenerate form.  The thin--thick distinction collapses to empty
versus non-empty arrival at a fresh layer: if the cluster does not enter the
layer, it has petered out, while the arrival of even a single vertex
brings the prospect of open escape.

Indeed, wherever such an arrival occurs, deep inside or against the irregular
boundary of $\Lambda$, the access assertion in part~\eqref{l.useful-layers.geom}
of Lemma~\ref{l.useful-layers} supplies a short path to the fresh
$\kappa$-array, and Lemma~\ref{l.layer-fabrics} gives a fixed positive chance
of joining the layer fabric.  A cluster of diameter $K/2$ must encounter
order $K/\kappa$ such layers and evade capture at every one if it is to remain
apart from the global cluster $\mc C$.

We run the sweep without conditioning on the global woven fabric~$\mc W$.  In each
layer the target is the layer fabric $\mc W_\Theta$, which is defined from the
fresh field in that layer.  Part~\eqref{l.layer-fabrics.weave} of Lemma~\ref{l.layer-fabrics}
then identifies all these fabrics as subsets of $\mc W$.  This separates the sequential
freshness estimate from the global weaving argument.  The edge configurations
of successive layers in a fixed stack are independent.  The capture events
need not be, since their entry vertices depend on the preceding exploration;
but layer-edge independence preserves their uniform conditional bounds, which
combine by the tower property.

\begin{proposition}[Absorption]
\label{p.nonaxial-absorption}
Except on an event of probability at most $Ce^{-c_{\rm abs}K}$, every open
path in $B$ of diameter at least $K/2$ is contained in $\mc C$, where
$c_{\rm abs}=c_{\rm abs}(d,p,\kappa)>0$.
\end{proposition}

\noindent{\bf Proof.} On $\mathsf{Weave}$, the fabric $\mc W$ is non-empty
and $\mc C$ is the open cluster of $B$ containing it.  Thus any open cluster
of diameter at least $K/2$ outside $\mc C$ must be disjoint from $\mc W$.
We bound the probability of such a cluster on $\mathsf{Weave}$, then add
$\PP(\mathsf{Weave}^c)$.  By part~\eqref{l.useful-layers.cross} of Lemma~\ref{l.useful-layers},
such a cluster has
$\varphi$-extent at least $c_{\rm span}K$ for some $\varphi\in\Xi$. For a vertex
$x\in V(B)$ and $\varphi\in\Xi$, let $E_{x,\varphi}$ be the event that the
open cluster of $x$ is disjoint from $\mc W$, has $\varphi$-extent at least
$c_{\rm span}K$, and has $x$ as the first, in a fixed lattice order, of its
$\varphi$-minimal vertices. The event to be bounded is contained in the union
of the $E_{x,\varphi}$ over at most $C_dK^d$ vertices and $2d$ signed
coordinates.

Fix $x$ and $\varphi$.  If
$\varphi(x)>\sup_\Lambda\varphi+1-c_{\rm span}K$, then
$E_{x,\varphi}$ is empty.  Otherwise, by part~\eqref{l.useful-layers.cross} of Lemma~\ref{l.useful-layers},
the usable layers in the deterministic $\varphi$-stack whose lower levels lie in
$[\varphi(x)+W_0\kappa,\ \varphi(x)+c_{\rm span}K/2]$
form a consecutive run of $m\ge c_0K/\kappa$ layers.
List them as $\Theta_1,\ldots,\Theta_m$ with lower levels
$t_1<\cdots<t_m$.  This list, and $m$, depend only on $x$ and
$\varphi$, not on the percolation configuration.  By the
definition of the stack,
$t_{j+1}=t_j+W_0\kappa$; put $t_{m+1}=t_m+W_0\kappa$.  Let
$\mathcal F_j$ be generated by the edges $e\in E(B)$ with
$\varphi^+(e)\le t_j$.  These $\sigma$-fields increase, the edges of
$\Theta_j$ are fresh given $\mathcal F_j$, and they are
$\mathcal F_{j+1}$-measurable.

Let $R_j$ be the event that the open cluster of $x$, explored using only
edges counted in $\mathcal F_j$, reaches a vertex incident to an edge of
$\Theta_j$.  On $R_j$, let $w_j$ be the first such vertex in the fixed lattice
order; on $R_j^c$, let $w_j$ be the first vertex in that order incident to any
edge of $\Theta_j$.  Thus $w_j$ is always defined, is
$\mathcal F_j$-measurable and is incident to an edge of $\Theta_j$.  Let
$T_j$ be the event that $w_j$ is joined to $\mc W_{\Theta_j}$ using edges of
$\Theta_j$; the capture is decided inside $\Theta_j$, whose upper interface
is at $t_{j+1}$, so $T_j\in\mathcal F_{j+1}$, while it is fresh given
$\mathcal F_j$. By part~\eqref{l.layer-fabrics.capture} of Lemma~\ref{l.layer-fabrics},
\begin{equation}
\label{e.capture-conditional}
   \PP\bigl(T_j\bigm|\mathcal F_j\bigr)\, \ge \, c_{\rm cap} .
\end{equation}

On $E_{x,\varphi}$ every $R_j$ occurs. Indeed, the cluster of $x$ contains a
vertex of $\varphi$-level exceeding $t_j$; on a path in the cluster from $x$
to it, let $u$ be the vertex preceding the first vertex of $\varphi$-level
exceeding $t_j$. Every edge of the initial segment has both endpoints at
$\varphi$-level at most $t_j$ and so is counted in $\mathcal F_j$, whence $u$
is reached by the restricted exploration; and the edge $e$ at $u$ that leaves
that level has, by part~\eqref{l.layer-box-geometry.lip} of Lemma~\ref{l.layer-box-geometry},
$t_j<\varphi^+(e)\le\varphi(u)+1\le t_j+1<t_{j+1}$, so $e$ is an edge of
$\Theta_j$ and $u$ is incident to it. On $\mathsf{Weave}$ we have
$\mc W_{\Theta_j}\subseteq\mc W$, so no $T_j$ can occur on
$E_{x,\varphi}\cap \mathsf{Weave}$, the cluster being disjoint from~$\mc W$.
Hence
\[
   E_{x,\varphi}\cap \mathsf{Weave}
   \, \subseteq \, \bigcap_{j=1}^{m}\bigl(R_j\cap T_j^c\bigr).
\]
Peel from the top: for each $j$ the events $R_i\cap T_i^c$ with $i<j$, and
$R_j$, are $\mathcal F_j$-measurable, because
$T_i\in\mathcal F_{i+1}\subseteq\mathcal F_j$; so by the tower property
and~\eqref{e.capture-conditional},
\[
   \PP\Bigl(\bigcap_{j\le n}R_j\cap T_j^c\Bigr)
   \, \le \, (1-c_{\rm cap})\,
     \PP\Bigl(\bigcap_{j\le n-1}R_j\cap T_j^c\Bigr),
   \qquad n\le m,
\]
and therefore, since $m\ge c_0K/\kappa$, there is
$c_{\rm sweep}=c_{\rm sweep}(d,p,\kappa)>0$ such that
\[
   \PP(E_{x,\varphi}\cap \mathsf{Weave})
   \, \le \, (1-c_{\rm cap})^m \, \le \, e^{-c_{\rm sweep}K}.
\]
Summing over the at most $C_dK^d$ vertices and $2d$ coordinates costs a
polynomial factor, which is absorbed by reducing the exponential rate.
part~\eqref{l.layer-fabrics.weave} of Lemma~\ref{l.layer-fabrics} gives an error
$Ce^{-cK/\kappa}$, also exponential in $K$ because $\kappa$ is fixed.  This
gives the stated $Ce^{-c_{\rm abs}K}$ bound. \qed

\section{Obtaining the layer-box estimate}
\label{ss.nonaxial-conclusion}

We may now complete the proof.
Propositions~\ref{p.nonaxial-woven} and~\ref{p.nonaxial-absorption}
supply the two aptness requirements; it remains to combine their bounds and
verify measurability, one-dependence and the offset-layer statement.

\noindent\textbf{Proof of Proposition~\ref{p.nonaxial-ap}.}
Let $\mathsf{Bad}(B)$ be the union of the failure events of
Proposition~\ref{p.nonaxial-woven} and
Proposition~\ref{p.nonaxial-absorption}.  Since $\kappa$ is now fixed,
$\PP(\mathsf{Bad}(B))\le Ce^{-cK}$. Off $\mathsf{Bad}(B)$ the open
cluster $\mc C$ contains, in every boxlet
of $B$, an open path of diameter at least~$3K/4$, and every open path in $B$
of diameter at least $K/2$ is contained in $\mc C$: that is, $B$ is apt in
the sense of Definition~\ref{d.apt-box}, with woven cluster $\mc C$.

By Definition~\ref{d.apt-box},
$\{B\text{ is apt}\}\in\sigma(E(B))$.

Finally, boxes whose full indices are at $\ell^\infty$-distance greater than
one have disjoint edge sets: in a transverse coordinate the ranges
$(K(v_i-2),Kv_i]$ and $(K(v'_i-2),Kv'_i]$ are disjoint when $|v_i-v'_i|\ge2$,
and in the layer index the upper-endpoint ranges $(a_k,b_k]$ and
$(a_{k'},b_{k'}]$ are disjoint when $|k-k'|\ge2$. The inapt-box field is
therefore one-dependent.

For an offset layer $N(a)$, replacing the upper-endpoint interval
$(a_k,b_k]$ by $(a,a+2K]$ leaves all preceding geometric bounds unchanged.
The same construction therefore gives, uniformly in $a$, the aptness estimate
and measurability with respect to the edge configuration on $E(B)$.  Boxes in
$\overline N(a)$ whose transverse indices have $\ell^\infty$-distance greater
than one have disjoint edge sets, so the corresponding inapt-box field on
$\Z^{d-1}$ is one-dependent. \qed

The Antal--Pisztora estimate derived here is uniform over the
one-directional shear family we have considered.  We expect the estimate to
hold likewise for cubes in general alignment, with the argument requiring
some modification of the section-array and access-path constructions.

\clearpage
\chapter{Glossary}
\label{s.notation-guide}

There is quite a lot of notation in this paper.  Here we recall some of its
principal instances, each with a brief summarising phrase and a link to a
defining or principal passage.  The work has a natural conceptual division,
reflected in its four parts: single-slide analysis and an overview of the
iterated slide; the elements needed for the iterated-slide analysis; renewal
production and survival via this analysis; and the Ornstein--Zernike
consequences via classical renewal theory.  We organise the glossary under the corresponding broad headings, with
the introductory chapter included in the~first.

\medskip
\noindent\textit{Chapters~1--5: Introduction, single-slide analysis, and the iterated-slide overview}
\begin{center}
\footnotesize
\renewcommand{\arraystretch}{1.12}
\begin{tabular}{@{}p{0.29\textwidth}p{0.61\textwidth}r@{}}
\textit{Symbol or term}&\textit{Meaning}&\textit{p.}\\ \hline
\(p>p_c(d),\ \PP;\ \mc C(0)\) & supercritical parameter and law; open cluster of the origin \fff{percolation-law}
\(\vecell,\ x\cvl,\ r=r(\vecell)\) & height direction, coordinate, and maximal one-edge increment \fff{height-axis}
\(\base,\head,\cluster_I\) & finite-cluster extreme vertices and prescribed-height event \fff{base-head-cluster-event}
\(\xi,\ K_\xi,\ \mathcal W\) & inverse correlation length, equidecay body, and polar Wulff body \fff{xi-bodies}
\(\drift(\vecell)\) & conditioned cluster-head direction (Cor.~\ref{c.conditioned-head-drift}) \fff{drift}
\(F_h,B_h;\ \upspace_h,\downspace_h\) & half-spaces and their lattice subgraphs \fff{height-spaces}
\(\vertbdry_h,\slab_{y,z}\) & boundary band and slab \fff{height-spaces}
\(R(U),\ \rlevel(U)\) & renewal points of \(U\) and their heights \fff{renewal-point}
mass gap & exponential moment for the tilted sausage height &\hyperref[t.mass-gap]{\pageref*{t.mass-gap}}\\
\(\sting_u(x),\ \stingset_u(x)\) & height-\(u\) string event at \(x\), and its unique string \fff{string}
\(\clusteruthree,\ \Q\) & finite-cluster event and its conditioned law \fff{conditioned-law}
with overwhelming probability & informal usage: probability at least \(1-Ce^{-cu}\) under the law in force \fff{overwhelming-probability}
TipToe & twin-pressure mantra: \emph{thin: instantly perish; thick: openly escape} \fff{tiptoe}
Tether, \(\mathcal T\) & coarse-box path through the midpiece that anchors the upcluster shift \fff{tether}
\(B(v),\ \overline A,\ \partial_{\rm ext}\overline A\) & axial \(K\)-boxes, box support, and exterior boundary \fff{k-box}
\(\hadj\) & bounded-offset coordinate in which the upcluster is shift-invariant \fff{slide-decomposition}
\(\downcluster_h,\midpiece_h,\upcluster_h\) & down--mid--up pieces at slide height \(h\) \fff{slide-decomposition}
\(\fine_h,\ \veryfine_h\) & fine and stronger input very-fine regularity events \fff{fine-very-fine}
\(\success\) & two-edge reconstruction producing a renewal point at its midvertex \fff{simple-join}
\(\slide_h\) & single-slide kernel: resampling when fine (active), identity when dormant \fff{single-slide}
\(\Rlaw^v_{\uclust,\dclust}\) & fixed-shift reconstruction law for proposed head \(v\) \fff{fixed-shift-law}
\(\iter_{\overline h},\ \Qdyn\) & iterated-slide operator and cluster-state trajectory law \fff{dynamic-law}
\(\ivf,\ \llf\) & initially very-fine and limelight-fine sampled heights &\hyperref[ss.breaker-llf-points]{\pageref*{ss.breaker-llf-points}}\\
limelight & for a scheduled height, the moment its slide is attempted \fff{limelight}
IVF-to-LLF transition & positive fraction of IVF heights become LLF; fine delivery \fff{fine-delivery}
\(\Accept_v\) & fixed-shift event giving an admissible finite connection \fff{accept-event}
\end{tabular}
\end{center}

\Needspace{20\baselineskip}
\bigskip
\noindent\textit{Chapters~6--9: Elements for the iterated-slide analysis}
\begin{center}
\footnotesize
\renewcommand{\arraystretch}{1.12}
\begin{tabular}{@{}p{0.29\textwidth}p{0.61\textwidth}r@{}}
\textit{Symbol or term}&\textit{Meaning}&\textit{p.}\\ \hline
DMU decomposition & downcluster, midpiece, and upcluster at a given height \fff{dmu}
down-backtrack \(\dbt(v)\) & lowest ceiling beneath which \(v\) connects to the base \fff{backtrack-protection}
protection height \(\prot(A)\) & height by which every vertex of \(A\) has entered the downcluster \fff{backtrack-protection}
cut-it-and-lose-it & higher-slide exit of material entering at a lower slide \fff{cut-lose}
\(\overline N_k,\ N_k\) & transverse box lattice and its underlying physical \(\vecell\)-slab \fff{layer}
\(\Bbox{k,u}\) & axial-column boxlet \fff{overlapping-layer-box}
\(\Obox{k,v}\) & \(2K\)-wide overlapping layer box \fff{overlapping-layer-box}
apt box (in layer) & contains woven cluster crossing every boxlet and absorbing all \(K/2\)-paths \fff{apt-woven}
\(\mathcal W(B)\) & woven cluster of an apt box~\(B\) \fff{apt-woven}
\(\Mesh,\ \MeshBox z\) & disjoint $2K$-mesh of the layer and its individual mesh boxes &\hyperref[l.disjoint-twoK-mesh]{\pageref*{l.disjoint-twoK-mesh}}\\
\(\MeshSupp(H)\) & mesh boxes containing edges of~$H$, for $H$ an edge set or subgraph &\hyperref[e.mesh-support-bound]{\pageref*{e.mesh-support-bound}}\\
\(A^{\mathrm{ov}}\) & overlapping boxes contained edgewise in a mesh region~$A$ &\hyperref[l.overlap-mesh-region]{\pageref*{l.overlap-mesh-region}}\\
\(A(D)\) & exterior mesh boundary of a finite mesh region~$D$ &\hyperref[l.mesh-boundary-capture]{\pageref*{l.mesh-boundary-capture}}\\
\(\MeshInt(S)\) & captured mesh boxes lying strictly behind~$S$ &\hyperref[e.mesh-edge-partition]{\pageref*{e.mesh-edge-partition}}\\
quilt, \(\mathscr Q\); \(\mathcal W(\mathscr Q)\) & finite connected collection of apt boxes and its connected woven fabric \fff{quilt}
\(\Capt_{\overline N}(S),\ \Ext_{\overline N}(S)\) & capture and exterior of~$S$ in the layer-box lattice \fff{exterior-capture}
\(\MCap(S),\ \MExt(S)\) & captured and exterior mesh sets of~$S$ &\hyperref[e.mesh-edge-partition]{\pageref*{e.mesh-edge-partition}}\\
\(\ECap(S),\ \mathcal U(S)\) & captured and exterior edge sets, the latter unread for a seal &\hyperref[e.mesh-edge-partition]{\pageref*{e.mesh-edge-partition}}\\
\(\Halo(S)\) & exterior collar of boxes meeting the capture &\hyperref[ss.visible-welcoming-exteriors]{\pageref*{ss.visible-welcoming-exteriors}}\\
welcoming & apt boundary with halo boxes likely apt given captured-side data &\hyperref[ss.visible-welcoming-exteriors]{\pageref*{ss.visible-welcoming-exteriors}}\\
seed & vertex supplied as a starting point for further layer growth \fff{seed}
stopping domain & adaptively revealed edge set whose layer complement has product law &\hyperref[l.freshness]{\pageref*{l.freshness}}\\
seal, \(\mathscr S\) & welcoming stopping-domain quilt collection \fff{seal}
virgin or sealed seed & seed with access to the unread seal exterior, or with none &\hyperref[d.virgin-seed]{\pageref*{d.virgin-seed}}\\
activated seal component & seal component whose known woven fabric has been reached &\hyperref[d.virgin-seed]{\pageref*{d.virgin-seed}}\\
seal enlargement & same-configuration outward seal deformation capturing finite material &\hyperref[ss.efficient-extension]{\pageref*{ss.efficient-extension}}\\
initial seal & one-off enclosure seal around a prescribed finite open-and-closed record &\hyperref[ss.initial-seal]{\pageref*{ss.initial-seal}}\\
Encloser control problem & adaptive seeding that seeks finite volume without open escape &\hyperref[ss.one-life-percolation-control]{\pageref*{ss.one-life-percolation-control}}\\
wavebreaker & thin row giving a uniform chance to arrest a layer exploration \fff{wavebreaker}
\(\mathscr B_h;\ \mathscr B_h^0,\mathscr B_h^1,\mathscr B_h^2\) & breaker window and its three breaker rows \fff{breaker-window}
\end{tabular}
\end{center}

\Needspace{27\baselineskip}
\bigskip
\noindent\textit{Chapters~10--13: Renewal production and survival via the iterated-slide analysis}
\begin{center}
\footnotesize
\renewcommand{\arraystretch}{1.12}
\begin{tabular}{@{}p{0.29\textwidth}p{0.61\textwidth}r@{}}
\textit{Symbol or term}&\textit{Meaning}&\textit{p.}\\ \hline
\(\ell,\ k,\ \cS(k)\) & target height, depth below it, and associated partition layer \fff{target-layer}
advance phase & schedule phase before passage through the layer \fff{two-phases}
post-passage phase & schedule phase after passage, up to target limelight \fff{two-phases}
passage & schedule step across a monitored layer, from \(s^-\) to \(s^+\) \fff{two-phases}
host \(\RetRec\) and guest \(\MidRec\) & retained and midpiece records forming the passage record \fff{retained-midpiece}
multi-layer paused exploration & fixed-shift cluster exploration interleaved across monitored layers &\hyperref[d.multi-layer-paused]{\pageref*{d.multi-layer-paused}}\\
\(\SlideCost_b,\ \SealCost_b\) & one-slide cost and its seal-related part &\hyperref[d.ssr-one-slide-charges]{\pageref*{d.ssr-one-slide-charges}}\\
\(\Seal_j^0,\ \PassPrep_j,\ \PassExc_j\) & passage seal, selected prospective excess, and remaining seal-size excess &\hyperref[p.ssr-passage-seal]{\pageref*{p.ssr-passage-seal}}\\
depth regularity & linear count of layer--target pairs of bounded depth and breaker count &\hyperref[d.ssr-depth-regular]{\pageref*{d.ssr-depth-regular}}\\
post-passage seal replacement & fresh seal around preceding capture at each active post-passage slide &\hyperref[ss.post-passage-resealing]{\pageref*{ss.post-passage-resealing}}\\
reference dynamics & unconditioned fixed-shift dynamics, sent to a cemetery state at rejection \fff{reference-dynamics}
augmented reference dynamics & also carries sealing and paused-exploration data \fff{augmented-reference-dynamics}
marked trajectory comparison & likelihood-ratio transfer from reference to accepted trajectories &\hyperref[l.marked-whole-block-change]{\pageref*{l.marked-whole-block-change}}\\
replacement front & outward propagation of seal boundaries through successive replacements &\hyperref[ss.rebuilt-seals-to-limelight]{\pageref*{ss.rebuilt-seals-to-limelight}}\\
reference-law multi-layer accrual & joint exponential tail for accrual across monitored layers &\hyperref[c.reference-multilayer-accrual]{\pageref*{c.reference-multilayer-accrual}}\\
towers, shadows, sunny targets & layers cast shadows on targets when positive-depth fine conditions fail \fff{towers-shadows}
\(\kappa_j;\ \mathcal G_{j-},\mathcal G_j\) & attempted index of the active slide and its pre- and post-slide fields &\hyperref[ss.slide-index-time-queue]{\pageref*{ss.slide-index-time-queue}}\\
\(\mathfrak b_{j,\rho};\ \mathsf U_j,\beta_j\) & labelled breaker unit, and two-row non-entry event and indicator &\hyperref[ss.usable-breaker-birth]{\pageref*{ss.usable-breaker-birth}}\\
\(J_j,\ \varrho_j\) & SimpleJoin queue-birth indicator and its renewal point &\hyperref[e.queue-start-probability]{\pageref*{e.queue-start-probability}}\\
\(\mathbf U_n^{(\sigma)},\ Q_n^{(\sigma)}\) & ordered protection queue of single-use breaker units and its current size \fff{protection-queue}
\(D_n^{(\sigma)},\ \tau_\sigma\) & conservative queue loss and abandonment time \fff{protection-queue}
\end{tabular}
\end{center}

\Needspace{20\baselineskip}
\bigskip
\noindent\textit{Chapters~14--15: Ornstein--Zernike consequences via classical renewal theory}
\begin{center}
\footnotesize
\renewcommand{\arraystretch}{1.12}
\begin{tabular}{@{}p{0.29\textwidth}p{0.61\textwidth}r@{}}
\textit{Symbol or term}&\textit{Meaning}&\textit{p.}\\ \hline
sausage & string with no renewal points other than its endpoints \fff{sausage}
\(\mathsf O(S)\) & occurrence event for a specified string or sausage \(S\) \fff{occurrence-event}
\(\zeta,\ \mu,\ \sausageset\) & string decay rate, tilted sausage law, and origin-based sausages \fff{sausage-law}
\(\mu^{\otimes\N},\ \rlevel,\ \mathsf M_{\rm str}\) & independent-sausage chain, renewal levels, and terminal-height mass \fff{renewal-chain}
\(\culdesac,\ \icds\) & lower and inverted endpoint cul-de-sac events \fff{culdesac}
\(T_n,\ \eta^-,\eta^+\) & renewal levels and tilted endpoint measures \fff{renewal-measure}
\(U\) & renewal measure, including the level at zero \fff{renewal-chain}
\(\mathcal D,\ \mathcal W,\ K_\xi,\ \xi\) & moment domain, its closure, and polar support-function pair \fff{moment-body}
\(X(S),\ \hgt(S),\ \pi(S),\operatorname{tsp}(S)\) & increment, height, transverse displacement, and transverse span \fff{vector-increment}
\(\Lambda_{\mathrm{vec}},\ \mathcal W_{\mathrm{str}}\) & vector sausage generating function and its unit sublevel set &\hyperref[e.lambda-vec]{\pageref*{e.lambda-vec}}\\
\(\mathcal N,\ t_{\mathrm{str}},\ \xi_{\mathrm{str}}\) & target cone, conjugate tilt, and string inverse correlation length &\hyperref[l.conjugate-tilt-interior]{\pageref*{l.conjugate-tilt-interior}}\\
\(\Psi^{\mathrm{str}}\) & positive analytic string Ornstein--Zernike amplitude &\hyperref[e.oz-string-llt]{\pageref*{e.oz-string-llt}}\\
\(W^{\base},W^{\head}\) & base- and head-rooted endpoint weights &\hyperref[l.jell.cds]{\pageref*{l.jell.cds}}\\
\(0\fconnext_{\vecell_0}x\) & \(0\) and \(x\) are the cluster base and head in the fixed \(\vecell_0\)-frame &\hyperref[t.oz-rooted]{\pageref*{t.oz-rooted}}\\
\(B^{\base},B^{\head}\) & marked boundary kernels for finite connectivity &\hyperref[e.marked-boundary-kernels]{\pageref*{e.marked-boundary-kernels}}\\
\(H^{\base},H^{\head}\) & weighted one-sided endpoint occupation factors &\hyperref[e.environment-amplitude]{\pageref*{e.environment-amplitude}}\\
\end{tabular}
\end{center}

\clearpage
\backmatter

\end{document}